# Categorical Reconstruction Theory

## Dissertation

zur Erlangung des akademischen Grades
Doctor rerum naturalium (Dr. rer. nat.)

vorgelegt
dem Bereich Mathematik und Naturwissenschaften
der Technischen Universität Dresden
von

Tony Zorman, M. Sc.

am 01. April 2025

begutachtet von
Prof. Dr. rer. nat. Ulrich Krähmer
Prof. Dr. rer. nat. Catharina Stroppel

angefertigt von November 2021 bis September 2025
am Institut für Geometrie

# INFORMATION ABOUT THE ARXIV VERSION

Due to compilation time constraints, this is a PDF only version of the thesis. The source code is nevertheless available at the following URL:

codeberg.org/slotThe/dissertation/src/branch/arXiv

Für Isabell
und meine Eltern

I almost wish I hadn't gone down that rabbit-hole—and yet—and yet—it's rather curious, you know, this sort of life!

Lewis Caroll, Alice in Wonderland

Maybe I'll start with the acknowledgements...

MIKE SLACKENERNY; PhD Comics #149

# ACKNOWLEDGEMENTS

It is surely impossible to write a dissertation without coming into contact with—and indeed becoming indebted to—a large number of wonderful people. This non-exhaustive list must certainly include

Uli, for being the best advisor, teacher, and life coach one could hope for—thank you for tolerating all of my stupid questions;

Isabell, Mum, and Dad, for giving me strength throughout, and for being patient with me when I was thinking too much about maths;

Catharina, for agreeing to be the second examiner of this thesis, and for organising the never-ending Bonn–Dresden seminar;

Sebastian, for not only putting up with me as an office mate, but also as my first coauthor—I'm not sorry for having inflicted Emacs upon you;

Matti, for being the best company in the world, and for talking to me about Hopf monads, even though you would much rather think about socles;

Zbiggi, for being a patient teacher of all things representation theory, and for not only making our office 85% more flammable, but also 100% more fun;

Benny, Florian, Julius, and Lukas, for squatting in our office, asking—and answering—lots of questions about maths and life;

Anna, for learning everything about duoidal categories in two hours flat;

Marcel and Sven, for always being up for some good old brotherly rivalry;

Florian, for encouraging me to start writing up early, because "even just streamlining notation takes more time than you think"—you were right;

Ivan and Myriam, for being the grownups in the office in times of need;

Philip, for—being knowledgeable about approximately everything—many enlightening conversation about mathematics, philosophy, music, movies, typesetting, and even Emacs; and

Yue, for—even now—never failing to make me smile.

I gracefully acknowledge the financial support of the Deutsche Forschungsgemeinschaft, who supported me from October 2022 to September 2025 via the grant KR 5036/2–1 "Cocommutative Comonoids".

# CONTENTS

That depends on where you want to get to;
"Well it really doesn't matter—";
Then it doesn't really matter which way you go.

Lewis Caroll, Alice in Wonderland

# INTRODUCTION

1

This dissertation generalises several reconstruction results in classical algebra to the language of monoidal categories, module categories, and monads. We start with some brief historical remarks, and summarise our main contributions in Section 1.1 below.

Monoidal and module categories are ubiquitous in many fields of mathematics. Briefly, a monoidal category is a category $\mathscr{C}$ equipped with a tensor product functor $\otimes\colon \mathscr{C}\times\mathscr{C} \longrightarrow \mathscr{C}$ and a unit $1 \in \mathscr{C}$, satisfying associativity and unitality conditions. A (left) module category $\mathscr{M}$ over $\mathscr{C}$, in turn, is endowed with an associative and unital action functor $\triangleright\colon \mathscr{C}\times\mathscr{M} \longrightarrow \mathscr{M}$. One motivating example of these structures comes from linear algebra. The category $\mathsf{Vect}$ of vector spaces over a field $\Bbbk$ is monoidal, and the right modules over a $\Bbbk$-algebra $A$ form a left $\mathsf{Vect}$-module category. Given $V \in \mathsf{Vect}$ and a right $A$-module $M$, the tensor product $V \otimes_{\Bbbk} M$ is again a right $A$-module, with action given by $(v \otimes m).a \coloneqq v \otimes m.a$, for all $a \in A$, $m \in M$, and $v \in V$.

The areas in which the language and theory of monoidal and module categories has been applied include quantum field theories [FRS02; Day07], categorification in representation theory [Str23; LMGRSW24; SW24], algebraic geometry [BZFN10; BZBJ18; Pas24], and many aspects of the theory of Hopf algebras and tensor categories [Sch00a; KK14; EGNO15; FGJS22; Shi23b; Str24b]. For a survey on the role of module categories in applied category theory and computer science, better known as actegories therein, see [CG22].

Just as in the category of vector spaces, given an algebra object $A \in \mathscr{C}$, the category $\mathrm{mod}_{\mathscr{C}}(A)$ of right $A$-modules in $\mathscr{C}$ is naturally a left $\mathscr{C}$-module category. Conversely, when considering a given $\mathscr{C}$-module category $\mathscr{M}$, it is often useful to find and study an algebra object $A$ in $\mathscr{C}$ such that there is an equivalence $\mathscr{M} \simeq \mathrm{mod}_{\mathscr{C}}(A)$ of $\mathscr{C}$-module categories. This kind of *reconstruction* process lies at the heart of Chapters 8 and 9 of the present thesis. An early result of this kind is [Ost03, Theorem 1], in which a "nice" module category over a finite tensor category $\mathscr{C}$—one that is equivalent to the finite-dimensional modules over a finite-dimensional algebra—is expressed as an algebra object in $\mathscr{C}$, see Theorem A of Chapter 8. Many generalisations and variants have appeared since; for instance [EGNO15, Theorem 7.10.1], [DSPS19, Theorem 2.24], [MMMT19, Theorem 4.7], and [BZBJ18, Theorem 4.6].

All of the reconstruction results cited above involve rigid monoidal categories: those in which every object has a dual. For example, to every finite-dimensional vector space $V$, there exist canonical evaluation and coevaluation maps $\mathrm{ev}\colon V^* \otimes_{\Bbbk} V \longrightarrow \Bbbk$ and $\mathrm{coev}\colon \Bbbk \longrightarrow V \otimes_{\Bbbk} V^*$, where $V^* := \mathrm{Hom}_{\Bbbk}(V, \Bbbk)$ is the dual space. However, in [DSPS19, Example 2.20] it is shown that, in the absence of rigidity, there exist $\mathscr{C}$-modules categories from which it is impossible to reconstruct an algebra object in $\mathscr{C}$. This forces us to consider a further generalisation—monads.

The term "monad" first appears in print in Section 5.4 of [Bén67]. Other names include "triple", "standard construction", or "triad".

Historically, the study of (co)monads arose out of (co)homology and homotopy theory, where one can use comonads to construct simplicial resolutions, see [God58; Hub61]. Already noted in [Hub61] is the deep connection of monads with the theory of adjunctions: every adjunction gives rise to a monad and, in turn, every monad arises from an adjunction. In fact, there exists a category of all adjunctions that realise a given monad. The initial and terminal object therein—the Kleisli and Eilenberg–Moore category of the monad, [EM65; Kle65]—play a pivotal role in this work.

Further applications are manifold and include category theory [Bec69; Str72; BKP89; LS02; AM24], categorical representation theory [Moe02; BLV11; Str24b; BHV24], Hopf algebras and quantum groups [BV07; AC12; KKS15; HL18], universal algebra [Lin66; Wal70; HP07], formal semantics [Mog89; Mog91; TP97], functional programming [Wad92; Wad93; PJW93], and algebraic effects [BHM02; PP02; PP09; BP15].

For this thesis, the most important aspect of monads lies in monadic reconstruction theory—relating additional structure on a monad $T$ on $\mathscr{C}$ to structure on its Eilenberg–Moore category $\mathscr{C}^T$. The connection to the theory of adjunctions yields a canonical forgetful functor $U^T : \mathscr{C}^T \longrightarrow \mathscr{C}$. This relationship can thus be seen as a generalisation of *Tannaka–Krein duality* [Tan38; Kre49], in which one reconstructs a compact group from its category of representations. Generalisations of this type of result have been proposed in a variety of different contexts, see [SR72; Wor88; Del90; Lur04; Szl09; Sch13].

A crucial feature of all of these statements is the involvement of a *fibre functor*: a forgetful strict monoidal functor to a nice underlying category, like the category of vector spaces or bimodules over a ring.

In the monadic case, Tannaka–Krein-type reconstruction results were obtained for bimonads [Moe02; McC02], Hopf monads [BV07], as well as linearly distributive and $*$-autonomous monads [PS09; Pas12]. We enrich this already broad landscape of results with new insights from the module category-theoretic perspective, offering a new reconstruction result for comodule monads over bimonads in the spirit of the algebraic version of considering coideal subalgebras over bialgebras.

## 1.1 SUMMARY

We shall now provide a short outline; a more detailed description and introduction will be provided at the start of each respective chapter. All original contributions of this thesis are already contained in either of the following articles or preprints.

The order in which the articles are listed here is that in which they were uploaded to the arXiv.

- [HZ24b]: Sebastian Halbig and Tony Zorman. Pivotality, twisted centres, and the anti-double of a Hopf monad. In: Theory Appl. Categ. 41 (2024), pp. 86–149. issn: 1201-561x.
- [HZ24a]: Sebastian Halbig and Tony Zorman. Diagrammatics for Comodule Monads. In: Appl. Categ. Struct. 32 (2024). Id/No 27, p. 17. issn: 0927-2852. doi: 10.1017/cbo9781139542333.
- [HZ23]: Sebastian Halbig and Tony Zorman. Duality in Monoidal Categories; 2023. arXiv: 2301.03545.
- [SZ24]: Mateusz Stroiński and Tony Zorman. Reconstruction of module categories in the infinite and non-rigid settings; 2024. arXiv: 2409.00793.
- [Zor25]: Tony Zorman. Duoidal R-Matrices; 2025. arXiv: 2503.03445.

To keep a common thematic focal point, the article [CSZ25]—although created during the author's PhD studies—will not be part of this work.

In the interest of brevity, we shall refrain from citing either of the articles [HZ23; HZ24a; HZ24b; SZ24; Zor25] beyond this introduction.

*Chapter 2: Preliminaries*

Being an amalgamate of the preliminaries of all articles, this chapter fixes notation and recalls basic properties of, for example, 2-categories, monoidal and module categories, the Drinfeld centre, linear and abelian categories, coends, as well as monads. We furthermore introduce one of the main tools for computation employed throughout—the graphical calculus.

*Chapter 3: Duality theory for monoidal categories*

Being mainly based on [HZ23], in this chapter we study various dualities for monoidal categories and how they interact: rigidity, tensor representability, Grothendieck–Verdier duality, and closedness. These notions form a hierarchy, with closedness—the least restrictive type of duality—on one end, and rigidity on the other. However, it is not obvious whether all inclusions are strict. Taken together, Propositions 3.7, 3.12, and 3.16; Examples 3.15 and 3.17; and Theorem 3.23 yield

$$\text{Rigid} \subsetneq \text{Tensor representable} \subsetneq \text{Grothendieck–Verdier} \subsetneq \text{Closed},$$

answering a question of Heunen.

We then more closely investigate duality structures on finite-dimensional functor categories endowed with Day convolution as its tensor product. Intuitively, this generalises the tensor product of modules over a commutative algebra, see Example 2.123. Given a nice base category, closedness and Grothendieck–Verdier duality lift to this setting.

**Proposition 3.27.** *Let $\mathscr{C}$ be a $\Bbbk$-linear hom-finite Grothendieck–Verdier category. Under the assumptions of Hypothesis 3.25, the finite-dimensional functor category $[\mathscr{C}, \mathsf{vect}]$ is a Grothendieck–Verdier category.*

Being motivated by representation theoretic questions, the Cauchy completion of a category also plays an important role for us. For example, the full subcategory of projective modules over a ring is the Cauchy completion of the full subcategory of all free modules. We find that the Cauchy completion completely mirrors the duality behaviour of its base category.

**Corollary 3.43.** *Let $\mathscr{C}^{\mathrm{op}}$ be a $\Bbbk$-linear right closed monoidal category. Then $\mathscr{C}^{\mathrm{op}}$ is right rigid (tensor representable) if and only if its Cauchy completion $\overline{\mathscr{C}^{\mathrm{op}}}$ is right rigid (tensor representable). Further, $(\mathscr{C}^{\mathrm{op}}, d)$ is a right Grothendieck–Verdier category if and only if $(\overline{\mathscr{C}^{\mathrm{op}}}, \text{よ}_d)$ is.*

*Chapter 4: Twisted centres*

In this chapter, mainly based on [HZ24b], we hoist Theorem 4.1 from the Hopf-theoretic world into that of monoidal categories. It lays the groundwork for later monadic generalisations.

**Theorem 4.23.** *Let $\mathcal{C}$ be a rigid category. There is a bijection between*

(i) *equivalence classes of quasi-pivotal structures on $\mathcal{C}$,*
(ii) *the Picard heap of the anti-centre, and*
(iii) *isomorphism classes of equivalences of module categories between the centre and the anti-centre.*

We then more closely investigate the heap structure of the anti-centre, answering a question of Shimizu about the relationship to quasi-pivotal structures on the underlying category.

**Theorem 4.50.** *There exists a category $\mathcal{C}$ on which there exists a pivotal structure that is not induced by any element of the Picard heap of the anti-centre of $\mathcal{C}$.*

*Chapter 5: Monadic Tannaka–Krein reconstruction*

In this chapter, which is mainly based on [HZ24a; HZ24b; SZ24], we study comodule monads in the sense of [AC12]. These structures can be seen as generalising comodule algebras over a bialgebra, see Example 5.14. We then prove our main Tannaka–Krein reconstruction result for comodule monads in the spirit of [Moe02, Theorem 7.1] and [McC02, Corollary 3.13].

**Theorem 5.31.** *Let $B$ be a bimonad on the monoidal category $\mathcal{C}$, and $T$ a monad on a right $\mathcal{C}$-module category $\mathcal{M}$. Coactions of $B$ on $T$ are in bijection with right actions of $\mathcal{C}^B$ on $\mathcal{M}^T$ such that $U^T$ is a strict comodule functor over $U^B$.*

In particular, we prove statements akin to Kelly's doctrinal adjunctions result in the case of comodule and $\mathcal{C}$-module monads—the latter play an important role in Chapters 8 and 9.

**Theorem 5.28 and Porism 5.29.** *Let $F\colon \mathcal{C} \rightleftarrows \mathcal{D} :U$ be an oplax monoidal adjunction. Lifts of an ordinary adjunction $G\colon \mathcal{M} \rightleftarrows \mathcal{N} :V$ to a comodule adjunction are in bijection with lifts of $V\colon \mathcal{N} \longrightarrow \mathcal{M}$ to a strong comodule functor.*

*An adjunction $F\colon \mathcal{M} \rightleftarrows \mathcal{N} :U$ between $\mathcal{C}$-module categories yields a bijection of oplax $\mathcal{C}$-module structures on $F$ and lax $\mathcal{C}$-module structures on $U$.*

*Chapter 6: Monadic twisted centres*

This chapter, mainly based on [HZ24b], provides a monadic interpretation of Chapter 4. Based on the Drinfeld double of a Hopf monad, we develop the notion of an anti-double, which realises the anti-centre of *ibid* as its Eilenberg–Moore category, and present a monadic variant of Theorem 4.23.

**Theorem 6.44.** *Let $\mathscr{C}$ be a rigid monoidal category. For a Hopf monad $H$ on $\mathscr{C}$ that admits a double and anti-double, the following statements are equivalent:*

(*i*) *the monoidal unit of $\mathscr{C}$ lifts to a module over the anti-double,*
(*ii*) *the double and anti-double of $H$ are isomorphic as comodule monads, and*
(*iii*) *the double and anti-double of $H$ are isomorphic as monads.*

*If $\mathscr{C}$ is pivotal, the above statements hold if and only if $H$ admits a pair in involution.*

Investigating the interplay between the double and anti-double, we obtain a new criterion for when a rigid monoidal category is pivotal.

**Corollary 6.45.** *Let $\mathscr{C}$ be a rigid monoidal category. If $\mathscr{C}$ admits a central Hopf monad $\mathfrak{D}(\mathscr{C})$ and an anti-central comodule monad $\mathfrak{Q}(\mathscr{C})$, then it is pivotal if and only if $\mathfrak{D}(\mathscr{C}) \cong \mathfrak{Q}(\mathscr{C})$ as monads.*

*Chapter 7: Duoidal R-matrices*

Being mainly based on [Zor25], this chapter introduces R-matrices in duoidal categories, generalising the well-known R-matrices for bialgebras, and those for bimonads [BV07]. As it turns out, R-matrices on a suitable monad correspond to duoidal structures its Eilenberg–Moore category.

**Theorem 7.21.** *Let $\mathfrak{D}$ be a category with monoidal structures $\circ$ and $\bullet$, and $T$ a monad on $\mathfrak{D}$ that has a $\circ$-oplax monoidal and a $\bullet$-oplax monoidal structure. Then quasitriangular structures on $T$ are in bijection with duoidal structures on $\mathfrak{D}^T$.*

We then explore the connection between normal duoidal categories, and (non-planar) linearly distributive categories [CS97; GLF16]. The latter can be seen as an analogue of Grothendieck–Verdier categories without an explicit notion of dual. In Proposition 7.31, we relate a cocommutative version of duoidal monads to the linearly distributive monads of [Pas12].

*Chapter 8: Infinite and non-rigid reconstruction*

This chapter, mainly based on [SZ24], studies reconstruction without a fibre functor, as done for example in [Ost03, Theorem 1].

Since every object $m$ in a module category $\mathcal{M}$ over $\mathcal{C}$ gives rise to a strong $\mathcal{C}$-module functor $- \triangleright m \colon \mathcal{C} \longrightarrow \mathcal{M}$, and in good cases this functor has a right adjoint, the resulting monad is naturally a lax $\mathcal{C}$-module monad. However, it is a priori not clear whether the Eilenberg–Moore category of a (right exact) lax $\mathcal{C}$-module monad even has a canonical $\mathcal{C}$-module structure. Using multicategorical techniques, we establish when such an *extension* exists.

**Theorems 8.25 and 8.26.** *Let $\mathcal{C}$ be an abelian monoidal and $\mathcal{M}$ an abelian $\mathcal{C}$-module category. Suppose that the action functor $\triangleright \colon \mathcal{C} \otimes_{\Bbbk} \mathcal{M} \longrightarrow \mathcal{M}$ is right exact in both variables, and let $T \colon \mathcal{M} \longrightarrow \mathcal{M}$ be a right exact lax $\mathcal{C}$-module monad. Then there exists an essentially uniquely $\mathcal{C}$-module structure on $\mathcal{M}^T$, such that the canonical inclusion $\iota \colon \mathcal{M}_T \longrightarrow \mathcal{M}^T$ is a strong $\mathcal{C}$-module functor.*

Using the constructed $\mathcal{C}$-module structure, we then establish a reconstruction and classification result under mild additional assumptions.

**Theorems 8.48 to 8.50.** *Assume that $\mathcal{C}$ and $\mathcal{M}$ have enough projectives (injectives) and that there exists an object $\ell \in \mathcal{M}$ such that:*

- *there is a right adjoint $\lfloor \ell, - \rfloor$ (left adjoint $\lceil \ell, - \rceil$) to $- \triangleright \ell$;*
- *for $x \in \mathcal{C}$ projective (injective), the object $x \triangleright \ell$ is projective (injective);*
- *any projective (injective) object of $\mathcal{M}$ is a direct summand of an object of the form $x \triangleright \ell$, for $x$ projective (injective).*

*Let $T$ be the monad $\lfloor \ell, - \triangleright \ell \rfloor$. Then there is an equivalence $\mathcal{M} \simeq \mathcal{C}^T$ of $\mathcal{C}$-module categories, where the $\mathcal{C}$-module structure of the category of $T$-modules is extended from the Kleisli category. Furthermore, this gives rise to a bijection*

$$\{\, (\mathcal{M}, \ell) \text{ as above} \,\} \big/ (\mathcal{M} \simeq \mathcal{N}) \;\leftrightarrow\; \left\{ \begin{array}{l} \text{Right exact lax } \mathcal{C}\text{-module} \\ \text{monads on } \mathcal{C} \end{array} \right\} \big/ \, (\mathcal{C}^T \simeq \mathcal{C}^S)$$

$$(\mathcal{M}, \ell) \longmapsto \lfloor \ell, - \triangleright \ell \rfloor$$

$$(\mathcal{C}^T, T1) \longleftarrow T$$

*Chapter 9: Hopf trimodules*

The final chapter, mainly based on [SZ24], further develops the non-rigid reconstruction results of Chapter 8 in the case of several examples.

For a Hopf algebra $B$, the *Fundamental Theorem of Hopf Modules* says that the category ${}^{B}_{B}\mathsf{Vect}$ of Hopf modules is naturally equivalent to the category of vector spaces [LS69]. Likewise, the category of *Hopf trimodules* ${}^{B}_{B}\mathsf{Vect}^{B}$ is equivalent to that of right comodules. In comparison to plain Hopf modules, one can generalise Hopf trimodules to the quasi-Hopf algebraic world.

If $B$ is just a bialgebra—i.e., does not have an antipode—then we can still extract interesting structure from such modules.

**Theorem 9.2.** *For a bialgebra $B$ and $\mathcal{V} := {}^{B}\mathsf{Vect}$, there is a monoidal equivalence*

$${}^{B}_{B}\mathsf{Vect}^{B} \simeq \mathsf{LexfLax}\mathcal{V}\mathsf{Mod}(\mathcal{V}, \mathcal{V})$$

*between the category of Hopf trimodules, and the category of left exact finitary lax $\mathcal{V}$-module endofunctors on $\mathcal{V}$.*

This allows us to give a categorical interpretation of the aforementioned Fundamental Theorem of Hopf Modules.

**Proposition 9.38 and Corollary 9.45.** *The functor ${}^{B}\mathsf{Vect} \longrightarrow {}^{B}_{B}\mathsf{Vect}^{B}$ corresponds to the inclusion of strong ${}^{B}\mathsf{Vect}$-module endofunctors into the category of lax ${}^{B}\mathsf{Vect}$-module endofunctors. It is an equivalence if and only if ${}^{B}\mathsf{vect}$ is left rigid, which is the case if and only if $B$ has a twisted antipode.*

We also give applications in the setting of semigroup algebras, see Sections 9.4 and 9.5. This can be seen as an example of a non-rigid category where the reconstruction procedure described in [EGNO15, Chapter 7] fails.

Lastly, the fusion operators of a bimonad may be expressed in terms of its canonical module structure.

**Lemma 9.47 and Proposition 9.49.** *Let $F\colon \mathcal{C} \rightleftarrows \mathcal{D} :U$ be an oplax monoidal adjunction. The strong monoidal structure of $U$ turns $\mathcal{C}$ into a $\mathcal{D}$-module category, by defining $- \triangleright = := U(-) \otimes =$. With respect to this $\mathcal{D}$-module structure, the bimonad $T := UF$ on $\mathcal{C}$ becomes an oplax $\mathcal{D}$-module monad. The right fusion operator is given by evaluating the $\mathcal{D}$-module coherence morphism at a free $T$-module, and furthermore $T$ is right Hopf if and only if it is a strong $\mathcal{C}^{T}$-module monad.*

Das Rad gehört nicht zur Maschine,
das man drehen kann, ohne daß
etwas anderes sich mitbewegt.

LUDWIG WITTGENSTEIN;
Philosophische Untersuchungen

# PRELIMINARIES

# 2

WE ASSUME THE READER'S FAMILIARITY with basic category theoretical concepts, as well as the main features of the theories of monoidal and 2-categories, as discussed for example in [ML98; EGNO15; Rie17; JY21]. Nevertheless, we start with a brief review of the most important terminology and notation.

## 2.1 BICATEGORIES

THE THEORY OF BICATEGORIES IS ESSENTIAL to the graphical calculus we will use throughout, see Section 2.3, as well as to some formal arguments we shall employ. Following [JY21], we recall the most important definitions.

**Definition 2.1.** A *bicategory* $\mathbb{B}$ has as data:

- A collection of *objects* $\operatorname{Ob}\mathbb{B}$, where we often write $x \in \mathbb{B}$ for $x \in \operatorname{Ob}\mathbb{B}$;
- for all $x, y \in \mathbb{B}$ a *hom-category* $\mathbb{B}(x, y)$, whose objects are called *1-cells* (or 1-morphisms) and whose morphisms are called *2-cells* (or 2-morphisms);
- for all $x, y, z \in \mathbb{B}$ a *horizontal composition* functor

$$\otimes_y : \mathbb{B}(y, z) \times \mathbb{B}(x, y) \longrightarrow \mathbb{B}(x, z);$$

- for all $x \in \mathbb{B}$ an *identity* $1_x \in \mathbb{B}(x, x)$;
- a natural isomorphism $\alpha_{M,N,P} : (M \otimes_y N) \otimes_x P \longrightarrow M \otimes_y (N \otimes_x P)$;
- natural isomorphisms $\lambda_M : 1_y \otimes_y M \longrightarrow M$ and $\rho_M : M \longrightarrow M \otimes_x 1_x$.

This data is subject to the commutativity of the following diagrams, for all suitable 1-cells $M, N, P, R$:

$$\begin{array}{ccc}
(M \otimes_z (N \otimes_y P)) \otimes_x R & \xrightarrow{\alpha_{M,N\otimes_y P,R}} & M \otimes_z ((N \otimes_y P) \otimes_x R) \\
\uparrow{\scriptstyle \alpha_{MNP}\otimes R} & & \downarrow{\scriptstyle M\otimes\alpha_{NPR}} \\
((M \otimes_z N) \otimes_y P) \otimes_x R & & M \otimes_z (N \otimes_y (P \otimes_x R)) \\
& \searrow{\scriptstyle \alpha_{M\otimes_z N,P,R}} \quad \nearrow{\scriptstyle \alpha_{M,N,P\otimes_x R}} & \\
& (M \otimes_z N) \otimes_y (P \otimes_x R) &
\end{array}$$

$$\begin{array}{ccc}
(M \otimes_x 1_x) \otimes_x N & \xrightarrow{\alpha_{M,1_x,N}} & M \otimes_x (1_x \otimes_x N) \\
\uparrow{\scriptstyle \rho_M\otimes N} & & \downarrow{\scriptstyle M\otimes\lambda_N} \\
M \otimes_x N & = & M \otimes_x N
\end{array}$$

If $\alpha$, $\lambda$, and $\rho$ are all identities, we call $\mathbb{B}$ a *2-category*.

The composition of morphisms inside of the category $\mathbb{B}(x, y)$ will be called the *vertical composition*. We denote the terminal bicategory consisting of a single object, a single 1-, and a single 2-morphism by $\heartsuit$.

**Example 2.2.** The category $\mathbb{C}$at of (small) categories, functors, and natural transformations is a 2-category. Horizontal composition is given by horizontal composition of functors and natural transformations.

Examples of bicategories are plentiful throughout this work. We refrain from spelling them out here before having introduced more of the necessary notation, but see Examples 2.103, 2.111, and 2.113 and Remarks 4.8 and 5.2

**Remark 2.3.** The use of $\otimes$ for the horizontal composition is non-standard, and inspired by [Gar22]. We use it here to highlight the parallels to the definition of monoidal categories, see Definition 2.31.

**Notation 2.4.** To syntactically differentiate bicategories from ordinary categories, we will generally start the name of a bicategory with a blackboard bold letter. For example, when talking about the 1-category of all (small) categories and functors, we write Cat instead of $\mathbb{C}$at.

Much like monoidal categories, see Theorem 2.40, bicategories admit a *coherence* and *strictification* result, in that every bicategory is biequivalent to a 2-category, see for example [JY21, Theorem 8.4.1 and Corollary 8.4.2].

**Definition 2.5.** A *lax functor* $F\colon \mathbb{B} \longrightarrow \mathbb{C}$ of bicategories $\mathbb{B}$, $\mathbb{C}$ has as data:

- an object assignment $F\colon \operatorname{Ob}\mathbb{B} \longrightarrow \operatorname{Ob}\mathbb{C}$;
- for all $x, y \in \mathbb{B}$, a functor $F_{x,y}\colon \mathbb{B}(x,y) \longrightarrow \mathbb{C}(Fx, Fy)$;
- a family of natural transformations $F_2\colon F\otimes F \Longrightarrow F\circ\otimes$, with components $F_{2,M,N}\colon F_{yz}M \otimes_{Fy} F_{xy}N \longrightarrow F_{xz}(M\otimes_y N)$, for all $x, y, z \in \mathbb{B}$, $M \in \mathbb{B}(y,z)$, and $N \in \mathbb{B}(x,y)$; and
- for all $x \in \mathbb{B}$, arrows $F_{0,x}\colon 1_{Fx} \longrightarrow F1_x$.

This data is subject to the following axioms in $\mathbb{C}(Fw, Fz)$ and $\mathbb{C}(Fx, Fy)$, respectively, for all admissible $M, N, P$:

$$\begin{array}{ccc}
(FM \otimes_{Fy} FN) \otimes_{Fx} FP & \xrightarrow{\alpha} & FM \otimes_{Fy} (FN \otimes_{Fx} FP) \\
\downarrow{\scriptstyle F_{2;M,N}\otimes \mathrm{id}} & & \downarrow{\scriptstyle \mathrm{id}\otimes F_{2;N,P}} \\
F(M \otimes_y N) \otimes_{Fx} FP & & FM \otimes_{Fy} F(N \otimes_x P) \\
\downarrow{\scriptstyle F_{2;M\otimes N,P}} & & \downarrow{\scriptstyle F_{2;M,N\otimes P}} \\
F((M \otimes_y N) \otimes_x P) & \xrightarrow[F\alpha]{} & F(M \otimes_y (N \otimes_x P))
\end{array}$$

$$\begin{array}{ccccc}
1_{Fy} \otimes_{Fy} FM & \xrightarrow{\lambda} & FM & \xrightarrow{\rho} & FM \otimes_{Fx} 1_{Fx} \\
\downarrow{\scriptstyle F_{0;y}\otimes \mathrm{id}} & & \Vert & & \downarrow{\scriptstyle \mathrm{id}\otimes F_{0;x}} \\
F1_y \otimes_{Fy} FM & & \Vert & & FM \otimes_{Fx} F1_x \\
\downarrow{\scriptstyle F_{2;1_y,M}} & & \Vert & & \downarrow{\scriptstyle F_{2;M,1_x}} \\
F(1_y \otimes_y M) & \xrightarrow[F\lambda]{} & FM & \xrightarrow[F\rho]{} & F(M \otimes_x 1_x)
\end{array}$$

A lax functor is called a *pseudofunctor* if all $F_2$'s and $F_0$'s are invertible, and a *2-functor* if they are all identities.

Analogously, one defines *oplax functors* between 2-categories.

**Definition 2.6.** An *oplax transformation* $\omega\colon F \Longrightarrow G$ between lax functors $(F, F_2, F_0), (G, G_2, G_0)\colon \mathbb{B} \longrightarrow \mathbb{C}$ consists of a 1-cell $\omega_x\colon Fx \longrightarrow Gx$ in $\mathbb{C}$ for every $x \in \mathbb{B}$, and for all $x, y \in \mathbb{B}$ a natural transformation

$$\omega\colon \omega_y \otimes F(-) \Longrightarrow G(-) \otimes \omega_x\colon \mathbb{B}(x,y) \longrightarrow \mathbb{C}(Fx, Gy)$$

with component 2-cells $\omega_f\colon \omega_y \otimes_{Fy} Ff \Longrightarrow Gf \otimes_{Gx} \omega_x$ for all $f\colon x \longrightarrow y$ in $\mathbb{B}$. Graphically, we draw these 2-cells like the following:

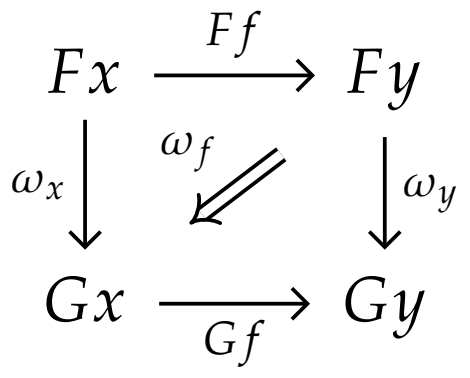

An oplax transformation is called a *pseudonatural transformation*[1] if all $\omega_f$'s are isomorphisms, and a *2-natural transformation* if they are all identities.

[1] Pseudonatural transformations are called *strong transformations* in [JY21].

Analogously, one defines *lax* transformations between lax functors, as well as lax and oplax transformations between oplax functors.

**Definition 2.7** ([Lac10])**.** Let $F, G\colon \mathbb{B} \longrightarrow \mathbb{C}$ be lax functors between bicategories. An *icon* between $F$ and $G$ consists of an assertion that $Fx = Gx$, for all $x \in \mathbb{B}$; and an oplax transformation $\omega\colon F \Longrightarrow G$, such that for all $x \in \mathbb{B}$, the 1-cell $\omega_x\colon Fx \longrightarrow Gx$ is the identity.

Definition 2.7 might seem a bit contrived at first, but it is an important concept: there is no 2-category with (small) bicategories as 0-cells, lax (even pseudo-) functors as 1-cells, and 2-natutral transformations as 2-cells. Instead, icons yield the "right" kind of 2-cells for this construction, see [JY21, Section 4.6] for a more extensive account.

**Definition 2.8.** Let $\mathbb{B}$ and $\mathbb{C}$ be bicategories, $F, G\colon \mathbb{B} \longrightarrow \mathbb{C}$ lax functors, and suppose that $\alpha, \beta\colon F \Longrightarrow G$ are oplax transformations. A *modification* $\Upsilon\colon \alpha \Rrightarrow \beta$ between $\alpha$ and $\beta$ consists of a 2-cell $\Upsilon_x\colon \alpha_x \Longrightarrow \beta_x$ for every $x \in \mathbb{B}$, such that the following diagram commutes for all 1-cells $f \in \mathbb{B}(x, y)$:

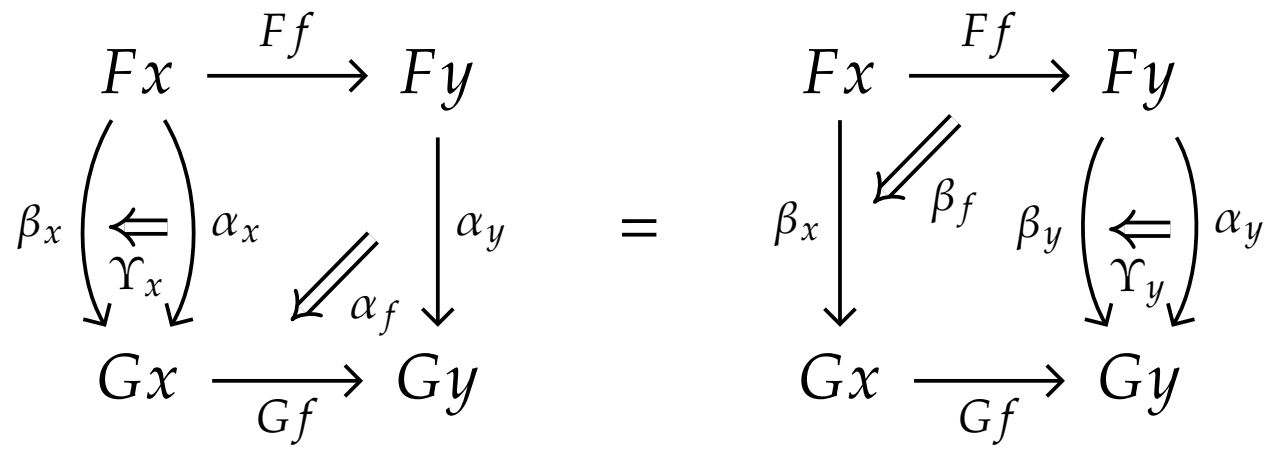

In fact, bicategories, lax functors, lax natural transformations, and modifications between them form a *tricategory*. We shall not give the formal definition here, but see [GPS95; Gur06] and [JY21, Chapter 11].

## 2.2 MONADS AND ADJUNCTIONS

**Definition 2.9.** A *monad* on a category $\mathscr{C}$ consists of an endofunctor $T$ on $\mathscr{C}$ and two natural transformations $\mu\colon T^2 \Longrightarrow T$ and $\eta\colon \mathrm{Id}_{\mathscr{C}} \Longrightarrow T$, called the *multiplication* and *unit* of $T$, satisfying *associativity* and *unitality* axioms:

$$
\begin{array}{ccc}
T^3 & \xrightarrow{T\mu} & T^2 \\
{\scriptstyle \mu T}\downarrow & & \downarrow{\scriptstyle \mu} \\
T^2 & \xrightarrow[\mu]{} & T
\end{array}
\qquad\qquad
\begin{array}{ccccc}
T^2 & \xleftarrow{\eta T} & T & \xrightarrow{T\eta} & T \\
& {\scriptstyle \mu}\searrow & \Big| & \swarrow{\scriptstyle \mu} & \\
& & T & &
\end{array}
$$

**Definition 2.10.** An *adjunction* consists of a pair of functors $F\colon \mathscr{C} \rightleftarrows \mathscr{D} :U$, together with two natural transformations, the *unit* $\eta\colon \mathrm{Id}_{\mathscr{C}} \Longrightarrow UF$ and the *counit* $\varepsilon\colon FU \Longrightarrow \mathrm{Id}_{\mathscr{D}}$, satisfying the *snake* or *triangle* identities; see for example [Rie17, Definition 4.2.5].

Alternatively, one could define an adjunction between $F\colon \mathscr{C} \longrightarrow \mathscr{D}$ and $G\colon \mathscr{D} \longrightarrow \mathscr{C}$ to require the existence of a natural isomorphism

$$\mathscr{D}(F(-), =) \cong \mathscr{C}(-, U(=)),$$

from which one recovers the unit and counit, see [Rie17, Proposition 4.2.6]. Note that, a priori, different natural isomorphisms lead to different adjunctions between the same functors. However, fixing for example $G$, then for two left adjoints $F$ and $F'$ there exists a unique isomorphism $\Theta\colon F \Longrightarrow F'$ that commutes with the respective units and counits; see [Rie17, Proposition 4.4.1].

**Example 2.11.** If $F\colon \mathscr{C} \rightleftarrows \mathscr{D} :U$ is an adjunction with unit $\eta$ and counit $\varepsilon$, then $UF\colon \mathscr{C} \longrightarrow \mathscr{C}$ is a monad; the multiplication $\mu^T\colon UFUF \Longrightarrow UF$ is given by $U\varepsilon F$, and the unit $\eta^T\colon \mathrm{Id}_{\mathscr{C}} \Longrightarrow UF$ is equal to $\eta$.

**Definition 2.12.** Given a monad $(T, \mu, \eta)$ on a category $\mathscr{C}$, a $T$-*algebra* consists of an object $x \in \mathscr{C}$ and a morphism $\alpha\colon Tx \longrightarrow x$, satisfying

$$\alpha \circ T\alpha = \alpha \circ \mu_x \qquad \text{and} \qquad \alpha \circ \eta_x = \mathrm{id}_x.$$

Given two $T$-algebras $(x, \alpha)$ and $(y, \beta)$, a *morphism* between them consists of a morphism $f\colon x \longrightarrow y$ in $\mathscr{C}$, such that $\beta \circ Tf = f \circ \alpha$.

**Remark 2.13.** For any monad $T$, the $T$-algebras and their morphisms form a category: the *Eilenberg–Moore category of* $T$. We shall denote it by $\mathscr{C}^T$.

Our notation for the Eilenberg–Moore category already ascribes a certain "terminal" quality to it; see Remark 2.23.

The Eilenberg–Moore category of $T$ is also often called the *category of* $T$*-algebras* or, following for example [BV07], the *category of modules* over $T$. We use all three terminologies interchangeably.

A monad is intimately connected to its Eilenberg–Moore category.

**Example 2.14.** There is a 2-category $\mathbb{M}\mathsf{on}(\mathbb{C}\mathsf{at})$ of monads in $\mathbb{C}\mathsf{at}$, [Str72, § 1]. The inclusion 2-functor maps a category to its identity monad:

$$\mathbb{C}\mathsf{at} \longrightarrow \mathbb{M}\mathsf{on}(\mathbb{C}\mathsf{at}), \qquad \mathscr{C} \longmapsto (\mathrm{Id}_{\mathscr{C}}, \mathrm{id}_{\mathrm{Id}_{\mathscr{C}}}, \mathrm{id}_{\mathrm{Id}_{\mathscr{C}}}).$$

By assumption, $\mathbb{C}\mathsf{at}$ *admits the construction of algebras*: there exists a right adjoint to the above functor:

$$\mathbb{M}\mathsf{on}(\mathbb{C}\mathsf{at}) \longrightarrow \mathbb{C}\mathsf{at}, \qquad (T\colon \mathscr{C} \longrightarrow \mathscr{C}, \mu, \eta) \longmapsto \mathscr{C}^T,$$

where $\mathscr{C}^T \in \mathbb{C}\mathsf{at}$ is the Eilenberg–Moore category of $T$. Using the previous 2-adjunction, one proves that to every monad $(T, \mu, \eta)$ on $\mathscr{C}$ there exist an *Eilenberg–Moore adjunction* $F^T\colon \mathscr{C} \longrightarrow \mathscr{C}^T$ and $U^T\colon \mathscr{C}^T \longrightarrow \mathscr{C}$, such that

$$T = U^T F^T, \qquad \mu = F^T \varepsilon U^T, \qquad \eta = \eta,$$

where $\eta\colon 1_{\mathscr{C}} \Longrightarrow U^T F^T$ and $\varepsilon\colon F^T U^T \Longrightarrow 1_{\mathscr{C}^T}$ are the unit and counit of the Eilenberg–Moore adjunction. We shall call $F^T$ and $U^T$ the *free* and *forgetful* functor *associated to* $T$, respectively.

For the following definition, we follow [BV07; TV17].

**Definition 2.15.** Suppose that $T$ and $S$ are two monads on the category $\mathscr{C}$. A *morphism of monads* between $T$ and $S$ is a natural transformation $f\colon T \Longrightarrow S$, such that the following diagrams commute

$$\begin{array}{ccc} \mathrm{Id} & \xrightarrow{\eta^T} & T \\ & \searrow^{\eta^S} & \downarrow f \\ & & S \end{array} \qquad\qquad \begin{array}{ccccc} TT & \xrightarrow{Tf} & TS & \xrightarrow{fS} & SS \\ \mu^T \downarrow & & & & \downarrow \mu^S \\ T & & \xrightarrow{\quad f \quad} & & S \end{array} \tag{2.2.1}$$

**Remark 2.16.** The terminology of Definition 2.15 is slightly non-standard. What we call a morphism of monads is often called a *oplax* (or colax) monad morphism, see for example [Str72, § 1].

**Remark 2.17.** One can define a monad in any bicategory $\mathbb{B}$ by considering $(C, t, \eta, \mu)$, where $C \in \mathbb{B}$ is an object, $t\colon C \longrightarrow C$ is a 1-cell, and $\eta\colon \mathrm{Id}_C \Longrightarrow t$ and $\mu\colon tt \Longrightarrow t$ are 2-cells, satisfying relations analogous to Definition 2.9. An oplax morphism of monads from $(C, t, \eta^t, \mu^t)$ to $(D, s, \eta^s, \mu^s)$ then consists of a 1-cell $u\colon C \longrightarrow D$ and a 2-cell $\phi\colon ut \Longrightarrow su$, subject to identities reminiscent of Diagram (2.2.1). A *lax morphism of monads* involves a 1-cell $u\colon C \longrightarrow D$ and a 2-cell $\phi\colon su \Longrightarrow ut$, satisfying similar properties.

The following example sheds some additional light on this terminology.

**Example 2.18.** Monads can alternatively be defined as lax functors—in the sense of Definition 2.5—from the terminal 2-category $\heartsuit$ to $\mathbb{C}\mathsf{at}$. Unravelling this definition, a lax functor $\mathfrak{T}\colon \heartsuit \longrightarrow \mathbb{C}\mathsf{at}$ consists of:

- an object assignment $\mathfrak{T}\colon \operatorname{Ob}\heartsuit \longrightarrow \operatorname{Ob}\mathbb{C}\mathsf{at}$, sending the unique object $*$ to a category $\mathscr{C}$;
- a functor $\mathfrak{T}(*,*)\colon \heartsuit(*,*) \longrightarrow \mathbb{C}\mathsf{at}(\mathscr{C},\mathscr{C})$ from the terminal 1-category $\heartsuit(*,*)$ to the category of endofunctors on $\mathscr{C}$, sending the unique 1-morphism $\mathrm{id}_*\colon * \longrightarrow *$ to $T\colon \mathscr{C} \longrightarrow \mathscr{C}$ and the unique 2-morphism $1_{\mathrm{id}_*}\colon \mathrm{id}_* \Longrightarrow \mathrm{id}_*$ to the identity natural transformation $T \Longrightarrow T$;
- a 2-cell $\mathfrak{T}_2\colon \mathfrak{T}\mathrm{id}_* \otimes \mathfrak{T}\mathrm{id}_* \Longrightarrow \mathfrak{T}\mathrm{id}_*$, which we write as $\mu\colon TT \Longrightarrow T$; and
- a 2-cell $\mathfrak{T}_0\colon 1_{\mathfrak{T}(*)} \Longrightarrow \mathfrak{T}\mathrm{id}_*$ that we write as $\eta\colon \mathrm{Id}_{\mathscr{C}} \Longrightarrow T$.

Recall that $\otimes$ is the horizontal composition in $\mathbb{B}$.

The properties of Definition 2.5 for $\mathfrak{T}_2$ and $\mathfrak{T}_0$ translate to the associativity and unitality properties of $\mu$ and $\eta$. In this setting, a morphism of monads becomes an oplax transformation in the sense of Definition 2.6.

**Example 2.19.** A monad $T$ on $\mathscr{C}$ has another canonical category associated to it: its *Kleisli category* $\mathscr{C}_T$. On objects, it is given by $\operatorname{Ob}(\mathscr{C}_T) := \operatorname{Ob}(\mathscr{C})$, and for $x, y \in \mathscr{C}_T$ we have $\mathscr{C}_T(x,y) := \mathscr{C}(x, Ty)$. Composition is defined by

$$\begin{aligned} \circ\colon \mathscr{C}_T(y,z) \times \mathscr{C}_T(x,y) &\longrightarrow \mathscr{C}_T(x,z) \\ (g,f) &\longmapsto \big(x \xrightarrow{f} Ty \xrightarrow{Tg} T^2 z \xrightarrow{\mu_z} Tz\big). \end{aligned}$$

**Proposition 2.20.** *Let $T$ be a monad on a category $\mathscr{C}$. There exists an adjunction*

$$\mathscr{C} \underset{U_T}{\overset{F_T}{\rightleftarrows}} \mathscr{C}_T \quad (F_T \dashv U_T)$$

*where $F_T$ is identity on objects and sends $f \in \mathscr{C}(x,y)$ to $\eta_y \circ f$, and $U_T$ sends $x$ to $Tx$ and $f \in \mathscr{C}(x,y)$ to $\mu_y \circ Tf$.*

### 2.2.1 *Comparison functors*

Given a monad $T$ coming from the adjunction $F\colon \mathscr{C} \rightleftarrows \mathscr{D} : U$, we might ask how much the functors $F$ and $U$ "differ" from the free and forgetful functors $F^T\colon \mathscr{C} \longrightarrow \mathscr{C}^T$ and $U^T\colon \mathscr{C}^T \longrightarrow \mathscr{C}$ of $T$. Roughly summarised we are interested in the following:

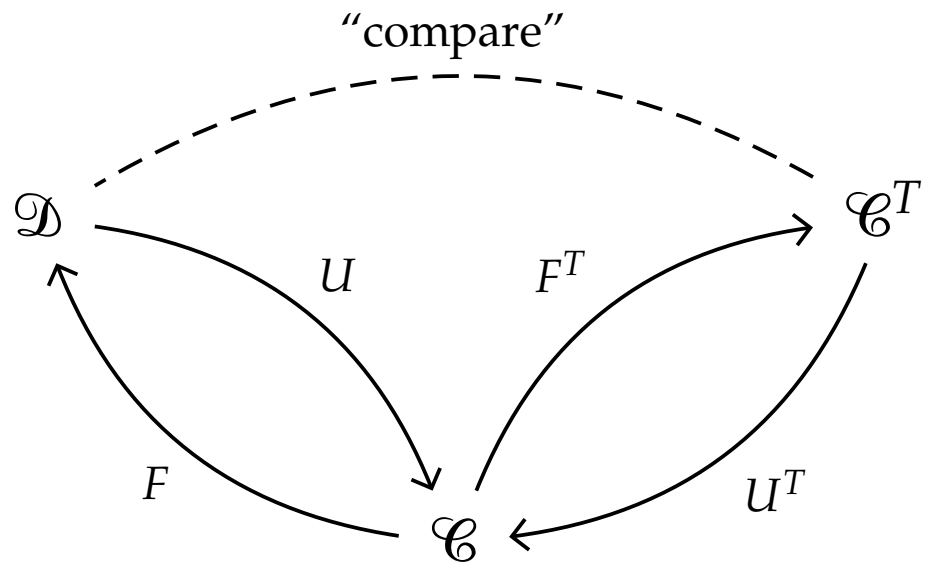

**Proposition 2.21** ([Fre69, Theorem 2])**.** *Let $T$ and $S$ be monads on a category $\mathscr{C}$. Then there exists a bijection between monad morphisms from $T$ to $S$ and functors $F\colon \mathscr{C}^S \longrightarrow \mathscr{C}^T$ between their categories of algebras, such that $U^T F = U^S$.*

Note that, if we were to define morphisms of monads as lax transformations, Proposition 2.21 would yield a covariant assignment.

**Lemma 2.22** ([Str72, Theorem 3])**.** *Let $F\colon \mathscr{C} \rightleftarrows \mathscr{D} : U$ be an adjunction and set $T := UF$. There exists a unique functor $K^T\colon \mathscr{D} \longrightarrow \mathscr{C}^T$ satisfying $K^T F = F^T$ and $U^T K^T = U$. On objects it is given by $K^T d = (Ud, U\varepsilon_d)$, for all $d \in \mathscr{D}$.*

We call the unique functor from Lemma 2.22 the *comparison functor*. An adjunction is called *monadic* if its comparison functor is an equivalence.

There exists an analogous version of Lemma 2.22 for the Kleisli category of a monad $T$. The comparison functor in this situation is given by

$$K_T\colon \mathscr{C}_T \longrightarrow \mathscr{D}, \qquad x \longmapsto Fx, \qquad f \in \mathscr{D}(x, Ty) \longmapsto \varepsilon_{Fy} \circ Ff.$$

**Remark 2.23.** In fact, the Eilenberg–Moore and the Kleisli category are the terminal and initial objects in the suitably defined category of adjunctions $F\colon \mathscr{C} \rightleftarrows \mathscr{D} : U$ producing the monad $T$. Thus, the diagram

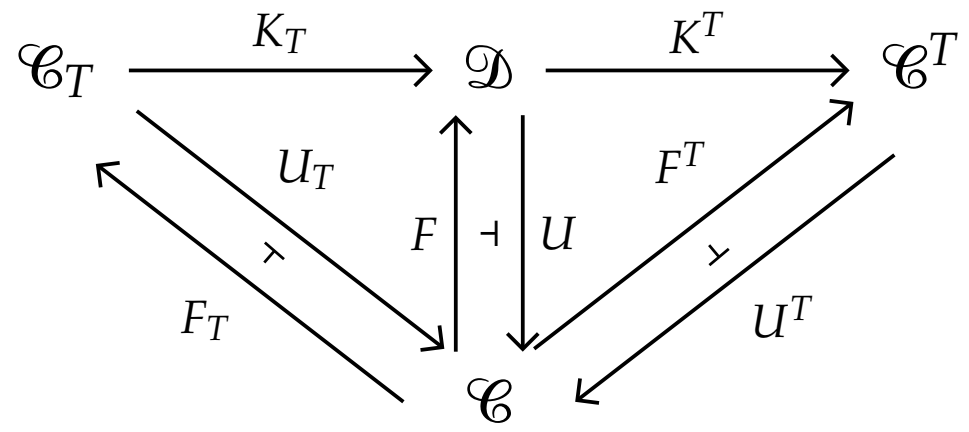

commutes, and its commutativity characterises $K^T$ and $K_T$ completely.

The composite $K^T \circ K_T$ is fully faithful, its full image consisting of *free modules*. We will denote the canonical inclusion of $\mathscr{C}_T$ into $\mathscr{C}^T$ by

$$\iota \colon \mathscr{C}_T \hookrightarrow \mathscr{C}^T, \qquad x \longmapsto (Tx, \mu_x), \qquad f \in \mathscr{C}(x, Ty) \longmapsto \mu_y \circ Tf. \tag{2.2.2}$$

**Example 2.24.** A *comonad* on a category $\mathscr{C}$ consists of an endofunctor $S$ on $\mathscr{C}$ together with natural transformations $\Delta \colon S \Longrightarrow S^2$ and $\varepsilon \colon S \Longrightarrow \mathrm{Id}_{\mathscr{C}}$, satisfying axioms dual to those of Definition 2.9.

For a comonad $S$ on $\mathscr{C}$ an *$S$-comodule* or *$S$-coalgebra* is an object $x \in \mathscr{C}$ together with a morphism $x \longrightarrow Sx$ satisfying axioms analogous to Definition 2.12. We will denote the category of $S$-comodules by $\mathscr{C}^S$ and the Kleisli category of $S$ by $\mathscr{C}_S$. They have analogously defined adjunctions

$$\mathscr{C}_S \underset{U_S}{\overset{F_S}{\underset{\perp}{\rightleftarrows}}} \mathscr{C} \qquad \text{and} \qquad \mathscr{C}^S \underset{U^S}{\overset{F^S}{\underset{\perp}{\rightleftarrows}}} \mathscr{C}.$$

The Eilenberg–Moore category and the Kleisli category for the comonad $S := FU$ can be characterised as a terminal and initial object, respectively. This yields functors $K_S$ and $K^S$ such that the following diagram commutes

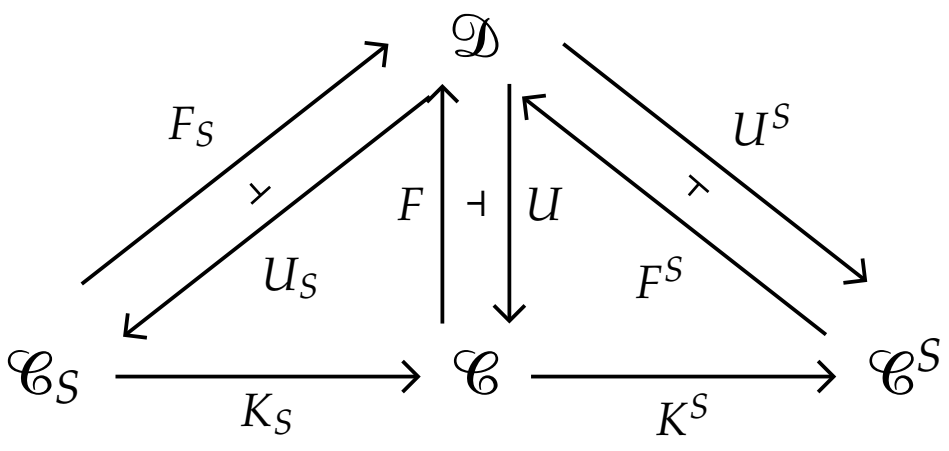

Similarly to Example 2.11, given an adjunction $F \colon \mathscr{C} \rightleftarrows \mathscr{D} : U$, we obtain a comonad $FU$ on $\mathscr{D}$.

It may be that a monad $T$ on $\mathscr{C}$ has a right adjoint $G$—in this case, one can automatically equip $G$ with a canonical comonad structure.

**Proposition 2.25** ([MLM92, Theorem V.8.2])**.** *If $T$ is a monad on a category $\mathscr{C}$ and $G$ is a right adjoint to $T$, then $G$ is a comonad and there is a canonical isomorphism between the Eilenberg–Moore categories: $\mathscr{C}^T \cong \mathscr{C}^G$.*

**Proposition 2.26** ([Kle90, Theorem 3])**.** *If $T$ is a monad on $\mathscr{C}$ and $S$ is a left adjoint to $T$, then $S$ is a comonad and there is a canonical isomorphism $\mathscr{C}_T \cong \mathscr{C}_S$.*

In general, it is *not* true that $\mathscr{C}^T \cong \mathscr{C}^S$ in the setting of Proposition 2.26.

### 2.2.2 *Distributive laws*

Beck's theory of distributive laws concerns itself with the question when, given a monad $T$ on a category $\mathscr{C}$, a functor $S\colon \mathscr{C} \longrightarrow \mathscr{C}$ *lifts* to a functor on the Eilenberg–Moore category of $T$, [Bec69]. That is, there is some $\tilde{S}\colon \mathscr{C}^T \longrightarrow \mathscr{C}^T$, that satisfies the lifting condition

$$\begin{array}{ccc} & & \mathscr{C}^T \\ & \overset{\tilde{S}}{\nearrow} & \big\downarrow U^T \\ \mathscr{C}^T & \xrightarrow[S\circ U^T]{} & \mathscr{C} \end{array}$$

If $S$ is a monad this in turn is equivalent to $ST$ being a monad itself. Distributive laws thus yield a witness for the composability of two monads. Street further developed the theory of distributive laws intrinsic to certain well-behaved bicategories, see [Str72].

**Definition 2.27.** Let $T$ and $S$ be two monad on a category $\mathscr{C}$. A *distributive law* of $T$ over $S$ consists of a natural transformation $\Omega\colon TS \Longrightarrow ST$ such that the following diagrams commute:

$$\begin{array}{ccccc} TSS & \xrightarrow{\Omega S} & STS & \xrightarrow{S\Omega} & SST \\ {\scriptstyle T\mu^S}\big\downarrow & & & & \big\downarrow{\scriptstyle \mu^S T} \\ TS & & \xrightarrow[\Omega]{} & & ST \end{array} \qquad \begin{array}{ccc} T & \xrightarrow{T\eta^S} & TS \\ & {\scriptstyle \eta^S T}\searrow & \big\downarrow{\scriptstyle \Omega} \\ & & ST \end{array}$$

$$\begin{array}{ccccc} TTS & \xrightarrow{T\Omega} & TST & \xrightarrow{\Omega T} & STT \\ {\scriptstyle \mu^T S}\big\downarrow & & & & \big\downarrow{\scriptstyle S\mu^T} \\ TS & & \xrightarrow[\Omega]{} & & ST \end{array} \qquad \begin{array}{ccc} S & \xrightarrow{\eta^T S} & TS \\ & {\scriptstyle S\eta^T}\searrow & \big\downarrow{\scriptstyle \Omega} \\ & & ST \end{array}$$

One analogously defines a distributive law between two comonads, and *mixed distributive laws* between a monad and a comonad.

**Theorem 2.28** ([Bec69])**.** *Let $S$ and $T$ be monads on a category $\mathscr{C}$. Then the following are equivalent:*

(*i*) *distributive laws of $T$ over $S$;*
(*ii*) *monad structures $S \circ_\Omega T := (ST, \mu, \eta)$ on $ST$, such that $S\eta^T$ and $\eta^S T$ are morphisms of monads and $1_{ST} = \mu \circ S\eta^T\eta^T T$; and*

(*iii*) *lifts of $S$ to the Eilenberg–Moore category of $T$, such that $S$ is a monad on $\mathscr{C}^T$.*

*Sketch of proof.* We outline the main constructions.

Given a distributive law $\Omega\colon TS \Longrightarrow ST$, for two monads on $\mathscr{C}$, the monad structure on $ST$ is given by

$$\mu\colon STST \xrightarrow{S\Omega T} SSTT \xrightarrow{\mu^S\mu^T} ST \qquad \text{and} \qquad \eta\colon \mathrm{Id}_{\mathscr{C}} \xrightarrow{\eta^T} T \xrightarrow{\eta^S T} ST.$$

The lift $\tilde{S}$ of $S$ to $\mathscr{C}^T$ is defined by $\tilde{S}(x, \nabla_x) := (Sx, TSx \xrightarrow{\Omega_x} STx \xrightarrow{S\nabla_x} Sx)$. Its multiplication and unit is given by those of $S$.

Given a monad structure $(ST, \mu, \eta)$, one defines a distributive law by

$$TS \xrightarrow{\eta^S TS\eta^T} STST \xrightarrow{\mu} ST.$$

Finally, if $\tilde{S}$ is a lift of $S$, we obtain a distributive law via

$$TS \xrightarrow{TS\eta^T} TSU^TF^T = U^TF^TU^T\tilde{S}F^T \xrightarrow{U^T\varepsilon^T\tilde{S}F^T} U^T\tilde{S}F^T = ST. \qquad \square$$

## 2.3 STRING DIAGRAMS

RECALL FROM EXAMPLE 2.2 THAT $\mathbb{C}\mathrm{at}$ is the 2-category of all (small) categories, functors, and natural transformations. We use juxtaposition for the horizontal composition of $\mathbb{C}\mathrm{at}$, or, in case we want to emphasise the direction of composition,

$$- \odot = \colon \mathbb{C}\mathrm{at}(\mathscr{B}, \mathscr{C}) \times \mathbb{C}\mathrm{at}(\mathscr{A}, \mathscr{B}) \longrightarrow \mathbb{C}\mathrm{at}(\mathscr{A}, \mathscr{C}), \qquad \mathscr{A}, \mathscr{B}, \mathscr{C} \in \mathbb{C}\mathrm{at}.$$

String diagrams will serve as an important tool for doing computations. In the case of $\mathbb{C}\mathrm{at}$, a *string diagram* consists of regions labelled with categories, strings labelled with functors, and vertices between the strings labelled with natural transformations. If two string diagrams can be transformed into each other—that is, if they are isotopic—then the natural transformations they represent are equal. A more detailed description is given in [JS91; Sel11]. Our convention is to read diagrams from bottom to top and right to left. Horizontal and vertical composition are given by horizontal and vertical gluing of diagrams, respectively. Identity natural transformations are given by unlabelled vertices. The edge for the identity functor of a category will not be drawn. If the involved categories are clear from the context, we omit writing them explicitly. Figure 2.1 details our conventions.

More generally, our applications of the graphical calculus can be formulated in any bicategory that admits the construction of algebras in the sense of Example 2.14. For ease of presentation we shall stick with $\mathbb{C}\mathrm{at}$.

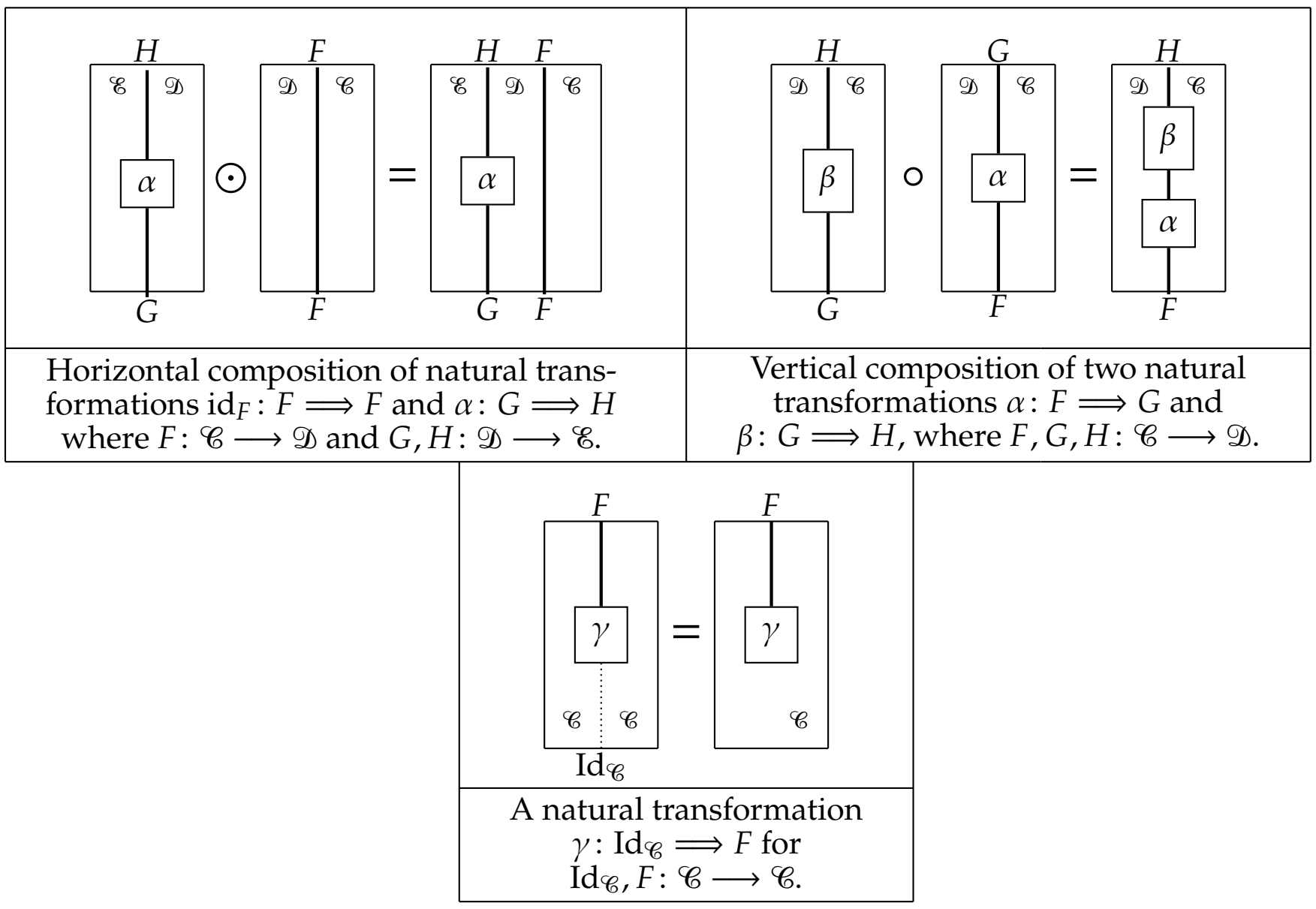

Figure 2.1: Basic string diagrammatic conventions.

**Example 2.29.** Let $T$ be a monad on $\mathscr{C}$. We represent the multiplication and unit of $T$ in terms of string diagrams:

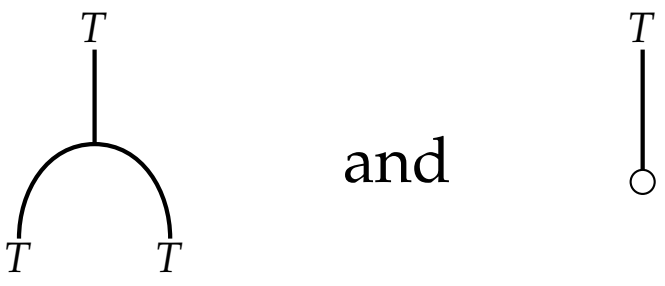

Their associativity and unitality then equate to

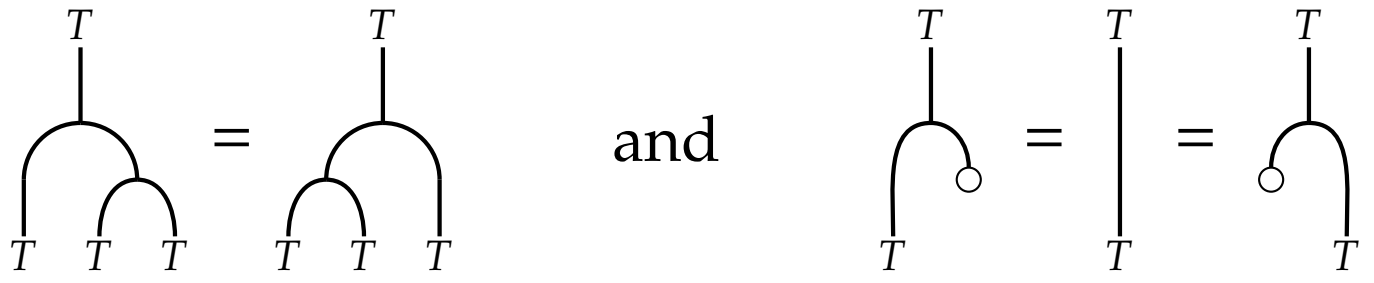

We depict the unit and counit of the Eilenberg–Moore adjunction as

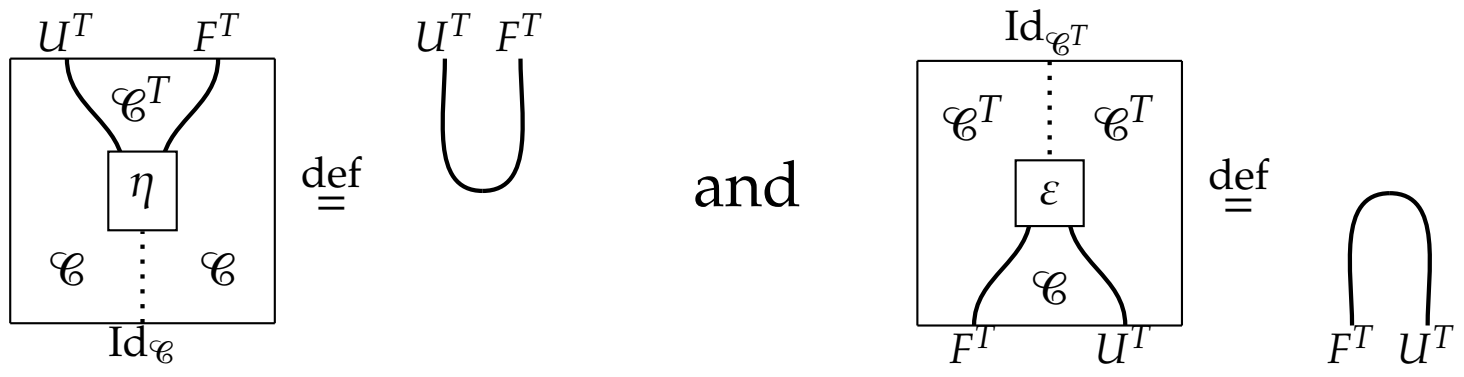

The defining equations of adjunctions translate to the snake equations

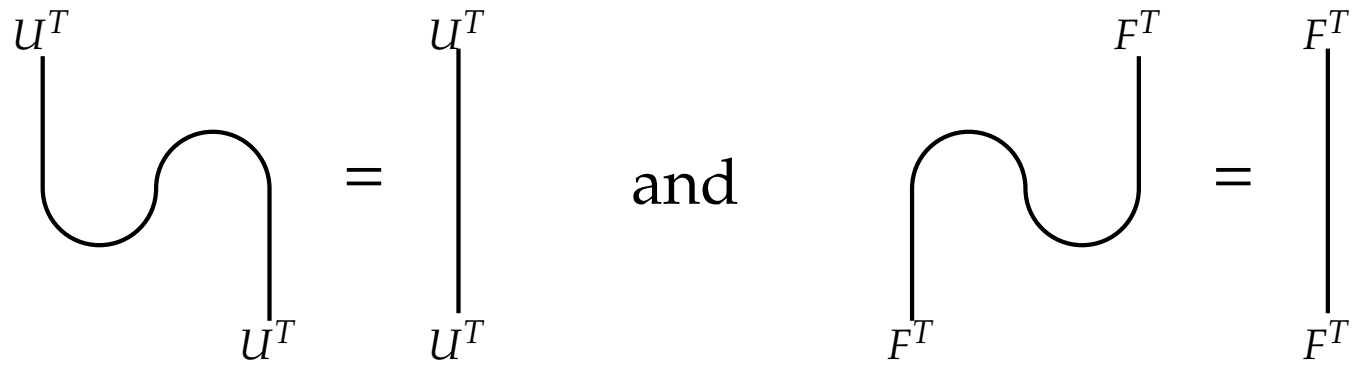

**Remark 2.30.** Given a monad $T$, its category of algebras can be incorporated into the graphical calculus, see [Wil08]. Define the natural transformation

$$\nabla\colon TU^T = U^T F^T U^T \xrightarrow{U^T \varepsilon} U^T. \tag{2.3.1}$$

Now, one may write

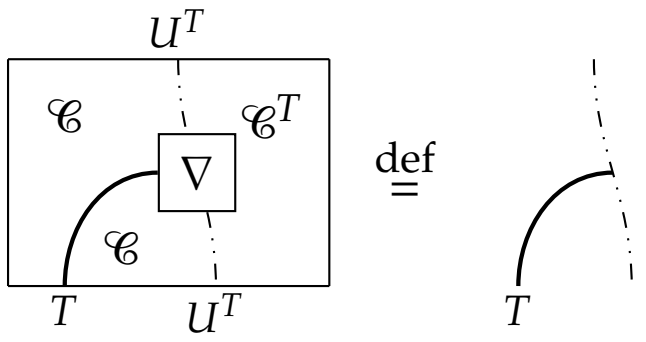

We think of $\nabla\colon TU^T \longrightarrow U^T$ as an action of $T$ on $U^T$ due to the identities

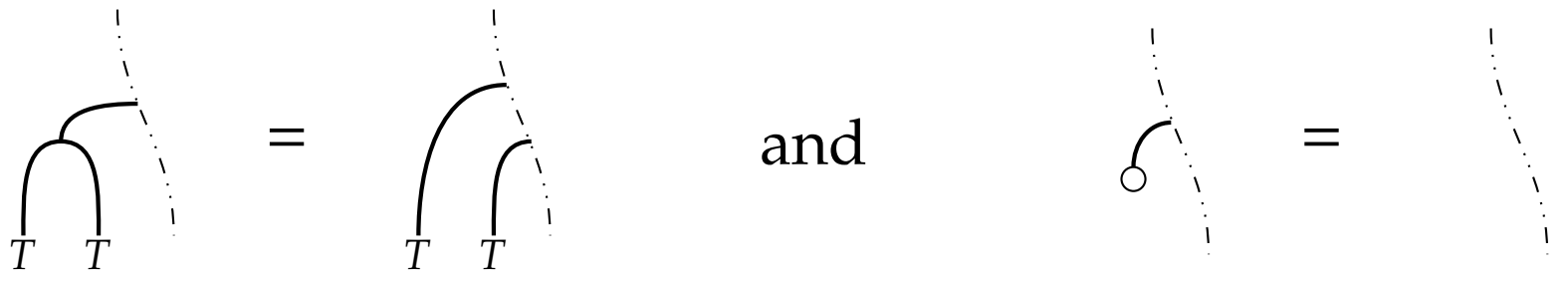

By abuse of notation, we shall often simply speak of an object $x \in \mathscr{C}^T$, leaving the action $\nabla_x$ implicit.

## 2.4 MONOIDAL AND MODULE CATEGORIES

FOR A MORE EXTENSIVE ACCOUNT concerning the notions of monoidal and module categories, we refer the reader to [EGNO15, Chapter 2].

**Definition 2.31.** A *monoidal category* consists of a category $\mathscr{C}$ together with a *tensor product* functor $\otimes\colon \mathscr{C} \times \mathscr{C} \longrightarrow \mathscr{C}$, a *unit* $1 \in \mathscr{C}$, and for all $x, y, z \in \mathscr{C}$, natural isomorphisms

$$\alpha_{x,y,z}\colon (x \otimes y) \otimes z \xRightarrow{\sim} x \otimes (y \otimes z), \qquad \lambda_x\colon 1 \otimes x \xRightarrow{\sim} x, \qquad \rho_x\colon x \otimes 1 \xRightarrow{\sim} x,$$

satisfying *coherence* axioms; see [EGNO15, Definition 2.1.1].

A monoidal category is called *strict* if $\alpha$, $\lambda$, and $\rho$ are identities.

**Example 2.32.** Given a monoidal category $(\mathscr{C}, \otimes, 1)$, we denote by $\mathscr{C}^{\mathrm{rev}}$ the monoidal category $(\mathscr{C}, \otimes^{\mathrm{op}}, 1)$, where for all $x, y \in \mathscr{C}$ one sets $x \otimes^{\mathrm{op}} y := y \otimes x$.

Recall *Sweedler notation*: if $B$ is a bialgebra in $\mathsf{Vect}$ with comultiplication $\Delta\colon B \longrightarrow B \otimes_{\Bbbk} B$, then for $b \in B$ we write $b_{(1)} \otimes b_{(2)} := \Delta(b) \in B \otimes_{\Bbbk} B$. For example, in this notation coassociativity translates to

$$b_{(1)} \otimes (b_{(2)})_{(1)} \otimes (b_{(2)})_{(2)} = b_{(1)} \otimes b_{(2)} \otimes b_{(3)} = (b_{(1)})_{(1)} \otimes (b_{(1)})_{(2)} \otimes b_{(2)}.$$

This naturally extends to left and right coactions on $B$, which we respectively denote by $x \longmapsto x_{(-1)} \otimes x_{(0)}$ and $x \longmapsto x_{(0)} \otimes x_{(1)}$.

**Example 2.33** ([Yet90]). Let $H \in \mathsf{Vect}$ be a bialgebra The category ${}^H_H\mathcal{YD}$ of (left-left) *Yetter–Drinfeld modules* is defined as the category of those $Y \in \mathsf{Vect}$ that are left $H$-modules with action $\triangleright$ and also left $H$-comodules with coaction $y \longmapsto y_{(-1)} \otimes y_{(0)}$, such that the following compatibility condition holds for all $h \in H$ and $y \in Y$:

$$h_{(1)}y_{(-1)} \otimes (h_{(2)} \triangleright y_{(0)}) = (h_{(1)} \triangleright y)_{(-1)}h_{(2)} \otimes (h_{(1)} \triangleright y)_{(0)} \tag{2.4.1}$$

Then ${}^H_H\mathcal{YD}$ is monoidal: given Yetter–Drinfeld modules $Y$ and $Z$, their tensor product $Y \otimes_{\Bbbk} Z$ can be equipped with the structure of a Yetter–Drinfeld module. For all $h \in H$, $y \in Y$, and $z \in Z$, define

$$h \triangleright (y \otimes z) := (h_{(1)} \triangleright y) \otimes (h_{(2)} \triangleright z) \quad \text{and} \quad y \otimes z \longmapsto y_{(-1)}z_{(-1)} \otimes y_{(0)} \otimes z_{(0)}.$$

In fact, if $H$ is a Hopf algebra with invertible antipode $S$—for example, when $H$ is finite-dimensional—then Equation (2.4.1) is equivalent to

$$(h \triangleright y)_{(-1)} \otimes (h \triangleright y)_{(0)} = h_{(1)}y_{(-1)}S(h_{(3)}) \otimes (h_{(2)} \triangleright y_{(0)}), \tag{2.4.2}$$

see [Mon93, Proposition 10.6.16].

**Definition 2.34.** A *lax monoidal functor* between monoidal categories $\mathcal{C}$ and $\mathcal{D}$ consists of a functor $F\colon \mathcal{C} \longrightarrow \mathcal{D}$ together with a morphism $F_0\colon 1 \longrightarrow F1$ and a natural transformation $F_2\colon F(-) \otimes F(=) \Longrightarrow F(- \otimes =)$, satisfying *associativity* and *unitality* conditions, see [EGNO15, Definition 2.4.1].

Analogously, an *oplax monoidal functor* is a functor $G\colon \mathcal{C} \longrightarrow \mathcal{D}$ together with morphisms $G_0\colon G1 \longrightarrow 1$ and $G_2\colon G(- \otimes =) \Longrightarrow G(-) \otimes G(=)$, making $G^{\mathrm{op}}\colon \mathcal{C}^{\mathrm{op}} \longrightarrow \mathcal{D}^{\mathrm{op}}$ lax monoidal. We call $F$ *strong monoidal* if $F_0$ and $F_2$ and invertible, and *strict monoidal* if they are identities.

Observe that an oplax monoidal functor with invertible coherence morphisms canonically defines a strong monoidal functor.

**Remark 2.35.** Let $T$ be an oplax monoidal endofunctor on a monoidal category $\mathcal{C}$. By coassociativity of $T_2$, there exists a well-defined natural transformation

$$T_3\colon T(- \otimes = \otimes \equiv) \Longrightarrow T(-) \otimes T(=) \otimes T(\equiv),$$

which, for all $x, y, z \in \mathcal{C}$, has components

$$T_{3;x,y,z} := (T_{2;x,y} \otimes Tz) \circ T_{2;x \otimes y,z} = (Tx \otimes T_{2;y,z}) \circ T_{2;x,y \otimes z}.$$

**Definition 2.36.** Let $\mathscr{C}$ and $\mathscr{D}$ be two monoidal categories and suppose that $F, G\colon \mathscr{C} \longrightarrow \mathscr{D}$ are two lax monoidal functors. A *(lax) monoidal natural transformation* is a natural transformation $\varphi\colon F \Longrightarrow G$ such that the following two diagrams commute for all $x, y \in \mathscr{C}$:

$$\begin{array}{ccc} Fx \otimes Fy & \xrightarrow{F_{2;x,y}} & F(x \otimes y) \\ {\scriptstyle \varphi_x \otimes \varphi_y}\downarrow & & \downarrow{\scriptstyle \varphi_{x \otimes y}} \\ Gx \otimes Gy & \xrightarrow[G_{2;x,y}]{} & G(x \otimes y) \end{array} \qquad \begin{array}{ccc} 1 & \xrightarrow{F_0} & F1 \\ {\scriptstyle G_0}\downarrow & \swarrow{\scriptstyle \varphi_1} & \\ G1 & & \end{array}$$

Dually, given two oplax monoidal functors $F, G\colon \mathscr{C} \longrightarrow \mathscr{D}$, one defines an *oplax monoidal natural transformation* as an arrow $\varphi\colon F \Longrightarrow G$ making the following diagrams commute for all $x, y \in \mathscr{C}$:

$$\begin{array}{ccc} F(x \otimes y) & \xrightarrow{F_{2;x,y}} & Fx \otimes Fy \\ {\scriptstyle \varphi_{x \otimes y}}\downarrow & & \downarrow{\scriptstyle \varphi_x \otimes \varphi_y} \\ G(x \otimes y) & \xrightarrow[G_{2;x,y}]{} & Gx \otimes Gy \end{array} \qquad \begin{array}{ccc} F1 & \xrightarrow{F_0} & 1 \\ {\scriptstyle \varphi_1}\downarrow & \nearrow{\scriptstyle G_0} & \\ G1 & & \end{array}$$

**Definition 2.37.** An adjunction $F\colon \mathscr{C} \rightleftarrows \mathscr{D} :U$ is called *lax monoidal* if $F$ and $U$ are lax monoidal functors, and the unit and counit of the adjunction are monoidal natural transformations.

Analogously to Definition 2.37 one defines *oplax monoidal adjunctions*.

**Definition 2.38.** Let $F\colon \mathscr{C} \rightleftarrows \mathscr{D} :U$ be an ordinary adjunction between monoidal categories. A *lift of $F \dashv U$ to a lax monoidal adjunction* consists of choices of lax monoidal structures on $F$ and $U$, together with assertions that the unit and counit of the adjunction are lax monoidal natural transformations. In other words, it is the structure necessary for a lax monoidal adjunction whose left adjoint is $F$ and whose right adjoint is $U$.

Similarly, one can define a lift of $F \dashv U$ to an oplax monoidal adjunction.

**Remark 2.39.** A strong monoidal functor is called a *monoidal equivalence* if it additionally is an equivalence of categories. The quasi-inverse $F^{-1}\colon \mathscr{D} \longrightarrow \mathscr{C}$ of a monoidal equivalence $F\colon \mathscr{C} \longrightarrow \mathscr{D}$ is again monoidal, and there are monoidal natural isomorphisms $F \circ F^{-1} \xrightarrow{\sim} \mathrm{Id}_{\mathscr{D}}$ and $F^{-1} \circ F \xrightarrow{\sim} \mathrm{Id}_{\mathscr{C}}$.

Unless otherwise specified, all monoidal categories in the rest of this thesis are assumed to be strict. This is justified by the following *strictification* result of Mac Lane, [ML63], see also [EGNO15, Theorem 2.8.5].

**Theorem 2.40.** *Every monoidal category is monoidally equivalent to a strict monoidal category.*

**Definition 2.41.** A *(left) $\mathscr{C}$-module category* of a monoidal category $\mathscr{C}$ consists of a category $\mathscr{M}$ together with a functor $\triangleright\colon \mathscr{C} \times \mathscr{M} \longrightarrow \mathscr{M}$ and isomorphisms

$$\mathscr{M}_{\mathsf{a};x,y,m}\colon (x \otimes y) \triangleright m \xRightarrow{\sim} x \triangleright (y \triangleright m),$$

natural in $m \in \mathscr{M}$ and $x, y \in \mathscr{C}$, satisfying *coherence* axioms similar to those of a monoidal category; see [EGNO15, Definition 7.1.1].

Analogously, one can define *right* module categories over $\mathscr{C}$, which involve a functor $\triangleleft\colon \mathscr{M} \times \mathscr{C} \longrightarrow \mathscr{M}$ and isomorphisms $m \triangleleft (x \otimes y) \xRightarrow{\sim} (x \triangleleft y) \triangleleft m$. For brevity, we may sometimes speak of just "$\mathscr{C}$-module categories"; this will always mean "left $\mathscr{C}$-module categories".

**Example 2.42** ([HKRS04])**.** Let $H$ be a Hopf algebra with invertible antipode $S$. Recall from Example 2.33 that the category ${}^{H}_{H}\mathscr{YD}$ of Yetter–Drinfeld modules is a monoidal category. A left-left *anti-Yetter–Drinfeld module* is some $M \in \mathsf{Vect}$ equipped with a left $H$-module structure $\triangleright$ and a left $H$-comodule structure $m \longmapsto m_{(-1)} \otimes m_{(0)}$ such that the following holds for all $h \in H$ and $m \in M$:

Notice that the only difference to Equation (2.4.2) is that we use $S^{-1}$ instead of $S$.

$$(h \triangleright m)_{(-1)} \otimes (h \triangleright m)_{(0)} = h_{(1)} m_{(-1)} S^{-1}(h_{(3)}) \otimes (h_{(2)} \triangleright m_{(0)})$$

The category of anti-Yetter–Drinfeld modules do not form a monoidal category, but rather a $\mathscr{C}$-module category over the Yetter–Drinfeld modules. Given $M \in {}^{H}_{H}\mathsf{a}\mathscr{YD}$ and $Y \in {}^{H}_{H}\mathscr{YD}$, define $Y \triangleright M \in {}^{H}_{H}\mathsf{a}\mathscr{YD}$ as the vector space $Y \otimes_{\Bbbk} M$, equipped with the following action and coaction:

$$h \triangleright (y \otimes m) := (h_{(1)} \triangleright y) \otimes (h_{(2)} \triangleright m), \qquad y \otimes m \longmapsto y_{(-1)} m_{(-1)} \otimes y_{(0)} \otimes m_{(0)}.$$

**Definition 2.43.** Let $\mathscr{C}$ and $\mathscr{D}$ be monoidal categories. A $(\mathscr{C}, \mathscr{D})$*-bimodule category* consists of a left $\mathscr{C}$-module category $\mathscr{M}$ that is simultaneously a right $\mathscr{D}$-module category, such that there exists a natural isomorphism

$$(x \triangleright m) \triangleleft y \xRightarrow{\sim} x \triangleright (m \triangleleft y),$$

for all $x \in \mathscr{C}$, $y \in \mathscr{D}$, and $m \in \mathscr{M}$, called the *middle interchange*, which satisfies appropriate *associativity* constraints; see [EGNO15, Definition 7.1.7].

**Example 2.44.** Every monoidal category $\mathscr{C}$ is a $(\mathscr{C},\mathscr{C})$-bimodule category, whose actions are both given by the tensor product of $\mathscr{C}$. We call this module category the *regular bimodule*, and denote it by ${}_{\mathscr{C}}\mathscr{C}_{\mathscr{C}}$. Analogously, one defines the left and right regular $\mathscr{C}$-module categories ${}_{\mathscr{C}}\mathscr{C}$ and $\mathscr{C}_{\mathscr{C}}$.

**Example 2.45.** Let $\mathscr{C}$ be a monoidal category. The regular action is not the only way in which we can consider $\mathscr{C}$ as a bimodule over itself. Suppose that $R\colon \mathscr{C} \longrightarrow \mathscr{C}$ is a strong monoidal functor. The action

$$\triangleleft\colon \mathscr{C}\times\mathscr{C} \xrightarrow{\mathscr{C}\times R} \mathscr{C}\times\mathscr{C} \xrightarrow{\otimes} \mathscr{C}$$

endows $\mathscr{C}$ with the structure of a right module category over itself, which we shall denote by $\mathscr{C}_R$. In other words, we have

$$x \triangleright y := Rx \otimes y, \qquad\qquad f \triangleright g := Rf \otimes g,$$

for $v, w, x, y \in \mathscr{C}$ and $f\colon v \longrightarrow w$, $g\colon x \longrightarrow y$.

We can also twist the action from the left with another monoidal functor $L$. The *twisted* bimodule obtained in this manner is denoted by ${}_L\mathscr{C}_R$.

**Remark 2.46.** One can easily imagine a more involved setting than Example 2.45 by twisting with an oplax monoidal functor $L\colon \mathscr{C} \longrightarrow \mathscr{C}$ from the left and a lax monoidal functor $R$ from the right. In this setting, ${}_L\mathscr{C}_R$ is a lax left and oplax right $\mathscr{C}$-bimodule category, see for example [Szl12, Section 2]. Since in most of this thesis we only consider the strong case, we refrain from more formally treating this case; however, see Remark 6.13.

**Definition 2.47.** Let $\mathscr{M}$ and $\mathscr{N}$ be left module categories over a monoidal category $\mathscr{C}$. A *lax $\mathscr{C}$-module functor* from $\mathscr{M}$ to $\mathscr{N}$ is a functor $F\colon \mathscr{M} \longrightarrow \mathscr{N}$, together with a collection of morphisms

$$F_{\mathsf{a};x,m}\colon x \triangleright Fm \longrightarrow F(x \triangleright m), \qquad \text{natural for all } x \in \mathscr{C} \text{ and } m \in \mathscr{M},$$

satisfying *associativity* and *unitality* conditions, see [EGNO15, Definition 7.2.1].

*Oplax* and *strong* $\mathscr{C}$-module functors are defined similarly; in the former case, one considers arrows $F_{\mathsf{a};x,m}\colon F(x \triangleright m) \longrightarrow x \triangleright Fm$, while for the latter $F_{\mathsf{a}}$ should be invertible.

**Definition 2.48.** Let $F, G\colon \mathscr{M} \longrightarrow \mathscr{N}$ be $\mathscr{C}$-module functors between the left $\mathscr{C}$-module categories $\mathscr{M}$ and $\mathscr{N}$. A natural transformation $\phi\colon F \Longrightarrow G$ is called a *$\mathscr{C}$-module transformation* if

$$G_{\mathsf{a};x,m} \circ (\phi_x \triangleright m) = \phi_{x\triangleright m} \circ F_{\mathsf{a};x,m}, \qquad\qquad \text{for all } x \in \mathscr{C} \text{ and } m \in \mathscr{M}.$$

**Notation 2.49.** For $\mathscr{C}$-module categories $\mathscr{M}$ and $\mathscr{N}$, we obtain the following categories of $\mathscr{C}$-module functors:

- $\mathsf{\mathscr{C}Mod}(\mathscr{M}, \mathscr{N}) := \mathsf{Str\mathscr{C}Mod}(\mathscr{M}, \mathscr{N})$: strong $\mathscr{C}$-module functors from $\mathscr{M}$ to $\mathscr{N}$,
- $\mathsf{Lax\mathscr{C}Mod}(\mathscr{M}, \mathscr{N})$: lax $\mathscr{C}$-module functors from $\mathscr{M}$ to $\mathscr{N}$,
- $\mathsf{Oplax\mathscr{C}Mod}(\mathscr{M}, \mathscr{N})$: oplax $\mathscr{C}$-module functors from $\mathscr{M}$ to $\mathscr{N}$.

**Example 2.50.** For an object $x$ is a (non-strict) monoidal category $\mathscr{C}$, the functor $F := - \otimes x$ is a strong left $\mathscr{C}$-module functor. The transformation $F_{\mathsf{a}}$ is given by the associator of $\mathscr{C}$. Further, a morphism $f\colon x \longrightarrow y$ in $\mathscr{C}$ yields a $\mathscr{C}$-module transformation $- \otimes f\colon - \otimes x \Longrightarrow - \otimes y$.

In fact, Example 2.50 generalises and completely describes module functors from the regular module. We emphasise that this result is a consequence of the bicategorical Yoneda lemma, see [JY21, Lemma 8.3.16].

**Proposition 2.51.** *Let $\mathscr{M}$ be a left module category over the monoidal category $\mathscr{C}$. Then there is an equivalence of module categories*

$$\mathsf{Str\mathscr{C}Mod}(\mathscr{C}, \mathscr{M}) \xrightarrow{\sim} \mathscr{M}, \qquad F \longmapsto F1, \qquad - \triangleright m \longleftarrow m.$$

In particular, this yields a monoidal equivalence $\mathsf{Str\mathscr{C}Mod} \simeq \mathscr{C}^{\mathrm{rev}}$.

### 2.4.1 *Braidings*

Braidings are natural transformations relating the tensor product to its opposite. They where introduced by Joyal and Street in [JS85; JS86; JS93], and build on the notion of symmetries studied in [ML63; EK66].

**Definition 2.52.** A *braiding* on a monoidal category $\mathscr{C}$ is a natural isomorphism $\sigma_{-,=}\colon - \otimes = \Longrightarrow = \otimes -$, satisfying $\sigma_{x,1} = \mathrm{id}_x$ and the *hexagon axioms*[2].
The pair $(\mathscr{C}, \sigma)$ will be called a *braided monoidal category*.

[2] The name "hexagon axioms" is due to the fact that in the non-strict setting, the defining equations can be organised as a hexagon-shaped diagram, see [JS93].

**Remark 2.53.** For a braiding $\sigma$, the condition that $\sigma_{x,1} = \mathrm{id}_x$ is actually implied by the hexagon identities: The calculation

$$\sigma_{x,1} = \sigma_{x,1\otimes 1} = (\mathrm{id}_1 \otimes \sigma_{x,1}) \circ (\sigma_{x,1} \otimes \mathrm{id}_1) = \sigma_{x,1} \circ \sigma_{x,1}$$

shows that $\sigma_{x,1}$ is an invertible idempotent, hence the claim follows.

**Example 2.54** ([Yet90, Theorem 7.2]). The category of left-left Yetter–Drinfeld modules of Example 2.33 is braided monoidal when equipped with

$$\sigma := \{\, \sigma_{Y,Z} \colon Y \otimes Z \longrightarrow Z \otimes Y, \quad y \otimes z \longmapsto (y_{(-1)} \triangleright z) \otimes y_{(0)} \,\}_{Y,Z \in {}^{H}_{H}\mathcal{YD}},$$

called the *Yetter–Drinfeld braiding*. Note that $\sigma$ is invertible because $H$ has an invertible antipode. Its inverse is given for all $Y, Z \in {}^{H}_{H}\mathcal{YD}$ by

$$\sigma^{-1} := \{\, \sigma^{-1}_{Y,Z} \colon Z \otimes Y \longrightarrow Y \otimes Y, \quad z \otimes y \longmapsto y_{(0)} \otimes (S^{-1}(y_{(-1)}) \triangleright z) \,\}_{Y,Z \in {}^{H}_{H}\mathcal{YD}}.$$

**Definition 2.55.** A braiding $\sigma$ on a monoidal category $\mathscr{C}$ is called *symmetric* if $\sigma^{-1}_{x,y} = \sigma_{y,x}$ for all $x, y \in \mathscr{C}$. A monoidal category equipped with a symmetric braiding will be referred to as a *symmetric monoidal category*.

Braidings are depicted in the graphical calculus by crossings of strings subject to Reidemeister-esque identities, see [Sel11]. The following figure shows a braiding, its inverse, the hexagon identity, and the naturality of the braiding in its first argument:

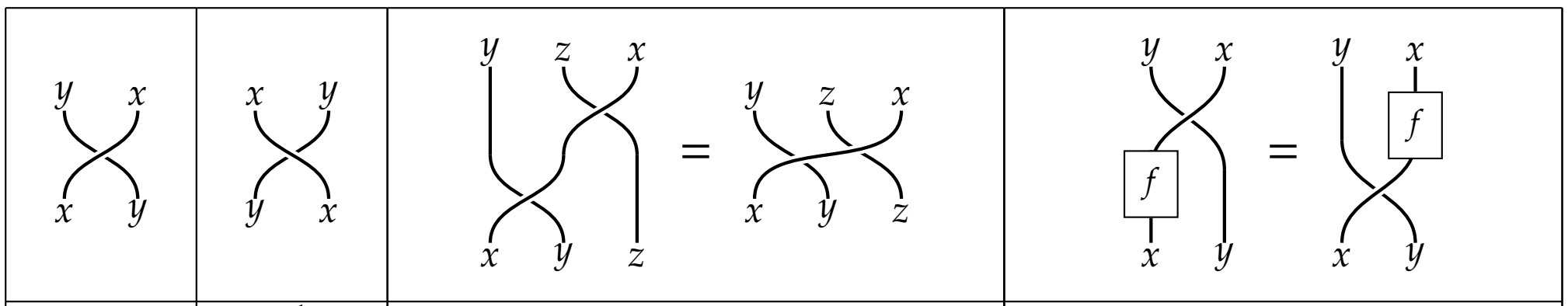

| $\sigma_{x,y}$ | $\sigma^{-1}_{x,y}$ | $(\mathrm{id}_y \otimes \sigma_{x,z}) \circ (\sigma_{x,y} \otimes \mathrm{id}_z) = \sigma_{x,y\otimes z}$ | $\sigma_{w,y} \circ (f \otimes \mathrm{id}_y) = (y \otimes f) \circ \sigma_{x,y}$ |
|---|---|---|---|

### 2.4.2 *Closedness*

**Definition 2.56.** Let $\mathscr{C}$ be a monoidal category. An object $m$ in a $\mathscr{C}$-module category $\mathscr{M}$ is called *closed* if the functor $- \triangleright m \colon \mathscr{C} \longrightarrow \mathscr{M}$ has a right adjoint $\lfloor m, - \rfloor \colon \mathscr{M} \longrightarrow \mathscr{C}$, called the *internal hom from* $m$. We refer to

$$- \triangleright m \colon \mathscr{C} \rightleftarrows \mathscr{M} : \lfloor m, - \rfloor$$

as the *internal tensor–hom adjunction*.

Similarly to Definition 2.56, an object $m \in \mathscr{M}$ is called *coclosed* if the functor $- \triangleright m \colon \mathscr{C} \longrightarrow \mathscr{M}$ has a left adjoint; we denote this *internal cohom from* $m$ by $\lceil m, - \rceil \colon \mathscr{M} \longrightarrow \mathscr{C}$, and call $\lceil m, - \rceil \colon \mathscr{M} \rightleftarrows \mathscr{C} : - \triangleright m$ as the *internal cohom–tensor adjunction*. If every object of $\mathscr{M}$ is closed, we call $\mathscr{M}$ a *closed $\mathscr{C}$-module category*.

**Definition 2.57.** A monoidal category $\mathscr{C}$ is called *left closed* if the left regular $\mathscr{C}$-module category ${}_{\mathscr{C}}\mathscr{C}$ is closed, *right closed* if the right regular $\mathscr{C}$-module category $\mathscr{C}_{\mathscr{C}}$ is closed, and *closed* or *biclosed* if it is both left and right closed.

In these cases, we write $[x,-]_\ell$ and $[x,-]_r$ for $\lfloor x,-\rfloor$, respectively.

**Example 2.58.** A symmetric monoidal category is left closed if and only if it is right closed if and only if it is closed.

The object-wise tensor–hom adjunctions of a closed $\mathscr{C}$-module category specify a unique functor $\lfloor -,=\rfloor\colon \mathscr{M}^{\mathrm{op}} \times \mathscr{M} \longrightarrow \mathscr{C}$ such that we have

$$\mathscr{M}(c \triangleright m, n) \cong \mathscr{C}(c, \lfloor m, n\rfloor),$$

natural in all three variables, see [ML98, Section IV.7] for the monoidal case.

**Example 2.59.** For any object $x \in \mathscr{C}$ of a right closed monoidal category, a combination of the unit and counit

$$\eta_y^{(x)}\colon y \longrightarrow [x, x \otimes y]_r, \qquad \varepsilon_y^{(x)}\colon x \otimes [x,y]_r \longrightarrow y, \qquad \text{for all } y \in \mathscr{C}$$

gives rise to the morphism

$$\phi_y^{(x)} := [x,1]_r \otimes y \xrightarrow{\eta^{(x)}_{[x,1]_r\otimes y}} [x, x\otimes [x,1]_r \otimes y]_r \xrightarrow{[x,\varepsilon_1^{(x)}\otimes y]} [x,y]_r. \tag{2.4.3}$$

### 2.4.3 *Rigidity and pivotality*

Rigid monoidal categories are also called *autonomous*, or, in the symmetric case, *compact closed*.

Rigidity in the context of monoidal categories refers to a concept of duality similar to that of finite-dimensional vector spaces. Importantly, notions like dual bases and evaluations have their analogues in this setting.

Pivotal categories—also known as *balanced* or *sovereign* categories—are those where there exists an identification between objects and their double duals that is compatible with the tensor product.

Recall from Equation (2.4.3) there is a canonical morphism

$$\Phi_y^{(x)}\colon [x,1]_r \otimes y \longrightarrow [x,y]_r$$

for all objects $x, y$ in the right closed monoidal category $\mathscr{C}$. The next results links the invertibility of this morphism to conditions that $\mathscr{C}$ is rigid monoidal; see for example [Kel72], and [NW17, Proposition 2.1] for a concise proof.

**Proposition 2.60.** *Let $\mathscr{C}$ be a monoidal category. For any object $x \in \mathscr{C}$, the following statements are equivalent:*

(*i*) *the right internal hom of $x$ exists, and the canonical arrows $\phi_y^{(x)}$ are invertible for all $y \in \mathscr{C}$;*

(*ii*) *the right internal hom of $x$ exists and $\phi_x^{(x)} \colon [x, 1]_r \otimes x \longrightarrow [x, x]_r$ is an isomorphism; and*

(*iii*) *there exists an object $x^\vee \in \mathscr{C}$ together with morphisms $\mathrm{ev}_x^r \colon x \otimes x^\vee \longrightarrow 1$ and $\mathrm{coev}_x^r \colon 1 \longrightarrow x^\vee \otimes x$, satisfying the snake identities*

$$\begin{aligned} \mathrm{id}_x &= (\mathrm{ev}_x^r \otimes \mathrm{id}_x) \circ (\mathrm{id}_x \otimes \mathrm{coev}_x^r), \\ \mathrm{id}_{x^\vee} &= (\mathrm{id}_{x^\vee} \otimes \mathrm{ev}_x^r) \circ (\mathrm{coev}_x^r \otimes \mathrm{id}_{x^\vee}). \end{aligned} \tag{2.4.4}$$

**Definition 2.61.** We call $x$ *right (rigidly) dualisable* if any of the equivalent conditions of Proposition 2.60 are met. In this case, we have $[x, y]_r \cong x^\vee \otimes y$ for all $x, y \in \mathscr{C}$ and in particular $x^\vee \cong [x, 1]_r$. The *left (rigid) dualisability* of an object $x \in \mathscr{C}$ can be defined similarly by either the invertibility of a canonical morphism $\psi_y^{(x)} \colon y \otimes [x, 1]_\ell \longrightarrow [x, y]_\ell$ or by the existence of an object ${}^\vee x \in \mathscr{C}$ endowed with morphisms $\mathrm{ev}_x^\ell \colon {}^\vee x \otimes x \longrightarrow 1$ and $\mathrm{coev}_x^\ell \colon 1 \longrightarrow x \otimes {}^\vee x$, subject to suitable variants of the snake identities.

When the context leaves the choice unambiguous, we—in the interest of brevity—omit the left and right superscripts for the (co)evaluation morphisms. The pair $(y, x, \mathrm{ev}, \mathrm{coev})$, where $y$ is a left dual of $x$ and, thus, $x$ is a right dual of $y$ is called a *dual pair*.

**Definition 2.62.** A monoidal category $\mathscr{C}$ is called *left (right) rigid* if every object is left (right) rigidly dualisable, and *rigid* if it is left and right rigid.

We denote the left and right *dualising functors* of a rigid category $\mathscr{C}$ by

$${}^\vee(-) \cong [-, 1]_\ell \colon \mathscr{C}^{\mathrm{op,rev}} \longrightarrow \mathscr{C} \qquad \text{and} \qquad (-)^\vee \cong [-, 1]_r \colon \mathscr{C}^{\mathrm{op,rev}} \longrightarrow \mathscr{C}.$$

**Example 2.63.** Let $\mathscr{C}$ be a left rigid monoidal category. Given a morphism $f \colon x \longrightarrow y$ in $\mathscr{C}$, the induced morphism ${}^\vee f \colon {}^\vee y \longrightarrow {}^\vee x$ is given by

$$(\mathrm{ev}_y^\ell \otimes \mathrm{id}_{{}^\vee x}) \circ (\mathrm{id}_{{}^\vee y} \otimes f \otimes \mathrm{id}_{{}^\vee x}) \circ (\mathrm{id}_{{}^\vee y} \otimes \mathrm{coev}_x^\ell).$$

There are categories with only one-sided closedness or rigidity, see [Lor21, Theorem 6.3.3] and [TV17, Example 1.6.2]. Taking the opposite of a rigid category reverses the roles of evaluation and coevaluation. Thus the opposite $\mathscr{C}^{\mathrm{op}}$ of the left rigid category $\mathscr{C}$ is right rigid.

**Lemma 2.64** ([EGNO15, Exercise 2.10.6]). *Let $F\colon \mathscr{C} \longrightarrow \mathscr{D}$ be a strong monoidal functor. The image $Fx$ of any (rigidly) dualisable object $x \in \mathscr{C}$ is dualisable.*

*Sketch of proof.* If $x \in \mathscr{C}$ is an object with right dual $x^\vee$, define $(Fx)^\vee := F(x^\vee)$; the evaluation map is given by

$$Fx \otimes F(x^\vee) \xrightarrow{F_2} F(x \otimes x^\vee) \xrightarrow{F\mathrm{ev}_x^{(r)}} F1 \xrightarrow{F_0^{-1}} 1,$$

and the coevaluation is defined analogously. A straightforward calculation shows the snake identities to be satisfied. □

Graphically, evaluations and coevaluations will be represented by semicircles, possibly decorated with arrows to emphasise whether we consider their left or right version.

| | | | |
|---|---|---|---|
| 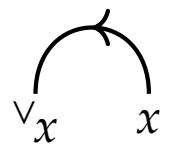 | 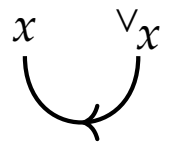 | 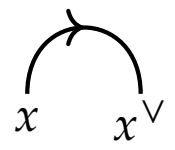 | 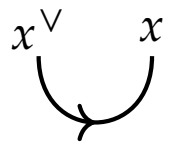 |
| $\mathrm{ev}_x^\ell\colon {}^\vee x \otimes x \longrightarrow 1$ | $\mathrm{coev}_x^\ell\colon 1 \longrightarrow x \otimes {}^\vee x$ | $\mathrm{ev}_x^r\colon x \otimes x^\vee \longrightarrow 1$ | $\mathrm{coev}_x^r\colon 1 \longrightarrow x^\vee \otimes x$ |

**Definition 2.65.** An object $x \in \mathscr{C}$ in a rigid monoidal category $\mathscr{C}$ is called *invertible* if its (left) evaluation and coevaluation are isomorphisms.

It follows that the right evaluations and coevaluations of an invertible objects are isomorphisms as well. Further, tensor products and duals of invertible objects are invertible, such that the full subcategory $\mathsf{Inv}(\mathscr{C})$ of invertible object of $\mathscr{C}$ is rigid monoidal.

**Definition 2.66** ([Cas05]). The *Picard group* $\mathrm{Pic}\,\mathscr{C}$ of a rigid monoidal category $\mathscr{C}$ is the group of isomorphism classes of invertible objects in $\mathscr{C}$. Its multiplication is induced by the tensor product; i.e., $[\alpha]\cdot[\beta] := [\alpha \otimes \beta]$ for $\alpha, \beta \in \mathsf{Inv}(\mathscr{C})$. The unit of $\mathrm{Pic}\,\mathscr{C}$ is $[1]$ and for any $\alpha \in \mathsf{Inv}(\mathscr{C})$ we have that $[\alpha]^{-1} = [{}^\vee\alpha]$.

**Proposition 2.67.** *For every object $x \in \mathscr{C}$ in a rigid category $\mathscr{C}$ we obtain two chains of adjoint endofunctors of $\mathscr{C}$:*

$$\ldots \dashv - \otimes x^\vee \dashv - \otimes x \dashv - \otimes {}^\vee x \dashv \ldots \qquad \ldots \dashv {}^\vee x \otimes - \dashv x \otimes - \dashv x^\vee \otimes - \dashv \ldots$$

*Further, $- \otimes x$ and $x \otimes -$ are equivalences of categories if and only if $x$ is invertible.*

*Proof.* The existence of the stated chains of adjunctions follows from [EGNO15, Proposition 2.10.8]. A straightforward calculation shows that tensoring (from the left or right) with an invertible object yields an equivalence of categories.

Conversely, suppose that $x \in \mathscr{C}$ is such that $F := - \otimes x$ is an equivalence of categories. The functor $F$ and its quasi-inverse $U$ are part of an adjunction with invertible unit $\eta\colon \mathrm{Id}_{\mathscr{C}} \Longrightarrow UF$ and counit $\varepsilon\colon FU \Longrightarrow \mathrm{Id}_{\mathscr{D}}$, see for example [Rie17, Proposition 4.4.5]. By [Rie17, Proposition 4.4.1], there exists a natural isomorphism $\theta\colon U \Longrightarrow - \otimes {}^{\vee}x$ that commutes with the respective counits and units. Applied to the monoidal unit $1 \in \mathscr{C}$, we obtain

$$\mathrm{coev}_x^{\ell} = \theta_x \circ \eta_1 \qquad \text{and} \qquad \mathrm{ev}_x^{\ell} \circ (\theta_1 \otimes \mathrm{id}_x) = \varepsilon_1.$$

It follows that $x$ is invertible. An analogous argument shows that $x \otimes -$ being an equivalence of categories also entails $x$ being invertible. □

**Remark 2.68.** Given a monoidal category $\mathscr{C}$, the mere presence of adjunctions

$$- \otimes x\colon \mathscr{C} \rightleftarrows \mathscr{C} : - \otimes Lx \qquad \text{and} \qquad x \otimes -\colon \mathscr{C} \rightleftarrows \mathscr{C} : Rx \otimes -$$

for objects $Lx, Rx \in \mathscr{C}$ does not lead to rigidity, but to the weaker notion of *tensor representability*, see Definitions 3.1 and 3.2.

To elucidate the underlying problem, let us assume for a moment that we are given objects $x, Rx \in \mathscr{C}$, such that $x \otimes -\colon \mathscr{C} \rightleftarrows \mathscr{C} : Rx \otimes -$. One can show that any right rigid dual of $x$—if it exists—has to be isomorphic to $Rx$; i.e., one has $Rx \cong x^{\vee}$ on objects. Hence, the unit $\eta_z^{(x)}\colon z \longrightarrow Rx \otimes x \otimes z$ and counit $\varepsilon_z^{(x)}\colon x \otimes Rx \otimes z \longrightarrow z$ of the adjunction provide us with natural candidates for the coevaluation and evaluation morphisms:

$$\mathrm{coev}_x := \eta_1^{(x)}\colon 1 \longrightarrow Rx \otimes x \qquad \text{and} \qquad \mathrm{ev}_x := \varepsilon_1^{(x)}\colon x \otimes Rx \longrightarrow 1.$$

Evaluating the triangle identities of the adjunction at the monoidal unit yields

$$\mathrm{id}_x = \varepsilon_x^{(x)} \circ (x \otimes \eta_1^{(x)}) \qquad \text{and} \qquad \mathrm{id}_{Rx} = (Rx \otimes \varepsilon_1^{(x)}) \circ \eta_{Rx}^{(x)}.$$

However, if $Rx$ is to be a dual of x in the rigid sense, the snake identities (2.4.4) must hold. For this, we ought to require the stronger condition that

$$\varepsilon_x^{(x)} = \varepsilon_1^{(x)} \otimes \mathrm{id}_x \qquad \text{and} \qquad \eta_{Rx}^{(x)} = \eta_1^{(x)} \otimes Rx.$$

**Remark 2.69.** A left closed monoidal such that for every $x \in \mathscr{C}$, the functor $-\otimes x$ admits a right dual $[x,-]_\ell$ as a $\mathscr{C}$-module functor—see Example 5.26—is already rigid. Explicitly, this involves the following map natural in $y, z \in \mathscr{C}$:

$$[x,z]_\ell \otimes y \xrightarrow{\sim} [x, z\otimes y]_\ell,$$

the resulting compatibility conditions of which imply the snake identities.

We will study these questions in greater details in Chapter 3; see also Remark 9.39 for a characterisation of rigidity in terms of $\mathscr{C}$-module functors.

**Definition 2.70.** A rigid monoidal category $\mathscr{C}$ is called *strict rigid* if the dual functors ${}^\vee(-)$, $(-)^\vee\colon \mathscr{C}^{\mathrm{op,rev}} \longrightarrow \mathscr{C}$ are strict and

$${}^\vee\big((-)^\vee\big) = \mathrm{Id}_{\mathscr{C}} = \big({}^\vee(-)\big)^\vee.$$

**Notation 2.71.** Let $\mathscr{C}$ be a rigid monoidal category; for $x \in \mathscr{C}$ and $n \in \mathbb{Z}$, write

$$(x)^n := \begin{cases} \text{The } n\text{-fold left dual of } x, & \text{if } n > 0; \\ x, & \text{if } n = 0; \\ \text{The } n\text{-fold right dual of } x, & \text{if } n < 0. \end{cases}$$

Our next result was conjectured in [Sch01, Section 5], and is a slight variation of [NS07, Theorem 2.2]. It shows that every rigid category admits a *rigid strictification*; i.e., a monoidally equivalent strict rigid category. The compatibility between the respective left and right duality functors is an immediate consequence of the fact that for any strong monoidal functor $F\colon \mathscr{C} \longrightarrow \mathscr{D}$ between rigid categories there are natural monoidal isomorphisms

$$\phi_x\colon F({}^\vee x) \xrightarrow{\sim} {}^\vee F x, \qquad \text{and} \qquad \psi_x\colon F(x^\vee) \xrightarrow{\sim} F x^\vee, \qquad \text{for all } x \in \mathscr{C}.$$

**Theorem 2.72.** *Every rigid category admits a rigid strictification.*

*Proof.* Suppose that $\mathscr{C}$ is a rigid strict monoidal category $\mathscr{C}$. Build a monoidally equivalent strict rigid category $\mathscr{D}$ as follows: the objects of $\mathscr{D}$ are (possibly empty) finite sequences $(x_1^{n_1}, \ldots, x_i^{n_i})$ of objects $x_1, \ldots, x_i \in \mathscr{C}$, adorned with integers $n_1, \ldots, n_i \in \mathbb{Z}$.

To define the morphisms of $\mathscr{D}$, recall Notation 2.71 and set

$$\mathscr{D}((x_1^{n_1}, \ldots, x_i^{n_i}), (y_1^{m_1}, \ldots, y_j^{m_j})) := \mathscr{C}((x_1)^{n_1} \otimes \cdots \otimes (x_i)^{n_i}, (y_1)^{m_1} \otimes \cdots \otimes (y_j)^{m_j}).$$

The category $\mathscr{D}$ is strict monoidal when equipped with the concatenation of sequences as tensor product and the empty sequence as unit. By construction, there exists a strict monoidal equivalence of categories $F\colon \mathscr{D} \longrightarrow \mathscr{C}$, which maps any object $(x_1^{n_1}, \dots, x_i^{n_i}) \in \mathscr{D}$ to $(x_1)^{n_1} \otimes \cdots \otimes (x_i)^{n_i} \in \mathscr{C}$, as well as every morphism to itself; the unit of $\mathscr{C}$ is the empty tensor product.

Define the left dual of an object $x := (x_1^{n_1}, \dots, x_i^{n_i}) \in \mathscr{D}$ as

$$^{\vee}x := (x_i^{n_i+1}, \dots, x_1^{n_1+1}),$$

with evaluation and coevaluation morphisms given by

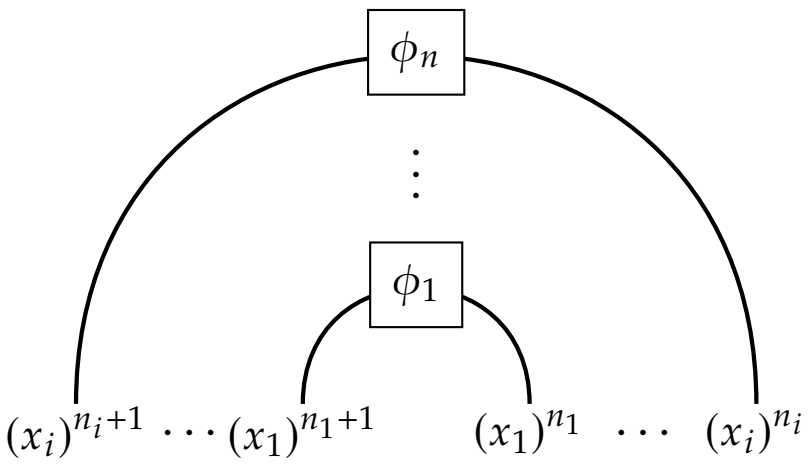

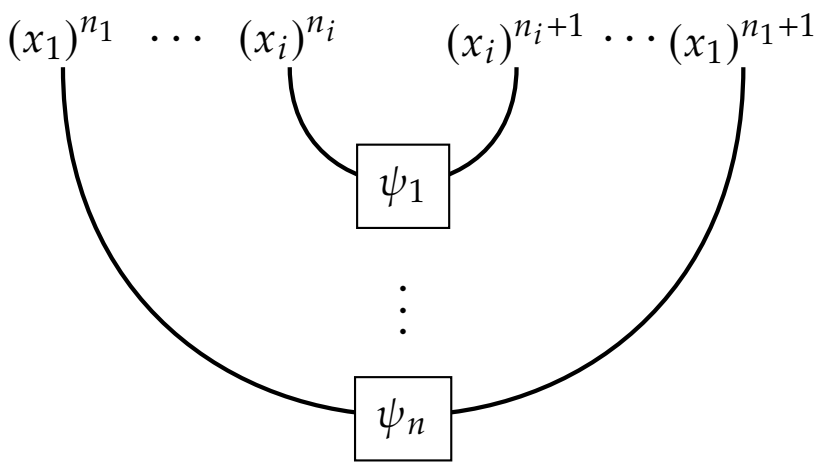

For all $1 \leq k \leq i$, we set

$$\phi_k := \begin{cases} \mathrm{ev}^{\ell}_{(x_k)^{n_k}}, & \text{if } n_k \geq 0; \\ \mathrm{ev}^{r}_{(x_k)^{n_k+1}}, & \text{if } n_k < 0; \end{cases} \qquad \text{and} \qquad \psi_k := \begin{cases} \mathrm{coev}^{\ell}_{(x_k)^{n_k}}, & \text{if } n_k \geq 0; \\ \mathrm{coev}^{r}_{(x_k)^{n_k+1}}, & \text{if } n_k < 0. \end{cases}$$

The right dual of $x$ if defined similarly as $x^{\vee} := (x_i^{n_i-1}, \dots, x_1^{n_1-1})$. Hence, $\mathscr{D}$ is strict rigid, which completes the proof. □

Module functors over rigid categories behave in particularly nice ways.

**Proposition 2.73** ([Ost03, Remark 4], [DSPS19, Lemma 2.10]). *If $\mathscr{C}$ is a left rigid monoidal, then every lax $\mathscr{C}$-module functor is strong. Dually, if $\mathscr{C}$ is a right rigid monoidal category, then every oplax $\mathscr{C}$-module functor is strong.*

Many applications require that the objects of a rigid monoidal category are isomorphic to their double duals in a way which is compatible with the monoidal structure. In Corollary 6.45 we prove a result for detecting this.

**Definition 2.74.** A *pivotal category* is a rigid monoidal category $\mathscr{C}$ equipped with a monoidal natural isomorphism $\rho\colon \mathrm{Id}_{\mathscr{C}} \xRightarrow{\sim} {}^{\vee\vee}-$, a *pivotal structure* of $\mathscr{C}$.

Rigid monoidal categories do not have to admit a pivotal structure and, if they do, it need not be unique. Examples coming from Hopf algebra theory are given in [KR93; HK19; Hal21]. However, every rigid monoidal category admits a universal pivotal category, called its *pivotal cover*, see [Shi15].

### 2.4.4 *The Drinfeld centre*

Classically, the centre construction is used to build a braided monoidal category from a monoidal one, see for example [EGNO15, Chapter 7]. Throughout especially Chapters 4 and 5, we work in a slightly more general setting, see for example [GNN09; BV12; HKS19; FH23; Kow24].

Analogously to Remark 2.53, one can show that $\sigma_{M,1} = \mathrm{id}_M$ is implied by the hexagon axiom.

**Definition 2.75.** Let $\mathcal{C}$ be a monoidal category, $\mathcal{M}$ a $\mathcal{C}$-bimodule category, and $m \in \mathcal{M}$. A *half-braiding* on $m$ is a natural isomorphism $\sigma_{-,=} \colon - \triangleleft = \Longrightarrow - \triangleright =$, such that for all $x, y \in \mathcal{C}$ we have $\sigma_{M,1} = \mathrm{id}_M$ and

$$\sigma_{m,x\otimes y} = (\mathrm{id}_x \triangleright \sigma_{m,y}) \circ (\sigma_{m,x} \triangleleft \mathrm{id}_y).$$

**Definition 2.76.** The *centre* of a $\mathcal{C}$-bimodule category $\mathcal{M}$ is the category $\mathsf{Z}(\mathcal{M})$ defined as follows:

- Objects are pairs $(m, \sigma_{m,-})$ of an object $m \in \mathcal{M}$ and a half-braiding $\sigma_{m,-}$.
- A morphism $f \colon (m, \sigma_{m,-}) \longrightarrow (n, \sigma_{n,-})$ consists of an $f \in \mathcal{M}(m,n)$ that commutes with the half-braidings:

$$(\mathrm{id}_x \triangleright f) \circ \sigma_{m,x} = \sigma_{n,x} \circ (f \triangleleft \mathrm{id}_x), \qquad \text{for all } x \in \mathcal{C}.$$

There is a canonical forgetful functor $U^{(M)} \colon \mathsf{Z}(\mathcal{M}) \longrightarrow \mathcal{M}$. Unlike classical representation theory where the centre of a bimodule is a subset of the bimodule, $U^{(M)}$ need not be injective on objects in general.

**Example 2.77.** The centre $\mathsf{Z}(\mathcal{C})$ of the regular bimodule of a monoidal category $\mathcal{C}$, see Example 2.44, is the *Drinfeld centre* of $\mathcal{C}$. The tensor product is defined by $(x, \sigma_{x,-}) \otimes (y, \sigma_{y,-}) := (x \otimes y, \sigma_{x\otimes y,-})$, with

$$\sigma_{x\otimes y,z} := (\sigma_{x,z} \otimes \mathrm{id}_y) \circ (\mathrm{id}_x \otimes \sigma_{y,z}), \qquad \text{for all } z \in \mathcal{C}.$$

The centre is braided monoidal, with braiding given by gluing together the respective half-braidings. The hexagon axioms follow from the definition of the half-braiding and the tensor product of $\mathsf{Z}(\mathcal{C})$.

**Example 2.78.** Let $H \in \mathsf{Vect}$ be a finite-dimensional Hopf algebra. Then the category of left-left Yetter–Drinfeld modules of Examples 2.33 and 2.54 is, as a braided monoidal category, equivalent to the Drinfeld centre of the category ${}_H\mathsf{Vect}$ of left $H$-modules, [Dri87; Yet90].

Further, by [Dri87] there exists another Hopf algebra $D(H)$, the *Drinfeld double* of $H$, that has $\mathsf{Z}({}_H\mathsf{Vect}) \simeq {}^H_H\mathcal{YD}$ as its category of modules. We refer to [Kas98, Definition IX.4.1] for a complete definition, and to [Kas98, Theorem XIII.5.1] for a proof of the above equivalence.

**Remark 2.79.** The category of anti-Yetter–Drinfeld modules of Example 2.42 is also the category of modules over a certain Hopf algebra: the *anti-Drinfeld double* [CMZ97; Sch99]. In analogy to how the anti-Yetter–Drinfeld modules are a module category over the Yetter–Drinfeld modules, the anti-Drinfeld double is a comodule algebra over the Drinfeld double.

For the next result, recall Notation 2.71 for iterated duals.

**Proposition 2.80** ([JS91, Lemma 7]). *The centre of a (strict) rigid category* $\mathscr{C}$ *is (strict) rigid: for all* $(x, \sigma_{x,-}) \in \mathsf{Z}(\mathscr{C})$, *we have that* $U^{(\mathsf{Z})}\big({}^{\vee}(x, \sigma_{x,-})\big) = {}^{\vee}x$ *and* $U^{(\mathsf{Z})}\big((x, \sigma_{x,-})^{\vee}\big) = x^{\vee}$. *Moreover, for every* $n \in \mathbb{Z}$ *and* $x \in \mathsf{Z}(\mathscr{C})$, *the equality* $\sigma_{(x)^n,(y)^n} = (\sigma_{x,y})^n$ *holds for all* $y \in \mathsf{Z}(\mathscr{C})$.

## 2.5 Linear and abelian categories

We shall pay special attention to $\Bbbk$-linear categories, for some field $\Bbbk$. Indeed, $\Bbbk$-linear functor categories encompass many phenomena in representation theory; see for example Section 3.2 and Chapter 9. We shall denote the category $\Bbbk$-vector spaces by $\mathsf{Vect}$, and write $\mathsf{vect}$ for the full subcategory of finite-dimensional $\Bbbk$-vector spaces.[3]

[3] This will be a general theme throughout the thesis: whenever there exists a full subcategory of "finite-dimensional objects" inside of a lager category of all objects, the finite-dimensional subcategory will start with a lower case letter, and we use a capital for the larger category.

**Definition 2.81.** A $\Bbbk$-*linear category* is a category enriched in $\mathsf{Vect}$

**Example 2.82.** Let $\mathscr{C}$ be an ordinary category. Its *linearisation* $\Bbbk\mathscr{C}$ has the same objects, and the space of morphisms between two objects $a, b \in \Bbbk\mathscr{C}$ is $\Bbbk\mathscr{C}(a, b) := \operatorname{span}_{\Bbbk}\mathscr{C}(a, b)$. To define the composition in $\Bbbk\mathscr{C}$, note that for any $a, b, c \in \Bbbk\mathscr{C}$ there is a unique linear map

$$- \circ =: \Bbbk\mathscr{C}(b, c) \otimes_{\Bbbk} \Bbbk\mathscr{C}(a, b) \longrightarrow \Bbbk\mathscr{C}(a, c),$$

such that $g \circ f = gf$ for all $g \in \mathscr{C}(b, c)$ and $f \in \mathscr{C}(a, b)$.

Let $\iota\colon \mathscr{C} \longrightarrow \Bbbk\mathscr{C}$ be the functor that is the identity on objects and maps any morphism to the corresponding basis vector. For any functor $F\colon \mathscr{C} \longrightarrow \mathscr{D}$ whose codomain is a $\Bbbk$-linear category, we obtain a commuting triangle:

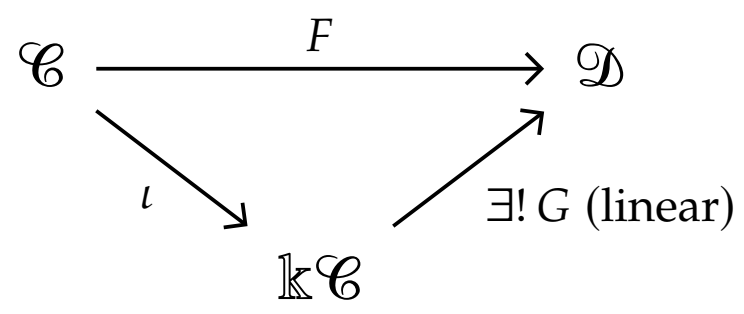

As a consequence, if we endow the ordinary functor category $[\mathscr{C}, \mathsf{Vect}]$ with the pointwise $\Bbbk$-linear structure, we obtain an isomorphism of linear categories between it and the category of $\Bbbk$-linear functors from $\Bbbk\mathscr{C}$ to $\mathsf{Vect}$.

### 2.5.1 *Finite and locally finite abelian categories*

In this section, we assume all categories and functors to be $\Bbbk$-linear, for some field $\Bbbk$. We call a functor $F\colon \mathcal{A} \longrightarrow \mathcal{B}$ between abelian categories $\mathcal{A}$ and $\mathcal{B}$ *left exact* if it preserves finite limits, *right exact* if it preserves finite colimits, and *exact* if it is right and left exact. The category of right exact functors from $\mathcal{A}$ to $\mathcal{B}$ will be denoted by $\mathsf{Rex}(\mathcal{A}, \mathcal{B})$, and we write $\mathsf{Lex}(\mathcal{A}, \mathcal{B})$ for the category of left exact functors. Further, let $\mathcal{A}$-proj denote the full subcategory of $\mathcal{A}$ consisting of its projective objects, and similarly for $\mathcal{A}$-inj.

We recall basic properties of projective and injective objects in this setting.

**Proposition 2.83** ([GKKP22, Lemma 3.3])**.** *Let $\mathcal{A}$ and $\mathcal{B}$ be abelian categories, and assume that $\mathcal{B}$ has enough projectives. Let $\iota\colon \mathcal{A}\text{-proj} \hookrightarrow \mathcal{A}$ be the inclusion functor. Then the restriction functor $\mathsf{Rex}(\mathcal{A}, \mathcal{B}) \xrightarrow{-\circ\iota} [\mathcal{A}\text{-proj}, \mathcal{B}]$ is an equivalence.*

**Proposition 2.84** ([GKKP22, Lemma 3.3])**.** *Let $\mathcal{A}$ and $\mathcal{B}$ be abelian categories, and assume that $\mathcal{B}$ has enough injectives. Let $\iota\colon \mathcal{A}\text{-inj} \hookrightarrow \mathcal{A}$ be the inclusion functor. Then $\mathsf{Lex}(\mathcal{A}, \mathcal{B}) \xrightarrow{-\circ\iota} [\mathcal{A}\text{-inj}, \mathcal{B}]$ is an equivalence.*

We refer to [EGNO15, Chapter 1] for the following definitions and results. An object of a $\Bbbk$-linear category is said to be *simple* if any non-zero endomorphism thereof is an isomorphism. A category $\mathcal{C}$ is said to be *semisimple* if any of its objects is a direct sum of simple objects.

**Proposition 2.85** ([TV17, Theorem C.6])**.** *A semisimple category is abelian. Any epimorphism or monomorphism in a semisimple category splits. In particular, any functor from a semisimple category to an abelian category is exact.*

An object $V$ in an abelian category $\mathcal{A}$ has *finite length* if there is a filtration

$$0 = V_0 \subseteq V_1 \subseteq \cdots \subseteq V_n = V,$$

such that $V_i/V_{i-1}$ is simple, for all $i \in \{1, \ldots, n\}$.

**Definition 2.86.** A category $\mathcal{C}$ is said to be *finite abelian* if it is abelian, hom-finite[4], has enough projectives, only finitely many isomorphism classes of simple objects, and all of its objects are of finite length.

[4] A ($\Bbbk$-linear) category $\mathcal{C}$ is called *hom-finite* if the hom space $\mathcal{C}(x, y)$ is finite-dimensional for all $x, y \in \mathcal{C}$.

**Lemma 2.87.** *A category $\mathcal{A}$ is finite abelian if and only if there is a finite-dimensional $\Bbbk$-algebra $A$ such that $\mathcal{A} \simeq A$-mod.*

**Definition 2.88.** A category is said to be *locally finite abelian* if it is abelian, hom-finite, and all of its objects are of finite length.

**Lemma 2.89.** *A category $\mathscr{C}$ is locally finite abelian if and only if there exists a $\mathbb{k}$-coalgebra $C$ such that $\mathscr{C} \simeq {}^{C}\mathsf{vect}$.*

Locally finite abelian categories satisfy a Krull–Schmidt-type theorem.

**Proposition 2.90** ([EGNO15, Theorem 1.5.7])**.** *Every object of finite length admits a unique decomposition into a direct sum of indecomposable objects.*

**Definition 2.91.** An abelian category $\mathscr{C}$ that satisfies Proposition 2.90 is called a *Krull–Schmidt category*.

General abelian monoidal categories admit a very useful special case of *Beck's monadicity theorem*, see for example [BZBJ18, Section 4.1].

**Theorem 2.92.** *Let $F\colon \mathscr{A} \rightleftarrows \mathscr{B} : U$ be an adjunction of abelian categories. If $U$ is right exact and reflects zero objects[5], the comparison functor $K\colon \mathscr{B} \longrightarrow \mathscr{A}^{UF}$ is an equivalence. Likewise, if $F$ is left exact and reflects zero objects then the comparison functor $K\colon \mathscr{A} \longrightarrow \mathscr{B}^{FU}$ is an equivalence.*

[5] A functor $F$ *reflects zero objects* if $Fx \cong 0 \implies x \cong 0$.

### 2.5.2 *Tensor and ring categories*

We refer the reader to [EGNO15, Chapter 4] for a comprehensive account of the theory of (multi)tensor and ring categories. All categories and functors in this section are again assumed to be $\mathbb{k}$-linear.

**Definition 2.93.**

- A *tensor category* is a locally finite abelian rigid monoidal category $\mathscr{C}$ such that $\mathrm{End}_{\mathscr{C}}(\mathbb{1}) \cong \mathbb{k}$.
- A *finite tensor category* is a finite abelian rigid monoidal category $\mathscr{C}$ such that $\mathrm{End}_{\mathscr{C}}(\mathbb{1}) \cong \mathbb{k}$.
- A *ring category* is a locally finite abelian separately exact monoidal category $\mathscr{C}$ such that $\mathrm{End}_{\mathscr{C}}(\mathbb{1}) \cong \mathbb{k}$.
- A *finite ring category* is a finite abelian separately exact monoidal category $\mathscr{C}$ such that $\mathrm{End}_{\mathscr{C}}(\mathbb{1}) \cong \mathbb{k}$.

Despite apparent similarity in terminology, the notion of a ring category is not related to that of a rig category, such as those considered in [JY21].

**Definition 2.94.** Let $\mathscr{C}$ be a ring category. A *fibre functor* for $\mathscr{C}$ is a faithful and exact monoidal functor $U\colon \mathscr{C} \longrightarrow \mathsf{vect}$.

The next result is a variant of Tannaka–Krein duality for ring categories.

**Theorem 2.95** ([EGNO15, Theorem 5.4.1]). *Let $\mathscr{C}$ be a ring category, and assume it admits a fibre functor $U\colon \mathscr{C} \longrightarrow \mathsf{vect}$. Then there exists a bialgebra $B$ such that there is a monoidal equivalence $K_B\colon \mathscr{C} \longrightarrow B\text{-comod}$ and a monoidal natural isomorphism $U_B \circ K_B \cong U$, where $U_B\colon B\text{-comod} \longrightarrow \mathsf{vect}$ is the forgetful functor. Further, $B$ is unique up to isomorphism, Hopf if and only if $\mathscr{C}$ is a tensor category, and finite-dimensional if and only if $\mathscr{C}$ is a finite ring category.*

## 2.6 Algebra and module objects

Let $\mathscr{C}$ be a monoidal category. An *algebra object in* $\mathscr{C}$ consists of an object $A \in \mathscr{C}$ together with a *multiplication* $\mu\colon A \otimes A \longrightarrow A$ and a *unit* $\eta\colon \mathbb{1} \longrightarrow A$, satisfying associativity and unitality axioms analogous to those for a $\Bbbk$-algebra. In fact, a $\Bbbk$-algebra is the same as an algebra object in $\mathsf{Vect}$.

A *coalgebra object* in $\mathscr{C}$ is an object $C \in \mathscr{C}$ equipped with a comultiplication $\Delta\colon C \longrightarrow C \otimes C$ and a counit morphism $\varepsilon\colon C \longrightarrow \mathbb{1}$, satisfying coassociativity and counitality axioms, such that $C$ becomes an algebra object in $\mathscr{C}^{\mathrm{op}}$.

**Example 2.96.** An object $A$ in a monoidal category $\mathscr{C}$ is an algebra object in $\mathscr{C}$ if and only if $- \otimes A\colon \mathscr{C} \longrightarrow \mathscr{C}$ is a monad.

**Example 2.97.** Let $\mathscr{C}$ be a monoidal category, and suppose that $x \in \mathscr{C}$ has a left dual ${}^\vee x \in \mathscr{C}$. Then this induces an algebra object

$$(x \otimes {}^\vee x,\ x \otimes \mathrm{ev}_x^\ell \otimes {}^\vee x\colon x \otimes {}^\vee x \otimes x \otimes {}^\vee x \longrightarrow x \otimes {}^\vee x,\ \mathrm{coev}_x^\ell\colon 1 \longrightarrow x \otimes {}^\vee x)$$

and, dually, a coalgebra object $({}^\vee x \otimes x,\ {}^\vee x \otimes \mathrm{coev}_x^\ell \otimes X,\ \mathrm{ev}_x^\ell)$ in $\mathscr{C}$.

A special case of Example 2.97 is when we consider the monoidal category $[\mathscr{C}, \mathscr{C}]$ of endofunctors. Then a left dual of some $U \in [\mathscr{C}, \mathscr{C}]$ is a left adjoint, with $\mathrm{ev}_U^\ell$ and $\mathrm{coev}_U^\ell$ giving the coevaluation and evaluation, respectively. the resulting algebra and coalgebra structures are the monad and comonad corresponding to the adjunction.

**Definition 2.98.** Let $(\mathcal{M}, \triangleright)$ be a left $\mathcal{C}$-module category and $(A, \mu, \eta)$ an algebra object in $\mathcal{C}$. A *left A-module in* $\mathcal{M}$ is an object $M \in \mathcal{M}$ together with an *action* morphism $\alpha \colon A \triangleright M \longrightarrow M$ such that the following diagrams commute

$$\begin{array}{ccccc} (A \otimes A) \triangleright M & \xrightarrow{\mathcal{M}_{\mathrm{a}}} & A \triangleright (A \triangleright M) & \xrightarrow{A \triangleright \alpha} & A \triangleright M \\ \downarrow{\scriptstyle \mu \triangleright m} & & & & \downarrow{\scriptstyle \alpha} \\ A \triangleright M & & \xrightarrow{\quad\alpha\quad} & & M \end{array} \qquad \begin{array}{ccc} \mathbb{1} \triangleright M & \xrightarrow{\mathcal{M}_{\mathrm{u}}} & M \\ \downarrow{\scriptstyle \eta \triangleright M} & \nearrow{\scriptstyle \alpha} & \\ A \triangleright M & & \end{array}$$

A *morphism of modules* is an arrow $f \colon M \longrightarrow N$ in $\mathcal{M}$ commuting with the respective actions in the sense that $f \circ \alpha_m = \alpha_n \circ (A \triangleright f)$.

**Example 2.99.** Clearly, left $A$-modules in $\mathcal{M}$ and their morphisms form a category, which we shall denote by $A\text{-}\mathsf{Mod}_{\mathcal{M}}$. The obvious forgetful functor $A\text{-Forget}_{\mathcal{M}} \colon A\text{-}\mathsf{Mod}_{\mathcal{M}} \longrightarrow \mathcal{M}$ admits a left adjoint

$$\begin{aligned} A\text{-Free}_{\mathcal{M}} \colon \mathcal{M} &\longrightarrow A\text{-}\mathsf{Mod}_{\mathcal{M}} \\ M &\longmapsto \left(A \triangleright M,\ A \triangleright (A \triangleright M) \xrightarrow{\mathcal{M}_{\mathrm{a}}^{-1}} (A \otimes A) \triangleright M \xrightarrow{\mu \triangleright M} A \triangleright M\right) \end{aligned}$$

The unit of this adjunction is given by $M \xrightarrow{\mathcal{M}_{\mathrm{u}}} \mathbb{1} \triangleright M \xrightarrow{\eta \triangleright M} A \triangleright M$, and the counit is $\alpha_M \colon A \triangleright M \longrightarrow M$.

Given a coalgebra object $C$ in $\mathcal{C}$, a *left C-comodule in* $\mathcal{M}$ is an object $N \in \mathcal{C}$ satisfying axioms dual to those for a module in $\mathcal{M}$, so that a left $C$-comodule in $\mathcal{M}$ is the same as a left $C$-module in the $\mathcal{C}^{\mathrm{op}}$-module category $\mathcal{M}^{\mathrm{op}}$. Comodule morphisms are defined similarly. We denote the category of $C$-comodules in $\mathcal{M}$ by $C\text{-}\mathsf{Comod}_{\mathcal{M}}$. There is an adjunction analogous to that for modules:

$$C\text{-Coforget}_{\mathcal{M}} \colon \mathcal{M} \rightleftarrows C\text{-}\mathsf{Comod}_{\mathcal{M}} \colon C\text{-Cofree}_{\mathcal{M}}.$$

Similarly, one can define right modules and comodules in a right module category, and bimodules and bicomodules in a bimodule category.

**Definition 2.100.** Let $A$ be an algebra in $\mathcal{C}$. A (left) *A-comodule algebra* is an algebra $C$, together with a left $A$-comodule structure $\rho \colon C \longrightarrow A \otimes C$, such that $\rho$ is a morphism of $A$-algebras.

Analogously to Definition 2.100 one defines right comodule algebras, as well as comodule coalgebras and module algebras.

**Notation 2.101.** If the underlying category is the category of vector spaces, then we use the following notation, to differentiate it from the general case:

$$^{C}\mathsf{Vect} := C\text{-}\mathsf{Comod}_{\mathsf{Vect}}.$$

for the category of left $C$-comodules of a coalgebra $C$.

Analogously, we for example define $_{A}\mathsf{Vect}$, $^{C}\mathsf{Vect}_{A}$; or $^{C}\mathsf{vect}$, $_{A}\mathsf{vect}$, and $^{B}_{B}\mathsf{vect}^{B}_{B}$ for the finite-dimensional case.

**Remark 2.102.** Let $A$ be an algebra object in $\mathscr{C}$ and let $\mathscr{M}$ and $\mathscr{N}$ be left $\mathscr{C}$-module categories. A lax $\mathscr{C}$-module functor $F\colon \mathscr{N} \longrightarrow \mathscr{M}$ can be lifted to

$$A\text{-}\mathsf{Mod}_F\colon A\text{-}\mathsf{Mod}_{\mathscr{N}} \longrightarrow A\text{-}\mathsf{Mod}_{\mathscr{M}}$$

on the respective categories of $A$-modules by sending $N \in A\text{-}\mathsf{Mod}_{\mathscr{N}}$ to the left $A$-module $FN$, whose action is given by

$$\alpha_{FN}\colon A \triangleright_{\mathscr{M}} FN \xrightarrow{F_{\mathrm{a};A,N}} F(A \triangleright_{\mathscr{N}} N) \xrightarrow{F\alpha_N} FN.$$

This functor is a lift of $F$—it satisfies the following relation:

$$\begin{array}{ccc} A\text{-}\mathsf{Mod}_{\mathscr{N}} & \xrightarrow{A\text{-}\mathsf{Mod}_F} & A\text{-}\mathsf{Mod}_{\mathscr{M}} \\ {\scriptstyle A\text{-}\mathrm{Forget}_{\mathscr{N}}}\downarrow & & \downarrow{\scriptstyle A\text{-}\mathrm{Forget}_{\mathscr{M}}} \\ \mathscr{N} & \xrightarrow[F]{} & \mathscr{M} \end{array}$$

**Example 2.103.** Let $\mathscr{C}$ be a cocomplete symmetric monoidal category, such that its is right exact in both variables. Then one may define the bicategory $\mathbb{B}\mathsf{imod}(\mathscr{C})$ of bimodules in $\mathscr{C}$ as follows:

- Objects are algebra objects in $\mathscr{C}$.
- For $A$ and $B$ in $\mathscr{C}$, the hom-category $\mathbb{B}\mathsf{imod}(\mathscr{C})(A, B)$ is given by $B$-$A$-bimodules in $\mathscr{C}$ and their homomorphisms. Horizontal composition is given by the balanced tensor product of bimodules.

Analogously to Example 2.103, there is a bicategory $\mathbb{B}\mathsf{icomod}(\mathscr{C})$ of $\mathscr{C}$-bicomodules, for $\mathscr{C}$ a complete symmetric monoidal category, with the tensor product being left exact in both variables.

**Definition 2.104.** A category $\mathscr{I}$ is *filtered* if every finite diagram has a cocone.

Equivalently, Definition 2.104 could be stated in the following way: a category $\mathcal{I}$ is filtered if the following two conditions hold:

- for all $x, y \in \mathcal{I}$ there exist arrows $k_x \colon x \longrightarrow k$ and $k_y \colon y \longrightarrow k$ in $\mathcal{I}$; and
- for all parallel arrows $f, g \colon x \longrightarrow y$ in $\mathcal{I}$ there exists a $k_{fg} \colon y \longrightarrow k$, such that $k_{fg} \circ g = k_{fg} \circ f$.

**Definition 2.105.** A diagram $F \colon \mathcal{I} \longrightarrow \mathcal{C}$ is called *filtered* if the domain $\mathcal{I}$ is a filtered category. A *filtered colimit* is a colimit of a filtered diagram, and a functor is called *finitary* if it preserves filtered colimits.

Ordinarily, limits and colimits do not commute. For example, for sets $X, Y$, $W$, and $V$, the canonical morphism $(X \times Y) \coprod (W \times V) \longrightarrow (X \coprod V) \times (Y \coprod V)$ is not isomorphism in general. However, it turns out that filtered colimits do commute with finite limits, see for example [ML98, Theorem IX.2.1].

**Proposition 2.106.** *Filtered colimits commute with finite limits in* $\mathsf{Set}$: *given a functor* $F \colon \mathcal{J} \times \mathcal{I} \longrightarrow \mathsf{Set}$ *of two diagram categories, where* $\mathcal{J}$ *is finite and* $\mathcal{I}$ *is small and filtered, the following canonical morphism is an isomorphism*:

$$\operatorname*{colim}_i \lim_j F(j, i) \xrightarrow{\sim} \lim_j \operatorname*{colim}_i F(j, i).$$

**Lemma 2.107.** *Let* $\mathcal{A}$ *and* $\mathcal{B}$ *be additive categories. Finite biproducts endow* $\mathcal{A}$ *and* $\mathcal{B}$ *with the structure of* $\mathsf{vect}$*-module categories. Any functor* $F \colon \mathcal{A} \longrightarrow \mathcal{B}$ *is a* $\mathsf{vect}$*-module functor with respect to these* $\mathsf{vect}$*-actions.*

*If* $\mathcal{B}$ *admits filtered colimits, F extends essentially uniquely to a finitary* $\mathsf{Vect}$*-module functor* $\mathsf{Ind}(\mathcal{A}) \longrightarrow \mathcal{B}$*, and any finitary* $\mathsf{Vect}$*-module functor is of this form.*

*Proof.* Let $\mathsf{add}(\mathcal{A})$ denote the cocompletion of $\mathcal{A}$ under finite biproducts—the *additive closure* of $\mathcal{A}$. As finite biproducts are absolute colimits, the left adjoint $\mathsf{add}(\mathcal{A}) \longrightarrow \mathcal{A}$ of the canonical inclusion $\mathcal{A} \hookrightarrow \mathsf{add}(\mathcal{A})$, coming from the universal property of Equation (2.9.1), is an equivalence of $\mathsf{vect}$-module categories. Hence, so is the inclusion $\mathcal{A} \hookrightarrow \mathsf{add}(\mathcal{A})$. We find the equivalences

$$\mathsf{Cat}_{\Bbbk}(\mathcal{A}, -) \simeq \mathsf{vect}\text{-}\mathsf{Mod}(\mathsf{add}(\mathcal{A}), -) \simeq \mathsf{vect}\text{-}\mathsf{Mod}(\mathcal{A}, -),$$

establishing the first claim. The second claim follows in a similar way, using that $\mathsf{Vect} \simeq \mathsf{Ind}(\mathsf{vect})$ and

$$\mathsf{vect}\text{-}\mathsf{Mod}(\mathcal{A}, \mathcal{B}) \simeq \mathsf{Vect}\text{-}\mathsf{Mod}_{\mathrm{filt}}(\mathsf{Ind}(\mathcal{A}), \mathcal{B}). \qquad \square$$

**Example 2.108.** Let $C$ be a coalgebra in $\mathsf{Vect}$. Given two $C$-comodules $(X, \rho)$ and $(Y, \lambda)$, one may form their *cotensor product*—the vector space $X \,\square\, Y$ given by the equaliser

$$X \,\square\, Y \hookrightarrow X \otimes_{\Bbbk} Y \underset{X \otimes_{\Bbbk} \lambda}{\overset{\rho \otimes_{\Bbbk} Y}{\rightrightarrows}} X \otimes_{\Bbbk} C \otimes_{\Bbbk} Y.$$

Sometimes, the formalism of module categories—even those over the seemingly trivial monoidal categories $\mathsf{vect}$ and $\mathsf{Vect}$—provides additional clarity to classical statements about algebraic and coalgebraic $\Bbbk$-linear structures. We give a brief proof of the result below as an example of this.

**Proposition 2.109** ([Tak77, Proposition 2.1])**.** *Let $\mathfrak{D}$ be an abelian category admitting filtered colimits, and let $C$ be a coalgebra over $\Bbbk$. There is an equivalence*

$$\begin{aligned} \mathsf{Lexf}({}^{C}\mathsf{Vect}, \mathfrak{D}) &\xrightarrow{\sim} \mathsf{Comod}_{\mathfrak{D}} C \\ F &\longmapsto FC \\ N \,\square\, - &\longleftarrow N, \end{aligned}$$

*where the left-hand side is the category of left exact finitary functors from ${}^{C}\mathsf{Vect}$, and the right-hand side—following the notation of Example 2.99—denotes the category of right $C$-comodules in $\mathfrak{D}$. Here, $\mathfrak{D}$ is endowed with the $\mathsf{Vect}$-module structure described in Proposition 2.131.*

*In particular, for a coalgebra $D$ over $\Bbbk$, we have*

$$\mathsf{Lexf}({}^{C}\mathsf{Vect}, {}^{D}\mathsf{Vect}) \simeq {}^{D}\mathsf{Vect}^{C}.$$

*Proof.* Since ${}^{C}\mathsf{Vect} \simeq \mathsf{Ind}({}^{C}\mathsf{vect})$, the functor $F$ can be seen as induced to filtered colimits from its restriction to ${}^{C}\mathsf{vect}$, and as such is a $\mathsf{Vect}$-module functor, following Lemma 2.107. Observe that $C$, being a bicomodule over itself, is a right $C$-comodule in ${}^{C}\mathsf{Vect}$. Thus, $FC$ is a right $C$-comodule in $\mathfrak{D}$.

To see that $F \cong FC \,\square\, -$, let $N \in {}^{C}\mathsf{Vect}$ and recall that the bar construction provides a functorial injective resolution

$$N \cong \mathrm{eq}\big(C \otimes N \underset{C \otimes \Delta_N}{\overset{\Delta_C \otimes N}{\rightrightarrows}} C \otimes C \otimes N\big).$$

Thus,

$$\begin{aligned} FN &\cong F\big(\mathrm{eq}(C \otimes N \underset{C \otimes \Delta_N}{\overset{\Delta_C \otimes N}{\rightrightarrows}} C \otimes C \otimes N)\big) \cong \mathrm{eq}\left(F(C \otimes N) \underset{F(C \otimes \Delta_N)}{\overset{F(\Delta_C \otimes N)}{\rightrightarrows}} F(C \otimes C \otimes N)\right) \\ &\cong \mathrm{eq}\left(FC \otimes N \underset{FC \otimes \Delta_N}{\overset{F(\Delta_C) \otimes N}{\rightrightarrows}} FC \otimes C \otimes N\right) = FC \,\square\, N, \end{aligned}$$

where the first isomorphism is using the bar construction, the second follows from left exactness of $F$, the third uses the fact that $F$ is a right Vect-module functor, and the fourth follows by observing that $F(\Delta_C) = \Delta_{F(C)}$. □

## 2.7 MONOIDAL BICATEGORIES

WE BRIEFLY RECALL THE NOTION OF A MONOIDAL BICATEGORY, which will play an important role in the string diagrams used throughout. More detailed accounts are given [JY21, Chatper 12], see also [GPS95; Gur06; GS16].

**Definition 2.110.** A *monoidal bicategory* consists of a bicategory $\mathbb{B}$ together with the following data:

- A pseudofunctor $\boxtimes\colon \mathbb{B} \times \mathbb{B} \longrightarrow \mathbb{B}$.
- A pseudofunctor $\mathbb{1}_{\mathbb{B}}\colon \heartsuit \longrightarrow \mathbb{B}$. The image of the unique object $\heartsuit$ in $\mathbb{B}$ will be referred to as the *identity object* of $\mathbb{B}$, and, abusing notation, also denoted by $\mathbb{1}_{\mathbb{B}}$. The pseudofunctoriality of $\mathbb{1}_{\mathbb{B}}$ yields additional structure: a 1-morphism $\mathrm{I}_{\mathbb{B}}\colon \mathbb{1}_{\mathbb{B}} \longrightarrow \mathbb{1}_{\mathbb{B}}$ and invertible 2-morphisms $\nabla^{\mathrm{I}}\colon \mathrm{I}_{\mathbb{B}} \circ \mathrm{I}_{\mathbb{B}} \xrightarrow{\sim} \mathrm{I}_{\mathbb{B}}$ and $\eta^{\mathrm{I}}\colon \mathrm{Id}_{\mathbb{1}_{\mathbb{B}}} \xrightarrow{\sim} \mathrm{I}_{\mathbb{B}}$, such that $(\mathbb{1}_{\mathbb{B}}, \nabla^{\mathrm{I}}, \eta^{\mathrm{I}})$ form an idempotent monad in $\mathbb{B}$.
- An adjoint pseudonatural equivalence $(\mathfrak{A}, \mathfrak{A}^{\blacklozenge}, \eta^{\mathfrak{A}}, \varepsilon^{\mathfrak{A}})$:

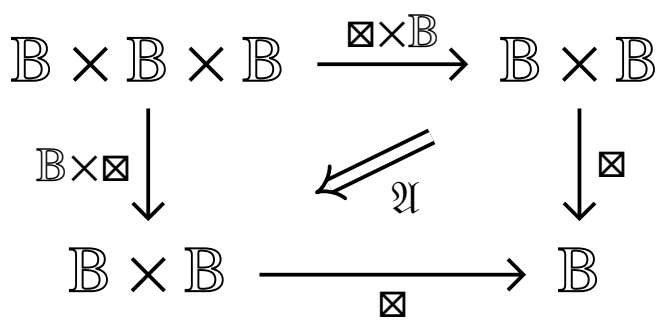

  which will be referred to as the *associator* for $\mathbb{B}$. Its components are 1-morphisms of the form $\mathfrak{A}_{\mathcal{A},\mathcal{B},\mathcal{C}}\colon (\mathcal{A} \boxtimes \mathcal{B}) \boxtimes \mathcal{C} \xrightarrow{\sim} \mathcal{A} \boxtimes (\mathcal{B} \boxtimes \mathcal{C})$.
- Pseudonatural adjoint equivalences $(\mathfrak{L}, \mathfrak{L}^{\blacklozenge}, \eta^{\mathfrak{L}}, \varepsilon^{\mathfrak{L}})$ and $(\mathfrak{R}, \mathfrak{R}^{\blacklozenge}, \eta^{\mathfrak{R}}, \varepsilon^{\mathfrak{R}})$, referred to as *left and right unitors*, respectively:

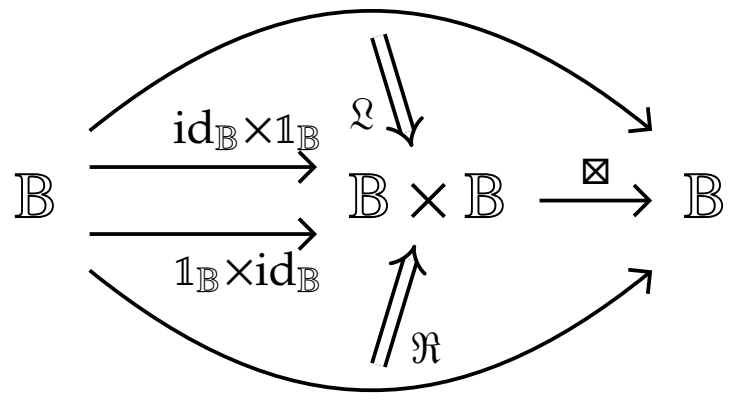

  Its components are 1-morphisms of $\mathbb{B}$ of the form $\mathfrak{R}\colon \mathbb{1}_{\mathbb{B}} \boxtimes \mathcal{A} \xrightarrow{\sim} \mathcal{A}$ and $\mathfrak{L}\colon \mathcal{A} \boxtimes \mathbb{1}_{\mathbb{B}} \xrightarrow{\sim} \mathcal{A}$, for $\mathcal{A} \in \mathrm{Ob}\,\mathbb{B}$.

- An invertible modification $\mathfrak{p}$, as indicated by the components of $\mathfrak{p}$ displayed below, where $\mathscr{A}, \mathscr{B}, \mathscr{C}, \mathscr{D} \in \operatorname{Ob} \mathbb{B}$:

$$
\begin{array}{ccccc}
(\mathscr{D} \boxtimes (\mathscr{C} \boxtimes \mathscr{B})) \boxtimes \mathscr{A} & & \longrightarrow & & \mathscr{D} \boxtimes ((\mathscr{C} \boxtimes \mathscr{B}) \boxtimes \mathscr{A}) \\
\uparrow & & \Big\Downarrow \mathfrak{p}_{\mathscr{D},\mathscr{C},\mathscr{B},\mathscr{A}} & & \downarrow \\
((\mathscr{D} \boxtimes \mathscr{C}) \boxtimes \mathscr{B}) \boxtimes \mathscr{A} & & & & \mathscr{D} \boxtimes (\mathscr{C} \boxtimes (\mathscr{B} \boxtimes \mathscr{A})) \\
& \searrow & & \swarrow & \\
& & (\mathscr{D} \boxtimes \mathscr{C}) \boxtimes (\mathscr{B} \boxtimes \mathscr{A}) & &
\end{array}
\tag{2.7.1}
$$

The arrows in the diagram are formed using the associator $\mathfrak{A}$ as indicated. We refer to $\mathfrak{p}$ as the *pentagonator*.

- Invertible modifications $\mathfrak{l}$, $\mathfrak{m}$, and $\mathfrak{r}$; the *left, middle, and right 2-unitors*:

$$
\begin{array}{ccc}
(\mathscr{B} \boxtimes \mathbb{1}_{\mathbb{B}}) \boxtimes \mathscr{A} & \xrightarrow{\mathfrak{A}_{\mathscr{B},\mathbb{1}_{\mathbb{B}},\mathscr{A}}} & \mathscr{B} \boxtimes (\mathbb{1}_{\mathbb{B}} \boxtimes \mathscr{A}) \\
{\scriptstyle \mathfrak{R}^{\blacklozenge} \boxtimes \mathscr{A}} \uparrow & \Downarrow \mathfrak{m}_{\mathscr{B},\mathscr{A}} & \downarrow {\scriptstyle \mathscr{B} \boxtimes \mathfrak{L}_{\mathscr{A}}} \\
\mathscr{B} \boxtimes \mathscr{A} & = & \mathscr{B} \boxtimes \mathscr{A}
\end{array}
\qquad
\begin{array}{ccc}
(\mathscr{B} \boxtimes \mathscr{A}) \boxtimes \mathbb{1}_{\mathbb{B}} & \xrightarrow{\mathfrak{A}_{\mathscr{B},\mathscr{A},\mathbb{1}_{\mathbb{B}}}} & \mathscr{B} \boxtimes (\mathscr{A} \boxtimes \mathbb{1}_{\mathbb{B}}) \\
{\scriptstyle \mathfrak{R}^{\blacklozenge}_{\mathscr{B} \boxtimes \mathscr{A}}} \uparrow \; \mathfrak{r}_{\mathscr{B},\mathscr{A}} & \Leftarrow & \nearrow \\
\mathscr{B} \boxtimes \mathscr{A} & \xrightarrow{\mathscr{B} \boxtimes \mathfrak{R}^{\blacklozenge}_{\mathscr{A}}} &
\end{array}
$$

$$
\begin{array}{ccc}
(\mathbb{1}_{\mathbb{B}} \boxtimes \mathscr{B}) \boxtimes \mathscr{A} & \xrightarrow{\mathfrak{A}_{\mathbb{1}_{\mathbb{B}},\mathscr{B},\mathscr{A}}} & \mathbb{1}_{\mathbb{B}} \boxtimes (\mathscr{B} \boxtimes \mathscr{A}) \\
{\scriptstyle \mathfrak{L}_{\mathscr{B}} \boxtimes \mathscr{A}} \downarrow \; \mathfrak{l}_{\mathscr{B},\mathscr{A}} & \Rightarrow & \swarrow {\scriptstyle \mathfrak{L}_{\mathscr{B} \boxtimes \mathscr{A}}} \\
\mathscr{B} \boxtimes \mathscr{A} & &
\end{array}
$$

This data is subject to coherence axioms, see [JY21, Explanation 12.1.3].

**Example 2.111.** The 2-category $\mathbb{C}\mathrm{at}_{\mathbb{k}}$ of $\mathbb{k}$-linear categories, $\mathbb{k}$-linear functors, and natural transformations between them is a monoidal bicategory when endowed with the *tensor product* of $\mathbb{k}$-linear categories: given $\mathscr{A}$ and $\mathscr{B}$, define the $\mathbb{k}$-linear category $\mathscr{A} \otimes_{\mathbb{k}} \mathscr{B}$ by

- $\operatorname{Ob}(\mathscr{A} \otimes_{\mathbb{k}} \mathscr{B}) := \operatorname{Ob} \mathscr{A} \times \operatorname{Ob} \mathscr{B}$,
- $(\mathscr{A} \otimes_{\mathbb{k}} \mathscr{B})((a,b),(a',b')) := \mathscr{A}(a,a') \otimes_{\mathbb{k}} \mathscr{B}(b,b')$.

Given $\mathbb{k}$-linear functors $F\colon \mathscr{A} \to \mathscr{B}$ and $G\colon \mathscr{C} \to \mathscr{D}$, the tensor product $F \otimes_{\mathbb{k}} G$ is defined by $\operatorname{Ob}(F \otimes_{\mathbb{k}} G) = \operatorname{Ob}(F) \times \operatorname{Ob}(G)$ and

$$
(F \otimes_{\mathbb{k}} G)_{(a,b),(c,d)}\colon \mathscr{A}(a,c) \otimes_{\mathbb{k}} \mathscr{B}(b,d) \xrightarrow{F_{a,c} \otimes_{\mathbb{k}} G_{b,d}} \mathscr{C}(Fa,Fc) \otimes_{\mathbb{k}} \mathscr{D}(Gb,Gd).
$$

The tensor product of natural transformations extends similarly.

The associators and unitors are lifted from $\mathsf{Set}$ and $\mathsf{Vect}$, and the pentagonators and 2-unitors are trivial. For example, given $\Bbbk$-linear categories $\mathscr{A}$, $\mathscr{B}$, and $\mathscr{C}$, the associator $\mathfrak{A}_{\mathscr{A},\mathscr{B},\mathscr{C}}$ is the equivalence

$$\begin{aligned}(\mathscr{A}\otimes_{\Bbbk}\mathscr{B})\otimes_{\Bbbk}\mathscr{C} &\xrightarrow{\sim} \mathscr{A}\otimes_{\Bbbk}(\mathscr{B}\otimes_{\Bbbk}\mathscr{C})\\ ((a,b),c) &\longmapsto (a,(b,c))\\ (f\otimes_{\Bbbk}g)\otimes_{\Bbbk}h &\longmapsto f\otimes_{\Bbbk}(g\otimes_{\Bbbk}h).\end{aligned}$$

Diagram (2.7.1) commutes in this case, and so the pentagonator may be chosen to be the identity.

Generalising the above example, for a symmetric closed monoidal category $(\mathscr{V},\otimes)$, the bicategory $\mathscr{V}\text{-}\mathbb{C}\mathsf{at}$ is monoidal, with the monoidal structure $\mathscr{A}\otimes\mathscr{B}$ defined analogously to $\mathscr{A}\otimes_{\Bbbk}\mathscr{B}$ in the case of $\mathscr{V} \coloneqq \mathsf{Vect}$.

**Example 2.112.** The bicategory $\mathbb{B}\mathsf{imod}(\mathscr{C})$ of Example 2.103 is a monoidal bicategory with respect to the underlying tensor product of $\mathscr{C}$. Tensoring bimodules ${}_BM_A$ and ${}_DN_C$ is given by the $(B\otimes D)$-$(A\otimes C)$-bimodule ${}_BM_A\otimes{}_DN_C$, and this extends to bimodule homomorphisms in the obvious way.

For details on the following example, together with its more general variant where $\mathscr{V}$ is only required to be braided, see [DS97, Section 7].

**Example 2.113.** The bicategory $\mathbb{P}\mathsf{rof}_{\Bbbk}$ has as objects $\Bbbk$-linear categories, and the hom-categories are given by

$$\mathbb{P}\mathsf{rof}_{\Bbbk}(\mathscr{C},\mathscr{D}) \coloneqq \mathbb{C}\mathsf{at}_{\Bbbk}(\mathscr{D}^{\mathrm{op}}\otimes_{\Bbbk}\mathscr{C},\mathsf{Vect}).$$

We call $F \in \mathbb{P}\mathsf{rof}_{\Bbbk}(\mathscr{C},\mathscr{D})$ a *profunctor*, and write $F\colon \mathscr{C} \nrightarrow \mathscr{D}$. Composition of profunctors $F\colon \mathscr{C} \nrightarrow \mathscr{D}$ and $G\colon \mathscr{D} \nrightarrow \mathscr{E}$ is given by the coend

$$(G\diamond F)(e,c) \coloneqq \int^{d\in\mathscr{D}} G(e,d)\otimes_{\Bbbk} F(d,c).$$

Further, $\mathbb{P}\mathsf{rof}_{\Bbbk}$ is monoidal when endowed the monoidal structure $\otimes_{\Bbbk}$. On objects, it coincides with that in Example 2.111: given categories $\mathscr{A}$ and $\mathscr{B}$, their tensor product is $\mathscr{A}\otimes_{\Bbbk}\mathscr{B}$. Given profunctors $\Phi\colon \mathscr{A} \nrightarrow \mathscr{A}'$ and $\Psi\colon \mathscr{B} \nrightarrow \mathscr{B}'$, the profunctor $\Phi\otimes_{\Bbbk}\Psi\colon \mathscr{A}\otimes_{\Bbbk}\mathscr{A}' \nrightarrow \mathscr{B}\otimes_{\Bbbk}\mathscr{B}'$ is defined via

$$(\Phi\otimes_{\Bbbk}\Psi)((b,b'),(a,a')) = \Phi(a',a)\otimes_{\Bbbk}\Psi(b',b).$$

This extends similarly to natural transformations between profunctors.

### 2.7.1 *String diagrams in monoidal bicategories*

This section extends the graphical language of Section 2.3 by incorporating the tensor product of $\mathbb{C}\mathrm{at}$, which turns it into a monoidal 2-category. In our presentation, we closely follow Willerton [Wil08], see [DS03; Str03; BMS24; DS25] for other examples of "sheet" or "surface" diagrams.

As before, we consider strings and vertices between them. These are labelled with functors and natural transformations, respectively. The strings and vertices are embedded into bounded rectangles, which we will call sheets. Each (connected) region of a sheet is decorated with a category. The same mechanics as for ordinary 2-dimensional string diagrams apply: horizontal and vertical gluing represents composition of functors and natural transformations. On top of these operations, we add stacking sheets behind each other to depict the monoidal product of $\mathbb{C}\mathrm{at}$. Our convention is to read diagrams from front to back, right to left, and bottom to top.

Given a monoidal category $\mathscr{C}$, two of the most vital building blocks in this new graphical language are its tensor product $\otimes\colon \mathscr{C}\times\mathscr{C} \longrightarrow \mathscr{C}$ and unit $1 \in \mathscr{C}$; they are shown in Figure 2.2. On the left, there are two sheets equating to two copies of $\mathscr{C}$ joined by a line: the tensor product of $\mathscr{C}$. On the right, we have the unit of $\mathscr{C}$ considered as a functor from the terminal category to $\mathscr{C}$. We represent $\mathbb{1}$ by the empty sheet, and the unit of $\mathscr{C}$ by a dashed line.

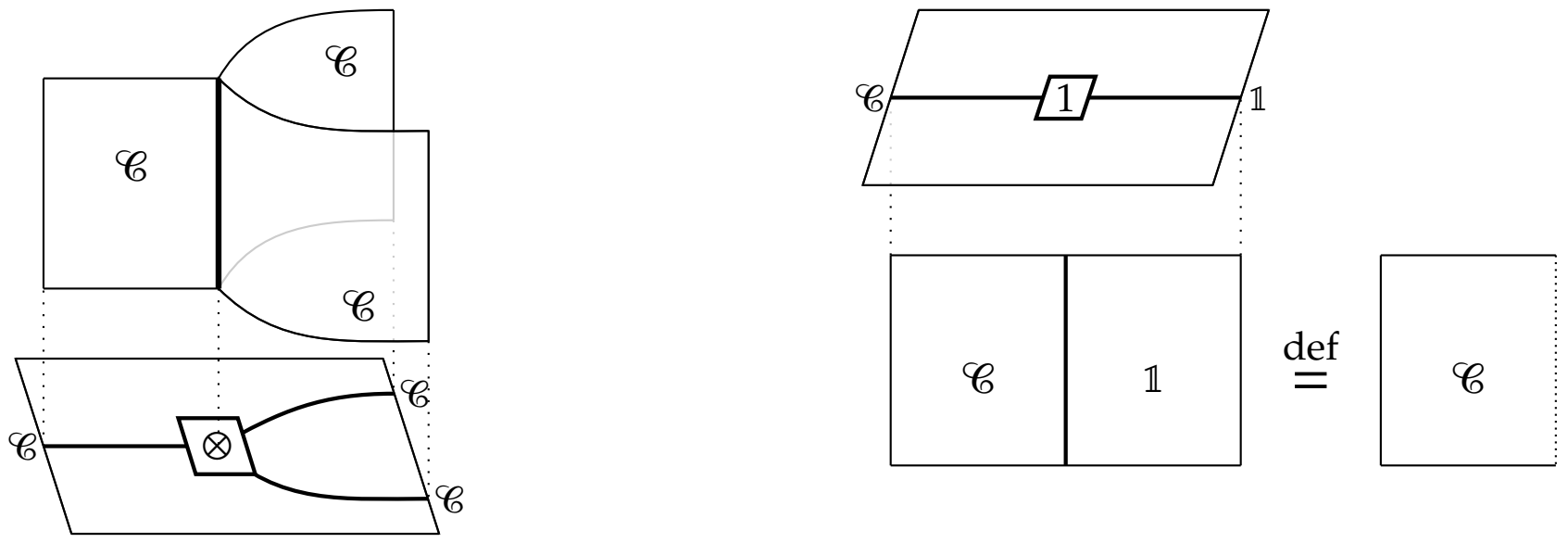

Figure 2.2: Graphical representation of the tensor product and unit of a monoidal category $(\mathscr{C}, \otimes, 1)$.

**Example 2.114.** Consider an oplax monoidal functor $(F, F_2, F_0)\colon \mathscr{C} \longrightarrow \mathscr{D}$. Figure 2.3 depicts the comultiplication $F_2$ as a "time evolution", where we start with $F(- \otimes =)$ on the bottom, and end up with $F(-) \otimes F(=)$ at the top. The coassociativity and counitality conditions are depicted in Figure 2.4.

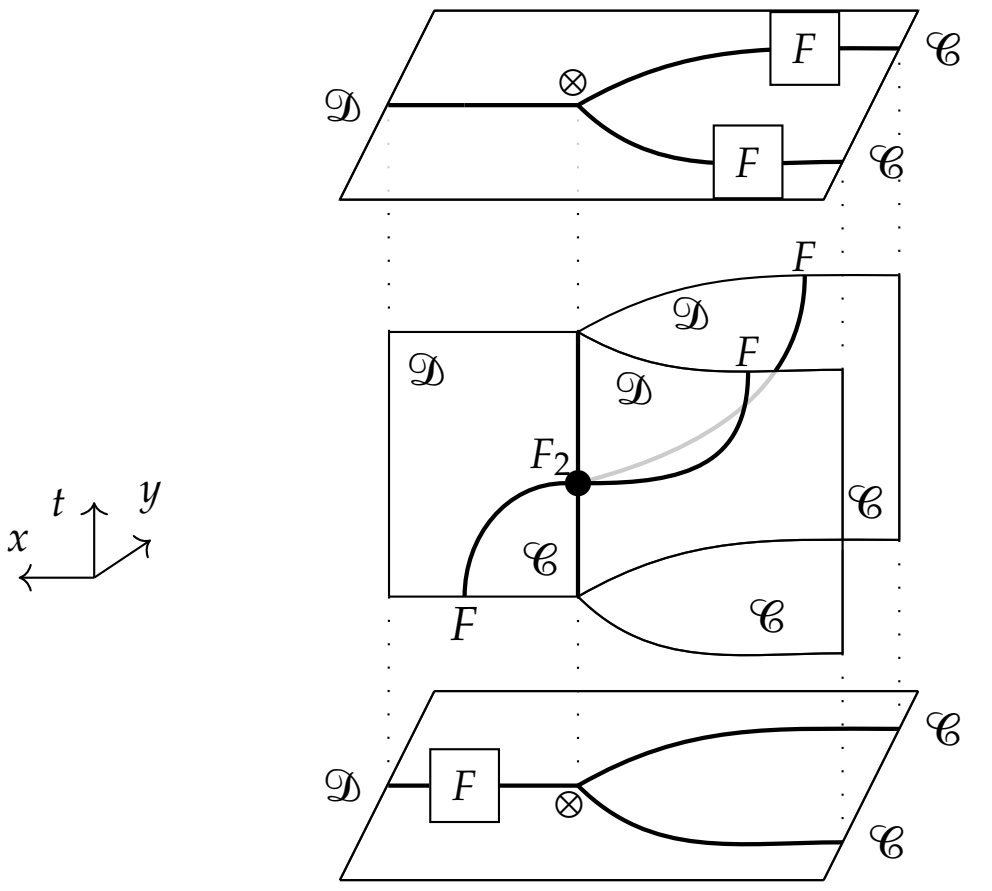

Figure 2.3: The comultiplication of an oplax monoidal functor $F\colon \mathscr{C} \longrightarrow \mathscr{D}$.

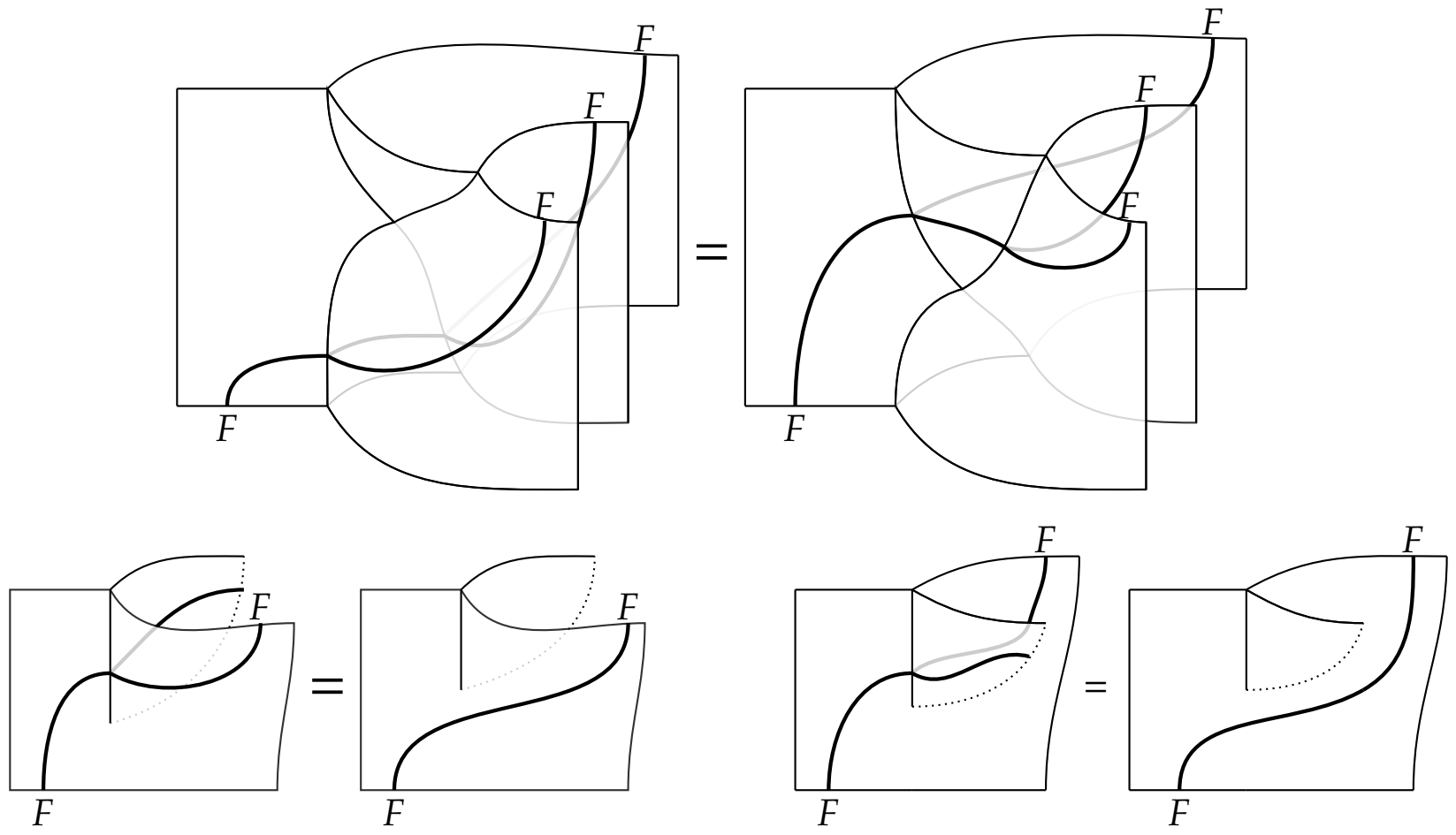

Figure 2.4: Coassociativity and counitality conditions of an oplax monoidal functor.

**Remark 2.115.** Given two oplax monoidal functors $F, G\colon \mathscr{C} \longrightarrow \mathscr{D}$, the string diagrammatic conditions for a natural transformation $\eta\colon F \Longrightarrow G$ to be one of oplax monoidal functors look like

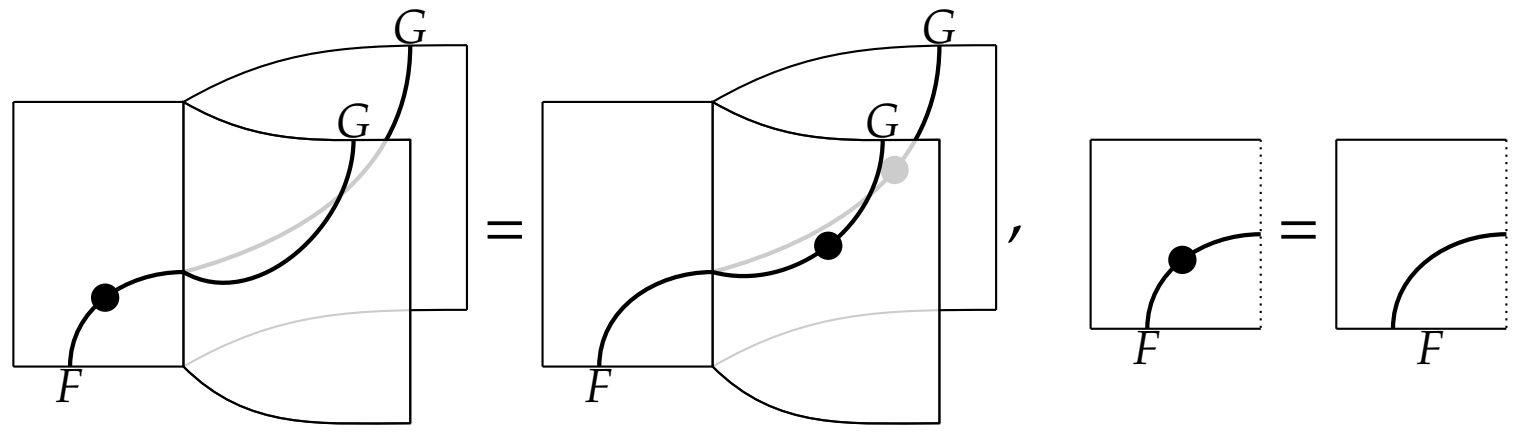

**Example 2.116.** Suppose that $F\colon \mathscr{C} \rightleftarrows \mathscr{D} :U$ is an oplax monoidal adjunction. In string diagrams, the conditions that the unit $\eta$ and the counit $\varepsilon$ are oplax monoidal natural transformation is given in Figure 2.5.

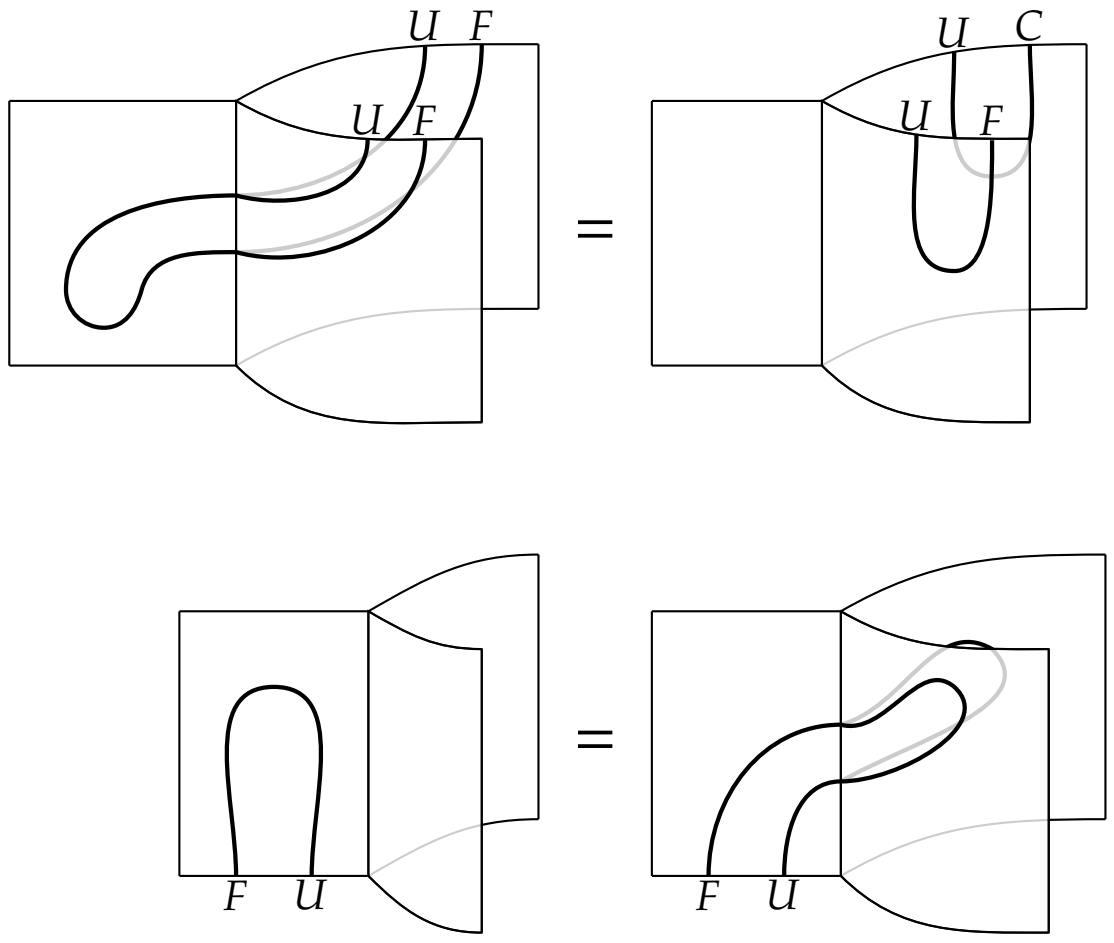

Figure 2.5: Conditions for $\eta$ and $\varepsilon$ to be oplax transformations.

## 2.8 COENDS

IN THIS SECTION WE GIVE THE DEFINITION OF an enriched coend in the special case of $\Bbbk$-linear categories. More general accounts can be found for example in [Gra80, Section 2.3], [Kel05, Chapter 3], and [Lor21, Chapter 4].

**Definition 2.117.** Let $\mathscr{C}$ be an essentially small $\Bbbk$-linear category. The *coend* of a $\Bbbk$-linear functor $P\colon \mathscr{C}^{\mathrm{op}} \otimes_{\Bbbk} \mathscr{C} \longrightarrow \mathsf{Vect}$ is following the coequaliser in $\mathsf{Vect}$:

$$\coprod_{a,b\in\mathscr{C}} \mathscr{C}(a,b) \otimes_{\Bbbk} P(b,a) \rightrightarrows \coprod_{a\in\mathscr{C}} P(a,a) \longrightarrow\!\!\!\!\rightarrow \int^{a\in\mathscr{C}} P(a,a).$$

For $a, b \in \mathscr{C}$ the two parallel morphisms are given by

$$\begin{aligned} \mathscr{C}(a,b) \otimes_{\Bbbk} P(b,a) &\longrightarrow P(a,a), & f \otimes x &\longmapsto P(f,a)x, \\ \mathscr{C}(a,b) \otimes_{\Bbbk} P(b,a) &\longrightarrow P(b,b), & g \otimes y &\longmapsto P(b,g)y. \end{aligned}$$

The *end* $\int_{a\in\mathscr{C}} P(a,a)$ of $P$ is defined analogously as an equaliser in $\mathsf{Vect}$.

**Remark 2.118.** A consequence of $\mathscr{C}$ being essentially small is that we can rewrite the (co)equalisers defining the respective (co)ends in Definition 2.117 to index over a set of objects, instead of a proper class. Since $\mathsf{Vect}$ is complete and cocomplete, the end and coend of any $P\colon \mathscr{C}^{\mathrm{op}} \otimes_{\Bbbk} \mathscr{C} \longrightarrow \mathsf{Vect}$ exists.

In case the indexing category $\mathscr{C}$ can be inferred unambiguously from the context, we will omit writing it explicitly.

**Remark 2.119.** Ends and coends have a number of useful properties, two of which we are going to frequently use in this thesis.

First, they are functorial: Any natural transformation $\alpha\colon P \Longrightarrow P'$ between functors $P, P'\colon \mathscr{C}^{\mathrm{op}} \times \mathscr{C} \longrightarrow \mathsf{Vect}$ induces a unique morphism

$$\int \mu\colon \int P \longrightarrow \int P'$$

by the universal property of (co-)equalisers.

Second, for appropriate functors $P$, the *Fubini–Tonelli interchange law* holds: if any of the following coends exist:

$$\int^{a\in\mathscr{C}} \int^{b\in\mathscr{C}} P(a,a,b,b), \qquad \int^{(a,b)\in\mathscr{C}\times\mathscr{C}} P(a,a,b,b), \qquad \int^{b\in\mathscr{C}} \int^{a\in\mathscr{C}} P(a,a,b,b),$$

then they all do and are isomorphic; we often simply write $\int^{a,b} P(a,a,b,b)$.

The following lemma is also known as the *enriched Yoneda lemma*.

**Lemma 2.120.** *Let $\mathscr{C}$ be a category, and let $F\colon \mathscr{C} \longrightarrow \mathsf{Vect}$ and $G\colon \mathscr{C}^{\mathrm{op}} \longrightarrow \mathsf{Vect}$ be functors. Then there are natural isomorphisms*

$$\int^{c} \mathscr{C}(c,-) \otimes_{\Bbbk} Fc \;\cong\; F \;\cong\; \int_{c} \mathsf{Vect}(\mathscr{C}(-,c), Fc), \tag{2.8.1}$$

*and*

$$\int^{c} \mathscr{C}(-,c) \otimes_{\Bbbk} Gc \;\cong\; G \;\cong\; \int_{c} \mathsf{Vect}(\mathscr{C}(c,-), Gc). \tag{2.8.2}$$

As the name suggests, the above lemma is an analogue of the Yoneda lemma in the enriched case; as such, we will only refer to it as "the Yoneda lemma" in subsequent considerations. For a proof, see [Kel05, Section 2.4].

**Definition 2.121.** Let $\mathscr{C}$ be a $\Bbbk$-linear monoidal category. Given two $\Bbbk$-linear functors $F, G\colon \mathscr{C} \longrightarrow \mathsf{Vect}$, the *Day convolution* of $F$ and $G$ is the $\Bbbk$-linear functor $F * G\colon \mathscr{C} \longrightarrow \mathsf{Vect}$ that is given by the coend

$$F * G := \int^{a,b} \mathscr{C}(a \otimes b, -) \otimes_{\Bbbk} Fa \otimes_{\Bbbk} Gb.$$

**Remark 2.122.** If $\mathscr{C}$ is left or right closed monoidal, we may simplify the Day convolution product using the Yoneda lemma for coends:

$$\begin{aligned} F * G &= \int^{a,b} \mathscr{C}(a \otimes b, -) \otimes_{\Bbbk} Fa \otimes_{\Bbbk} Gb \cong \int^{a,b} \mathscr{C}(a, [b, -]_\ell) \otimes_{\Bbbk} Fa \otimes_{\Bbbk} Gb \\ &\overset{(2.8.1)}{\cong} \int^{b} F([b, -]_\ell) \otimes_{\Bbbk} Gb, \end{aligned} \tag{2.8.3}$$

$$\begin{aligned} F * G &= \int^{a,b} \mathscr{C}(a \otimes b, -) \otimes_{\Bbbk} Fa \otimes_{\Bbbk} Gb \cong \int^{a,b} \mathscr{C}(b, [a, -]_r) \otimes_{\Bbbk} Fa \otimes_{\Bbbk} Gb \\ &\overset{(2.8.1)}{\cong} \int^{a} Fa \otimes G([a, -]_r). \end{aligned} \tag{2.8.4}$$

**Example 2.123.** Consider a $\Bbbk$-linear monoidal category $\mathcal{X}$ with a single object $x$. Then $A := \mathrm{End}_{\mathcal{X}}(x)$ is a commutative algebra and $[\mathcal{X}, \mathsf{Vect}] \cong \text{Mod-}A$. Let $F, G\colon \mathcal{X} \longrightarrow \mathsf{Vect}$ be two functors. Writing $M := Fx$ and $N := Gx$ for the corresponding modules over $A$, and $an := na$ for all $a \in A$ and $n \in N$, we obtain the following using Equation (2.8.3) and the definition of coends:

$$(F * G)x \cong {M \otimes_{\Bbbk} N}\big/{\langle ma \otimes_{\Bbbk} n - m \otimes_{\Bbbk} an \mid m \in M, n \in N, a \in A \rangle} = M \otimes_A N.$$

Thus, one recovers the tensor product of modules over commutative algebras.

**Theorem 2.124** ([Day71])**.** *For any monoidal category $\mathscr{C}$, the category $[\mathscr{C}^{\mathrm{op}}, \mathsf{Vect}]$ is closed monoidal with Day convolution as its tensor product, $\mathscr{C}(1, -)$ as its unit, and the internal homs given for all $F, G\colon \mathscr{C}^{\mathrm{op}} \longrightarrow \mathsf{Vect}$ by*

$$[F, G]_\ell := \int_{a,b} \mathsf{Vect}(\mathscr{C}(- \otimes a, b), \mathsf{Vect}(Fa, Gb)), \tag{2.8.5}$$

$$[F, G]_r := \int_{a,b} \mathsf{Vect}(\mathscr{C}(a \otimes -, b), \mathsf{Vect}(Fa, Gb)). \tag{2.8.6}$$

**Remark 2.125.** If the monoidal category $\mathscr{C}$ is closed, then the formulas for the internal homs of $[\mathscr{C}, \mathsf{Vect}]$ may be simplified by means of the Yoneda lemma:

$$\begin{aligned} [F, G]_r &= \int_{a,b} \mathsf{Vect}(\mathscr{C}(a \otimes -, b), \mathsf{Vect}(Fa, Gb)) \\ &\cong \int_{a,b} \mathsf{Vect}(\mathscr{C}(a \otimes -, b) \otimes_{\Bbbk} Fa, Gb) \cong \int_b \mathsf{Vect}\Big(\int^a \mathscr{C}(a \otimes -, b) \otimes_{\Bbbk} Fa, Gb\Big) \\ &\cong \int_b \mathsf{Vect}\Big(\int^a \mathscr{C}(a, [-, b]_\ell) \otimes_{\Bbbk} Fa, Gb\Big) \overset{(2.8.1)}{\cong} \int_b \mathsf{Vect}(F[-, b]_\ell, Gb), \end{aligned} \tag{2.8.7}$$

$$[F,G]_\ell \cong \int_b \mathsf{Vect}\Big(\int_a \mathscr{C}(a,[-,b]_r)\otimes_{\Bbbk} Fa, Gb\Big) \overset{(2.8.1)}{\cong} \int_a \mathsf{Vect}(F[-,b]_r, Gb). \tag{2.8.8}$$

**Remark 2.126.** Whenever we treat $\widehat{\mathscr{C}^{\mathrm{op}}} := [\mathscr{C}, \mathsf{Vect}]$ as a (closed) monoidal category, we implicitly equip with the convolution tensor product. Analogously, one may define a closed monoidal structure on $\widehat{\mathscr{C}} := [\mathscr{C}^{\mathrm{op}}, \mathsf{Vect}]$. Note, however, that in this case we cannot simplify the internal hom in the same way as in Remark 2.125; this would require the category $\mathscr{C}^{\mathrm{op}}$, as opposed to $\mathscr{C}$ itself, to be closed monoidal.

The convolution structure is particularly well-behaved on representables:

$$\mathscr{C}(-,x) * \mathscr{C}(-,y) = \int^{a,b\in\mathscr{C}} \mathscr{C}(-,a\otimes b)\otimes_{\Bbbk}\mathscr{C}(a,x)\otimes_{\Bbbk}\mathscr{C}(b,y) = \mathscr{C}(-,x\otimes y)$$

This connection extends to the entire functor category, see [Day71; IK86].

**Proposition 2.127.** *For a monoidal category $\mathscr{C}$, the Yoneda embedding*

$$よ\colon \mathscr{C} \longrightarrow \widehat{\mathscr{C}} = [\mathscr{C}^{\mathrm{op}}, \mathsf{Vect}], \qquad x \longmapsto \mathscr{C}(-,x) \tag{2.8.9}$$

*is a strong monoidal functor.*

## 2.9 (Co)completions

In this subsection we give a brief—informal—account of the results regarding the monoidal pseudofunctoriality of cocompletions and the resulting (co)completion operations for monoidal and module categories. We refer to [Kel05; KS06] for generalities on (co)limits and (co)completions. We implicitly assume all categories and functors to be $\Bbbk$-linear.

Let $\Phi$ be a class of diagrams. We say that a category is $\Phi$-*cocomplete* if it admits colimits of functors with domain in $\Phi$, and we say that a functor is $\Phi$-*cocontinuous* if it preserves such colimits.

**Definition 2.128.** A monoidal category $\mathscr{C}$ is called *separately $\Phi$-cocontinuous* if $\mathscr{C}$ is $\Phi$-cocomplete and its tensor product is separately $\Phi$-cocontinuous.

Similarly, for a $\Phi$-cocomplete monoidal category $\mathscr{D}$, a $\mathscr{D}$-module category $\mathscr{M}$ is said to be separately $\Phi$-cocontinuous if $\mathscr{M}$ is $\Phi$-cocomplete and the action $-\triangleright_{\mathscr{M}} =$ is separately $\Phi$-cocontinuous.

Let $\mathbb{C}\mathrm{at}_{\Phi\text{-Cocts}}$ be the 2-category of $\Phi$-cocomplete categories, $\Phi$-cocontinuous functors, and transformations between them. Let $\mathscr{C}$ be a small category. By [Kel82, Section 5.7], there is a category $\Phi(\mathscr{C})$, the $\Phi$-*cocompletion* of $\mathscr{C}$, and an embedding $\iota\colon \mathscr{C} \hookrightarrow \Phi(\mathscr{C})$, such that for any $\Phi$-cocomplete category $\mathscr{D}$ the restriction functor $\mathbb{C}\mathrm{at}_{\Phi\text{-Cocts}}(\mathscr{C}, \mathscr{D}) \xrightarrow{-\circ\iota} \mathbb{C}\mathrm{at}(\mathscr{C}, \mathscr{D})$ is an equivalence.

By the main results of [Zöb76; Koc95], the inclusion 2-functor of $\mathbb{C}\mathrm{at}_{\Phi\text{-Cocts}}$ into $\mathbb{C}\mathrm{at}$ is 2-monadic, with a left adjoint given by the pseudofunctor of $\Phi$-cocompletion, $\Phi(-)\colon \mathbb{C}\mathrm{at} \longrightarrow \mathbb{C}\mathrm{at}_{\Phi\text{-Cocts}}$. This 2-monad is a lax-idempotent 2-monad in the sense of [KL97].

By [KL00], there is a similar 2-monad $\Phi_c$ on $\mathbb{C}\mathrm{at}$, whose algebras are categories with a *prescribed choice* of $\Phi$-colimits. Strict morphisms are functors strictly preserving the chosen $\Phi$-colimits, and pseudomorphisms are functors preserving $\Phi$-colimits in the ordinary sense.

By [LF11, Theorem 6.2], this yields a pseudoclosed structure the 2-category $\mathbb{C}\mathrm{at}_{\Phi_c\text{-Cocts}}$ of categories with a choice of $\Phi$-colimits, $\Phi$-cocontinuous functors, and transformations between them. The category underlying the internal hom from $\mathscr{A}$ to $\mathscr{B}$ is given by $\mathbb{C}\mathrm{at}_{\Phi\text{-Cocts}}(\mathscr{A}, \mathscr{B})$, and the functors comprising the 2-monad are closed. Further, in some cases, for instance if $\Phi$ is the class of finite colimits, the pseudoclosed structure becomes pseudoclosed monoidal.

Since monoidal categories and module categories can be formulated as pseudomonoids internally to a pseudoclosed monoidal structure, see [HP02], we find the following result.

**Proposition 2.129.** *Let $\mathscr{C}$ be a monoidal category. Then $\Phi(\mathscr{C})$ is separately $\Phi$-cocontinuous monoidal, the inclusion $\iota\colon \mathscr{C} \hookrightarrow \Phi(\mathscr{C})$ is strong monoidal, and induces the following equivalence for any separately $\Phi$-cocontinuous monoidal category $\mathscr{D}$:*

$$\mathsf{StrMon}_{\Phi\text{-Cocts}}(\Phi(\mathscr{C}), \mathscr{D}) \xrightarrow{-\circ\iota} \mathsf{StrMon}(\mathscr{C}, \mathscr{D}).$$

*Similar equivalences are induced for lax and oplax monoidal functors.*

*Likewise, given a $\mathscr{C}$-module category $\mathscr{M}$, the $\Phi$-cocompletion $\Phi(\mathscr{M})$ is a separately $\Phi$-cocontinuous $\Phi(\mathscr{C})$hypmodule category, the inclusion $\iota\colon \mathscr{M} \hookrightarrow \Phi(\mathscr{M})$ is a strong $\mathscr{C}$-module functor, and the functor*

$$\mathsf{Str}\Phi(\mathscr{C})\mathsf{Mod}_{\Phi\text{-Cocts}}(\Phi(\mathscr{M}), \mathscr{N}) \xrightarrow{-\circ\iota} \mathsf{Str}\mathscr{C}\mathsf{Mod}(\mathscr{M}, \mathscr{N}), \tag{2.9.1}$$

*is an equivalence; similar equivalences exist for lax and oplax module functors.*

We remark that explicit constructions of the monoidal and module structures obtained in Proposition 2.129 and direct proofs of their universal properties are given in a variety of cases and ways in the literature, see for instance [MMMT19, Section 3] and [CG22, Example 3.2.9].

**Proposition 2.130.** *Let $\mathscr{A}$ be an additive category. The copower of $\mathscr{A}$ over* $\mathsf{vect}$*—see Example 5.22—is the essentially unique* $\mathsf{vect}$*-module structure on $\mathscr{A}$.*

*Proof.* Let $\mathsf{add}$ denote the collection of finite discrete categories. An $\mathsf{add}$-cocomplete category is simply an additive category, and an $\mathsf{add}$-cocontinuous functor is an additive functor.

Since any $\Bbbk$-linear functor preserves direct sums, an additive $\mathsf{vect}$-module category is necessarily separately additive in the sense of Definition 2.128. The result then follows by observing that $\mathsf{vect} \simeq \mathsf{add}(\{\circledast\})$, where $\{\circledast\}$ is the terminal monoidal category, so

$$\mathsf{StrMon}_{\mathsf{add}}(\mathsf{vect}, \mathbb{C}\mathsf{at}_{\mathsf{add}}(\mathscr{A}, \mathscr{A})) \simeq \mathsf{StrMon}(\{*\}, \mathbb{C}\mathsf{at}_{\mathsf{add}}(\mathscr{A}, \mathscr{A})) \simeq \{*\},$$

where $\{*\}$ is the terminal category and $\mathsf{StrMon}_{\mathsf{add}}$ denotes the category of strong monoidal additive functors. □

**Proposition 2.131.** *Let $\mathscr{A}$ be an additive category admitting filtered colimits. The action of* $\mathsf{vect}$ *on $\mathscr{A}$ extends uniquely to a separately finitary* $\mathsf{Vect}$*-module category structure on $\mathscr{A}$.*

*Proof.* Let $\{\mathsf{add}, \mathsf{filt}\}$ denote the collection of diagrams which are filtered or finite discrete. Then we have $\mathsf{Vect} \simeq \{\mathsf{add}, \mathsf{filt}\}(\{\circledast\})$, and

$$\begin{aligned} &\mathbf{StrMon}_{\mathsf{filt}}(\mathsf{Vect}, \mathbb{C}\mathsf{at}_{\mathsf{add},\mathsf{filt}}(\mathscr{A}, \mathscr{A})) = \mathbf{StrMon}_{\mathsf{add},\mathsf{filt}}(\mathsf{Vect}, \mathbb{C}\mathsf{at}_{\mathsf{add},\mathsf{filt}}(\mathscr{A}, \mathscr{A})) \\ &\simeq \mathbf{StrMon}(\{\circledast\}, \mathbb{C}\mathsf{at}_{\mathsf{add},\mathsf{filt}}(\mathscr{A}, \mathscr{A})) = \mathbf{StrMon}(\{\circledast\}, \mathbb{C}\mathsf{at}_{\mathsf{filt}}(\mathscr{A}, \mathscr{A})) \simeq \{*\}. \end{aligned}$$ □

In view of Lemma 2.89, filtered colimits and cocompletions under them, which we shall call *ind-completions*, are particularly important in the study of locally finite abelian categories. For example, the ind-completion $\mathsf{Ind}(\mathscr{C})$ of a locally finite abelian category $\mathscr{C}$ always has enough injectives. Another reason for our study of these completions is that categories of comodules are locally finitely presentable, whence we can make use of the characterisation of adjoint functors, see Propositions 2.133 and 2.134 below.

**Definition 2.132.** An object $x$ in a category $\mathscr{C}$ is called *compact* or *finitely presented* if $\mathscr{C}(x, -)$ is finitary; i.e., preserves filtered colimits. We denote the full subcategory of compact objects by $\mathscr{C}_c$. A category $\mathscr{C}$ is called *locally finitely presentable* if every object in $\mathscr{C}$ is a filtered colimit of compact objects.

Examples of ind-completions arise from the fact that $\mathsf{Ind}({}^C\mathsf{vect}) \simeq {}^C\mathsf{Vect}$.

**Proposition 2.133** (Special adjoint functor theorem for right adjoints)**.** *If $\mathscr{C}$ and $\mathscr{D}$ are locally finite abelian categories, then there are equivalences*

$$\mathsf{Rex}(\mathscr{C}, \mathscr{D}) \xrightarrow{\mathsf{Ind}} \mathsf{Cocont}(\mathsf{Ind}(\mathscr{C}), \mathsf{Ind}(\mathscr{D})) \xrightarrow{\sim} \mathsf{Map}(\mathsf{Ind}(\mathscr{C}), \mathsf{Ind}(\mathscr{D})),$$

*where $\mathsf{Map}(\mathscr{A}, \mathscr{B})$ denotes the category of functors admitting a right adjoint.*

*In particular, a functor from $\mathsf{Ind}(\mathscr{C})$ to $\mathsf{Ind}(\mathscr{D})$ admits a right adjoint if and only if it is cocontinuous, and the extension of a functor $F\colon \mathscr{C} \longrightarrow \mathscr{D}$ to $\mathsf{Ind}$-cocompletions is cocontinuous if and only if $F$ is right exact.*

**Proposition 2.134** (Special adjoint functor theorem for left adjoints)**.** *Let $\mathscr{C}$ and $\mathscr{D}$ be locally finite abelian categories. If a functor $G\colon \mathsf{Ind}(\mathscr{C}) \longrightarrow \mathsf{Ind}(\mathscr{D})$ is continuous and preserves filtered colimits, then it admits a left adjoint.*

*If a functor $G\colon \mathsf{Ind}(\mathscr{C}) \longrightarrow \mathsf{Ind}(\mathscr{D})$ admits a left adjoint $F$, then $G$ preserves filtered colimits if and only if $F$ preserves compact objects in $\mathscr{C}$.*

The latter part of Proposition 2.133 follows from $\mathsf{Ind}(F)$ preserving both finite and filtered colimits, and thus preserving all colimits. Since filtered colimits in a locally finitely presentable category are exact, the extension $\mathsf{Ind}(G)$ of a left exact functor $G$ is left exact. However, it is not clear that $\mathsf{Ind}(G)$ preserves arbitrary products.

**Definition 2.135.** An object $X \in \mathsf{Ind}(\mathscr{C})$ for a locally finite abelian category $\mathscr{C}$ is *finitely cogenerated* if the functor $\mathsf{Ind}(\mathscr{C})(-, X)$ preserves arbitrary products.

One can verify that for a coalgebra $D$ such that $\mathsf{Ind}(\mathscr{C}) \simeq$ Comod-$D$, a comodule $M$ is finitely cogenerated if and only if $M$ embeds into a comodule of the form $D^{\oplus m}$ for some $m \in \mathbb{N}$. Observe that an object $M \in \mathsf{Ind}(\mathscr{D})$ is compact if and only if it lies in the essential image of the embedding $\mathscr{D} \hookrightarrow \mathsf{Ind}(\mathscr{D})$.

**Definition 2.136.** Let $\mathscr{C}$ and $\mathscr{D}$ be locally finite abelian categories. A functor $F\colon \mathsf{Ind}(\mathscr{C}) \longrightarrow \mathsf{Ind}(\mathscr{D})$ is called *quasi-finite* if for any finitely cogenerated injective object $I$ in $\mathsf{Ind}(\mathscr{C})$ and compact object $M$ of $\mathscr{D}$, the $\mathbb{k}$-vector space $\mathsf{Ind}(\mathscr{D})(M, FI)$ is finite-dimensional.

Using the equivalence $\mathscr{C} \simeq {}^C\mathsf{vect}$ of Lemma 2.89, the following result is a direct consequence of [Tak77, Proposition 1.3].

**Proposition 2.137.** *Let $G\colon \mathscr{C} \longrightarrow \mathscr{D}$ be a left exact functor of locally finite abelian categories. The functor $\mathsf{Ind}(G)$ has a left adjoint if and only if it is quasi-finite.*

A stronger result holds for finite abelian categories, as a direct consequence of the Eilenberg–Watts theorem for categories of bimodules, see for example [FSS20, Lemma 2.1].

**Proposition 2.138.** *A functor $F\colon \mathscr{A} \longrightarrow \mathscr{B}$ of finite abelian categories is right exact if and only if it has a right adjoint, and left exact if and only if it has a left adjoint.*

Recall that an object $A$ in a category $\mathscr{A}$ that admits $\Phi$-colimits is said to be $\Phi$-*small* if the functor $\mathscr{A}(A, -)$ preserves them.

**Lemma 2.139.** *Let $G\colon \mathscr{C} \longrightarrow \mathscr{D}$ be a $\Phi$-cocontinuous functor, and assume $G$ has a left adjoint $F$. Then $F$ sends $\Phi$-small objects to $\Phi$-small objects.*

*Proof.* Let $x \in \mathscr{D}$ be a $\Phi$-small object. Then

$$\mathscr{C}(Fx, -) \cong \mathscr{D}(x, G(-)) = \mathscr{D}(x, -) \circ G.$$

The functor $\mathscr{C}(Fx, -)$ is thus naturally isomorphic to a composition of two $\Phi$-cocontinuous functors, and hence it itself is a $\Phi$-cocontinuous functor. □

May it be a light to you in dark places,
when all other lights go out.

J. R. R Tolkien; The Fellowship of the Ring

# DUALITY THEORY FOR MONOIDAL CATEGORIES

3

The aim of this chapter is to compare several notions of "duality" in monoidal categories. On the one hand, we have closedness, exemplifying the tensor–hom adjunction in the category of sets; on the other lies rigidity, generalising the duals of finite-dimensional vector spaces. Grothendieck–Verdier duality [BD13], also called $*$-autonomy [Bar79], describes a duality theory between the strict confinements of rigidity and the very general notion of monoidal closedness. Roughly speaking, Grothendieck–Verdier categories consist of a closed monoidal category and a chosen dualising object, such that the internal hom into this object induces an anti-equivalence of categories. While closed and rigid categories have left and right variants, Grothendieck–Verdier duality is an inherently ambidextrous concept.

Our first goal is to more closely investigate the converse of Proposition 2.67, see also Remark 2.68. That is, we study whether one can decide if a closed monoidal category is rigid solely by verifying that the internal hom is given by tensoring with an object. We thank Chris Heunen for suggesting this may not hold in general at the BCQT2022 summer school.

**Theorem 3.23.** *There exists a non-rigid tensor representable category.*

In particular, taking Theorem 3.23, Proposition 3.16 and Example 3.17, as well as Proposition 3.12 and Example 3.15 together yields a strict hierarchy:

$$\text{Rigid} \subsetneq \text{Tensor representable} \subsetneq \text{Grothendieck–Verdier} \subsetneq \text{Closed}.$$

For the second inclusion we assume that the natural transformations "comparing" the left and right-handed version of tensor representability are invertible, see Propositions 3.7 and 3.16. For rigid categories this is always the case.

We then investigate closed and Grothendieck–Verdier structures on functor categories endowed with Day convolution as its tensor product. Given a sufficiently nice base category, the functor category inherits these dualities.

**Proposition 3.27.** *Let $\mathscr{C}$ be a $\Bbbk$-linear hom-finite Grothendieck–Verdier category. If Hypothesis 3.25 holds—all relevant coends exist and are finite-dimensional vector spaces—then $[\mathscr{C}, \mathsf{vect}]$ is a Grothendieck–Verdier category.*

We furthermore investigate the duality properties of Cauchy completions. Since we are working with the functor category $[\mathscr{C}, \mathsf{vect}]$ instead of the presheaves $[\mathscr{C}^{\mathrm{op}}, \mathsf{vect}]$, we work with the Cauchy completion $\overline{\mathscr{C}^{\mathrm{op}}}$ of $\mathscr{C}^{\mathrm{op}}$, it being a subcategory of the free cocompletion $\widehat{\mathscr{C}^{\mathrm{op}}} = [\mathscr{C}, \mathsf{vect}]$ of $\mathscr{C}^{\mathrm{op}}$. Recall the definition of the Yoneda embedding from Equation (2.8.9).

**Corollary 3.43.** *Let $\mathscr{C}^{\mathrm{op}}$ be a $\Bbbk$-linear right closed monoidal category.*

*(i) $\mathscr{C}^{\mathrm{op}}$ is right rigid if and only if its Cauchy completion $\overline{\mathscr{C}^{\mathrm{op}}}$ is.*
*(ii) $\mathscr{C}^{\mathrm{op}}$ is right tensor representable if and only if $\overline{\mathscr{C}^{\mathrm{op}}}$ is.*
*(iii) $(\mathscr{C}^{\mathrm{op}}, d)$ is a right Grothendieck–Verdier category if and only if $(\overline{\mathscr{C}^{\mathrm{op}}}, よ_d)$ is.*

The chapter concludes with an investigation of three examples of abelian closed monoidal functor categories. The first is finite Boolean algebras, which historically arose in the study of order theory and universal algebra; the connection to Grothendieck–Verdier duality was already noted in [Bar79].

See Definition 3.53 for when an algebra is QF-2.

**Proposition 3.54.** *The path algebra of a finite Boolean algebra is* QF-2.

We subsequently study the category of finite-dimensional Mackey functors, objects abstracting and generalising the induction, restriction, and conjugation functors in classical representation theory [Lin76; TW95; Web00].

**Proposition 3.57.** *Let $G$ be a finite group and suppose $\Bbbk$ is a field. The category $\mathsf{mky}_{\Bbbk}(G)$ is a Grothendieck–Verdier category, and its rigid objects are precisely the finitely-generated projective Mackey functors. Furthermore, $\mathsf{mky}_{\Bbbk}(G)$ is rigid if and only if the Mackey algebra $\mathbb{M}G$ is semisimple, which is equivalent to $\operatorname{char}\Bbbk$ not dividing the order of $G$.*

Finally, we concern ourselves with strict 2-groups and their equivalent formulation as crossed modules [Wag21]. These are related to *r-categories*, a special case of Grothendieck–Verdier duality in which the dualising object is isomorphic to the monoidal unit.

**Proposition 3.65.** *Let $G$ and $H$ be finite groups, and suppose Suppose that*

$$\big(G,\, H,\, t\colon H \longrightarrow G,\, \alpha\colon G \longrightarrow \mathrm{Aut}(H)\big)$$

*is a crossed module. Let $\mathcal{G}$ be its associated strict 2-group. The category $\mathsf{rep}_{\mathcal{G}}(\ker t)$ of finite-dimensional representations of $\ker t$ is a right $r$-category.*

## 3.1 TENSOR REPRESENTABILITY

**Definition 3.1.** We call an endofunctor $F\colon \mathcal{C} \longrightarrow \mathcal{C}$ of a monoidal category *right tensor representable* if there exists an $x \in \mathcal{C}$ such that $F \cong x \otimes -$.

Proposition 2.67 states that the internal hom of a rigid monoidal category is tensor-representable at every object. From now on, for the sake of brevity, we will simply speak of "(monoidal) tensor representable categories" instead of "closed monoidal categories with tensor representable internal homs". As shown in Lemma 3.5 below, this is equivalent to the following definition, which does not rely on closedness.

**Definition 3.2.** A monoidal category $\mathcal{C}$ is called *right tensor representable* if for all $x \in \mathcal{C}$ there exists an object $Rx \in \mathcal{C}$, called the *right tensor dual* of $x$, such that $x \otimes - \dashv Rx \otimes -$.

A monoidal category $\mathcal{C}$ is called *left tensor representable* if $\mathcal{C}^{\mathrm{rev}}$ is right tensor representable. A left and right tensor representable category will be referred to as *tensor representable*.

Tensor representable categories encompass all rigid monoidal categories. As such, some results in the theory of the latter carry over to the former. For example, the following result is well-known in the case of rigid monoidal abelian categories, see for example [EGNO15, Corollary 4.2.13].

**Lemma 3.3.** *All objects of an abelian monoidal right tensor representable category are projective if and only if its unit is projective.*

*Proof.* As the first claim implies the second, we only have to show its converse. Fix an object $x$ in an abelian right tensor representable category $\mathcal{A}$. There exists a chain of adjunctions

$$x \otimes - \dashv Rx \otimes - \dashv R^2x \otimes - \dashv \dots .$$

Hence, the functor $Rx \otimes - \colon \mathcal{A} \longrightarrow \mathcal{A}$ has left and right adjoints and is therefore exact. Consequently, $\mathcal{A}(x, -) \cong \mathcal{A}(1, Rx \otimes -)$ is—as a composite of exact functors—exact itself. □

While Lemma 2.64 holds in the rigid case, it need not be true for tensor representable categories: given a strong monoidal functor $F$, the image $Fx$ of an object $x$ admitting a (right) tensor dual may not have a (right) tensor dual itself. However, the property is preserved by monoidal equivalences.

**Lemma 3.4.** *Let $\mathscr{C}$ be a right tensor representable category, and suppose that $\mathscr{D}$ is monoidal. Any strong monoidal adjoint equivalence $F\colon \mathscr{C} \rightleftarrows \mathscr{D} : F^{-1}$ equips $\mathscr{D}$ with a right tensor representable structure.*

*Proof.* For any $x \in \mathscr{D}$, define $Dx \coloneqq FRF^{-1}x \in \mathscr{D}$, where $RF^{-1}x$ is a right tensor dual of $F^{-1}x \in \mathscr{C}$. Then $Dx \otimes -$ is right adjoint to $x \otimes -$:

$$\begin{aligned}\mathscr{D}(x \otimes y, z) &\cong \mathscr{C}(F^{-1}(x \otimes y), F^{-1}z) \cong \mathscr{C}(F^{-1}x \otimes F^{-1}y, F^{-1}z)\\ &\cong \mathscr{C}(F^{-1}y, RF^{-1}x \otimes F^{-1}z) \cong \mathscr{D}(FF^{-1}y, F(RF^{-1}x \otimes F^{-1}z))\\ &\cong \mathscr{D}(y, FRF^{-1}x \otimes z) = \mathscr{D}(y, Dx \otimes z).\end{aligned}$$

Thus, $\mathscr{D}$ is right tensor representable. □

Just like in the rigid case, the relationship between closedness and tensor representability is governed by a family of natural isomorphisms.

**Lemma 3.5.** *A monoidal category $\mathscr{C}$ is right tensor representable if and only if it is right closed monoidal and for all $x \in \mathscr{C}$ there are isomorphisms*

(3.1.1) $$\varphi_y^{(x)}\colon [x, y]_r \longrightarrow [x, 1]_r \otimes y \qquad \textit{natural in } y \in \mathscr{C}.$$

*Proof.* For every $x \in \mathscr{C}$ in a right tensor representable category $\mathscr{C}$, by definition there exists an object $Rx \in \mathscr{C}$ such that $Rx \otimes -$ is right adjoint to $x \otimes -$. Hence, the category $\mathscr{C}$ is right closed monoidal, and the internal hom from $x \in \mathscr{C}$ is given by $Rx \otimes -$. The claim follows by noticing that Equation (3.1.1) is now even the functorial equality $Rx \otimes - = Rx \otimes 1 \otimes -$.

Conversely, for all $x \in \mathscr{C}$ there is an isomorphism

$$\mathscr{C}(y, [x, z]_r) \cong \mathscr{C}(x \otimes y, z) \cong \mathscr{C}(y, [x, 1]_r \otimes z)$$

natural in $y, z \in \mathscr{C}$; the result follows by the Yoneda lemma. □

**Definition 3.6.** Assume $\mathscr{C}$ to be a right tensor representable category. By Lemma 3.5, it is right closed monoidal. Then $R \coloneqq [-, 1]_r\colon \mathscr{C}^{\mathrm{op}} \longrightarrow \mathscr{C}$ is called the *right (tensor-)dualising functor* of $\mathscr{C}$. Analogously, $L \coloneqq [-, 1]_\ell\colon \mathscr{D}^{\mathrm{op}} \longrightarrow \mathscr{D}$ is the *left (tensor-)dualising functor* of a left tensor representable category $\mathscr{D}$.

A right rigid monoidal category is simultaneously left rigid if and only if its right dualising functor is an equivalence of categories. This is more complicated in the tensor representable case. Assume $\mathscr{C}$ to be of this type.

Write $R\colon \mathscr{C}^{\mathrm{op}} \longrightarrow \mathscr{C}$ and $L\colon \mathscr{C}^{\mathrm{op}} \longrightarrow \mathscr{C}$ for the right and left tensor dualising functors, respectively. Set

$$\eta_y^{(x)}\colon y \longrightarrow Rx \otimes x \otimes y, \qquad \varepsilon_y^{(x)}\colon x \otimes Rx \otimes y \longrightarrow y,$$
$$u_y^{(x)}\colon y \longrightarrow y \otimes x \otimes Lx, \qquad c_y^{(x)}\colon y \otimes Lx \otimes x \longrightarrow y$$

for the units and counits of the correspondig adjunctions. The left and right tensor dualising functors are related to each other by

$$x \xrightarrow{u_x^{(Rx)}} x \otimes Rx \otimes LRx \xrightarrow{\varepsilon_1^{(x)} \otimes LRx} LRx, \tag{3.1.2}$$
$$x \xrightarrow{\eta_x^{(Lx)}} RLx \otimes Lx \otimes x \xrightarrow{RLx \otimes c_1^{(x)}} RLx. \tag{3.1.3}$$

**Proposition 3.7.** *The right dualising functor $R\colon \mathscr{C}^{\mathrm{op}} \longrightarrow \mathscr{C}$ of a tensor representable category $\mathscr{C}$ is right adjoint to $L^{\mathrm{op}}\colon \mathscr{C} \longrightarrow \mathscr{C}^{\mathrm{op}}$. Further, $R$ is an equivalence if and only if the canonical morphisms of Equations (3.1.2) and (3.1.3) are invertible. In this case $L^{\mathrm{op}}$ is a quasi-inverse of $R$.*

*Proof.* For any two objects $x, y \in \mathscr{C}$, we have

$$\mathscr{C}^{\mathrm{op}}(Lx, y) = \mathscr{C}(y, Lx) \cong \mathscr{C}(y \otimes x, 1) \cong \mathscr{C}(x, Ry).$$

A direct computation shows that the unit and counit of this adjunction are given by the natural transformations displayed in Equations (3.1.2) and (3.1.3). If these morphisms are invertible, $R$ and $L^{\mathrm{op}}$ are quasi-inverse to each other.

Conversely, assume $R$ to be an equivalence of categories. Since any equivalence can be bettered to an adjoint equivalence and adjoints are unique up to unique natural isomorphism, $L^{\mathrm{op}}$ must be quasi-inverse to $R$ and the unit and counit morphisms of Equations (3.1.2) and (3.1.3) are invertible. □

### 3.1.1 *Grothendieck–Verdier duality*

We are going to loosely follow the notation and terminology of [BD13].

**Definition 3.8.** A *(right) Grothendieck–Verdier category* consists of a pair $(\mathscr{C}, d)$ of a monoidal category $\mathscr{C}$ and an object $d \in \mathscr{C}$, such that for all $x \in \mathscr{C}$ there exists a representing object $D_r x \in \mathscr{C}$ for $\mathscr{C}(x \otimes -, d)$; i.e.,

$$\mathscr{C}(x \otimes -, d) \cong \mathscr{C}(-, D_r x), \tag{3.1.4}$$

and the induced *right dualising functor* $D_r\colon \mathscr{C}^{\mathrm{op}} \longrightarrow \mathscr{C}$ is an equivalence of categories. If $d = 1$ is the monoidal unit, one speaks of an *r-category*.

Because $D_r$ is an equivalence, one also obtains

$$\mathscr{C}(x \otimes y, d) \cong \mathscr{C}(y, D_r x) \cong \mathscr{C}^{\mathrm{op}}(D_\ell^{-1} y, x) = \mathscr{C}(x, D_r^{-1} y), \quad \text{for all } x, y \in \mathscr{C}. \tag{3.1.5}$$

**Remark 3.9.** In [BD13], the authors instead consider a pair $(\mathscr{C}, d)$ subject to the natural isomorphism $\mathscr{C}(- \otimes y, d) \cong \mathscr{C}(-, D_\ell y)$, such that the induced functor $D_\ell \colon \mathscr{C}^{\mathrm{op}} \longrightarrow \mathscr{C}$ is an equivalence. In our framework, we shall call this a *left* Grothendieck–Verdier category structure on $\mathscr{C}$, with *left* dualising functor $D_\ell$. Equations (3.1.4) and (3.1.5) become

$$\mathscr{C}(x \otimes y, d) \cong \mathscr{C}(x, D_\ell y) \cong \mathscr{C}(y, D_\ell^{-1} x). \tag{3.1.6}$$

**Remark 3.10.** Notice that the nomenclature of Remark 3.9 might be misleading at first sight. Because the dualising functor of a Grothendieck–Verdier category is assumed to be an equivalence, Grothendieck–Verdier structures are always two-sided. In view of Equation (3.1.5), a right Grothendieck–Verdier category $(\mathscr{C}, d)$ is also a left Grothendieck–Verdier category with the same dualising object $d$, as well as $D_r^{-1} = D_\ell$. Likewise, all left Grothendieck–Verdier categories are also right-sided.

For the sake of clarity and to denote the directional bias, Grothendieck–Verdier structures will often nevertheless explicitly carry a prefix. In particular, in Section 3.2 we consider the induced Grothendieck–Verdier structure of the category of copresheaves $\widehat{\mathscr{C}^{\mathrm{op}}}$ over a Grothendieck–Verdier base category $\mathscr{C}$. In this case, the *left* dualising functor of $\mathscr{C}$ naturally gives rise to a *right* dualising functor for $\widehat{\mathscr{C}^{\mathrm{op}}}$; see Proposition 3.27.

The symmetric nature of Grothendieck–Verdier categories is in stark contrast to the one-sidedness of rigidity and tensor representability.

**Example 3.11.** Consider a Hopf algebra $H$ whose antipode is not invertible, see [Sch00b]. Its finite-dimensional right comodules are a left-rigid category comod-$H$ which, as discussed in [Ulb90, Remark 2], is not right rigid. If comod-$H$ was a Grothendieck–Verdier category, there would as a consequence of Remark 3.13 exist a comodule $M \in$ comod-$H$, such that ${}^{\vee}(-) \otimes M \colon (\text{comod-}H)^{\mathrm{op}} \longrightarrow \text{comod-}H$ is an equivalence of categories. This would imply that $M$ is one-dimensional, and therefore the left dualising functor ${}^{\vee}(-) \colon \text{comod-}H^{\mathrm{op}} \longrightarrow \text{comod-}H$ would be an equivalence—contradiction.

The following proposition, also discussed in [BD13, Section 2], is a direct consequence of Equations (3.1.4) and (3.1.5).

**Proposition 3.12.** *Every (right) Grothendieck–Verdier category* $(\mathscr{C}, d)$ *is closed monoidal; for all* $x, y \in \mathscr{C}$*, the left and right internal homs are given by*

$$\begin{aligned} [x, y]_\ell &:= D_r^{-1}(x \otimes D_r y) = D_\ell(x \otimes D_\ell^{-1} y), \\ [x, y]_r &:= D_r(D_r^{-1} y \otimes x) = D_\ell^{-1}(D_\ell y \otimes x). \end{aligned} \tag{3.1.7}$$

**Remark 3.13.** If $(\mathscr{C}, d)$ is a Grothendieck–Verdier category, one recovers its dualising object by applying the dualising functor to the monoidal unit of $\mathscr{C}$. In fact, as shown in [BD13, Remark 2.1.(4)], this relationship is bidirectional:

$$D_r 1 \cong d, \quad D_r d \cong 1, \quad D_\ell 1 \cong d, \quad D_\ell d \cong 1, \quad D_r^2 1 \cong 1, \quad D_\ell^2 d \cong d. \tag{3.1.8}$$

In combination with Proposition 3.12, one obtains an explicit description of the dualising functors of $\mathscr{C}$:

$$D_r x \cong [x, d]_r \qquad \text{and} \qquad D_\ell x = D_r^{-1} x \cong [x, d]_\ell, \qquad \text{for all } x \in \mathscr{C}. \tag{3.1.9}$$

While Grothendieck–Verdier duality is an extra bit of structure on a monoidal category—a priori making it quite distinct from the other dualities—in practise it still behaves very similarly to how a property would, as the invertible objects govern the possible choices for dualising objects.

**Lemma 3.14** ([BD13, Proposition 2.3]). *Given a right Grothendieck–Verdier category* $(\mathscr{C}, d)$*, there is an anti-equivalence between its full subcategories of invertible and dualising objects given by*

$$\mathrm{Inv}(\mathscr{C}) \longrightarrow \mathrm{Dual}(\mathscr{C}), \qquad \alpha \longmapsto \alpha^{-1} \otimes d.$$

**Example 3.15.** The category $(\mathsf{Set}, \times, \mathbb{1})$ of sets with its usual Cartesian monoidal structure is closed monoidal. If some set $X$ is invertible, then it must be isomorphic to the terminal set $\mathbb{1}$. Thus, to deduce that there is no Grothendieck–Verdier structure on $\mathsf{Set}$, it is sufficient to see that the functor $\mathsf{Set}(-, \mathbb{1})$ is not an anti-equivalence, which follows by $\mathbb{1}$ being terminal.

Let us now characterise the relationship between Grothendieck–Verdier and tensor representable categories.

**Proposition 3.16.** *Let* $\mathscr{C}$ *be a tensor representable category. Then the following statements are equivalent:*

(*i*) *the right dualising functor admits a quasi-inverse,*

*(ii)* *the category $\mathscr{C}$ is an $r$-category, and*
*(iii)* *the category $\mathscr{C}$ is a Grothendieck–Verdier category.*

*Proof.* By definition we have that (ii) $\Longrightarrow$ (iii), and (i) $\Longrightarrow$ (ii) follows by Definition 3.6 and Remark 3.13.

To prove (iii) $\Longrightarrow$ (ii), assume that $\mathscr{C}$ is Grothendieck–Verdier with dualising object $d \in \mathscr{C}$. Write $L, R \colon \mathscr{C}^{\mathrm{op}} \longrightarrow \mathscr{C}$ for its left and right tensor dualising functors and $D_r \colon \mathscr{C}^{\mathrm{op}} \longrightarrow \mathscr{C}$ for the dualising functor induced by $d \in \mathscr{C}$. Then

$$
\begin{aligned}
Rd \otimes d &\overset{(3.1.1)}{\cong} [d,d]_r \overset{(3.1.7)}{\cong} D_r(D_r^{-1}d \otimes d) \overset{(3.1.8)}{\cong} D_r d \overset{(3.1.8)}{\cong} 1,\\
d \otimes Ld &\overset{(3.1.1)}{\cong} [d,d]_\ell \overset{(3.1.7)}{\cong} D_r^{-1}(d \otimes D_r d) \overset{(3.1.8)}{\cong} D_r^{-1} d \overset{(3.1.8)}{\cong} 1,
\end{aligned}
$$

and therefore $d \otimes Rd \cong d \otimes Rd \otimes d \otimes Ld \cong 1$. In other words, $d$ is invertible and $1 \cong Rd \otimes d$ is another dualising object due to Lemma 3.14.

To complete the proof notice that the right dualising functor of $\mathscr{C}$ considered as an $r$-category coincides with its right tensor dualising functor. □

**Example 3.17.** The converse of Proposition 3.7 is not true. Let $A$ be a commutative Frobenius algebra over a field $\Bbbk$. We write $A\text{-mod}^{\mathrm{f.g.}}$ for the category of finitely-generated left $A$-modules. One can show, see for example [Lam99, Sections 15 and 16], that for all $M \in A\text{-mod}^{\mathrm{f.g.}}$ the canonical morphism

$$
\phi_M \colon M \longrightarrow M^{\vee\vee} = \mathrm{Hom}_A(\mathrm{Hom}_A(M,A),A), \qquad \phi_M(x)\alpha = \alpha(x)
$$

is an isomorphism; the module $M$ is called *reflexive*.

Put differently, $A\text{-mod}^{\mathrm{f.g.}}$ is an $r$-category with

$$
\mathrm{Hom}_A(-,A) \colon (A\text{-mod}^{\mathrm{f.g.}})^{\mathrm{op}} \longrightarrow A\text{-mod}^{\mathrm{f.g.}}
$$

as dualising functor. It is well-known that for any $A$-module $M$ there exists an $A$-module $N$ such that $M \otimes - \dashv N \otimes -$ if and only if $M$ is finitely-generated projective; see for example [NW17, Proposition 2.1]. Thus, $A\text{-mod}^{\mathrm{f.g.}}$ is tensor representable if and only if all (finitely-generated) $A$-modules are projective, which is equivalent to $A$ being semisimple.

However, not all commutative Frobenius algebras are semisimple. For example, the group algebra $\Bbbk[G]$ of a finite abelian group $G$ is Frobenius. By Maschke's theorem, it is semisimple if and only if the characteristic of $\Bbbk$ does not divide the order of $G$.

**Remark 3.18.** Let $\mathscr{C}$ be a right tensor representable category whose dualising functor $R\colon \mathscr{C}^{\mathrm{op}} \longrightarrow \mathscr{C}$ is an equivalence. Following [BD13, Section 4], for all $x, y \in \mathscr{C}$ there exists a canonical isomorphism

$$\begin{aligned}\tau_{x,y}\colon \mathscr{C}(y \otimes x \otimes Rx \otimes Ry, 1) &\xrightarrow{\sim} \mathscr{C}(Rx \otimes Ry, R(y \otimes x)) \\ &\xrightarrow{\sim} \mathscr{C}^{\mathrm{op}}(y \otimes x, R^{-1}(Rx \otimes Ry)).\end{aligned}$$

As shown in [BD13, Lemma 4.1 and Remark 4.2], the morphisms

$$f := \tau_{x,y}\big(\varepsilon_1^{(y)} \circ (y \otimes \varepsilon_1^{(x)} \otimes Ry)\big)\colon R^{-1}(Rx \otimes Ry) \longrightarrow y \otimes x$$

$$\mu_{x,y} := Rf\colon Rx \otimes Ry \longrightarrow R(y \otimes x) = R(x \otimes^{\mathrm{op}} y)$$

endow $R\colon \mathscr{C}^{\mathrm{op,rev}} \longrightarrow \mathscr{C}$ with the structure of a lax monoidal functor:

$$\begin{aligned}&\mu_{x,y\otimes^{\mathrm{op}}z} \circ (Rx \otimes \mu_{y,z}) = \mu_{x\otimes^{\mathrm{op}}y,z} \circ (\mu_{x,y} \otimes Rz),\\ &\mu_{1,x} \circ (R1 \otimes Rx) = \mathrm{id}_{Rx} = \mu_{x,1} \circ (Rx \otimes R1), \qquad \text{for all } x, y, z \in \mathscr{C}.\end{aligned}$$

In the situation of the above remark and as a consequence of [BD13, Proposition 4.4], we obtain the following characterisation of rigidity.

**Proposition 3.19.** *Let the dualising functor $R\colon \mathscr{C}^{\mathrm{op,rev}} \longrightarrow \mathscr{C}$ of a right tensor representable category be an equivalence. Then $\mathscr{C}$ is rigid if and only if the morphism $\mu_{x,y}\colon R(-) \otimes R(=) \Longrightarrow R(- \otimes^{\mathrm{op}} =)$ is an isomorphism, making $R$ strong monoidal.*

### 3.1.2 *The free tensor representable category is not rigid*

The aim of this section is to concretise Remark 2.68. That is, we explicitly construct a "free" tensor representable category $\mathscr{W}$ and show that it is not rigid; see Theorem 3.23. Due to the reliance of the proof on some technical arguments, we will first provide an overview of the underlying strategy.

As an intermediate step in the construction of $\mathscr{W}$, we define a strict monoidal category $\mathscr{V}$, which has the same generators as $\mathscr{W}$ but less relations. In other words, $\mathscr{W}$ itself will be a quotient of $\mathscr{V}$. The category $\mathscr{V}$ will have the property that all presentations of a given morphism have the same length. In order to show that $\mathscr{W}$ is not rigid, we will fix an object $t \in \mathscr{W}$ and assume it admits a rigid dual. In particular, $t$ would have to admit evaluation and coevaluation morphisms $\mathrm{ev}_t\colon t \otimes t^\vee \longrightarrow 1$ and $\mathrm{coev}_t\colon 1 \longrightarrow t^\vee \otimes t$. By showing that any morphism in $\mathscr{V}$ representing the class $(\mathrm{ev}_t \otimes t) \circ (t \otimes \mathrm{coev}_t) \in \mathscr{W}(t,t)$ has to have length at least 2, we conclude that the snake identities cannot hold. Hence, $\mathscr{W}$ is not rigid.

We define the strict monoidal category $\mathcal{V}$ in terms of generators and relations. For details of this type of construction we refer to [Kas98, Chapter xii].

The monoid of objects of the category $\mathcal{V}$ is freely generated by the alphabet $\{\dots, R^{-1}t, t, Rt, \dots\}$. That is, its elements are words of finite length whose letters are $R^n t$ for some integer $n \in \mathbb{Z}$. Write $1 \in \mathrm{Ob}(\mathcal{V})$ for the monoidal unit of $\mathcal{V}$, given by the empty word. For an object $w := (R^{k_1}t, \dots, R^{k_n}t) \in \mathrm{Ob}(\mathcal{V})$, set $Rw := (R^{k_n+1}t, \dots, R^{k_1+1}t) \in \mathrm{Ob}(\mathcal{V})$. This assignment extends to a group action of the free abelian group $\langle R \rangle \cong (\mathbb{Z}, +)$ on $\mathrm{Ob}(\mathcal{V})$. For any $i \in \mathbb{Z}$, we write $R^i w$ for the action of $R^i$ on $w \in \mathrm{Ob}(\mathcal{V})$.

The morphisms of $\mathcal{V}$ are tensor products and compositions of identities and the *generating morphisms*. These are given for all $x \neq 1, y \in \mathrm{Ob}(\mathcal{V})$ by

$$\eta_y^{(x)} \colon y \longrightarrow Rx \otimes x \otimes y, \qquad \varepsilon_y^{(x)} \colon x \otimes Rx \otimes y \longrightarrow y,$$
$$u_y^{(x)} \colon y \longrightarrow y \otimes x \otimes R^{-1}x, \qquad c_y^{(x)} \colon y \otimes R^{-1}x \otimes x \longrightarrow y.$$

These relations are tailored to implement natural transformations

$$\eta^{(x)} \colon \mathrm{Id}_{\mathcal{V}} \Longrightarrow Rx \otimes x \otimes -, \qquad \varepsilon^{(x)} \colon x \otimes Rx \otimes - \Longrightarrow \mathrm{Id}_{\mathcal{V}},$$
$$u^{(x)} \colon \mathrm{Id}_{\mathcal{V}} \Longrightarrow - \otimes x \otimes R^{-1}x, \qquad c^{(x)} \colon - \otimes R^{-1}x \otimes x \Longrightarrow \mathrm{Id}_{\mathcal{V}}.$$

Explicitly, for any generating morphism $g \colon c \longrightarrow d$, any $y = a \otimes c \otimes b$, $z = a \otimes d \otimes b$, any arrow $f = (\mathrm{id}_a \otimes g \otimes \mathrm{id}_b) \colon y \longrightarrow z$, and all $x \in \mathrm{Ob}(\mathcal{V}) \setminus \{1\}$, we require the following conditions to hold:

$$\eta_z^{(x)} \circ f = (Rx \otimes x \otimes f) \circ \eta_y^{(x)}, \qquad \varepsilon_z^{(x)} \circ (x \otimes Rx \otimes f) = f \circ \varepsilon_y^{(x)},$$
$$u_z^{(x)} \circ f = (f \otimes x \otimes R^{-1}x) \circ u_y^{(x)}, \qquad c_z^{(x)} \circ (f \otimes R^{-1}x \otimes x) = f \circ c_y^{(x)}.$$

We obtain the tensor representable category $\mathcal{W}$ by quotienting out the triangle identities.[6] Explained in more detail, the strict monoidal category $\mathcal{W}$ has the same objects and generating morphisms as $\mathcal{V}$, and the same identities hold. In addition, for any $x, y \in \mathrm{Ob}(\mathcal{W})$ with $x \neq 1$, we require

[6] The category is obtained by *monoidal localisation*, see [Day73].

(3.1.10)
$$\varepsilon_{x\otimes y}^{(x)} \circ (x \otimes \eta_y^{(x)}) = \mathrm{id}_x, \qquad (Rx \otimes \varepsilon_y^{(x)}) \circ \eta_{Rx\otimes y}^{(x)} = \mathrm{id}_{Rx},$$
$$c_{y\otimes x}^{(x)} \circ (u_y^{(x)} \otimes x) = \mathrm{id}_x, \qquad (c_y^{(x)} \otimes R^{-1}x) \circ u_{y\otimes R^{-1}x}^{(x)} = \mathrm{id}_{R^{-1}x},$$
$$(\varepsilon_1^{(x)} \otimes x) \circ u_x^{(Rx)} = \mathrm{id}_x, \qquad (x \otimes c_1^{(x)}) \circ \eta_x^{(R^{-1}x)} = \mathrm{id}_x.$$

The next result succinctly summarises the observations made so far concerning the internal hom of $\mathcal{W}$.

**Lemma 3.20.** *The category $\mathcal{W}$ is tensor representable. Furthermore, its right dualising functor $R\colon \mathcal{W}^{\mathrm{op}} \longrightarrow \mathcal{W}$ is an equivalence of categories.*

*Proof.* By construction, we have $x \otimes - \dashv Rx \otimes -$ and $- \otimes x \dashv - \otimes R^{-1}x$ for all $x \in \mathcal{W}$. Thus, $\mathcal{W}$ is tensor representable. Equation (3.1.10) and Proposition 3.7 imply that $R\colon \mathcal{W}^{\mathrm{op}} \longrightarrow \mathcal{W}$ and $R^{-1}\colon \mathcal{W}^{\mathrm{op}} \longrightarrow \mathcal{W}$ are quasi-inverses. □

To analyse the morphisms in $\mathcal{W}$ and show that it is not rigid monoidal, we will rely on two tools. The first is the *length* of an arrow $f \in \mathcal{W}(x, y)$. It is defined as the minimal number of generating morphisms needed to present $f$. The second tool will be given by invariants for morphisms in $\mathcal{W}$ arising from functors into the rigid monoidal category $\mathsf{vect}$ of finite-dimensional vector spaces over a field $\Bbbk$. We will write $(-)^*\colon \mathsf{vect}^{\mathrm{op}} \longrightarrow \mathsf{vect}$ for the right rigid dualising functor of $\mathsf{vect}$. The $n$-fold dual of a finite-dimensional vector space $V$ will be denoted by $V^{(n)}$.

**Lemma 3.21.** *For all $V \in \mathsf{vect}$ there is a strong monoidal functor $F_V\colon \mathcal{W} \longrightarrow \mathsf{vect}$ satisfying $F_V(R^n t) \cong V^{(n)}$ for all $n \in \mathbb{Z}$ and*

$$F_V(\eta_y^{(x)}) = \mathrm{coev}^r_{F_V x} \otimes F_V y, \qquad F_V(\varepsilon_y^{(x)}) = \mathrm{ev}^r_{F_V x} \otimes F_V y,$$
$$F_V(u_y^{(x)}) = F_V y \otimes \mathrm{coev}^\ell_{F_V x}, \qquad F_V(c_y^{(x)}) = F_V y \otimes \mathrm{ev}^\ell_{F_V x},$$

*for suitably chosen evaluation and coevaluation morphisms.*

*Proof.* Theorem 2.72 yields a monoidal equivalence $G\colon \mathsf{vect} \rightleftarrows \mathcal{A} : H$ between the category of finite-dimensional vector spaces and a rigid monoidal category $\mathcal{A}$ whose dualising functors $D_\ell, D_r\colon \mathcal{A}^{\mathrm{op}} \longrightarrow \mathcal{A}$ satisfy $D_\ell D_r^{\mathrm{op}} = \mathrm{Id}_{\mathcal{A}} = D_r^{\mathrm{op}} D_\ell$. Fix a quasi-inverse $H\colon \mathcal{A} \longrightarrow \mathsf{vect}$ and choose an object $v \in \mathcal{A}$ with $H(v) \cong V$, for some $V \in \mathsf{vect}$. A direct computation using the universal property of $\mathcal{W}$, see for example [Kas98, Proposition XII.1.4], shows that for any object $v \in \mathcal{A}$ there is a unique strict monoidal functor $F'_v\colon \mathcal{W} \longrightarrow \mathcal{A}$ such that $F'_v(t) = v$ and for all $x, y \in \mathcal{W}$, $x \neq 1$, we have

$$F'_v \eta_y^{(x)} = \mathrm{coev}^r_{F'_v x} \otimes F'_v y, \qquad F'_v \varepsilon_y^{(x)} = \mathrm{ev}^r_{F'_v x} \otimes F'_v y,$$
$$F'_v u_y^{(x)} = F'_v y \otimes \mathrm{coev}^\ell_{F'_v x}, \qquad F'_v c_y^{(x)} = F'_v y \otimes \mathrm{ev}^\ell_{F'_v x}.$$

The claim now follows by setting $F = HF'_{G(V)}$. □

**Corollary 3.22.** *The following arrows cannot be isomorphisms, for any* $1 \neq x \in \mathcal{W}$ *and any morphism* $g \in \mathcal{W}(a_1 \otimes y \otimes b_1, a_2 \otimes z \otimes b_2)$:

$$(\mathrm{id}_{a_2} \otimes u_z^{(x)} \otimes \mathrm{id}_{b_2}) \circ g, \qquad g \circ (\mathrm{id}_{a_1} \otimes \varepsilon_y^{(x)} \otimes \mathrm{id}_{b_1}).$$

*Proof.* Suppose that $g \in \mathcal{W}(a_1 \otimes y \otimes b_1, a_2 \otimes z \otimes b_2)$ and consider

$$f := g \circ (\mathrm{id}_{a_1} \otimes \varepsilon_y^{(x)} \otimes \mathrm{id}_{b_1}).$$

Let $V \in \mathsf{vect}$ be a vector space of dimension at least two. Applying $F_V$ to $f$, we get $F_V f = F_V g \circ F_V(\mathrm{id}_{a_1} \otimes \varepsilon_y^{(x)} \otimes \mathrm{id}_{b_1})$. However, due to the difference in the dimensions of its source and target, $F_V(\mathrm{id}_{a_1} \otimes \varepsilon_y^{(x)} \otimes \mathrm{id}_{b_1})$ must have a non-trivial kernel and thus $f$ cannot be an isomorphism.

A similar argument involving the cokernel proves that the composition $(\mathrm{id}_{a_2} \otimes u_z^{(x)} \otimes \mathrm{id}_{b_2}) \circ g$ is not invertible. □

**Theorem 3.23.** *The category* $\mathcal{W}$ *is not rigid.*

*Proof.* Assume that $t \in \mathcal{W}$ admits a right rigid dual $x \in \mathcal{W}$. By the uniqueness of adjoints, there exist isomorphisms $\vartheta\colon Rt \otimes t \longrightarrow x \otimes t$ and $\theta\colon x \longrightarrow Rt$ such that the evaluation and coevaluation morphisms are given by

$$\mathrm{coev} := \vartheta \circ \eta_1^{(t)}\colon 1 \longrightarrow x \otimes t, \qquad \mathrm{ev} := \varepsilon_1^{(t)} \circ (t \otimes \theta)\colon t \otimes x \longrightarrow 1.$$

Let $\pi\colon \mathcal{V} \longrightarrow\!\!\!\!\!\rightarrow \mathcal{W}$ be the projection functor. We consider the set

$$S := \left\{ (\varepsilon_1^{(t)} \otimes t)\, \phi\, (t \otimes \eta_1^{(t)}) \in \mathcal{V}(t,t) \;\middle|\; \begin{array}{c} \phi \in \mathcal{V}(t \otimes Rt \otimes t, t \otimes Rt \otimes t) \\ \text{such that } \pi(\phi) \text{ is invertible} \end{array} \right\}.$$

Notice that $S \subseteq \mathrm{Mor}(\mathcal{V})$. By construction, there exists a morphism $s \in S$ such that the composite arrow $(\mathrm{ev} \otimes t) \circ (t \otimes \mathrm{coev}) = \pi(s)$ corresponds to one of the snake identities. Furthermore, every element of $S$ has length at least two.[7] Thus, by proving that $S$ is closed under the relations arising from Equation (3.1.10), it follows that $\pi(s) \neq \mathrm{id}_t$, which concludes the proof.

[7] Note that the relations of $\mathcal{V}$ leave the number of generating morphisms in any presentation of a given arrow invariant.

Let us consider an arbitrary element $f := (\varepsilon_1^{(t)} \otimes t) \circ \phi \circ (t \otimes \eta_1^{(t)}) \in S$. There are two types of "moves" we have to study. First, suppose we expand an identity into one of the morphisms displayed in Equation (3.1.10). This equates to either pre or postcomposing $\phi$ with an arrow $\psi \in \mathcal{V}(t \otimes Rt \otimes t, t \otimes Rt \otimes t)$ that projects onto an isomorphism in $\mathcal{W}$, leading to another element in $S$.

Second, any of the morphisms of Equation (3.1.10) might be contracted to an identity. A priori, there are three ways in which this might occur:

$$f = (\varepsilon_1^{(t)} \otimes t) \circ \phi' \circ \varepsilon_t^{(t)} \circ (t \otimes \eta_1^{(t)}), \quad \text{with } \phi = \phi' \circ \varepsilon_t^{(t)}, \text{ or} \tag{3.1.11}$$

$$f = (\varepsilon_1^{(t)} \otimes t) \circ u_t^{(Rt)} \circ \phi'' \circ (t \otimes \eta_1^{(t)}), \quad \text{with } \phi = u_t^{(Rt)} \circ \phi'', \text{ or} \tag{3.1.12}$$

$$f = (\varepsilon_1^{(t)} \otimes t) \circ \phi_2 \circ g \circ \phi_1 \circ (t \otimes \eta_1^{(t)}), \text{ with } \phi = \phi_2 \circ g \circ \phi_1 \text{ and } \pi(g) = \mathrm{id}. \tag{3.1.13}$$

Due to Corollary 3.22, neither $\pi(\phi') \circ \pi(\varepsilon_t^{(t)})$ nor $\pi(u_t^{(Rt)}) \circ \pi(\phi'')$ are isomorphisms, contradicting Cases (3.1.11) and (3.1.12).

Now, assume $f = (\varepsilon_1^{(t)} \otimes t) \circ \phi_2 \circ g \circ \phi_1 \circ (t \otimes \eta_1^{(t)})$ and $\phi = \phi_2 \circ g \circ \phi_1$. Using the functoriality of $\pi \colon \mathcal{V} \longrightarrow \mathcal{W}$, we get

$$\pi(\phi) = \pi(\phi_2 \circ g \circ \phi_1) = \pi(\phi_2) \circ \pi(g) \circ \pi(\phi_1) = \pi(\phi_2) \circ \pi(\phi_1) = \pi(\phi_2 \circ \phi_1).$$

Thus, $\pi(\phi_2 \circ \phi_1)$ is an isomorphism and $(\varepsilon_1^{(t)} \otimes t) \circ \phi_2 \circ \phi_1 \circ (t \otimes \eta_1^{(t)})$ is an element of $S$. □

## 3.2 FUNCTOR CATEGORIES

WE WILL NOW INVESTIGATE the previously discussed types of dualities in the context of functor categories. The starting point is the fact that Day convolution is closed monoidal, see Theorem 2.124. Examples of such categories of functors are abound, and include modules over path algebras, Mackey functors, and group-graded representations arising from crossed modules; we shall more closely investigate these examples in Section 3.3.

**Hypothesis 3.24.** Fix a field $\Bbbk$, and write $(-)^* \colon \mathsf{vect}^{\mathrm{op}} \longrightarrow \mathsf{vect}$ for the $\Bbbk$-linear dualising functor. Unless explicitly stated, for the remainder of this section all (monoidal) categories and functors will be enriched over $\mathsf{Vect}$, and $\mathscr{C}$ will denote an essentially small monoidal $\Bbbk$-linear category. For the sake of brevity, the prefixes "enriched", "$\Bbbk$-linear", and "(essentially) small" will often be omitted. The category of functors from $\mathscr{C}$ to $\mathsf{Vect}$ will be denoted by $\widehat{\mathscr{C}^{\mathrm{op}}}$.

Instead of studying functors $F \in \widehat{\mathscr{C}^{\mathrm{op}}} = [\mathscr{C}, \mathsf{Vect}]$ from a hom-finite category $\mathscr{C}$ to *all* vector spaces, we are interested in the full subcategory $\widehat{\mathscr{C}^{\mathrm{op}}}_{\mathrm{fin}}$ of (point-wise) *finite-dimensional* functors. That is, $F \in \widehat{\mathscr{C}^{\mathrm{op}}}$ such that $Fx$ is finite-dimensional for all $x \in \mathscr{C}$.

**Hypothesis 3.25.** Within this section, we assume that all relevant ends and coends are finite-dimensional vector spaces.

In particular, Hypothesis 3.25 implies that $\widehat{\mathscr{C}^{\mathrm{op}}}_{\mathrm{fin}}$ is closed monoidal; that is, the tensor product and the internal hom are finite-dimensional at every point. These assumptions impose a certain finiteness condition on $\mathscr{C}$ itself, or—as the next example will show—on a dense[8] subcategory of it.

[8] A subcategory $\mathfrak{S}$ of $\mathscr{C}$ is *dense* if the restricted Yoneda embedding from $\mathscr{C}$ to $[\mathfrak{S}^{\mathrm{op}}, \mathsf{Vect}]$ is fully faithful.

**Example 3.26.** Let $\mathscr{C}$ be a hom-finite category, and assume that there is a full subcategory $\mathcal{S}$ such that every object of $\mathscr{C}$ may be written as a finite direct sum of objects in $\mathcal{S}$. In this setting, one can reformulate Lemma 2.120 to index over this finite set: for example, given $F\colon \mathscr{C} \longrightarrow \mathsf{vect}$, we have

$$Fx \cong F\Big(\bigoplus_{i=1}^{n} s_i\Big) \cong \bigoplus_{i=1}^{n} Fs_i \cong \bigoplus_{i=1}^{n} \int^{s\in\mathcal{S}} \mathcal{S}(s, s_i) \otimes Fs \cong \bigoplus_{i=1}^{n} \int^{s\in\mathcal{S}} \mathscr{C}(s, s_i) \otimes Fs$$
$$\cong \int^{s\in\mathcal{S}} \mathscr{C}\Big(s, \bigoplus_{i=1}^{n} s_i\Big) \otimes Fs \cong \int^{s\in\mathcal{S}} \mathscr{C}(s, x) \otimes Fs.$$

Panchadcharam and Street used this to show that finite-dimensional Mackey functors are a Grothendieck–Verdier category, see [PS07, Section 9].

Hypothesis 3.25 does not only imply a closed structure on $\widehat{\mathscr{C}^{\mathrm{op}}}_{\mathrm{fin}}$, but yields a canonical notion of a dual.

**Proposition 3.27.** *Let $(\mathscr{C}, d)$ be a hom-finite left Grothendieck–Verdier category with dualising functor $D_\ell$. Then $\widehat{\mathscr{C}^{\mathrm{op}}}_{\mathrm{fin}}$ is a right Grothendieck–Verdier category with dualising object $\mathscr{C}(-, d)^* \in \widehat{\mathscr{C}^{\mathrm{op}}}_{\mathrm{fin}}$ and dualising functor*

$$\mathfrak{D}_r \colon \widehat{\mathscr{C}^{\mathrm{op}}}_{\mathrm{fin}}^{\mathrm{op}} \longrightarrow \widehat{\mathscr{C}^{\mathrm{op}}}_{\mathrm{fin}}, \qquad F \longmapsto F(D_\ell -)^*.$$

*Proof.* Using Equation (3.1.9), one can recover the right dualising functor of a Grothendieck–Verdier category from its right internal hom by evaluating it at the dualising object; i.e., $\mathfrak{D}_r \cong [-, \mathscr{C}(-, d)^*]_r$. In our case, we have

$$\begin{aligned}
[F, \mathscr{C}(-, d)^*]_r x &\overset{(2.8.7)}{\cong} \int_b \mathsf{vect}\big(F[x,b]_\ell, \mathsf{vect}(\mathscr{C}(b,d), \Bbbk)\big) \\
&\cong \int_b \mathsf{vect}\big(\mathscr{C}(b,d) \otimes F[x,b]_\ell, \Bbbk\big) \\
&\cong \mathsf{vect}\Big(\int^b \mathscr{C}(b,d) \otimes F[x,b]_\ell, \Bbbk\Big) \overset{(2.8.1)}{\cong} \mathsf{vect}\big(F[x,d]_\ell, \Bbbk\big) \\
&\overset{(3.1.9)}{\cong} F(D_\ell x)^*.
\end{aligned}$$

To conclude the proof, notice that $\mathfrak{D}_r \colon \widehat{\mathscr{C}^{\mathrm{op}}}_{\mathrm{fin}}^{\mathrm{op}} \longrightarrow \widehat{\mathscr{C}^{\mathrm{op}}}_{\mathrm{fin}}$ has a quasi-inverse

$$\mathfrak{D}_r^{-1} \colon \widehat{\mathscr{C}^{\mathrm{op}}}_{\mathrm{fin}} \longrightarrow \widehat{\mathscr{C}^{\mathrm{op}}}_{\mathrm{fin}}^{\mathrm{op}}, \qquad F \longmapsto F(D_\ell^{-1}-)^{*}. \qquad \square$$

While we will develop efficient means to detect tensor representability and rigidity based on properties which exist in the abelian case, it is worthwhile to explore how such structures arise as an interplay between the duality type of the base category and that of vect.

**Corollary 3.28.** *Let $(\mathscr{C}, d)$ be a hom-finite left Grothendieck–Verdier category with dualising functor $D_\ell$. Then the category $\widehat{\mathscr{C}^{\mathrm{op}}}_{\mathrm{fin}}$ is right tensor representable if there are isomorphisms*

Not all assumptions are strictly needed, see Remark 3.34.

$$D_\ell^2 a \cong a, \qquad D_\ell(a \otimes x) \cong D_\ell x \otimes D_\ell a, \qquad \textit{natural in } a, x \in \mathscr{C}, \tag{3.2.1}$$

*and for all objects $F, G \in \widehat{\mathscr{C}^{\mathrm{op}}}_{\mathrm{fin}}$ we have*

$$\int^{a} \mathsf{vect}(F[-,a]_\ell, Ga) \cong \int_{a} \mathsf{vect}(F[-,a]_\ell, Ga). \tag{3.2.2}$$

*Proof.* Let $\mathfrak{D}_r$ be the dualising functor of Proposition 3.27, and fix functors $F, G \in \widehat{\mathscr{C}^{\mathrm{op}}}_{\mathrm{fin}}$. Then, by the following calculation, we obtain an isomorphism $(\mathfrak{D}_r F * G)x \cong [F,G]_r x$, and hence the claim follows:

$$\begin{aligned}
(\mathfrak{D}_r F * G)x &\overset{(2.8.3)}{\cong} \int^{a} \mathfrak{D}_r F([a,x]_\ell) \otimes Ga \overset{(3.1.7)}{\cong} \int^{a} \mathfrak{D}_r F(D_\ell(a \otimes D_\ell x)) \otimes Ga \\
&\overset{3.27}{\cong} \int^{a} F(D_\ell^2(a \otimes D_\ell x))^{*} \otimes Ga \overset{(3.2.1)}{\cong} \int^{a} \mathsf{vect}\big(F(a \otimes D_\ell x), Ga\big) \\
&\overset{(3.2.1)}{\cong} \int^{a} \mathsf{vect}\big(F(D_\ell^2 a \otimes D_\ell x), Ga\big) \overset{(3.2.1)}{\cong} \int^{a} \mathsf{vect}\big(F(D_\ell(x \otimes D_\ell a)), Ga\big) \\
&\overset{(3.1.7)+(3.2.1)}{\cong} \int^{a} \mathsf{vect}\big(F([x,a]_\ell), Ga\big) \overset{(3.2.2)}{\cong} \int_{a} \mathsf{vect}\big(F([x,a]_\ell), Ga\big) \\
&\overset{(2.8.8)}{=} [F,G]_r x. \qquad \square
\end{aligned}$$

**Remark 3.29.** Before we discuss conditions for the interchangeability of ends and coends, let us briefly mention some cases where the dualising functor of a left Grothendieck–Verdier category $(\mathscr{C}, d)$ admits natural isomorphisms as stated in Equation (3.2.1).

The requirement $D_\ell(x \otimes y) \cong D_\ell y \otimes D_\ell x$ implies that $\mathscr{C}$ must be left tensor representable, since then for all $x, y, z \in \mathscr{C}$ we have

$$\mathscr{C}(x \otimes y, z) \overset{(3.1.6)}{\cong} \mathscr{C}(x \otimes y \otimes D_\ell^{-1} z, d) \overset{(3.1.6)}{\cong} \mathscr{C}(x, D_\ell(y \otimes D_\ell^{-1} z)) \overset{(3.2.1)}{\cong} \mathscr{C}(x, z \otimes D_\ell y).$$

In the absence of further coherence assumptions, this does not imply that $\mathscr{C}$ must be rigid. The condition $D_\ell^2 \cong \mathrm{Id}_{\mathscr{C}}$ is met for example if $\mathscr{C}$ is braided. For all $x, y \in \mathscr{C}$ we have

$$\mathscr{C}(x, y) \overset{(3.1.6)}{\cong} \mathscr{C}(D_\ell y \otimes x, d) \cong \mathscr{C}(x \otimes D_\ell y, d) \overset{(3.1.6)}{\cong} \mathscr{C}(x, D_\ell^2 y),$$

and the Yoneda lemma implies that $\mathrm{Id}_{\mathscr{C}} \cong D_\ell^2$.

The following constructions are an adaptation of [Day06]. For simplicity, we replaced the promonoidal structure in *ibid* with the one induced from a monoidal structure. That is, $J := \mathscr{C}(1, -)$, $P(x, y, z) := \mathscr{C}(x \otimes y, z)$, for $J$ and $P$ the promonoidal unit and multiplication.

**Remark 3.30.** Let $\mathscr{C}$ be a category. Given any functor $T\colon \mathscr{C}^{\mathrm{op}} \otimes_{\Bbbk} \mathscr{C} \longrightarrow \mathsf{Vect}$, there exists a canonical map

$$\int^{a \in \mathscr{C}} \int_{b \in \mathscr{C}} \mathscr{C}(b, a) \otimes T(a, b) \longrightarrow \int_{b \in \mathscr{C}} \int^{a \in \mathscr{C}} \mathscr{C}(b, a) \otimes T(a, b). \tag{3.2.3}$$

This follows from the fact that ends and coends are functorial; in particular, for any arrow $f\colon x \longrightarrow y$, the following diagram commutes:

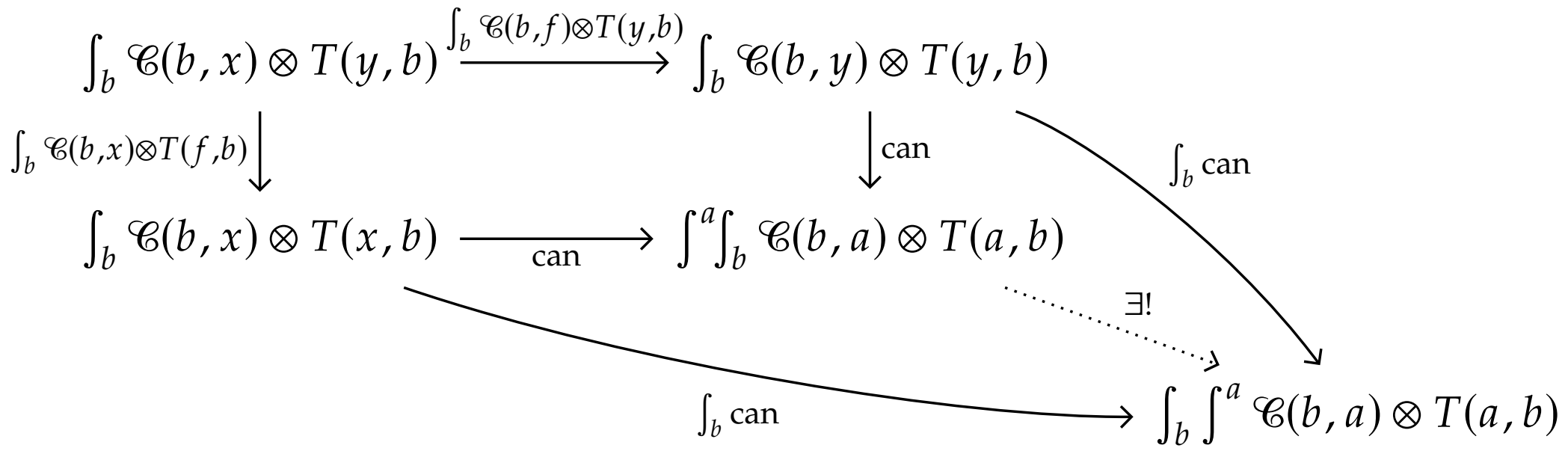

**Lemma 3.31** ([Day06, p. 1])**.** *Let $\mathscr{C}$ be a hom-finite category, $P\colon \mathscr{C}^{\mathrm{op}} \otimes_{\Bbbk} \mathscr{C} \longrightarrow \mathsf{vect}$ a functor, and suppose that the map given in Equation (3.2.3) is invertible. If for all $a, b \in \mathscr{C}$ there is a natural isomorphism $\varphi\colon \mathscr{C}(b, a) \xrightarrow{\sim} \mathscr{C}(a, b)^*$, then*

$$\int^a P(a, a) \xrightarrow{\sim} \int_b P(b, b).$$

*Proof.* We calculate

$$\begin{aligned}\int^a P(a,a) &\cong \int^a\int_b \mathsf{vect}(\mathscr{C}(a,b), P(a,b)) \cong \int^a\int_b \mathscr{C}(a,b)^* \otimes P(a,b) \\ &\overset{\varphi^{-1}}{\cong} \int^a\int_b \mathscr{C}(b,a)\otimes P(a,b) \overset{(3.2.3)}{\cong} \int_b\int^a \mathscr{C}(b,a)\otimes P(a,b) \cong \int_b P(b,b).\end{aligned}$$

The first and last isomorphisms follow from Lemma 2.120 and the second one is a consequence of $\mathscr{C}(a,b)$ being finite-dimensional and $\mathsf{vect}$ being rigid. □

In general, it seems difficult to decide whether a natural isomorphism as that in Lemma 3.31 exists. However, for certain classes of examples, trace maps provide us with viable candidates.

**Remark 3.32.** Let $\mathscr{C}$ be a hom-finite and pivotal category with pivotal structure $\psi\colon (-)^\vee \longrightarrow {}^\vee(-)$, and suppose that $\mathscr{C}(1,1)\cong \Bbbk$. For all $a,b\in\mathscr{C}$, consider

$$\begin{aligned}\varphi_{a,b}\colon \mathscr{C}(b,a) &\longrightarrow \mathscr{C}(a,b)^* \\ f &\longmapsto \big(g \longmapsto \mathrm{tr}(f\circ g) := \mathrm{ev}^\ell_a \circ (\psi_a \otimes (f\circ g))\circ \mathrm{coev}^r_a\big).\end{aligned} \tag{3.2.4}$$

The dual of $\varphi$ is given by

$$\begin{aligned}(\varphi_{a,b})^*\colon \mathscr{C}(a,b)^{**}\cong \mathscr{C}(a,b) &\longrightarrow \mathscr{C}(b,a)^* \\ g &\longmapsto \big(f \longmapsto \mathrm{tr}(f\circ g) = \mathrm{tr}(g\circ f) = \varphi_{b,a}(g)(f)\big).\end{aligned}$$

Thus, $\varphi_{a,b}$ is injective if and only if $\varphi_{b,a}$ is surjective.

Suppose $\Bbbk$ has characteristic zero and $\mathscr{C} := H\text{-mod}$ is the category of finite-dimensional modules of a semisimple Hopf algebra $H$. As a consequence of [LR88, Theorem 4], left and right duals coincide and we can chose the "quantum trace" of Equation (3.2.4) to agree with the usual trace of endomorphisms between finite-dimensional vector spaces. Let $f\colon M\longrightarrow N$ be a morphism in $\mathscr{C}$. By semisimplicity, the following short exact sequence splits:

$$0 \longrightarrow \ker f \longrightarrow M \underset{\iota}{\overset{f}{\rightleftarrows}} \operatorname{im} f \longrightarrow 0$$

The arrow $g = \iota\oplus 0\colon N\cong \operatorname{im} f\oplus N/\operatorname{im} f \longrightarrow M$ satisfies

$$\mathrm{tr}(fg) = \mathrm{tr}(f\iota) = \mathrm{tr}(\mathrm{id}_{\operatorname{im} f}) = \dim \operatorname{im} f,$$

implying that $\varphi_{M,N}$ is injective.

**Example 3.33.** Consider a field $\Bbbk$ of characteristic $p$, as well as $q \in \mathbb{N}$ a multiple of $p$. The group algebra $H := \Bbbk[\mathrm{GL}_q(p)]$ of the group of invertible $q \times q$-matrices over $\Bbbk$ provides us with an example of a pivotal hom-finite category where the morphism of Equation (3.2.4) is not invertible. Matrix-vector multiplication turns $\Bbbk^q$ into a simple $H$-module whose endomorphism algebra is one-dimensional. This implies that for all $f, g \in \mathrm{End}_{\mathrm{GL}_q(p)}(\Bbbk^p)$ there exists some $\lambda \in \Bbbk$ such that $\mathrm{tr}(fg) = \lambda\, \mathrm{tr}(\mathrm{id}_{\Bbbk^q}) = \lambda \cdot 0 = 0$.

**Remark 3.34.** One can use a combination of Corollary 3.28 and Lemma 3.31 to deduce that a given category of finite-dimensional functors is right tensor representable. To reduce the number of axioms to be checked, note that the square of the left dualising functor being isomorphic to the identity is implied by the other assumptions. Let $(\mathscr{C}, d)$ be a hom-finite left Grothendieck–Verdier category with dualising functor $D_\ell$, such that there are isomorphisms $\mathscr{C}(a, b) \cong \mathscr{C}(b, a)^*$ and $D_\ell(a \otimes b) \cong D_\ell b \otimes D_\ell a$, natural in $a, b \in \mathscr{C}$. Using that

$$\begin{aligned}\mathscr{C}(d, a \otimes D_\ell b) &\overset{(3.1.6)}{\cong} \mathscr{C}(d \otimes D_\ell^{-1}(a \otimes D_\ell b), d) \cong \mathscr{C}(d \otimes b \otimes D_\ell^{-1} a, d)\\ &\overset{(3.1.6)}{\cong} \mathscr{C}(D_\ell^{-1} a, D_\ell^{-1}(d \otimes b)) \overset{(3.1.8)}{\cong} \mathscr{C}(D_\ell^{-1} a, D_\ell^{-1} b) \cong \mathscr{C}(b, a),\end{aligned}$$

one obtains

$$\mathscr{C}(a, D_\ell^2 b) \overset{(3.1.6)}{\cong} \mathscr{C}(a \otimes D_\ell b, d) \cong \mathscr{C}(d, a \otimes D_\ell b)^* \cong \mathscr{C}(b, a)^* \cong \mathscr{C}(a, b)^{**} \cong \mathscr{C}(a, b).$$

The claim follows from the Yoneda lemma.

**Lemma 3.35.** *Assume $\Bbbk$ to be a perfect field. If the hom-finite category $\mathscr{C}$ has finitely many objects and $\widehat{\mathscr{C}^{\mathrm{op}}}_{\mathrm{fin}}$ is semisimple, then the canonical map*

$$\int^a \int_b T(b, b, a, a) \cong \int_b \int^a T(a, a, b, b)$$

*is invertible for all functors $T\colon \mathscr{C}^{\mathrm{op}} \otimes_\Bbbk \mathscr{C} \otimes_\Bbbk \mathscr{C}^{\mathrm{op}} \otimes_\Bbbk \mathscr{C} \longrightarrow \mathsf{vect}$.*

*Proof.* We endow the vector space $A := \bigoplus_{a,b \in \mathscr{C}} \mathscr{C}(a, b)$ with the structure of an associative unital algebra via the multiplication following multiplication, for $f \in \mathscr{C}(a, b)$, $g \in \mathscr{C}(c, d)$:

$$f \cdot g := \begin{cases} g \circ f, & \text{if } \mathrm{b} = \mathrm{c},\\ 0, & \text{otherwise.}\end{cases}$$

The unit $\sum_{a\in\mathscr{C}} \mathrm{id}_a$ of $A$ is a sum of orthogonal idempotents. Thus, any right module $M$ of $A$ decomposes as a vector space into a direct sum $M \cong \bigoplus_{a\in\mathscr{C}} M_a$, where $M_a = M \cdot \mathrm{id}_a$, and the action of any $f \in \mathscr{C}(a,b)$ defines a linear map $M_a \longrightarrow M_b$. Accordingly, a morphism of right $A$-modules corresponds to a collection of homomorphisms $\{ \phi_a : M_a \longrightarrow N_a \}_{a\in\mathscr{C}}$ such that

$$\phi_b(m \cdot f) = \phi_a(m) \cdot f \qquad \text{for all } m \in M_a \text{ and } f \in \mathscr{C}(a,b).$$

This defines a functor $\Theta\colon \text{Mod-}A \longrightarrow \widehat{\mathscr{C}^{\mathrm{op}}}$. Its quasi-inverse $\Omega$ maps any $F \in \widehat{\mathscr{C}^{\mathrm{op}}}$ to the module $\bigoplus_{a\in\mathscr{C}} Fa$, and any arrow $\{ \psi_a : Fa \longrightarrow Ga \}_{a\in\mathscr{C}}$ to the module homomorphism $\bigoplus_{a\in\mathscr{C}} \psi_a \colon \bigoplus_{a\in\mathscr{C}} Fa \longrightarrow \bigoplus_{a\in\mathscr{C}} Ga$. Thus, $\widehat{\mathscr{C}^{\mathrm{op}}}_{\mathrm{fin}}$ corresponds to the category mod-$A$ of finite-dimensional right $A$-modules.

Since $\widehat{\mathscr{C}^{\mathrm{op}}}_{\mathrm{fin}}$ is semisimple, so is $A$ and $A^{\mathrm{op}}$. Furthermore, $\Bbbk$ being perfect implies that the tensor product of any two finite-dimensional semisimple algebras is semisimple, see [FD93, Section 3]. In particular, we have that the algebra $B := A^{\mathrm{op}} \otimes_\Bbbk A \otimes_\Bbbk A^{\mathrm{op}} \otimes_\Bbbk A$ is semisimple, and

$$[\mathscr{C}^{\mathrm{op}} \otimes_\Bbbk \mathscr{C} \otimes_\Bbbk \mathscr{C}^{\mathrm{op}} \otimes_\Bbbk \mathscr{C}, \mathsf{vect}] \cong \text{mod-}B.$$

Write $\mathscr{D} := [\mathscr{C}^{\mathrm{op}} \otimes_\Bbbk \mathscr{C}, [\mathscr{C}^{\mathrm{op}} \otimes_\Bbbk \mathscr{C}, \mathsf{vect}]]$; then the functor

$$[\mathscr{C}^{\mathrm{op}} \otimes_\Bbbk \mathscr{C} \otimes_\Bbbk \mathscr{C}^{\mathrm{op}} \otimes_\Bbbk \mathscr{C}, \mathsf{vect}] \longrightarrow \mathscr{D}, \qquad F \longmapsto \widehat{F},$$

where $\widehat{F}(x,y)[u,v] = F(u,v,x,y)$, is a $\Bbbk$-linear equivalence of categories.

Further, since limits commute with limits, the functor

$$\mathrm{end}\colon \mathscr{D} \longrightarrow [\mathscr{C}^{\mathrm{op}} \otimes_\Bbbk \mathscr{C}, \mathsf{vect}], \qquad F \longmapsto (u,v \longmapsto \Big(\int_b F(b,b)\Big)(u,v))$$

is left exact. By semisimplicity, any short exact sequence in $\mathscr{D}$ splits and split epimorphisms are preserved by all functors. Thus, $\mathrm{end}\colon \mathscr{D} \longrightarrow [\mathscr{C}^{\mathrm{op}} \otimes_\Bbbk \mathscr{C}, \mathsf{vect}]$ is exact and must preserve colimits.

Given $F \in [\mathscr{C}^{\mathrm{op}} \otimes_\Bbbk \mathscr{C} \otimes_\Bbbk \mathscr{C}^{\mathrm{op}} \otimes_\Bbbk \mathscr{C}, \mathsf{vect}]$, we now compute

$$\begin{aligned} \int^a \int_b F(a,a)(b,b) &\cong \int^a \int_b \widehat{F}(b,b)[a,a] \cong \int_b \int^a \widehat{F}(b,b)(a,a) \\ &\cong \int_b \int^a F(a,a,b,b). \end{aligned}$$

□

### 3.2.1 *Cauchy completions*

Having assembled all of the necessary tools, we can now state explicit criteria for the Cauchy completion of (the opposite of) a $\Bbbk$-linear category to carry certain duality structures. As with modules over (commutative) rings, these notions are closely connected with objects being finitely-generated projective. Again, we implicitly assume all categories and functors to be $\Bbbk$-linear, for a field $\Bbbk$. However, we do not require finiteness of hom-spaces and work solely under the assumptions made in Hypothesis 3.24.

**Definition 3.36.** Suppose that $\mathscr{C}$ is a $\Bbbk$-linear category. The *Cauchy completion* of $\mathscr{C}$ is the full subcategory $\overline{\mathscr{C}}$ of $\widehat{\mathscr{C}}$ that consists of all presheaves $F$, such that $\widehat{\mathscr{C}}(F, -)$ commutes with all small limits.

**Remark 3.37.** Definition 3.36 differs from the usual definition of the Cauchy completion of a $\Bbbk$-linear category $\mathscr{C}$ as the additive and idempotent completion of $\mathscr{C}$. However, by [BD86, Proposition 2] and [LT22, Corollary 4.22], these two notions are equivalent. We choose the former because it is more convenient to work with for the purposes of this thesis.

**Example 3.38.** Let $\iota\colon \mathscr{C} \hookrightarrow \overline{\mathscr{C}}$ denote the inclusion of a category $\mathscr{C}$ into its Cauchy completion, sending $x \in \mathscr{C}$ to $\mathscr{C}(-, x)$. There is an adjoint equivalence

$$\left(\overline{\mathscr{C}}(-, \iota =),\ \overline{\mathscr{C}}(\iota -, =),\ \eta,\ \varepsilon\right),$$

between $\mathscr{C}$ and its Cauchy completion in the monoidal bicategory of profunctors, see Example 2.113. The unit

$$\eta_{x,y}\colon \mathscr{C}(x, y) \longrightarrow \int^{\overline{c} \in \overline{\mathscr{C}}} \overline{\mathscr{C}}(\iota_x, \overline{c}) \times \overline{\mathscr{C}}(\overline{c}, \iota_y)$$

is defined by

$$\mathscr{C}(x, y) \xrightarrow{\iota_{x,y}} \overline{\mathscr{C}}(\iota_x, \iota_y) \cong \int^{\overline{c} \in \overline{\mathscr{C}}} \overline{\mathscr{C}}(\iota_x, \overline{c}) \times \overline{\mathscr{C}}(\overline{c}, \iota_y),$$

and for the counit we have

$$\begin{aligned} \varepsilon_{\overline{x},\overline{y}}\colon \int^{c \in \mathscr{C}} \overline{\mathscr{C}}(\overline{x}, \iota_c) \times \overline{\mathscr{C}}(\iota_c, \overline{y}) &\longrightarrow \overline{\mathscr{C}}(\overline{x}, \overline{y}) \\ [(f, g)] &\longmapsto g \circ f. \end{aligned}$$

Since $\iota$ is fully faithful, it immediately follows that $\eta_{x,y}$ is an isomorphism, for all $x, y \in \mathscr{C}$, and the counit $\varepsilon$ is one by [AEBSJ01, Theorem 1.1].

In principle, the absence of free objects in the category of copresheaves necessitates a development of these notions in abstract categorical terms; this is discussed extensively for example in [Pre09]. For our purposes, the following characterisation in terms of the Cauchy completion of $\mathscr{C}^{\mathrm{op}}$ suffices.

**Lemma 3.39.** *Let $\mathscr{C}$ be a category. For any $F \in \widehat{\mathscr{C}^{\mathrm{op}}}$, the following are equivalent:*

(*i*) *$F$ is finitely-generated projective,*
(*ii*) *$F$ is a direct summand of a finite direct sum of representable functors, and*
(*iii*) *the functor $\widehat{\mathscr{C}^{\mathrm{op}}}(F,-)$ commutes with small colimits.*

*Proof.* The equivalence between (i) and (ii) is proven in Corollary 10.1.14 of [Pre09]. In order to show that (ii) and (iii) are equivalent, observe that the full subcategory $\overline{\mathscr{C}^{\mathrm{op}}}$ of $\widehat{\mathscr{C}^{\mathrm{op}}}$ consisting of direct summands of finite direct sums of representable functors is Cauchy complete by [LT22, Corollary 4.22]. The claim now follows by proceeding analogous to Proposition 2 of [BD86]. □

**Notation 3.40.** In view of Lemma 3.39, we shall adopt the following notation to emphasise that $X \in \overline{\mathscr{C}^{\mathrm{op}}}$ is a direct summand of a direct sum of representables:

$$X \underset{\pi_X}{\overset{\iota_X}{\rightleftarrows}} \bigoplus_{i=1}^{n} よ_{u_i},$$

where $\pi_X \circ \iota_X = \mathrm{id}_X$, and

$$よ : \mathscr{C}^{\mathrm{op}} \longrightarrow [(\mathscr{C}^{\mathrm{op}})^{\mathrm{op}}, \mathsf{Vect}] = [\mathscr{C}, \mathsf{Vect}], \qquad x \longmapsto \mathscr{C}^{\mathrm{op}}(-,x) = \mathscr{C}(x,-)$$

is the contravariant Yoneda embedding as in Equation (2.8.9).

**Proposition 3.41.** *For all objects $X \underset{\pi_X}{\overset{\iota_X}{\rightleftarrows}} \bigoplus_{i=1}^{n} よ_{u_i}$ and $Y \underset{\pi_Y}{\overset{\iota_Y}{\rightleftarrows}} \bigoplus_{j=1}^{m} よ_{v_j}$, the right internal hom of the Cauchy completion $\overline{\mathscr{C}^{\mathrm{op}}}$ of a right closed monoidal category $\mathscr{C}^{\mathrm{op}}$ exists and satisfies*

$$[X,Y]_r \underset{[\iota_X,\pi_Y]_r}{\overset{[\pi_X,\iota_Y]_r}{\rightleftarrows}} [\oplus_{i=1}^{n} よ_{u_i}, \oplus_{j=1}^{m} よ_{v_j}] \cong \bigoplus_{i=1}^{n} \bigoplus_{j=1}^{m} よ_{[u_i,v_j]_r}. \tag{3.2.5}$$

*Proof.* The category $\overline{\mathscr{C}^{\mathrm{op}}} \subseteq \widehat{\mathscr{C}^{\mathrm{op}}}$ is closed under taking tensor products since the Yoneda embedding is strong monoidal and $\overline{\mathscr{C}^{\mathrm{op}}}$ is the full subcategory of $\widehat{\mathscr{C}^{\mathrm{op}}}$ consisting of direct summands of finite direct sums of representables.

Let $X, Y \in \overline{\mathscr{C}^{\mathrm{op}}}$. Write $X = \operatorname{colim}_i \mathscr{C}(u_i, -)$ and $Y = \operatorname{colim}_j \mathscr{C}(v_j, -)$ as direct summands of finite direct sums of representables. We compute

$$
\begin{aligned}
[X,Y]_r &\overset{(2.8.6)}{=} \int_{a,b} \mathsf{Vect}(\mathscr{C}(a \otimes -, b), \mathsf{Vect}(Xa, Yb)) \overset{(2.8.1)}{\cong} \int_a \mathsf{Vect}(Xa, Y(a \otimes -)) \\
&\cong \int_a \mathsf{Vect}(\operatorname{colim}_i \mathscr{C}(u_i, a), \operatorname{colim}_j \mathscr{C}(v_j, a \otimes -)) \\
&\cong \lim_i \operatorname{colim}_j \int_a \mathsf{Vect}(\mathscr{C}(u_i, a), \mathscr{C}(v_j, a \otimes -)) \cong \lim_i \operatorname{colim}_j \mathscr{C}(v_j, u_i \otimes -) \\
&\cong \lim_i \operatorname{colim}_j \mathscr{C}([u_i, v_j]_r, -).
\end{aligned}
$$

It follows that $\overline{\mathscr{C}^{\mathrm{op}}}$ is right closed monoidal and that Equation (3.2.5) holds. □

**Example 3.42.** Proposition 3.41 can be understood as a variation of the fact that a direct summand or direct sum of (rigidly) dualisable objects is dualisable. Indeed, suppose $\mathscr{C}^{\mathrm{op}}$ to be a closed monoidal category and consider three direct summands $X \underset{\pi_X}{\overset{\iota_X}{\rightleftarrows}} U$, $Y \underset{\pi_Y}{\overset{\iota_Y}{\rightleftarrows}} V$, and $Z \underset{\pi_W}{\overset{\iota_W}{\rightleftarrows}} W$ of objects in $\overline{\mathscr{C}^{\mathrm{op}}}$. The following diagram, whose horizontal arrows are the isomorphisms of the tensor–hom adjunction of $\overline{\mathscr{C}^{\mathrm{op}}}$, commutes:

$$
\begin{array}{ccc}
\overline{\mathscr{C}^{\mathrm{op}}}(X * Y, Z) & \xrightarrow{\quad \phi_{X,Y,Z} \quad} & \overline{\mathscr{C}^{\mathrm{op}}}(Y, [X, Z]_r) \\
{\scriptstyle \overline{\mathscr{C}^{\mathrm{op}}}(\pi_X * \pi_Y, \iota_W)} \downarrow \; \uparrow {\scriptstyle \overline{\mathscr{C}^{\mathrm{op}}}(\iota_X * \iota_Y, \pi_Z)} & & {\scriptstyle \overline{\mathscr{C}^{\mathrm{op}}}(\pi_V, [\pi_X, \iota_W]_r)} \downarrow \; \uparrow {\scriptstyle \overline{\mathscr{C}^{\mathrm{op}}}(\iota_V, [\iota_X, \pi_Z]_r)} \\
\overline{\mathscr{C}^{\mathrm{op}}}(U * V, W) & \xrightarrow{\quad \phi_{U,V,W} \quad} & \overline{\mathscr{C}^{\mathrm{op}}}(V, [U, W]_r)
\end{array}
$$

Thus the unit and counit of the adjunction $X * - : \overline{\mathscr{C}} \rightleftarrows \overline{\mathscr{C}} : [X, -]_r$ satisfy

$$
\begin{aligned}
\eta_Y^{(X)} &= Y \xrightarrow{\iota_Y} V \xrightarrow{\eta_V^{(U)}} [U, U * V] \xrightarrow{[\iota_X, \pi_X * \pi_Y]} [X, X * V], \\
\varepsilon_Y^{(X)} &= X * [X, Y] \xrightarrow{\iota_X * [\pi_X, \iota_Y]} U * [U, V] \xrightarrow{\varepsilon_V^{(U)}} V \xrightarrow{\pi_Y} Y.
\end{aligned}
\tag{3.2.6}
$$

**Corollary 3.43.** *Let $\mathscr{C}^{\mathrm{op}}$ be a right closed monoidal category. We have:*

(i) *$\mathscr{C}^{\mathrm{op}}$ is right rigid if and only if $\overline{\mathscr{C}^{\mathrm{op}}}$ is.*
(ii) *$\mathscr{C}^{\mathrm{op}}$ is right tensor representable if and only if $\overline{\mathscr{C}^{\mathrm{op}}}$ is.*
(iii) *$(\mathscr{C}^{\mathrm{op}}, d)$ is a right Grothendieck–Verdier category if and only if $(\overline{\mathscr{C}^{\mathrm{op}}}, よ_d)$ is.*

*Proof.* That right rigidity and tensor representability of $\mathscr{C}^{\mathrm{op}}$ imply the same property for $\overline{\mathscr{C}^{\mathrm{op}}}$ follows from the description of the internal hom of $\overline{\mathscr{C}^{\mathrm{op}}}$ and the units and counits of the tensor-hom adjunctions, see Equation (3.2.6).

If $(\mathscr{C}^{\mathrm{op}}, d)$ is a right Grothendieck–Verdier category, then its right dualising functor is, up to natural isomorphism, given by $[-, d]_r \colon \mathscr{C} \longrightarrow \mathscr{C}^{\mathrm{op}}$. A direct computation using Proposition 3.41 shows that $[-, \text{よ}_d]_r \colon \overline{\mathscr{C}^{\mathrm{op}}}^{\mathrm{op}} \longrightarrow \overline{\mathscr{C}^{\mathrm{op}}}$ is an equivalence, and therefore $\left(\overline{\mathscr{C}^{\mathrm{op}}}, \text{よ}_d\right)$ is right Grothendieck–Verdier.

By Proposition 3.41 the right internal hom of two representable functors $\text{よ}_x$ and $\text{よ}_y$ is $[\text{よ}_x, \text{よ}_y]_r \cong \text{よ}_{[x,y]_r}$. Thus, the converse of any of the three statements is a consequence of the fact that, via the Yoneda embedding, $\mathscr{C}^{\mathrm{op}}$ is equivalent as a right closed monoidal category to the full subcategory of $\overline{\mathscr{C}^{\mathrm{op}}}$ whose objects are representable functors. □

**Remark 3.44.** The linearisation $\Bbbk\mathscr{W}$ of the tensor representable category $\mathscr{W}$ discussed in Theorem 3.23 is tensor representable, with left and right dualising functors $L$ and $R$, but not rigid. The latter follows from the fact that by the proof of Theorem 3.23 there exists an object $t \in \mathscr{W}$ whose canonical morphism $[t, 1]_r \otimes t \longrightarrow [t, t]_r$ is not invertible.

Since the right dualising functor $R$ of $\mathscr{W}$ is an anti-equivalence, and so $L \cong R^{-1}$, there also is a tensor representable, but not rigid, structure on $\Bbbk\mathscr{W}^{\mathrm{op}}$. One takes the left dualising functor to be $R$, and the right one to be $L$. By Corollary 3.43, $\overline{\Bbbk\mathscr{W}^{\mathrm{op}}}$ is tensor representable but not rigid.

Notice that this is in stark contrast to many cases arising in representation theory, like modules over commutative rings or $\Bbbk$-algebras, where rigidity and tensor representability are equivalent; see for example [NW17, Proposition 2.1] for a slightly more general statement. Since a $\Bbbk$-algebra can be interpreted as a $\Bbbk$-linear category with one object, the following proposition may be seen as a "many object" version of the classical case.

**Proposition 3.45.** *Let $\mathscr{C}$ be a left rigid monoidal category. the following are equivalent for some $X \in \widehat{\mathscr{C}^{\mathrm{op}}}$:*

(i) *it has a right rigid dual,*
(ii) *there exists a $\mathbb{D}X \in \widehat{\mathscr{C}^{\mathrm{op}}}$ such that $X * - \colon \widehat{\mathscr{C}^{\mathrm{op}}} \rightleftarrows \widehat{\mathscr{C}^{\mathrm{op}}} \colon \mathbb{D}X * -$, and*
(iii) *it is finitely-generated projective.*

*Proof.* By definition, (i) $\Longrightarrow$ (ii). To show that (ii) $\Longrightarrow$ (iii), we assume that there exists a $\mathbb{D}X \in \widehat{\mathscr{C}^{\mathrm{op}}}$ such that $X * - \colon \widehat{\mathscr{C}^{\mathrm{op}}} \rightleftarrows \widehat{\mathscr{C}^{\mathrm{op}}} \colon \mathbb{D}X * -$, and fix a small colimit $\operatorname{colim}_{i\in I} Fi \in \widehat{\mathscr{C}^{\mathrm{op}}}$ of some diagram $F\colon I \longrightarrow \widehat{\mathscr{C}^{\mathrm{op}}}$. For every $i \in I$, write $\iota_i\colon Fi \longrightarrow \operatorname{colim}_{i\in I} Fi \in \widehat{\mathscr{C}^{\mathrm{op}}}$ for its structure morphisms. Now consider the commutative diagram:

$$\begin{array}{ccc}
\operatorname{colim}_{i\in I} \widehat{\mathscr{C}^{\mathrm{op}}}(X, Fi) & \xrightarrow{\operatorname{colim}_{i\in I}(\widehat{\mathscr{C}^{\mathrm{op}}}(X,\iota_i))} & \widehat{\mathscr{C}^{\mathrm{op}}}(X, \operatorname{colim}_{i\in I} Fi) \\
\Big\downarrow{\scriptstyle \operatorname{colim}_{i\in I}\phi_{X,Fi}} & & \Big\downarrow{\scriptstyle \phi_{X,\operatorname{colim}_{i\in I} Fi}} \\
\operatorname{colim}_{i\in I} \widehat{\mathscr{C}^{\mathrm{op}}}(1, \mathbb{D}X * Fi) & \xrightarrow{\operatorname{colim}_{i\in I}(\widehat{\mathscr{C}^{\mathrm{op}}}(1,\mathbb{D}X*\iota_i))} & \widehat{\mathscr{C}^{\mathrm{op}}}(1, \mathbb{D}X * \operatorname{colim}_{i\in I} Fi)
\end{array}$$

Its horizontal arrows are induced by the universal property of the colimit and the vertical arrows are due to the tensor–hom adjunction of $\widehat{\mathscr{C}^{\mathrm{op}}}$. The functors $\widehat{\mathscr{C}^{\mathrm{op}}}(1,-)$ and $\mathbb{D}X \otimes -$ commute with all small colimits; the first one due to the fact that $\widehat{\mathscr{C}^{\mathrm{op}}}(1,-)$ is finitely-generated projective, see Lemma 3.39, and the second one due to Day convolution being a colimit. Therefore, the horizontal arrow at the bottom is invertible. Since the vertical arrows are also invertible, the canonical arrow $\operatorname{colim}_{i\in I} \widehat{\mathscr{C}^{\mathrm{op}}}(X, Fi) \longrightarrow \widehat{\mathscr{C}^{\mathrm{op}}}(X, \operatorname{colim}_{i\in I} Fi)$ displayed at the top of the diagram must be an isomorphism. Again, using Lemma 3.39, $X$ must be finitely-generated projective. Thus, (ii) implies (iii).

Finally, if $X \in \widehat{\mathscr{C}^{\mathrm{op}}}$ is finitely-generated projective, then it is contained in the Cauchy-completion of $\mathscr{C}^{\mathrm{op}}$, which is right rigid. Then, by Corollary 3.43, so is $\overline{\mathscr{C}^{\mathrm{op}}}$. Therefore $X$ admits a rigid dual and (iii) $\Longrightarrow$ (i). □

## 3.3 Applications

We conclude the chapter by discussing several examples.

We investigate Boolean algebras, their applications in group and ring theory, and how they induce abelian $\Bbbk$-linear Grothendieck–Verdier categories.

Next, we focus on Mackey functors. These are, roughly speaking, collections of vector spaces indexed by all subgroups of a fixed finite group, together with morphisms subject to relations resembling the behaviour of induction, restriction, and conjugation operations, including the eponymous Mackey identity; see [Lin76; TW95]. Finite-dimensional Mackey functors form a Grothendieck–Verdier category [PS07]; we show in Proposition 3.57 that it is rigid if and only if it is semisimple.

The last example arises in the study of crossed modules, which—in categorical terms—correspond to strict 2-groups. The functors from any finite strict 2-group to vect form an abelian monoidal category that is equivalent to a direct sum of representation categories of the isotropy-group of the monoidal unit of the 2-group. In Proposition 3.66 we prove the rigidity of this category to be equivalent to the semisimplicity of a certain group algebra.

### 3.3.1 *Boolean algebras*

**Definition 3.46.** A *lattice* consists of a set $L$ together with two associative and commutative operations

$$\wedge\colon L\times L\longrightarrow L,\quad (a,b)\longmapsto a\wedge b\quad \text{and}\quad \vee\colon L\times L\longrightarrow L,\quad (a,b)\longmapsto a\vee b,$$

called *meet* and *join*, which satisfy the *absorption laws*, for all $a,b\in L$:

$$a\vee(a\wedge b)=a,\qquad a\wedge(a\vee b)=a.$$

**Remark 3.47.** Any lattice $(L,\wedge,\vee)$ defines a poset via the relation

$$a\leq b\iff b=a\vee b\iff a=a\wedge b,\qquad a,b\in L.$$

Conversely, a poset $(P,\leq)$ that admits for any pair of objects $a,b\in P$ a least upper bound $a\vee b\in P$ and a greatest lower bound $a\wedge b\in P$ is a lattice.

A direct computation shows that an element $1\in L$ is maximal with respect to the partial order of the lattice $L$ if and only if $1\vee a=1$ for all $a\in L$. Analogously, $0\in L$ being minimal equates to $0\wedge a=0$ for any $a\in L$. Minimal and maximal elements are unique. In case they exist, we call $L$ *bounded*.

**Definition 3.48.** A *Boolean algebra* is a bounded lattice $(L,\wedge,\vee)$ satisfying the following *distributivity* condition for all $a,b,c\in L$:

$$a\wedge(b\vee c)=(a\wedge b)\vee(a\wedge c),\qquad a\vee(b\wedge c)=(a\wedge b)\vee(a\wedge c),\tag{3.3.1}$$

and every $a\in L$ admits a *complement* $a^{\perp}\in L$ in the sense that

$$a\vee a^{\perp}=1\qquad\text{and}\qquad a\wedge a^{\perp}=0.$$

**Remark 3.49.** For any Boolean algebra $(L,\wedge,\vee)$, any of the two distributivity requirements of Equation (3.3.1) implies the other. Further, a direct computation shows that complements are unique. We obtain an involutive map $(-)^{\perp}\colon L\longrightarrow L$ that maps $a$ to $a^{\perp}$. The maximal element $1\in L$ is a unit for $\wedge$:

$$a\wedge 1=a\wedge(a\vee a^{\perp})=a,\qquad\text{for all } a\in L.$$

An analogous computation shows that the minimal element $0$ in a Boolean algebra satisfies $a\vee 0=a$.

**Example 3.50.**

- *Central idempotents*: the set $C := \{ e \in Z(R) \mid e^2 = e \}$ of central idempotents of a ring $R$ is a Boolean algebra when endowed with

$$\begin{aligned} \wedge\colon C \times C &\longrightarrow C & e \wedge f &= ef, \\ \vee\colon C \times C &\longrightarrow C & e \vee f &= e + f - ef, \\ (-)^{\perp}\colon C &\longrightarrow C & e^{\perp} &= 1 - e. \end{aligned}$$

- *Annihilators in semiprime rings*: consider a commutative semiprime ring $R$ with 1. That is, its Jacobson radical is trivial. In case $R$ is furthermore Artinean, this is equivalent to it being semisimple. The *annihilator* of an ideal $I \subset R$ is $I^{\perp} := \{ x \in R \mid xI = 0 \}$; it is a radical ideal. We define on the set $\mathrm{Ann}(R)$ of annihilators the maps

$$\begin{aligned} \wedge\colon \mathrm{Ann}(R) \times \mathrm{Ann}(R) &\longrightarrow \mathrm{Ann}(R), & I \wedge J &= I \cap J, \\ \vee\colon \mathrm{Ann}(R) \times \mathrm{Ann}(R) &\longrightarrow \mathrm{Ann}(R), & I \vee J &= (I + J)^{\perp}. \end{aligned}$$

  This defines the structure of a Boolean algebra on $\mathrm{Ann}(R)$. In fact, the (right) annihilators of a not necessarily commutative ring form a Boolean algebra if and only if the ring is semiprime, see [DT21].

- *The subgroup lattice*: let $H \subseteq G$ be a subgroup of a finite group. We write $H^{\perp}$ for the intersection of all maximal subgroups that do not contain $H$. The minimal group constructed in this manner is the Frattini subgroup $\Phi(G) = G^{\perp}$ of $G$. Deaconescu, Isaacs, and Wall showed that the set

$$\{ \Phi(G) \subseteq H \subseteq G \mid G = HH^{\perp} \},$$

  partially ordered under inclusion forms a Boolean algebra, see [DIW11].

Barr showed in [Bar79] that the duality of Boolean algebras fits within the framework of Grothendieck–Verdier categories.

**Proposition 3.51.** *Let $(L, \wedge, \vee)$ be a Boolean algebra. The associated poset-category $\mathfrak{L}$ is a left Grothendieck–Verdier category with meet as tensor product, 1 as monoidal unit, and 0 as dualising object. The dualising functor is induced by $(-)^{\perp}$.*

*If $L$ is finite, the category $\langle \mathfrak{L}, \mathsf{vect} \rangle$ of ordinary functors from $\mathfrak{L}$ to the category of finite-dimensional $\Bbbk$-vector spaces can be equipped with a right Grothendieck–Verdier structure, with $\Bbbk\mathfrak{L}(-, s)^*$ as dualising object and dualising functor*

$$\mathfrak{D}_r\colon \langle \mathfrak{L}, \mathsf{vect} \rangle^{\mathrm{op}} \longrightarrow \langle \mathfrak{L}, \mathsf{vect} \rangle, \qquad F \longmapsto F(-^{\perp})^*.$$

*Proof.* A direct consequence of Equation (3.3.1) is that for $a \leq c$ and $b \leq d$, we have $a \wedge b \leq c \wedge d$. Thus, since $\wedge \colon L \times L \longrightarrow L$ is associative, commutative, and unital, it induces a (symmetric) monoidal structure on $\mathscr{L}$.

To see that $(\mathscr{L}, 0)$ is a Grothendieck–Verdier category, we compute its left internal hom. Let us fix a triple of elements $a, b, c \in L$. We have

$$\begin{aligned} b \wedge a \leq c \iff (b \wedge a) \vee a^{\perp} \leq c \vee a^{\perp} \iff (b \vee a^{\perp}) \wedge (a \vee a^{\perp}) \leq c \vee a^{\perp} \\ \iff b \vee a^{\perp} \leq c \vee a^{\perp} \iff b \leq c \vee a^{\perp}, \end{aligned}$$

where the last equivalence is due to the fact that $b = b \vee 0 \leq b \vee a^{\perp}$ and $b \leq c \vee a^{\perp}$, implying that $b \vee a^{\perp} \leq (c \vee a^{\perp}) \vee a^{\perp} = c \vee a^{\perp}$. It follows that

$$[-,=]_{\ell} \colon L \times L \longrightarrow L, \qquad [a,b]_{\ell} \coloneqq (a \wedge b^{\perp})^{\perp} = b \vee a^{\perp},$$

defines the left internal hom of $\mathscr{L}$, and the order reversing involution

$$[-,0]_{\ell} = (-)^{\perp} \colon L \longrightarrow L$$

induces an equivalence of categories $D_{\ell} \colon \mathscr{L}^{\mathrm{op}} \longrightarrow \mathscr{L}$. As discussed in Remark 3.13, $(\mathscr{L}, 0)$ is a left Grothendieck–Verdier category.

By Example 2.82, we have an equivalence $\langle \mathscr{L}, \mathsf{vect} \rangle \cong [\Bbbk\mathscr{L}, \mathsf{vect}]$. Under the assumption that $L$ is finite, Proposition 3.27 shows that $\langle \mathscr{L}, \mathsf{vect} \rangle$ has the structure of a right Grothendieck–Verdier category with the specified dualising object and dualising functor. □

Suppose $L$ is a finite Boolean algebra and $\mathscr{L}$ is its poset-category. As discussed in the first steps of the proof of Lemma 3.35, the category $[\Bbbk\mathscr{L}, \mathsf{vect}]$ can be identified with the finite-dimensional right modules of the path algebra $A \cong \bigoplus_{e,f \in L} \Bbbk\mathscr{L}(e,f)$ of the poset associated to $L$.

**Example 3.52.** The partial order on the set $L \coloneqq \{0, a, b, 1\}$ displayed in the following diagram defines a Boolean algebra:

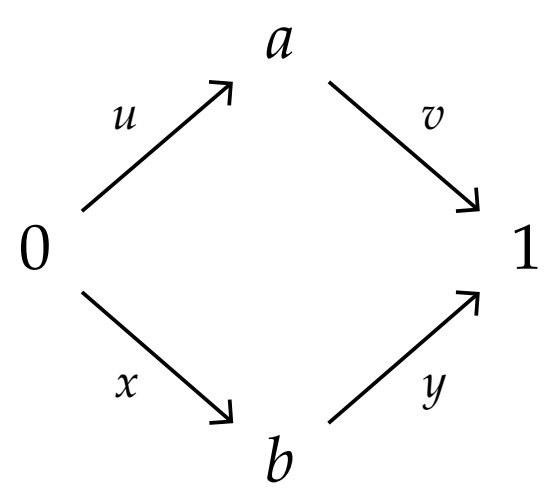

A direct computation shows that its path algebra $A$ is isomorphic to a subalgebra of the $\Bbbk$-valued $4 \times 4$ upper triangular matrices:

$$A \cong \left\{ \begin{pmatrix} 0 & u & x & z \\ 0 & a & 0 & v \\ 0 & 0 & b & y \\ 0 & 0 & 0 & 1 \end{pmatrix} \middle| \; 0, u, x, z, a, v, b, y, 1 \in \Bbbk \right\}.$$

Subalgebras of upper triangular matrices were studied by Thrall in [Thr48] to provide generalisations of quasi-Frobenius algebras.[9] The following definition generalises this property.

[9] A $\Bbbk$-algebra is *quasi-Frobenius* if all projective right modules are injective.

**Definition 3.53.** A $\Bbbk$-algebra is called QF-2 if all indecomposable projective right and all indecomposable projective left modules have simple socles.

Using the Grothendieck–Verdier structure of Proposition 3.51, we will show that the path algebra of any finite Boolean algebra $L$ is QF-2. Hereto we need the following observation. Taking complements in $L$ extends to an anti-algebra isomorphism $\phi\colon A \longrightarrow A$; its pushforward $(-)_\phi\colon \text{mod-}A \longrightarrow A\text{-mod}$ is an involutive equivalence of categories. It maps any left module $M$ to the right module $M_\phi$ that has the same underlying vector space and is endowed with the action $a \triangleright m := m \triangleleft \phi(a)$ for all $m \in M$ and $a \in A$.

**Proposition 3.54.** *Let $A$ be the path algebra of a finite Boolean algebra $L$. Then*

$$(-_\phi)^*\colon (\text{mod-}A)^{\text{op}} \longrightarrow \text{mod-}A$$

*is an equivalence of categories. In particular, $M \in \text{mod-}A$ is projective if and only if ${M_\phi}^*$ is injective. Further, $A$ is* QF-2, *and quasi-Frobenius if and only if $|L| = 1$.*

*Proof.* The first statement follows directly from the equivalence of categories between $[\Bbbk\mathcal{L}, \mathsf{vect}]$ and mod-$A$ given in the proof of Lemma 3.35, and the definition of the dualising functor of $[\Bbbk\mathcal{L}, \mathsf{vect}]$ stated in Proposition 3.51.

In order to show the second claim, we observe that the finite-dimensional algebra $A = \bigoplus_{n \geq 0} A_n$ is graded by the path lengths. The elements of $L$ form a basis of $A_0$ and correspond to the *primitive idempotents*—non-zero idempotents, which cannot be written as a sum of two non-zero orthogonal idempotents. Furthermore, the Jacobson radical of $A$ is $J(A) = \bigoplus_{n \geq 1} A_n$. The indecomposable projectives of $A$ are of the form $eA$, for $e \in L$. Note that $eA$ has a vector space basis by paths $[e, y]$ for $e \leq y \leq 1$. The socle of $eA$ is

$$\text{soc}(eA) = \{\, m \in eA \mid mJ(A) = 0 \,\} = \text{span}_\Bbbk\{[e, 1]\},$$

which is one-dimensional and therefore simple. Using that $A$-mod $\cong$ mod-$A$, it follows that $A$ is a QF-2 algebra.

For $|L| = 1$, we obtain the trivial (quasi-)Frobenius-algebra $A = \Bbbk$. Otherwise, there exists an $e \in L$ such that $\dim eA \geq 2$. A direct computation shows that $\mathrm{soc}(eA) \cong 1A$. If $eA$ was injective, the inclusion $1A \cong \mathrm{soc}(eA) \hookrightarrow eA$ would have a retraction, contradicting the indecomposability of $eA$. □

### 3.3.2 *Mackey functors*

Let $G$ be a finite group and write $\mathsf{Sp}_G$ for the category of isomorphism classes of spans of finite $G$-sets. One can show that each hom-set is a free and finitely-generated commutative monoid, see for example the discussion preceding Lemma 2.1 of [TW95]. Besides the definition of Mackey functors sketched at the beginning of Section 3.3, there is a more succinct formulation due to Lindner [Lin76].

In this section we again implicitly assume all categories and functors to be $\Bbbk$-linear.

**Definition 3.55.** Let $\Bbbk$ be a field and $G$ a finite group. The category $\mathsf{Mky}_{\Bbbk}(G)$ of *Mackey functors* of $G$ is given by the $\Bbbk$-linear functor category $[\Bbbk\mathsf{Sp}_G, \mathsf{Vect}]$.

Examples of Mackey functors are plentiful: any representation of $G$ defines one, see [Thé95, Example 53.1] as well as [PS07, Propostion 10.1]. For an extensive overview, we refer the reader to [Thé95, Chapter 53].

**Example 3.56.** To any finite group $G$, one can associate a finite-dimensional algebra $\mathbb{M}G$—the *Mackey algebra of* $G$—whose category of modules is equivalent to $\mathsf{Mky}_{\Bbbk}(G)$, see [TW95, Propositions 3.1 and 3.2]. Its finite-dimensional modules are in correspondence with the (pointwise) finite-dimensional Mackey functors, which we denote by $\mathsf{mky}_{\Bbbk}(G)$.

Given that $\Bbbk\mathsf{Sp}_G$ is symmetric monoidal, $\mathsf{Mky}_{\Bbbk}(G)$ is closed symmetric monoidal when equipped with Day convolution as its tensor product.[10] Furthermore, $\Bbbk\mathsf{Sp}_G$ has a finite dense subcategory whose objects form a complete set of representatives of transitive $G$-sets. The arguments of Example 3.26 may now be used to show that the tensor product and internal hom of finite-dimensional Mackey functors is finite-dimensional, see [PS07, Section 9].

The following results gives an alternative proof for the sketched argument in [Bou05, Lemma 2.2] that the classes of rigid and finitely-generated projective Mackey functors coincide.

[10] As such, we refrain from adding "left" or "right" prefixes to the different duality notions when talking about Mackey functors.

**Proposition 3.57.** *Let $G$ be a finite group and suppose $\Bbbk$ is a field. The category* $\mathsf{mky}_{\Bbbk}(G)$ *is a Grothendieck–Verdier category with* $\mathsf{mky}_{\Bbbk}(-,1)^*$ *as dualising object. Its rigid objects are precisely the finitely-generated projective Mackey functors.*

*Furthermore,* $\mathsf{mky}_{\Bbbk}(G)$ *is rigid itself if and only if* $\mathbb{M}G$ *is semisimple, which is equivalent to* $\operatorname{char}\Bbbk$ *not dividing the order of $G$.*

*Proof.* The fact that $\mathsf{mky}_{\Bbbk}(G)$ is a Grothendieck–Verdier category was proven in [PS07, Theorem 9.2]. Alternatively, we may use that $\Bbbk\mathsf{Sp}_G$ is a rigid category, see [PS07, Section 2], and apply Proposition 3.27 to obtain that $\mathsf{mky}_{\Bbbk}(G)$ is a Grothendieck–Verdier category with $\mathsf{mky}_{\Bbbk}(-,1)^* \in \mathsf{mky}_{\Bbbk}(G)$ as its dualising object. Due to Proposition 3.45, a Mackey functor is rigid dualisable if and only if it is finitely-generated projective.

Now let us assume $\mathsf{mky}_{\Bbbk}(G)$ is rigid and recall that it is equivalent to the category $\mathbb{M}G$-mod of finite-dimensional modules of the Mackey algebra. Since $\mathbb{M}G$ is finite-dimensional and every object in $\mathsf{mky}_{\Bbbk}G$ is projective by Proposition 3.45, any submodule of $\mathbb{M}G$ must have a complement. In particular, $\mathbb{M}G$ is semisimple. Conversely, in case $\mathbb{M}G$ is semisimple, all objects of $\mathbb{M}G$-mod, and therefore also $\mathsf{mky}_{\Bbbk}(G)$, are projective. As they are furthermore finitely-generated, Proposition 3.45 implies the rigidity of $\mathsf{mky}_{\Bbbk}(G)$.

Corollary 14.4 of [TW95], shows that $\mathbb{M}G$ is semisimple if $\operatorname{char}\Bbbk$ does not divide $|G|$. On the other hand, the category $\operatorname{Rep}_{\Bbbk}(G)$ of finite-dimensional representations of $G$ over $\Bbbk$ is a full subcategory of $\mathsf{mky}_{\Bbbk}(G)$, see [PS07, Proposition 10.1]. Thus, if $\mathsf{mky}_{\Bbbk}(G)$ is semisimple, so is $\operatorname{Rep}_{\Bbbk}(G)$, implying that $\operatorname{char}\Bbbk$ does not divide $|G|$ by Maschke's theorem. □

### 3.3.3 *Crossed modules*

Another example for Grothendieck–Verdier structures on abelian $\Bbbk$-linear functor categories arises from studying (strict) 2-groups, which can be identified with crossed modules. Our exposition follows [Wag21].

**Definition 3.58.** A *strict 2-group* is a (small) groupoid $\mathscr{G}$ endowed with a strict monoidal structure, such that the monoid of objects $(\operatorname{Ob}(\mathscr{G}), \otimes, 1)$ is a group.

Every strict 2-group is a rigid category; the left and right dual of an object $g \in \mathscr{G}$ is given by its inverse $g^{-1} \in \mathscr{G}$.

Let $G$ be the group of objects of $\mathscr{G}$ and $H = \mathscr{G}(1,-)$ the set of arrows which start at the monoidal unit $1 \in G$. Writing $t\colon H \longrightarrow G$ for the *target map*, we define a group structure on $H$ via the multiplication

$$h'h := (h' \otimes \mathrm{id}_{t(h)}) \circ h, \qquad \text{for all } h, h' \in H.$$

For any $g \in G$ and $h \in H$ there is a unique $h' \in H$ such that $\mathrm{id}_g \otimes h = h' \otimes \mathrm{id}_g$. This induces a map $\alpha\colon G \longrightarrow \mathrm{Aut}(H)$ and turns the quadruple $(G, H, t, \alpha)$ into a crossed module as defined below.

**Definition 3.59.** A *crossed module* is a quadruple $(G, H, t, \alpha)$ consisting of two groups $G, H$ and two group homomorphisms $t\colon H \longrightarrow G$, $\alpha\colon G \longrightarrow \mathrm{Aut}(H)$, such that for all $g \in G$ and $h, l \in H$, we have

$$t(\alpha(g)h) = g \circ t(h) \circ g^{-1}, \qquad \alpha(t(l))h = l \circ h \circ l^{-1}. \tag{3.3.2}$$

**Remark 3.60.** Any crossed module $(G, H, t, \alpha)$ defines a strict 2-group $\mathcal{G}$, with $\mathrm{Ob}\,\mathcal{G} = G$, and $\mathcal{G}(g, g') = \{ (h, g) \in H \times G \mid t(h)g = g' \}$. Composition is given by $(h', t(h)g) \circ (h, g) := (h'h, g)$, and for the tensor product we have $(h, g) \otimes (h', g') := (h\,\alpha(g)h', gg')$. The left and right dualising functors coincide; they map any object $b \in G$ to ${}^{\vee}b = b^{-1} = b^{\vee}$. For a morphism $f = (h, b)\colon b \longrightarrow c$, we have

$${}^{\vee}f = (\alpha(c^{-1})h, c^{-1}) = (\alpha(b^{-1})h, c^{-1}) = f^{\vee}.$$

There is an equivalence of categories between crossed modules and strict 2-groups as is shown for example in [Wag21, Theorem 1.9.2].

**Example 3.61.**

- *Hopf–Galois extensions and skew braces*: the *holomorph* of a group $H$ is the semidirect product $H \rtimes \mathrm{Aut}(H)$. The two group homomorphisms $t\colon H \longrightarrow \mathrm{Aut}(H)$, $t(l)h = lhl^{-1}$, and $\alpha = \mathrm{id}\colon \mathrm{Aut}(H) \longrightarrow \mathrm{Aut}(H)$ turn $(\mathrm{Aut}(H), H, t, \mathrm{id}_{\mathrm{Aut}(H)})$ into a crossed module. As explained in the Introduction of [Byo24] and [Bac16, Section 2], holomorphs can be used to classify certain Hopf–Galois extensions as well as skew-braces.

- *Representations of finite groups*: a short exact sequence of groups

  $$0 \longrightarrow C \xrightarrow{\ \iota\ } E \xrightarrow{\ t\ } G \longrightarrow 0$$

  such that $C$ embeds into the centre of $E$ is called a *central extension* of $G$. Let $\kappa\colon G \longrightarrow E$ be a set-theoretical section of $t\colon E \longrightarrow G$. A direct computation shows that the map $\alpha\colon G \longrightarrow \mathrm{Aut}(E)$ defined by $\alpha(g)e := \kappa(g)e\kappa(g)^{-1}$ is independent of the choice of the section and a homomorphism of groups. The tuple $(E, G, t, \alpha)$ forms a crossed module, see [Wag21, Section 1.3].

As discussed in [HH92, Chapter 1], if $G$ is a finite group, then any projective representation $\rho\colon G \longrightarrow \mathrm{PGL}_{\mathbb{C}}(V)$ of $G$ can be lifted to a linear representation $\varrho\colon E \longrightarrow \mathrm{GL}(V)$ of a certain central extension $E$ of $G$. Projective representations themselves are studied in the context of representation theory of semidirect products, see for example [CSST22].

Let $(G, H, t, \alpha)$ be a crossed module whose associated 2-group we denote by $\mathscr{G}$. Using the results of Section 3.2, we determine a Grothendieck–Verdier structure on the category $\langle \mathscr{G}, \mathsf{vect}\rangle$ of ordinary functors between $\mathscr{G}$ and $\mathsf{vect}$.

Hereto, we need to analyse the structure of $\mathscr{G}$ in more detail. By Equation (3.3.2), the image $K := \operatorname{im} t \subset G$ is a normal subgroup of $G$ and the kernel $L := \ker t$ is normal and central in $H$. The latter follows from the identity

$$lhl^{-1} = \alpha(t(l))h = \alpha(e)h = h \qquad \text{for all } l \in L, h \in H.$$

For any $k = t(h) \in K$ and $l \in L$, we get $\alpha(k)l = \alpha(t(h))l = hlh^{-1} = l$. Hence, there is a unique group homomorphism $\overline{\alpha}\colon Q := G/N \longrightarrow \mathrm{Aut}(L)$, with

$$\overline{\alpha}([g])l = \alpha(g)l \in L \subset H.$$

for all $g \in G$ and $l \in L$.

The connected components of $\mathscr{G}$ are in bijection with the elements of $Q$. Given any element $g \in G$, we write $\mathscr{G}_g$ for the maximal full connected pointed subgroupoid of $\mathscr{G}$ whose distinguished object is $g$. The set of objects of $\mathscr{G}_g$ is $Kg$ and its morphisms correspond to $H \times Kg$.

**Lemma 3.62.** *Let $g \in G$ and write $\mathbf{B}L$ for the delooping of $L$ with single object $g$. The canonical inclusion $\Theta_g\colon \mathbf{B}L \hookrightarrow \mathscr{G}_g$ sending $g$ to $g$ and $l\colon g \longrightarrow g$ to $(l, g)\colon g \longrightarrow g$ is essentially surjective and fully faithful—an equivalence.*

To specify a quasi-inverse, fix a set-theoretical section $\iota\colon K \longrightarrow H$ of the surjective homomorphisms of groups $t\colon H \longrightarrow K$. Define the quasi-inverse $\Psi_g\colon \mathscr{G}_g \longrightarrow \mathbf{B}L$ by mapping each object to $g$ and any morphism $(h, kg)$ to $\iota((t(h)k)^{-1})h\iota(k)$. Using that $\langle \mathbf{B}L, \mathsf{vect}\rangle$ can be identified with the category $\mathsf{rep}(L)$ of finite-dimensional representations of $L$, the pushforwards of $\Theta_g$ and $\Psi_g$ establish an equivalence of categories

$$\Theta_g^*\colon \langle \mathscr{G}_g, \mathsf{vect}\rangle \rightleftarrows \mathsf{rep}(L) : \Psi_g^*$$

$$F \longmapsto (Fg, \rho_F) \quad (\text{where } \rho_F(l) = F(l, g))$$

(3.3.3) $$\left\{ \begin{aligned} kg &\longmapsto M \\ (hk, g) &\longmapsto \rho\left(\iota((t(h)k)^{-1})h\iota(k)\right) \end{aligned} \right\} \longleftarrow (M, \rho).$$

**Notation 3.63.** Define $\mathsf{rep}_{\mathcal{G}}(L) := \bigoplus_{q\in Q} \mathsf{rep}(L)_q$, where $\mathsf{rep}(L)_q = \mathsf{rep}(L)$. Given some $M \in \mathsf{rep}_{\mathcal{G}}(L)$, we write $M_q$ for its the homogeneous component of degree $q \in Q$. Note that in case $K$ has finite index in $G$, we have

$$\langle \mathcal{G}, \mathsf{vect}\rangle \cong \bigoplus_{[g]\in Q} \langle \mathcal{G}_g, \mathsf{vect}\rangle.$$

**Lemma 3.64.** *Let $\mathcal{G}$ be a strict 2-group, and write $(G, H, t, \alpha)$ for the corresponding crossed module. If $K = \operatorname{im} t \subset G$ has finite-index, then $\mathsf{rep}_{\mathcal{G}}(L) \simeq \langle \mathcal{G}, \mathsf{vect}\rangle$.*

The tools developed in Section 3.2 allow us to endow $\langle \mathcal{G}, \mathsf{vect}\rangle$ with the structure of a Grothendieck–Verdier category which we will transfer to $\mathsf{rep}_{\mathcal{G}}(L)$. The tensor product obtained in this manner will permute the homogeneous components and "twist" the action of $L$ by virtue of $\overline{\alpha}$. In this regard, we introduce the following notation: for any $q \in Q$ and $M \in \mathsf{rep}(L)$, we denote by $M^{\overline{\alpha}(q)}$ the representation of $L$ whose underlying vector space is $M$, endowed with the action $l \blacktriangleright m = \overline{\alpha}(q)l \triangleright m$ for all $l \in L$ and $m \in M$.

**Proposition 3.65.** *Suppose $(G, H, t, \alpha)$ is a crossed module with $G$ and $H$ finite and let $\mathcal{G}$ be its associated strict 2-group. The category $\mathsf{rep}_{\mathcal{G}}(L)$ is a right r-category. Its tensor product and dualising functor are defined by the assignments*

$$M \otimes N := (M \otimes_{\Bbbk L} N^{\overline{\alpha}(p^{-1})}) \in \mathsf{rep}_{\mathcal{G}}(L)_{pq}, \qquad DM := \left(M^{\overline{\alpha}(p)}\right)^* \in \mathsf{rep}_{\mathcal{G}}(L)_{p^{-1}}$$

*for $p, q \in Q$ and $M \in \mathsf{rep}_{\mathcal{G}}(L)_p$, $N \in \mathsf{rep}_{\mathcal{G}}(L)_q$.*

*Proof.* By Proposition 3.27, $\langle \mathcal{G}, \mathsf{vect}\rangle \cong [\Bbbk\mathcal{G}, \mathsf{vect}]$ can be endowed with the structure of a right Grothendieck–Verdier category with $\Bbbk\mathcal{G}(-, 1)^*$ as dualising object; write $R\colon \langle \mathcal{G}, \mathsf{vect}\rangle^{\mathrm{op}} \longrightarrow \langle \mathcal{G}, \mathsf{vect}\rangle$ for its dualising functor.

By Lemma 3.64, there are $\Bbbk$-linear equivalences

$$\widehat{\Theta}\colon [\Bbbk\mathcal{G}, \mathsf{vect}] \rightleftarrows \mathsf{rep}_{\mathcal{G}}(L) : \widehat{\Psi}$$

that are determined by the pushforwards of the equivalences

$$\Theta_g\colon \langle \mathcal{G}_g, \mathsf{vect}\rangle \rightleftarrows \mathbf{B}L : \Psi_g,$$

for each homogeneous component $[\Bbbk\mathcal{G}_{[g]}, \mathsf{vect}]$. Using the formulas of Equation (3.3.3), we can therefore explicitly transfer the Grothendieck–Verdier structure of $[\Bbbk\mathcal{G}, \mathsf{vect}]$ to $\mathsf{rep}_{\mathcal{G}}(L)$.

The dualising object is mapped to $\Theta(\Bbbk\mathcal{G}(-,1)^*) = \Bbbk L^* \cong \Bbbk L \in \mathsf{rep}_{\mathcal{G}}(L)_{[1]}$, where we used that $\Bbbk L$ is a Frobenius algebra for the last equality. The tensor product and right dualising functor of $\mathsf{rep}_{\mathcal{G}}(L)$ are given by

$$\otimes := \widehat{\Theta} * (\widehat{\Psi} \times \widehat{\Psi})\colon \mathsf{rep}_{\mathcal{G}}(L) \times \mathsf{rep}_{\mathcal{G}}(L) \longrightarrow \mathsf{rep}_{\mathcal{G}}(L),$$
$$D := \widehat{\Theta}\mathsf{R}\widehat{\Psi}^{\mathrm{op}}\colon \mathsf{rep}_{\mathcal{G}}(L)^{\mathrm{op}} \longrightarrow \mathsf{rep}_{\mathcal{G}}(L).$$

In order to compute it explicitly, we consider two elements $p, q \in Q$ as well as two representations $(M,\rho) \in \mathsf{rep}_{\mathcal{G}}(L)_p$ and $(N,\tau) \in \mathsf{rep}_{\mathcal{G}}(L)_q$. The Day convolution of the functors corresponding to $M$ and $N$ is computed via a certain colimit and $\widehat{\Theta}$, as an equivalence of categories, commutes with this colimit. Thus, by Equation (2.8.4) and the definition of coends, the homogeneous component $(M \otimes N)_x$ of degree $x \in Q$ is the coequaliser of

$$\coprod_{r,s\in Q} \Bbbk L_{r,s} \otimes_{\Bbbk} M_r \otimes_{\Bbbk} N_{s^{-1}x} \rightrightarrows \coprod_{r\in Q} M_r \otimes_{\Bbbk} N_{r^{-1}x},$$

where $L_{r,s} = L$ if $r = s$ and $\varnothing$ otherwise. Its parallel morphisms are

$$l \otimes_{\Bbbk} m \otimes_{\Bbbk} n \longmapsto F\rho(l)m \otimes_{\Bbbk} n,$$
$$l \otimes_{\Bbbk} m \otimes_{\Bbbk} n \longmapsto m \otimes \tau(l^{\vee})n = m \otimes \tau(\overline{\alpha}(s^{-1})l)n,$$

where $l \in L_{r,s}$, $m \in M_r$ and $n \in N_{s^{-1}x}$.

In order to compute $D$, we note that $\mathsf{R}(F) = F({}^{\vee}-)^*$ for any $F \in [\Bbbk\mathcal{G}, \mathsf{vect}]$. Thus, given $x \in Q$, we compute

$$DM_x = \begin{cases} M^* & x^{-1} = p, \\ \{0\} & \text{otherwise.} \end{cases}$$

Hence, the action of any $l \in L$ on $M^*$ is given by $\big(\rho(\overline{\alpha}(p)l)\big)^*\colon M_p^* \longrightarrow M_p^*$. □

**Proposition 3.66.** *Let $(G, H, \alpha, t)$ be a crossed module with $G$ and $H$ finite and write $\mathcal{G}$ for its associated strict 2-group. The category $\mathsf{rep}_{\mathcal{G}}(L)$ is rigid if and only if $\operatorname{char}\Bbbk$ does not divide the order of $L = \ker t$.*

*Proof.* It follows from Proposition 3.45 that $\mathsf{rep}_{\mathcal{G}}(L)$ is rigid if and only if all of its objects are finitely-generated and projective. As $\mathsf{rep}_{\mathcal{G}}(L)$ is a direct sum of finitely many copies of $\mathsf{rep}(L)$ and every object of $\mathsf{rep}(L)$ is finitely-generated, this is the case if and only if all objects of $\mathsf{rep}(L)$ are projective. The latter is equivalent to $\Bbbk L$ being semisimple, which corresponds to $\operatorname{char}\Bbbk$ and $|L|$ being coprime by Maschke's theorem. □

Ich gehe am Strand, dudududu;
ich gehe am Strand, rududu;
ich gehe am Strand, dudududuu;
ich gehe am Strand, ruududuu.

SpongeBob Schwammkopf, S01E18b

# 4 TWISTED CENTRES

A peculiarity of Hopf-cyclic cohomology in the sense of Connes and Moscovici is the lack of "canonical" coefficients [CM99]. Originally, modular pairs in involution were considered [CM00]. These consist of a group-like element and a character of the Hopf algebra under consideration implementing the square of the antipode by their respective adjoint actions. Later, Hajac, Khalkhali, Rangipour, and Sommerhäuser obtained a more general source of coefficients in the category of anti-Yetter–Drinfeld modules, [HKRS04]. As mentioned in Example 2.42, anti-Yetter–Drinfeld modules do not form a monoidal category, but rather a module category over the Yetter–Drinfeld modules. This is reflected by the fact that they can be identified with the modules over the anti-Drinfeld double, a comodule algebra over the Drinfeld double, see Remark 2.79. The special role of pairs in involution is captured by the following theorem due to Hajac and Sommerhäuser.

**Theorem 4.1** ([Hal21, Theorem 3.4])**.** *For a finite-dimensional Hopf algebra $H$, the following statements are equivalent*:

(*i*) *The Hopf algebra $H$ admits a pair in involution.*
(*ii*) *There exists a one-dimensional anti-Yetter–Drinfeld module over $H$.*
(*iii*) *The Drinfeld double and anti-Drinfeld double of $H$ are isomorphic algebras.*

Pairs in involution are of categorical interest because they give rise to pivotal structures on the Yetter–Drinfeld modules. Using the theory of heaps, certain algebraic structures equipped with a ternary operation, we can lift the Picard group of a space into the Picard heap of a category. These can be seen as an analogue of the classical pairs in involution. Further abstracting the anti-Yetter–Drinfeld modules into the anti-centre, we can reformulate Theorem 4.1 in more categorical terms, emphasising the role of pivotal structures.

In Theorem 6.44 we shall give a monadic interpretation of this result.

**Theorem 4.23.** *Let $\mathscr{C}$ be a rigid category. There is a bijection between*

(*i*) *equivalence classes of quasi-pivotal structures on $\mathscr{C}$,*
(*ii*) *the Picard heap of the anti-centre, and*
(*iii*) *isomorphism classes of equivalences of module categories between the centre and the anti-centre.*

Further investigating the heap structure of the anti-centre of $\mathscr{C}$, we study the natural injective map $\kappa\colon \operatorname{Pic} \mathsf{A}(\mathscr{C}) \longrightarrow \operatorname{Piv} \mathsf{Z}(\mathscr{C})$, see Theorem 4.33, first introduced by Shimizu [Shi19]. By [Shi23a, Theorem 4.1], this arrow is always bijective if $\mathscr{C}$ is a finite tensor category. However, we show that this fails to be true in general, confirming a conjecture raised in the introduction of *ibid*.

**Theorem 4.50.** *There exists a pivotal structure on a category $\mathscr{C}$ that is not induced by any element of the Picard heap of the anti-centre of $\mathscr{C}$. In other words, the map $\kappa\colon \operatorname{Pic} \mathsf{A}(\mathscr{C}) \longrightarrow \operatorname{Piv} \mathsf{Z}(\mathscr{C})$ is not surjective.*

## 4.1 Heaps

Heaps can be thought of as groups without a fixed neutral element. They were first studied under the name *Schar*, [Prü24; Bae29]. Recently, their homological properties were studied in [ESZ21], and a generalisation towards a "quantum version" is hinted at in [Ško07]. Heaps are equivalent to *G-torsors*—nonempty sets with a freely transitive $G$-action—for a group $G$, see for example [Gir71]. We follow Section 2 of [Brz20] for our exposition.

**Definition 4.2.** A *heap* is a set $G$ together with a ternary operation

$$\langle -, =, \equiv \rangle \colon G \times G \times G \longrightarrow G,$$

which we call the *heap operation*, satisfying a generalised associativity axiom and the *Mal'cev identities*, which we think of as unitality axioms:

$$\langle g, h, \langle i, j, k \rangle \rangle = \langle \langle g, h, i \rangle, j, k \rangle, \qquad \text{for all } g, h, i, j, k \in G, \tag{4.1.1}$$

$$\langle g, g, h \rangle = h = \langle h, g, g \rangle, \qquad \text{for all } g, h \in G. \tag{4.1.2}$$

There are two peculiarities we would like to point out. First, our definition does intentionally not exclude the empty set from being a heap. Second, due to a slightly different setup, an additional "middle associativity axiom" is required in [HL17]. However, as noted in [Brz20, Lemma 2.3], this is implied by Equations (4.1.1) and (4.1.2).

**Definition 4.3.** A map $f\colon G \longrightarrow H$ between heaps is a *morphism of heaps* if

$$f\big(\langle g,h,i\rangle\big) = \langle f(g), f(h), f(i)\rangle, \qquad \text{for all } g,h,i \in G.$$

The next lemma can be shown analogously to its group-theoretical version.

**Lemma 4.4.** *A morphism of heaps is an isomorphism if and only if it is bijective.*

By forgetting its unit, any group defines a heap, Conversely, any non-empty heap can be turned into a group by choosing a unit, see [Cer43].

**Lemma 4.5.** *Every group $(G,\cdot,e)$ is a heap via*

$$\langle -,=,\equiv\rangle\colon G\times G\times G \longrightarrow G, \qquad \langle g,h,i\rangle := g\cdot h^{-1}\cdot i.$$

*A morphism of groups becomes a morphism of the induced heaps.*

**Lemma 4.6.** *A non-empty heap $G$ with a fixed element $e \in G$ can be considered as a group with unit $e$ via the multiplication*

$$-\cdot_e = \colon G\times G \longrightarrow G, \qquad g\cdot_e h := \langle g,e,h\rangle.$$

*The inverse of an element $g \in G$ with respect to $\cdot_e$ is given by $g^{-1} := \langle e,g,e\rangle$. A morphism of heaps is a morphism of the induced groups, provided it maps the fixed element of its source to the fixed element of its target.*

More generally, if $\mathscr{G}$ is a groupoid and $x,y \in \mathscr{G}$, then the set of morphisms $\mathscr{G}(x,y)$ becomes a heap with heap operation given by $\langle f,g,h\rangle := fg^{-1}h$. The next example is a special case of this construction, which will play a prominent role in this chapter.

**Example 4.7.** Let $F,G\colon \mathscr{C} \longrightarrow \mathscr{C}$ be oplax monoidal endofunctors. The set

$$\mathrm{Iso}_\otimes(F,G) := \big\{\text{oplax monoidal natural isomorphisms from } F \text{ to } G\big\}$$

bears a heap structure with heap operation

$$\langle -,=,\equiv\rangle\colon \mathrm{Iso}_\otimes(F,G)^3 \longrightarrow \mathrm{Iso}_\otimes(F,G), \qquad \langle \phi,\psi,\xi\rangle := \phi\psi^{-1}\xi.$$

## 4.2 Pivotal structures and twisted centres

### 4.2.1 *Twisted centres and their Picard heaps*

Recall from Example 2.45 that we can twist the regular action of a monoidal category on itself with strong monoidal functors. If $\mathscr{C}$ is a rigid monoidal category, then we refer to the centre $\mathsf{Z}({}_L\mathscr{C}_R)$ as a *twisted centre*, and to $\mathsf{Z}(\mathscr{C}_R)$ and $\mathsf{Z}({}_L\mathscr{C})$ as *right* and *left* twisted centres.

**Remark 4.8.** A natural description of twisted centres is obtained from the perspective of bicategories. Given a monoidal category $\mathscr{C}$, let $\mathbf{B}\mathscr{C}$ denote its *delooping*; that is, $\mathbf{B}\mathscr{C}$ is a bicategory with a single object $\bullet$, and $\mathbf{B}\mathscr{C}(\bullet, \bullet) \simeq \mathscr{C}$. In this setting, strong monoidal functors are identified with pseudofunctors between one-object bicategories, and strong monoidal transformations correspond to pseudo-icons, see for example [JY21, Proposition 4.6.9].

Let $\mathbb{P}\mathsf{s}^{\mathsf{ps}}(\mathbf{B}\mathscr{C})$ be the bicategory of endopseudofunctors of $\mathbf{B}\mathscr{C}$, their pseudonatural transformations, and modifications. Then $\mathbb{P}\mathsf{s}^{\mathsf{ps}}(\mathbf{B}\mathscr{C})(L, R) \simeq \mathsf{Z}({}_L\mathscr{C}_R)$. See [FH23, Proposition 3.6] for a proof of this fact.

For the rest of this section, fix a strict rigid category $\mathscr{C}$, see Theorem 2.72. The forgetful functor from the centre of a twisted bimodule category to the underlying monoidal category is faithful, so we can use a variant of the graphical calculus for monoidal categories, as long as we pay special attention to the half-braidings. Whence, we introduce a colouring scheme to help us keep track of the various categories.

(i) Red for objects in the right twisted centre $\mathsf{Z}(\mathscr{C}_R)$,
(ii) blue for objects in the left twisted centre $\mathsf{Z}({}_L\mathscr{C})$, and
(iii) black for objects in the Drinfeld centre $\mathsf{Z}(\mathscr{C})$ or in $\mathscr{C}$.

For example, the half-braidings of objects $a \in \mathsf{Z}(\mathscr{C}_R)$ and $q \in \mathsf{Z}({}_L\mathscr{C})$ are:

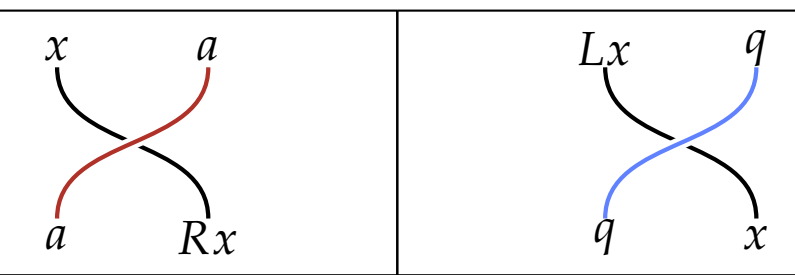

| $\sigma_{a,x} : a \otimes Rx \longrightarrow x \otimes a$ | $\sigma_{q,x} : q \otimes x \longrightarrow Lx \otimes q$ |
|---|---|

**Hypothesis 4.9.** In the rest of this chapter, we are predominantly interested in twisting with the same strict monoidal functor from the left or right. For the purpose of brevity, we therefore fix such a functor $L = R : \mathscr{C} \longrightarrow \mathscr{C}$ and consider the categories ${}_L\mathscr{C}$ and $\mathscr{C}_R$.

Suppose we are given three objects

$$(a, \sigma_{a,-}) \in \mathsf{Z}(\mathscr{C}_R), \qquad (q, \sigma_{q,-}) \in \mathsf{Z}({}_L\mathscr{C}), \qquad \text{and} \qquad (x, \sigma_{x,-}) \in \mathsf{Z}(\mathscr{C}).$$

**Proposition 4.10.** *The tensor product of $\mathscr{C}$ extends to a left action of $\mathsf{Z}(\mathscr{C})$ on $\mathsf{Z}(\mathscr{C}_R)$ and a right action of $\mathsf{Z}(\mathscr{C})$ on $\mathsf{Z}({}_L\mathscr{C})$. The half-braidings are as defined in Figure 4.1.*

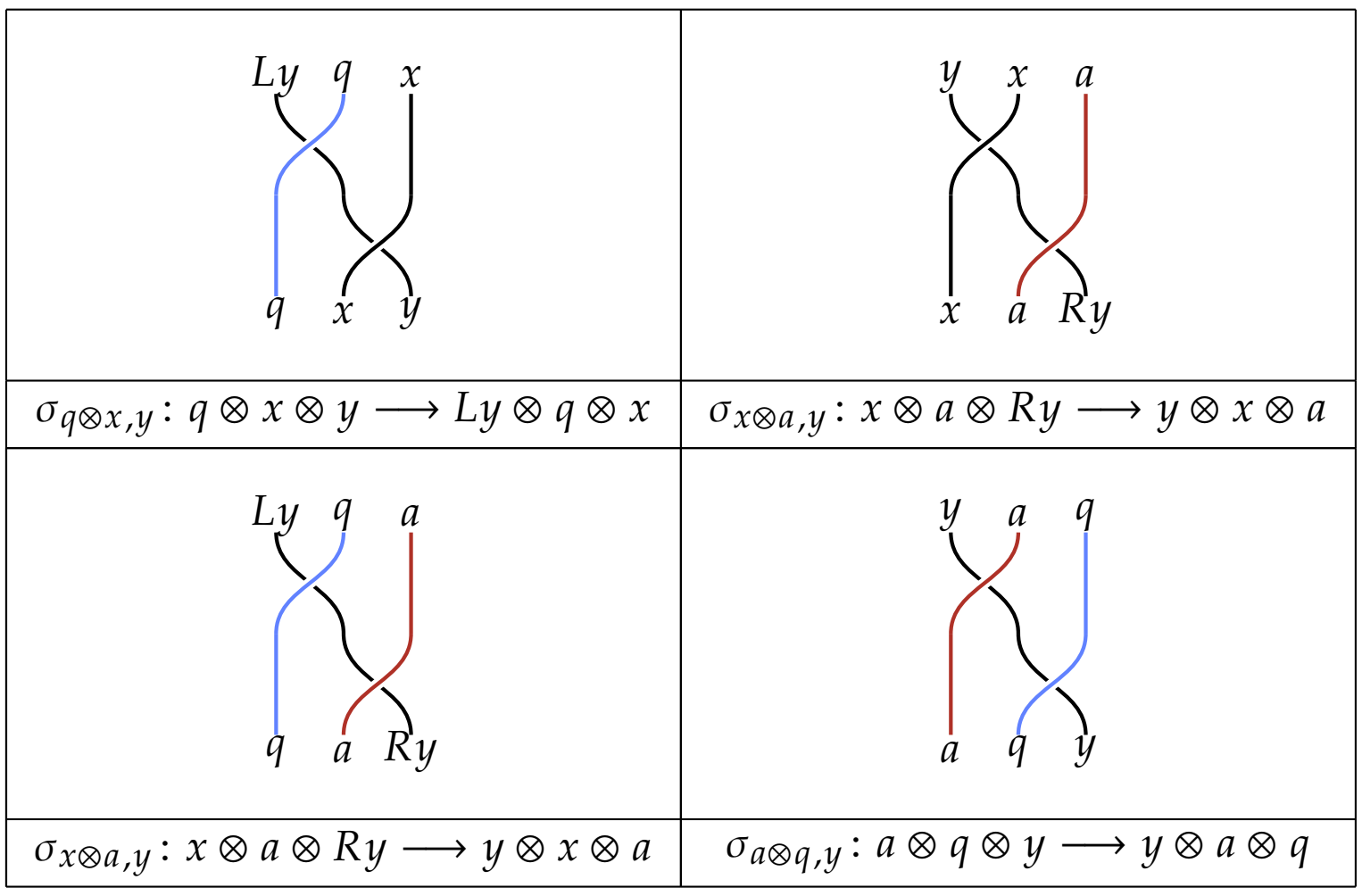

Figure 4.1: Canonical actions of $\mathsf{Z}(\mathscr{C})$ on $\mathsf{Z}({}_L\mathscr{C})$ and $\mathsf{Z}(\mathscr{C}_R)$.

**Remark 4.11.** A direct computation proves the categories $\mathsf{Z}(\mathscr{C}_R^{\mathrm{op,rev}})$ and $\mathsf{Z}({}_R\mathscr{C})^{\mathrm{op}}$ to be equivalent. Furthermore, this identification is compatible with the module structure: for all $x \in \mathsf{Z}(\mathscr{C}^{\mathrm{op,rev}})$ and $a \in \mathsf{Z}(\mathscr{C}_R^{\mathrm{op,rev}})$, we have $x \otimes^{\mathrm{op}} R(a) = R(a) \otimes x$ and $\sigma_{x\otimes^{\mathrm{op}} a,-} = \sigma_{a\otimes x,-}$. Hence, we restrict ourselves to the study of right twisted centres.

Recall from Proposition 2.80 that, since $\mathscr{C}$ is a rigid monoidal category, the centre $\mathsf{Z}(\mathscr{C})$ is rigid as well. While the same cannot be said of the twisted centre, the left dual ${}^\vee a$ of any object $(a, \sigma_{a,-}) \in \mathsf{Z}(\mathscr{C}_R)$ can be turned into an object of $\mathsf{Z}({}_R\mathscr{C})$ if we equip it with the half-braiding

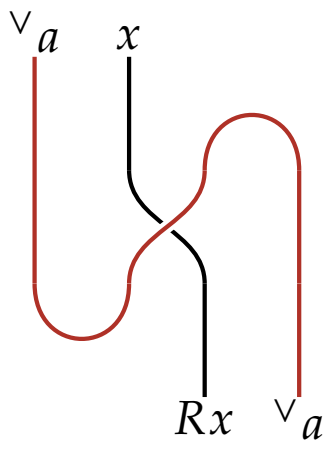

This suggests that, in analogy with Proposition 2.80, we may lift the dualising functor of $\mathscr{C}$ to the level of twisted centres by interchanging right with left twists. A more conceptual description is provided in [FH23, Proposition 3.15].

**Proposition 4.12.** *The left dualising functor* ${}^{\vee}(-)\colon \mathscr{C} \longrightarrow \mathscr{C}^{\mathrm{op,rev}}$ *lifts to a functor between right and left twisted centres*

$$ {}^{\vee}(-)\colon \mathsf{Z}(\mathscr{C}_R) \longrightarrow \mathsf{Z}({}_R\mathscr{C})^{\mathrm{op}}. $$

The half-braidings displayed in the right column of Figure 4.1 show that every object $a \in \mathsf{Z}(\mathscr{C}_R)$ gives rise to two functors of left $\mathsf{Z}(\mathscr{C})$-modules:

$$ - \otimes a\colon \mathsf{Z}(\mathscr{C}) \longrightarrow \mathsf{Z}(\mathscr{C}_R) \qquad \text{and} \qquad - \otimes {}^{\vee}a\colon \mathsf{Z}(\mathscr{C}_R) \longrightarrow \mathsf{Z}(\mathscr{C}). \tag{4.2.1} $$

Let $(a, \sigma_{a,-}) \in \mathsf{Z}(\mathscr{C}_R)$. In view of Proposition 4.12, we use the following notation for the evaluation and coevaluation morphisms of twisted centres:

(4.2.2)

| | |
|---|---|
| 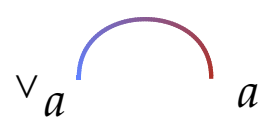 | 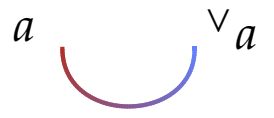 |
| $\mathrm{ev}^{\ell}_a\colon {}^{\vee}a \otimes a \longrightarrow 1$ | $\mathrm{coev}^{\ell}_a\colon 1 \longrightarrow a \otimes {}^{\vee}a$ |

**Proposition 4.13.** *Every object* $a \in \mathsf{Z}(\mathscr{C}_R)$ *induces adjoint* $\mathsf{Z}(\mathscr{C})$*-module functors*

$$ - \otimes a\colon \mathsf{Z}(\mathscr{C}) \rightleftarrows \mathsf{Z}(\mathscr{C}_R) : - \otimes {}^{\vee}a. \tag{4.2.3} $$

*Proof.* We know that for any object $(a, \sigma_{a,-}) \in \mathsf{Z}(\mathscr{C}_R)$, Equation (4.2.3) is an ordinary adjunction, whose unit $\eta$ and counit $\varepsilon$ are given by

$$ \begin{aligned} \eta_y &:= \mathrm{id}_y \otimes \mathrm{coev}^{\ell}_a\colon y \longrightarrow y \otimes a \otimes {}^{\vee}a, && \text{for all } y \in \mathsf{Z}(\mathscr{C}), \\ \varepsilon_x &:= \mathrm{id}_x \otimes \mathrm{ev}^{\ell}_a\colon x \otimes {}^{\vee}a \otimes a \longrightarrow x, && \text{for all } x \in \mathsf{Z}(\mathscr{C}_R). \end{aligned} $$

The next diagram shows that $\varepsilon_x$ is a morphism in $\mathsf{Z}(\mathscr{C}_R)$ for every $x \in \mathsf{Z}(\mathscr{C}_R)$:

(4.2.4)

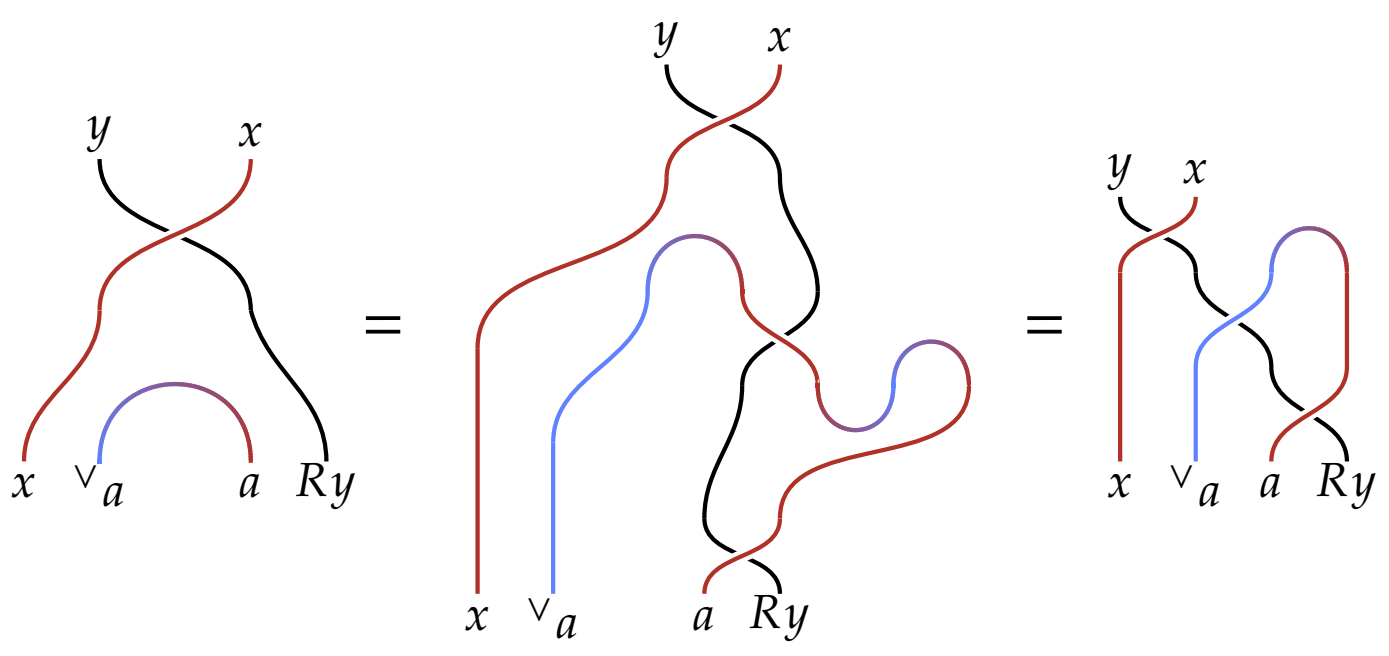

Furthermore, $\varepsilon_{w\otimes x} = \mathrm{id}_w \otimes \varepsilon_x$ for all $w \in \mathsf{Z}(\mathscr{C})$. A similar argument shows that the unit of the adjunction is a natural transformation of module functors. □

Since the forgetful functor from the (twisted) centre to its underlying category is faithful, it is also conservative[11], which allows us to characterise equivalences of module categories between $\mathsf{Z}(\mathscr{C})$ and right twisted centres.

[11] A functor $F\colon \mathscr{C} \longrightarrow \mathscr{D}$ is called *conservative* if it reflects isomorphisms; i.e., if for all $f \in \mathscr{C}(x,y)$, the fact that $Ff$ is an isomorphism implies that $f$ is one.

**Remark 4.14.** Recall from Proposition 2.51 that for a monoidal category $\mathscr{C}$ and a $\mathscr{C}$-module category $\mathscr{M}$, there is an equivalence of $\mathscr{C}$-module categories

$$\mathsf{Str}\mathscr{C}\mathsf{Mod}(\mathscr{C},\mathscr{M}) \xrightarrow{\sim} \mathscr{M}, \qquad F \longmapsto F1, \qquad - \triangleright m \longleftarrow m.$$

In particular, any functor of left module categories $F\colon \mathsf{Z}(\mathscr{C}) \longrightarrow \mathsf{Z}(\mathscr{C}_R)$ is naturally isomorphic to $- \otimes F1\colon \mathsf{Z}(\mathscr{C}) \longrightarrow \mathsf{Z}(\mathscr{C}_R)$. By Proposition 2.67, $F$ is an equivalence if and only if $F1 \in \mathscr{C}$ is invertible, and two left $\mathsf{Z}(\mathscr{C})$-module functors $F, G\colon \mathsf{Z}(\mathscr{C}) \longrightarrow \mathsf{Z}(\mathscr{C}_R)$ are isomorphic if and only if $F1 \cong G1$.

**Definition 4.15.** An object $(\alpha, \sigma_{\alpha,-}) \in \mathsf{Z}(\mathscr{C}_R)$ in is called *$\mathscr{C}$-invertible* if the image $U\alpha \in \mathscr{C}$ of $\alpha$ under the forgetful functor $U^{(R)}\colon \mathsf{Z}(\mathscr{C}_R) \longrightarrow \mathscr{C}$ is invertible.

**Notation 4.16.** Let $(\alpha, \sigma_{(\alpha,-)}) \in \mathsf{Z}(\mathscr{C}_R)$ be a $\mathscr{C}$-invertible element. Analogously to Equation (4.2.2), we use the following notation for the inverses $\mathrm{ev}^{-\ell}$ and $\mathrm{coev}^{-\ell}$ of $\mathrm{ev}^{\ell}$ and $\mathrm{coev}^{\ell}$, respectively:

| 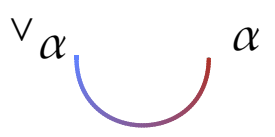 | 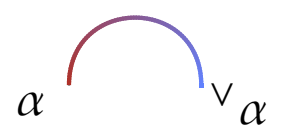 |
|---|---|
| $\mathrm{ev}_\alpha^{-\ell}\colon 1 \longrightarrow {}^\vee\alpha \otimes \alpha$ | $\mathrm{coev}_\alpha^{-\ell}\colon \alpha \otimes {}^\vee\alpha \longrightarrow 1$ |

The notion of heaps allows us to define an algebraic structure on the isomorphism classes of objects implementing module equivalences between the Drinfeld centre $\mathsf{Z}(\mathscr{C})$ and its twisted relative $\mathsf{Z}(\mathscr{C}_R)$. In analogy with the Picard group, we call this the Picard heap of a twisted centre.

**Definition 4.17.** The *Picard heap* of the right twisted centre $\mathsf{Z}(\mathscr{C}_R)$ is the set of isomorphism classes

$$\operatorname{Pic}\mathsf{Z}(\mathscr{C}_R) := \left\{[\alpha] \;\middle|\; \alpha \in \mathsf{Z}(\mathscr{C}_R) \text{ is } \mathscr{C}\text{-invertible}\right\}$$

together with the heap operation defined for $[\alpha], [\beta], [\gamma] \in \operatorname{Pic}\mathsf{Z}(\mathscr{C}_R)$ by

$$\langle [\alpha], [\beta], [\gamma] \rangle = [\alpha \otimes {}^\vee\beta \otimes \gamma].$$

**Lemma 4.18.** *The Picard heap defined in Definition 4.17 forms a heap.*

*Proof.* The generalised associativity, see Equation (4.1.1), follows from the associativity of the tensor product of $\mathscr{C}$ and its compatibility with the "gluing" of half-braidings. To show that the Mal'cev identities hold, fix $\mathscr{C}$-invertible objects $\alpha, \beta \in \mathsf{Z}(\mathscr{C}_R)$. Proposition 2.80 and Equation (4.2.4) imply that

$$\alpha \otimes {}^{\vee}\alpha \otimes \beta \xrightarrow{\mathrm{coev}_\alpha^{\ell\ -1} \otimes \mathrm{id}_\beta} \beta \qquad \text{and} \qquad \beta \otimes {}^{\vee}\alpha \otimes \alpha \xrightarrow{\mathrm{id}_\beta \otimes \mathrm{ev}_\alpha^\ell} \beta$$

are isomorphisms in $\mathsf{Z}(\mathscr{C}_R)$, hence $\langle [\alpha], [\alpha], [\beta] \rangle = [\beta] = \langle [\beta], [\alpha], [\alpha] \rangle$. □

In general, the twisted centre $\mathsf{Z}(\mathscr{C}_R)$ does not inherit a monoidal structure from $\mathscr{C}$. Lemma 4.18, however, hints towards a slight generalisation where the tensor product is replaced by a trivalent functor, essentially categorifying heaps (without the Mal'cev identities). The well-definedness of this concept was hinted at in [Ško07] under the name of *heapy categories*.

### 4.2.2 *Quasi-pivotality*

A particularly interesting consequence of our previous findings arises in the case of $R := {}^{\vee\vee}(-)$; the centre of the regular bimodule twisted on the right by $R$ can be understood as a generalisation of anti-Yetter–Drinfeld modules, see [HKS19, Theorem 2.3]. As before, fix a strict rigid category $\mathscr{C}$ and consider the twisted bimodules categories $\mathscr{C}_{{}^{\vee\vee}(-)}$ and ${}_{{}^{\vee\vee}(-)}\mathscr{C}$.

**Notation 4.19.** We write $\mathsf{A}(\mathscr{C}) := \mathsf{Z}(\mathscr{C}_{{}^{\vee\vee}(-)})$ and $\mathsf{Q}(\mathscr{C}) := \mathsf{Z}({}_{{}^{\vee\vee}(-)}\mathscr{C})$, and call the former the *anti-Drinfeld centre* of $\mathscr{C}$.

We have already mentioned the connection between the twisted centre $\mathsf{A}(\mathscr{C})$ and anti-Yetter–Drinfeld modules over Hopf algebras given in [HKS19]. The case where $\mathscr{C}$ is the category of modules over a Hopf algebroid was investigated by Kowalzig in [Kow24]. The counterpart $\mathsf{Q}(\mathscr{C})$ of the generalised anti-Yetter–Drinfeld modules is less common in the literature, but plays a crucial role in the monadic investigations of Chapters 5 and 6.

**Definition 4.20** ([Shi23a, Section 4])**.** A *quasi-pivotal structure* on a rigid monoidal category $\mathscr{C}$ is a pair $(\beta, \rho_\beta)$ consisting of an invertible object $\beta \in \mathscr{C}$ and a monoidal natural isomorphism

$$\rho_\beta \colon (-) \xRightarrow{\sim} \beta \otimes {}^{\vee\vee}(-) \otimes {}^{\vee}\beta.$$

The tuple $(\mathscr{C}, (\beta, \rho_\beta))$ will be called a *quasi-pivotal* category.

If $\mathscr{C}$ is the category of finite-dimensional modules over a finite-dimensional Hopf algebra, quasi-pivotal structures have a well-known interpretation—they translate to pairs in involution. This can be deduced from a slight variation of [Hal21, Lemma 5.6], the main observation being that the invertible object $\beta$ of a quasi-pivotal structure $(\beta, \rho_\beta)$ on $\mathscr{C}$ corresponds to a character, and that $\rho_\beta$ determines a group-like element. The fact that $\rho_\beta$ is a natural transformation from the identity to a conjugate of the double dual functor is captured by the character and group-like implementing the square of the antipode. We study a monadic analogue of this in Section 6.5.

**Remark 4.21.** Every pivotal category is quasi-pivotal; the converse does not hold. A counterexample are the finite-dimensional modules over the generalised Taft algebras discussed in [HK19]. Any of these Hopf algebras admit pairs in involution but in general neither the character nor the group-like can be trivial. The previous discussion and [Hal21, Lemma 5.6] show that Mod-$H$ is quasi-pivotal but not pivotal—in contrast to its Drinfeld centre $\mathsf{Z}$(Mod-$H$), which admits a pivotal structure by [Hal21, Lemma 5.5].

Let $(\beta, \rho_\beta)$ be a quasi-pivotal structure on $\mathscr{C}$ and $\phi\colon \beta' \longrightarrow \beta$ an isomorphism in $\mathscr{C}$. Clearly, the pair $(\beta', (\phi^{-1} \otimes \mathrm{id} \otimes {}^\vee\phi) \circ \rho_\beta)$ is another quasi-pivotal structure on $\mathscr{C}$. This defines an equivalence relation and we write

$$\mathrm{QPiv}(\mathscr{C}) := \{\, [(\beta, \rho_\beta)] \mid (\beta, \rho_\beta) \text{ is a quasi-pivotal structure on } \mathscr{C} \,\}$$

for the set of equivalence classes of quasi-pivotal structures on $\mathscr{C}$.

**Lemma 4.22.** *Let $\mathscr{C}$ be a strict rigid category. The Picard heap* Pic $\mathsf{A}(\mathscr{C})$ *and the set of equivalence classes of quasi-pivotal structures* $\mathrm{QPiv}(\mathscr{C})$ *are in bijection.*

*Proof.* Let $(\beta, \rho_\beta)$ be a quasi-pivotal structure on $\mathscr{C}$. Writing $\mathrm{ev}_\beta^{-\ell} := (\mathrm{ev}_\beta^{\ell})^{-1}$, we define the half-braiding

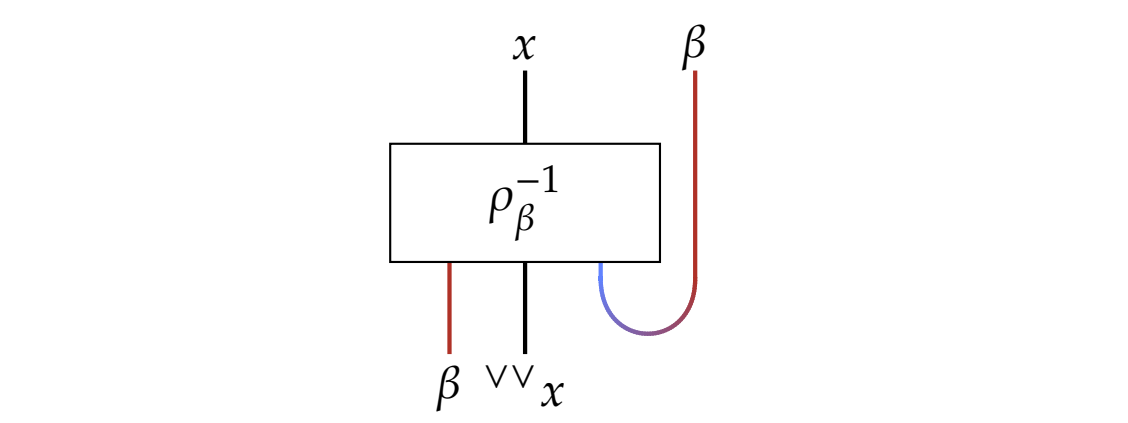

$\sigma_{\beta,x} := (\rho_{\beta,x}^{-1} \otimes \beta) \circ (\mathrm{id} \otimes \mathrm{ev}_\beta^{-\ell})\colon \beta \otimes {}^{\vee\vee}x \longrightarrow x \otimes \beta$

Then, as $\rho_\beta$ is monoidal, $\sigma_{\beta,x}$ satisfies the hexagon identity. This defines

$$\phi\colon\ \mathrm{QPiv}(\mathscr{C}) \longrightarrow \operatorname{Pic}\mathsf{A}(\mathscr{C}), \qquad [(\beta,\rho_\beta)] \longmapsto [(\beta,\sigma_{\beta,-})].$$

Conversely, let $(\alpha,\sigma_{\alpha,-}) \in \mathsf{A}(\mathscr{C})$ be $\mathscr{C}$-invertible. From its half-braiding we obtain a monoidal natural transformation

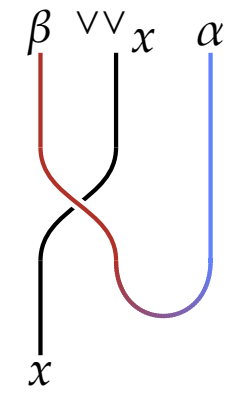

$\sigma_{\beta,x} = (\rho^{-1}_{\beta,x} \otimes \beta) \circ (\mathrm{id} \otimes \mathrm{ev}^{-\ell}_{\beta})\colon \beta \otimes {}^{\vee\vee}x \longrightarrow x \otimes \beta$

Due to the snake identities, the following map is the inverse of $\phi$:

$$\psi\colon\ \operatorname{Pic}\mathsf{A}(\mathscr{C}) \longrightarrow \mathrm{QPiv}(\mathscr{C}), \qquad [(\alpha,\sigma_{\alpha,-})] \longmapsto [(\alpha,\rho_\alpha)]. \qquad \square$$

**Theorem 4.23.** *Let $\mathscr{C}$ be a strict rigid category. Then there is a bijection between*

(i) *equivalence classes of quasi-pivotal structures on $\mathscr{C}$,*
(ii) *the Picard heap of $\mathsf{A}(\mathscr{C})$, and*
(iii) *isomorphism classes of $\mathsf{Z}(\mathscr{C})$-module equivalences between $\mathsf{Z}(\mathscr{C})$ and $\mathsf{A}(\mathscr{C})$.*

*Proof.* The equivalence of (ii) and (iii) follows by Remark 4.14, and for (i) being equivalent to (ii) we invoke Lemma 4.22. $\square$

### 4.2.3 *Pivotality of the Drinfeld centre*

In this section we shall examine the relationship between pairs in involution and pivotal structures, see for example Remark 4.21, from a categorical perspective. Let $a := (\alpha,\sigma_{\alpha,-}) \in \mathsf{A}(\mathscr{C})$ be $\mathscr{C}$-invertible and write $\Omega := (\omega,\sigma_{\omega,-}) \in \mathsf{Q}(\mathscr{C})$ for its left dual. Figure 4.2 collects some useful properties of the coevaluation of $\alpha$ and its inverse $\mathrm{coev}^{-\ell}_{\alpha}$.

**Definition 4.24.** Let $a := (\alpha,\sigma_{\alpha,-}) \in \mathsf{A}(\mathscr{C})$ be a $\mathscr{C}$-invertible element with left dual $\Omega := (\omega,\sigma_{\omega,-}) \in \mathsf{Q}(\mathscr{C})$. For any $x \in \mathsf{Z}(\mathscr{C})$, define the morphism $\rho_{a,x}$ by

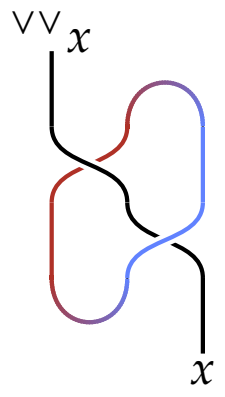

$\rho_{a,x} := ({}^{\vee\vee}x \otimes \mathrm{coev}^{-\ell}_{\alpha}) \circ (\sigma^{-1}_{{}^{\vee\vee}x,\alpha} \otimes \omega) \circ (\alpha \otimes \sigma_{\omega,x}) \circ (\mathrm{coev}^{\ell}_{\alpha} \otimes x)\colon x \longrightarrow {}^{\vee\vee}x$

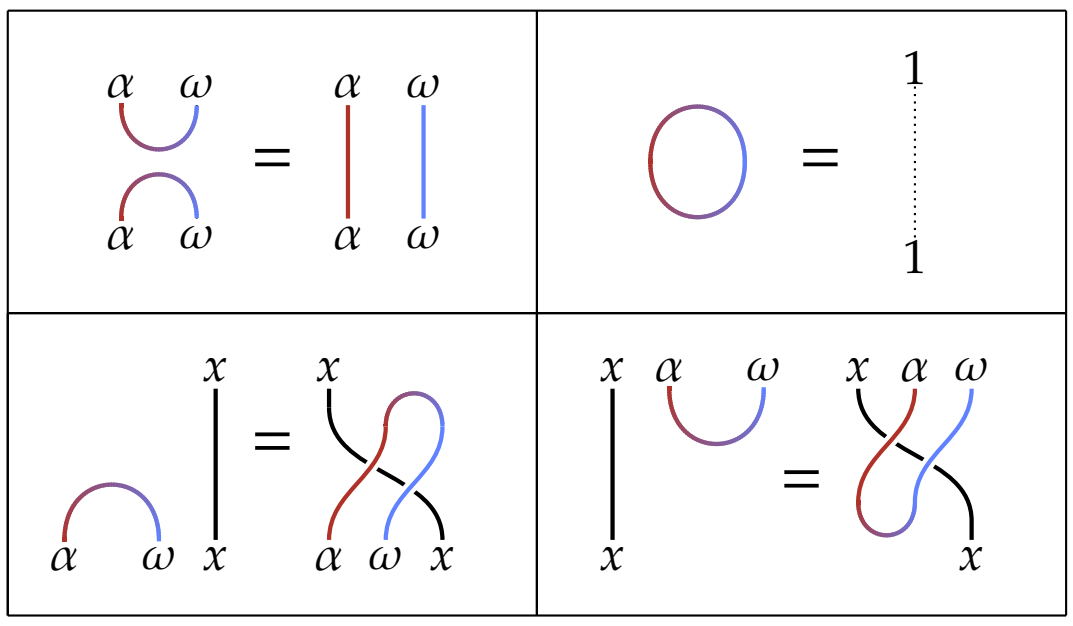

Figure 4.2: Properties of $\mathrm{coev}^{\ell}_{\alpha}$ and $\mathrm{coev}^{-\ell}_{\alpha}$.

As mentioned in Remark 4.8, objects in twisted centres correspond to pseudonatural transformations between deloopings of monoidal functors. The following lemma shows that if these objects are $\mathscr{C}$-invertible, then one can reconstruct monoidal natural isomorphisms from them.

**Lemma 4.25** ([Shi23a, Section 4.4]). *For any $\mathscr{C}$-invertible object $a \in \mathsf{A}(\mathscr{C})$ the map $\rho_{a,-}\colon (-) \Longrightarrow {}^{\vee\vee}(-)$ defines a pivotal structure on $\mathsf{Z}(\mathscr{C})$.*

*Proof.* Fix an object $a := (\alpha, \sigma_{\alpha,-}) \in \mathsf{A}(\mathscr{C})$, such that $\alpha$ is invertible in $\mathscr{C}$, and write $\Omega := (\omega, \sigma_{\omega,-}) \in \mathsf{Q}(\mathscr{C})$ for its left dual. Furthermore, we assume $x \in \mathsf{Z}(\mathscr{C})$ to be any object in the Drinfeld centre of $\mathscr{C}$. We note that for any $y \in \mathscr{C}$ a variant of the Yang–Baxter identity holds:

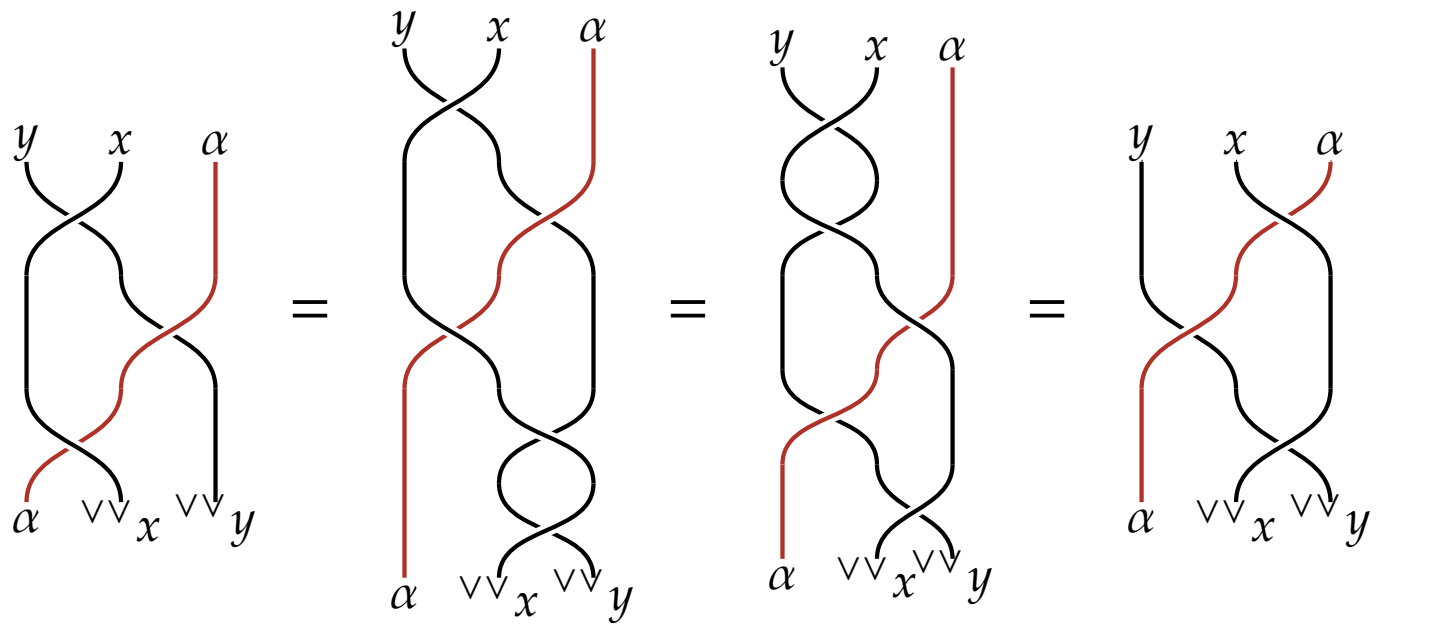

Combining this fact with Figure 4.2, one obtains that $\rho_{a,x}\colon x \longrightarrow {}^{\vee\vee}x$ is a morphism in the Drinfeld centre of $\mathscr{C}$:

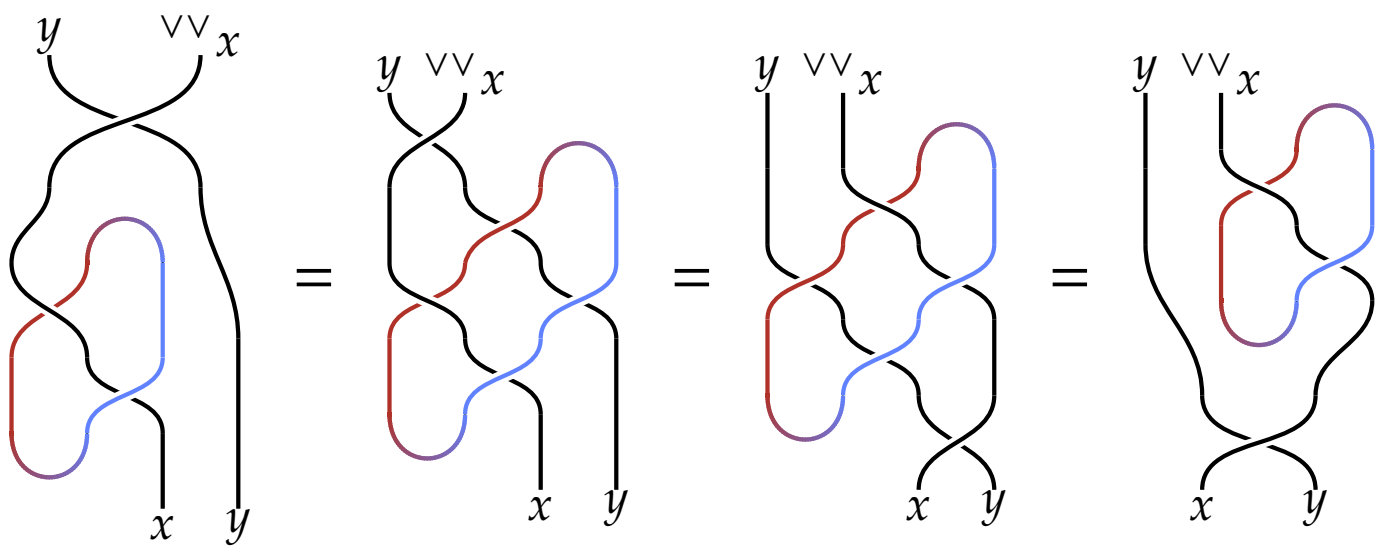

Since the forgetful functor $U^{(Z)}\colon \mathsf{Z}(\mathscr{C}) \longrightarrow \mathscr{C}$ is conservative and $\rho_{a,x}$ is a composite of isomorphisms in $\mathscr{C}$, it is an isomorphism in the centre $\mathsf{Z}(\mathscr{C})$.

The naturality of the half-braidings implies that $\rho_a$ is natural as well:

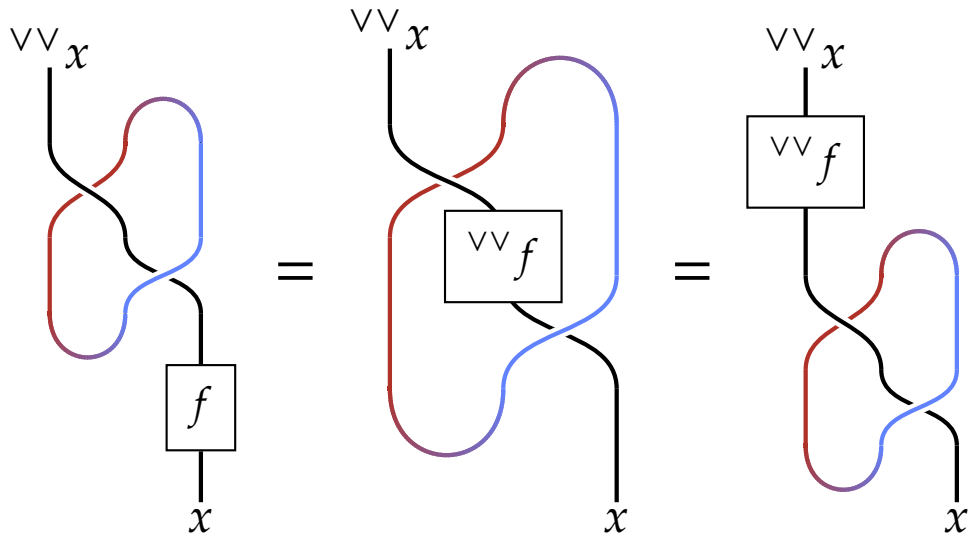

Lastly, the natural isomorphism $\rho_a\colon \mathrm{Id}_{\mathsf{Z}(\mathscr{C})} \longrightarrow {}^{\vee\vee}(-)$ being monoidal is established by the hexagon identities.

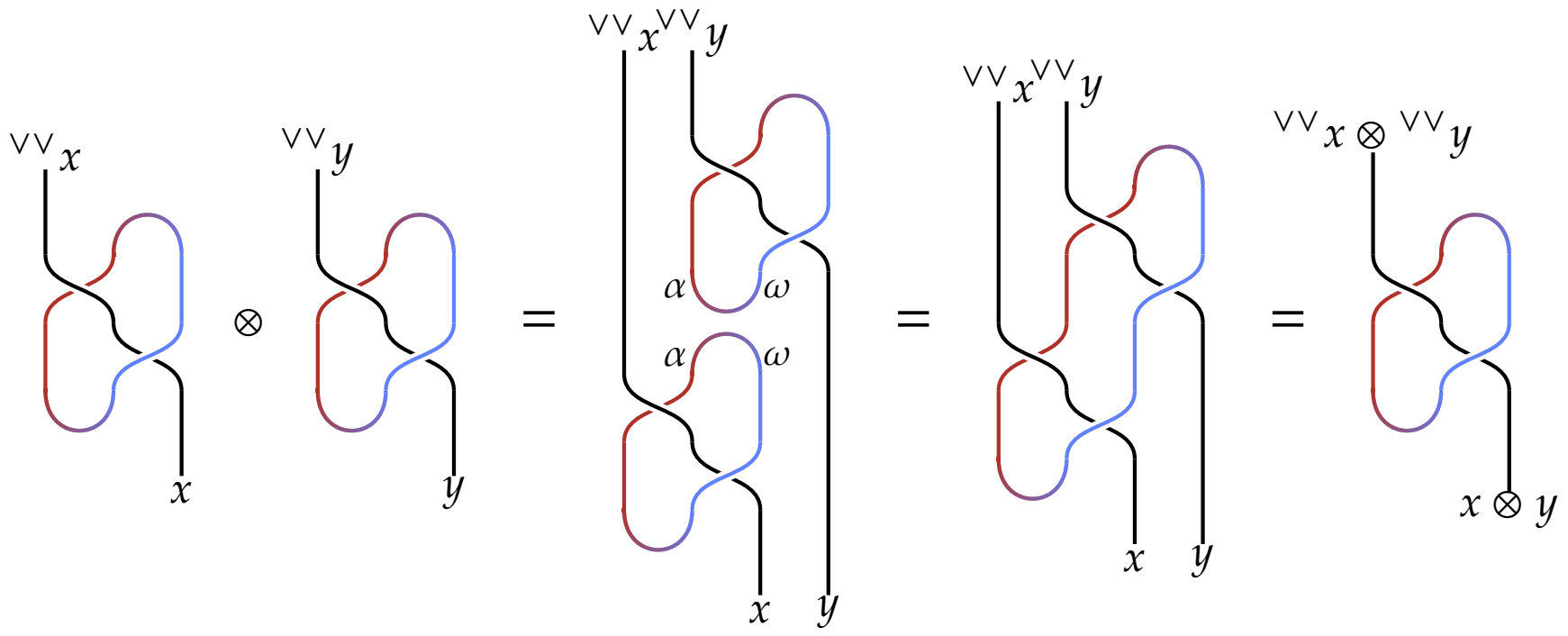

□

Given different pivotal structures of $\mathsf{Z}(\mathscr{C})$, induced by $\mathscr{C}$-invertible objects in $\mathsf{A}(\mathscr{C})$, it is a priori unclear whether they coincide. The following lemma is a first step in this direction. It shows that the induced pivotal structures only depend on the isomorphism classes of $\mathscr{C}$-invertible objects in $\mathsf{A}(\mathscr{C})$.

**Lemma 4.26.** *Suppose that $a_1, a_2 \in \mathsf{A}(\mathscr{C})$ are two representatives of the equivalence class $[a_1] = [a_2] \in \operatorname{Pic}\mathsf{A}(\mathscr{C})$. Then $\rho_{a_1} = \rho_{a_2}$.*

*Proof.* Suppose that there exists an isomorphism $\phi\colon a_1 \longrightarrow a_2$ in the anti-Drinfeld centre. Then by Figure 4.3 the induced pivotal structures $\rho_{a_1}$ and $\rho_{a_2}$ are the same, for any $x \in \mathsf{Z}(\mathscr{C})$. □

**Definition 4.27.** We call an object $x \in \mathsf{Z}(\mathscr{C})$ *symmetric* if we have

$$\sigma_{x,y}^{-1} = \sigma_{y,x}, \qquad \text{for all } y \in \mathsf{Z}(\mathscr{C}).$$

We call the full (symmetric) monoidal subcategory $\mathsf{SZ}(\mathscr{C})$ of $\mathsf{Z}(\mathscr{C})$ whose objects are symmetric the *symmetric* or *Müger centre* of $\mathsf{Z}(\mathscr{C})$, see [Müg13].

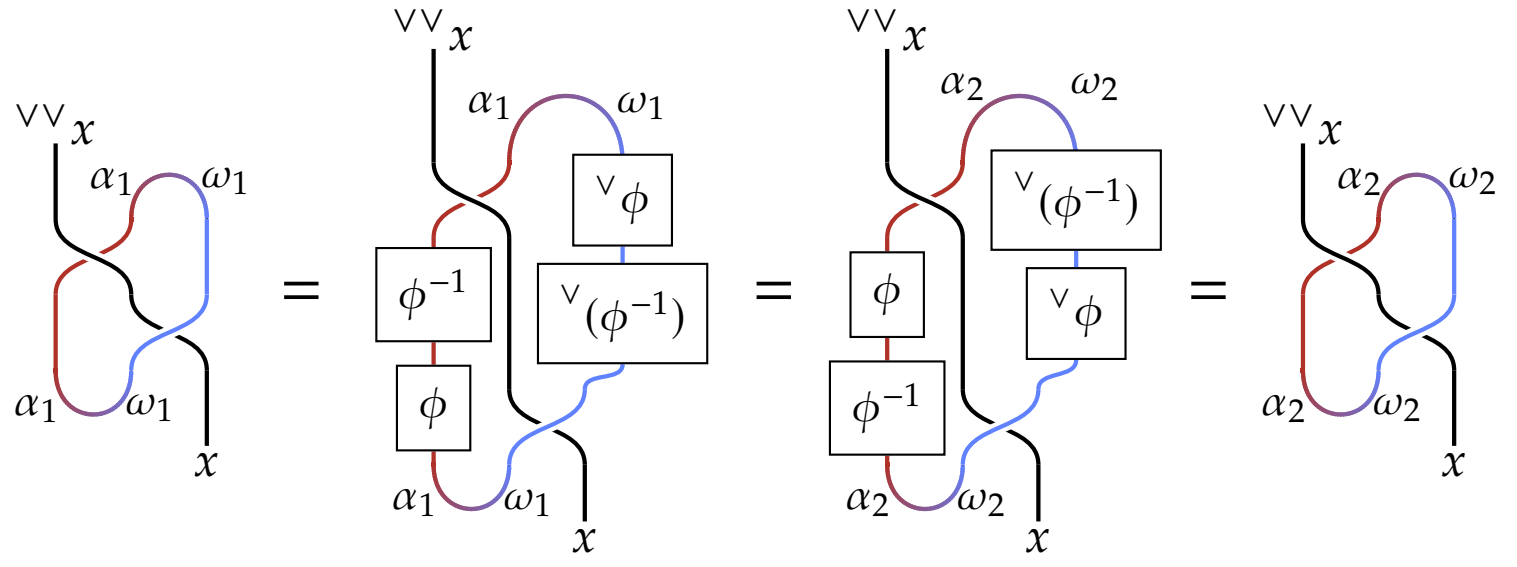

Figure 4.3: The induced pivotal structures of $\rho_{a_1}$ and $\rho_{a_2}$ coincide.

**Lemma 4.28.** *If $\mathscr{C}$ is a rigid monoidal category, then $\mathsf{SZ}(\mathscr{C})$ is rigid monoidal.*

*Proof.* Suppose $x \in \mathsf{Z}(\mathscr{C})$ to be symmetric and let $y \in \mathsf{Z}(\mathscr{C})$. We compute

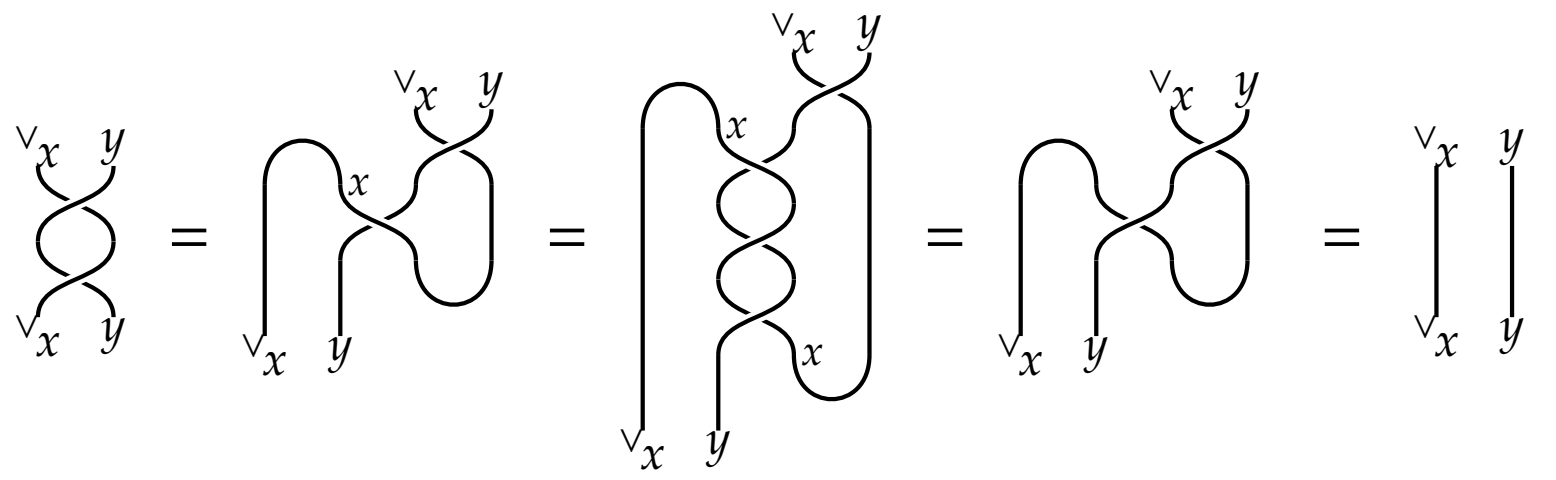

This implies that $\sigma^{-1}_{^\vee x,y} = \sigma_{y,^\vee x}$. Since the left dual of any $x \in \mathsf{SZ}(\mathscr{C}) \subseteq \mathsf{Z}(\mathscr{C})$ can be equipped with the structure of a right dual and $\mathsf{SZ}(\mathscr{C})$ is a full subcategory of $\mathsf{Z}(\mathscr{C})$, it must be rigid. □

Let us now consider the Picard group $\operatorname{Pic}\mathsf{SZ}(\mathscr{C})$. It acts on $\operatorname{Pic}\mathsf{A}(\mathscr{C})$ via tensoring from the left, as shown in Figure 4.1. Two elements $a, c \in \operatorname{Pic}\mathsf{A}(\mathscr{C})$ are equivalent if they are contained in the same orbit; that is,

$$[a] \sim [c] \iff \text{there exists a } [b] \in \operatorname{Pic}\mathsf{SZ}(\mathscr{C}) \text{ such that } [b \otimes a] = [c].$$

To show that two objects are equivalent if and only if they induce the same pivotal structure on $\mathsf{Z}(\mathscr{C})$, we need two technical observations. First, an alternate description of symmetric invertible objects; and second, a more detailed investigation into the inverse of an induced pivotal structure.

**Lemma 4.29.** *An invertible object $(\beta, \sigma_{\beta,-}) \in \mathsf{Z}(\mathscr{C})$ is symmetric if and only if, for all $x \in \mathsf{Z}(\mathscr{C})$, it satisfies*

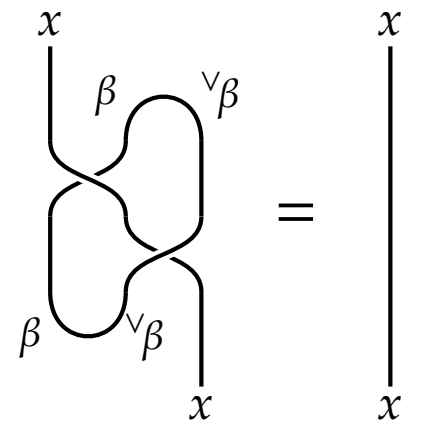

(4.2.5)

*Proof.* Let $B = (\beta, \sigma_{\beta,-}) \in \mathsf{Z}(\mathscr{C})$ be invertible and $x \in \mathsf{Z}(\mathscr{C})$. The left-hand side of Equation (4.2.5) can be rephrased as:

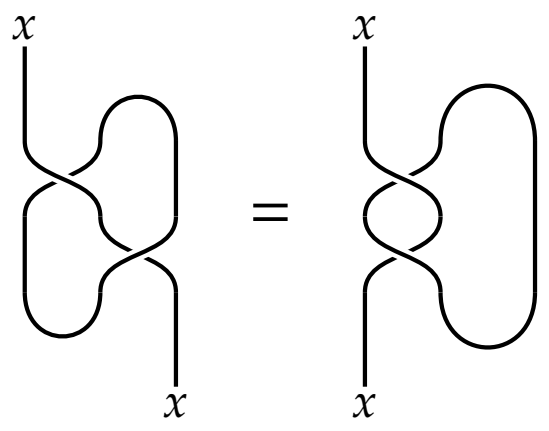

Define the morphism $f := \mathrm{id}_x \otimes \mathrm{coev}^{\ell}_{\beta} \colon x \longrightarrow x \otimes \beta \otimes {}^{\vee}\beta$ and observe that Equation (4.2.5) is identical to

$$f^{-1} \circ ((\sigma_{\beta,x} \circ \sigma_{x,\beta})^{-1} \otimes \mathrm{id}_{{}^{\vee}\beta}) \circ f = \mathrm{id}_x.$$

This is equivalent to $(\sigma_{\beta,x} \circ \sigma_{x,\beta}) \otimes \mathrm{id}_{{}^{\vee}\beta} = \mathrm{id}_{x\otimes\beta} \otimes \mathrm{id}_{{}^{\vee}\beta}$. As the functor $- \otimes {}^{\vee}\beta$ is conservative, the claim follows. □

**Lemma 4.30.** *Let $a := (\alpha, \sigma_{\alpha,-}) \in \mathsf{A}(\mathscr{C})$ be $\mathscr{C}$-invertible, and write $\Omega := (\omega, \sigma_{\omega,-})$ for its dual in $\mathsf{Q}(\mathscr{C})$. For any $x \in \mathsf{Z}(\mathscr{C})$, the inverse of $\rho_{a,x}$ is*

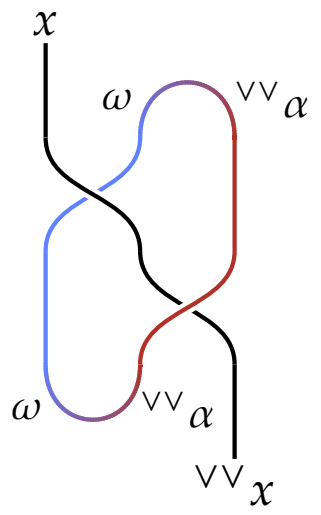

*Proof.* For $x \in \mathsf{Z}(\mathscr{C})$, the snake identities and a variant of Figure 4.2 imply Figure 4.4. Thus, for $\Omega := ({}^{\vee}\alpha, \sigma_{{}^{\vee}\alpha,-}) \in \mathsf{Z}(\mathscr{C})$, we have $\rho_{a,x} \circ \rho_{\Omega,x} = \mathrm{id}_x$. □

**Proposition 4.31.** *Two elements $[a], [c] \in \operatorname{Pic} \mathsf{A}(\mathscr{C})$ induce the same pivotal structure on $\mathsf{Z}(\mathscr{C})$ if and only if there exists a $[b] \in \operatorname{Pic} \mathsf{SZ}(\mathscr{C})$ such that $[b \otimes a] = [c]$.*

*Proof.* Let $[a], [c] \in \operatorname{Pic} \mathsf{A}(\mathscr{C})$. Suppose there exists a $[b] \in \operatorname{Pic} \mathsf{SZ}(\mathscr{C})$ such that $[b \otimes a] = [c]$. For any $x \in \mathsf{Z}(\mathscr{C})$, we compute:

(4.2.6)

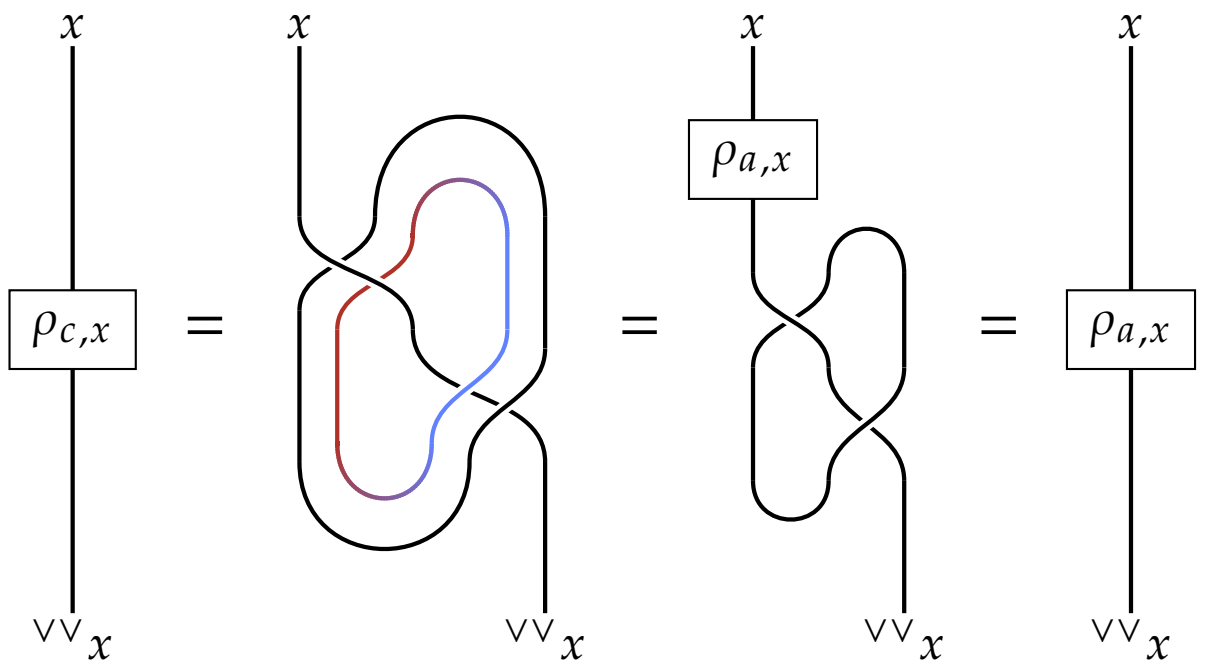

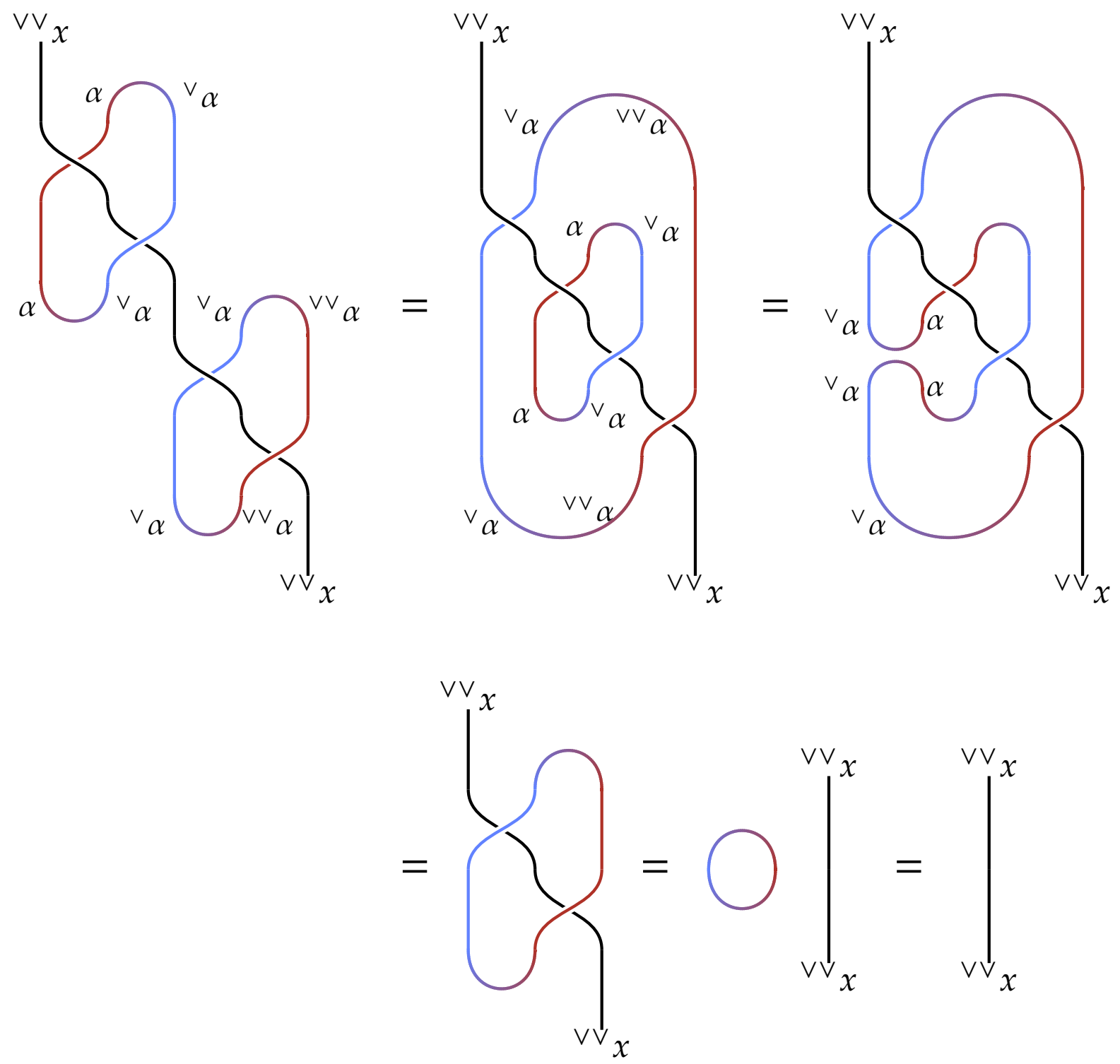

Figure 4.4: Verification of Lemma 4.30.

If conversely $\rho_a = \rho_c$, we claim that $c \otimes {}^\vee a$ is symmetric. By Lemma 4.29 we have to show that for every $x \in \mathsf{Z}(\mathscr{C})$ the "entwinement" $\rho_{c\otimes {}^\vee a}$ of $c \otimes {}^\vee a$ with $x$ is the identity. Indeed,

$$\rho_{c\otimes {}^\vee a,x} = \rho_{{}^\vee a,x} \circ \rho_{c,x} = \rho_{a,x}^{-1} \circ \rho_{c,x} = \mathrm{id}_x.$$

For the first equality we used the hexagon identities as in Equation (4.2.6) to separate $\rho_{c\otimes {}^\vee a,x}$ into two parts. The second one follows from the description of the inverse of $\rho_{a,x}$ given in Lemma 4.30. Finally, since $\mathrm{id}_c \otimes \mathrm{ev}_a^\ell \colon c \otimes {}^\vee a \otimes a \longrightarrow c$ is an isomorphism in $\mathsf{A}(\mathscr{C})$, we have $[(c \otimes {}^\vee a) \otimes a] = [c]$. □

The isomorphism classes of $\mathscr{C}$-invertible objects in $\mathsf{A}(\mathscr{C})$ are not just a set, but form the Picard heap $\mathrm{Pic}\,\mathsf{A}(\mathscr{C})$.

**Lemma 4.32.** *The canonical projection* $\pi\colon \mathrm{Pic}\,\mathsf{A}(\mathscr{C}) \longrightarrow \mathrm{Pic}\,\mathsf{A}(\mathscr{C})/\mathrm{Pic}\,\mathsf{SZ}(\mathscr{C})$ *induces a heap structure on the set of equivalence classes* $\mathrm{Pic}\,\mathsf{A}(\mathscr{C})/\mathrm{Pic}\,\mathsf{SZ}(\mathscr{C})$.

*Proof.* The claim follows from a general observation. Let $x \in \mathsf{Z}(\mathscr{C})$ and $a \in \mathsf{A}(\mathscr{C})$. The half-braiding $\sigma_{x,a} \colon x \otimes a \longrightarrow a \otimes x$ is an isomorphism in $\mathsf{A}(\mathscr{C})$:

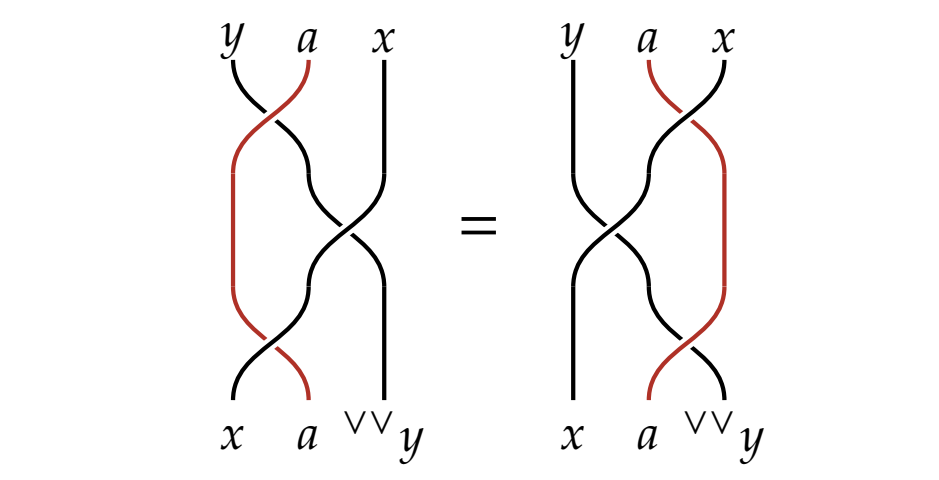

$\sigma_{a\otimes x,y} \circ (\sigma_{x,a} \otimes {}^{\vee\vee}y) = (y \otimes \sigma_{x,a}) \circ \sigma_{x\otimes a,y}$

Likewise, $\sigma_{x,{}^{\vee}a} \colon x \otimes {}^{\vee}a \longrightarrow {}^{\vee}a \otimes x$ is an isomorphism in $\mathsf{Q}(\mathscr{C})$. For all $[a], [a'], [a''] \in \operatorname{Pic} \mathsf{A}(\mathscr{C})$ and $[b], [b'], [b''] \in \operatorname{Pic} \mathsf{SZ}(\mathscr{C})$, one calculates

$$\begin{aligned} \pi\left(\langle [a], [a'], [a''] \rangle\right) &= \pi\left([a \otimes {}^{\vee}a' \otimes a'']\right) = \pi\left([b \otimes {}^{\vee}b' \otimes b'' \otimes a \otimes {}^{\vee}a' \otimes a'']\right) \\ &= \pi\left([b \otimes a \otimes {}^{\vee}(b' \otimes a') \otimes b'' \otimes a'']\right) = \pi\left(\langle [b \otimes a], [b' \otimes a'], [b'' \otimes a''] \rangle\right). \end{aligned} \qquad \square$$

Recall that, by Example 4.7, the pivotal structures on $\mathsf{Z}(\mathscr{C})$ form a heap. Using Definition 4.24, we can relate this structure to that of the anti-centre.

**Theorem 4.33.** *The morphism of heaps*

$$\kappa \colon \operatorname{Pic} \mathsf{A}(\mathscr{C}) \longrightarrow \operatorname{Piv} \mathsf{Z}(\mathscr{C}), \qquad [a] \longmapsto \rho_a$$

*induces a unique injective heap morphism* $\iota \colon \operatorname{Pic} \mathsf{A}(\mathscr{C})/\operatorname{Pic} \mathsf{SZ}(\mathscr{C}) \longrightarrow \operatorname{Piv} \mathsf{Z}(\mathscr{C})$, *such that the following diagram commutes*:

$$\begin{array}{ccc} \operatorname{Pic} \mathsf{A}(\mathscr{C}) & \xrightarrow{\quad\kappa\quad} & \operatorname{Piv} \mathsf{Z}(\mathscr{C}) \\ {\scriptstyle\pi}\downarrow & \nearrow{\scriptstyle \exists!\iota} & \\ \operatorname{Pic} \mathsf{A}(\mathscr{C})/\operatorname{Pic} \mathsf{SZ}(\mathscr{C}) & & \end{array} \tag{4.2.7}$$

*Proof.* Lemmas 4.25 and 4.26 show that $\kappa$ is well-defined. Given three elements $[a], [b], [c] \in \operatorname{Pic} \mathsf{A}(\mathscr{C})$, we compute

$$\kappa(\langle [a], [b], [c] \rangle) = \rho_{a \otimes {}^{\vee}b \otimes c} = \rho_a \circ \rho_{{}^{\vee}b} \rho_c = \rho_a \circ \rho_b^{-1} \circ \rho_c = \langle \rho_a, \rho_b^{-1}, \rho_c \rangle.$$

The second equality follows from the hexagon identities as Equation (4.2.6), and the third one by Lemma 4.30.

Proposition 4.31 states that for any two elements $[a], [b] \in \operatorname{Pic} \mathsf{A}(\mathscr{C})$ we have $\kappa([a]) = \kappa([b])$ if and only if $\pi([a]) = \pi([b])$. It follows from Lemma 4.32 that the unique injective map $\iota\colon \operatorname{Pic} \mathsf{A}(\mathscr{C})/\operatorname{Pic} \mathsf{SZ}(\mathscr{C}) \longrightarrow \operatorname{Piv} \mathsf{Z}(\mathscr{C})$ that lets Diagram (4.2.7) commute is a morphism of heaps. □

**Example 4.34.** By [Shi19, Theorem 1.1], the Picard group $\operatorname{Pic} \mathsf{SZ}(\mathscr{C})$ of the Müger centre of a finite tensor category $\mathscr{C}$ over an algebraically closed field is trivial. In this setting, the induced pivotal structures depend only on the Picard heap $\operatorname{Pic} \mathsf{A}(\mathscr{C})$ and not on a quotient thereof.

On the other side of the spectrum, one might consider the discrete category $\mathbf{E}G$ of an abelian group $G$; its set of objects is $G$ and all morphisms are identities.[12] The category $\mathbf{E}G$ is rigid monoidal; the tensor product given by the multiplication of $G$ and the left and right duals given by the respective inverses. A direct computation shows that $\mathsf{SZ}(\mathbf{E}G) = \mathsf{Z}(\mathbf{E}G) \cong \mathbf{E}G$. Since $\mathbf{E}G$ is skeletal and every object is invertible, $\operatorname{Pic} \mathsf{SZ}(\mathbf{E}G) \cong G$. As the double dual and identity functor on $\mathbf{E}G$ coincide, the same argument implies that $\operatorname{Pic} \mathsf{A}(\mathbf{E}G) \cong G$, whence it follows that $\operatorname{Pic} \mathsf{A}(\mathbf{E}G)/\operatorname{Pic} \mathsf{SZ}(\mathbf{E}G) \cong \{1\}$.

[12] In analogy with the delooping $\mathbf{B}G$ of $G$, the category $\mathbf{E}G$ is also called the *total space* of the group.

In good cases, all pivotal structures on the centre of $\mathscr{C}$ are induced by the quasi-pivotal structures of $\mathscr{C}$.

**Proposition 4.35** ([Shi23a, Theorem 4.1])**.** *For a finite tensor category $\mathscr{C}$, the map*

$$\iota\colon \operatorname{Pic} \mathsf{A}(\mathscr{C}) \big/ \operatorname{Pic} \mathsf{SZ}(\mathscr{C}) \longrightarrow \operatorname{Piv} \mathsf{Z}(\mathscr{C})$$

*is bijective.*

However, in the introduction of [Shi23a] the author states that it is not to be expected that this holds true in general. In the remainder of this section, we will construct an explicit counterexample.

Let us sketch our general approach: suppose there is an object $x \in \mathscr{C}$ that can be endowed with two different half-braidings $\sigma_{x,-}$ and $\chi_{x,-}$. Assume furthermore that there is a pivotal structure $\zeta\colon \mathrm{Id}_{\mathsf{Z}(\mathscr{C})} \longrightarrow {}^{\vee\vee}(-)$ on $\mathsf{Z}(\mathscr{C})$ such that $\zeta_{(x,\sigma_{x,-})} \neq \zeta_{(x,\chi_{x,-})}$. If the unit of $\mathscr{C}$ is the only invertible object, there is no (quasi-)pivotal structure inducing $\zeta$ and therefore $\iota$ cannot be surjective.

Let us define such a category $\mathscr{C}$ in terms of generators and relations. Our approach is similar to Section 3.1.2 and again follows [Kas98, Chapter XII]. Consider the free monoidal category $\mathscr{C}^{\text{free}}$ generated by a single object $x$.

Their tensor product is given by $x^n \otimes x^m := x^{n+m}$. The morphisms of $\mathscr{C}^{\text{free}}$ are identities on objects plus the set $M$ of generating morphisms:

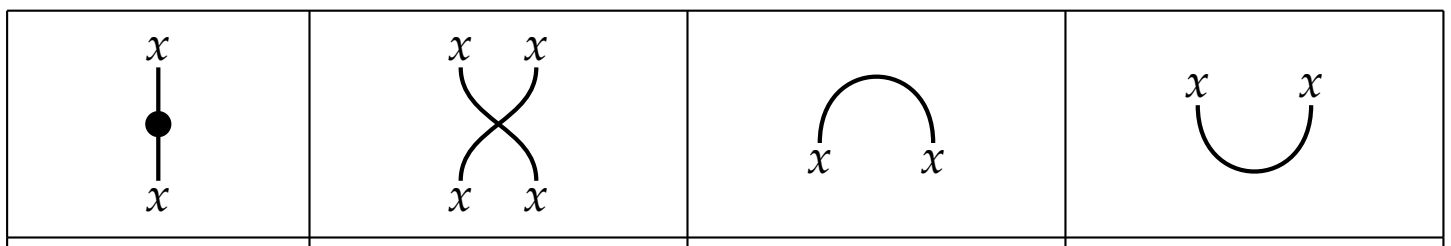

| | | | |
|---|---|---|---|
| $\rho_x : x \longrightarrow x$ | $\sigma_{x,x} : x^2 \longrightarrow x^2$ | $\mathrm{ev}_x : x^2 \longrightarrow 1$ | $\mathrm{coev}_x : 1 \longrightarrow x^2$ |

**Remark 4.36.** By [Kas98, Lemma XII.1.2], every morphism $f : x^n \longrightarrow x^m$ in $\mathscr{C}^{\text{free}}$ is either the identity or can be written as

$$f = (\mathrm{id}_{x^{j_l}} \otimes f_l \otimes \mathrm{id}_{x^{i_l}}) \circ \cdots \circ (\mathrm{id}_{x^{j_2}} \otimes f_2 \otimes \mathrm{id}_{x^{i_2}}) \circ (\mathrm{id}_{x^{j_1}} \otimes f_1 \otimes \mathrm{id}_{x^{i_1}}),$$

where $i_1, j_1, \ldots, i_l, j_l \in \mathbb{N}$ and $f_1, \ldots, f_l \in M$. Such a presentation is not unique, but the number $l \in \mathbb{N}$ of generating morphisms needed to write $f$ in such a manner is. We call it the *degree* of $f$ and write $\deg(f) := l$.

**Definition 4.37.** The category $\mathscr{C}$ is defined as the quotient of $\mathscr{C}^{\text{free}}$ by the relations depicted below. This turns $\mathscr{C}$ into a pivotal (strict) rigid monoidal category, and allows us to extend $\sigma$ to a symmetric braiding. To increase readability, we omit labelling the strings with $x$.

(4.2.8)

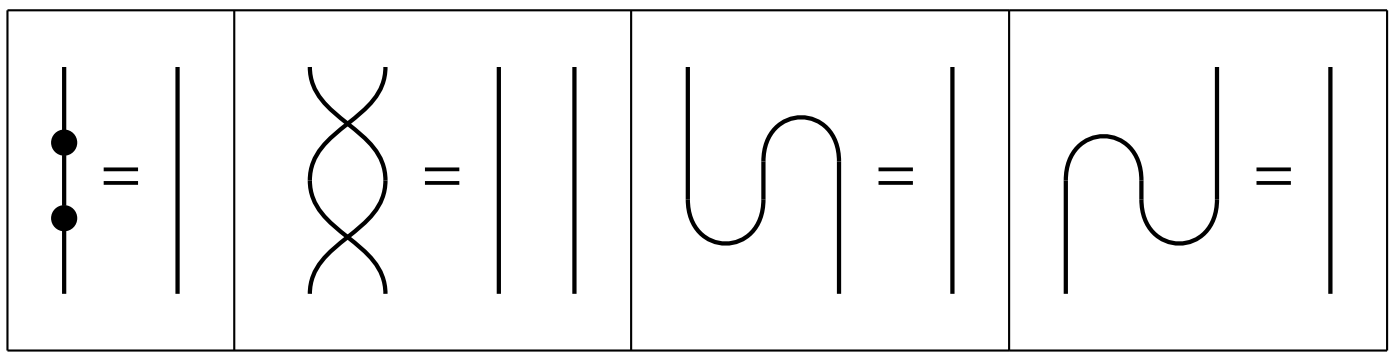

(4.2.9)

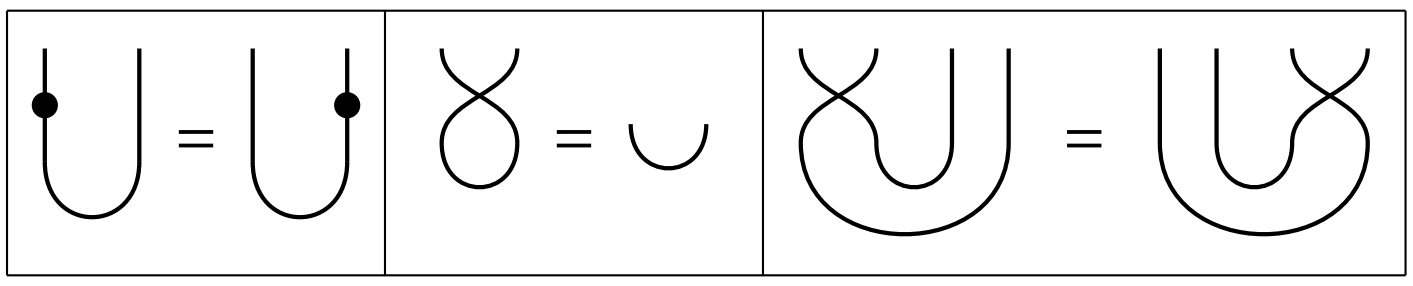

(4.2.10)

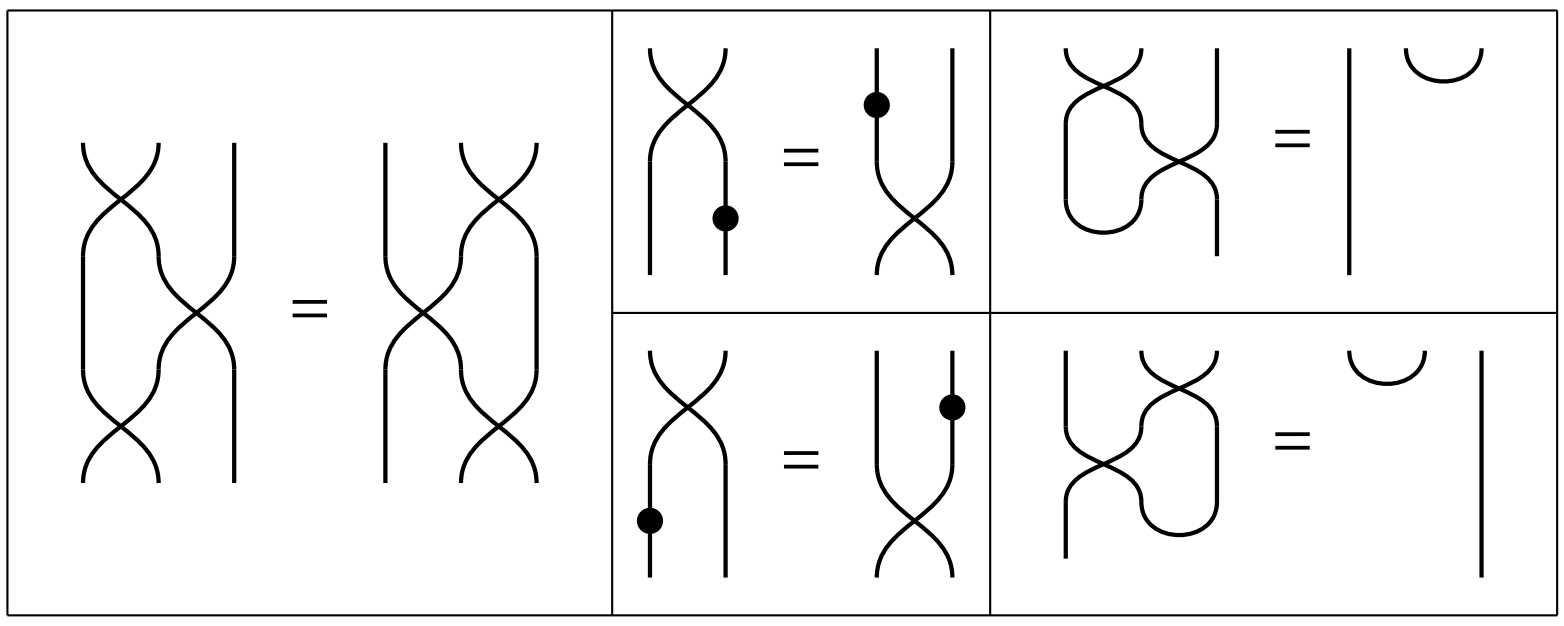

By [Kas98, Proposition XII.1.4] there is a unique functor $P\colon \mathscr{C}^{\mathrm{free}} \longrightarrow \mathscr{C}$ that maps objects to themselves and the generating morphisms to their respective equivalence classes.

**Definition 4.38.** Consider a morphism $f \in \mathscr{C}(x^n, x^m)$. A *presentation* of $f$ is a morphism $g \in \mathscr{C}^{\mathrm{free}}(x^n, y^n)$ such that $f = Pg$. If the degree of $g$ is minimal amongst the presentations of $f$, we call it a *minimal presentation*.

**Theorem 4.39.** *The category $\mathscr{C}$ of Definition 4.37 is strict rigid and the double dual functor is the identity. Further, $\mathrm{id}_x$, $\rho_x\colon x \longrightarrow x$ can be extended to pivotal structures, and $\sigma_{x,x}\colon x^2 \longrightarrow x^2$ may be extended to a symmetric braiding.*

*Proof.* The evaluation and coevaluation morphisms plus their snake identities make $x \in \mathscr{C}$, and by extension every object of $\mathscr{C}$, its own left and right dual. Using Equation (4.2.9), we compute

$$
{}^{\vee}\rho_x = \rho_x = \rho_x{}^{\vee}, \qquad\qquad {}^{\vee}\sigma_{x,x} = \sigma_{x,x} = \sigma_{x,x}{}^{\vee},
$$
$$
{}^{\vee}\mathrm{ev}_x = \mathrm{coev}_x = \mathrm{ev}_x{}^{\vee}, \qquad\qquad {}^{\vee}\mathrm{coev}_x = \mathrm{ev}_x = \mathrm{coev}_x{}^{\vee}.
$$

Hence, $\mathscr{C}$ is strict rigid and its double dual functor is equal to the identity.

Our candidate for a non-trivial pivotal structure on $\mathscr{C}$ is

$$
\rho\colon \mathrm{Id}_{\mathscr{C}} \longrightarrow \mathrm{Id}_{\mathscr{C}}, \qquad \rho_{x^n} := \rho_x \otimes \cdots \otimes \rho_x\colon x^n \longrightarrow x^n, \text{ for } n \in \mathbb{N}.
$$

By construction, this family of isomorphisms is compatible with the monoidal structure of $\mathscr{C}$, so we only have to prove naturality, for which it suffices to check the generators. Equation (4.2.10) implies that $\rho_{x^2}$ commutes with $\sigma_{x,x}$. For the evaluation of $x \in \mathscr{C}$ we use the dual of Equation (4.2.9) to compute

$$
\mathrm{ev}_x \circ \rho_{x^2} = \mathrm{ev}_x \circ (\rho_x \otimes \rho_x) = \mathrm{ev}_x \circ ({}^{\vee}\rho_x \otimes \rho_x) = \mathrm{ev}_x \circ (\mathrm{id}_x \otimes \rho_x^2) = \rho_1 \circ \mathrm{ev}_x.
$$

Applying the left dualising functor, one obtains $\mathrm{coev}_x \circ \rho_1 = \rho_{x^2} \circ \mathrm{coev}_x$, whence $\rho$ defines a pivotal structure.

Lastly, $\sigma_{x,x}$ implements a symmetry $\sigma\colon \otimes \Longrightarrow \otimes^{\mathrm{op}}$ on $\mathscr{C}$. Set

$$
\sigma_{x,x^m} := (\mathrm{id}_x \otimes \sigma_{x,x^{m-1}}) \circ (\sigma_{x,x} \otimes \mathrm{id}_{x^{m-1}}), \quad \text{for } m \in \mathbb{N},
$$

and extend this to arbitrary objects by

$$
\sigma_{x^n,x^m} := (\sigma_{x^{n-1,m}} \otimes \mathrm{id}_x) \circ (\mathrm{id}_{x^{n-1}} \otimes \sigma_{x,x^m}), \quad \text{for } n, m \in \mathbb{N}.
$$

As this family of isomorphisms is constructed according to the hexagon axioms, we only have to prove its naturality. Again, it suffices to consider the generating morphisms. Equation (4.2.10) implies that $\sigma$ is natural with respect to $\rho_x$, $\sigma_{x,x}$, and $\mathrm{coev}_x$. The self-duality of $\sigma_{x,x}$ and ${}^{\vee}\mathrm{coev}_x = \mathrm{ev}_x$ imply the desired commutation between $\sigma$ and $\mathrm{ev}_x$. Thus $\sigma$ is a braiding on $\mathscr{C}$, which is symmetric by Equation (4.2.8). □

We think of a generic morphism of $\mathscr{C}$ to be of the form

(4.2.11)

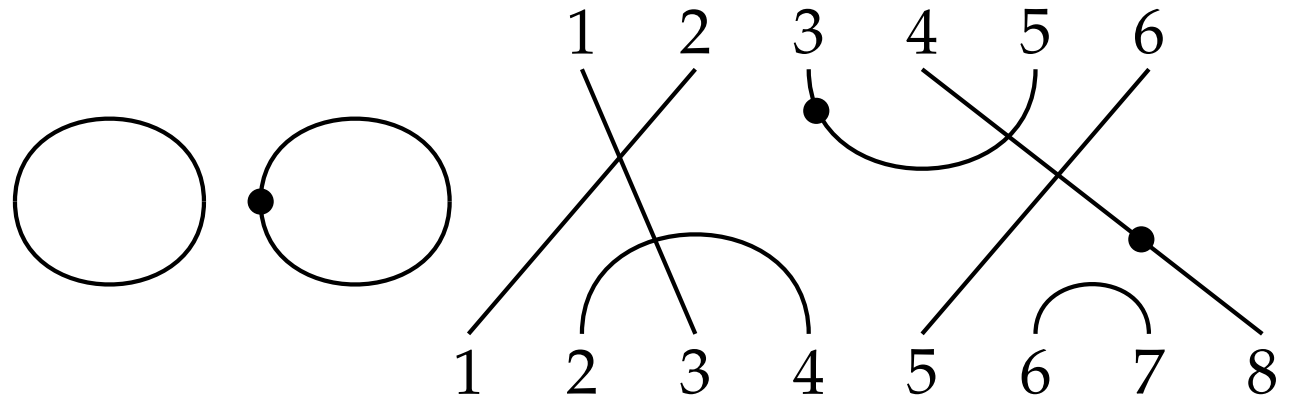

This diagram suggests a distinction between different kinds of morphisms: there are *connectors*, which link an input to an output vertex, *closed loops*, and *half-circles* of *evaluation* and *coevaluation-type*. Connectors induce a permutation on a subset of $\mathbb{N}$. For example, the permutation arising from Equation (4.2.11) can be identified with (1 2)(3 4).

As usual, $S_n$ denotes the symmetric group on $\{1, 2, \dots, n\}$.

Conversely, let $s := t_{i_1} \dots t_{i_l} \in S_n$ be a permutation written as a product of elementary transpositions and set $f_s := f_{t_{i_1}} \dots f_{t_{i_l}} : x^n \longrightarrow x^n$, where

$$f_{t_i} := \mathrm{id}_{x^{i-1}} \otimes \sigma_{x,x} \otimes \mathrm{id}_{x^{n-(i+1)}} : x^n \longrightarrow x^n, \qquad \text{for} \qquad 1 \leq i \leq n-1.$$

As the braiding $\sigma$ is symmetric, $f_s$ does not depend on the presentation of $s$. If the presentation of $s$ is minimal, however, so is the corresponding one of $f_s$.

To derive a normal form of the automorphisms of $\mathscr{C}$ and turn our previously explained thoughts into precise mathematical statements, we need to study the "topological" features of the morphisms in $\mathscr{C}$.

**Remark 4.40.** We recall the *category of tangles* $\mathfrak{T}$, a close relative to the string diagrams arising from $\mathscr{C}$, based on [Kas98, Chapter XII.2]. Its objects are finite sequences in $\{+, -\}$ and its morphisms are isotopy classes of oriented tangles. A detailed discussion of tangles is given in [Kas98, Definition X.5.1]. For us, it suffices to think of an oriented tangle $L$ of type $(n, m)$ as a finite disjoint union of embeddings of either the unit circle $S^1$ or the interval $[0, 1]$ into $\mathbb{R}^2 \times [0, 1]$, such that

$$\partial L = L \cap \left(\mathbb{R}^2 \times \{0, 1\}\right) = ([n] \times \{(0, 0)\}) \cup ([l] \times \{(0, 1)\}),$$

where $[n] = \{1, \ldots, n\}$ and $[l] = \{1, \ldots, l\}$. The orientation on each of the connected components of $L$ is induced by the counter-clockwise orientation of $S^1$ and the (ascending) orientation of $[0,1]$. The tensor product of tangles is given by pasting them next to each other. Their composition is implemented, by appropriate gluing and rescaling.

To distinguish isotopy classes of tangles, one can study their images under the projection $\mathbb{R}^2 \times [0,1] \twoheadrightarrow \mathbb{R} \times [0,1]$. This leads to a combinatorial description of $\mathfrak{T}$, see for example [Kas98, Theorem XII.2.2].

**Proposition 4.41.** *The strict monoidal category $\mathfrak{T}$ is generated by the morphisms*

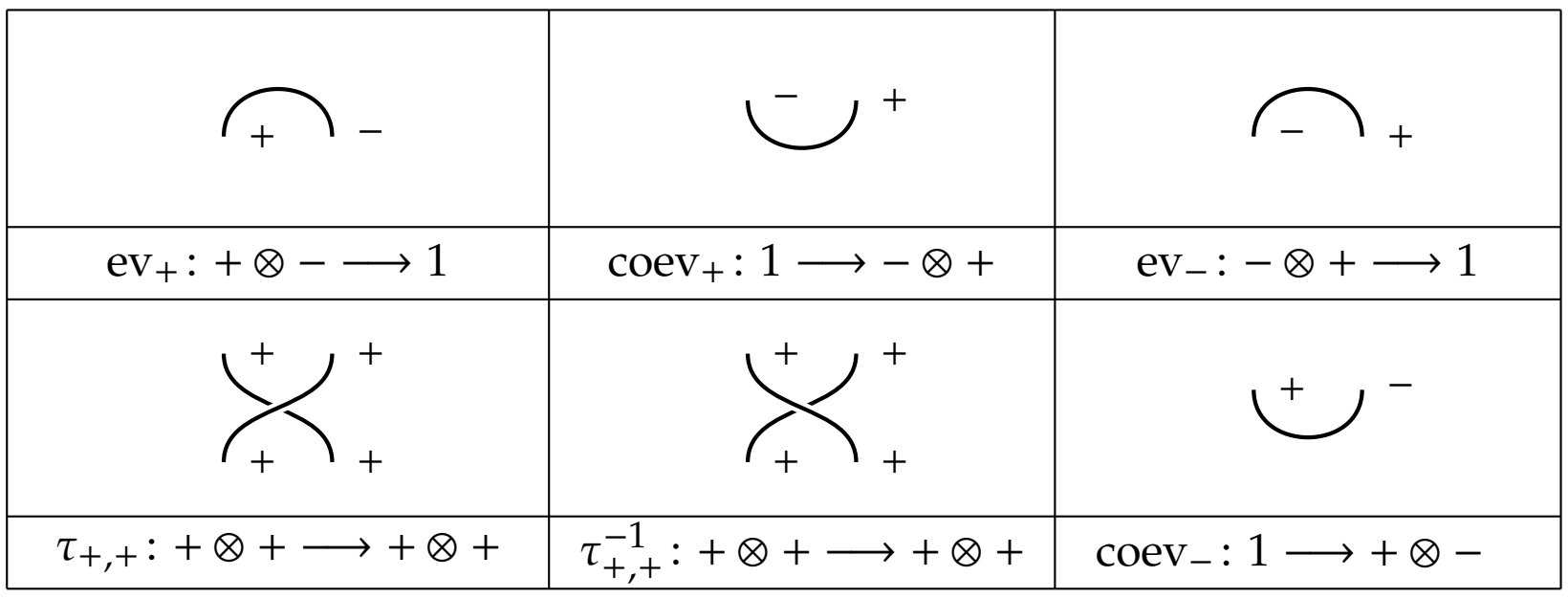

*They are subject to the following relations:*

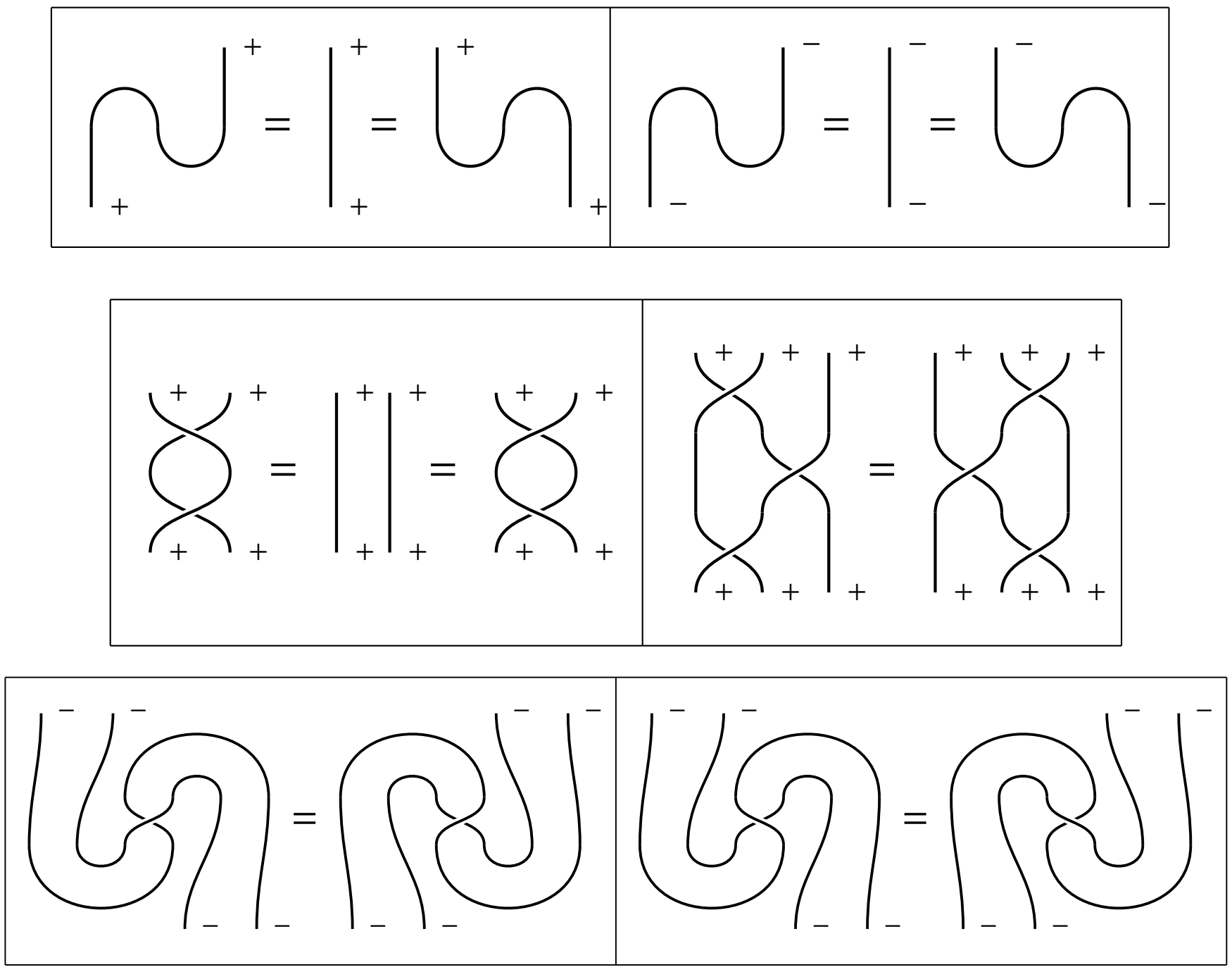

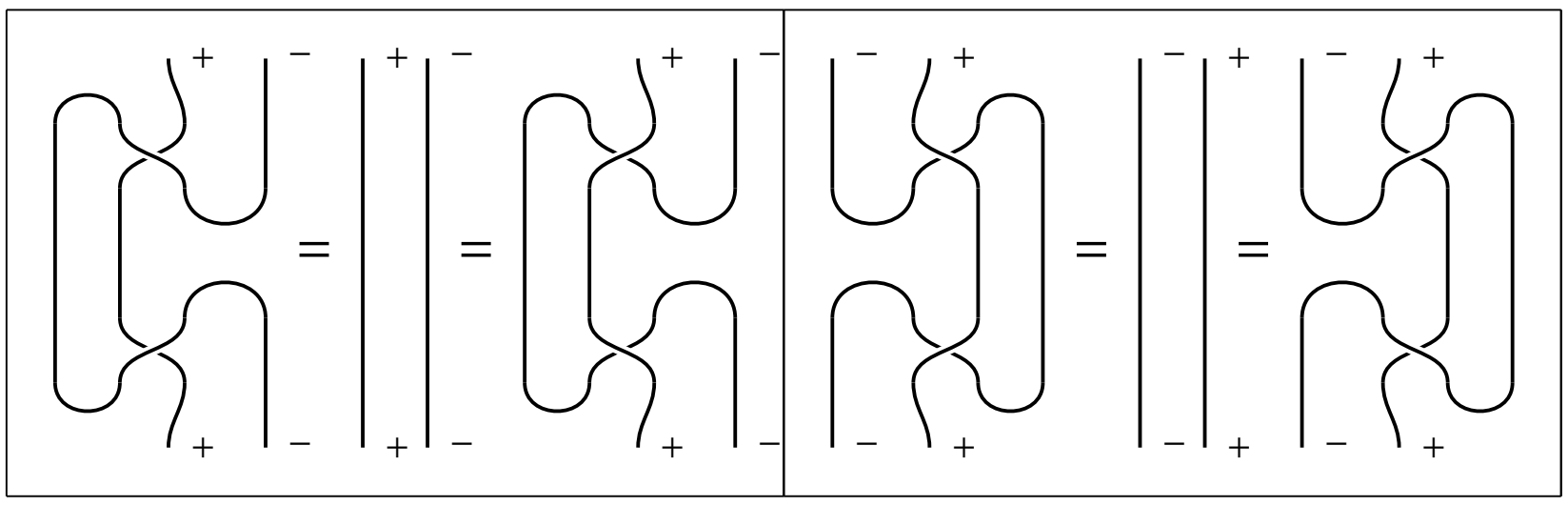

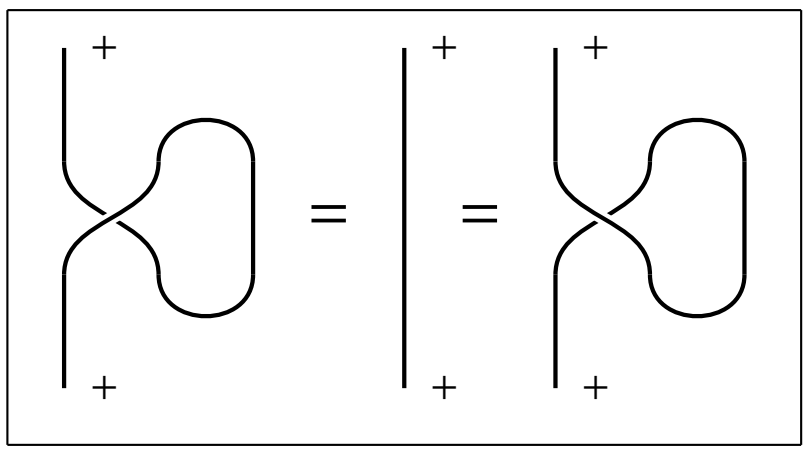

The connection between tangles and the category $\mathscr{C}$ is attained by applying [Kas98, Proposition xii.1.4] in a similar way as we did in Lemma 3.21.

**Lemma 4.42.** *There exists a strict monoidal functor $S\colon \mathfrak{T} \longrightarrow \mathscr{C}$ that is uniquely determined by $S(+) = x = S(-)$ and*

$$S(\mathrm{ev}_\pm) = \mathrm{ev}_x, \qquad S(\mathrm{coev}_\pm) = \mathrm{coev}_x, \qquad S(\tau^{\pm}_{+,+}) = \sigma_{x,x}.$$

We now want to lift the morphisms of $\mathscr{C}$ to $\mathfrak{T}$. Hereto we want to "trivialise" the generator $\rho_{x,x}\colon x \longrightarrow x$. Set $\mathscr{C}/\langle\rho_x\rangle$ to be the category obtained from $\mathscr{C}$ by identifying $\rho_x$ with $\mathrm{id}_x$. The projection functor $\mathrm{Pr}\colon \mathscr{C} \longrightarrow \mathscr{C}/\langle\rho_x\rangle$ allows us to define an equivalence relation on the morphisms of $\mathscr{C}$:

$$f \sim g \qquad \Longleftrightarrow \qquad \mathrm{Pr}(f) = \mathrm{Pr}(g).$$

For example the endomorphisms ○: $1 \longrightarrow 1$ and ◐: $1 \longrightarrow 1$ of the monoidal unit of $\mathscr{C}$ would be equivalent with respect to this relation.

**Proposition 4.43.** *Every $f \in \mathrm{Aut}_{\mathscr{C}}(x^n, x^n)$ can be uniquely written as $f = f_s \circ f_\phi$, where $f_s\colon x^n \longrightarrow x^n$ is the automorphism induced by a permutation $s \in S_n$ and $f_\phi = \rho_x^{\phi_1} \otimes \cdots \otimes \rho_x^{\phi_n}$, with $\phi_1, \ldots, \phi_n \in \mathbb{Z}_2$. Furthermore, if a minimal presentation $s := t_{i_1} \ldots t_{i_l}$ is fixed, the resulting presentation of $f$ is minimal as well.*

*Proof.* For any $f \in \mathrm{Aut}_{\mathscr{C}}(x^n)$ there exists another automorphism $g \in \mathrm{Aut}_{\mathscr{C}}(x^n)$, such that $\mathrm{Pr}\, f = \mathrm{Pr}\, g$ and $g$ has a presentation in which no copies of $\rho$ occur.

By proceeding analogous to [Kas98, Lemma x.3.3], we construct a tangle $L_g$ out of $g$ such that $S(L_g) = g$, and it is isotopic to a tangle $L'_g$ whose images of its connected components under the projection $\mathbb{R}^2 \times [0,1] \longrightarrow \mathbb{R} \times [0,1]$ are either closed loops, half-circles of evaluation or coevaluation-type, or straight lines. Write $L_n^{\mathrm{triv}}$ for a tangle that projects to $n$ parallel straight lines:

$$\{ (k,t) \mid t \in [0,1] \text{ and } k \in \{1,\dots,n\} \}.$$

Since $g$ is invertible by assumption, we can lift its inverse $g^{-1}$ to a tangle $L_{g^{-1}}$ with $[L_g] \circ [L_{g^{-1}}] = [L_n^{\mathrm{triv}}] = [L_{g^{-1}}] \circ [L_g]$. This equation implies that $L'_g$ could not have contained any loops or half-circles. In other words $g = f_s$, where $f_s$ is the morphism obtained from the permutation $s \in S_n$, induced by the projection of $L'_g$ onto $\mathbb{R} \times [0,1]$. Due to the naturality of $\sigma_{x,x}$, the equivalence between $f$ and $g$ implies $f = f_s \circ f_\phi$, with $f_\phi$ being a tensor product of identities and copies of $\rho_x$. Consequentially, a minimal representation of $s$ induces a minimal representation of $f$. ◻

The first step in showing that the $\iota$ defined in Theorem 4.33 cannot be surjective is to prove that the $\operatorname{Pic} \mathsf{A}(\mathscr{C})$ contains at most two elements.

**Corollary 4.44.** *The only quasi-pivotal structures on the category $\mathscr{C}$ of Definition 4.37 are* $\mathrm{id}\colon \mathrm{Id}_{\mathscr{C}} \longrightarrow \mathrm{Id}_{\mathscr{C}}$ *and* $\rho\colon \mathrm{Id}_{\mathscr{C}} \longrightarrow \mathrm{Id}_{\mathscr{C}}$.

*Proof.* The only invertible object of $\mathscr{C}$ is its monoidal unit, which implies that any quasi-pivotal structure on $\mathscr{C}$ is pivotal. By Proposition 4.43, these are determined by their value on $x$ and $\mathrm{Aut}_{\mathscr{C}}(x) = \{ \mathrm{id}_x, \rho_x \}$. ◻

Let us now focus on the various ways in which we can equip an object $y \in \mathscr{C}$ with a half-braiding. Our classification of automorphisms in $\mathscr{C}$ allows us to easily verify that on $x \in \mathscr{C}$ there are four different half-braidings, which are determined by

| 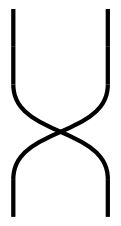 | 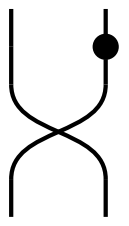 | 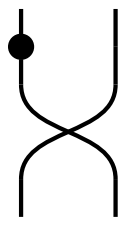 | 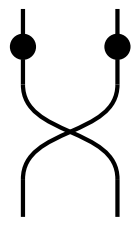 |
| --- | --- | --- | --- |
| $\sigma^{\circ,\circ}_{x,x}\colon x^2 \longrightarrow x^2$ | $\sigma^{\circ,\bullet}_{x,x}\colon x^2 \longrightarrow x^2$ | $\sigma^{\bullet,\circ}_{x,x}\colon x^2 \longrightarrow x^2$ | $\sigma^{\bullet,\bullet}_{x,x}\colon x^2 \longrightarrow x^2$ |

The fact that these braidings are distinguished by the appearances of $\rho$ on the respective strings motivates our next definition.

**Definition 4.45.** Let $f \coloneqq f_s f_\phi\colon x^n \longrightarrow x^n$ be an automorphism in $\mathscr{C}$. Its *characteristic sequence* is $\phi \coloneqq (\phi_1, \dots, \phi_n) \in (\mathbb{Z}_2)^n$ with

$$f_\phi = \rho_x^{\phi_1} \otimes \cdots \otimes \rho_x^{\phi_n}.$$

Indeed, it is the interplay between instances of $\rho$ and the underlying permutation that determine whether an automorphism $\chi_{y,x}\colon y\otimes x \longrightarrow x\otimes y$ can be lifted to a half-braiding.

**Lemma 4.46.** *Any automorphism $\chi_{y,x}\colon y \otimes x \longrightarrow x \otimes y$ extends to a half-braiding on $y$ if and only if there exists an $f \in \mathrm{Aut}_{\mathscr{C}}(y)$ with characteristic sequence $(\phi_1, \dots, \phi_n)$ and underlying permutation $s \in S_n$, such that $s^2(i) = i$ and $\phi_{s(i)} = \phi_i$ for all $1 \le i \le n$, and $\chi_{y,x} = \sigma_{y,x} \circ (f \otimes \rho_x^j)$ for an integer $j \in \mathbb{Z}_2$.*

*Proof.* Assume $\chi_{y,x}\colon y\otimes x \longrightarrow x\otimes y$ to induce a half-braiding on $y := x^n$. Due to Proposition 4.43, we can write $\chi_{y,x} = \sigma_{y,x} \circ (f \otimes \rho_x^j)$, where $f\colon y \longrightarrow y$ is an automorphism of $y$ and $j \in \mathbb{Z}_2$. Let $\phi = (\phi_1, \dots, \phi_n)$ be the characteristic sequence of $f$ and $s \in S_n$ its underlying permutation. Write $f_s\colon y \longrightarrow y$ for the morphism induced by $s$ and set

$$f_\phi = \rho_x^{\phi_1} \otimes \cdots \otimes \rho_x^{\phi_n}, \qquad f_{s^{-1}(\phi)} = \rho_x^{\phi_{s^{-1}(1)}} \otimes \cdots \otimes \rho_x^{\phi_{s^{-1}(n)}}.$$

Using that $f = f_s \circ f_\phi$, the naturality of $\chi_{y,-}$, and Equation (4.2.9), we compute

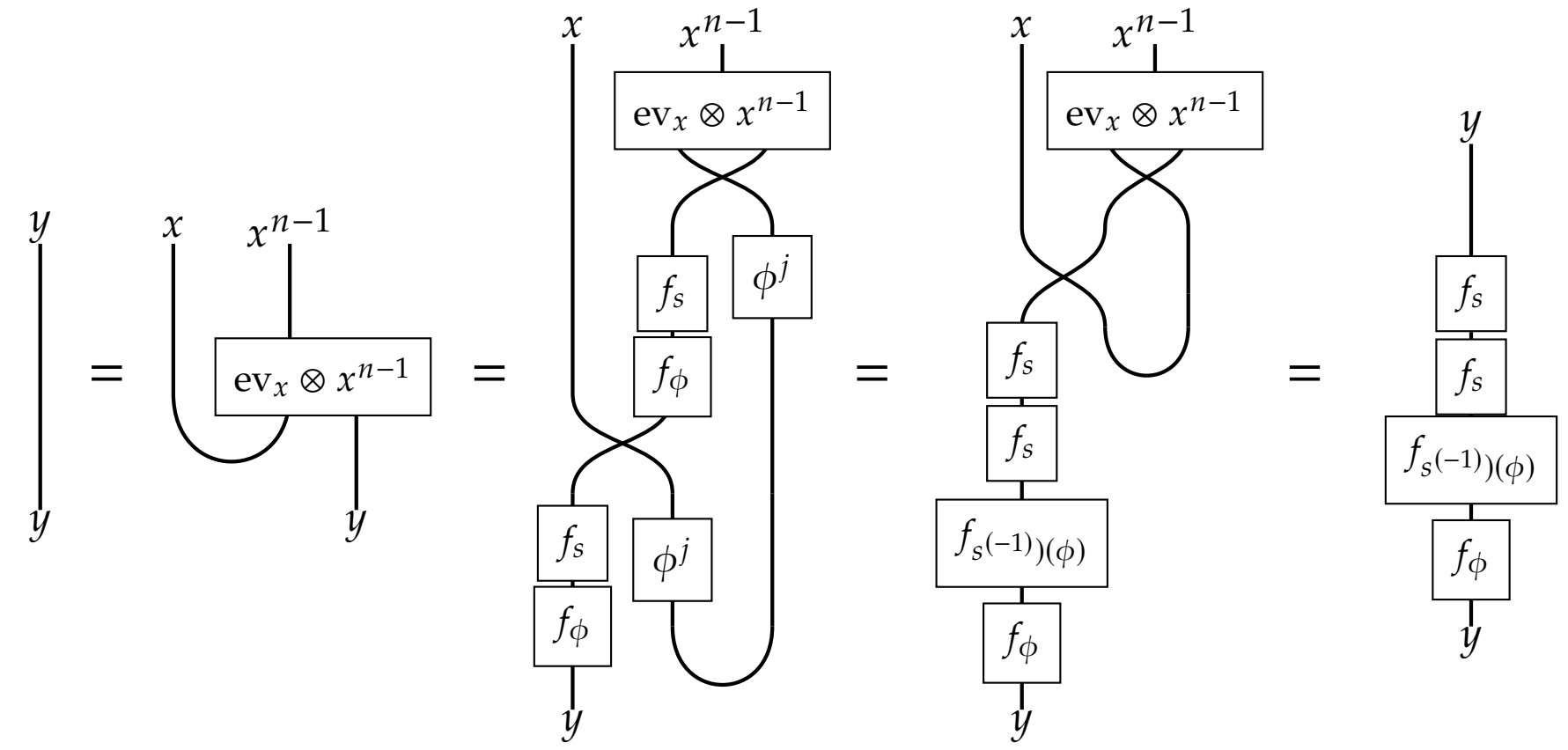

This is equivalent to $s$ being an involution and $\phi$ being invariant under $s$.

Conversely, let $\chi_{y,x} = \sigma_{y,x} \circ (f \otimes \rho_x^j)\colon y \otimes x \longrightarrow x \otimes y$, where $f$ is an automorphism satisfying the assumptions of the lemma. We extend it to a family of automorphisms $\chi_{y,-}\colon y \otimes - \longrightarrow - \otimes y$ according to the hexagon axioms and verify its naturality on the generators of $\mathscr{C}$. For $\rho_x$ and $\sigma_{x,x}$ this is immediate consequence of their respective naturality conditions. To prove the commutation relations between $\chi_{y,-}$, $\mathrm{coev}_x$, and $\mathrm{ev}_x$, argue as above. □

The previous lemma severely restricts the number of possibilities in which an automorphism of $\mathscr{C}$ can lift to the centre $\mathsf{Z}(\mathscr{C})$.

**Corollary 4.47.** *Consider an object $x^n \in \mathscr{C}$ equipped with two half-braidings*

$$\chi_{x^n,x} = \sigma_{x^n,x} \circ ((f_s \circ f_\phi) \otimes \rho_x^j), \qquad \theta_{x^n,x} = \sigma_{x^n,x} \circ ((f_t \circ f_\psi) \otimes \rho_x^k).$$

*If $g = g_r \circ g_\lambda \in \mathrm{Aut}_{\mathscr{C}}(x^n)$ lifts to a morphism $g\colon (x^n, \chi_{x^n,-}) \longrightarrow (x^n, \theta_{x^n,-})$ of objects in the centre of $\mathscr{C}$, then*

$$\phi_i \circ \lambda_{sr(i)} = \psi_{r(i)} \circ \lambda_{r(i)} \qquad \text{for all } 1 \leq i \leq n.$$

*Proof.* For $g = f_r \circ f_\lambda \in \mathrm{Aut}_{\mathscr{C}}(x^n)$ to lift to the centre it must satisfy

$$\begin{aligned} \sigma_{x^n,x} \circ ((f_s \circ f_\phi \circ g) \otimes \rho_x^j) = \chi_{x^n,x} \circ (g \otimes \mathrm{id}_x) &= (\mathrm{id}_x \otimes g) \circ \theta_{x^n,x} \\ &= \sigma_{x^n,x} \circ ((g \circ f_t \circ f_\psi) \otimes \rho_x^k). \end{aligned}$$

This implies that $f_s \circ f_\phi \circ g = g \circ f_t \circ f_\psi$, and therefore $\phi_{s(i)} \circ \lambda_{sr(i)} = \lambda_{r(i)} \circ \psi_{rt(i)}$ for all $1 \leq i \leq n$. Since $\mathbb{Z}_2$ is abelian and $\phi_{s(i)} = \phi_i$ as well as $\psi_{t(i)} = \psi_i$, the claim follows. □

In view of Lemma 4.46, we state a refined version of Definition 4.45.

**Definition 4.48.** Consider an object $y = (x^n, \chi_{x^n,x}) \in \mathsf{Z}(\mathscr{C})$ whose half-braiding is defined by $\chi_{x^n,x} = \sigma_{x^n,x} \circ (f \otimes \rho_x^j)$ for an integer $j \in \mathbb{Z}_2$. We call the characteristic sequence $\phi$ of $f$ the *signature* of $y$.

We now construct a pivotal structure on the centre of $\mathscr{C}$ that differs from the lifts of id and $\rho$ from $\mathscr{C}$ to $\mathsf{Z}(\mathscr{C})$.

**Theorem 4.49.** *The Drinfeld centre $\mathsf{Z}(\mathscr{C})$ of $\mathscr{C}$ admits a pivotal structure $\zeta$ with*

$$\zeta_{(x,\sigma_{x,-}^{\circ,\circ})} = \mathrm{id}_x, \qquad \zeta_{(x,\sigma_{x,-}^{\circ,\bullet})} = \mathrm{id}_x, \qquad \zeta_{(x,\sigma_{x,-}^{\bullet,\circ})} = \rho_x, \qquad \zeta_{(x,\sigma_{x,-}^{\bullet,\bullet})} = \rho_x.$$

*Proof.* For any object $y \in \mathsf{Z}(\mathscr{C})$, define

$$\zeta_y = \rho_x^{\phi_1} \otimes \cdots \otimes \rho_x^{\phi_n}, \qquad \text{where } \phi = (\phi_1, \ldots, \phi_n) \text{ is the signature of } y.$$

Since the signature of a tensor product of objects in $\mathsf{Z}(\mathscr{C})$ is given by concatenating the signatures of its components, this defines a family of isomorphisms $\zeta\colon \mathrm{Id}_{\mathsf{Z}(\mathscr{C})} \longrightarrow \mathrm{Id}_{\mathsf{Z}(\mathscr{C})}$ that is compatible with the monoidal structure.

It remains to prove the naturality of $\zeta$, which can be verified by considering all possible lifts of identities and generators of $\mathscr{C}$ to its Drinfeld centre. For $\mathrm{id}_x, \rho_x\colon x \longrightarrow x$ and $\sigma_{x,x}\colon x^2 \longrightarrow x^2$, this follows from Corollary 4.47. To study the coevaluation of $x$, fix a half-braiding $\chi_{x^2,-}\colon x^2 \otimes - \longrightarrow - \otimes x^2$ on $x^2$. Due to Lemma 4.46, it is determined by

$$\chi_{x^2,x} = \sigma_{x^2,x} \circ \left((\sigma_{x,x}^i \circ (\rho_x^j \otimes \rho_x^k)) \otimes \rho_x^l\right), \qquad \text{where } i,j,k,l \in \mathbb{Z}_2.$$

Now suppose $\mathrm{coev}_x\colon 1 \longrightarrow x^2$ lifts to a morphism in $\mathsf{Z}(\mathscr{C})$, where $x^2$ is equipped with this half-braiding. Equation (4.2.9) and the self-duality of $\sigma_{x,x}$ imply that $\sigma_{x,x} \circ \mathrm{coev}_x = \mathrm{coev}_x$ and $\mathrm{ev}_x \circ \sigma_{x,x} = \mathrm{ev}_x$; we compute

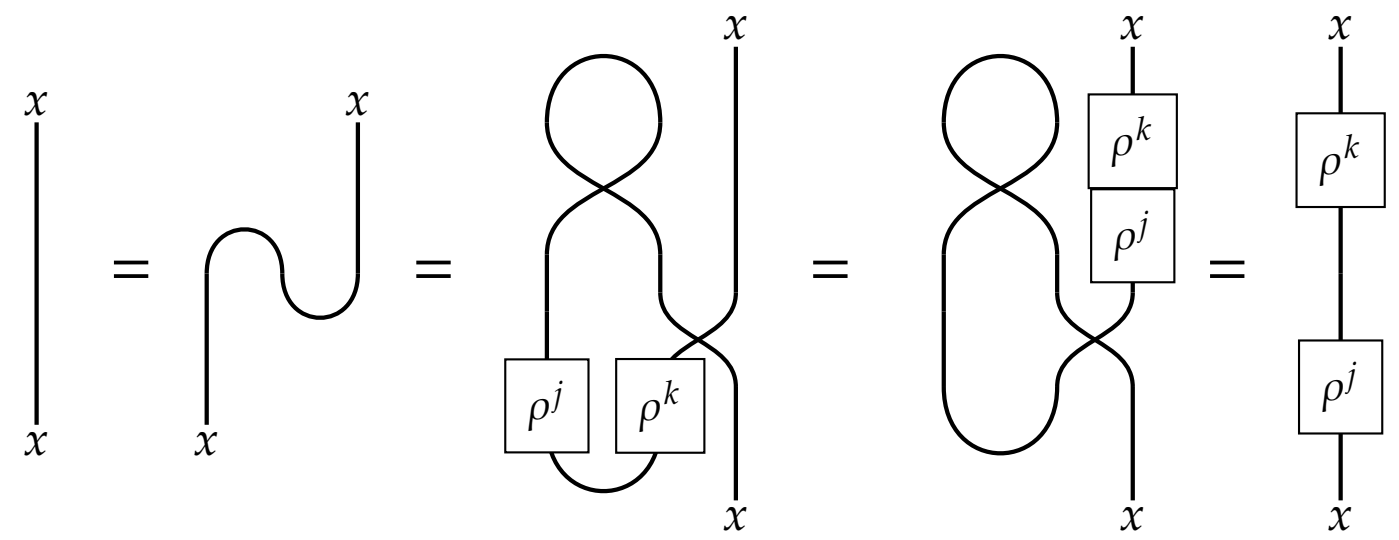

Therefore, $j = k$ and $\zeta_{(x^2,\chi_{x^2,-})} = \mathrm{id}_x^2$ or $\zeta_{(x^2,\chi_{x^2,-})} = \rho_x^2$, implying naturality. A similar argument for the evaluation of $x$ concludes the proof. □

By Corollary 4.44, the Picard heap of $\mathsf{A}(\mathscr{C})$ can have at most two elements. However, the above theorem constructs a third pivotal structure on $\mathsf{Z}(\mathscr{C})$; this implies our desired result.

**Theorem 4.50.** *The pivotal structure $\zeta$ of $\mathsf{Z}(\mathscr{C})$ is not induced by the Picard heap of $\mathsf{A}(\mathscr{C})$. In particular, the map $\iota\colon \operatorname{Pic}\mathsf{A}(\mathscr{C})/\operatorname{Pic}\mathsf{SZ}(\mathscr{C}) \longrightarrow \operatorname{Piv}\mathsf{Z}(\mathscr{C})$ of Theorem 4.33 is not surjective.*

[D]ie Monade, das
unvollkommenste aller Wesen.

Arthur Schopenhauer;
Die Kunst, Recht zu behalten

# 5 MONADIC TANNAKA–KREIN RECONSTRUCTION

Bimonads and Hopf monads are a vast generalisation of bialgebras and Hopf algebras. They naturally arise in the study of (rigid) monoidal categories and topological quantum field theories, see amongst others [KL01; Moe02; BV07; BLV11; TV17]. Recall that the definition of a bialgebra object necessarily requires a braided monoidal category as a base, in order to write down what it means for the multiplication to be a morphism of comonoids. However, in general the category of endofunctors is not braided, so the naïve notion of bialgebras does not generalise to the monadic setting and needs to be adjusted. One possible way of overcoming this problem was introduced and studied by Moerdijk under the name *Hopf monad*, [Moe02].[13] Here, the structure morphisms of an oplax monoidal functor serve as a substitute of the comultiplication and counit. There are other, sometimes non-equivalent, notions of Hopf monad, see [Boa95; MW11]. We follow [Moe02], with a slight terminology change due to [BV07; BLV11]. This definition aims to generalise the paradigmatic example of the free–forgetful adjunction of a Hopf algebra. Thus, equipping a monad with the prefixes "bi-" or "Hopf" refers to additional structure or properties put on its Eilenberg–Moore category.

[13] As remarked in [Moe02], this concept is strictly dual to that of monoidal comonads, which are studied in [Boa95].

More generally, a monadic interpretation of module categories was given by Aguiar and Chase under the name *comodule monad*, [AC12]. A comodule monad over a bimonad generalises the notion of a comodule algebra over a bialgebra, see Example 5.14 below. The main result of this chapter extends the reconstruction results of [AC12, Proposition 4.1.2] and [TV17, Lemma 7.10]:

**Theorem 5.28.** *Let $\mathscr{C}$ and $\mathscr{D}$ be monoidal categories, and suppose that $\mathscr{M}$ and $\mathscr{N}$ are right $\mathscr{C}$- and $\mathscr{D}$-module categories, respectively. Let $F\colon \mathscr{C} \rightleftarrows \mathscr{D} : U$ be an oplax monoidal adjunction. Lifts of an adjunction $G\colon \mathscr{M} \rightleftarrows \mathscr{N} : V$ to a comodule adjunction are in bijection with lifts of $V\colon \mathscr{N} \longrightarrow \mathscr{M}$ to a strong comodule functor.*

From the proof of this result one immediately obtains an analogue of Kelly's doctrinal adjunction result, [Kel74], for $\mathscr{C}$-module categories.

**Porism 5.29.** *An adjunction $G\colon \mathcal{M} \rightleftarrows \mathcal{N} : V$ of $\mathcal{C}$-module categories yields a bijection between oplax $\mathcal{C}$-module structures on $F$ and lax $\mathcal{C}$-module structures on $U$.*

We then obtain a Tannaka–Krein reconstruction result for comodule monads in the spirit of [Moe02, Theorem 7.1] and [McC02, Corollary 3.13].

**Theorem 5.31.** *Let $B$ be a bimonad on the monoidal category $\mathcal{C}$, and $T$ a monad on a right $\mathcal{C}$-module category $\mathcal{M}$. Coactions of $B$ on $T$ are in bijection with right actions of $\mathcal{C}^B$ on $\mathcal{M}^T$, such that $U^T$ is a strict comodule functor over $U^B$.*

We can subsequently apply these results in order to study the comparison functor of Section 2.2.1. Namely, this functor always inherits the strong module structure of the "strong adjoint"; i.e., the right adjoint for oplax $\mathcal{C}$-module adjunctions, and the left one in the lax case.

**Proposition 5.35.** *Let $\mathcal{M}$ and $\mathcal{N}$ be left $\mathcal{C}$-module categories.*

- *The comparison functor $K^{VG}$ of an oplax $\mathcal{C}$-module adjunction $G\colon \mathcal{M} \rightleftarrows \mathcal{N} : V$ is a strong $\mathcal{C}$-module functor.*
- *The comparison functor $K_{VG}$ of a lax $\mathcal{C}$-module adjunction $G\colon \mathcal{M} \rightleftarrows \mathcal{N} : V$ is a strong $\mathcal{C}$-module functor.*

## 5.1 Bimonads

**Definition 5.1.** A *bimonad* on a monoidal category $\mathcal{C}$ consists of an oplax monoidal endofunctor $B$ on $\mathcal{C}$, as well as oplax monoidal natural transformations $\mu\colon B^2 \Longrightarrow B$, $\eta\colon \mathrm{Id}_{\mathcal{C}} \Longrightarrow B$ such that $(B, \mu, \eta)$ is a monad on $\mathcal{C}$.

A *morphism of bimonads* is a natural transformation $f\colon B \Longrightarrow H$ between bimonads that is oplax monoidal as well as a morphism of monads.

Let us make Definition 5.1 explicit. In order for a monad $(B, \mu, \eta)$ to be a bimonad, there has to exist a natural transformation $B_2\colon B \circ \otimes \Longrightarrow B \otimes B$ and a morphism $B_0\colon B1 \longrightarrow 1$, such that all diagrams in Figure 5.1 commute. Further, what it means for $\mu$ and $\eta$ to be oplax monoidal transformations in the diagrammatic calculus of Section 2.7.1 is displayed in Figure 5.2.

**Remark 5.2.** We can also express Definition 5.1 in more bicategorical language. Then, a bimonad may be defined as a monad[14] in the bicategory $\mathbb{O}\mathsf{plMon}^{\mathsf{opl}}$ where objects are monoidal categories, 1-cells are oplax monoidal functors, and 2-cells are oplax monoidal transformations.

[14] A monoid in the monoidal category $\mathbb{O}\mathsf{plMon}^{\mathsf{opl}}(\mathcal{C}, \mathcal{C})$.

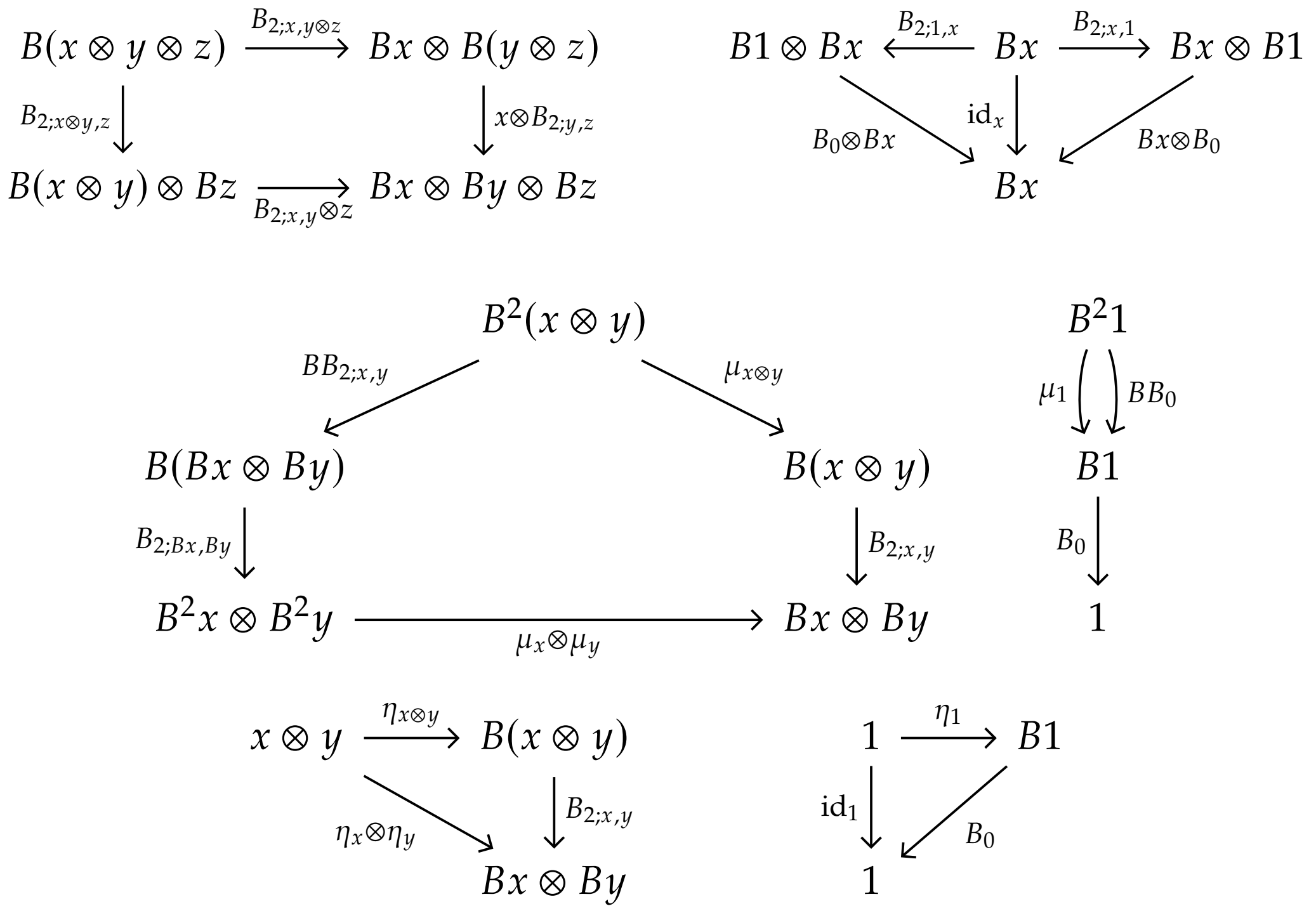

Figure 5.1: Explicit conditions for $B$ to be a bimonad.

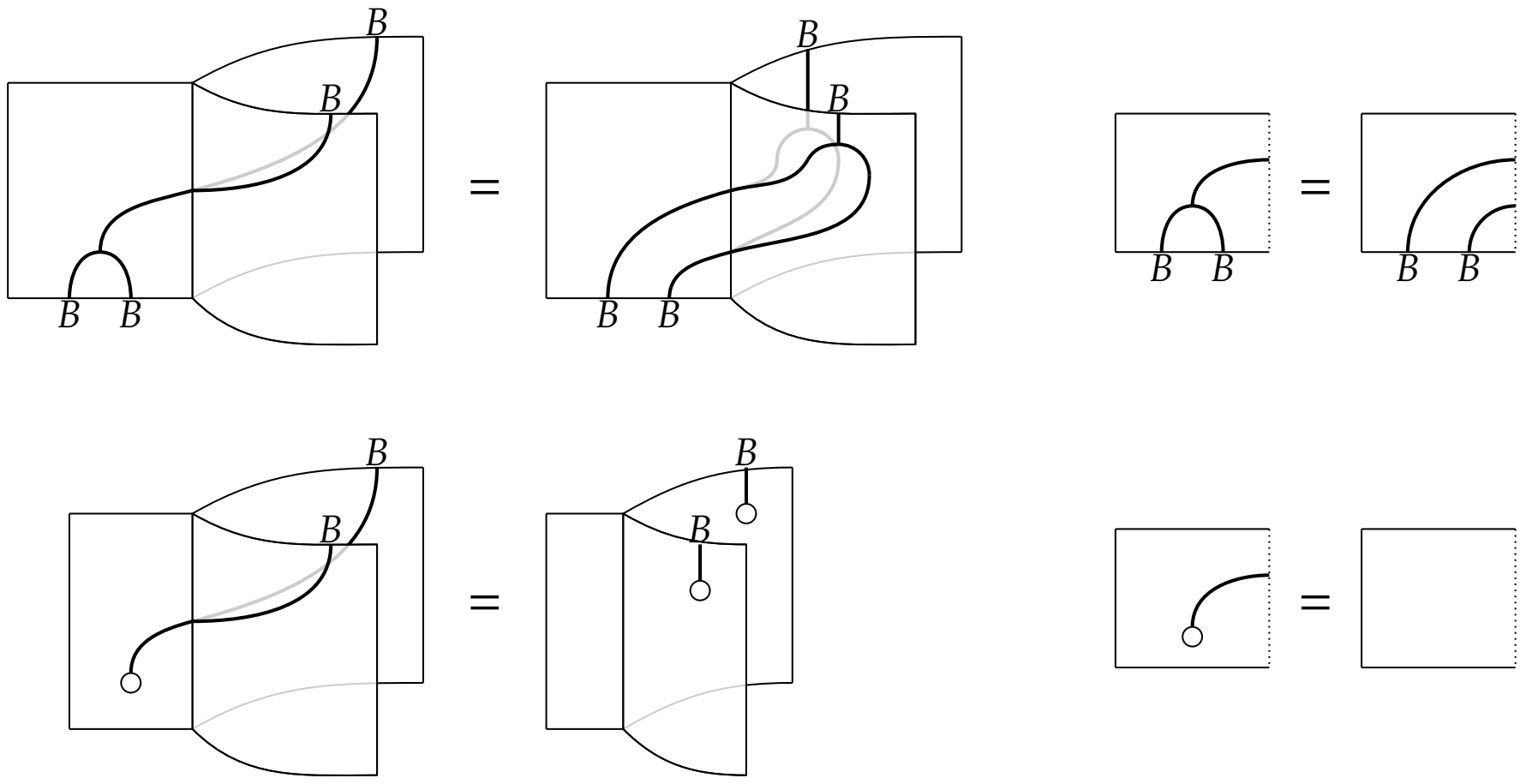

Figure 5.2: Conditions for the multiplication and unit of a bimonad to be oplax monoidal natural transformations.

**Example 5.3.** Let $A \in \mathscr{C}$ be an object in a monoidal category. In Example 2.96 we saw that the functor $A \otimes -$ is a monad on $\mathscr{C}$ if and only if $A$ is an algebra object in $\mathscr{C}$. This relationship extends to bialgebras: let $B \in \mathscr{C}$ be an algebra object, and suppose $T := B \otimes -$ is the associated monad on $\mathscr{C}$. Then $T$ is a bimonad if and only if $B$ is a bialgebra. More precisely, starting with a bialgebra structure $(\Delta, \varepsilon)$ on $B$, set

$$T_2 \colon B \otimes - \xrightarrow{\Delta \otimes -} B \otimes B \otimes - \qquad \text{and} \qquad T_0 \colon B \otimes - \xrightarrow{\varepsilon \otimes -} -.$$

This turns $T$ into a bimonad.

Conversely, assume $(T, T_2, T_0)$ is an oplax monoidal functor. One defines a bialgebra structure on $B$ by

$$\Delta \colon B = B \otimes 1 \otimes 1 \xrightarrow{T_{2;1,1}} B \otimes 1 \otimes B \otimes 1 = B \otimes B, \qquad \varepsilon \colon B = B \otimes 1 \xrightarrow{T_0} 1.$$

The fact that these arrows satisfy all of the required identities follows from a straightforward calculation.

**Example 5.4.** The intricate interplay between monads and adjunctions of Example 2.11 transcends to monoidal categories and bimonads. Given an oplax monoidal adjunction $F \colon \mathscr{C} \rightleftarrows \mathscr{D} :U$, the monad $UF \colon \mathscr{C} \longrightarrow \mathscr{C}$ is a bimonad, whose comultiplication is defined for every $x, y \in \mathscr{C}$ as the composition

$$UF(x \otimes y) \xrightarrow{UF_{2;x,y}} U(Fx \otimes Fy) \xrightarrow{U_{2;Fx,Fy}} UFx \otimes UFy.$$

Its counit is given by $UF1 \xrightarrow{UF_0} U1 \xrightarrow{U_0} 1$.

The following statement is a special case of [Kel74, Theorem 1.2].

**Proposition 5.5.** *Given monoidal categories $\mathscr{C}$ and $\mathscr{D}$, as well as an adjunction $F \colon \mathscr{C} \rightleftarrows \mathscr{D} :U$, there is a bijection between oplax monoidal structures on $F$ and lax monoidal structures on $U$.*

The proof of [Kel74, Theorem 1.2] is given constructively in terms of mates. Given an oplax monoidal structure $(F_2, F_0)$ on $F$, the lax monoidal natural transformation $U_2$ of $U$ is, for all $x, y \in \mathscr{D}$, given by

$$Ux \otimes Uy \xrightarrow{\eta_{Ux \otimes Uy}} UF(Ux \otimes Uy) \xrightarrow{UF_{2;Ux,Uy}} U(FUx \otimes FUy) \xrightarrow{U(\varepsilon_x \otimes \varepsilon_y)} U(x \otimes y),$$

where $\eta$ and $\varepsilon$ are the unit and counit of the adjunction. Constructing an oplax monoidal structure on $F$ from a lax monoidal structure on $U$ is similar.

Let us now highlight an important special case of Proposition 5.5.

**Lemma 5.6** ([TV17, Lemma 7.10]). *Any adjunction between two monoidal categories is oplax monoidal if and only if the right adjoint is strong monoidal.*

**Example 5.7.** Let $\Bbbk$ be a field, and $B \in \mathsf{Vect}$ a bialgebra. Then, by Example 5.4, the monad $B \otimes_{\Bbbk} -\colon \mathsf{Vect} \longrightarrow \mathsf{Vect}$ is a bimonad, and hence the adjunction $B \otimes_{\Bbbk} -\colon \mathsf{Vect} \rightleftarrows B\text{-Mod} : U$ is oplax monoidal, where $U\colon B\text{-Mod} \longrightarrow \mathsf{Vect}$ is the canonical forgetful functor that simply forgets the $B$-module structure. By Lemma 5.6, we obtain that $U$ is even strong monoidal.

More generally, let $B\colon \mathscr{C} \longrightarrow \mathscr{C}$ be the bimonad arising from an oplax monoidal adjunction $F\colon \mathscr{C} \rightleftarrows \mathscr{D} : U$. As the forgetful functor $U^B\colon \mathscr{C}^B \longrightarrow \mathscr{C}$ is strict monoidal, the Eilenberg–Moore adjunction $F^B \dashv U^B$ is, oplax monoidal.

The following result is due to [Kel74], see also [BV07, Theorem 2.6].

**Lemma 5.8.** *Let $F\colon \mathscr{C} \rightleftarrows \mathscr{D} : U$ be an oplax monoidal adjunction. Then the comparison functor $K^{UF}\colon \mathscr{D} \longrightarrow \mathscr{C}^{UF}$ of the bimonad $UF$ is strong monoidal, and $U^{UF}K^{UF} = U$ and $K^{UF}F = F^{UF}$ as strong respectively oplax monoidal functors.*

The question to which extent the monoidal structure on $\mathscr{C}^B$ is unique was answered by Moerdijk [Moe02, Theorem 7.1] and McCrudden [McC02, Corollary 3.13]. In particular, this kind of Tannaka–Krein reconstruction gives another justification for the name "bimonad".

**Proposition 5.9.** *Let $(B, \mu, \eta)$ be a monad on a monoidal category $\mathscr{C}$. There exists a one-to-one correspondence between bimonad structures on $B$ and monoidal structures on $\mathscr{C}^B$ such that the forgetful functor $U^B$ is strict monoidal.*

*Sketch of proof.* Given a bimonad $(B, \mu, \eta, B_2, B_0)\colon \mathscr{C} \longrightarrow \mathscr{C}$, as well as two modules $(m, \nabla_m)$ and $(n, \nabla_n) \in \mathscr{C}^B$, we set

$$(m, \nabla_m) \otimes (n, \nabla_n) := \big(m \otimes n,\ B(m \otimes n) \xrightarrow{B_{2;m,n}} Bm \otimes Bn \xrightarrow{\nabla_m \otimes \nabla_n} m \otimes n\big).$$

Moreover, define $\nabla_1\colon B1 \xrightarrow{B_0} 1$. The coassociativity and counitality of the comultiplication of $B$ imply that the above construction implements a monoidal structure on $\mathscr{C}^B$, parallel to that on the modules over a bialgebra.

Conversely, let $(B, \mu, \eta)$ be a monad on the monoidal category $\mathscr{C}$. Suppose $\mathscr{C}^T$ is monoidal such that the forgetful functor $U^T$ is strict monoidal. Consider the free modules $(Bm, \mu_m)$ and $(Bn, \mu_n)$, and denote their tensor product by

$$(Bm, \mu_m) \otimes (Bn, \mu_n) = (Bm \otimes Bn, \mu_m \cdot \mu_n\colon B(Bm \otimes Bn) \longrightarrow Bm \otimes Bn).$$

The comultiplication of $B$ is then given by

$$B(m \otimes n) \xrightarrow{B(\eta_m \otimes \eta_n)} B(Bm \otimes Bn) \xrightarrow{\mu_m \cdot \mu_n} Bm \otimes Bn.$$

For the counit, take the action $B1 \longrightarrow 1$ of the unit of $\mathscr{C}^B$. □

## 5.2 Hopf monads

We can now incorporate rigidity into the framework of Section 5.1.

**Definition 5.10** ([BLV11, Section 2.6]). Let $H$ be a bimonad on a monoidal category $\mathscr{C}$. The *right fusion operator* $H_{\mathsf{rf}}$ of $H$ is the natural transformation

$$H_{\mathsf{rf};x,y}\colon H(Hx \otimes y) \xrightarrow{H_{2;Hx,y}} H^2x \otimes Hy \xrightarrow{\mu_x \otimes Hy} Hx \otimes Hy, \qquad \text{for } x, y \in \mathscr{C}.$$

The bimonad $H$ is called *right Hopf* if its right fusion operator is invertible.

Left fusion operators and left Hopf monads are defined dually, and a bimonad is said to be a *Hopf monad* if it is both left and right Hopf.

**Theorem 5.11** ([BLV11, Theorem 3.6]). *Let $\mathscr{C}$ be a right closed monoidal category and let $H$ be a bimonad on $\mathscr{C}$. Then $H$ is a right Hopf monad if and only if $\mathscr{C}^H$ is right closed and the forgetful functor $U^H\colon \mathscr{C}^H \longrightarrow \mathscr{C}$ is right closed.*

A Hopf monad $H\colon \mathscr{C} \longrightarrow \mathscr{C}$ can be defined on any monoidal category. If $\mathscr{C}$ is rigid then one may use [BLV11, Lemma 3.4 and Theorem 3.6] to obtain a rigid—not just closed—monoidal structure on the category of $H$-modules. This is reflected by the existence of two natural transformations

(5.2.1) $$s^\ell\colon H({}^\vee H) \Longrightarrow {}^\vee H \qquad \text{and} \qquad s^r\colon H(H^\vee) \Longrightarrow H^\vee, \qquad \text{for all } x \in \mathscr{C},$$

called the *left* and *right antipode* of $H$. In Example 2.4 of [BV12] it is explained how these generalise the antipode of a Hopf algebra; an analogous result to Example 5.3 holds in the Hopf case.

## 5.3 (co)module monads

Alongside the theory of comodule monads of Aguiar and Chase, [AC12], we will also develop the special case of (op)lax $\mathscr{C}$-module monads in this section. The latter will play an important role in Chapters 8 and 9.

**Definition 5.12** ([AC12, Definition 3.3.1]). Let $(F, F_2, F_0)\colon \mathscr{C} \longrightarrow \mathscr{D}$ be an oplax monoidal functor, $\mathscr{M}$ a right $\mathscr{C}$-module, and $\mathscr{N}$ a right $\mathscr{D}$-module category. A *(right) comodule functor over* $F$[15] is a pair $(G, G_{\mathsf{a}})$ consisting of a functor $G\colon \mathscr{M} \longrightarrow \mathscr{N}$ together with a natural transformation

$$G_{\mathsf{a};m,x}\colon G(m \triangleleft x) \longrightarrow Gm \triangleleft Fx, \qquad \text{for all } x \in \mathscr{C} \text{ and } m \in \mathscr{M},$$

called the *coaction* of $G$, which is *coassociative* and *counital*, in the sense that the following diagrams commute for all $x, y \in \mathscr{C}$ and $m \in \mathscr{M}$:

$$\begin{array}{ccc}
G((m \triangleleft x) \triangleleft y) & \xrightarrow{\cong} & G(m \triangleleft (x \otimes y)) \\
{\scriptstyle G_{\mathsf{a};m \triangleleft x,y}}\downarrow & & \downarrow{\scriptstyle G_{\mathsf{a};m,x\otimes y}} \\
G(m \triangleleft x) \triangleleft Fy & & Gm \triangleleft F(x \otimes y) \\
{\scriptstyle G_{\mathsf{a};m,x} \triangleleft Fy}\downarrow & & \downarrow{\scriptstyle Gm \triangleleft F_{2;x,y}} \\
(Gm \triangleleft Fx) \triangleleft Fy & \xrightarrow[\cong]{} & Gm \triangleleft (Fx \otimes Fy)
\end{array}
\qquad
\begin{array}{ccc}
G(m \triangleleft 1) & \xrightarrow{\cong} & Gm \\
{\scriptstyle G_{\mathsf{a};m,1}}\downarrow & & \uparrow{\scriptstyle \cong} \\
Gm \triangleleft F1 & \xrightarrow[Gm \triangleleft F_0]{} & Gm \triangleleft 1
\end{array}$$

A comodule functor is called *strong* if its coaction is an isomorphism, and *strict* if it is the identity.

[15] Alternatively, we will say that $G$ is a *(right) $F$-comodule functor*.

**Example 5.13.** A right comodule functor over the identity functor $\mathrm{Id}_{\mathscr{C}}$ of a monoidal category $\mathscr{C}$ is nothing more than an oplax $\mathscr{C}$-module functor in the sense of Definition 2.47, considering right instead of left $\mathscr{C}$-module categories.

**Example 5.14** ([AC12, Section 6.1]). If $B$ is a bialgebra over a commutative ring $\Bbbk$, then $B \otimes_{\Bbbk} -\colon \Bbbk\text{-Mod} \longrightarrow \Bbbk\text{-Mod}$ becomes an oplax monoidal functor. Let $A$ be a right $B$-comodule algebra with coaction $\nu\colon A \longrightarrow A \otimes_{\Bbbk} B$. Recall that the subalgebra of *$B$-coinvariants* of $A$ is given by

$$A^{\mathrm{co}B} := \{\, a \in A \mid \nu(a) = a \otimes 1 \,\}.$$

Suppose that $C$ is a subalgebra of $A^{\mathrm{co}B}$. Then $\nu$ is a map of $C$-$C$-bimodules, which turns $A \otimes_C -\colon C\text{-Mod} \longrightarrow C\text{-Mod}$ into a right comodule functor over $B \otimes_{\Bbbk} -$. The action of $\Bbbk$-Mod on $C$-Mod is given by tensoring over $\Bbbk$.

Recall the constructions of the centre and twisted centre from Sections 2.4.4 and 4.2.1. The canonical functor $U^{(\mathsf{Z})}\colon \mathsf{Z}(\mathscr{C}) \longrightarrow \mathscr{C}$ that forgets the half braiding is strict monoidal. Similarly, given a strong monoidal functor $L$, the forgetful functor $U^{(L)}\colon \mathsf{Z}({}_L\mathscr{C}) \longrightarrow \mathscr{C}$ is a strict comodule functor over $U^{(\mathsf{Z})}$.

**Remark 5.15.** The diagrammatic calculus of Section 2.7.1 can be adapted to the setting of comodule functors. In addition to combining sheets using tensor products, we now additionally consider actions to do so as well. In order to keep track of which splitting occurred, functors between the module categories shall be coloured blue. For example, a coaction will be drawn as:

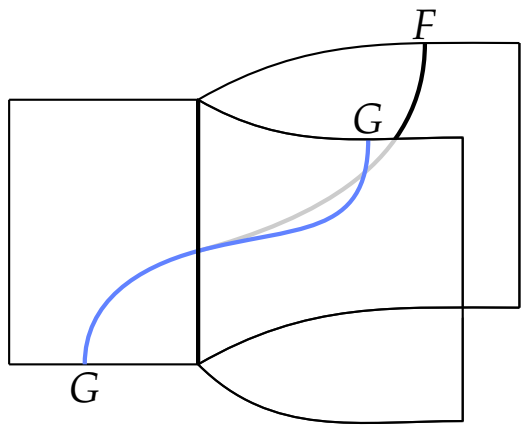

The coassociativity and counitality of $G_{\mathsf{a}}$ is displayed in Figure 5.3; notice the analogous nature of the diagrams to Figure 2.4.

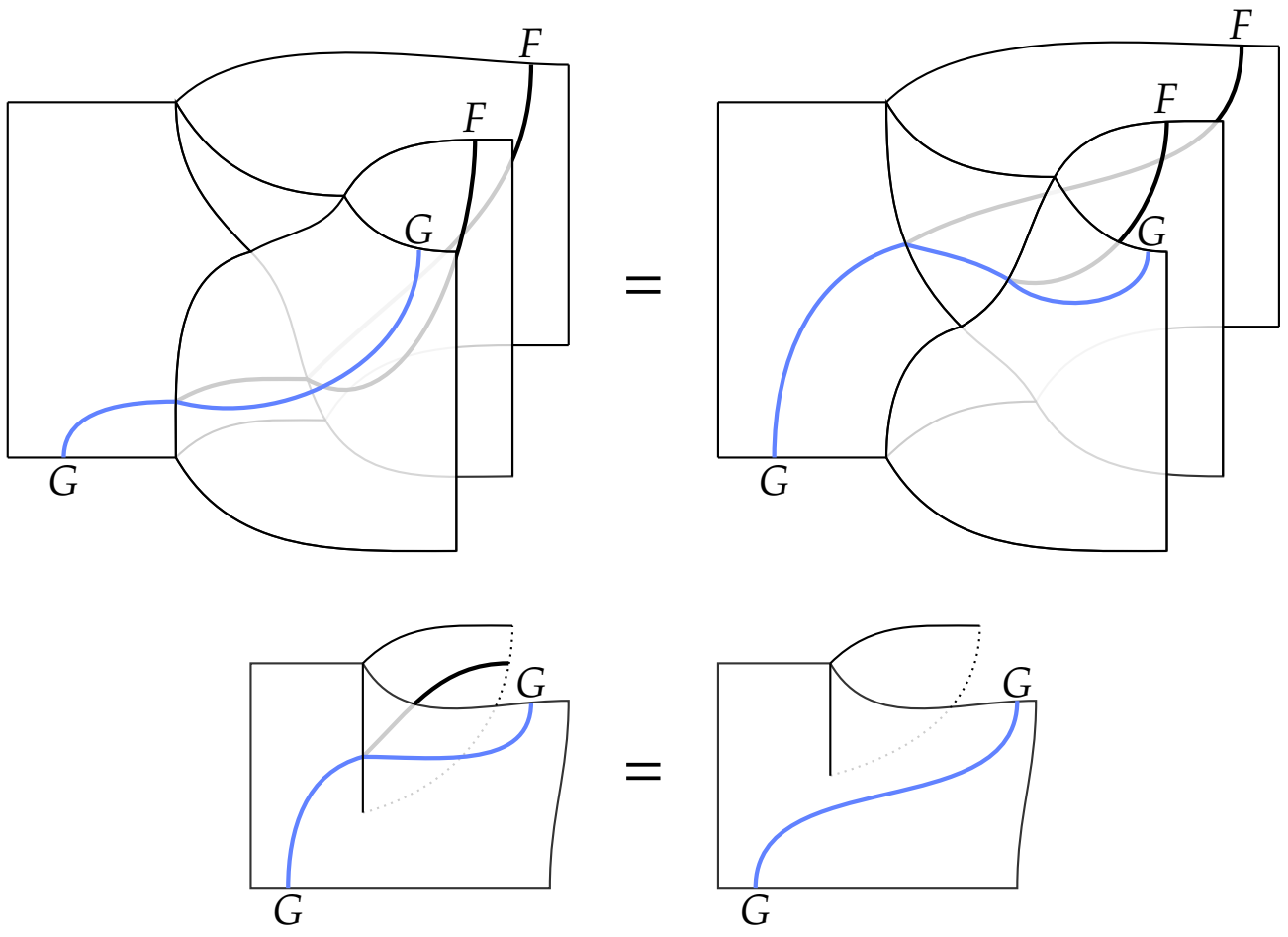

Figure 5.3: Coassociativity and counitality conditions of the coaction of a comodule functor.

**Definition 5.16.** Suppose that $C, G\colon \mathcal{M} \longrightarrow \mathcal{N}$ are comodule functors over $B, F\colon \mathcal{C} \longrightarrow \mathcal{D}$. A *comodule (natural) transformation* from $C$ to $G$ comprises a pair of natural transformations $\phi\colon C \Longrightarrow G$ and $\psi\colon B \Longrightarrow F$, such that the following diagram commutes for all $x \in \mathcal{C}$ and $m \in \mathcal{M}$:

(5.3.1)
$$\begin{array}{ccc} C(m \triangleleft x) & \xrightarrow{\phi_{m \triangleleft x}} & G(m \triangleleft x) \\ {\scriptstyle C_{\mathsf{a};m,x}}\downarrow & & \downarrow{\scriptstyle G_{\mathsf{a};m,x}} \\ Cm \triangleleft Bx & \xrightarrow[\phi_m \triangleleft \psi_x]{} & Gm \triangleleft Fx \end{array}$$

We call $(\phi, \psi)$ a *morphism of comodule functors* if $B = F$ and $\psi = \mathrm{id}_B$.

**Example 5.17.** If $\mathscr{C} = \mathscr{D}$ and $B = F = \mathrm{Id}_{\mathscr{C}}$, then $C$ and $G$ are oplax $\mathscr{C}$-module functors, see Example 5.13. A comodule natural transformation from $C$ to $G$ is exactly a $\mathscr{C}$-module transformation in the sense of Definition 2.48.

Given a comodule natural transformation ($\phi\colon G \Longrightarrow C$, $\psi\colon B \Longrightarrow F$), the graphical version of Diagram (5.3.1) is displayed in our next picture, where the blue dot represents $\phi$ and the black dot represents $\psi$.

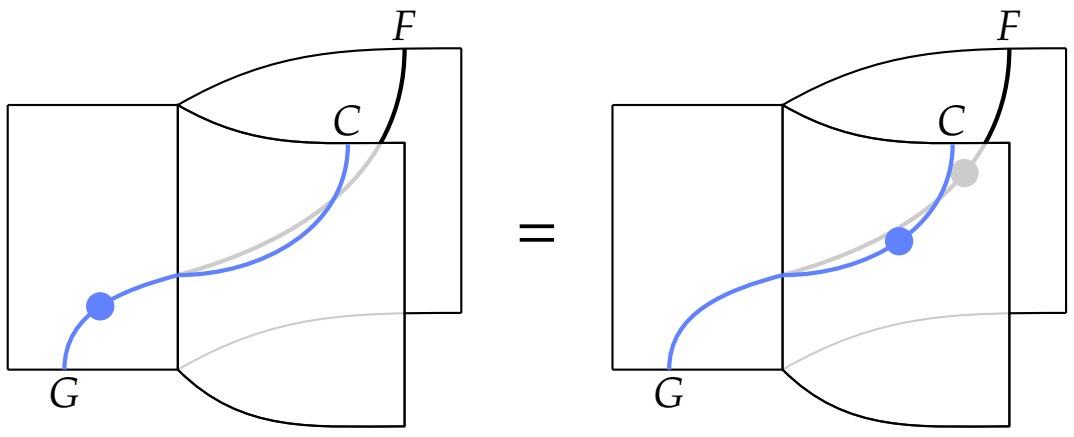

The notion of a morphism of comodule functors might seem restrictive; however, it is actually sufficient to only consider arrows of this form.

**Example 5.18.** Let the pair $\phi\colon C \Longrightarrow G$ and $\psi\colon B \Longrightarrow F$ be a comodule transformation. We can view $\phi\colon C \Longrightarrow G$ as a morphism of comodule functors over $F$ if we equip $C$ with a new coaction. It is given by

$$C(m \triangleleft x) \xrightarrow{C_{\mathsf{a};m,x}} Cm \triangleleft Bx \xrightarrow{Cm \triangleleft \psi_x} Cm \triangleleft Fx, \qquad \text{for all } x \in \mathscr{C} \text{ and } m \in \mathscr{M}.$$

Thus, by suitably altering the involved coactions, comodule natural transformations and morphisms of comodule functors can be identified.

Let $B\colon \mathscr{C} \longrightarrow \mathscr{C}$ be a bimonad and $\mathscr{M}$ a module category over $\mathscr{C}$. The unit $\eta\colon \mathrm{Id}_{\mathscr{C}} \Longrightarrow B$ implements a coaction on $\mathrm{Id}_{\mathscr{M}}\colon \mathscr{M} \Longrightarrow \mathscr{M}$ via

$$\mathrm{id}_m \triangleleft \eta_x\colon \mathrm{Id}_{\mathscr{M}}(m \triangleleft x) \longrightarrow \mathrm{Id}_{\mathscr{M}} m \triangleleft Bx, \qquad \text{for all } x \in \mathscr{C}, m \in \mathscr{M}.$$

Using the multiplication $\mu\colon B^2 \Longrightarrow B$, we can equip the composition $CG$ of two comodule functors $C, G\colon \mathscr{M} \longrightarrow \mathscr{M}$ with a comodule structure:

$$(CG)_{\mathsf{a}}\colon CG(- \triangleleft =) \xrightarrow{CG_{\mathsf{a}}} C(G(-) \triangleleft B(=)) \xrightarrow{C_{\mathsf{a}}} CG(-) \triangleleft B^2(=) \xrightarrow{\mathrm{id} \triangleleft \mu} CG(-) \triangleleft B(=).$$

Due to the associativity and unitality of the multiplication of $B$, in this way the category of comodule endofunctors on $\mathscr{M}$ over $B$ becomes monoidal.

**Definition 5.19.** Consider a bimonad $B\colon \mathscr{C} \longrightarrow \mathscr{C}$ and a module category $\mathscr{M}$ over $\mathscr{C}$. A *comodule monad* over $B$ on $\mathscr{M}$ comprises a comodule endofunctor $(C, C_{\mathsf{a}})\colon \mathscr{M} \longrightarrow \mathscr{M}$ together with comodule morphisms $\mu\colon C^2 \Longrightarrow C$ and $\eta\colon \mathrm{Id}_{\mathscr{M}} \Longrightarrow C$ such that $(C, \mu, \eta)$ is a monad.

A *morphism of comodule monads* is a natural transformation of comodule functors $f\colon C \Longrightarrow G$ that is also a morphism of monads.

The conditions for the multiplication and unit of a comodule monad $C$ on $\mathscr{M}$ over a bimonad $B$ on $\mathscr{C}$ to be morphisms of comodule functors is given in Figure 5.4. Notice how the conditions are analogous to Figure 5.2.

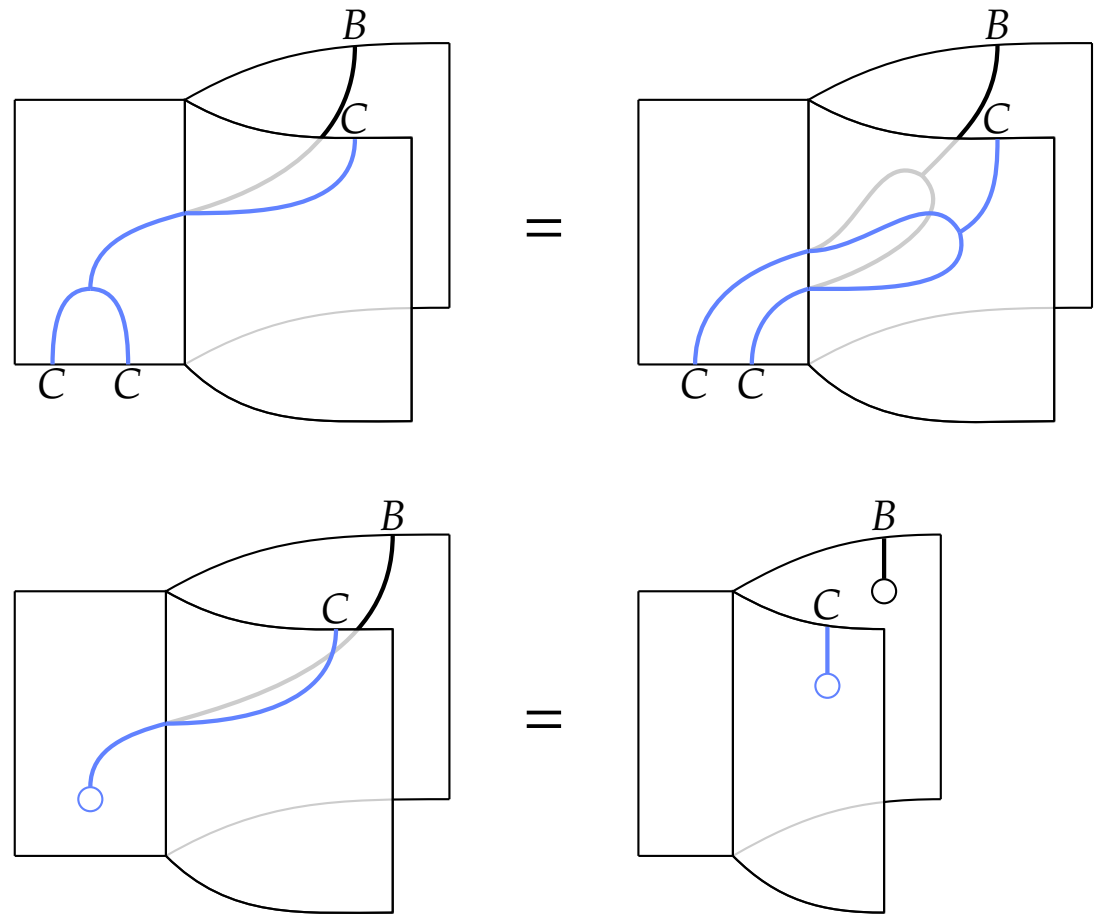

Figure 5.4: Conditions for the multiplication and unit of a comodule monad to be morphisms of comodule functors.

**Example 5.20.** An *oplax $\mathscr{C}$-module monad* is a comodule monad over the identity monad on $\mathscr{C}$. Alternatively, it is monoid in the category $\mathsf{Oplax}\mathscr{C}\mathsf{Mod}(\mathscr{M}, \mathscr{M})$ of oplax $\mathscr{C}$-module endofunctors on a (left) $\mathscr{C}$-module category $\mathscr{M}$.

Analogously, one can define the notion of a *lax $\mathscr{C}$-module comonad*.

**Example 5.21.** Consider the poset $(\mathbb{R}, \leq)$. It is a monoidal category with addition as tensor product and $0 \in \mathbb{R}$ as unit. In [HL18], Hasegawa and Lemay defined a bimonad using the ceiling function. Let

$$H\colon \mathbb{R} \longrightarrow \mathbb{R}, \qquad x \longmapsto \lceil x \rceil := \min\{\, n \in \mathbb{Z} \mid n \geq x \,\}.$$

As $\lceil x + y \rceil \leq \lceil x \rceil + \lceil y \rceil$ for all $x, y \in \mathbb{R}$ and $\lceil 0 \rceil = 0$, the functor $H$ is oplax monoidal. The comultiplication and counit are given by the unique arrows

$$H_0\colon H0 = 0 \longrightarrow 0 \qquad \text{and} \qquad H_{2;x,y}\colon H(x+y) \longrightarrow Hx + Hy,$$

for all $x, y \in \mathbb{R}$. The idempotence of the ceiling function implies that the identity natural transformation defines a multiplication $\mu\colon H^2 \Longrightarrow H$. Its unit corresponds to $\{\, x \longrightarrow Hx \,\}_{x\in\mathbb{R}}$.

Given a number $z$ in the half-open interval $[0,1)$, let

$$C^{(z)}\colon \mathbb{R} \longrightarrow \mathbb{R}, \qquad x \longmapsto \min\{\, z+n \mid n \in \mathbb{Z}, z+n \geq x \,\}.$$

In case $z = 0$, we have $C^{(z)} = H$. Otherwise $C^{(z)}0 = z > 0$ and $C^{(z)}$ cannot be oplax monoidal. Nonetheless, $C^{(z)}(x+y) \leq C^{(z)}x + \lceil y \rceil = C^{(z)}x + Hy$ holds, and the unique natural arrow

$$C^{(z)}_{\mathsf{a};x,y}\colon C^{(z)}(x+y) \longrightarrow C^{(z)}x + Hy, \qquad \text{for all } x, y \in \mathbb{R},$$

defines a coaction of $H$ on $C$. Again, we have that $\left(C^{(z)}\right)^2 = C^{(z)}$ is idempotent and $x \leq C^{(z)}x$, for all $x \in \mathbb{R}$. Thus, it is a comodule monad over $H$.

**Example 5.22.** Let $(\mathcal{V}, \otimes, 1)$ be a closed symmetric monoidal category. A $\mathcal{V}$-category $\mathcal{C}$ is said to be *copowered* over $\mathcal{V}$, see [Kel05, Section 3.7], if there exists a functor $- \cdot = \colon \mathcal{C} \times \mathcal{V} \longrightarrow \mathcal{C}$, such that for all $c \in \mathcal{C}$ we have

$$c \cdot - \colon \mathcal{V} \rightleftarrows \mathcal{C} : \mathcal{C}(c, -).$$

One obtains $c \cdot (v \otimes w) \cong (c \cdot v) \cdot w$ by the Yoneda lemma:

$$\begin{aligned} \mathcal{C}\big(c \cdot (v \otimes w), x\big) &\cong \mathcal{V}\big(v \otimes w, \mathcal{C}(c, x)\big) \cong \mathcal{V}\big(w, \mathcal{V}(v, \mathcal{C}(c, x))\big) \\ &\cong \mathcal{V}\big(w, \mathcal{C}((c \cdot v), x)\big) \cong \mathcal{C}\big((c \cdot v) \cdot w, x\big). \end{aligned}$$

In fact, this turns $\mathcal{C}$ into a $\mathcal{V}$-module category, and therefore the identity monad on $\mathcal{C}$ into a comodule monad over $\mathrm{Id}_{\mathcal{V}}$ with trivial coaction.

Suppose that $B \coloneqq UF$ is a bimonad on $\mathcal{V}$. The unit $\eta\colon \mathrm{Id}_{\mathcal{V}} \Longrightarrow B$ is a morphism of bimonads and we may extend the coaction of $\mathrm{Id}_{\mathcal{C}}$ to

$$\delta\colon - \cdot = \xrightarrow{\;- \cdot \eta\;} - \cdot B=,$$

turning $B$ into a $\mathcal{C}$-module monad.

The following proposition connects the notion of a $\mathcal{C}$-module monad to that of an algebra object in $\mathcal{C}$ in the sense of Section 2.6.

**Proposition 5.23.** *Let $\mathcal{C}$ be a monoidal category. There is a bijection between strong $\mathcal{C}$-module monads on $\mathcal{C}$ and algebra objects in $\mathcal{C}$.*

*Proof.* By Proposition 2.51, there is an equivalence of left module categories

$$\mathsf{Str}\mathscr{C}\mathsf{Mod}(\mathscr{C},\mathscr{C}) \xrightarrow{\sim} \mathscr{C}^{\mathrm{rev}}, \qquad T \longmapsto T1, \qquad -\otimes A \longleftarrow A,$$

and by Example 2.96, $A$ is an algebra object if and only if $-\otimes A$ is a monad.

In particular, if $(A,\tau,\upsilon)$ is an algebra object in $\mathscr{C}$, then $(-\otimes A,\ -\otimes\tau,\ -\otimes\upsilon)$ becomes a strong $\mathscr{C}$-module monad on $\mathscr{C}$. The strong $\mathscr{C}$-module structure is given by the associator of $\mathscr{C}$:

$$(-\otimes A)_{\mathsf{a};X} \coloneqq X\otimes(-\otimes A) \xrightarrow{\sim} (X\otimes -)\otimes A, \qquad \text{for all } X\in\mathscr{C}.$$

If $(T,\mu,\eta)$ is a strong $\mathscr{C}$-module monad, then evaluating at the monoidal unit yields an algebra object in $\mathscr{C}$, with multiplication $\mu_1$ and unit $\eta_1$. □

In particular, the category of right $A$-modules and the Eilenberg–Moore category of $-\otimes A$ coincide as $\mathscr{C}$-module categories.

**Remark 5.24.** Let $B\colon \mathscr{C}\longrightarrow\mathscr{C}$ be a bimonad and $(C,C_{\mathsf{a}})\colon \mathscr{M}\longrightarrow\mathscr{M}$ a comodule monad over it. The coaction of $C$ allows us to define an action

$$\triangleleft\colon \mathscr{M}^C\times\mathscr{C}^B\longrightarrow\mathscr{M}^C.$$

For any two modules $(m,\nabla_m)\in\mathscr{M}^C$ and $(x,\nabla_x)\in\mathscr{C}^B$, it is given by

$$(m,\nabla_m)\triangleleft(x,\nabla_x) \coloneqq \big(m\triangleleft x, (\nabla_m\triangleleft\nabla_x)\circ C_{\mathsf{a};m,x}\big).$$

The axioms of the coaction of $B$ on $C$ translate precisely to the compatibility of the action of $\mathscr{C}^B$ on $\mathscr{M}^C$ with the tensor product and unit of $\mathscr{C}^B$.

**Definition 5.25.** Let $F\colon \mathscr{C}\rightleftarrows\mathscr{D} : U$ be an oplax monoidal adjunction. Suppose that $G\colon \mathscr{M}\rightleftarrows\mathscr{N} : V$ is an adjunction such that $G$ is an $F$-comodule functor, and $V$ is a $U$-comodule functor. The pair $(G\dashv V,\ F\dashv U)$ is a *comodule adjunction* if the following identities hold:

$$\begin{array}{ccc} m\triangleleft x & \xrightarrow{\eta^{(G\dashv V)}_{m\triangleleft x}} & VG(m\triangleleft x) \\ {\scriptstyle \eta^{(G\dashv V)}_m\triangleleft\eta^{(F\dashv U)}_x}\big\downarrow & & \big\downarrow{\scriptstyle VG_{\mathsf{a};m,x}} \\ VGm\triangleleft UFx & \xleftarrow[V_{\mathsf{a};Gm,Fx}]{} & V(Gm\triangleleft Fx) \end{array} \qquad \begin{array}{ccc} GV(m\triangleleft x) & \xrightarrow{GV_{\mathsf{a};m,x}} & G(Vm\triangleleft Ux) \\ {\scriptstyle \varepsilon^{(G\dashv V)}_{m\triangleleft x}}\big\downarrow & & \big\downarrow{\scriptstyle G_{\mathsf{a};Vm,Ux}} \\ m\triangleleft x & \xleftarrow[\varepsilon^{(G\dashv V)}_m\triangleleft\varepsilon^{(F\dashv U)}_x]{} & GVm\triangleleft FUx \end{array}$$

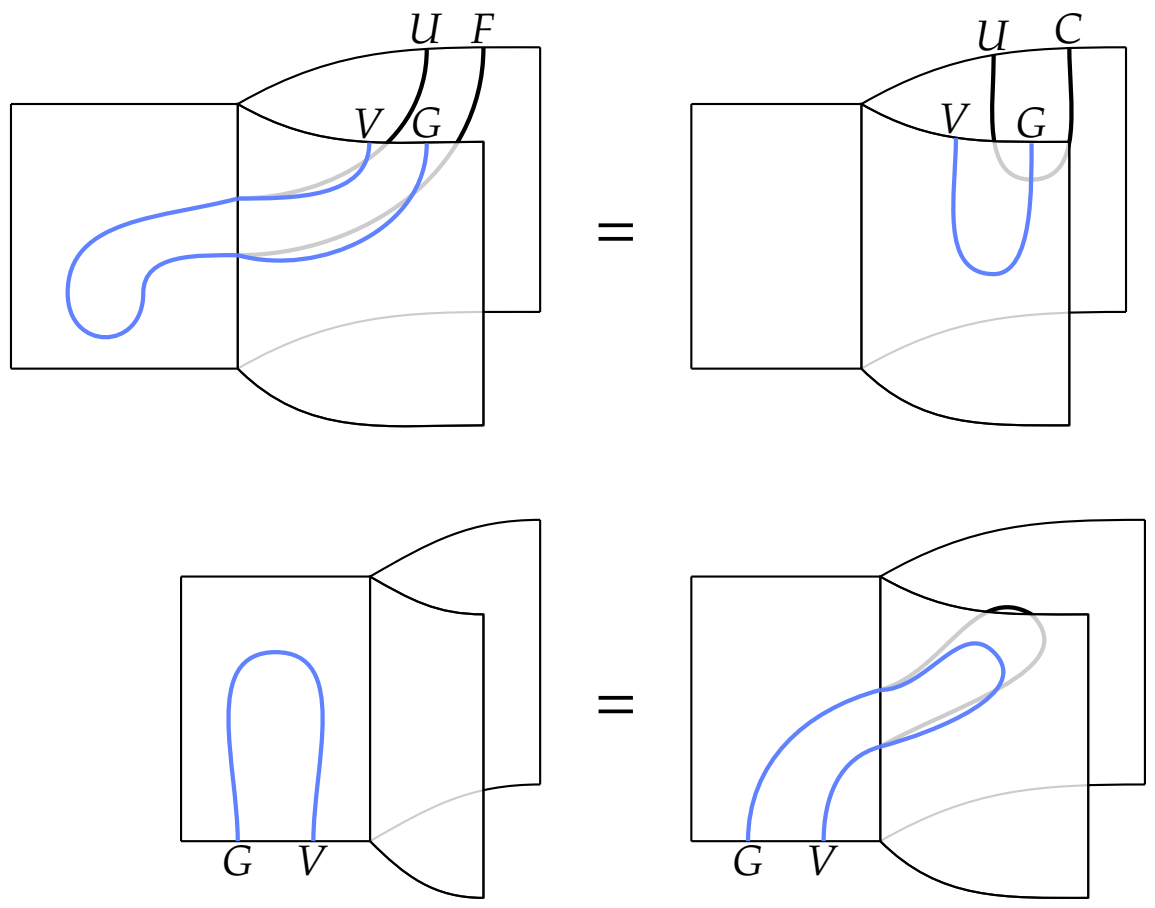

Figure 5.5: String diagrammatic conditions for a comodule adjunction.

**Example 5.26.** If $F = U = \mathrm{Id}_{\mathcal{M}}$ for a (left) $\mathcal{C}$-module category $\mathcal{M}$, then we say that $G \dashv V$ is an *oplax $\mathcal{C}$-module adjunction*. More explicitly, we have that $G$ and $V$ are oplax $\mathcal{C}$-module functors such that the unit and counit of the adjunction are $\mathcal{C}$-module transformations.

Analogously to this, one can define the notion of a *lax $\mathcal{C}$-module adjunction*.

Figure 5.5 makes plain that the conditions required by Definition 5.25 are analogous to those stated in Figure 2.5.

**Example 5.27.** The philosophy that monads and adjunctions are two sides of the same coin extends to comodule functors. Suppose that we have an oplax monoidal adjunction $F\colon \mathcal{C} \rightleftarrows \mathcal{D} : U$ and over it a comodule adjunction $G\colon \mathcal{M} \rightleftarrows \mathcal{N} : V$. By [AC12, Proposition 4.3.1], the bimonad $B := UF$ admits a coaction on the monad $C := VG$; for all $m \in \mathcal{M}$ and $x \in \mathcal{C}$ it is given by

$$VG(m \triangleleft x) \xrightarrow{VG\mathsf{a}_{;m,x}} V(Gm \triangleleft Fx) \xrightarrow{V\mathsf{a}_{;Gm,Fx}} VGm \triangleleft UFx = Cm \triangleleft Bx.$$

### 5.3.1 *Reconstruction for comodule monads*

The next theorem offers a way to reconstruct comodule monads from their module categories. It extends [AC12, Proposition 4.1.2] and [TV17, Lemma 7.10]. Recall Definition 2.38, and note that one can analogously define lifts of adjunctions to comodule adjunctions.

**Theorem 5.28.** *Let $\mathcal{C}$ and $\mathcal{D}$ be monoidal categories, and suppose that $\mathcal{M}$ and $\mathcal{N}$ are right $\mathcal{C}$- and $\mathcal{D}$-module categories, respectively. Let $F\colon \mathcal{C} \rightleftarrows \mathcal{D} : U$ be an oplax monoidal adjunction. Lifts of an adjunction $G\colon \mathcal{M} \rightleftarrows \mathcal{N} : V$ to a comodule adjunction are in bijection with lifts of $V\colon \mathcal{N} \longrightarrow \mathcal{M}$ to a strong $U$-comodule functor.*

*Proof.* Let $G \dashv V$ be a comodule adjunction over $F \dashv U$. Define the inverse $V_{\mathsf{a}}^{-1} \colon V - \triangleleft U = \Longrightarrow V(- \triangleleft =)$ of $V_{\mathsf{a}}$ by

(5.3.2) 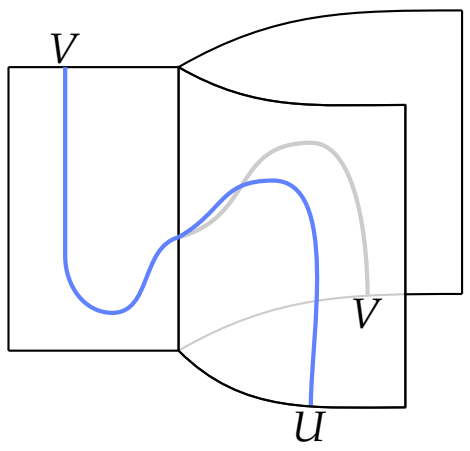

Using that $G$ and $V$ are part of a comodule adjunction, a straightforward computation proves $V_{\mathsf{a}}^{-1} \circ V_{\mathsf{a}} = \mathrm{id}_{V(- \triangleleft =)}$:

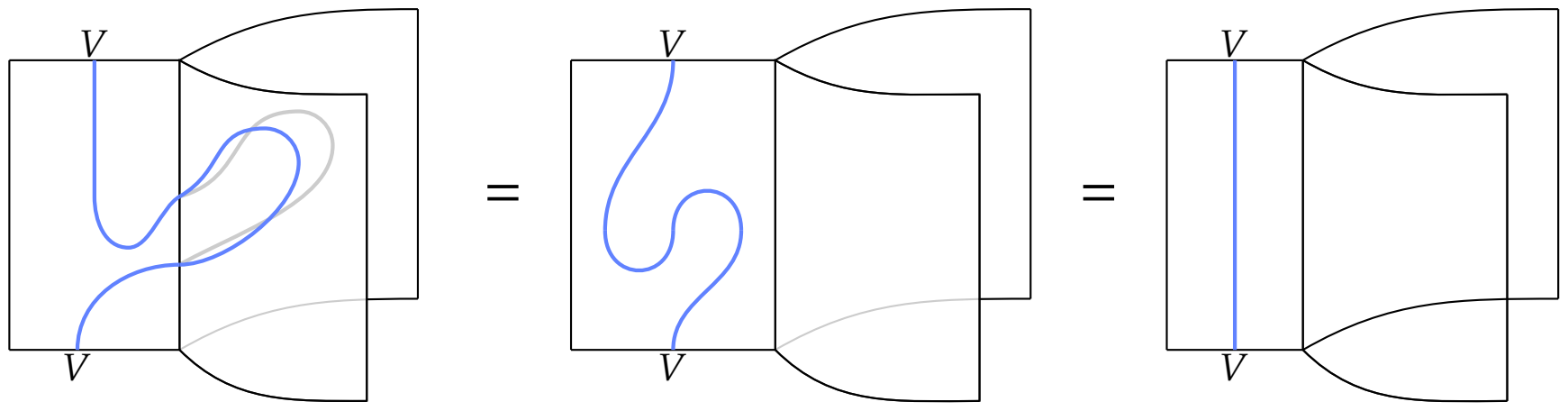

A similar strategy can be used to show that $V_{\mathsf{a}} \circ V_{\mathsf{a}}^{-1} = \mathrm{id}_{V - \triangleleft U =}$. Thus, $V_{\mathsf{a}}$ is a natural isomorphism and therefore $V$ is a strong comodule functor.

Conversely, suppose that $V \colon \mathcal{N} \longrightarrow \mathcal{M}$ is a strong comodule functor over $U$. Define an arrow $G_{\mathsf{a}} \colon G(- \triangleleft =) \Longrightarrow G(-) \triangleleft F(=)$ by

(5.3.3) 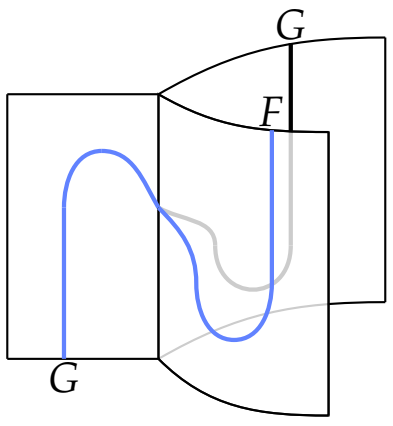

By [TV17, Lemma 7.10], the comultiplication and counit of $F$ are given by

$$F(x \otimes y) \xrightarrow{F(\eta_x \otimes \eta_y)} F(UFx \otimes UFy) \xrightarrow{FU_{-2;Fx,Fy}} FU(Fx \otimes Fy) \xrightarrow{\varepsilon_{Fx \otimes Fx}} Fx \otimes Fy,$$

$$F1 \xrightarrow{FU_{-0}} FU1 \xrightarrow{\varepsilon_1} 1, \qquad \text{for all } x, y \in \mathscr{C}.$$

Note that, graphically, $F_2$ looks just like Diagram (5.3.3), with black strings taking the place of blue ones.

We will show that $G$ is a comodule functor over $F$. Figure 5.6 shows that $G_{\mathsf{a}} \colon G(- \triangleleft =) \Longrightarrow G(-) \triangleleft F(=)$ is coassociative in the sense of Figure 5.3.

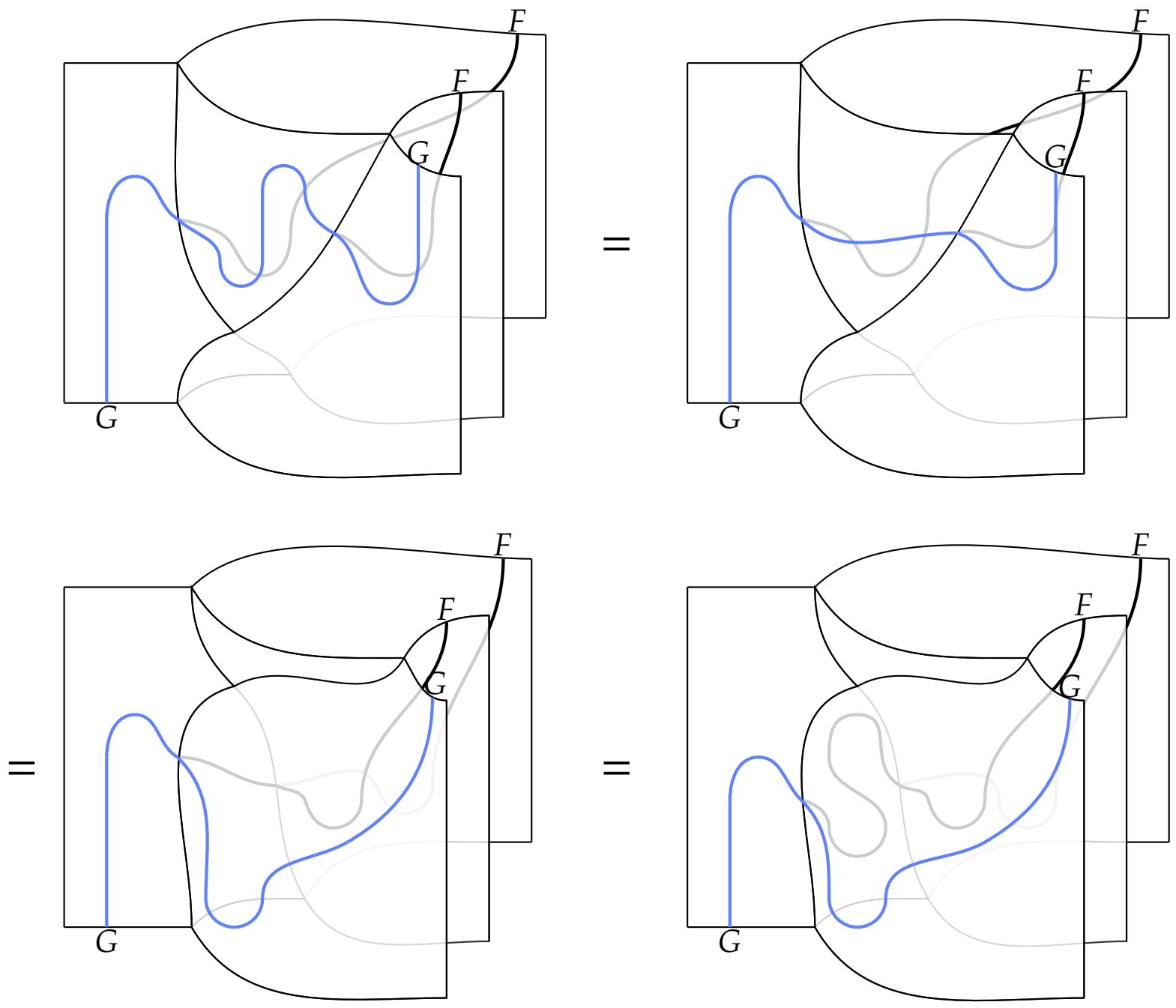

Figure 5.6: The coassociativity condition of $G$.

The fact that $G_{\mathsf{a}}$ is counital follows from

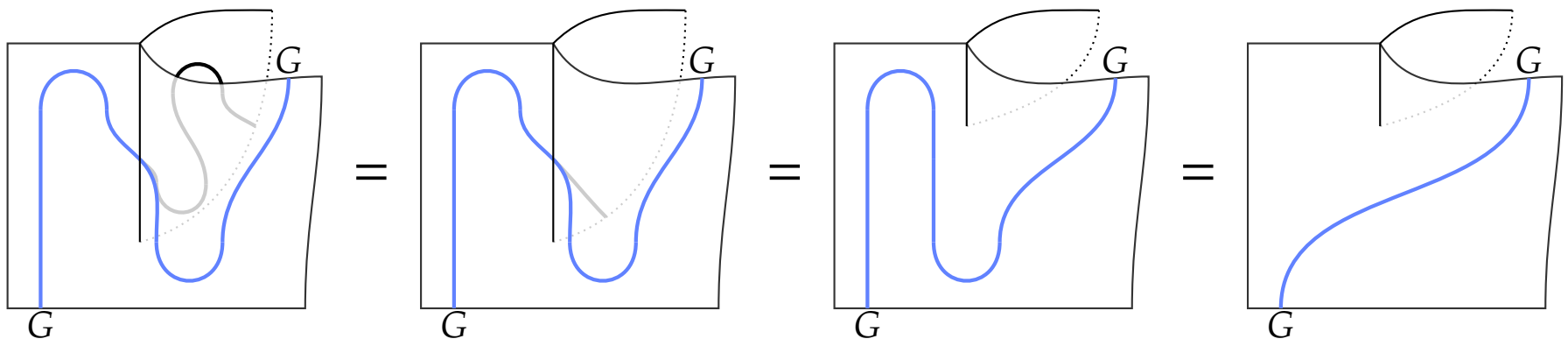

A straightforward computation proves that the unit of the adjunction $G \dashv V$ satisfies the axioms displayed in Figure 5.5:

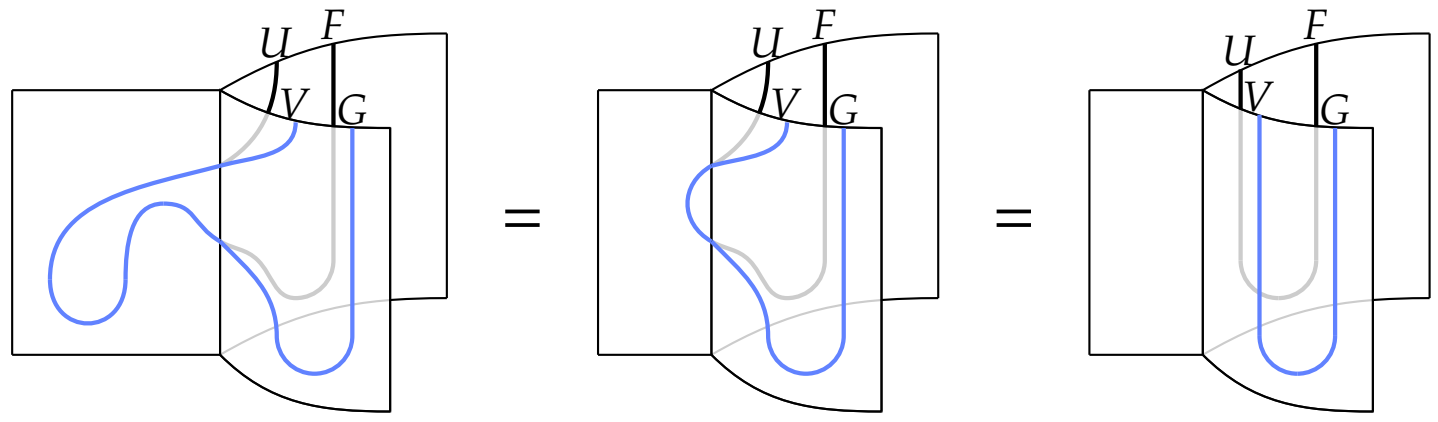

A similar argument for the counit shows that $G \dashv V$ is a comodule adjunction.

To see that these constructions are inverse to each other, first suppose that we have a comodule adjunction $G \dashv V$. By utilising $V_{\mathsf{a}}^{-1}$ as given in Diagram (5.3.2), we obtain another coaction $\lambda$ on $G$, see Diagram (5.3.3). Now, a direct computation shows that $G_{\mathsf{a}} = \lambda$:

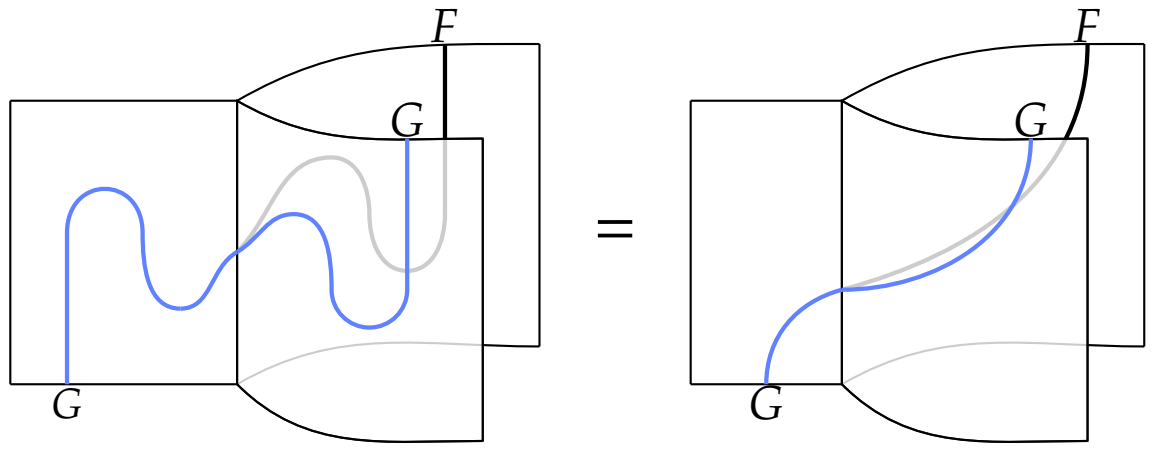

The converse is clear: the map that associates a comodule adjunction $G \dashv V$ to any strong comodule structure on $V$ preserves the coaction of $V$. □

The proof of Theorem 5.28 yields an analogue of Proposition 5.5, which we formulate in the language of $\mathscr{C}$-module functors.

**Porism 5.29.** *Let $\mathcal{M}$ and $\mathcal{N}$ be $\mathscr{C}$- and $\mathscr{D}$-module categories, respectively. For an adjunction $F\colon \mathcal{M} \rightleftarrows \mathcal{N} : U$, there is a bijection between oplax $\mathscr{C}$-module structures on $F$ and lax $\mathscr{C}$-module structures on $U$.*

**Corollary 5.30.** *Let $\mathcal{M}$ be a $\mathscr{C}$-module category, and let $T$ be a monad on $\mathcal{M}$. Then there exists a bijective correspondence*

$$\{\text{oplax } \mathscr{C}\text{-module monad structures on } T\} \leftrightarrow \left\{\begin{array}{l}\mathscr{C}\text{-module category structures on } \mathcal{M}^T \\ \text{such that } U^T \text{ is a strict module functor}\end{array}\right\}$$

$$\{\, T(- \triangleright =) \xrightarrow{T_{\mathsf{a}}} - \triangleright T(=) \,\} \longmapsto \{\, T(- \triangleright =) \xrightarrow{T_{\mathsf{a}}} - \triangleright T(=) \xrightarrow{\mathrm{id} \triangleright \mathrm{act}} - \triangleright = \,\}$$

$$\{\, U^T F^T(- \triangleright =) \xrightarrow{U^T_{\mathsf{a}} \circ U^T F^T_{\mathsf{a}}} - \triangleright U^T F^T(=) \,\} \longleftarrow \{\, U^T(- \triangleright =) \xrightarrow{U^T_{\mathsf{a}}} - \triangleright U^T(=) \,\}$$

Theorem 5.28 also yields a description of a comodule monad coaction in terms of its Eilenberg–Moore adjunction; i.e., a Tannaka–Krein reconstruction result in the spirit of Proposition 5.9.

**Theorem 5.31.** *Let $B$ be a bimonad on the monoidal category $\mathscr{C}$, and $T$ a monad on a right $\mathscr{C}$-module category $\mathcal{M}$. Coactions of $B$ on $T$ are in bijection with right actions of $\mathscr{C}^B$ on $\mathcal{M}^T$ such that $U^T$ is a strict comodule functor over $U^B$.*

*Proof.* Suppose $\mathscr{C}^B$ acts from the right on $\mathcal{M}^T$ such that $U^T$ is a strict comodule functor. Due to Theorem 5.28, $T = U^T F^T$ is a comodule monad via the coaction

$$T_{\mathsf{a}}\colon T(- \triangleleft =) \xrightarrow{U^T F^T_{\mathsf{a}}} U^T(F^T(-) \triangleleft F^B(=)) \xrightarrow{U^T_{\mathsf{a}}} T(-) \triangleleft B(=),$$

which is equal to $U^T F^T_{\mathsf{a}}$, as $U^T$ is a strict $U^B$-comodule functor.

Conversely, if $T$ is a comodule monad, then $\mathcal{M}^T$ becomes a right $\mathcal{C}^B$-module category, with action given as in Remark 5.24.

As $T_{\mathsf{a}}$ and the action of $\mathcal{C}^B$ on $\mathcal{M}^T$ determine the coaction $F^T_{\mathsf{a}}$ of $F^T$ uniquely, the two constructions are inverse to each other by Theorem 5.28. □

The following result is dual to Theorem 5.28 for the special case of $\mathcal{C}$-module adjunctions, where we study strong $\mathcal{C}$-module structures on the left rather than the right adjoint.

**Proposition 5.32.** *Let $\mathcal{C}$ be a monoidal category, and suppose that $\mathcal{M}$ and $\mathcal{N}$ are left $\mathcal{C}$-module categories. Given an adjunction $F\colon \mathcal{M} \rightleftarrows \mathcal{N} : U$, there is a one-to-one correspondence between lifts of $F$ to a strong $\mathcal{C}$-module functor, and lifts of $F \dashv U$ to a lax $\mathcal{C}$-module adjunction.*

*Proof.* Suppose that $F\colon \mathcal{M} \rightleftarrows \mathcal{N} : U$ is a lax $\mathcal{C}$-module adjunction. Define the inverse of the action $F_{\mathsf{a}}\colon - \triangleright F(=) \Longrightarrow F(- \triangleright =)$ by

$$F(- \triangleright =) \xrightarrow{F(\eta \triangleright =)} F(UF(-) \triangleright =) \xrightarrow{FU_{\mathsf{a}}} FU(F(-) \triangleright =) \xrightarrow{\varepsilon} F(-) \triangleright =.$$

Conversely, given an adjunction $F\colon \mathcal{M} \rightleftarrows \mathcal{N} : U$ such that $F$ is a strong $\mathcal{C}$-module functor, define

$$U(-) \triangleright = \xrightarrow{\eta} UF(U(-) \triangleright =) \xrightarrow{UF_{\mathsf{a}}^{-1}} U(FU(-) \triangleright =) \xrightarrow{U(\varepsilon \triangleright =)} U(- \triangleright =).$$

Reading the string diagrams in the proof of Theorem 5.28 upside down verifies the necessary the coherence conditions, and that these two constructions are inverses of each other. □

We obtain a version of Corollary 5.30 for the Kleisli category of a monad.

**Corollary 5.33.** *Let $\mathcal{M}$ be a left $\mathcal{C}$-module category, and suppose $T$ to be a monad on $\mathcal{M}$. There is a bijection between lax $\mathcal{C}$-module monad structures on $T$, and $\mathcal{C}$-module category structures on $\mathcal{M}_T$, such that $F_T$ is a strict $\mathcal{C}$-module functor.*

*Proof.* Let $\mathcal{M}_T$ be a $\mathcal{C}$-module category such that $F_T$ is a strict $\mathcal{C}$-module functor. Then $T \coloneqq U_T F_T\colon \mathcal{M} \longrightarrow \mathcal{M}$ is a lax $\mathcal{C}$-module monad by Proposition 5.32.

Conversely, if $T$ is a lax $\mathcal{C}$-module monad, $\mathcal{M}_T$ becomes a $\mathcal{C}$-module category as follows: for $x \in \mathcal{C}$ and $m, n \in \mathcal{M}$, set $x \triangleright_{\mathcal{M}_T} m \coloneqq x \triangleright_{\mathcal{M}} m$ on objects, and for a morphism $f\colon m \longrightarrow Tn$ in $\mathcal{M}_T$, we define the action of $x \in \mathcal{C}$ by

$$x \triangleright_{\mathcal{M}_T} f \coloneqq x \triangleright m \xrightarrow{x \triangleright f} x \triangleright Tn \xrightarrow{T_{\mathsf{a};x,n}} T(x \triangleright n).$$

It is easy to check that these assignments define a $\mathcal{C}$-module structure, for which $F_T$ is a strict $\mathcal{C}$-module functor.

The $\mathcal{C}$-module structure of $\mathcal{M}_T$ and the lax module monad structure on $T$ uniquely determine the $\mathcal{C}$-module structure of $F_T$. Hence, $F^{-1}_{T;\mathsf{a}}$ is given as in Proposition 5.32 and the two constructions are inverse of each other. □

Using these insights, we may now clarify the structure of comparison functors associated to comodule adjunctions.

**Proposition 5.34.** *Consider a comodule adjunction $G\colon \mathcal{M} \rightleftarrows \mathcal{N} : V$ over an oplax monoidal adjunction $F\colon \mathcal{C} \rightleftarrows \mathcal{D} : U$, and denote the associated comodule monad and bimonad by $C := VG\colon \mathcal{M} \longrightarrow \mathcal{M}$ and $B := UF\colon \mathcal{C} \longrightarrow \mathcal{C}$, respectively.*

*Then the comparison functor $K^C\colon \mathcal{N} \longrightarrow \mathcal{M}^C$ is a strong comodule functor over $K^B\colon \mathcal{D} \longrightarrow \mathcal{C}^B$, and the following identities of comodule functors hold:*

$$U^C K^C = V \qquad \textit{and} \qquad K^C G = F^C.$$

*Proof.* We proceed analogously to [BV07, Theorem 2.6]. For any $n \in \mathcal{N}$ we have $K^C n = (Vn, V\varepsilon_n)$ and a direct computation shows that the coaction of $V$ lifts to a coaction of $K^C$. That is,

$$U^C K^C_{\mathsf{a};n,y} = V_{\mathsf{a};n,y}, \qquad\qquad \text{for all } n \in \mathcal{N} \text{ and } y \in \mathcal{D}.$$

Using that $U^C\colon \mathcal{M}^C \longrightarrow \mathcal{M}$ is a faithful and conservative functor, one observes that $K^C$ becomes a strong comodule functor in this manner. Furthermore, as $U^C$ is a strict comodule functor, the coactions of $U^C K^C$ and $V$ coincide. Lastly, we compute the following, for any $x \in \mathcal{C}$ and $m \in \mathcal{M}$:

$$(K^C G)_{\mathsf{a};m,x} = (U^C K^C G)_{\mathsf{a};m,x} = (VG)_{\mathsf{a};m,x} = C_{\mathsf{a};m,x} = (U^C F^C)_{\mathsf{a};m,x} = F^C_{\mathsf{a};m,x}. \quad \square$$

**Proposition 5.35.** *Let $\mathcal{M}$ and $\mathcal{N}$ be left $\mathcal{C}$-module categories.*

- *Let $F\colon \mathcal{M} \rightleftarrows \mathcal{N} : U$ be an oplax $\mathcal{C}$-module adjunction and consider the corresponding $\mathcal{C}$-module category structure on $\mathcal{M}^T$ of Corollary 5.30. Then for $T := UF$ the comparison functor $K^T$ is a strong $\mathcal{C}$-module functor.*
- *Let $F\colon \mathcal{M} \rightleftarrows \mathcal{N} : U$ be a lax $\mathcal{C}$-module adjunctions and consider the corresponding $\mathcal{C}$-module category structure on $\mathcal{M}_T$ of Corollary 5.33. Then for $T := UF$ the comparison functor $K_T$ is a strong $\mathcal{C}$-module functor.*

**Proposition 5.36.** *Let $\mathcal{M}$ be a $\mathscr{C}$-module category, $T\colon \mathcal{M} \longrightarrow \mathcal{M}$ an oplax $\mathscr{C}$-module monad, and $S\colon \mathcal{M} \longrightarrow \mathcal{M}$ a right adjoint to $T$. Then the isomorphism $\mathcal{M}^T \cong \mathcal{M}^S$ of Proposition 2.25 is a strict $\mathscr{C}$-module isomorphism, where $\mathcal{M}^T$ is endowed with the $\mathscr{C}$-module structure of Proposition 5.30, and similarly for $\mathcal{M}^S$.*

*Proof.* Note that, by Theorem 5.28 the comonad $S\colon \mathcal{M} \longrightarrow \mathcal{M}$ comes equipped with the following lax $\mathscr{C}$-module structure, for all $x \in \mathscr{C}$ and $m \in \mathcal{M}$:

$$S_{\mathsf{a};x,m}\colon x \triangleright Sm \xrightarrow{\eta_{x\triangleright Sm}} ST(x \triangleright Sm) \xrightarrow{ST_{\mathsf{a};x,Sm}} S(x \triangleright TSm) \xrightarrow{S(x\triangleright\varepsilon_m)} S(x \triangleright m).$$

Further, the isomorphism $L\colon \mathcal{M}^T \xrightarrow{\sim} \mathcal{M}^S$ of Proposition 2.25 is given by

$$\big(m,\ \nabla_m\colon Tm \longrightarrow m\big) \longmapsto \big(m,\ m \xrightarrow{\eta_m} STm \xrightarrow{S\nabla_m} Sm\big).$$

We have to prove that $L(x \triangleright (m, \nabla_m)) = x \triangleright L((m, \nabla_m))$, for all $x \in \mathscr{C}$ and $(m, \nabla_m) \in \mathcal{M}^T$. This is equivalent to the equality of

$$\big(x \triangleright m,\ x \triangleright m \xrightarrow{\eta_{x\triangleright m}} ST(x \triangleright m) \xrightarrow{ST_{\mathsf{a};x,m}} S(x \triangleright Tm) \xrightarrow{S(x\triangleright\nabla_m)} S(x \triangleright m)\big)$$

and

$$\big(x \triangleright m,\ x \triangleright m \xrightarrow{x\triangleright\eta_m} x \triangleright STm \xrightarrow{x\triangleright S\nabla_m} x \triangleright Sm \xrightarrow{\eta_{x\triangleright Sm}} ST(x \triangleright Sm)$$
$$\xrightarrow{S(T_{\mathsf{a}})_{x,S(m)}} S(x \triangleright TSm) \xrightarrow{S(x\triangleright\varepsilon_m)} S(x \triangleright m)\big).$$

This is evidenced by the following commutative diagram:

$$\begin{array}{ccccccc}
x \triangleright m & \xrightarrow{\eta_{x\triangleright m}} & ST(x \triangleright m) & & \xrightarrow{ST_{\mathsf{a};x,m}} & & S(x \triangleright Tm) \\
\downarrow{\scriptstyle x\triangleright\eta_m} & \text{nat } \eta & \downarrow{\scriptstyle ST(x\triangleright\eta_m)} & \text{nat } T_{\mathsf{a}} & \swarrow{\scriptstyle S(x\triangleright T\eta_m)} & \text{adj} & \Vert \\
x \triangleright STm & \xrightarrow[\eta_{x\triangleright STm}]{} & ST(x \triangleright STm) & \xrightarrow{ST_{\mathsf{a};x,STm}} & S(x \triangleright TSTm) & \xrightarrow[S(x\triangleright\varepsilon_{Tm})]{} & S(x \triangleright Tm) \\
\downarrow{\scriptstyle x\triangleright S\nabla_m} & \text{nat } \eta & \downarrow{\scriptstyle ST(x\triangleright S\nabla_m)} & \text{nat } T_{\mathsf{a}} & \downarrow{\scriptstyle S(x\triangleright TS\nabla_m)} & \text{nat } \varepsilon & \downarrow{\scriptstyle S(x\triangleright\nabla_m)} \\
x \triangleright Sm & \xrightarrow[\eta_{x\triangleright Sm}]{} & ST(x \triangleright Sm) & \xrightarrow[ST_{\mathsf{a};x,Sm}]{} & S(x \triangleright TSm) & \xrightarrow[S(x\triangleright\varepsilon_m)]{} & S(x \triangleright m)
\end{array}$$

□

In an analogous way to Proposition 5.36, one obtains the following result.

**Proposition 5.37.** *Let $\mathcal{M}$ be a $\mathscr{C}$-module category, $T\colon \mathcal{M} \longrightarrow \mathcal{M}$ a lax $\mathscr{C}$-module monad, and $L\colon \mathcal{M} \longrightarrow \mathcal{M}$ a left adjoint to $T$. Then the isomorphism $\mathcal{M}_T \cong \mathcal{M}_L$ of Proposition 2.26 is a $\mathscr{C}$-module isomorphism, where $\mathcal{M}_T$ and $\mathcal{M}_L$ are endowed with the $\mathscr{C}$-module structure of Corollary 5.33.*

[Der Wert dieser Arbeit] wird umso grösser sein, je besser die Gedanken ausgedrückt sind. Je mehr der Nagel auf den Kopf getroffen ist.—Hier bin ich mir bewusst, weit hinter dem Möglichen zurückgeblieben zu sein. Einfach darum, weil meine Kraft zur Bewältigung der Aufgabe zu gering ist.

LUDWIG WITTGENSTEIN; Tractatus Logico-Philosophicus

# 6 MONADIC TWISTED CENTRES

THE ANTI-YETTER–DRINFELD MODULES of a finite-dimensional Hopf algebra are a module category over the Yetter–Drinfeld modules. Subsequently, they are implemented by a comodule algebra over the Drinfeld double, see [HKRS04]. As explained in Section 4.2, we find ourselves in a similar situation. The anti-Drinfeld centre, our replacement of the anti-Yetter–Drinfeld modules, is a module category over the Drinfeld centre.

Replacing finite-dimensional vector spaces by a rigid, possibly pivotal, category $\mathscr{C}$, and the underlying Hopf algebra with a Hopf monad $H$ on $\mathscr{C}$, this section serves to study a Hopf monad $D(H)\colon \mathscr{C} \longrightarrow \mathscr{C}$ and over it a comodule monad $Q(H)\colon \mathscr{C} \longrightarrow \mathscr{C}$, which realise the centre and its twisted cousin as their respective modules. Bruguières and Virelizier gave a transparent description of $D(H)$ in [BV12] by extending results of Day and Street, [DS07]. If $H$ is the identity functor, one defines the *central Hopf monad* $\mathfrak{D}(\mathscr{C}^H)$ on $\mathscr{C}^H$, with $\mathsf{Z}(\mathscr{C}^H)$ as its Eilenberg–Moore category. As an application of Beck's theory of distributive laws, one obtains the Drinfeld double $D(H)\colon \mathscr{C} \longrightarrow \mathscr{C}$. We apply the same techniques to define the anti-double $Q(H)$ of $H$, whose modules are isomorphic to the "dual" of the anti-Drinfeld centre $\mathsf{Q}(\mathscr{C}^H)$. This approach is best summarised by Figure 6.1. We obtain a monadic version of Theorem 4.1.

**Theorem 6.44.** *Let $\mathscr{C}$ be a rigid monoidal category, and suppose that $H\colon \mathscr{C} \longrightarrow \mathscr{C}$ is a Hopf monad that admits a double $D(H)$ and anti-double $Q(H)$. The following statements are equivalent*:

(*i*) *the monoidal unit $1 \in \mathscr{C}$ lifts to a module over $Q(H)$,*
(*ii*) *there is an isomorphism of comodule monads $D(H) \cong Q(H)$, and*
(*iii*) *there is an isomorphism of monads $Q(H) \cong D(H)$.*

*If $\mathscr{C}$ is pivotal with pivotal structure $\phi$, any of the above statements hold if and only if $H$ and $\phi$ admit a pair in involution.*

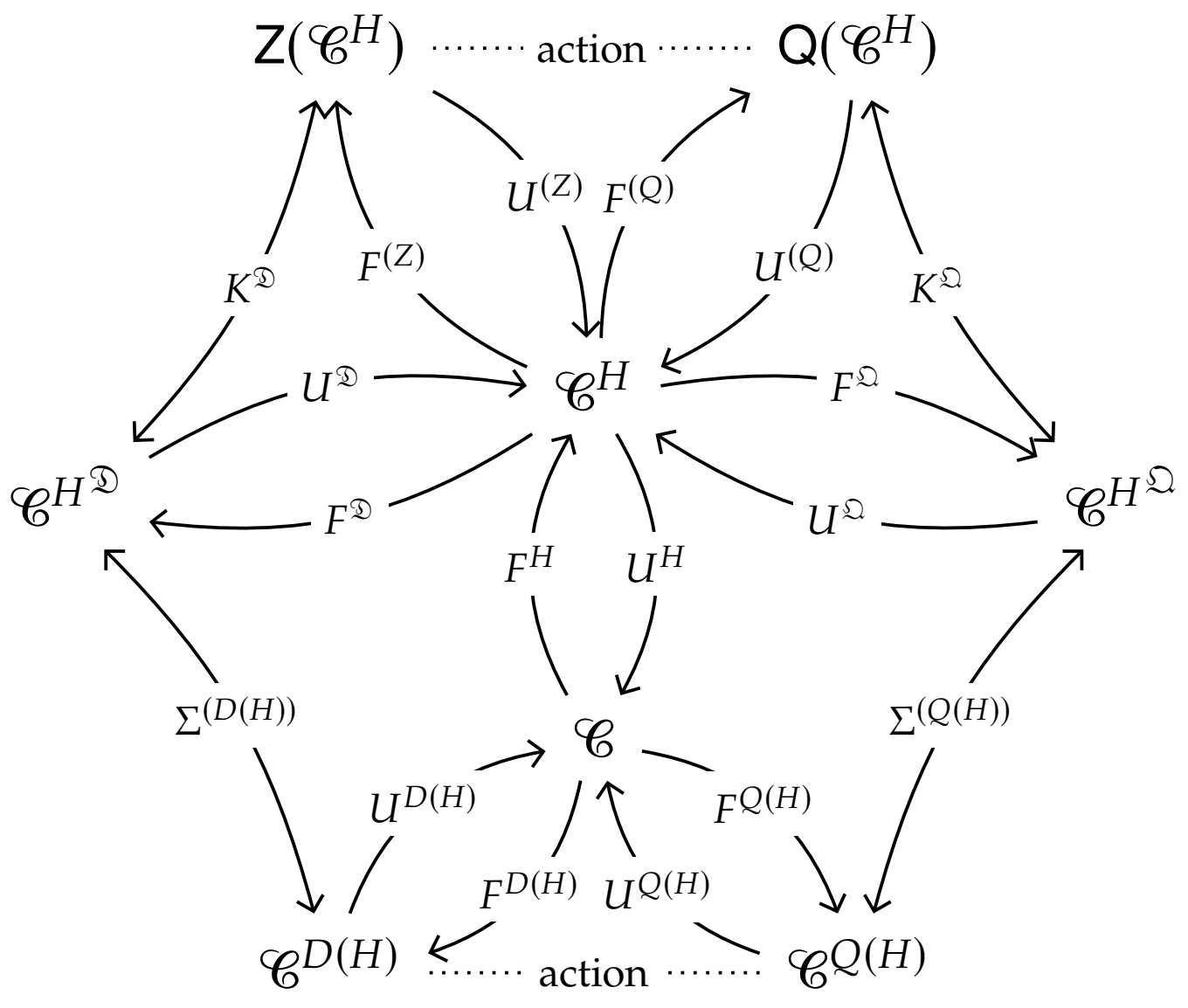

Figure 6.1: A cobweb of adjunctions, monads, and various versions of the centre and anti-centre.

From Theorem 6.44 we can deduce how pivotal structures on $\mathscr{C}^H$ arise from module morphisms between the central Hopf monad $\mathfrak{D}$ and the anti-central comodule monad $\mathfrak{Q}$.

**Corollary 6.45.** *Let $\mathscr{C}$ be a rigid monoidal category. If $\mathscr{C}$ admits a central Hopf monad $\mathfrak{D}(\mathscr{C})$ and an anti-central comodule monad $\mathfrak{Q}(\mathscr{C})$, then it is pivotal if and only if $\mathfrak{D}(\mathscr{C}) \cong \mathfrak{Q}(\mathscr{C})$ as monads.*

## 6.1 Cross products and distributive laws

The Hopf monadic description of the Drinfeld centre $\mathsf{Z}(\mathscr{C}^H)$, for a rigid category $\mathscr{C}$ and Hopf monad $H$ on $\mathscr{C}$, is achieved as a two-step process: one first finds a suitable monad on $\mathscr{C}^H$, which is then "lifted" to a monad on $\mathscr{C}$. We shall review this lifting process based on [BV12, Sections 3 and 4].

**Definition 6.1.** Let $H\colon \mathscr{C} \longrightarrow \mathscr{C}$ be a monad and $T\colon \mathscr{C}^H \longrightarrow \mathscr{C}^H$ a functor. The *cross product* $T \rtimes H$ of $T$ by $H$ is the endofunctor $U^H T F^H\colon \mathscr{C} \longrightarrow \mathscr{C}$.

If $(T, \mu, \eta)$ is a monad, then the cross product $T \rtimes H$ inherits this structure: the multiplication and unit are given by

Here $\varepsilon$ and $\eta$ are the unit and counit of the adjunction $F^H \dashv U^H$.

$$\mu^{(T\rtimes H)}\colon U^H T F^H U^H T F^H \xrightarrow{U^H T \varepsilon T F^H} U^H T T F^H \xrightarrow{U^H \mu^{(T)} F^H} U^H T F^H$$

and

$$\eta^{(T\rtimes H)}\colon \mathrm{Id}_{\mathscr{C}} \xrightarrow{\eta} U^H F^H \xrightarrow{U^H \eta^{(T)}} U^H T F^H.$$

Given two bimonads $H\colon \mathscr{C} \longrightarrow \mathscr{C}$ and $B\colon \mathscr{C}^H \longrightarrow \mathscr{C}^H$, iteratively applying the comultiplication and counit of $U^H$, $B$, and $F^H$ yields a bimonad structure on $B \rtimes H\colon \mathscr{C} \longrightarrow \mathscr{C}$. The comultiplication is given by

$$U^H_{2;TF^H(-),TF^H(=)} \circ U^H T_{2;F^H(-),F^H(=)} \circ U^H TF^H_{2;-,=},$$

and for the counit we have

$$\varepsilon^{(U^H)} \circ U^H \varepsilon^{(T)} \circ U^H T \varepsilon^{(F^H)}.$$

Similar considerations imply the following result.

**Lemma 6.2.** *Let $H\colon \mathscr{C} \longrightarrow \mathscr{C}$ and $B\colon \mathscr{C}^H \longrightarrow \mathscr{C}^H$ be bimonads which respectively coact on the comodule monads $G\colon \mathscr{M} \longrightarrow \mathscr{M}$ and $C\colon \mathscr{M}^G \longrightarrow \mathscr{M}^G$. The cross product $C \rtimes G\colon \mathscr{M} \longrightarrow \mathscr{M}$ is a comodule monad over $B \rtimes H$ via the coaction*

$$U^G_{\mathsf{a};CF^G(-),BF^H(=)} \circ U^G C_{\mathsf{a};F^G(-),F^H(=)} \circ U^G CF^G_{\mathsf{a};-,=}.$$

**Remark 6.3.** Let $H\colon \mathscr{C} \longrightarrow \mathscr{C}$ and $B\colon \mathscr{C}^H \longrightarrow \mathscr{C}^H$ be monads. The question under which conditions the modules $\mathscr{C}^{B\rtimes H}$ of $B\rtimes H$ are isomorphic to $(\mathscr{C}^H)^B$ is closely related to the theory of distributive laws of Section 2.2.2. We may apply the formal theory thereof, see [Str72], to the bicategory $\mathbb{O}\mathsf{plMon}^{\mathsf{opl}}$ of monoidal categories, oplax monoidal functors, and oplax monoidal natural transformations, to obtain a description of bimonads and *oplax monoidal distributive laws*—oplax monoidal natural transformations $\Lambda\colon HB \longrightarrow BH$ between bimonads $H, B\colon \mathscr{C} \longrightarrow \mathscr{C}$ that are moreover distributive laws, see [McC02].

Suppose $\Lambda\colon HB \longrightarrow BH$ to be an oplax monoidal distributive law. The comultiplication and counit of the underlying functor $BH\colon \mathscr{C} \longrightarrow \mathscr{C}$ turn the composite $B \circ_\Lambda H$ into a bimonad. Thus, comodule monads can be intrinsically described in the bicategory that has

- as objects, pairs $(\mathscr{M}, \mathscr{C})$ consisting of a right module category $\mathscr{M}$ over a monoidal category $\mathscr{C}$;
- as 1-morphisms, pairs $(G, F)$ of a comodule functor $G$ over an oplax monoidal functor $F$; and
- as 2-morphisms, pairs $(\phi, \psi)$ that form a comodule transformation.

The subsequent results arise immediately from [Str72].

**Definition 6.4.** Let $G, C\colon \mathcal{M} \longrightarrow \mathcal{M}$ be two comodule monads over the bimonads $H, B\colon \mathcal{C} \longrightarrow \mathcal{C}$, respectively. A *comodule distributive law* is a pair of distributive laws $\Omega\colon GC \longrightarrow CG$ and $\Lambda\colon HB \longrightarrow BH$ such that $(\Lambda, \Omega)$ is a comodule natural transformation.

**Proposition 6.5.** *Consider two comodule monads $G, C\colon \mathcal{M} \longrightarrow \mathcal{M}$ over the bimonads $H, B\colon \mathcal{C} \longrightarrow \mathcal{C}$. There exists a bijective correspondence between:*

(*i*) *comodule distributive laws $\big(GC \xrightarrow{\Omega} CG, HB \xrightarrow{\Lambda} BH\big)$; and*
(*ii*) *lifts of $B$ to a bimonad $\widetilde{B}\colon \mathcal{C}^H \longrightarrow \mathcal{C}^H$ together with lifts of $C$ to a comodule monad $\widetilde{C}\colon \mathcal{M}^G \longrightarrow \mathcal{M}^G$ over $\widetilde{B}$, such that $BU^H = U^H\widetilde{B}$ as oplax monoidal functors and $CU^G = U^G\widetilde{C}$ as comodule functors.*

Let $\big(GC \xrightarrow{\Omega} CG, HB \xrightarrow{\Lambda} BH\big)$ be a comodule distributive law. The coactions of $G$ and $C$ turn $C \circ_\Omega G$ into a comodule monad over $B \circ_\Lambda H$.

**Lemma 6.6.** *Suppose $\Omega\colon GC \longrightarrow CG$ and $\Lambda\colon HB \longrightarrow BH$ to form a comodule distributive law. Then there are equivalences $(\mathcal{C}^H)^{\widetilde{B}^\Lambda} \simeq \mathcal{C}^{B\circ_\Lambda H}$ as monoidal categories, and $(\mathcal{M}^G)^{\widetilde{C}^\Omega} \simeq \mathcal{M}^{C\circ_\Omega G}$ as $\mathcal{C}^{B\circ_\Lambda H}$-module categories.*

In fact, by [BV12, Section 4.5], the previous results transcend to the Hopf monadic setting. Let $B, H\colon \mathcal{C} \longrightarrow \mathcal{C}$ be Hopf monads. If $\Lambda\colon HB \longrightarrow BH$ is an oplax monoidal distributive law, then $B \circ_\Lambda H\colon \mathcal{C} \longrightarrow \mathcal{C}$ and the lift $\widetilde{B}^\Lambda\colon \mathcal{C}^H \longrightarrow \mathcal{C}^H$ are Hopf monads as well.

## 6.2 Centralisable functors and the central bimonad

The construction of the double of a Hopf monad $H\colon \mathcal{C} \longrightarrow \mathcal{C}$ given in [BV12] relies on explicitly describing the left dual of the forgetful functor $U^{(Z)}\colon \mathsf{Z}(\mathcal{C}^H) \longrightarrow \mathcal{C}^H$, which is realised as a coend.

**Definition 6.7.** Let $\mathcal{C}$ be a rigid category and $T\colon \mathcal{C} \longrightarrow \mathcal{C}$ an endofunctor. We call $T$ *centralisable* if the following coend exists for all $x \in \mathcal{C}$:

$$Z_T(x) := \int^{y\in\mathcal{C}} {}^{\vee}Ty \otimes x \otimes y.$$

**Remark 6.8.** Any centralisable functor $T\colon \mathcal{C} \longrightarrow \mathcal{C}$ admits a *universal coaction*

$$\chi_{x,y} := (\mathrm{id}_{Ty} \otimes \mathrm{copr}_{y,x}) \circ (\mathrm{coev}^l_{Ty} \otimes \mathrm{id}_{x\otimes y}), \qquad \text{for } x, y \in \mathcal{C},$$

which is natural in both variables. We call the pair $(Z_T, \chi)$ a *centraliser* of $T$.

**Remark 6.9.** Graphically, we represent the universal coaction as follows:

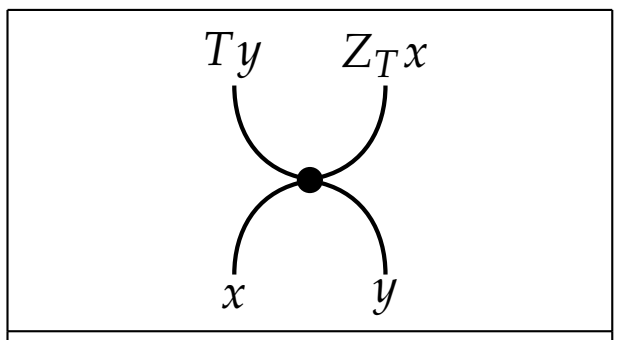

$\chi_{x,y} \colon x \otimes y \longrightarrow Ty \otimes Z_T x$

For all $f \colon x \longrightarrow x'$ and $g \colon y \longrightarrow y'$, the naturality condition equates to

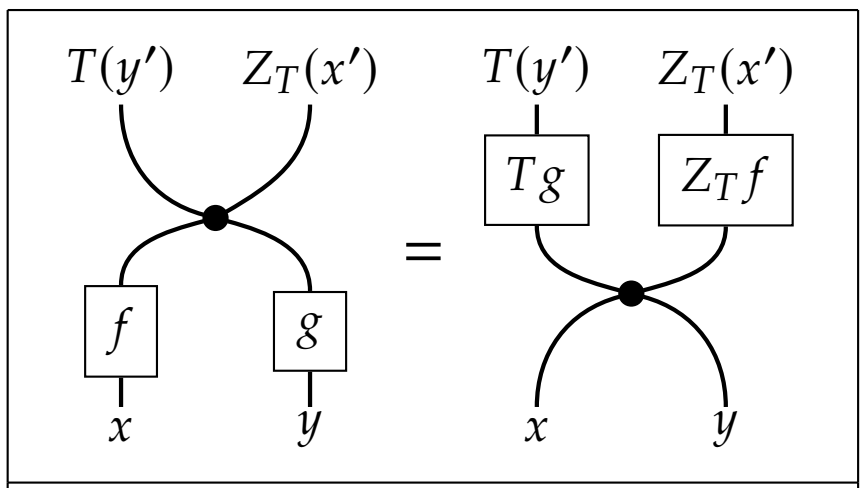

$\chi_{x',y'} \circ (f \otimes g) = (Tg \otimes Z_T f) \circ \chi_{x,y}$

Universal coactions have a certain universal property that will be vital in the rest of this chapter—the *extended factorisation property*. In particular, it provides us with a potent tool for constructing bi- and comodule monads.

**Lemma 6.10** ([BV12, Lemma 5.4])**.** *Let $(Z_T, \chi)$ be the centraliser of a functor $T \colon \mathscr{C} \longrightarrow \mathscr{C}$ and suppose that $L, R \colon \mathscr{D} \longrightarrow \mathscr{C}$ are two functors. For any $n \in \mathbb{N}$, $y \in \mathscr{D}$, and $x_1, \ldots, x_n \in \mathscr{C}$, and any natural transformation*

$$\phi_{y,x_1,\ldots,x_n} \colon Ly \otimes x_1 \otimes \cdots \otimes x_n \longrightarrow Tx_1 \otimes \cdots \otimes Tx_n \otimes Ry,$$

*there exists a unique natural transformation $\nu \colon Z_T^n L \Longrightarrow R$, satisfying*

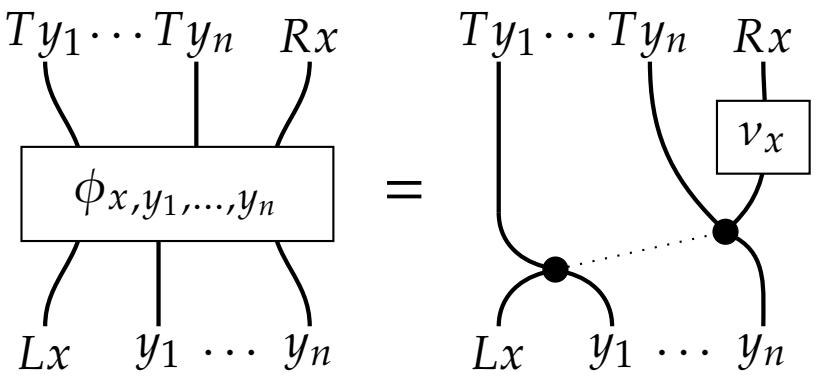

**Example 6.11.** Let $T \colon \mathscr{C} \longrightarrow \mathscr{C}$ be an oplax monoidal functor with centraliser $(Z_T, \chi)$. For all $x \in \mathscr{C}$, define the unit of $Z_T$ by

$$\eta_x^{(Z_T)} \colon x \xrightarrow{\chi_{x,1}} T1 \otimes Z_T x \xrightarrow{T_0 \otimes Z_T x} Z_T x.$$

By Lemma 6.10, we can derive a unique multiplication $\mu^{(Z_T)}\colon Z_T^2 \Longrightarrow Z_T$ from the comultiplication of $T$.

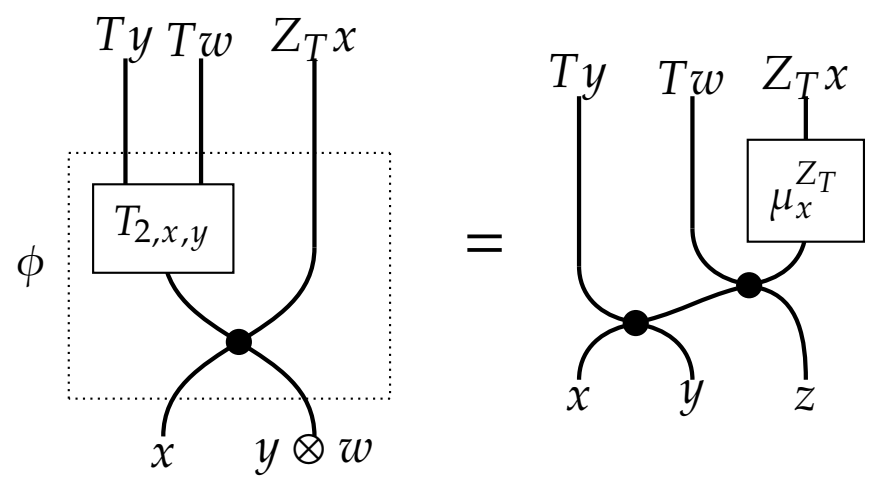

The following result is due to Day and Street, [DS07, pp. 191–192], see also [BV12, Theorem 5.6] for a proof.

**Lemma 6.12.** *The centraliser $(Z_T, \chi)$ of an oplax monoidal endofunctor $T$ on $\mathscr{C}$ is a monad with multiplication and unit as given in Example 6.11.*

In the proof of Lemma 6.12 given in [BV12, Theorem 5.6], the authors further consider $T\colon \mathscr{C} \longrightarrow \mathscr{C}$ to be equipped with a Hopf monad structure and show that in this case $Z_T$ is a Hopf monad as well. The extended factorisation property given in Lemma 6.10 reconstructs a comultiplication on $Z_T$ from a twofold application of the universal coaction and the multiplication of $T$:

(6.2.1)

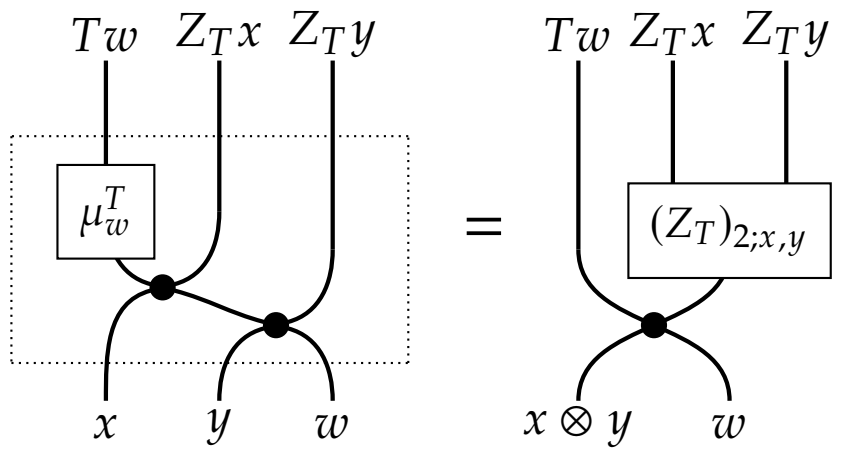

Likewise, the unit of $T$ induces a counit on $Z_T$ via

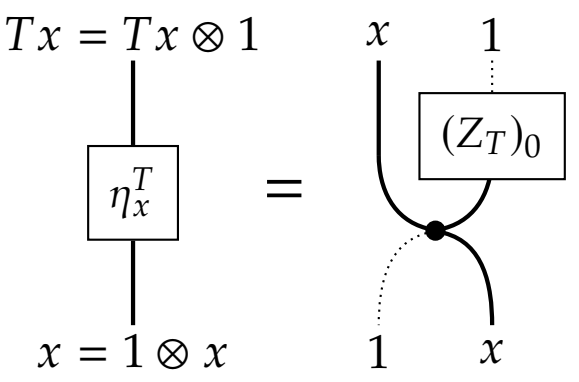

A direct computation verifies that the centraliser $Z_T$ is a bimonad as well. For the construction of left and right antipodes, see [BV12, Theorem 5.6].

**Remark 6.13.** We think of $\mathsf{Z}({}_H\mathscr{C})$ as the centre of an oplax bimodule category as stated in Remark 2.46, see also [BV07, Section 5.5]. Objects in $\mathsf{Z}({}_H\mathscr{C})$ are pairs $(x, \sigma_{x,-})$, where $x \in \mathscr{C}$ and $\sigma_{x,-}\colon x \otimes - \Longrightarrow H(-) \otimes x$ is a natural transformation, satisfying for all $x, y, z \in \mathscr{C}$

$$(H_{2;y,z} \otimes x) \circ \sigma_{x,y\otimes z} = (Hy \otimes \sigma_{x,z}) \circ (\sigma_{x,y} \otimes z),$$
$$(H_0 \otimes x) \circ \sigma_{x,1} = \mathrm{id}_x.$$

Analogous to the centres studied before, the arrows in $\mathsf{Z}({}_H\mathscr{C})$ are those morphisms of $\mathscr{C}$ that commute with the half-braidings. As shown in [BV12, Proposition 5.9], the structure morphisms of a Hopf monad $H\colon \mathscr{C} \longrightarrow \mathscr{C}$ can be used to define a rigid structure on $\mathsf{Z}({}_H\mathscr{C})$. For example, the tensor product of $(x, \sigma_{x,-}), (y, \sigma_{y,-}) \in \mathsf{Z}({}_H\mathscr{C})$ is $x \otimes y \in \mathscr{C}$, together with the half-braiding

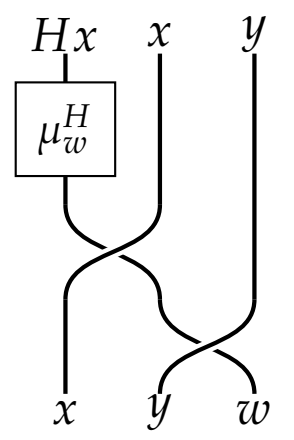

For a treatment of centres twisted by (op)lax monoidal functors, see [FH23].

Since centralisers of Hopf monads are Hopf monads themselves, their modules implement the twisted centres of Remark 6.13 as a rigid category. This is proven in [BV12, Theorem 5.12 and Corollary 5.14].

**Proposition 6.14.** *Suppose $H\colon \mathscr{C} \longrightarrow \mathscr{C}$ to be a centralisable Hopf monad. The modules $\mathscr{C}^{Z_H}$ of its centraliser $(Z_H, \chi)$ are isomorphic as a rigid category to $\mathsf{Z}({}_H\mathscr{C})$.*

Applying the above proposition to the identity functor $\mathrm{Id}\colon \mathscr{C} \longrightarrow \mathscr{C}$, we obtain a Hopf monadic description of the Drinfeld centre $\mathsf{Z}(\mathscr{C})$ of a rigid category $\mathscr{C}$. The terminology of our next definition is due to Shimizu, see [Shi17].

**Definition 6.15.** Let $\mathscr{C}$ be a monoidal category, and let the identity functor on $\mathscr{C}$ be centralisable with centraliser $(Z, \chi)$. The *central Hopf monad of* $\mathscr{C}$ is $\mathfrak{D}(\mathscr{C}) := Z_{\mathrm{Id}}\colon \mathscr{C} \longrightarrow \mathscr{C}$. We denote its Eilenberg–Moore by $\mathscr{C}^{\mathfrak{D}}$.

An important step in proving Proposition 6.14 is determining an inverse to the comparison functor $K^{Z_T}\colon \mathsf{Z}({}_T\mathscr{C}) \longrightarrow \mathscr{C}^{Z_T}$. This construction will also play a substantial role in our monadic description of the anti-Drinfeld centre, so we recall it in its full generality. Let $T\colon \mathscr{C} \longrightarrow \mathscr{C}$ be a centralisable oplax monoidal endofunctor with $(Z_T, \chi)$ as its centraliser. To every $Z_T$-module $(m, \nabla_m)$, we associate a half-braiding $\sigma_{m,-}\colon m \otimes - \Longrightarrow T(-) \otimes m$. For any $x \in \mathscr{C}$, it is given by the composition

(6.2.2) $$\sigma_{m,x}\colon m \otimes x \xrightarrow{\chi_{m,x}} Tx \otimes Z_T x \xrightarrow{Tx \otimes \nabla_m} Tx \otimes m.$$

This yields a functor $E^{Z_T}\colon \mathscr{C}^{Z_T} \longrightarrow \mathsf{Z}({}_T\mathscr{C})$, which is the identity on morphisms and on objects is given by

(6.2.3) $$E^{Z_T}(m, \nabla_m) = (m, \sigma_{m,-}), \qquad \text{for all } (m, \nabla_M) \in \mathscr{C}^{Z_T}.$$

Conversely, to every object $(m, \sigma_{m,-}) \in \mathsf{Z}({}_T\mathscr{C})$ we may assign a $Z_T$-module, whose action $\nabla_m$ is, due to Lemma 6.10, uniquely defined by

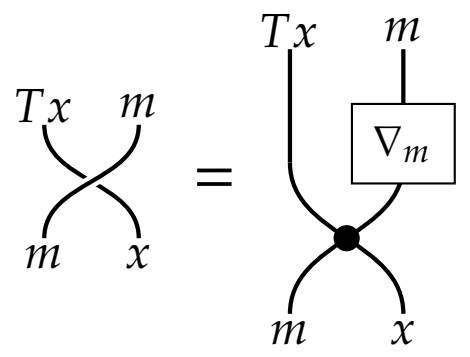

This yields the comparison functor $K^{Z_T}\colon \mathsf{Z}({}_T\mathscr{C}) \longrightarrow \mathscr{C}^{Z_T}$.

**Remark 6.16.** Suppose $T\colon \mathscr{C} \longrightarrow \mathscr{C}$ to be a centralisable oplax monoidal endofunctor with $(Z_T, \chi)$ as its centraliser. Denote the free functor of the Eilenberg–Moore adjunction of $Z_T$ by $F^{Z_T}\colon \mathscr{C} \longrightarrow \mathscr{C}^{Z_T}$. The composition

$$\mathscr{C} \xrightarrow{F^{Z_T}} \mathscr{C}^{Z_T} \xrightarrow{E^{Z_T}} \mathsf{Z}({}_T\mathscr{C})$$

defines a left adjoint of the forgetful functor $U^{(T)}\colon \mathsf{Z}({}_T\mathscr{C}) \longrightarrow \mathscr{C}$.

In fact, the central adjunction $F^{(T)} \dashv U^{(T)}$ is monadic.

**Proposition 6.17** ([BV12, Theorem 5.12])**.** *Let $(Z_T, \chi)$ be a centraliser of the oplax monoidal endofunctor $T\colon \mathscr{C} \longrightarrow \mathscr{C}$. The comparison functor $K^{Z_T}\colon \mathsf{Z}({}_T\mathscr{C}) \longrightarrow \mathscr{C}^{Z_T}$ is an isomorphism of categories with inverse $E^{Z_T}\colon \mathscr{C}^{Z_T} \longrightarrow \mathsf{Z}({}_T\mathscr{C})$.*

## 6.3 Centralisers and comodule monads

We will now apply the methods of Bruguières and Virelizier to twisted centres, for the purpose of obtaining a comodule monad that implements the anti-Drinfeld centre. Hereto, we need a generalised version of the concept of modules over a monad. Our approach is based on [MW11].

**Definition 6.18.** Let $(B, \mu, \eta)\colon \mathscr{C} \longrightarrow \mathscr{C}$ be a bimonad and $F\colon \mathscr{C} \longrightarrow \mathscr{D}$ an oplax monoidal functor. An *oplax monoidal right action* of $B$ on $F$ is an oplax natural transformation $\alpha\colon FB \Longrightarrow F$, such that the following diagrams commute:

$$\begin{array}{ccc} FBB & \xrightarrow{F\mu} & FB \\ {\scriptstyle \alpha B}\downarrow & & \downarrow{\scriptstyle \alpha} \\ FB & \xrightarrow[\alpha]{} & F \end{array} \qquad \begin{array}{ccc} F & \xrightarrow{F\eta} & FB \\ & {\scriptstyle \mathrm{id}_F}\searrow & \downarrow{\scriptstyle \alpha} \\ & & F \end{array}$$

Similarly, one defines oplax monoidal left actions.

**Remark 6.19.** Recall from Remark 5.2 that a bimonad $B$ is a monoid in the monoidal category $\mathbb{O}\mathsf{plMon}^{\mathsf{opl}}(\mathscr{C}, \mathscr{C})$ of oplax monoidal endofunctors on $\mathscr{C}$, with oplax monoidal natural transformations between them. There is an obvious right action of this category on $\mathbb{O}\mathsf{plMon}^{\mathsf{opl}}(\mathscr{C}, \mathscr{D})$ given by composition of functors. In this context, an oplax monoidal right action of a bimonad $B\colon \mathscr{C} \longrightarrow \mathscr{C}$ is the same as an $\mathbb{O}\mathsf{plMon}^{\mathsf{opl}}(\mathscr{C}, \mathscr{D})$-module of $B$.

A prime example of an oplax monoidal left action is given by the forgetful functor $U^B\colon \mathscr{C}^B \longrightarrow \mathscr{C}$ of a bimonad $B\colon \mathscr{C} \longrightarrow \mathscr{C}$ together with the action $\nabla = U^B\varepsilon\colon BU^B \longrightarrow U^B$, see Remark 2.30.

**Hypothesis 6.20.** To keep our notation concise, in the following we fix an oplax monoidal functor $L\colon \mathscr{C} \longrightarrow \mathscr{C}$ with an oplax right action $\alpha\colon LB \Longrightarrow L$ by a bimonad $B\colon \mathscr{C} \longrightarrow \mathscr{C}$ and assume that $L$ and $B$ are centralisable. Their centralisers will be denoted by $(Z_L, \xi)$ and $(Z_B, \chi)$, respectively.

We think of $\mathsf{Z}({}_B\mathscr{C})$ as a more general version of the Drinfeld centre which is supposed to act on $\mathsf{Z}({}_L\mathscr{C})$ from the right. To emphasise this, and in line with the colouring scheme of Section 4.2, we use black for objects in $\mathscr{C}$ or its generalised Drinfeld centre and blue for objects in $\mathsf{Z}({}_L\mathscr{C})$.

Consider objects $(m, \sigma_{m,-}) \in \mathsf{Z}({}_L\mathscr{C})$ and $(x, \sigma_{x,-}) \in \mathsf{Z}({}_B\mathscr{C})$. The action of $B$ on $L$ and the half-braidings of $m$ and $x$ yield a natural transformation

$$\sigma_{m\otimes x, y}\colon m \otimes x \otimes y \longrightarrow Ly \otimes m \otimes x,$$

which we can depict graphically by

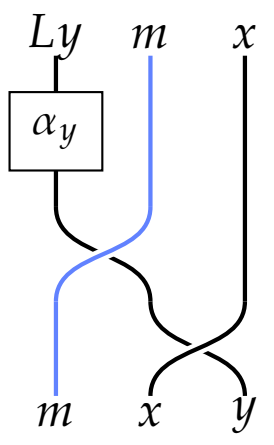

(6.3.1)

**Lemma 6.21.** *The centre* $\mathsf{Z}({}_B\mathscr{C})$ *acts on* $\mathsf{Z}({}_L\mathscr{C})$ *from the right by tensoring the underlying objects and gluing together the half-braidings as in Equation (6.3.1). With respect to this action, the forgetful functor* $U^{(L)}\colon \mathsf{Z}({}_L\mathscr{C}) \longrightarrow \mathscr{C}$ *is a strict comodule functor over* $U^{(B)}\colon \mathsf{Z}({}_B\mathscr{C}) \longrightarrow \mathscr{C}$.

*Proof.* We proceed as in [BV12, Proposition 5.9]. Fix objects $(m, \sigma_{m,-}) \in \mathsf{Z}({}_L\mathscr{C})$ and $(x, \sigma_{x,-}) \in \mathsf{Z}({}_B\mathscr{C})$. The compatibility of the half-braiding of $m \otimes x$ with the unit of $\mathscr{C}$ is a short computation:

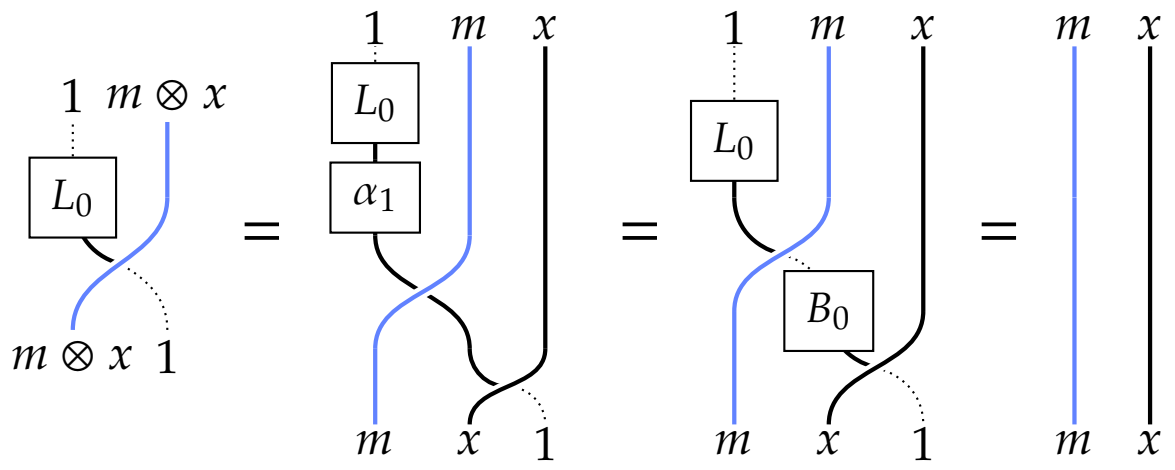

Similarly, we verify the hexagon axiom:

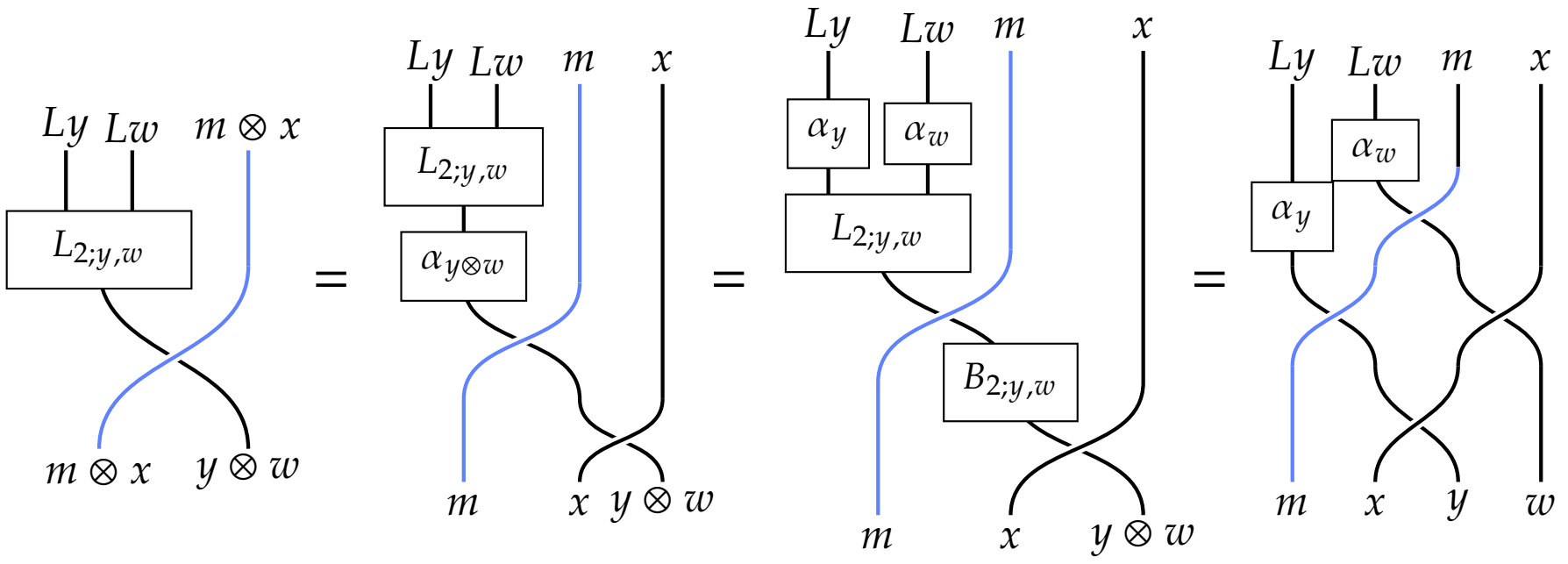

The compatibility of the action $\alpha\colon LB \Longrightarrow L$ with the multiplication and unit of $B$ asserts that $\mathsf{Z}({}_L\mathscr{C})$ is a right $\mathsf{Z}({}_B\mathscr{C})$-module category. By construction, for all $(m, \sigma_{m,-}) \in \mathsf{Z}({}_L\mathscr{C})$ and $(x, \sigma_{x,-}) \in \mathsf{Z}({}_B\mathscr{C})$ we have

$$U^{(L)}\big((m, \sigma_{m,-}) \triangleleft (x, \sigma_{x,-})\big) \;=\; m \otimes x \;=\; U^{(L)}(m, \sigma_{m,-}) \otimes U^{(B)}(x, \sigma_{x,-}).$$

Thus, $U^{(L)}$ is a strict comodule functor over $U^{(B)}$. □

**Notation 6.22.** We extend our colouring scheme to universal coactions:

| | |
|---|---|
| 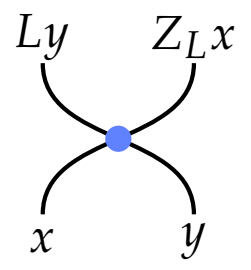 | 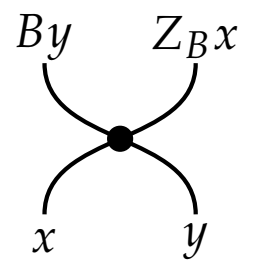 |
| $\xi_{x,y}\colon x \otimes y \longrightarrow Ly \otimes Z_L x$ | $\chi_{x,y}\colon x \otimes y \longrightarrow By \otimes Z_B x$ |

**Remark 6.23.** The identification of $\mathscr{C}^{Z_B}$ and $\mathscr{C}^{Z_L}$ with the generalised Drinfeld centre and its twisted cousin suggest that $Z_L$ is a comodule monad over $Z_B$. In analogy with Equation (6.2.1) we define the coaction $Z_{L;\mathsf{a}}$ of $Z_L$ by

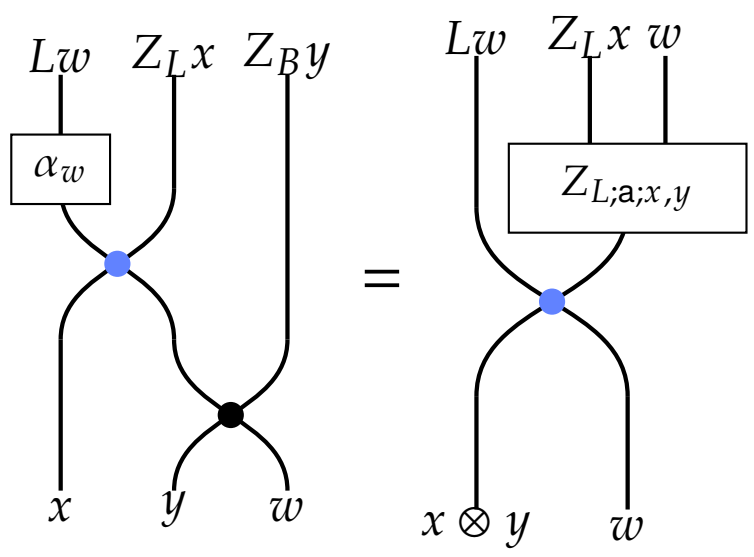

(6.3.2)

**Proposition 6.24.** *Let $\alpha\colon LB \Longrightarrow L$ be an oplax monoidal right action of a bimonad $B\colon \mathscr{C} \longrightarrow \mathscr{C}$ on an oplax monoidal functor $L\colon \mathscr{C} \longrightarrow \mathscr{C}$. Suppose furthermore that the centralisers $(Z_L, \xi)$ of $L$ and $(Z_B, \chi)$ of $B$ exist. The coaction of Equation (6.3.2) turns $Z_L$ into a comodule monad over $Z_B$ such that $\mathscr{C}^{Z_L}$ is isomorphic as a right module category over $\mathscr{C}^{Z_B}$ to $\mathsf{Z}({}_L\mathscr{C})$.*

*Proof.* By Remark 6.16 and Proposition 6.17, we have monadic adjunctions

$$F^{(B)}\colon \mathscr{C} \rightleftarrows \mathsf{Z}({}_B\mathscr{C}) : U^{(B)} \qquad \text{and} \qquad F^{(L)}\colon \mathscr{C} \rightleftarrows \mathsf{Z}({}_L\mathscr{C}) : U^{(L)}$$

that, due to [BV12, Remark 5.13], give rise to the bimonad $Z_B$ and monad $Z_L$. Lemma 6.21 shows that $U^{(L)}$ is a strict comodule functor over $U^{(B)}$.

By Theorem 5.31, see also Remark 5.24 and Example 5.27, $Z_L$ is a comodule monad over $Z_B$; the coaction $\lambda\colon Z_L(-\otimes =) \longrightarrow Z_L(-) \otimes Z_B(=)$ implementing the action of $\mathscr{C}^{Z_B}$ on $\mathscr{C}^{Z_L}$ is for all $x, y \in \mathscr{C}$ given by

$$\lambda_{x,y}\colon Z_L(x\otimes y) \xrightarrow{Z_L(\eta_x^{(Z_L)} \otimes \eta_y^{(Z_B)})} Z_L(Z_L x \otimes Z_B y) \xrightarrow{\nabla_{Z_L x \otimes Z_B y}} Z_L x \otimes Z_B y,$$

where $\nabla$ is defined as in Equation (2.3.1).

By using the relation between universal coactions and half-braidings, explained in Equation (6.2.2), and applying the hexagon identity we compute Figure 6.2. The uniqueness of universal coactions implies that $\lambda = Z_{L;\mathsf{a}}$.

It remains to show that $\mathscr{C}^{Z_L}$ and $\mathsf{Z}({}_L\mathscr{C})$ are isomorphic as modules over $\mathscr{C}^{Z_B}$. By Proposition 5.34, the comparison functor $K^{Z_B}\colon \mathsf{Z}({}_B\mathscr{C}) \longrightarrow \mathscr{C}^{Z_B}$ is strong monoidal and $K^{Z_L}\colon \mathsf{Z}({}_L\mathscr{C}) \longrightarrow \mathscr{C}^{Z_L}$ is a strong comodule functor over it. Furthermore, due to Proposition 6.17, both $K^{Z_B}$ and $K^{Z_L}$ admit inverses

$$E^{Z_B}\colon \mathscr{C}^{Z_B} \longrightarrow \mathsf{Z}({}_B\mathscr{C}) \qquad \text{and} \qquad E^{Z_L}\colon \mathscr{C}^{Z_L} \longrightarrow \mathsf{Z}({}_L\mathscr{C}).$$

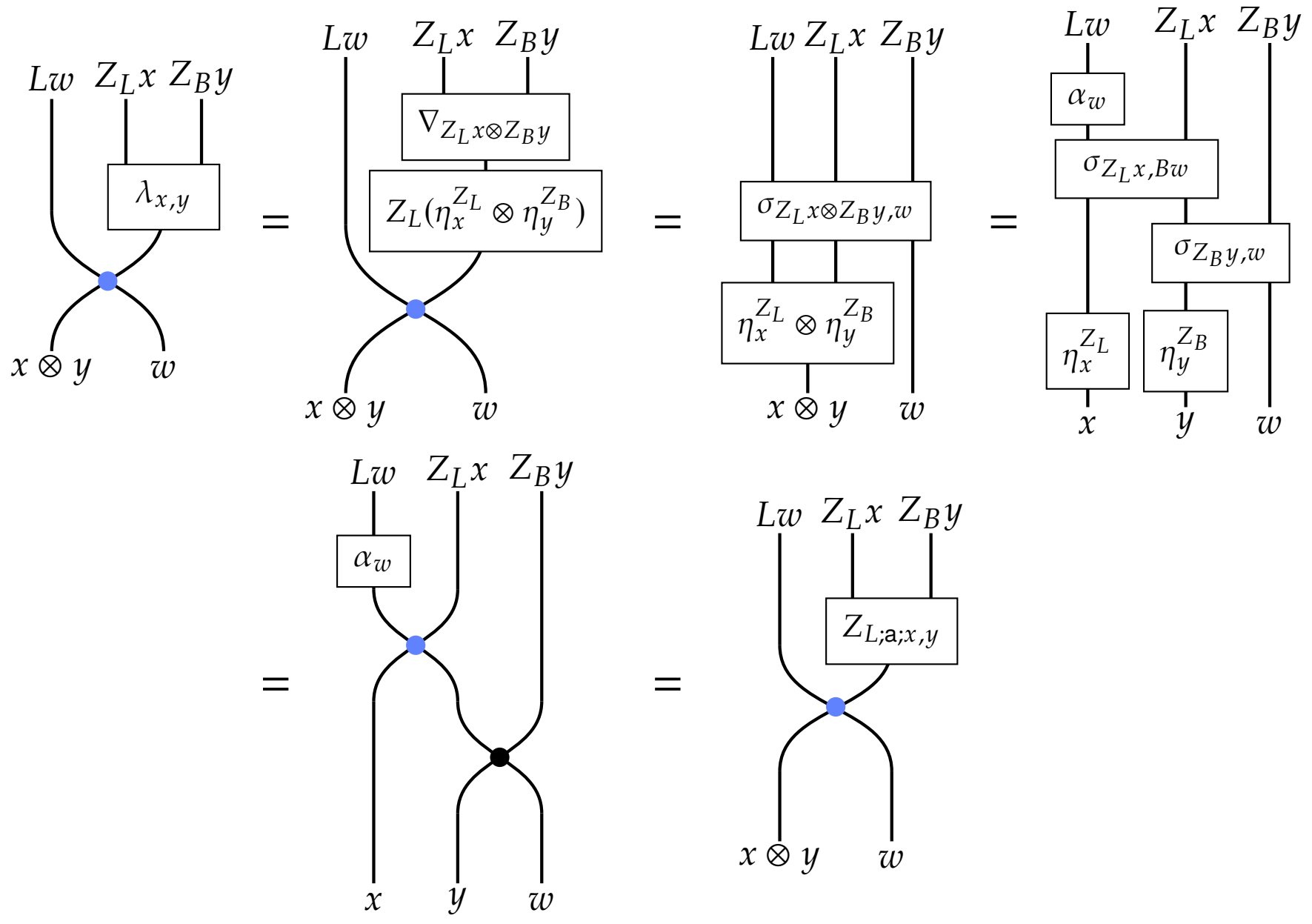

Figure 6.2: The arrows $\lambda$ and $Z_{L;\mathsf{a}}$ satisfy the same universal property.

Using that $E^{Z_B}$ is monoidal, we identify the right action of $\mathsf{Z}({}_B\mathscr{C})$ on $\mathsf{Z}({}_L\mathscr{C})$ with a right action $\blacktriangleleft\colon \mathsf{Z}({}_L\mathscr{C}) \times \mathscr{C}^{Z_B} \longrightarrow \mathsf{Z}({}_L\mathscr{C})$ of $\mathscr{C}^{Z_B}$ by setting

$$\mathsf{Z}({}_L\mathscr{C}) \times \mathscr{C}^{Z_B} \xrightarrow{\mathrm{id}\times E^{Z_B}} \mathsf{Z}({}_L\mathscr{C}) \times \mathsf{Z}({}_B\mathscr{C}) \xrightarrow{\triangleleft} \mathsf{Z}({}_L\mathscr{C}).$$

For any $m \in \mathsf{Z}({}_L\mathscr{C})$ and $x \in \mathsf{Z}({}_L\mathscr{C})$ we have

$$K^{Z_L}(m \blacktriangleleft x) = K^{Z_L}(m \triangleleft E^{Z_B}x) \xrightarrow{K^{Z_L}_{\mathsf{a};m,E^{Z_B}x}} K^{Z_L}m \triangleleft K^{Z_B}E^{Z_B}x = K^{Z_L}m \triangleleft x,$$

and hence $K^{Z_L}\colon \mathsf{Z}({}_L\mathscr{C}) \longrightarrow \mathscr{C}^{Z_L}$ is an isomorphism of module categories. □

**Example 6.25.** Let $\mathscr{C}$ be a rigid monoidal category, and let $(Z_L, \xi)$ and $(Z_B, \chi)$ be the centralisers of ${}^{\vee\vee}(-)\colon \mathscr{C} \longrightarrow \mathscr{C}$ and $\mathrm{Id}_{\mathscr{C}}$, respectively. Then there exists a trivial right action of $\mathrm{Id}_{\mathscr{C}}$ on ${}^{\vee\vee}(-)$:

$$\mathrm{id}_x\colon {}^{\vee\vee}(\mathrm{Id}_{\mathscr{C}}(x)) \longrightarrow {}^{\vee\vee}x, \qquad \text{for all } x \in \mathscr{C}.$$

This action turns $Z_L$ into a comodule monad over $Z_B$, and its modules $\mathscr{C}^{Z_L}$ become isomorphic to $\mathsf{Q}(\mathscr{C})$ as a $\mathscr{C}^{Z_B}$-module category.

Recall that, by Remark 4.11, we can identify $\mathsf{Q}(\mathscr{C})$ with $\mathsf{A}(\mathscr{C}^{\mathrm{op,rev}})^{\mathrm{op}}$.

**Definition 6.26.** Assume ${}^{\vee\vee}(-)$, $\mathrm{Id}_{\mathscr{C}}\colon \mathscr{C} \longrightarrow \mathscr{C}$ to admit centralisers $(Z_L, \xi)$ and $(Z_B, \chi)$. We call $\mathfrak{Q}(\mathscr{C}) := Z_L$ the *anti-central comodule monad of* $\mathscr{C}$.

## 6.4 The Drinfeld and anti-Drinfeld double of a Hopf monad

We are now able to untangle the relationship between the various different categories and adjunctions in Figure 6.1.

**Hypothesis 6.27.** For the rest of this chapter, fix a Hopf monad $H$ on a rigid category $\mathscr{C}$, together with an oplax monoidal endofunctor $L$ on $\mathscr{C}^H$, a bimonad $B$ on $\mathscr{C}^H$, an oplax monoidal right action $\alpha\colon LB \Longrightarrow B$, and assume that the cross products $B \rtimes H$ and $L \rtimes H$ have centralisers $(Z_{B\rtimes H}, \nu)$ and $(Z_{L\rtimes H}, \tau)$.

The following observation—which follows by a straightforward calculation—extends the action of $B$ on $L$ to an action of the respective cross products.

**Lemma 6.28.** *The oplax monoidal right action $\alpha\colon LB \Longrightarrow B$ induces an oplax monoidal action of $B \rtimes H$ on $L \rtimes H$ by*

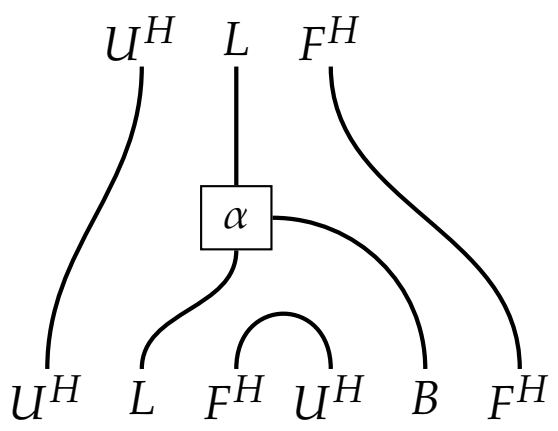

The next result is a variant of [BV12, Theorem 7.4].

**Theorem 6.29.** *The functors $B, L\colon \mathscr{C}^H \longrightarrow \mathscr{C}^H$ admit centralisers $(Z_B, \chi)$, $(Z_L, \xi)$, such that $Z_B$ lifts $Z_{B\rtimes H}$ as a bimonad and $Z_L$ lifts $Z_{L\rtimes H}$ as a comodule monad.*

*Proof.* By [BV12, Theorem 7.4(a)], there are centralisers $(Z_L, \xi)$ and $(Z_B, \chi)$ of $L$ and $B$ that, for all $(x, \nabla_x)$, $(y, \nabla_y) \in \mathscr{C}^H$, satisfy

$$U^H Z_L(x, \nabla_x) = Z_{L\rtimes H}x, \qquad U^H \xi_{(x,\nabla_x),(y,\nabla_y)} = (U^H L\nabla_y \otimes Z_{L\rtimes H}x) \circ \tau_{x,y},$$
$$U^H Z_B(x, \nabla_x) = Z_{B\rtimes H}x, \qquad U^H \chi_{(x,\nabla_x),(y,\nabla_y)} = (U^H B\nabla_y \otimes Z_{B\rtimes H}x) \circ \nu_{x,y}.$$

The second and third part of *ibid* state that $Z_L$ is a lift of the monad $Z_{L\rtimes H}$ and $Z_B$ is a lift of the bimonad $Z_{B\rtimes H}$. It remains for us to show that the coactions of $Z_L$ and $Z_{L\rtimes H}$ are compatible with the forgetful functor $U^H\colon \mathscr{C}^H \longrightarrow \mathscr{C}$. Fixing objects $(x, \nabla_x)$, $(y, \nabla_y) \in \mathscr{C}^H$ and $w \in \mathscr{C}$, this follows from Figure 6.3.

The uniqueness property of universal coactions as given in Lemma 6.10 then implies that $U^H Z_{L;\mathrm{a};(x,\nabla_x),(y,\nabla_y)} = Z_{L\rtimes H;\mathrm{a};x,y}$. Since $U^H\colon \mathscr{C}^H \longrightarrow \mathscr{C}$ is a strict comodule functor, the claim follows. □

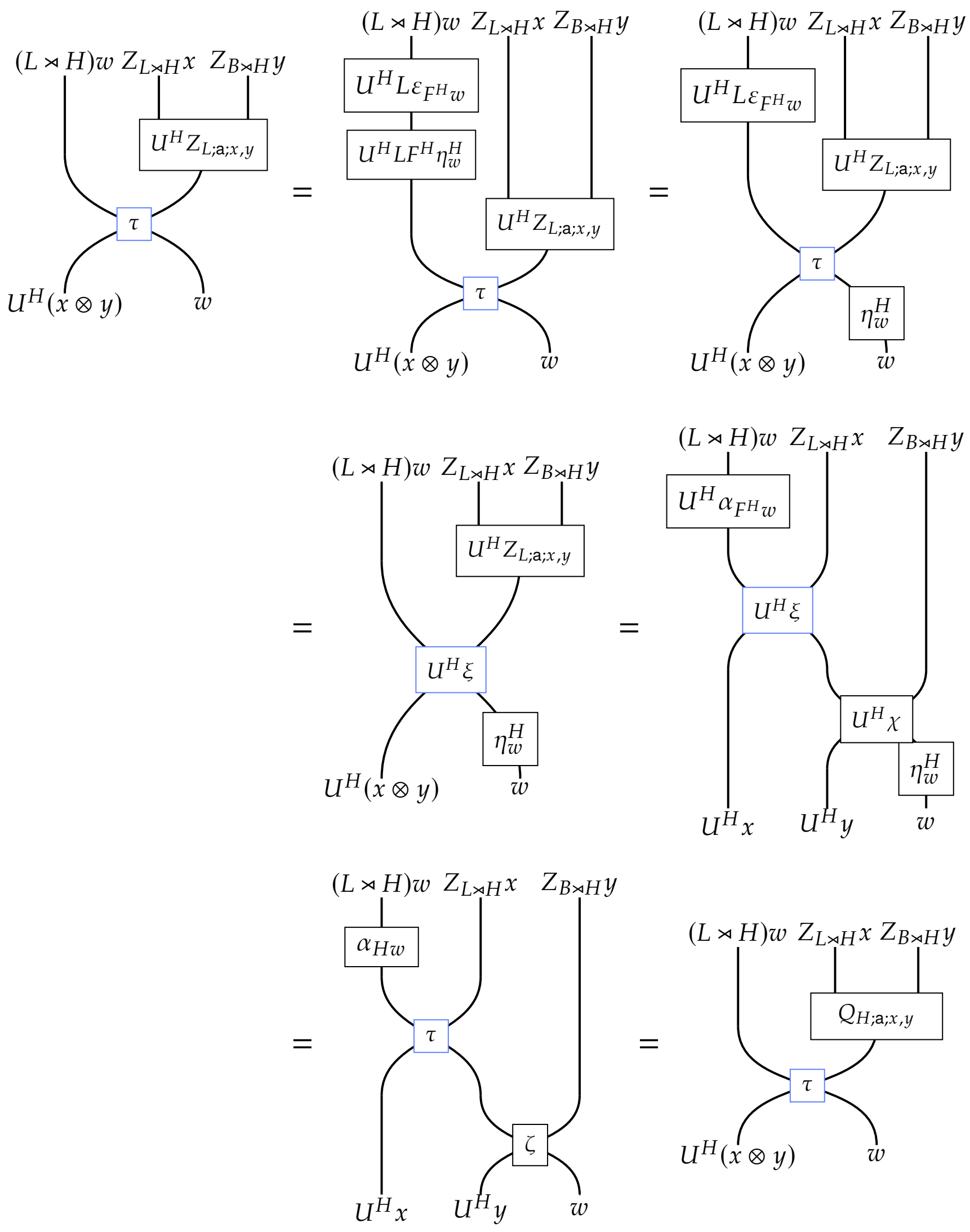

Figure 6.3: The coactions of $Z_L$ and $Z_{L\rtimes H}$ are compatible with the forgetful functor.

**Remark 6.30.** The previous theorem together with Lemma 6.2 imply that we obtain a comodule monad $D(L,H) := Z_L \rtimes H$ over $D(B,H) := Z_B \rtimes H$. The correspondence between lifts and monads given in Proposition 6.5 yields a unique comodule distributive law

$$\left(HZ_{L\rtimes H} \xrightarrow{\Omega} Z_{L\rtimes H}H,\ HZ_{B\rtimes H} \xrightarrow{\Lambda} Z_{B\rtimes H}H\right),$$

such that

$$D(L,H) = Z_{L\rtimes H} \circ_\Omega H \qquad \text{and} \qquad D(B,H) = Z_{B\rtimes H} \circ_\Lambda H.$$

**Definition 6.31.** We call $D(B,H)$ and $D(L,H)$ of Remark 6.30 the *double* and *twisted double* of the pairs $(B,H)$ and $(L,H)$, respectively.

The relationship between doubles and generalised Drinfeld centres is studied in [BV12, Proposition 7.5 and Theorem 7.6]. Our next result uses the same techniques to prove how twisted doubles parameterise twisted centres.

**Theorem 6.32.** *The twisted double $D(L,H)$ is a comodule monad over $D(B,H)$, and $\mathcal{C}^{D(L,H)}$ is isomorphic to $\mathsf{Z}({}_L\mathcal{C}^H)$ as a $\mathcal{C}^{D(B,H)}$-module category.*

*Proof.* Since $Z_L$ is a lift of $Z_{L\rtimes H}$ as a comodule monad, the twisted double $D(L,H)$ is a comodule monad over $D(B,H)$. By Lemma 6.6, this implies the existence of an isomorphism $I\colon \mathcal{C}^{D(L,H)} \xrightarrow{\sim} \left(\mathcal{C}^H\right)^{Z_L}$ of $\mathcal{C}^{D(B,H)}$-module categories. Due to the proof of Proposition 6.24, the comparison functor $K^{Z_L}\colon \mathsf{Z}({}_L\mathcal{C}^H) \longrightarrow \left(\mathcal{C}^H\right)^{Z_L}$ implements an isomorphism of module categories and the statement follows by considering

$$\mathcal{C}^{D(L,H)} \xrightarrow{I} \left(\mathcal{C}^H\right)^{Z_L} \xrightarrow{E^{Z_L}} \mathsf{Z}({}_L\mathcal{C}^H). \qquad \square$$

The following definition can be understood as an extension of the notion of the anti-Drinfeld double given by [HKRS04] to the monadic framework.

**Definition 6.33.** Let $H$ be a Hopf monad on a rigid category $\mathcal{C}$. For $B := \mathrm{Id}_{\mathcal{C}^H}$ and $L := {}^{\vee\vee}(-)\colon \mathcal{C}^H \longrightarrow \mathcal{C}^H$, we call $D(H) := D(B,H)$ and $Q(H) := D(L,H)$ the *Drinfeld* and *anti-Drinfeld double* of $H$, respectively.

## 6.5 Pairs in involution for Hopf monads

We now consider a Hopf monad that admits a double and anti-double, and develop the notion of pairs in involution in this setting. Classically, these consist of a group-like and character of a Hopf algebra that implement the square of its antipode by their adjoint actions.

**Definition 6.34.** Let $H$ be a Hopf monad on a rigid monoidal category $\mathcal{C}$. A *character* of $H$ is an $H$-algebra $\beta := (1, \nabla_\beta) \in \mathcal{C}^H$, whose underlying object is the monoidal unit $1 \in \mathcal{C}$.

**Remark 6.35.** Explicitly, Definition 6.34 says that a group-like $g$ of $H$ satisfies

$$H_{2;x,y} \circ g_{x\otimes y} = g_x \otimes g_y \qquad \text{and} \qquad H_0 \circ g_1 = \mathrm{id}_1, \qquad \text{for all } x,y \in \mathcal{C}.$$

The characters $\mathrm{Char}(H)$ of a Hopf monad $H$ on $\mathcal{C}$ form a monoid and, by [BV07, Lemma 3.21], the set $\mathrm{Gr}(H)$ of group-likes bears a group structure.

**Example 6.36.** The group-likes of a Hopf monad $H$ act on it by conjugation. We recall this construction based on [BV07, Section 1.4]. Given a natural transformation $g\colon \mathrm{Id}_{\mathscr{C}} \longrightarrow H$, define the *left* and *right regular action* of $g$ on $H$ to be the natural transformations defined by

$$L_g\colon H \xrightarrow{gH} H^2 \xrightarrow{\mu^{(H)}} H \qquad \text{and} \qquad R_g\colon H \xrightarrow{Hg} H^2 \xrightarrow{\mu^{(H)}} H.$$

**Definition 6.37.** Every pair $(g \in \mathrm{Gr}(H),\ \beta \in \mathrm{Char}(H))$ of a Hopf monad $H\colon \mathscr{C} \longrightarrow \mathscr{C}$ gives rise to natural transformations

Recall the definition of $H_3$ from Remark 2.35.

$$\begin{aligned}\mathrm{Ad}_g &:= L_g \circ R_{g^{-1}}\colon H \Longrightarrow H,\\ \mathrm{Ad}_\beta &:= (\nabla_\beta \otimes H(-) \otimes \nabla_{{}^\vee\beta}) \circ H_{3;1,-,1}\colon H \Longrightarrow H,\end{aligned}$$

called the *adjoint actions* of $g$ and $\beta$ on $H$, respectively.

To define pairs in involution, we need an analogue of the square of the antipode of a Hopf algebra. This notion is developed in [BV07, Section 7.3].

**Definition 6.38.** Suppose $\phi\colon \mathrm{Id}_{\mathscr{C}} \longrightarrow {}^{\vee\vee}(-)$ to be a pivotal structure on $\mathscr{C}$ and let $H\colon \mathscr{C} \longrightarrow \mathscr{C}$ be a Hopf monad. The *square of the antipode* of $H$ is a natural transformation $S^2\colon H \Longrightarrow H$, given for all $x \in \mathscr{C}$ by

$$S^2_x := \phi^{-1}_{Hx} \circ s^l_{{}^\vee Hx} \circ H{}^\vee s^l_x \circ H\phi_x,$$

where $s^l$ is the left antipode of $H$, see Equation (5.2.1) and [BV07, Section 3.3].

Analogous to the Hopf algebraic case, we state the following:

**Definition 6.39.** Let $H\colon \mathscr{C} \longrightarrow \mathscr{C}$ be a Hopf monad, and $\phi\colon \mathrm{Id}_{\mathscr{C}} \longrightarrow {}^{\vee\vee}(-)$ a pivotal structure. A *pair in involution* $(g,\beta) \in \mathrm{PI}^\phi_H$ of $H$ and $\phi$ consists of a group-like $g \in \mathrm{Gr}(H)$ and a character $\beta \in \mathrm{Char}(H)$, such that for all $x \in \mathscr{C}$

$$Ad_{g,x} = Ad_{\beta,x} \circ S^2_x.$$

To prove that pairs in involution correspond to certain pivotal structures on the Drinfeld centre of $\mathscr{C}^H$, we need two technical results. The first one is classical; for a proof see for example [BV07, Lemmas 1.2 and 1.3].

**Lemma 6.40.** *Let $H$ be a monad with canonical forgetful functor $U^H\colon \mathscr{C}^H \longrightarrow \mathscr{C}$. If $F, G\colon \mathscr{C} \longrightarrow \mathscr{D}$ are functors, for some category $\mathscr{D}$, then there is a canonical bijection*

$$(-)^\sharp\colon\ \mathrm{Nat}(F, GH) \longrightarrow \mathrm{Nat}(FU^H, GU^H), \qquad \alpha \longmapsto \alpha^\sharp,$$

*where* $\alpha^\sharp_{(m,\nabla_m)} := G\nabla_m \circ \alpha_m$.

The next lemma is a variant of [BV07, Lemma 7.5].

**Lemma 6.41.** *Let $\phi\colon \mathrm{Id}_{\mathscr{C}} \longrightarrow {}^{\vee\vee}(-)$ be a pivotal structure on $\mathscr{C}$ and $H\colon \mathscr{C} \longrightarrow \mathscr{C}$ a Hopf monad. For any group-like $g \in \mathrm{Gr}(H)$ and character $\beta \in \mathrm{Char}(H)$ the following are equivalent:*

(*i*) *the morphisms $g$ and $\beta$ form a pair in involution of $H$ and $\phi$; and*
(*ii*) *the arrow $\phi g^\sharp \in \mathrm{Nat}(U^H, {}^{\vee\vee}(-) \circ U^H)$ lifts to $\mathrm{Nat}(\mathrm{Id}_{\mathscr{C}^H}, \beta \otimes {}^{\vee\vee}(-) \otimes {}^{\vee}\beta)$.*

*Proof.* Consider a module $(m, \nabla_m) \in \mathscr{C}^H$. By [BV07, Theorem 3.8(a)] and the definition of $S^2$, the action on ${}^{\vee\vee}m$ is given by

$$\nabla_{{}^{\vee\vee}m} = {}^{\vee\vee}\nabla_m \circ s^l_{{}^{\vee\vee}Hm} \circ H({}^{\vee}s^l_m) = \phi_m \circ \nabla_m \circ S^2_m \circ H(\phi_m^{-1}),$$

and therefore we have

$$\nabla_{\beta\otimes{}^{\vee\vee}m\otimes{}^{\vee}\beta} = (\nabla_\beta \otimes \nabla_{{}^{\vee\vee}m} \otimes \nabla_{{}^{\vee}\beta})H_{3;1,m,1} = (\nabla_\beta \otimes \phi_m\nabla_m S^2_m H(\phi_m^{-1}) \otimes \nabla_{{}^{\vee}\beta})H_{3;1,m,1}.$$

By definition, $\phi g^\sharp$ lifts to a natural transformation from $\mathrm{Id}_{\mathscr{C}^H}$ to $\beta \otimes {}^{\vee\vee}(-) \otimes {}^{\vee}\beta$ if and only if for any $H$-module $(m, \nabla_m)$ we have

$$(\phi g^\sharp)_m \nabla_m = \nabla_{\beta\otimes{}^{\vee\vee}m\otimes{}^{\vee}\beta} H\big((\phi g^\sharp)_m\big). \tag{6.5.1}$$

Let us now successively simplify both sides of this equation. Using the naturality of $g\colon \mathrm{Id}_{\mathscr{C}} \Longrightarrow H$, the fact that $\nabla_m$ is an action, and the definition of $g^\sharp$ as given in Lemma 6.40, we can rewrite the left hand side as

$$(\phi g^\sharp)_m \nabla_m = \phi_m \nabla_m g_m \nabla_m = \phi_m \nabla_m H(\nabla_m) g_{Hm} = \phi_m \nabla_m \mu^{(H)}_m g_{Hm}.$$

Similarly, we simplify the right-hand side to

$$\begin{aligned}
\nabla_{\beta\otimes{}^{\vee\vee}m\otimes{}^{\vee}\beta} H\big((\phi g^\sharp)_m\big) &= (\nabla_\beta \otimes \phi_m\nabla_m S^2_m H(\phi_m^{-1}) \otimes \nabla_{{}^{\vee}\beta})H_{3;1,m,1} H\big((\phi g^\sharp)_m\big) \\
&= (\nabla_\beta \otimes \phi_m\nabla_m S^2_m H(\phi_m^{-1}) H\big((\phi g^\sharp)_m\big) \otimes \nabla_{{}^{\vee}\beta})H_{3;1,m,1} \\
&= (\nabla_\beta \otimes \phi_m\nabla_m S^2_m H(\nabla_m g_m) \otimes \nabla_{{}^{\vee}\beta})H_{3;1,m,1} \\
&= (\nabla_\beta \otimes \phi_m\nabla_m H(\nabla_m g_m) S^2_m \otimes \nabla_{{}^{\vee}\beta})H_{3;1,m,1} \\
&= \phi_m\nabla_m H(\nabla_m g_m)(\nabla_\beta \otimes \mathrm{id}_{Hm} \otimes \nabla_{{}^{\vee}\beta})H_{3;1,m,1} S^2_m \\
&= \phi_m\nabla_m \mu^{(H)}_m H(g_m)\, \mathrm{Ad}_{\beta,m}\, S^2_m.
\end{aligned}$$

Using the fact that $\phi$ is an isomorphism, Equation (6.5.1) can be restated as

$$\begin{aligned} & \nabla_m \mu_m^{(H)} g_{Hm} = \nabla_m \mu_m^{(H)} H(g_m)\,\mathrm{Ad}_{\beta,m}\, S_m^2 \\ \iff \quad & \nabla_m L_{g,m} = \nabla_m R_{g,m}\,\mathrm{Ad}_{\beta,m}\, S_m^2 \end{aligned}$$

By Lemma 6.40, the above is equivalent to $L_{g,m} = R_{g,m}\,\mathrm{Ad}_{\beta,m}\, S_m^2$. We conclude the proof by multiplying both sides with $R_{g^{-1},m}$. □

Lemma 6.41 leads to an identification of pairs in involution of $H$ and $\phi$ with certain quasi-pivotal structures on $\mathscr{C}^H$.

**Proposition 6.42.** *Let $H$ be a Hopf monad on $\mathscr{C}$ and $\phi\colon \mathrm{Id}_{\mathscr{C}} \longrightarrow {}^{\vee\vee}(-)$ a pivotal structure on $\mathscr{C}$. Then $(H,\phi)$ admits a pair in involution if and only if there exists a quasi-pivotal structure on $\mathscr{C}^H$ whose underlying invertible object is a character.*

*Proof.* We proceed analogous to [BV07, Proposition 7.6]. Suppose $(g,\beta) \in \mathrm{PI}_H^{\phi}$. By Lemma 6.41, the natural transformation $\phi g^\sharp$ lifts to a natural isomorphism

$$\rho_{\beta,x}\colon x \longrightarrow \beta \otimes {}^{\vee\vee}x \otimes {}^{\vee}\beta, \qquad \text{for all } x \in \mathscr{C}^H.$$

Since $\phi$ is monoidal by definition, and $g^\sharp$ is monoidal by virtue of $g$ being a group-like—see [BV07, Lemma 3.20]—we obtain a quasi-pivotal structure

$$\rho_\beta\colon \mathrm{Id}_{\mathscr{C}^H} \longrightarrow \beta \otimes {}^{\vee\vee}(-) \otimes {}^{\vee}\beta.$$

On the other hand, let $(\beta,\rho_\beta)$ be a quasi-pivotal structure, for $\beta \in \mathrm{Char}(H)$. Since the forgetful functor $U^H$ is strong monoidal and thus

$$U^H(\beta \otimes {}^{\vee\vee}(-) \otimes {}^{\vee}\beta) = U^H({}^{\vee\vee}(-)) = {}^{\vee\vee}(U^H(-)),$$

there exists a monoidal natural transformation

$$\phi^{-1}_{U^H x} \circ U^H(\rho_{\beta,x})\colon U^H x \longrightarrow U^H x, \qquad \text{for all } x \in \mathscr{C}^H.$$

Apply [BV07, Lemma 3.20] to obtain a unique group-like $g \in \mathrm{Gr}(H)$, with

$$g^\sharp = \phi^{-1}_{U^H(x)} \circ U^H(\rho_{\beta,x}).$$

As $\phi g^\sharp = U^H(\rho_\beta)$ lifts to the quasi-pivotal structure $(\beta,\rho_\beta)$ on $\mathscr{C}^H$, Lemma 6.41 implies that $(g,\beta) \in \mathrm{PI}_H^{\phi}$. □

Let us now study a variant of [BV07, Lemma 2.9].

**Proposition 6.43.** *Let $C, D\colon \mathcal{M} \longrightarrow \mathcal{M}$ be two comodule monads over a bimonad $B\colon \mathcal{C} \longrightarrow \mathcal{C}$. There is a bijective correspondence between morphisms of comodule monads $f\colon D \Longrightarrow C$ and strict module functors $F\colon \mathcal{M}^C \longrightarrow \mathcal{M}^D$ with $U^D F = U^C$.*

*Proof.* As shown in [BV07, Lemma 1.7], any functor $F\colon \mathcal{M}^C \longrightarrow \mathcal{M}^D$ with $U^D F = U^C$ is induced by a unique morphism of monads $f\colon D \longrightarrow C$. That is, $F$ is the identity on morphisms and on objects it is defined by

$$F(m, \nabla_m) = (m, \nabla_m f_m), \qquad \text{for all } (m, \nabla_m) \in m^C.$$

It remains to show that $f$ is a morphism of comodules if and only if $F$ is a strict $\mathcal{C}^B$-module functor. Let $(m, \nabla_m) \in m^C$ and $(x, \nabla_x) \in \mathcal{C}^B$. We compute

$$F\big((m, \nabla_m) \triangleleft (x, \nabla_x)\big) = (m \triangleleft x,\ (\nabla_m \triangleleft \nabla_x) \circ \mathsf{C}_{\mathsf{a};m,x} \circ f_{m \triangleleft x}),$$
$$F(m, \nabla_m) \triangleleft (x, \nabla_x) = (m \triangleleft x,\ (\nabla_m \triangleleft \nabla_x) \circ (f_m \triangleleft Bx) \circ \mathsf{D}_{\mathsf{a};m,x}).$$

According to [BV07, Lemma 1.4], these modules coincide if and only if

$$\mathsf{C}_{\mathsf{a};m,x} \circ f_{m \triangleleft x} = (f_m \triangleleft Bx) \circ \mathsf{D}_{\mathsf{a};m,x},$$

which is exactly the condition for $f$ to be a comodule morphism. □

The above result implies the desired monadic version of Theorem 4.1.

**Theorem 6.44.** *Let $\mathcal{C}$ be a rigid monoidal category, and suppose that $H\colon \mathcal{C} \longrightarrow \mathcal{C}$ is a Hopf monad that admits a double $D(H)$ and anti-double $Q(H)$. The following statements are equivalent:*

- *(i) the monoidal unit $1 \in \mathcal{C}$ lifts to a module over $Q(H)$;*
- *(ii) there is an isomorphism of comodule monads $D(H) \cong Q(H)$; and*
- *(iii) there is an isomorphism of monads $Q(H) \cong D(H)$.*

*Additionally, if $\mathcal{C}$ is pivotal with pivotal structure $\phi$, any of the above statements hold if and only if $H$ and $\phi$ admit a pair in involution.*

*Proof.* (i) $\Longrightarrow$ (ii): Suppose $\omega \in \mathsf{Q}(\mathcal{C}^H)$ with $U^{Q(H)}\omega = 1$. As shown in Equation (4.2.1), this induces a functor of module categories

$$\omega \otimes -\colon \mathcal{C}^{D(H)} \longrightarrow \mathcal{C}^{Q(H)}.$$

Since $U^{Q(H)}\omega = 1 \in \mathscr{C}$, we can apply Proposition 6.43, and obtain that $Q(H)$ and $D(H)$ are isomorphic as comodule monads.

It immediately follows that (ii) implies (iii); we proceed with (iii) $\Longrightarrow$ (i): consider an isomorphism of monads $g \colon Q(H) \Longrightarrow D(H)$. It gives rise to a functor $G \colon \mathscr{C}^{D(H)} \longrightarrow \mathscr{C}^{Q(H)}$ that, on objects, is defined by

$$G(m, \nabla_m) = (m, \nabla_m g_m), \qquad\qquad \text{for all } (m, \nabla_m) \in \mathscr{C}^{D(H)}.$$

We compose $G$ with $E^{Q(H)} \colon \mathscr{C}^{Q(H)} \longrightarrow \mathsf{Q}(\mathscr{C}^H)$—the inverse of the comparison functor defined in Equation (6.2.3)—and see that there exists an object

$$1^{(Q)} := E^{Q(H)}G(1) \in \mathsf{Q}(\mathscr{C}^H),$$

whose underlying object is the unit of $\mathscr{C}$. Now let $(\mathscr{C}, \phi)$ be pivotal. By Lemma 4.22, lifts of $1 \in \mathscr{C}$ to the dual of the anti-centre $\mathsf{Q}(\mathscr{C}^H)$ are in correspondence with quasi-pivotal structures $(\beta, \rho_\beta)$, where $\beta \in \mathsf{Char}(H)$. By Proposition 6.42, such a quasi-pivotal structure exists if and only if there exists a pair in involution for $H$ and $\phi$. □

As a corollary, we can determine whether a category is pivotal in terms of monad isomorphisms between the central and anti-central monad. For a category $\mathscr{C}$, recall Definition 6.15 of its central Hopf monad, and Definition 6.26 of its anti-central comodule monad.

**Corollary 6.45.** *Let $\mathscr{C}$ be a rigid monoidal category. If $\mathscr{C}$ admits a central Hopf monad $\mathfrak{D}(\mathscr{C})$ and an anti-central comodule monad $\mathfrak{Q}(\mathscr{C})$, then it is pivotal if and only if $\mathfrak{D}(\mathscr{C}) \cong \mathfrak{Q}(\mathscr{C})$ as monads.*

*Proof.* We consider the identity $\mathrm{Id}_{\mathscr{C}} \colon \mathscr{C} \longrightarrow \mathscr{C}$ as a Hopf monad. Its Drinfeld and anti-Drinfeld double are $D(\mathrm{Id}_{\mathscr{C}}) = \mathfrak{D}(\mathscr{C}) \rtimes \mathrm{Id}_{\mathscr{C}}$ and $Q(\mathrm{Id}_{\mathscr{C}}) = \mathfrak{Q}(\mathscr{C}) \rtimes \mathrm{Id}_{\mathscr{C}}$. From here it follows that $D(\mathrm{Id}_{\mathscr{C}}) = \mathfrak{D}(\mathscr{C})$, and similarly $Q(\mathrm{Id}_{\mathscr{C}}) = \mathfrak{Q}(\mathscr{C})$. The proof is concluded by Theorem 6.44. □

> The struggle itself toward the heights is enough to fill a man's heart. One must imagine Sisyphus happy.
>
> ALBERT CAMUS; The Myth of Sisyphus

# 7 DUOIDAL R-MATRICES

DUOIDAL CATEGORIES WERE INTRODUCED in [AM10] under the name *2-monoidal categories*, in order to study bilax monoidal functors and various constructions on linear species. They generalise both braided monoidal categories, by considering two monoidal structures that are connected by a non-invertible interchange law, as well as the 2-fold monoidal categories of [BFSV03] where the two tensor products are assumed to share a unit. Duoidal categories have since been used to study higher-dimensional Hopf theory [BCZ13; BS13; AHLF18; LFV20; Böh21], and have also found applications in various other fields of mathematics, [GLF16; SS22; Rom24; Tor24].

The aim of this chapter is to generalise a reconstruction-type result for R-matrices on bimonads, [BV07, Proposition 8.5], which in turn generalises the classical theory of R-matrices for bialgebras. The former has the additional advantage of not requiring a braided monoidal base category, as bimonads—in contrast with bialgebras—may be defined on any monoidal category.

**Theorem 7.21.** *Let $\mathscr{D}$ be a category with monoidal structures $\circ$ and $\bullet$, and $T$ a monad on $\mathscr{D}$ that has a $\circ$-oplax monoidal and a $\bullet$-oplax monoidal structure. Then quasitriangular structures on $T$ are in bijection with duoidal structures on $\mathscr{D}^T$.*

This result can be seen as a generalisation of the double comonoidal monads of [AM10, Section 7], analogous to how R-matrices for bialgebras generalise cocommutative bialgebras.

In Section 7.3 we study the relationship between normal duoidal and linearly distributive categories from this point of view. We see non-planar linearly distributive categories $\mathscr{L}$ as an analogue of preduoidal categories, in the sense that the additional structure trivialises in the monoidal case, see Example 7.27. Equipping $\mathscr{L}$ with a planar structure, we can relate double comonoidal monads to linearly distributive monads in the sense of [Pas12].

## 7.1 duoidal categories

We shall loosely follow the nomenclature of [BM12, Definition 3], who introduced the term *duoidal category*, and the notation of [BCZ13].

**Definition 7.1** ([AM10, Definition 6.1]). A *duoidal* category is a quintuple $(\mathcal{D}, \circ, \perp, \bullet, 1)$, consisting of the following data:

- monoidal categories $(\mathcal{D}, \circ, \perp)$ and $(\mathcal{D}, \bullet, 1)$;
- a not-necessarily invertible natural transformation

$$\zeta\colon (x \bullet y) \circ (a \bullet b) \Longrightarrow (x \circ a) \bullet (y \circ b),$$

  called the *middle interchange law*;
- three morphisms

$$\nu\colon \perp \longrightarrow \perp \bullet \perp, \qquad \varpi\colon 1 \circ 1 \longrightarrow 1, \qquad \iota\colon \perp \longrightarrow 1;$$

This data has to satisfy the following relations:

- $(1, \varpi, \iota)$ is a monoid in $(\mathcal{D}, \circ, \perp)$;
- $(\perp, \nu, \iota)$ is a comonoid in $(\mathcal{D}, \bullet, 1)$;
- the following diagrams commute, witnessing *associativity*:

(7.1.1)
$$\begin{array}{ccc} ((x \bullet y) \circ (a \bullet b)) \circ (c \bullet d) & \xrightarrow{\alpha} & (x \bullet y) \circ ((a \bullet b) \circ (c \bullet d)) \\ \downarrow{\scriptstyle \zeta \circ \mathrm{id}} & & \downarrow{\scriptstyle \mathrm{id} \circ \zeta} \\ ((x \circ a) \bullet (y \circ b)) \circ (c \bullet d) & & (x \bullet y) \circ ((a \circ c) \bullet (b \circ d)) \\ \downarrow{\scriptstyle \zeta} & & \downarrow{\scriptstyle \zeta} \\ ((x \circ a) \circ c) \bullet ((y \circ b) \circ d) & \xrightarrow[\alpha \bullet \alpha]{} & (x \circ (a \circ c)) \bullet (y \circ (b \circ d)) \end{array}$$

(7.1.2)
$$\begin{array}{ccc} ((x \bullet a) \bullet c) \circ ((y \bullet b) \bullet d) & \xrightarrow{\alpha \circ \alpha} & (x \bullet (a \bullet c)) \circ (y \bullet (b \bullet d)) \\ \downarrow{\scriptstyle \zeta} & & \downarrow{\scriptstyle \zeta} \\ ((x \bullet a) \circ (y \bullet b)) \bullet (c \circ d) & & (x \circ y) \bullet ((a \bullet c) \circ (b \bullet d)) \\ \downarrow{\scriptstyle \zeta \bullet \mathrm{id}} & & \downarrow{\scriptstyle \mathrm{id} \bullet \zeta} \\ ((x \circ y) \bullet (a \circ b)) \bullet (c \circ d) & \xrightarrow[\alpha]{} & (x \circ y) \bullet ((a \circ b) \bullet (c \circ d)) \end{array}$$

- the following diagrams commute, witnessing *unitality*:

$$
\begin{array}{ccc}
\bot \circ (a \bullet b) & \xrightarrow{\nu \circ \mathrm{id}} & (\bot \circ \bot) \circ (a \bullet b) \\
{\scriptstyle\lambda}\downarrow & & \downarrow{\scriptstyle\zeta} \\
a \bullet b & \xrightarrow[\lambda^{-1} \bullet \lambda^{-1}]{} & (\bot \circ a) \bullet (\bot \circ b)
\end{array}
\qquad
\begin{array}{ccc}
(a \bullet b) \circ \bot & \xrightarrow{\mathrm{id} \circ \nu} & (a \bullet b) \circ (\bot \bullet \bot) \\
{\scriptstyle\lambda}\downarrow & & \downarrow{\scriptstyle\zeta} \\
a \bullet b & \xrightarrow[\lambda^{-1} \bullet \lambda^{-1}]{} & (a \circ \bot) \bullet (b \circ \bot)
\end{array}
$$

$$
\begin{array}{ccc}
(1 \bullet a) \circ (1 \bullet b) & \xrightarrow{\zeta} & (1 \circ 1) \bullet (a \circ b) \\
{\scriptstyle\lambda \circ \lambda}\downarrow & & \downarrow{\scriptstyle\varpi \bullet \mathrm{id}} \\
a \circ b & \xrightarrow[\lambda]{} & 1 \bullet (a \circ b)
\end{array}
\qquad
\begin{array}{ccc}
(a \bullet 1) \circ (b \bullet 1) & \xrightarrow{\zeta} & (a \circ b) \bullet (1 \circ 1) \\
{\scriptstyle\lambda \circ \lambda}\downarrow & & \downarrow{\scriptstyle\mathrm{id} \bullet \varpi} \\
a \circ b & \xrightarrow[\lambda]{} & (a \circ b) \bullet 1
\end{array}
\tag{7.1.3}
$$

By abuse of notation, we shall often call $\mathfrak{D}$ a duoidal category, leaving the rest of the data implicit.

**Definition 7.2.** A duoidal category $\mathfrak{D}$ is called *normal* if $\bot \cong 1$.

Note that explicitly requiring the existence of $\iota\colon \bot \longrightarrow 1$ in Definition 7.1 is not strictly necessary, as it may be derived from the other specified data:

$$\iota\colon \bot \xrightarrow{\lambda^{-\bot}} \bot \circ \bot \xrightarrow{\lambda^{-1} \circ \rho^{-1}} (1 \bullet \bot) \circ (\bot \bullet 1) \xrightarrow{\zeta} (1 \circ \bot) \bullet (\bot \circ 1) \xrightarrow{\lambda^{\bot} \bullet \rho^{\bot}} 1 \bullet 1 \xrightarrow{\lambda^{1}} 1.$$

**Example 7.3.** By [AM10, Proposition 6.10], a braided monoidal category $(\mathscr{C}, \otimes, 1, \sigma)$ yields a duoidal category $(\mathscr{C}, \otimes, 1, \otimes, 1)$ with structure morphisms

$$\zeta := (a \otimes b) \otimes (c \otimes d) \cong a \otimes (b \otimes c) \otimes d \xrightarrow{a \otimes \sigma_{b,c} \otimes d} a \otimes (c \otimes b) \otimes d \cong (a \otimes c) \otimes (b \otimes d),$$

$$\varpi := 1 \otimes 1 \xrightarrow{\lambda} 1, \qquad \nu := 1 \xrightarrow{\lambda^{-1}} 1 \otimes 1, \qquad \iota := 1 \xrightarrow{\mathrm{id}_1} 1.$$

Note in particular that we have $\rho_1 = \lambda_1$.

**Example 7.4.** The converse of Example 7.3 also holds. If $\mathfrak{D}$ is a duoidal category, such that the interchange law and structure morphisms are isomorphisms, then [AM10, Proposition 6.11] yields a braiding on $(\mathfrak{D}, \circ, \bot)$ and $(\mathfrak{D}, \bullet, 1)$, such that they become isomorphic as braided monoidal categories, and the interchange law arises from the braiding.

Note, however, that there exist non-trivial duoidal structures on a monoidal category $(\mathscr{C}, \otimes, 1)$. Recall the definition of the category ${}^{H}_{H}\mathcal{YD}$ of (left-left) Yetter–Drinfeld modules from Example 2.54. If the Hopf algebra $H$ does not admit an invertible antipode, then ${}^{H}_{H}\mathcal{YD}$ is *lax braided*—the Yetter–Drinfeld braiding is a non-invertible natural transformation that satisfies the braid equations. This yields a duoidal structure on $(\mathscr{C}, \otimes, 1, \otimes, 1)$ that is not braided.

There are various equivalent definitions of duoidal categories. For example, as pseudomonoids in the monoidal 2-category of monoidal categories, oplax monoidal functors, and oplax monoidal natural transformations, see [AM10, Proposition 6.72] and [GLF16, Definition 1]. In particular, this means that $\bullet$ is a lax monoidal and that $\circ$ is an oplax monoidal functor;[16] from this characterisation, one may obtain a *coherence* result, see [Lew72], [AM10, Section 6.2], and [MP22, Theorem 5.9].

[16] This extends to normal duoidal categories, in which $\circ\colon \mathfrak{D}\times\mathfrak{D} \longrightarrow \mathfrak{D}$ and $\perp\colon \mathbb{1} \longrightarrow \mathfrak{D}$ are normal oplax monoidal functors, where $\mathbb{1}$ is the terminal category.

**Proposition 7.5.** *Any* ETC *diagram in a duoidal category commutes.*

Loosely speaking, an ETC diagram is a *formal* diagram $F\colon \mathfrak{J} \longrightarrow \mathfrak{D}$ in the sense of [MP22, p. 20], consisting of only structure morphisms of the duoidal category, such that for all $j \in \mathfrak{J}$ the object $Fj$ is not isomorphic to any of the two units. We refer to [MP22, Definition 5.8 and Theorem 5.9] for a precise definition and a proof of Proposition 7.5. A counterexample in the case of a formal diagram with parallel arrows $1 \bullet 1 \rightrightarrows 1$ is given in [Rom23, Proposition 3.1.6 and Example 3.1.7].

**Remark 7.6.** The tensor product and unit being normal monoidal functors, normal duoidal categories admit an analogue of the well-known coherence result for braided monoidal categories [JS93]. That is, any formal diagram comprised only of the structure morphisms in a normal duoidal category commutes, see [MP22, Theorem 5.18].

### 7.1.1 *Double opmonoidal monads*

**Definition 7.7** ([AM10, Definition 6.25])**.** Suppose that $\mathfrak{D}$ is a duoidal category. A *bimonoid* in $\mathfrak{D}$ is a quintuple $(B, \mu, \eta, \Delta, \varepsilon)$, consisting of a monoid $(B, \mu, \eta)$ in $(\mathfrak{D}, \circ, \perp)$, and a comonoid $(B, \Delta, \varepsilon)$ in $(\mathfrak{D}, \bullet, 1)$, such that the following compatibility conditions are satisfied:

$$
\begin{array}{ccccc}
B \circ B & \xrightarrow{\ \mu\ } & B & \xrightarrow{\ \Delta\ } & B \bullet B \\
{\scriptstyle \Delta\circ\Delta}\downarrow & & & & \uparrow{\scriptstyle \mu\bullet\mu} \\
(B \bullet B) \circ (B \bullet B) & & \xrightarrow[\ \zeta\ ]{} & & (B \circ B) \bullet (B \circ B)
\end{array}
$$

$$
\begin{array}{ccc}
B \circ B & \xrightarrow{\ \varepsilon\circ\varepsilon\ } & 1 \circ 1 \\
{\scriptstyle \mu}\downarrow & & \downarrow{\scriptstyle \varpi} \\
B & \xrightarrow[\ \varepsilon\ ]{} & 1
\end{array}
\qquad
\begin{array}{ccc}
\perp & \xrightarrow{\ \eta\ } & B \\
{\scriptstyle \nu}\downarrow & & \downarrow{\scriptstyle \Delta} \\
\perp \bullet \perp & \xrightarrow[\ \eta\bullet\eta\ ]{} & B \bullet B
\end{array}
\qquad
\begin{array}{ccc}
\perp & & \\
{\scriptstyle \eta}\downarrow & \searrow{\scriptstyle \iota} & \\
B & \xrightarrow[\ \varepsilon\ ]{} & 1
\end{array}
$$

**Proposition 7.8** ([BS13]). *For a monoid $b$ in a duoidal category $(\mathfrak{D}, \circ, \perp, \bullet, 1)$ there is a bijective correspondence between bimonoid structures on $b$, and bimonad structures on the monad $b \circ -$ on $(\mathfrak{D}, \bullet, 1)$.*

**Example 7.9.** A bimonoid in a braided monoidal category $\mathscr{C}$ is the same as a bimonoid in the duoidal category $\mathscr{C}$ from Example 7.3. In this way one recovers the fact that an object $b \in \mathscr{C}$ is a bimonoid if and only if the induced monad $b \otimes -$ is a bimonad on $\mathscr{C}$.

**Example 7.10.** Let $\Bbbk$ be a commutative ring, and $A$ a commutative $\Bbbk$-algebra. In [AM10, Example 6.18] it is shown that the category of $A$-bimodules is duoidal, with the two tensor products given by

$$M \bullet N := M \otimes_A N := {}^{M \otimes_k N}\!\big/_{\langle ma \otimes n - m \otimes an \rangle},$$

and

$$M \circ N := M \otimes_{A \otimes_k A} N := {}^{M \otimes_k N}\!\big/_{\langle amb \otimes n - m \otimes anb \rangle},$$

for all $a \in A$, $m \in M$, and $n \in N$. Furthermore, from [AM10, Example 6.44] we know that a bimonoid in this duoidal category is an $A$-bialgebroid in the sense of Ravenel, see [Rav86, Definition A1.1.1]. In this setting, Proposition 7.8 recovers a special case of [Szl03, Theorems 5.1 and 5.4].

**Definition 7.11** ([AHLF18, Section 7]). A *double opmonoidal monad* on a duoidal category $\mathfrak{D}$ consists of a monad $(T, \mu, \eta)$ on $\mathfrak{D}$, together with a bimonad structures $(T, B_2^\bullet, B_0^\bullet)$ on $(\mathfrak{D}, \bullet, 1)$ and $(T, B_2^\circ, B_0^\circ)$ on $(\mathfrak{D}, \circ, \perp)$, such that the following diagrams commute:

$$\begin{array}{ccc}
T(1\circ 1) & \xrightarrow{T\varpi} & T1 \\
{\scriptstyle T^\circ_{2,1,1}}\downarrow & & \\
T1 \circ T1 & & \downarrow{\scriptstyle T^\bullet_0} \\
{\scriptstyle T^\bullet_0 \circ T^\bullet_0}\downarrow & & \\
1 \circ 1 & \xrightarrow[\varpi]{} & 1
\end{array}
\qquad
\begin{array}{ccc}
T\perp & \xrightarrow{T\nu} & T(\perp \bullet \perp) \\
 & & \downarrow{\scriptstyle T^\bullet_{2,\perp,\perp}} \\
{\scriptstyle T^\circ_0}\downarrow & & T\perp \bullet T\perp \\
 & & \downarrow{\scriptstyle T^\circ_0 \bullet T^\circ_0} \\
\perp & \xrightarrow{\nu} & \perp \bullet \perp
\end{array}
\qquad
\begin{array}{ccc}
T\perp & \xrightarrow{T\iota} & T1 \\
{\scriptstyle T^\circ_0}\downarrow & & \downarrow{\scriptstyle T^\bullet_0} \\
\perp & \xrightarrow[\iota]{} & 1
\end{array}
\tag{7.1.4}$$

$$\begin{array}{ccc}
T((a\bullet b)\circ(c\bullet d)) & \xrightarrow{T\zeta_{a,b,c,d}} & T((a\circ c)\bullet(b\circ d)) \\
{\scriptstyle T^\circ_{2,a\bullet b,c\bullet d}}\downarrow & & \downarrow{\scriptstyle T^\bullet_{2,a\circ c,b\circ d}} \\
T(a\bullet b)\circ T(c\bullet d) & & T(a\circ c)\bullet T(b\circ d) \\
{\scriptstyle T^\bullet_{2,a,b}\circ T^\bullet_{2,c,d}}\downarrow & & \downarrow{\scriptstyle T^\circ_{2,a,c}\bullet T^\circ_{2,b,d}} \\
(Ta\bullet Tb)\circ(Tc\bullet Td) & \xrightarrow[\zeta_{Ta,Tb,Tc,Td}]{} & (Ta\circ Tc)\bullet(Tb\circ Td)
\end{array}
\tag{7.1.5}$$

**Example 7.12.** Let $(\mathcal{C}, \otimes, 1)$ be a braided monoidal category with braiding $\sigma$, seen as a duoidal category as in Example 7.3. A bimonad $B$ on $(\mathcal{C}, \otimes, 1)$ that additionally satisfies the equation $B_2 \circ B\sigma = \sigma \circ B_2$ is a double opmonoidal monad on $\mathcal{C}$, where the two oplax monoidal structures are the same, and the commutativity of Diagram (7.1.4) amounts to the fact that the monoidal structure morphisms of $\mathcal{C}$ lift to the category of $B$-algebras; see Proposition 5.9.

**Example 7.13.** For a bialgebra $B$ in $(\mathsf{Vect}, \otimes, \Bbbk)$, the endofunctor $B \otimes -$ is a double opmonoidal monad in $(\mathsf{Vect}, \otimes, \Bbbk, \otimes, \Bbbk)$. For $U, V, W, X \in \mathsf{Vect}$, as well as $b \in B$, $u \in U$, $v \in V$, $w \in W$, and $x \in X$, Diagram (7.1.5) simplifies to

$$b_{(1)} \otimes u \otimes b_{(3)} \otimes w \otimes b_{(2)} \otimes v \otimes b_{(4)} \otimes x \;=\; b_{(1)} \otimes u \otimes b_{(2)} \otimes w \otimes b_{(3)} \otimes v \otimes b_{(4)} \otimes x,$$

which is equivalent to $b_{(1)} \otimes b_{(2)} = b_{(2)} \otimes b_{(1)}$; i.e., $B$ has to be cocommutative.

As in the case of R-matrices for bialgebras and bimonads, requiring that the interchange morphism of a duoidal category $\mathfrak{D}$ lifts to the category of modules is a rather strong condition.

**Proposition 7.14** ([AHLF18, Theorem 7.2]). *Let $\mathfrak{D}$ be a duoidal category and $T\colon \mathcal{C} \longrightarrow \mathcal{C}$ a monad. Then the structure morphisms and interchange law of $\mathfrak{D}$ lift to $\mathfrak{D}^T$ if and only if $T$ is a double opmonoidal monad.*

*In particular, if $T$ is a double opmonoidal monad, then $\mathfrak{D}^T$ is a duoidal category.*

## 7.2 R-matrices

Instead of the situation of Proposition 7.14, we are instead interested in studying which additional structure one can impose on $T$ such that $\mathfrak{D}^T$ becomes duoidal, where the interchange morphism is instead giving by "twisting" that of $\mathfrak{D}$. This generalises so-called *R-matrices* for bialgebras and bimonads [Kas98, Part II] and [BV07, Section 8.2].

**Proposition 7.15** ([BV07, Theorem 8.5]). *Let $B$ be a bimonad on the monoidal category $\mathcal{C}$. Braidings on $\mathcal{C}^T$ are in bijective correspondence with R-matrices on $B$.*

A crucial feature of R-matrices for bimonads is that they can be defined on not necessarily braided monoidal categories. Our definition of duoidal R-matrices incorporates this feature.

**Definition 7.16.** A category $\mathfrak{D}$ is called *preduoidal* if it is equipped with two monoidal structures $(\circ, \perp)$ and $(\bullet, 1)$.

**Definition 7.17.** Let $\mathcal{D}$ be a preduoidal category. A monad $T$ on $\mathcal{D}$ that is equipped with two bimonad structures over $(\mathcal{D}, \circ, \perp)$ and $(\mathcal{D}, \bullet, 1)$ is called a *separately opmonoidal monad* on $\mathcal{D}$.

**Definition 7.18.** Let $\mathcal{D}$ be a preduoidal category and $T$ a separately opmonoidal monad on $\mathcal{D}$. An *R-matrix* on $T$ consists of a natural transformation

$$R := \{ R_{a,b,c,d} \colon (a \bullet b) \circ (c \bullet d) \Longrightarrow (Ta \bullet Tc) \circ (Tb \bullet Td) \}_{a,b,c,d \in \mathcal{D}},$$

as well as morphisms of $T$-algebras

$$\nu \colon (\perp, T_0^\circ) \longrightarrow (\perp, T_0^\circ) \bullet (\perp, T_0^\circ), \qquad \varpi \colon (1, T_0^\bullet) \circ (1, T_0^\bullet) \longrightarrow (1, T_0^\bullet),$$
$$\iota \colon (\perp, T_0^\circ) \longrightarrow (1, T_0^\bullet),$$

such that $(1, \varpi, \iota)$ is a monoid in $(\mathcal{D}^T, \circ, \perp)$; the tuple $(\perp, \nu, \iota)$ is a comonoid in $(\mathcal{D}^T, \bullet, 1)$; and the following diagrams commute for all $a, b, c, d, x, y \in \mathcal{D}$:

$$\begin{array}{ccc}
\perp \circ (a \bullet b) & \xrightarrow{\nu \circ \mathrm{id}} & (\perp \bullet \perp) \circ (a \bullet b) \\
 & & \downarrow{\scriptstyle R} \\
\downarrow{\scriptstyle \lambda} & & (T\perp \circ Ta) \bullet (T\perp \circ Tb) \\
 & & \downarrow{\scriptstyle (T_0^\circ \circ \alpha) \bullet (T_0^\circ \circ \beta)} \\
a \bullet b & \xrightarrow[\lambda^{-1} \bullet \lambda^{-1}]{} & (\perp \circ a) \bullet (\perp \circ b)
\end{array}
\qquad
\begin{array}{ccc}
(a \bullet b) \circ \perp & \xrightarrow{\mathrm{id} \circ \nu} & (a \bullet b) \circ (\perp \bullet \perp) \\
 & & \downarrow{\scriptstyle R} \\
\downarrow{\scriptstyle \lambda} & & (Ta \circ T\perp) \bullet (Tb \circ T\perp) \\
 & & \downarrow{\scriptstyle (\alpha \circ T_0^\circ) \bullet (\beta \circ T_0^\circ)} \\
a \bullet b & \xrightarrow[\lambda^{-1} \bullet \lambda^{-1}]{} & (a \circ \perp) \bullet (b \circ \perp)
\end{array}
\qquad (7.2.1)$$

$$\begin{array}{ccc}
(1 \bullet a) \circ (1 \bullet b) & \xrightarrow{R} & (T1 \circ T1) \bullet (Ta \circ Tb) \\
 & & \downarrow{\scriptstyle (T_0^\bullet \circ T_0^\bullet) \bullet (\alpha \circ \beta)} \\
\downarrow{\scriptstyle \lambda \circ \lambda} & & (1 \circ 1) \bullet (a \circ b) \\
 & & \downarrow{\scriptstyle \varpi \bullet \mathrm{id}} \\
a \circ b & \xrightarrow[\lambda]{} & 1 \bullet (a \circ b)
\end{array}
\qquad
\begin{array}{ccc}
(a \bullet 1) \circ (b \bullet 1) & \xrightarrow{R} & (Ta \circ Tb) \bullet (T1 \circ T1) \\
 & & \downarrow{\scriptstyle (\alpha \circ \beta) \bullet (T_0^\bullet \circ T_0^\bullet)} \\
\downarrow{\scriptstyle \lambda \circ \lambda} & & (a \circ b) \bullet (1 \circ 1) \\
 & & \downarrow{\scriptstyle \mathrm{id} \bullet \varpi} \\
a \circ b & \xrightarrow[\lambda]{} & (a \circ b) \bullet 1
\end{array}
\qquad (7.2.2)$$

$$\begin{array}{ccccc}
T((a \bullet b) \circ (c \bullet d)) & \xrightarrow{TR_{a,b,c,d}} & T((Ta \circ Tc) \bullet (Tb \circ Td)) & \xrightarrow{T^\bullet_{2,Ta \circ Tc, Tb \circ Td}} & T(Ta \circ Tc) \bullet T(Tb \circ Td) \\
\downarrow{\scriptstyle T^\circ_{2,a \bullet b, c \bullet d}} & & & & \downarrow{\scriptstyle T^\circ_{2,Ta,Tc} \bullet T^\circ_{2,Tb,Td}} \\
T(a \bullet b) \circ T(c \bullet d) & & & & (T^2 a \circ T^2 c) \bullet (T^2 b \circ T^2 d) \\
\downarrow{\scriptstyle T^\bullet_{2,a,b} \circ T^\bullet_{2,c,d}} & & & & \downarrow{\scriptstyle (\mu_a \circ \mu_c) \bullet (\mu_b \circ \mu_d)} \\
(Ta \bullet Tb) \circ (Tc \bullet Td) & \xrightarrow[R_{Ta,Tb,Tc,Td}]{} & (T^2 a \circ T^2 c) \bullet (T^2 b \circ T^2 d) & \xrightarrow[(\mu_a \circ \mu_c) \bullet (\mu_b \circ \mu_d)]{} & (Ta \circ Tc) \bullet (Tb \circ Td)
\end{array}
\qquad (7.2.3)$$

$$
\begin{array}{ccc}
(a \bullet b) \circ (c \bullet d) \circ (x \bullet y) & \xrightarrow{\mathrm{id} \circ R_{c,d,x,y}} & (a \bullet b) \circ ((Tc \circ Tx) \bullet (Td \circ Ty)) \\
\downarrow{\scriptstyle R_{a,b,c,d} \circ \mathrm{id}} & & \downarrow{\scriptstyle R_{a,b,Tc\circ Tx,Td\circ Ty}} \\
((Ta \circ Tc) \bullet (Tb \circ Td)) \circ (x \bullet y) & & (Ta \circ T(Tc \circ Tx)) \bullet (Tb \circ T(Td \circ Ty)) \\
\downarrow{\scriptstyle R_{Ta\circ Tc,Tb\circ Td,x,y}} & & \downarrow{\scriptstyle (Ta \circ T^{\circ}_{2,Tc,Tx}) \bullet (Tb \circ T^{\circ}_{2,Td,Ty})} \\
(T(Ta \circ Tc) \circ Tx) \bullet (T(Tb \circ Td) \circ Ty) & & (Ta \circ (T^2c \circ T^2x)) \bullet (Tb \circ (T^2d \circ T^2y)) \\
\downarrow{\scriptstyle (T^{\circ}_{2,Ta,Tc} \circ Tx) \bullet (T^{\circ}_{2,Tb,Td} \circ Ty)} & & \\
((T^2a \circ T^2c) \circ Tx) \bullet ((T^2b \circ T^2d) \circ Ty) & & \downarrow{\scriptstyle (Ta \circ \mu_c \circ \mu_x) \bullet (Tb \circ \mu_d \circ \mu_y)} \\
\downarrow{\scriptstyle (\mu_a \circ \mu_c \circ Tx) \bullet (\mu_b \circ \mu_d \circ Ty)} & & \\
((Ta \circ Tc) \circ Tx) \bullet ((Tb \circ Td) \circ Ty) & \xrightarrow[\cong]{} & (Ta \circ (Tc \circ Tx)) \bullet (Tb \circ (Td \circ Ty))
\end{array}
\tag{7.2.4}
$$

$$
\begin{array}{ccc}
((x \bullet a) \bullet c) \circ ((y \bullet b) \bullet d) & \xrightarrow{\cong} & (x \bullet (a \bullet c)) \circ (y \bullet (b \bullet d)) \\
\downarrow{\scriptstyle R_{x\bullet a,c,y\bullet b,d}} & & \downarrow{\scriptstyle R_{x,a\bullet c,y,b\bullet d}} \\
(T(x \bullet a) \circ T(y \bullet b)) \bullet (Tc \circ Td) & & (Tx \circ Ty) \bullet (T(a \bullet c) \circ T(b \bullet d)) \\
\downarrow{\scriptstyle T^{\bullet}_{2,x,a} \circ T^{\bullet}_{2,y,b} \bullet \mathrm{id}} & & \downarrow{\scriptstyle \mathrm{id} \bullet (T^{\bullet}_{2,a,c} \circ T^{\bullet}_{2,b,d})} \\
((Tx \bullet Ta) \circ (Ty \bullet Tb)) \bullet (Tc \circ Td) & & (Tx \circ Ty) \bullet ((Ta \bullet Tc) \circ (Tb \bullet Td)) \\
\downarrow{\scriptstyle R_{Tx,Ta,Ty,Tb} \bullet \mathrm{id}} & & \downarrow{\scriptstyle \mathrm{id} \bullet R_{Ta,Tc,Tb,Td}} \\
((T^2x \circ T^2y) \bullet (T^2a \circ T^2b)) \bullet (Tc \circ Td) & & (Tx \circ Ty) \bullet ((T^2a \circ T^2b) \bullet (T^2c \bullet T^2d)) \\
\downarrow{\scriptstyle (\mu_x \circ \mu_y) \bullet (\mu_a \circ \mu_b) \bullet \mathrm{id}} & & \downarrow{\scriptstyle \mathrm{id} \bullet ((\mu_a \circ \mu_b) \bullet (\mu_c \bullet \mu_d))} \\
((Tx \circ Ty) \bullet (Ta \circ Tb)) \bullet (Tc \circ Td) & \xrightarrow[\cong]{} & (Tx \circ Ty) \bullet ((Ta \circ Tb) \bullet (Tc \bullet Td))
\end{array}
\tag{7.2.5}
$$

A *quasitriangular* opmonoidal monad is one equipped with an R-matrix.

**Example 7.19.** Let $(\mathscr{C}, \otimes, 1)$ be a strict monoidal category, and $T$ a bimonad on $\mathscr{C}$. Let $R$ be an R-matrix on $T$ in the sense of [BV07, Section 8.2], and define

$$S := \{ \eta_a \otimes R_{b,c} \otimes \eta_d : a \otimes b \otimes c \otimes d \longrightarrow Ta \otimes Tc \otimes Tb \otimes Td \}_{a,b,c,d \in \mathscr{C}}.$$

Then $S$, together with $\nu$, $\varpi$, and $\iota$ being the identity, is an R-matrix on $T$, seen as a separately opmonoidal monad on the preduoidal category $\mathscr{C}$.

Diagrams (7.2.1) and (7.2.2) commute because $(\beta \otimes \alpha)R_{a,b}$ is a braiding by [BV07, Theorem 8.5]. Diagram (7.2.3) follows by Figure 7.1, where $T_3$ is defined as in Remark 2.35. The other diagrams are proved similarly.

By [BV07, Example 8.4] we also obtain that every R-matrix on a bialgebra $B$ yields an R-matrix on $B \otimes -$ in the sense of Definition 7.18.

**Remark 7.20.** Note that the converse of Example 7.19 is not necessarily true. Let $\mathscr{C}$ be a monoidal category seen as a preduoidal category, and assume that $T$ is a separately opmonoidal monad on $\mathscr{C}$ where the two oplax monoidal structures are the same. Then an R-matrix on $T$ does not necessarily yield an R-matrix in the sense of [BV07, Section 8.2], since we do not require $R$ to be $*$-*invertible*[17], which by [BV07, Theorem 8.5] corresponds bijectively to the braiding on $\mathscr{C}^T$ being invertible.

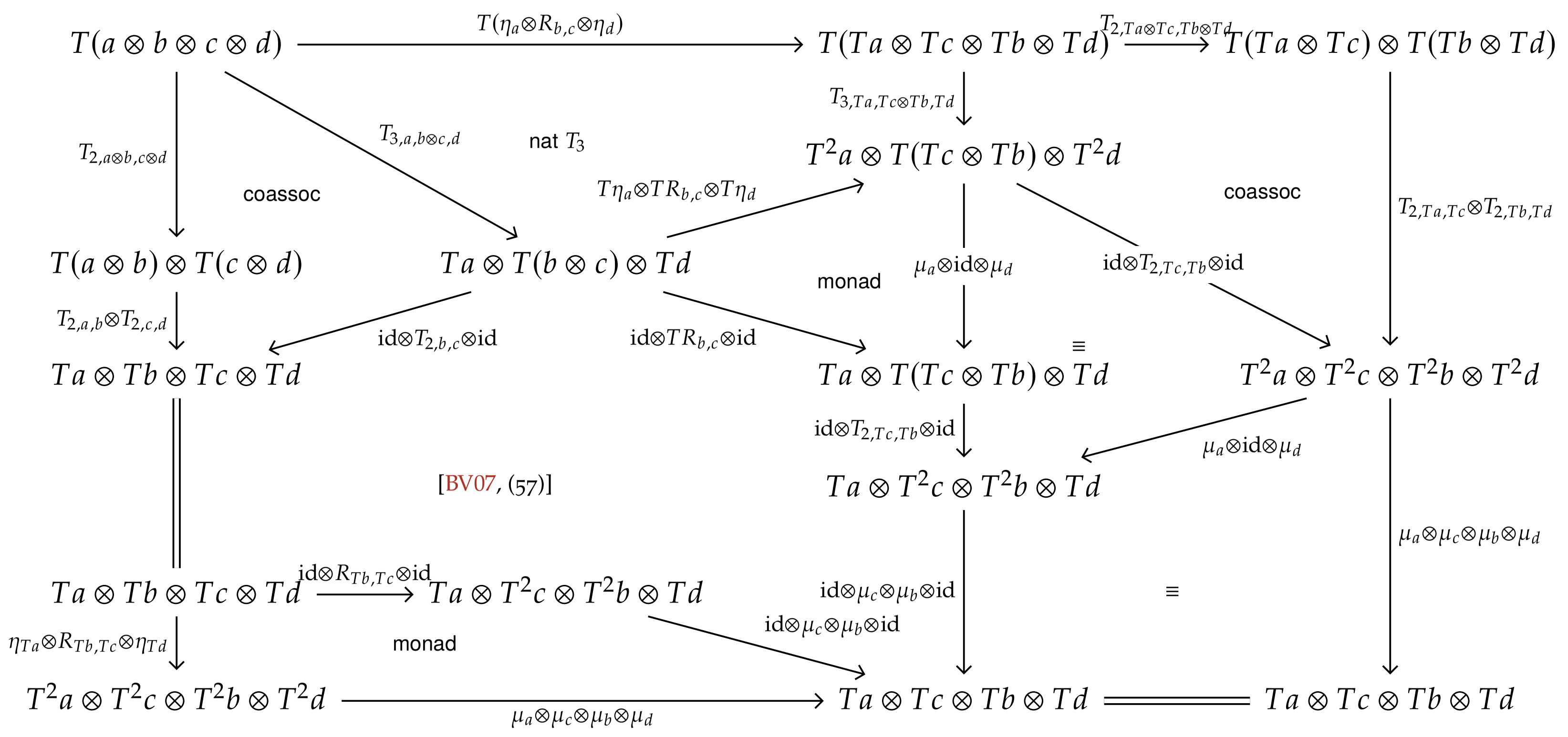

Figure 7.1: Verification that $S$ satisfies Diagram (7.2.3).

[17] A natural transformation $R\colon \otimes \Longrightarrow T \otimes^{\mathrm{op}} T$ is called $*$-*invertible* if there exists an "inverse" natural transformation $R^{-1}\colon \otimes^{\mathrm{op}} \Longrightarrow T \otimes T$, such that $(\mu \otimes \mu) \circ R^{-1} \circ R$ is equal to $\eta \otimes \eta$, and similarly for the other direction.

By Theorem 7.21 below, the R-matrices of Definition 7.18 correspond to duoidal structures on $\mathscr{C}^T$. Since the two tensor products on $\mathscr{C}^T$ agree, by arguments analogous to those in [AM10, Section 6.3], this forces the interchange law to come from a lax braiding. However, there is no a priori reason for this morphism to be invertible, see Example 7.4.

### 7.2.1 *From R-matrices to duoidal structures and back*

This section contains the main result of the chapter, which can be seen as an analogue of [BV07, Theorem 8.5], and a non-cocommutative counterpart to [AHLF18, Theorem 7.2].

**Theorem 7.21.** *Let $\mathfrak{D}$ be a preduoidal category and suppose that $T$ is a separately opmonoidal monad on $\mathfrak{D}$. For all $T$-algebras $(a,\alpha)$, $(b,\beta)$, $(c,\gamma)$, and $(d,\delta)$, a quasitriangular structure on $T$ yields an interchange law*

$$\xi := ((\alpha \circ \gamma) \bullet (\beta \circ \delta))R_{a,b,c,d}\colon (a \bullet b) \circ (c \bullet d) \longrightarrow (a \circ c) \bullet (b \circ d)$$

*on $\mathscr{C}^T$. Conversely, an interchange law $\xi$ on $\mathfrak{D}^T$ gives rise to an R-matrix*

$$R := \xi_{Ta,Tb,Tc,Td}((\eta_a \bullet \eta_b) \circ (\eta_c \bullet \eta_d))\colon (a \bullet b) \circ (c \bullet d) \longrightarrow (Ta \circ Tc) \bullet (Tb \circ Td)$$

*on $T$. These constructions are mutually inverse to each other.*

We split up the proof of Theorem 7.21 into its respective directions. Given these results, the rest of the proof is straightforward.

*Proof.* Combining Propositions 7.22 and 7.23 below, it is left to prove that the constructions are mutually inverse. In one direction we calculate

$$\begin{aligned} &((\alpha \circ \gamma) \bullet (\beta \circ \delta))\xi_{Ta,Tb,Tc,Td}((\eta_a \bullet \eta_b) \circ (\eta_c \bullet \eta_d)) \\ &\quad = ((\alpha \circ \gamma) \bullet (\beta \circ \delta))((\eta_a \circ \eta_c) \bullet (\eta_b \circ \eta_d))\xi_{a,b,c,d} && \text{naturality of } \xi \\ &\quad = ((\alpha\eta_a \circ \gamma\eta_c) \bullet (\beta\eta_b \circ \delta\eta_d))\xi_{a,b,c,d} && \text{functoriality of } \bullet \text{ and } \circ \\ &\quad = \xi_{a,b,c,d} && \text{monadicity of } T; \end{aligned}$$

and for the converse we have

$$\begin{aligned} &((\mu_a \circ \mu_c) \bullet (\mu_b \circ \mu_d))R_{Ta,Tb,Tc,Td}((\eta_a \bullet \eta_b) \circ (\eta_c \bullet \eta_d)) \\ &\qquad = ((\mu_a \circ \mu_c) \bullet (\mu_b \circ \mu_d))((\eta_{Ta} \circ \eta_{Tc}) \bullet (\eta_{Tb} \circ \eta_{Td}))R_{a,b,c,d} \\ &\qquad = R_{a,b,c,d}. \end{aligned}$$

□

**Proposition 7.22.** *Let $\mathcal{D}$ be a preduoidal category and $T$ a quasitriangular separately opmonoidal monad on $\mathcal{D}$ with R-matrix $(R, \nu, \varpi, \iota)$. Then $\mathcal{C}^T$ is a duoidal category, with structure morphisms $\nu$, $\varpi$, and $\iota$, and interchange law*

$$\xi := ((\alpha \circ \gamma) \bullet (\beta \circ \delta))R_{a,b,c,d} \colon (a \bullet b) \circ (c \bullet d) \longrightarrow (a \circ c) \bullet (b \circ d)$$

*for all $(a, \alpha)$, $(b, \beta)$, $(c, \gamma)$, and $(d, \delta) \in \mathcal{D}^T$.*

*Proof.* The claim that $\xi \in \mathcal{D}^T((a \bullet b) \circ (c \bullet d), (a \circ c) \bullet (b \circ d))$—that is, $\xi$ is a morphism of $T$-algebras—follows from Diagram (7.2.3), as seen in Figure 7.2. Diagram (7.1.1) follows by the commutativity of Figure 7.3, where we have left out the respective associators for readability; see Proposition 7.5. The proof of Diagram (7.1.2) is analogous. Lastly, Diagrams (7.2.1) and (7.2.2) immediately imply Diagram (7.1.3). □

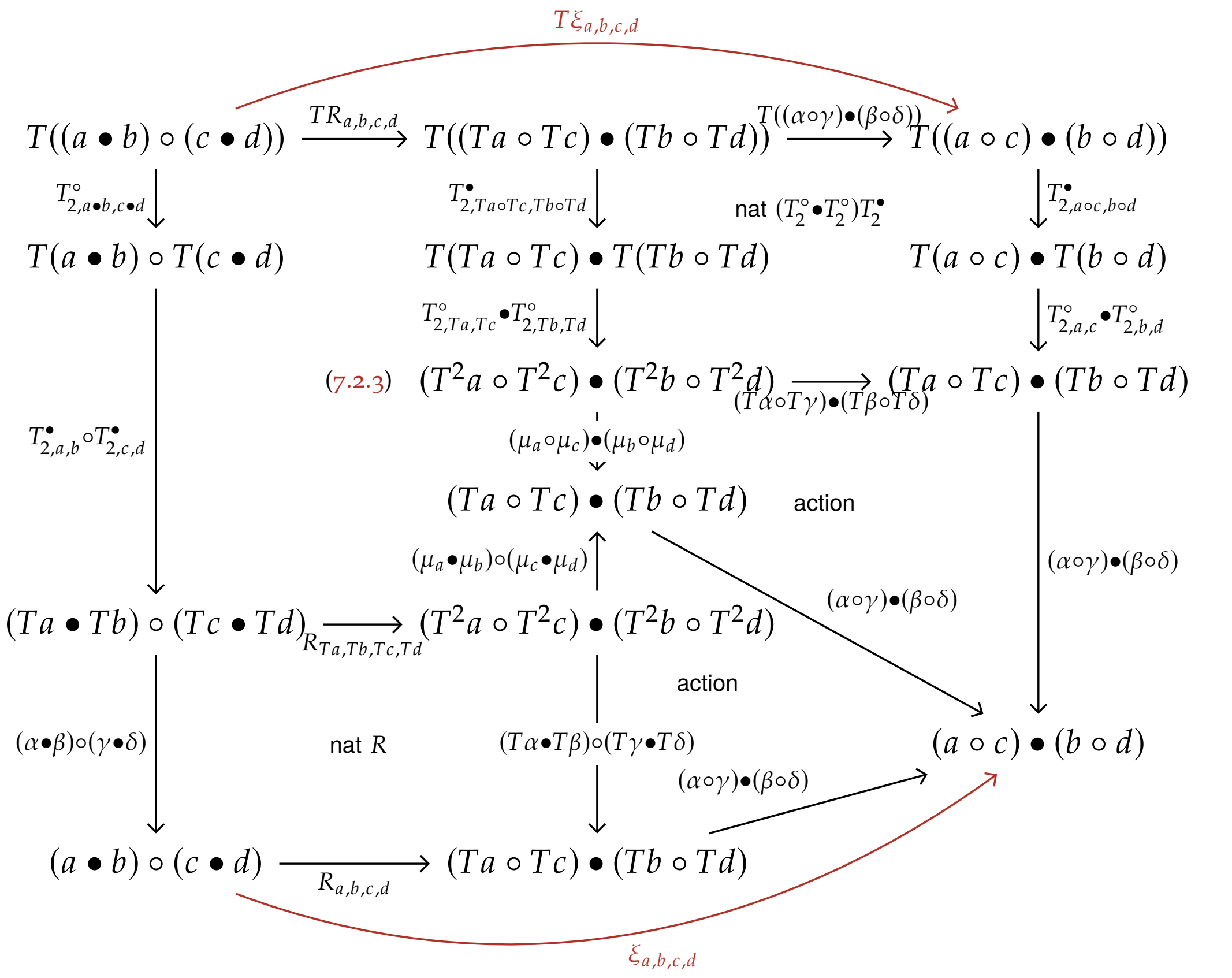

Figure 7.2: Proof that $\xi$ is a morphism of $T$-algebras.

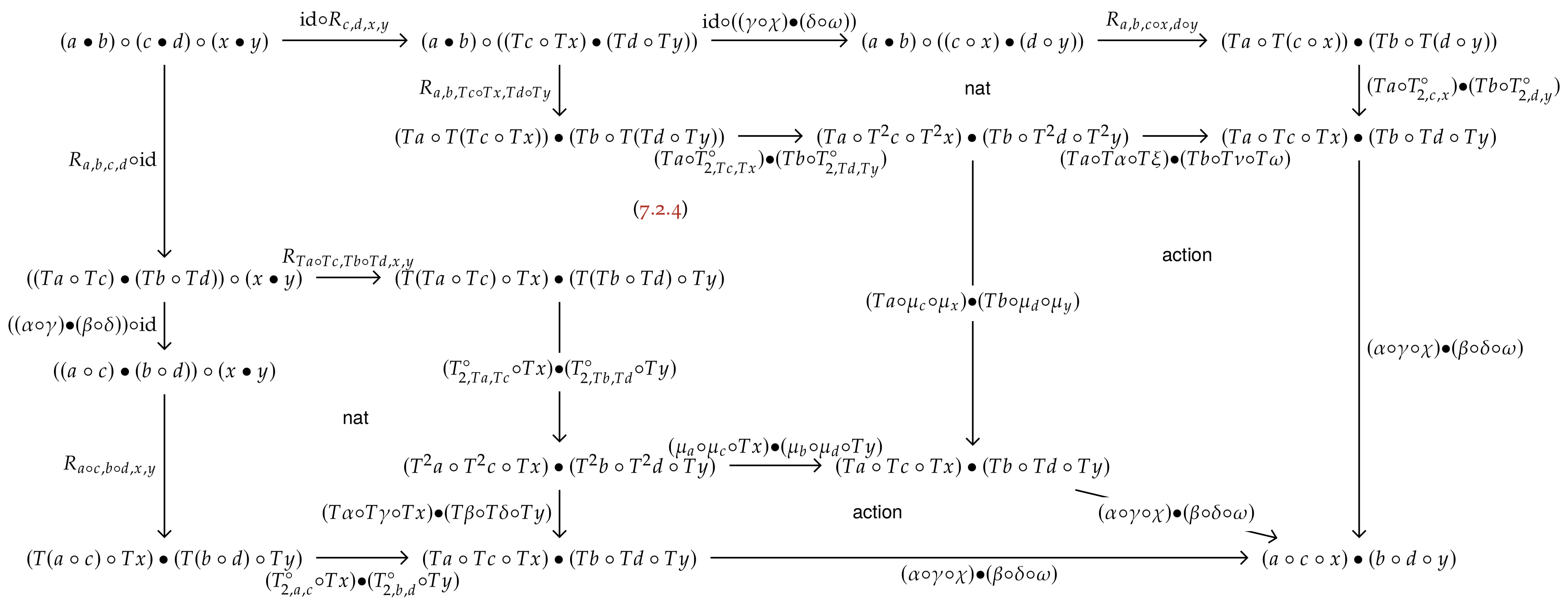

Figure 7.3: Proof that $\xi$ satisfies Diagram (7.1.1).

**Proposition 7.23.** *Let $\mathfrak{D}$ be a preduoidal category, $T$ a separately opmonoidal monad on $\mathfrak{D}$, and suppose that $\mathfrak{D}^T$ is a duoidal category with interchange law*

$$\xi_{a,b,c,d}\colon (a \bullet b) \circ (c \bullet d) \longrightarrow (a \circ c) \bullet (b \circ d).$$

*Then the structure morphisms of $\mathfrak{D}^T$, together with*

$$R := \xi_{Ta,Tb,Tc,Td}((\eta_a \bullet \eta_b) \circ (\eta_c \bullet \eta_d))\colon (a \bullet b) \circ (c \bullet d) \longrightarrow (Ta \circ Tc) \bullet (Tb \circ Td)$$

*yield an R-matrix for $T$.*

*Proof.* First, let us verify that $R$ satisfies Diagram (7.2.3). Let $a, b, c, d \in \mathfrak{D}^T$; then the claim follows from the commutativity of Figure 7.4. The fact that Diagram (7.2.4) holds is due to Figure 7.5, and Diagram (7.2.5) is similar.

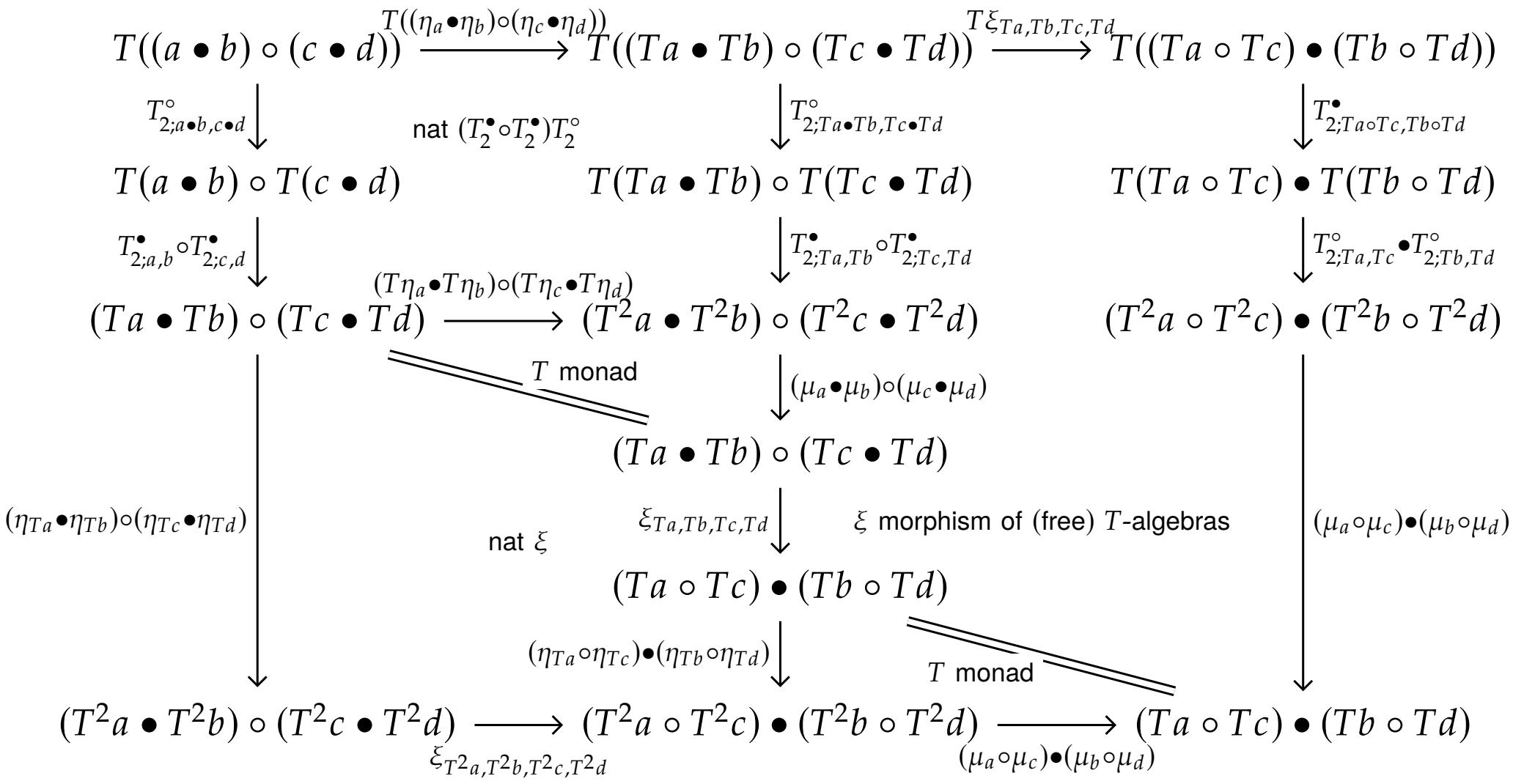

Figure 7.4: The map $R$ satisfies Diagram (7.2.3).

It is left to show the commutativity of Diagrams (7.2.1) and (7.2.2). For example, the first diagram in the former follows by the commutativity of

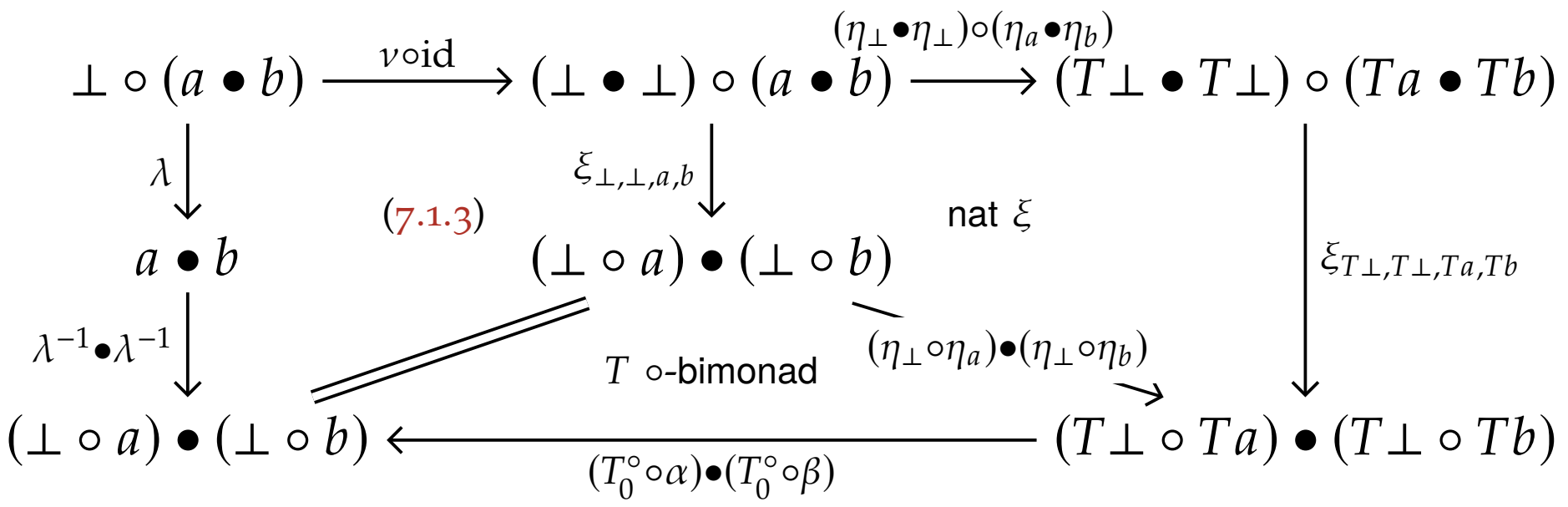

and the other diagrams are similar. ◻

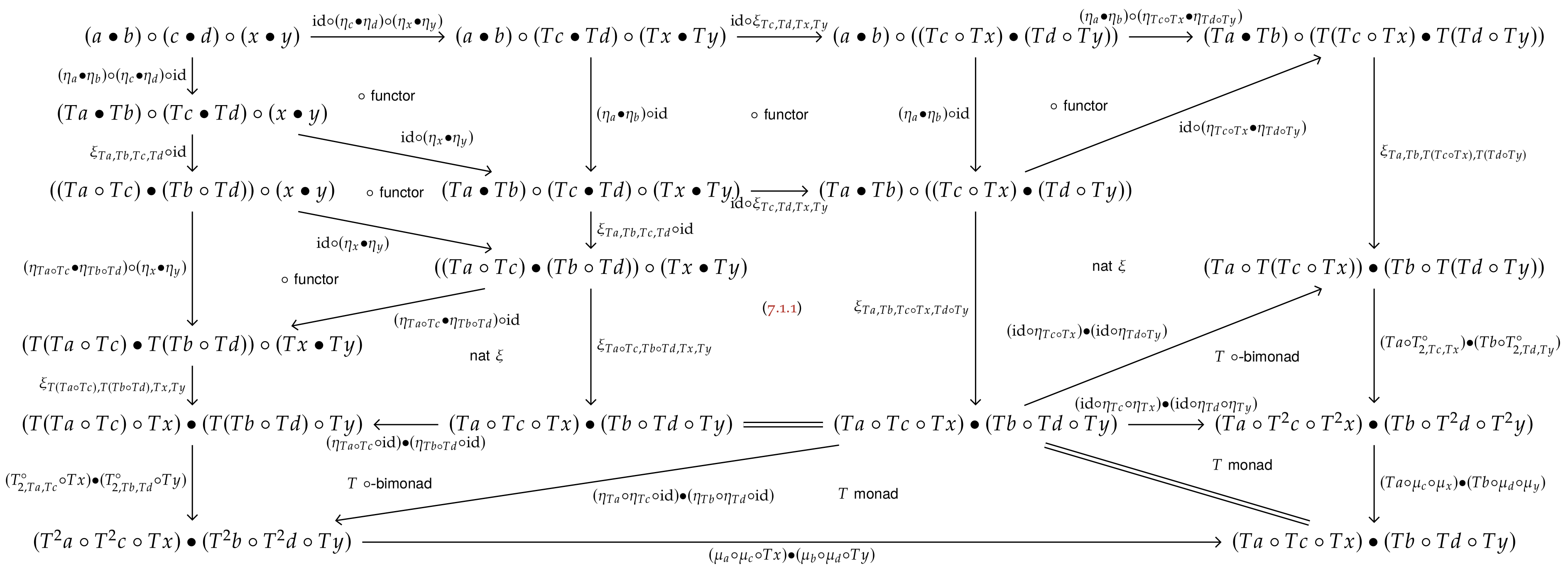

Figure 7.5: The R-matrix satisfies Diagram (7.2.4).

## 7.3 LINEARLY DISTRIBUTIVE MONADS

THE NORMAL DUOIDAL CATEGORIES of Definition 7.2 have connections to linear logic: in [GLF16, p. 7] it is shown that every normal duoidal category $\mathfrak{D}$ has the structure of a *linearly distributive category*; see [CS97]. In that case, the linear distributors are given by

$$\begin{aligned}
\partial^\ell_\ell &: a \circ (b \bullet c) \cong (a \bullet 1) \circ (b \bullet c) \xrightarrow{\zeta} (a \circ b) \bullet (1 \bullet c) \cong (a \circ b) \bullet c, \\
\partial^\ell_r &: a \circ (b \bullet c) \cong (1 \bullet a) \circ (b \bullet c) \xrightarrow{\zeta} (1 \circ b) \bullet (a \circ c) \cong b \bullet (a \circ c), \\
\partial^r_\ell &: (b \bullet c) \circ a \cong (b \bullet c) \circ (a \bullet 1) \xrightarrow{\zeta} (b \circ a) \bullet (c \circ 1) \cong (b \circ a) \bullet c. \\
\partial^r_r &: (b \bullet c) \circ a \cong (b \bullet c) \circ (1 \bullet a) \xrightarrow{\zeta} (b \circ 1) \bullet (c \circ a) \cong b \bullet (c \circ a).
\end{aligned} \tag{7.3.1}$$

**Remark 7.24.** By [MP22, Theorem 5.18], normal duoidal categories satisfy a much stronger form of coherence, so structures on them require fewer axioms to be fully specified. If $T$ is a double opmonoidal monad on a normal duoidal category $(\mathfrak{D}, \circ, \bot, \bullet, 1)$, then the following diagram commutes:

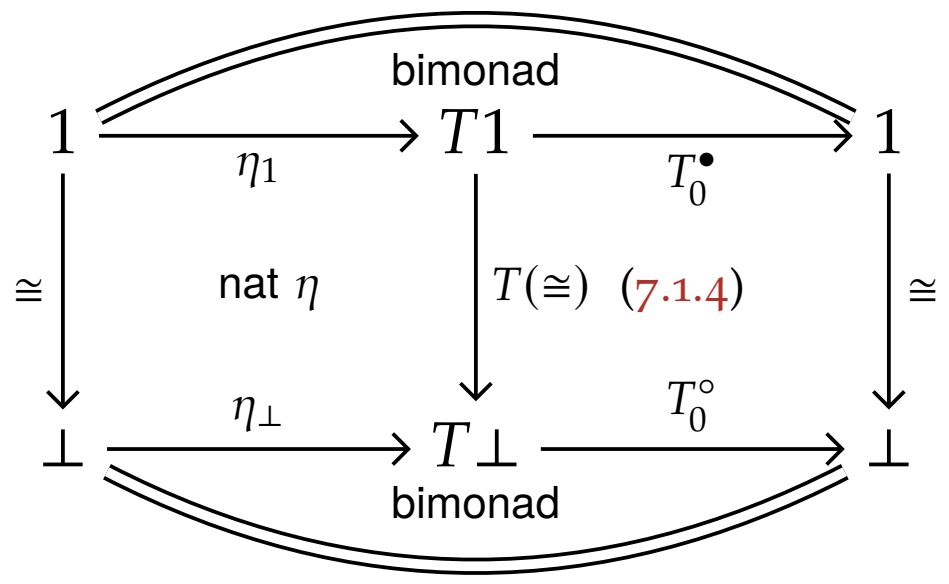

In particular, $T_0^\bullet$ and $T_0^\circ$ are conjugations of each other by isomorphisms:

$$(T1 \xrightarrow{T_0^\bullet} 1) = (T1 \xrightarrow{T(\cong)} T\bot \xrightarrow{T_0^\circ} \bot \xrightarrow{\cong^{-1}} 1).$$

In this setting, Diagram (7.1.4) automatically holds. For simplicity, assume $\mathfrak{D}$ to be strict, and write $T_0 := T_0^\bullet = T_0^\circ$. Then we for example have

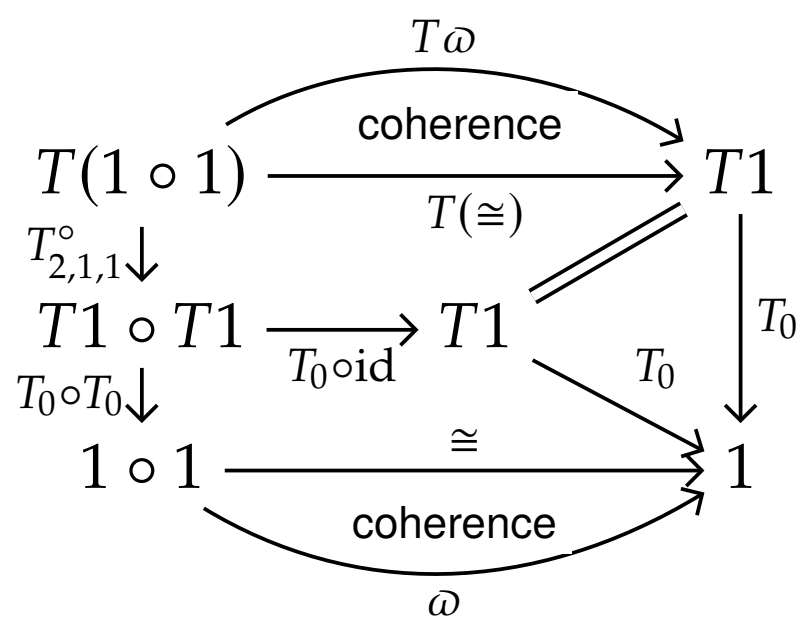

(7.3.2)

**Remark 7.25.** Sometimes, one considers only so-called *non-planar* linearly distributive categories, see [CS97, Section 2.1]. These are categories in which only $\partial_\ell^\ell$ and $\partial_r^r$ of Equation (7.3.1) exist. What we call a linearly distributive category is referred to as a *planar* linearly distributive category in *ibid*.

Conditions for a comonad to lift the (non-planar) linear distributive structure of its base category to its category of coalgebras were defined in [Pas12, Proposition 2.1]. For the convenience of the reader, the next proposition expresses this relation in terms of monads.

**Proposition 7.26.** *Let* $(\mathcal{L}, \otimes, \odot)$ *be a non-planar linearly distributive category, and suppose that the monad* $T$ *on* $\mathcal{L}$ *is separately opmonoidal. If the diagrams*

$$\begin{array}{ccccc} T(a\otimes(b\odot c)) & \xrightarrow{T^{\otimes}_{2,a,b\odot c}} & Ta\otimes T(b\odot c) & \xrightarrow{Ta\otimes T^{\odot}_{2,b,c}} & Ta\otimes(Tb\odot Tc) \\ \downarrow{\scriptstyle T\partial_l} & & & & \downarrow{\scriptstyle \partial_l} \\ T((a\otimes b)\odot c) & \xrightarrow[T^{\odot}_{2,a\otimes b,c}]{} & T(a\otimes b)\odot Tc & \xrightarrow[T^{\otimes}_{2,a,b}\odot Tc]{} & (Ta\otimes Tb)\odot Tc \end{array} \tag{7.3.3}$$

$$\begin{array}{ccccc} T((b\odot c)\otimes a) & \xrightarrow{T^{\otimes}_{2,b\odot c,a}} & T(b\odot c)\otimes Ta & \xrightarrow{T^{\odot}_{2,b,c}\otimes Ta} & (Tb\odot Tc)\otimes Ta \\ \downarrow{\scriptstyle T\partial_r} & & & & \downarrow{\scriptstyle \partial_r} \\ T(b\odot(c\otimes a)) & \xrightarrow[T^{\odot}_{2,b,c\otimes a}]{} & Tb\odot T(c\otimes a) & \xrightarrow[Tb\odot T^{\otimes}_{2,c,a}]{} & Tb\odot(Tc\otimes Ta) \end{array} \tag{7.3.4}$$

*commute for all* $T$*-algebras* $a$, $b$, *and* $c$, *then* $\mathcal{L}^T$ *is non-planar linearly distributive.*

**Example 7.27.** Setting $\otimes = \odot$, every monoidal category $\mathcal{C}$ is a linearly distributive category. The linear distributors are given by the associator and its inverse. Hence, a bimonad $B$ on $\mathcal{C}$ satisfies all assumptions of Proposition 7.26. Diagrams (7.3.3) and (7.3.4) reduce to the coassociativity of $B_2$.

Lifting the interchange morphism of a normal duoidal category is more involved than lifting only the non-planar linear distributors, much like lifting the preduoidal structure is easier than lifting the entire duoidal structure.

**Example 7.28.** Let $\mathcal{C}$ be a braided monoidal category, which is normal duoidal by Example 7.3. As such, the linear distributor $\partial_\ell^\ell$ is the isomorphism

$$\partial_\ell^\ell\colon x\otimes(y\otimes z)\cong x\otimes(1\otimes y)\otimes z \xrightarrow{x\otimes\sigma_{1,y}\otimes z} x\otimes(y\otimes 1)\otimes z\cong x\otimes(y\otimes z),$$

and $\partial_r^r$ is similar. By Proposition 7.26, this structure lifts to $\mathscr{C}^{(B\otimes -)}$, which is equal to the category of $B$-modules on $\mathscr{C}$. Analogously to Example 7.27, Diagrams (7.3.3) and (7.3.4) reduce to the coassociativity of $\Delta$. However, it is not true that the modules over an arbitrary bialgebra are braided monoidal; see for example [EGNO15, Example 8.3.5]. In other words, the planar structure

$$\partial_r^\ell \colon x\otimes(y\otimes z)\cong 1\otimes(x\otimes y)\otimes z \xrightarrow{1\otimes\sigma_{x,y}\otimes z} 1\otimes(y\otimes x)\otimes z\cong y\otimes(x\otimes z),$$
$$\partial_\ell^r \colon (x\otimes y)\otimes z\cong x\otimes(y\otimes z)\otimes 1 \xrightarrow{x\otimes\sigma_{y,z}\otimes 1} x\otimes(z\otimes y)\otimes 1\cong(x\otimes z)\otimes y,$$

does not lift to the category of $B$-modules.

As stated in the introduction, planar duoidal categories also capture and generalise the notion of a braiding, much like duoidal categories do. The following is a straightforward reformulation of Proposition 7.26.

**Proposition 7.29.** *Let $(\mathscr{L},\otimes,\odot)$ be a linearly distributive category with a separately opmonoidal monad $T$ on it. If, in addition to Diagrams (7.3.3) and (7.3.4), the following diagrams commute for all $T$-algebras $a$, $b$, and $c$:*

$$\begin{array}{ccccc} T(a\otimes(b\odot c)) & \xrightarrow{T^{\otimes}_{2;a,b\odot c}} & Ta\otimes T(b\odot c) & \xrightarrow{Ta\otimes T^{\odot}_{2;b,c}} & Ta\otimes(Tb\odot Tc) \\ {\scriptstyle T\partial_r^\ell}\big\downarrow & & & & \big\downarrow{\scriptstyle \partial_r^\ell} \\ T(b\odot(a\otimes c)) & \xrightarrow[T^{\odot}_{2;b,a\otimes c}]{} & Tb\odot T(a\otimes c) & \xrightarrow[Tb\odot T^{\otimes}_{2;a,c}]{} & Tb\odot(Ta\otimes Tc) \end{array} \tag{7.3.5}$$

$$\begin{array}{ccccc} T((a\odot b)\otimes c) & \xrightarrow{T^{\otimes}_{2;a\odot b,c}} & T(a\odot b)\otimes Tc & \xrightarrow{T^{\odot}_{2;a,b}\otimes Tc} & (Ta\odot Tb)\otimes Tc \\ {\scriptstyle T\partial_\ell^r}\big\downarrow & & & & \big\downarrow{\scriptstyle \partial_r^\ell} \\ T(a\odot(c\otimes b)) & \xrightarrow[T^{\odot}_{2;a,c\otimes b}]{} & Ta\odot T(c\otimes b) & \xrightarrow[Ta\odot T^{\otimes}_{2;c,b}]{} & Ta\odot(Tc\otimes Tb) \end{array} \tag{7.3.6}$$

*then $\mathscr{L}^T$ is linearly distributive.*

**Example 7.30.** Let $B\in\mathsf{Vect}$ be a bialgebra. Focusing on the planar linear distributor $\partial_r^\ell$, for all $b\in B$, $x\in a$, $y\in b$, and $z\in c$, Diagram (7.3.5) becomes

$$b_{(1)}\otimes y\otimes b_{(2)}\otimes x\otimes b_{(3)}\otimes z = b_{(2)}\otimes y\otimes b_{(1)}\otimes x\otimes b_{(3)}\otimes z,$$

which is easily seen to be equivalent to $b_{(2)}\otimes b_{(1)} = b_{(1)}\otimes b_{(2)}$.

Thus, linearly distributive monads seem to be connected to the double opmonoidal monads of Section 7.1.1.

**Proposition 7.31.** *Let* $(\mathfrak{D}, \bullet, \circ, 1)$ *be a normal duoidal category. Then double opmonoidal monads on* $\mathfrak{D}$ *are linear distributive bimonads on* $\mathfrak{D}$.

*Proof.* Let $T$ be a cocommutative bimonad on $\mathfrak{D}$ as a duoidal category. Then the left-left linear distributor $\partial_\ell^\ell$ is given by

$$a \circ (b \bullet c) \cong (a \bullet 1) \circ (b \bullet c) \xrightarrow{\zeta} (a \circ b) \bullet (1 \circ c) \cong (a \circ b) \bullet c.$$

Now, Diagram (7.3.3) is satisfied by the commutativity of Figure 7.6; Diagram (7.3.4) is similar. Diagram (7.3.5) is satisfied by

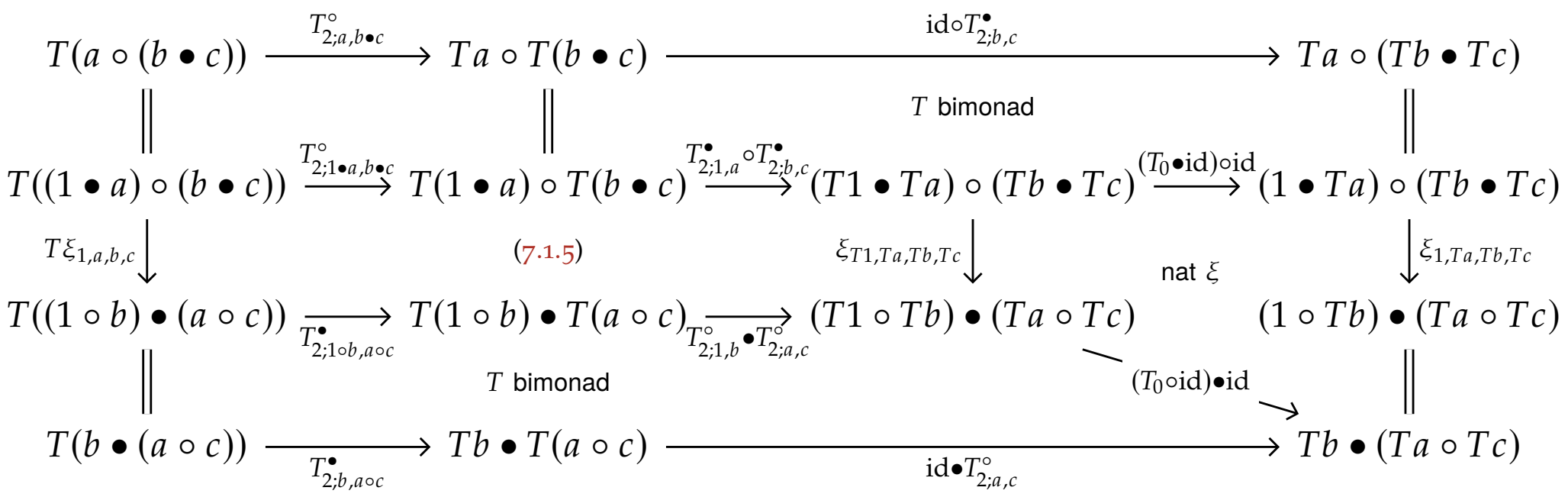

where we have assumed the normal duoidal structure to be strict for ease of readability. Diagram (7.3.6) follows similarly. □

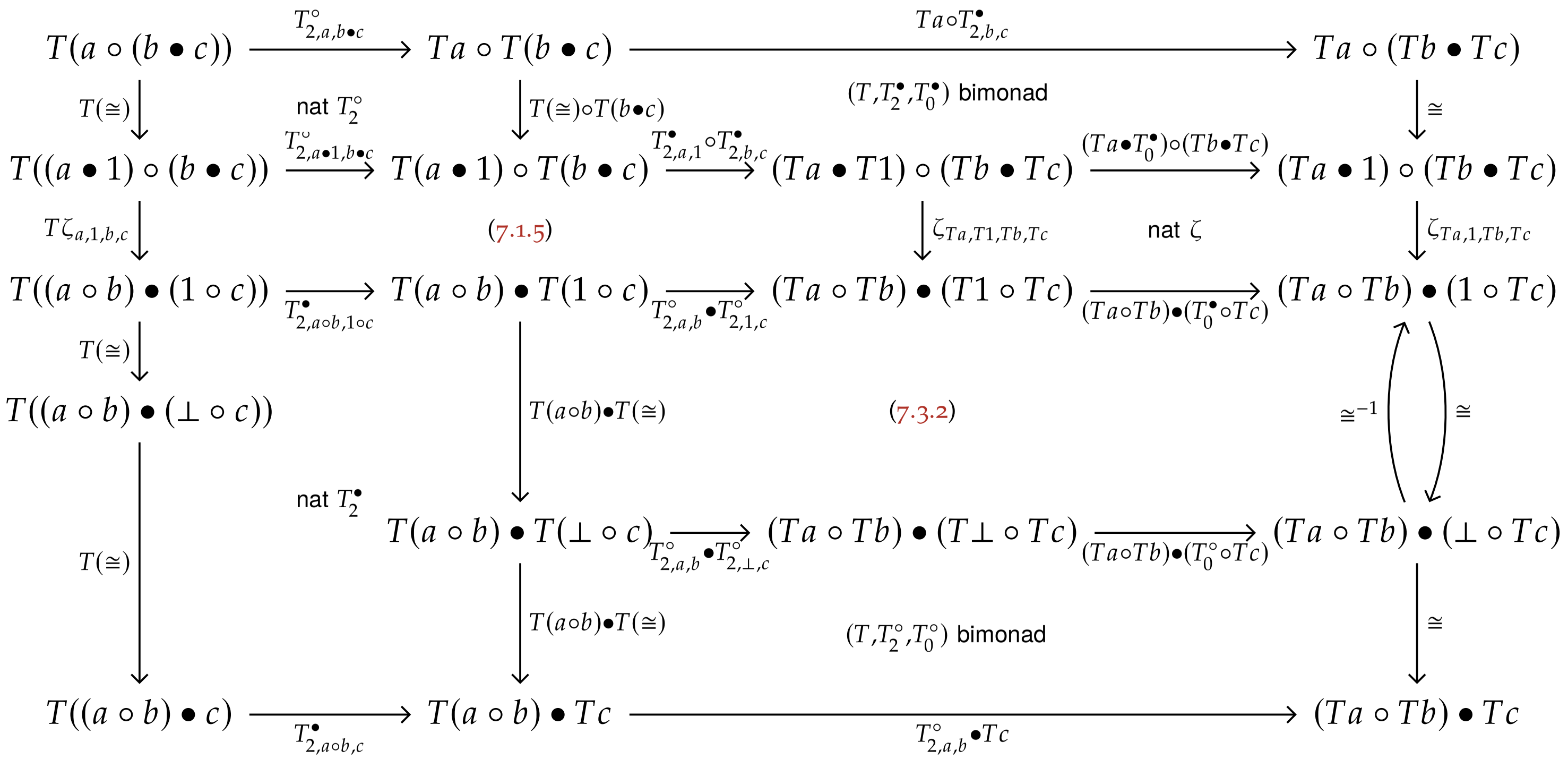

Figure 7.6: The left-left linear distributor satisfies Diagram (7.3.3).

Der Mensch, das tapferste und leidgewohnteste Thier, verneint an sich nicht das Leiden: er will es, er sucht es selbst auf, vorausgesetzt, dass man ihm einen Sinn dafür aufzeigt, ein Dazu des Leidens.

FRIEDRICH NIETZSCHE; Zur Genealogie der Moral.

# 8 INFINITE AND NON-RIGID RECONSTRUCTION

THE GOAL OF THIS CHAPTER IS to generalise the following theorem:

**Theorem A** ([Ost03, Theorem 1]). *Let $\mathscr{C}$ be a finite tensor category and let $\mathscr{M}$ be a finite abelian $\mathscr{C}$-module category, such that the evaluation functor $- \triangleright m : \mathscr{C} \longrightarrow \mathscr{M}$ is exact, for all $m \in \mathscr{M}$. Then there exists an algebra object $A \in \mathscr{C}$ such that there is an equivalence of $\mathscr{C}$-module categories $\mathrm{mod}_{\mathscr{C}}(A) \simeq \mathscr{M}$.*

This theorem makes several assumptions on $\mathscr{C}$ and $\mathscr{M}$:

- Finiteness assumptions: $\mathscr{C}$ and $\mathscr{M}$ are required to be finite abelian.
- The monoidal structure $- \otimes =$ of $\mathscr{C}$ and the action functor $- \triangleright =$ of $\mathscr{M}$ are both assumed to be exact in both variables.
- Rigidity assumptions on $\mathscr{C}$; i.e., that all of its objects admit left and right duals with respect to its monoidal structure.

We will generalise Theorem A in a way that greatly relaxes the first two kinds of assumptions, and removes the third one altogether. As mentioned before, by [DSPS19, Example 2.20], this reconstruction cannot solely rely on algebra objects in $\mathscr{C}$, and we have to approach Theorem A from a more monadic point of view. Our main result will be the following, and presents a certain *classification* of right exact lax $\mathscr{C}$-module monads.

**Theorems 8.48 to 8.50.** *Assume that $\mathscr{C}$ and $\mathscr{M}$ have enough projectives (injectives) and that there is an object $\ell \in \mathscr{M}$ such that:*

- *there is a right adjoint $\lfloor \ell, - \rfloor$ (left adjoint $\lceil \ell, - \rceil$) to $- \triangleright \ell$;*
- *for $x \in \mathscr{C}$ projective (injective), the object $x \triangleright \ell$ is projective (injective);*
- *any projective (injective) object of $\mathscr{M}$ is a direct summand of an object of the form $x \triangleright \ell$, for $x$ projective (injective).*

*Let $T$ be the monad $\lfloor \ell, - \triangleright \ell \rfloor$. Then $\mathcal{M} \simeq \mathcal{C}^T$, and the $\mathcal{C}$-module structure of the category of $T$-modules is extended from the Kleisli category. Furthermore, this extends to a bijection*

$$\{ (\mathcal{M}, \ell) \text{ as above} \} \big/ (\mathcal{M} \simeq \mathcal{N}) \leftrightarrow \left\{ \begin{matrix} \text{Right exact lax } \mathcal{C}\text{-module} \\ \text{monads on } \mathcal{C} \end{matrix} \right\} / \, (\mathcal{C}^T \simeq \mathcal{C}^S)$$

$$(\mathcal{M}, \ell) \longmapsto \lfloor \ell, - \triangleright \ell \rfloor$$

$$(\mathcal{C}^T, T1) \longleftarrow T$$

In the case of tensor categories, this statement may on the surface seem inapplicable: a tensor category that is not finite will generally not have enough projectives or injectives. However, since the category underlying a tensor category is equivalent to the category of finite-dimensional comodules for a coalgebra, the ind-completion of a tensor category has enough injectives. Further, since multitensor categories can be realised as categories of compact objects in locally finitely presentable categories, we may use an appropriate variant of the special adjoint functor theorem, see Propositions 2.133 and 2.134, characterising when the internal cohom exists for an object in $\mathsf{Ind}(\mathcal{M})$. As such, we are able to "pull back" the appropriate versions of Theorems 8.56 and 8.57 to the underlying categories.

**Theorem 8.59.** *Let $\mathcal{C}$ be a multitensor category, and let $\mathcal{M}$ be an abelian $\mathcal{C}$-module category such that the $\mathsf{Ind}(\mathcal{C})$-module category $\mathsf{Ind}(\mathcal{M})$ admits a coclosed $\mathsf{Ind}(\mathcal{C})$-injective $\mathsf{Ind}(\mathcal{C})$-cogenerator.*

*Then there is a coalgebra object $C$ in $\mathsf{Ind}(\mathcal{C})$ such that $\mathsf{Ind}(\mathcal{M}) \simeq \mathsf{Comod}_{\mathsf{Ind}(\mathcal{C})}C$. Further, $\mathcal{M}$ is the category of compact $C$-comodule objects..*

Consider the special case where a tensor category $\mathcal{C}$ admits a fibre functor to the category $\mathsf{vect}$ of finite-dimensional vector spaces, and is thus monoidally equivalent to the category ${}^H\mathsf{vect}$ of finite-dimensional comodules over a finite-dimensional Hopf algebra $H$, see [Ulb90]. Theorem A realises a finite $\mathcal{C}$-module category $\mathcal{M}$ as the category of finite-dimensional modules over a finite-dimensional $H$-comodule algebra. Since $H$ and $A$ are finite-dimensional, the category of $A$-modules is equivalent as an ${}^H\mathsf{vect}$-module category to the category of comodules for the $H$-comodule coalgebra $A^*$. In this sense, the following Hopf-theoretic corollary of the above theorem about ind-completions is an immediate generalisation of Theorem A.

**Corollary 8.60.** *Let $H$ be a Hopf algebra and let $\mathscr{C} = {}^{H}\mathsf{vect}$, meaning that there is a monoidal equivalence $\mathsf{Ind}(\mathscr{C}) \simeq {}^{H}\mathsf{Vect}$. Let $\mathscr{M}$ be an abelian $\mathscr{C}$-module category such that $\mathsf{Ind}(\mathscr{M})$ admits a coclosed $\mathsf{Ind}(\mathscr{C})$-injective $\mathsf{Ind}(\mathscr{C})$-cogenerator.*

*Then there exists a (possibly infinite-dimensional) $H$-comodule coalgebra $C$, such that there is an equivalence $\mathsf{Ind}(\mathscr{M}) \simeq \mathrm{Comod}_H C$ of $\mathsf{Ind}(\mathscr{C})$-module categories, restricting to a $\mathscr{C}$-module equivalence $\mathscr{M} \simeq \mathsf{comod}_H C$.*

Since we do not assume rigidity for our most general theorems, we instead invoke the bicategorical Yoneda lemma to realise the objects of $\mathscr{C}$ as the $\mathscr{C}$-module endofunctors of $\mathscr{C}$ itself. Then, we view the *lax* $\mathscr{C}$-module endofunctors of $\mathscr{C}$ as one appropriate generalisation of objects in $\mathscr{C}$, and the lax $\mathscr{C}$-module monads as the corresponding generalisation of algebra objects in $\mathscr{C}$. Using the reconstruction results of Section 5.3, the right adjoint of a (strong) $\mathscr{C}$-module functor is canonically lax, and thus the (co)monads we study are canonically (op)lax $\mathscr{C}$-module functors.

## 8.1 EXTENDING MODULE STRUCTURES

AN ADDITIONAL DIFFICULTY WE ENCOUNTER is that while the Kleisli category of a lax $\mathscr{C}$-module monad is naturally a $\mathscr{C}$-module category, the same is not true for the Eilenberg–Moore category. In this section, we shall investigate under which conditions the $\mathscr{C}$-module structure of the former lifts essentially uniquely to the latter. More precisely, we first establish uniqueness in Theorem 8.9, and then complement that with an existence result of such an extension in Theorem 8.25, under the assumption of right exactness (left exactness for comonads) of the (co)monad, using so-called Linton coequalisers and multicategorical techniques similar to those in [AHLF18].

**Hypothesis 8.1.** From now on until the end of this thesis, we implicitly assume all categories and functors to be $\Bbbk$-linear, for a field $\Bbbk$.

The fact that we only consider right exact lax $\mathscr{C}$-module monads is exemplified by the following fundamental fact.

**Proposition 8.2** ([BZBJ18, Proposition 3.2])**.** *Let $\mathscr{A}$ be abelian, $T$ a right exact monad on $\mathscr{A}$, and $S$ a left exact comonad on $\mathscr{A}$. Then $\mathscr{A}^T$ and $\mathscr{A}^S$ are abelian.*

**Proposition 8.3.** *Let $\mathscr{A}$ be a locally finite abelian category, and let $S$ be a left exact, finitary comonad on $\mathsf{Ind}(\mathscr{A})$. Then $\mathsf{Ind}(\mathscr{A})^S$ also is of the form $\mathsf{Ind}(\mathscr{E})$, for a locally finite abelian category $\mathscr{E}$. Further, $\mathscr{E}$ can be chosen to be the category of compact* *objects in* $\mathsf{Ind}(\mathcal{A})^S$*, and it can be characterised as the objects sent to compact objects under the forgetful functor* $F^S\colon \mathsf{Ind}(\mathcal{A})^S \longrightarrow \mathsf{Ind}(\mathcal{A})$*.*

*In particular, if* $S$ *is a left exact comonad on* $\mathcal{A}$*, we have* $\mathsf{Ind}(\mathcal{A})^{\mathsf{Ind}(S)} \simeq \mathsf{Ind}(\mathcal{A}^S)$ *and* $\mathcal{A}^S$ *is locally finite abelian.*

*Proof.* Let $D$ be the $\Bbbk$-coalgebra such that $\mathcal{A} \simeq {}^D\mathsf{vect}$. Then we also have that $\mathsf{Ind}(\mathcal{A}) \simeq {}^D\mathsf{Vect}$. In particular, $S$ is a comonad on ${}^D\mathsf{Vect}$. By Proposition 2.109, there is a $D$-$D$-bicomodule $C$ such that $S \cong C \,\square_D\, -$.

Similarly to [Tak77, Remark 2.4], under this isomorphism the comonad structure on $S$ corresponds to maps $C \longrightarrow C\square_D C$ and $C \longrightarrow D$, which endow $C$ with the structure of a $\Bbbk$-coalgebra, together with a coalgebra morphism $C \longrightarrow D$. Formally, a $S$-comodule in ${}^D\mathsf{Vect}$ is a $D$-comodule together with a $C$-comodule structure that restricts to the given $D$-comodule along the coalgebra morphism $C \longrightarrow D$. Morphisms of $S$-comodules in ${}^D\mathsf{Vect}$ are precisely $C$-comodule morphisms, and so we have

$$\mathsf{Ind}(\mathcal{A})^S \simeq ({}^D\mathsf{Vect})^S \simeq {}^C\mathsf{Vect}. \tag{8.1.1}$$

Thus we may set $\mathcal{C} = {}^C\mathsf{vect}$. The characterisation of compact objects in $({}^D\mathsf{Vect})^S$ in terms of the images of the functor $F^S$ follows immediately from observing that under the equivalence of Equation (8.1.1), the functor $F^S$ is simply the restriction functor ${}^C\mathsf{Vect} \longrightarrow {}^D\mathsf{Vect}$.

The second part of the statement follows by observing that $\mathcal{A}^S$ is the category of compact objects in $\mathsf{Ind}(\mathcal{A})^{\mathsf{Ind}(S)}$, by the first part of the statement, and by observing that $\mathcal{A}^S$ is also the category of compact objects in $\mathsf{Ind}(\mathcal{A}^S)$. Since both $\mathsf{Ind}(\mathcal{A}^S)$ and $\mathsf{Ind}(\mathcal{A})^{\mathsf{Ind}(S)}$ are locally finitely presentable, the equivalence between their respective categories of compact objects establishes an equivalence between the categories themselves. □

**Remark 8.4.** Suppose that $T$ is a monad on the category $\mathcal{C}$. For $T$-modules $(x, \alpha)$ and $(y, \beta)$, notice that there is a bijection

$$\begin{aligned} \mathcal{C}^T(Tx, y) &\cong \mathcal{C}(x, y) \\ (f\colon Tx \longrightarrow y) &\longmapsto (f \circ \eta_x) \\ (\beta \circ Tg) &\longleftarrow (g\colon x \longrightarrow y). \end{aligned} \tag{8.1.2}$$

This induces an isomorphism $\mathcal{C}^T(T(-), =) \cong \mathcal{C}(-, =)$. Using the definition of $\iota$ from Equation (2.2.2), the left-hand side is also obtained from the functor

$$\mathcal{C}^T \xrightarrow{\text{よ}} [(\mathcal{C}^T)^{\mathrm{op}}, \mathsf{Vect}] \xrightarrow{\iota^*} [\mathcal{C}_T{}^{\mathrm{op}}, \mathsf{Vect}].$$

The analogue of Proposition 8.3 for finite abelian categories is simpler.

**Proposition 8.5.** *Let $T$ be a right exact monad on a finite abelian category $\mathcal{A}$. Then $\mathcal{A}^T$ is a finite abelian category.*

*Further, any projective object in $\mathcal{A}^T$ is a direct summand of one of the form $\iota(P)$, where $\iota\colon \mathcal{A}_T \longrightarrow \mathcal{A}^T$ is the canonical embedding of Equation (2.2.2) and $P \in \mathcal{A}$-proj. Denoting the category of objects of this form by $\mathsf{Kl}_p(T)$, we obtain an equivalence*

$$\mathcal{A}^T\text{-proj} \simeq \overline{\mathsf{Kl}_p(T)}$$

*between the projectives of $\mathcal{A}^T$ and the Cauchy completion of $\mathsf{Kl}_p(T)$.*

*Proof.* Let $A$ be a finite-dimensional $\Bbbk$-algebra such that $\mathcal{A} \simeq A$-mod. Then there is an $A$-$A$-bimodule $B$ such that $T \cong B \otimes_A -$. Similarly to the proof of Proposition 8.3, monad structures on $T$ correspond to pairs consisting of a finite-dimensional $\Bbbk$-algebra structure on $B$ and an algebra homomorphism $A \longrightarrow B$. Under this correspondence we have $\mathcal{A}^T \simeq B$-mod. The latter category is clearly a finite abelian category.

Now suppose that $P \in \mathcal{A}$-proj. Applying Remark 8.4, one obtains an isomorphism $\mathcal{A}^T(\iota(P), -) \cong \mathcal{A}(P, -)$, proving the exactness of the left-hand side, and thus projectivity of $\iota(P)$. For $(X, \alpha) \in \mathcal{A}^T$, using the fact that $\mathcal{A}$ has enough projectives, we may fix an epimorphism $q\colon Q \twoheadrightarrow X$ in $\mathcal{A}$, where $Q \in \mathcal{A}$-proj. For the latter claim, notice that there is a composite epimorphism

$$\iota(P) \xrightarrow{\iota(q)} \iota(X) = TX \xrightarrow{\alpha} X,$$

the first part of which is epic by right exactness of $\iota$. □

**Lemma 8.6.** *Let $T$ be a monad on $\mathcal{A}$. The embedding $\mathcal{A}^T \hookrightarrow [\mathcal{A}_T^{\mathrm{op}}, \mathsf{Vect}]$ can be corestricted to an embedding $\mathcal{A}^T \hookrightarrow \mathrm{Fin}_{\mathrm{co}}(\mathcal{A}_T)$ into the finite cocompletion of $\mathcal{A}_T$.*

*Proof.* Let $(x, \alpha)$ and $(y, \beta)$ be two $T$-algebras. Recall that the coequaliser

$$T^2x \underset{T\alpha}{\overset{\mu_x}{\rightrightarrows}} Tx \xrightarrow{\alpha} x$$

in $\mathcal{A}^T$ is an absolute coequaliser[18] in $\mathcal{A}$; thus, its image under the Yoneda embedding よ yields a coequaliser

$$\mathcal{A}(-, T^2x) \underset{(T\alpha)_*}{\overset{(\mu_x)_*}{\rightrightarrows}} \mathcal{A}(-, Tx) \xrightarrow{\alpha_*} \mathcal{A}(-, x).$$

[18] A coequaliser is called *absolute* if it is preserved by any functor; see [Par69] for details.

Passing under the isomorphism of Equation (8.1.2), one obtains a coequaliser

$$\mathscr{A}^T(T(-),T^2x) \underset{(T\alpha)_*}{\overset{(\mu_x)_*}{\rightrightarrows}} \mathscr{A}^T(T(-),Tx) \xrightarrow{\alpha_*} \mathscr{A}^T(T(-),x),$$

which, in turn, is isomorphic to

$$\mathscr{A}_T(-,Tx) \underset{(T\alpha)_*}{\overset{(\mu_x)_*}{\rightrightarrows}} \mathscr{A}_T(-,x) \xrightarrow{\alpha_*} \mathscr{A}^T(T(-),x).$$

This proves the result. □

**Proposition 8.7.** *Let $T$ be a right exact monad on an abelian category $\mathscr{A}$. The inclusion functor $\mathscr{A}^T \hookrightarrow \mathrm{Fin}_{\mathrm{co}}(\mathscr{A}_T)$ has a left adjoint, and the counit of the adjunction is a natural isomorphism. Further, the left adjoint $\mathrm{Fin}_{\mathrm{co}}(\mathscr{A}_T) \longrightarrow \mathscr{A}^T$ is the right exact extension of the inclusion $\iota\colon \mathscr{A}_T \hookrightarrow \mathscr{A}^T$.*

*Proof.* By Proposition 8.2, the $\mathscr{A}^T$ is finitely cocomplete. Similarly to [KS06, Proposition 6.3.1], there is a functor $\mathrm{Fin}_{\mathrm{co}}(\mathscr{A}^T) \longrightarrow \mathscr{A}^T$ that extends the functor $\mathrm{Id}\colon \mathscr{A}^T \longrightarrow \mathscr{A}^T$, defined by sending some $x = \widetilde{\mathrm{colim}}_i\, x_i$ in $\mathrm{Fin}_{\mathrm{co}}(\mathscr{A}^T)$ to $\mathrm{colim}_i\, x_i$ in $\mathscr{A}^T$, where $\widetilde{\mathrm{colim}}_i$ denotes the formally added colimit. The functor $\mathrm{Fin}_{\mathrm{co}}(\mathscr{A}^T) \longrightarrow \mathscr{A}^T$ is left adjoint to the inclusion $\mathscr{A}^T \hookrightarrow \mathrm{Fin}_{\mathrm{co}}(\mathscr{A}^T)$.

On the other hand, there is an adjunction

(8.1.3) $$[\mathscr{A}_T^{\mathrm{op}}, \mathsf{Vect}] \underset{\iota^*}{\overset{\iota_!}{\underset{\longleftarrow}{\overset{\longrightarrow}{\perp}}}} [(\mathscr{A}^T)^{\mathrm{op}}, \mathsf{Vect}],$$

where $\iota_!$ is the left Kan extension of $\text{よ}_{\mathscr{A}^T} \circ \iota$ along $\text{よ}_{\mathscr{A}_T}$, see for example [Str24a, Corollary 3.3]. Equation (8.1.3) restricts to an adjunction

$$\mathrm{Fin}_{\mathrm{co}}(\mathscr{A}_T) \underset{\iota^*}{\overset{\iota_!}{\underset{\longleftarrow}{\overset{\longrightarrow}{\perp}}}} \mathrm{Fin}_{\mathrm{co}}(\mathscr{A}^T),$$

since, by definition, $\iota_!$ sends finite colimits of representables to finite colimits of representables, and $\iota^*$ has the same property, since it preserves colimits, and, by Lemma 8.6, sends representables to finite colimits of representables. Thus we have the following commutative diagram

$$\begin{array}{ccccc}
\mathscr{A}^T & \underset{\longrightarrow}{\overset{\longleftarrow}{\perp}} & \mathrm{Fin}_{\mathrm{co}}(\mathscr{A}^T) & \hookleftarrow & \mathscr{A}^T \\
 & \searrow & \iota_! \uparrow \dashv \downarrow \iota^* & & \uparrow \iota \\
 & & \mathrm{Fin}_{\mathrm{co}}(\mathscr{A}_T) & \hookleftarrow & \mathscr{A}_T
\end{array}$$

This realises the embedding $\mathscr{A}^T \hookrightarrow \mathrm{Fin}_{\mathrm{co}}(\mathscr{A}_T)$ as the composition of two right adjoints, and thus as a right adjoint. Lastly, the composite

$$\mathscr{A}_T \longrightarrow \mathscr{A}^T \longrightarrow \mathrm{Fin}_{\mathrm{co}}(\mathscr{A}^T) \longrightarrow \mathscr{A}^T$$

is naturally isomorphic to $\iota\colon \mathscr{A}_T \hookrightarrow \mathscr{A}^T$, which proves the latter statement. □

**Remark 8.8.** In other words, Proposition 8.7 says that for a right exact monad $T$ the category $\mathscr{A}^T$ is a *reflective subcategory*—one where the inclusion functor has a left adjoint—of $\mathrm{Fin}_{\mathrm{co}}(\mathscr{A}_T)$.

The next result proves the essential uniqueness of the module structure on the Eilenberg–Moore category, under the assumption that it is "induced" from the Kleisli category.

**Theorem 8.9.** *Let $\mathscr{C}_T$ denote the Kleisli category for a monad $T$ on $\mathscr{C}$, equipped with a fixed left $\mathscr{C}$-module structure. Equip $\mathscr{C}^T$ with two left $\mathscr{C}$-module category structures $\mathscr{C}_1^T$ and $\mathscr{C}_2^T$, such that the inclusion $\iota\colon \mathscr{C}_T \hookrightarrow \mathscr{C}^T$ gives strong $\mathscr{C}$-module functors $\mathscr{C}_T \longrightarrow \mathscr{C}_1^T$ and $\mathscr{C}_T \longrightarrow \mathscr{C}_2^T$.*

*Then there is an equivalence $\mathscr{C}_1^T \simeq \mathscr{C}_2^T$ of left $\mathscr{C}$-module categories.*

*Proof.* Consider the following diagram:

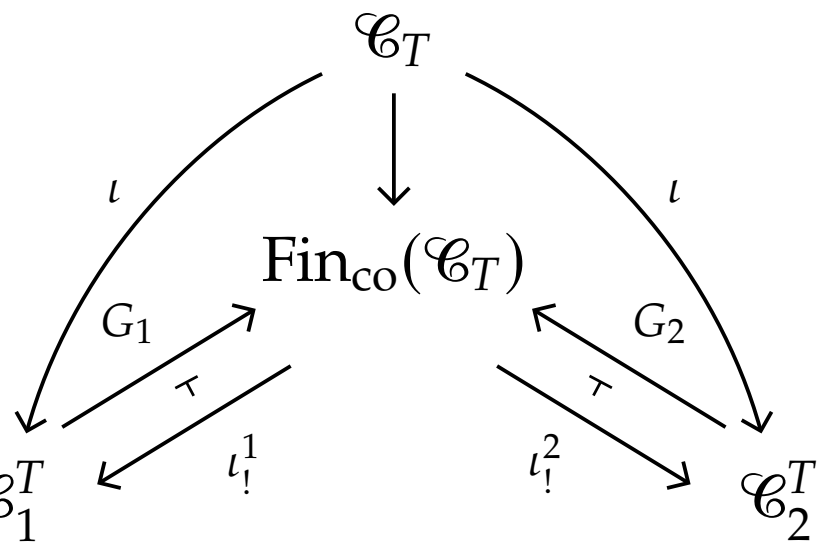

In particular, by Proposition 5.32, in this diagram every functor is a lax $\mathscr{C}$-module functor, and every adjunction is a $\mathscr{C}$-module adjunction. We claim that $\iota^2_! G_1\colon \mathscr{C}_1^T \longrightarrow \mathscr{C}_2^T$ and $\iota^1_! G_2\colon \mathscr{C}_2^T \longrightarrow \mathscr{C}_1^T$ define mutually quasi-inverse equivalences of $\mathscr{C}$-module categories. Indeed,

$$\iota^1_! G_2 \iota^2_! G_1 \xrightarrow{\varepsilon_2 \varepsilon_1} \mathrm{Id}_{\mathscr{C}_1^T}$$

is a $\mathscr{C}$-module isomorphism.

Similarly, there is a $\mathscr{C}$-module isomorphism $\mathrm{Id}_{\mathscr{C}_2^T} \cong \iota^2_! G_1 \iota^1_! G_2$. □

### 8.1.1 *Extendable monads*

**Definition 8.10.** We call a lax $\mathscr{C}$-module monad $T$ on a $\mathscr{C}$-module category $\mathscr{M}$ *extendable* if the category $\mathscr{M}^T$ is a $\mathscr{C}$-module category, for which the embedding $\iota\colon \mathscr{M}_T \hookrightarrow \mathscr{M}^T$ is a strong $\mathscr{C}$-module functor.

If $T$ is extendable, we refer to the essentially unique $\mathscr{C}$-module structure on $\mathscr{M}^T$ coming from $\mathscr{M}_T$ as the *extended* $\mathscr{C}$-module structure on $\mathscr{M}^T$.

Analogously to Definition 8.10, one defines extendable oplax $\mathscr{C}$-module comonads and extended module structures on their comodules.

**Proposition 8.11.** *Any strong $\mathscr{C}$-module monad $(T, \mu, \eta)\colon \mathscr{M} \longrightarrow \mathscr{M}$ is extendable.*

*Proof.* Since $T$ is in particular a lax $\mathscr{C}$-module monad, by Proposition 5.35, the comparison functor $K_T\colon \mathscr{M}_T \longrightarrow \mathscr{M}$ is a strong $\mathscr{C}$-module functor. Similarly, since $T$ is in particular an oplax $\mathscr{C}$-module monad, again by Proposition 5.35, we conclude that also the comparison functor $K^T\colon \mathscr{M} \longrightarrow \mathscr{M}^T$ is a strong $\mathscr{C}$-module functor. Thus, $\iota = K^T K_T$ is a strong $\mathscr{C}$-module functor. □

**Corollary 8.12.** *If $\mathscr{C}$ is left rigid, then any lax $\mathscr{C}$-module monad is extendable.*

*Proof.* By Proposition 2.73, a lax $\mathscr{C}$-module monad is automatically a strong $\mathscr{C}$-module monad. The result follows by Proposition 8.11. □

### 8.1.2 *Semisimple monads*

In this section we want to take a look at *semisimple* monads—those whose category of modules is semisimple. As it turns out, these monads are always extendable. Recall the definition of the Cauchy completion of a category from Section 3.2.1. The following result follows immediately from the fact that for a monad $T\colon \mathscr{C} \longrightarrow \mathscr{C}$, the forgetful functor $U^T\colon \mathscr{C}^T \longrightarrow \mathscr{C}$ creates limits.

**Lemma 8.13.** *The Eilenberg–Moore category of a monad on a Cauchy complete category is itself Cauchy complete.*

**Proposition 8.14.** *Let $T\colon \mathscr{A} \longrightarrow \mathscr{A}$ be a right exact monad on an abelian category $\mathscr{A}$. The category $\mathscr{A}^T$ is semisimple if and only if $\overline{\mathscr{A}_T}$ is. In this semisimple setting, the extension $\bar{\iota}\colon \overline{\mathscr{A}_T} \longrightarrow \mathscr{A}^T$ of the canonical embedding $\iota\colon \mathscr{A}_T \longrightarrow \mathscr{A}^T$ from Equation (2.2.2) to the Cauchy completion is an equivalence of categories.*

*Proof.* Since $\mathscr{A}$ is Cauchy complete, the adjunction $F_T\colon \mathscr{A} \rightleftarrows \mathscr{A}_T : U_T$ extends to an adjunction $\overline{F_T}\colon \mathscr{A} \rightleftarrows \overline{\mathscr{A}_T} : \overline{U_T}$, such that $\overline{U_T F_T} = T$.

It is easy to see that the resulting comparison functor $K\colon \overline{\mathscr{A}_T} \longrightarrow \mathscr{A}^T$ is naturally isomorphic to the extension $\bar{\iota}$ of $\iota$. Further, since $\iota$ is fully faithful, so is $\bar{\iota}$, hence it reflects zero objects and so do $U_T = U^T \circ \iota$ and $\overline{U_T}$.

Assume that $\mathscr{A}^T$ is semisimple. The monomorphisms and epimorphisms in $\mathscr{A}^T$ are split by Proposition 2.85, and hence they are reflected under the fully faithful functor $\bar{\iota}$. Thus, $\overline{\mathscr{A}_T}$ is an abelian category and $\bar{\iota}$ is an exact functor. In particular $\overline{U_T} \cong U_T \circ \bar{\iota}$ is exact and faithful, so it reflects zero objects. By Theorem 2.92, $K$ is an equivalence—in particular, $\overline{\mathscr{A}_T}$ is semisimple.

Assume now that $\overline{\mathscr{A}_T}$ is semisimple. Then $\overline{U_T}$ is exact by Proposition 2.85, and so it again satisfies the assumptions of Theorem 2.92. Thus, $K\colon \overline{\mathscr{A}_T} \xrightarrow{\sim} \mathscr{A}^T$ is an equivalence; in particular, $\mathscr{A}^T$ is semisimple. □

**Definition 8.15.** We say that a right exact monad $T\colon \mathscr{A} \longrightarrow \mathscr{A}$ on an abelian category $\mathscr{A}$ is *semisimple* if it satisfies the equivalent conditions of Proposition 8.14.

Since the opposite of an abelian category is abelian and the opposite of a semisimple category is semisimple, an analogous result to Proposition 8.14 holds for comonads, and we define *semisimple comonads* analogously.

**Proposition 8.16.** *Let $T\colon \mathscr{M} \longrightarrow \mathscr{M}$ be a semisimple lax $\mathscr{C}$-module monad on an abelian $\mathscr{C}$-module category. Then $T$ is extendable.*

*Proof.* By Proposition 8.14, the functor $\bar{\iota}\colon \overline{\mathscr{M}_T} \xrightarrow{\sim} \mathscr{M}^T$ is an equivalence. Further, following Proposition 2.129, $\overline{\mathscr{M}_T}$ has a canonical structure of a $\mathscr{C}$-module category, such that the inclusion $\mathscr{M}_T \hookrightarrow \overline{\mathscr{M}_T}$ is a strong $\mathscr{C}$-module functor. Transporting the $\mathscr{C}$-module structure along $\bar{\iota}$ endows $\iota$ with the structure of a strong $\mathscr{C}$-module functor, proving the claim. □

Semisimplicity also "transfers" to the left adjoint of a monad.

**Lemma 8.17.** *Let $T\colon \mathscr{A} \longrightarrow \mathscr{A}$ be a right exact monad on an abelian category, and $S\colon \mathscr{A} \longrightarrow \mathscr{A}$ a left adjoint to $T$. Then $T$ is semisimple if and only if $S$ is semisimple. In that case, there is an equivalence $\mathscr{A}^T \simeq \mathscr{A}^S$.*

*Proof.* By Proposition 8.14, $T$ is semisimple if and only if so is $\overline{\mathscr{A}_T}$. Moreover, the equivalence $\mathscr{A}_T \simeq \mathscr{A}_S$ of Proposition 2.26 extends to one on the Cauchy completions: $\overline{\mathscr{A}_T} \simeq \overline{\mathscr{A}_S}$. Thus, $\overline{\mathscr{A}_T}$ is semisimple if and only if $\overline{\mathscr{A}_S}$ is so, which is the case if and only if $S$ is semisimple. This establishes the first claim.

The latter claim follows from the equivalences

$$\mathcal{A}^T \simeq \overline{\mathcal{A}_T} \simeq \overline{\mathcal{A}_S} \simeq \mathcal{A}^S. \qquad \square$$

**Proposition 8.18.** *Let $T\colon \mathcal{M} \longrightarrow \mathcal{M}$ be a semisimple lax $\mathcal{C}$-module monad on an abelian $\mathcal{C}$-module category $\mathcal{M}$, and let $S\colon \mathcal{M} \longrightarrow \mathcal{M}$ be a left adjoint to $T$. The equivalence $\mathcal{M}^T \simeq \mathcal{M}^T$ of Lemma 8.17 is a $\mathcal{C}$-module equivalence, where $\mathcal{M}^T$ and $\mathcal{M}^S$ are endowed with the extended $\mathcal{C}$-module structures of $\mathcal{M}_T$ and $\mathcal{M}_S$, respectively.*

*Proof.* Following the proof of Proposition 8.16, the functor $\bar{\iota}\colon \overline{\mathcal{M}_T} \xrightarrow{\sim} \mathcal{M}^T$ is a $\mathcal{C}$-module equivalence, and similarly for $\mathcal{M}^S$ and $\overline{\mathcal{M}_S}$. Further, from Proposition 2.129 we obtain a $\mathcal{C}$-module equivalence $\overline{\mathcal{M}_T} \simeq \overline{\mathcal{M}_S}$ from the $\mathcal{C}$-module equivalence $\mathcal{M}_T \simeq \mathcal{M}_S$ of Proposition 5.37.

The result follows by the following chain of $\mathcal{C}$-module equivalences:

$$\mathcal{M}^T \simeq \overline{\mathcal{M}_T} \simeq \overline{\mathcal{M}_S} \simeq \mathcal{M}^S. \qquad \square$$

### 8.1.3 *Linton coequalisers via multiactegories*

There exists a quite general condition for extendability, in which one can even explicitly write down the resulting $\mathcal{C}$-module structure on the Eilenberg–Moore category. Our aim now is to establish Theorems 8.25 and 8.26. Thus, for the rest of this section, we make the following assumptions.

**Hypothesis 8.19.**

- Let $\mathcal{C}$ be an abelian monoidal and $\mathcal{M}$ an abelian $\mathcal{C}$-module category.
- The action functor $\triangleright\colon \mathcal{C} \otimes_{\Bbbk} \mathcal{M} \longrightarrow \mathcal{M}$ is right exact in both variables.
- Let $T\colon \mathcal{M} \longrightarrow \mathcal{M}$ be a right exact lax $\mathcal{C}$-module monad.

Note that, due to for example Proposition 8.2, many of these assumptions are already necessary for the theory we would like to develop in this chapter.

In our presentation, we closely follow [AHLF18], where an analogous result for lax monoidal monads—rather than lax module monads—can be found. As such, we shall often only describe the modifications necessary to make these results work in our case.

**Definition 8.20.** The *Linton coequaliser* $x \blacktriangleright m$ of $x \in \mathcal{C}$ and $(m, \nabla_m) \in \mathcal{M}^T$ is:

$$x \blacktriangleright m := \operatorname{coeq}\left( T(x \triangleright Tm) \underset{\mu_{x\triangleright m} \circ T\mathsf{T}_{\mathsf{a};x,m}}{\overset{T(x\triangleright\nabla_m)}{\rightrightarrows}} T(x \triangleright m) \right).$$

As $\mathcal{M}^T$ is abelian due to Proposition 8.2, by functoriality of colimits we obtain a functor $- \blacktriangleright = : \mathcal{C} \otimes_{\Bbbk} \mathcal{M}^T \longrightarrow \mathcal{M}^T$.

**Lemma 8.21.** *For $x \in \mathcal{C}$, the functor $x \blacktriangleright - : \mathcal{M}^T \longrightarrow \mathcal{M}^T$ is right exact.*

*Proof.* This immediately follows from the right exactness of $T(x \triangleright T(-))$ and $T(x \triangleright -)$, as well as the definition of the Linton coequaliser. □

**Lemma 8.22.** *There is a natural isomorphism $1 \blacktriangleright - \cong \mathrm{Id}_{\mathcal{M}^T}$.*

*Proof.* Recall that for all $T$-algebras $m \in \mathcal{M}^T$,

$$T^2 m \underset{T\nabla_m}{\overset{\mu_m}{\rightrightarrows}} Tm \xrightarrow{\nabla_m} m$$

is a coequaliser in $\mathcal{M}^T$, created by the underlying split coequaliser in $\mathcal{M}$. Now,

$$\begin{array}{ccc} T(1 \triangleright Tm) & \underset{\mu_{1 \triangleright m} \circ T T_{\mathsf{a};1,m}}{\overset{T(1 \triangleright \nabla_m)}{\rightrightarrows}} & T(1 \triangleright m) \\ {\scriptstyle T\lambda}\downarrow{\scriptstyle \cong} & & {\scriptstyle \cong}\downarrow{\scriptstyle T\lambda} \\ T^2 m & \underset{\mu_m}{\overset{T\nabla_m}{\rightrightarrows}} & Tm \end{array}$$

defines an isomorphism in the category of parallel pairs of morphisms in $\mathcal{M}^T$—the upper square commutes due to the naturality of the left unitor $\lambda$ of $\mathcal{M}$, and the lower by coherence for $- \triangleright =$. This isomorphism induces an isomorphism on the level of coequalisers, $m \cong 1 \blacktriangleright m$. Further, this construction of morphisms of parallel pairs is clearly natural in $m$, establishing the naturality of the isomorphism of coequalisers. □

**Lemma 8.23.** *For any $m \in \mathcal{M}$, the coequaliser of*

$$T(x \triangleright T^2 m) \underset{\mu_{x \triangleright Tm} \circ T T_{\mathsf{a};x,Tm}}{\overset{T(x \triangleright \mu_m)}{\rightrightarrows}} T(x \triangleright Tm)$$

*in $\mathcal{M}^T$ is natural in $m$, split, and isomorphic to $T(x \triangleright m)$. In other words the Linton coequaliser on a free $T$-module is given by $x \blacktriangleright Tm \cong T(x \triangleright m)$.*

*Proof.* The splitting data is given by

$$T(x \triangleright T^2 m) \underset{\mu_{x \triangleright Tm} \circ T T_{\mathsf{a};x,Tm}}{\overset{T(x \triangleright \mu_m)}{\rightrightarrows}} T(x \triangleright Tm) \xrightarrow{\mu_{x \triangleright m} \circ T T_{\mathsf{a};x,m}} T(x \triangleright m),$$

with sections $T(x \triangleright T\eta_m)\colon T(x \triangleright Tm) \to T(x \triangleright T^2 m)$ and $T(x \triangleright \eta_m)\colon T(x \triangleright m) \to T(x \triangleright Tm)$. □

**Lemma 8.24.** *For $x, y \in \mathscr{C}$, there is a natural isomorphism $(x \otimes y) \blacktriangleright - \cong x \blacktriangleright (y \blacktriangleright -)$.*

*Proof.* We calculate

$$\begin{aligned} x \blacktriangleright (y \blacktriangleright m) &= x \blacktriangleright \mathrm{coeq}(T(y \triangleright \nabla_m),\ \mu_{x \triangleright m} \circ T T_{\mathsf{a};y,m}) \\ &\cong \mathrm{coeq}(x \blacktriangleright T(y \triangleright \nabla_m),\ x \blacktriangleright \mu_{x \triangleright m} \circ T T_{\mathsf{a};y,m}) \\ &\cong \mathrm{coeq}(T(x \triangleright (y \triangleright \nabla_m)),\ \mu_{y \triangleright x \triangleright m} \circ T(T_{\mathsf{a};x,y \triangleright m} \circ y \triangleright T_{\mathsf{a};x,m})) \\ &\cong \mathrm{coeq}(T((x \otimes y) \triangleright \nabla_m),\ \mu_{(x \otimes y) \triangleright Tm} \circ T T_{\mathsf{a};x \otimes y, Tm}) \\ &= (x \otimes y) \blacktriangleright m, \end{aligned}$$

where the first isomorphism follows from Lemma 8.21, the second one is due to Lemma 8.23, and the third one follows by coherence of $- \triangleright =$. □

In view of Lemmas 8.22 and 8.24, all that remains in order to establish that $- \blacktriangleright =$ defines a $\mathscr{C}$-module category structure on $\mathscr{M}^T$ is the verification of coherence axioms. While this can be accomplished by diagram chasing—see [Sea13, Theorem 2.6.4]—we choose to give a more formal argument, mimicking the multicategorical approach used by [AHLF18] in the analogous case of monoidal monads. More specifically, we shall prove the following two results in the remainder of this section.

**Theorem 8.25.** *The functor $- \blacktriangleright = : \mathscr{C} \otimes_{\Bbbk} \mathscr{M}^T \longrightarrow \mathscr{M}^T$ defines a left $\mathscr{C}$-module category structure on $\mathscr{M}^T$.*

Assuming Theorem 8.25, we have $x \blacktriangleright Tm = T(x \triangleright m)$, and so the image of the inclusion $\iota \colon \mathscr{M}_T \longrightarrow \mathscr{M}^T$ of Equation (2.2.2) is a $\mathscr{C}$-module subcategory, where $\mathscr{M}_T$ is endowed with the action given in the proof of Corollary 5.33.

**Theorem 8.26.** *The unique embedding $\iota \colon (\mathscr{M}_T, \triangleright_{\mathscr{M}_T}) \longrightarrow (\mathscr{M}^T, \blacktriangleright)$ can be equipped with the structure of a strong $\mathscr{C}$-module functor.*

Let us now recall the theory of multicategories, see [Her00; Lei04].

**Definition 8.27.** A (locally small) *multicategory* $\mathbf{C}$ consists of

- a class $\mathrm{Ob}\,\mathbf{C}$ of objects of $\mathbf{C}$, where we write $x \in \mathbf{C}$ for $x \in \mathrm{Ob}\,\mathbf{C}$;
- for any finite sequence $(x_1, \ldots, x_n, y)$ of objects in $\mathbf{C}$, a set of *multimorphisms* $\mathbf{C}(x_1, \ldots, x_n; y)$ from $(x_1, \ldots, x_n)$ to $y$;
- for every $x \in \mathbf{C}$, an identity multimorphism $\mathrm{id}_x \in \mathbf{C}(x; x)$; and

- for any $y \in \mathrm{Ob}(\mathbf{C})$, $(x_i)_{i=1}^n \in \mathrm{Ob}(\mathbf{C})^n$, as well as

$$\big((u_1^1, \dots, u_1^{k_1}), \dots, (u_n^1, \dots, u_n^{k_n})\big) \in \mathrm{Ob}(\mathbf{C})^{k_1+\cdots+k_n},$$

  a composition operation

$$\begin{aligned}&\mathbf{C}(x_1, \dots, x_n; y) \times \mathbf{C}(u_1^1, \dots, u_1^{k_1}; x_1) \times \cdots \times \mathbf{C}(u_n^1, \dots, u_n^{k_n}; x_n)\\ &\quad \longrightarrow \mathbf{C}(u_1^1, \dots, u_n^{k_n}; y),\end{aligned}$$

subject to natural associativity and unitality conditions; see for example [Lei04, Section 2.1] or [Her00, Definition 2.1].

Since we want our multicategories and multiactegories to correspond to $\Bbbk$-linear monoidal and module categories, we should replace sets of multimorphisms with vector spaces, and the Cartesian products with tensor products. However, since our aim is to show Theorem 8.25, which can be verified on the level of the underlying (ordinary) categories, this is not essential. The following definition is a non-skew special case of [AM24, Definition 6.9].

**Definition 8.28.** Let $\mathbf{C}$ be a (locally small) multicategory. A (locally small) *left multiactegory* $\mathbf{M}$ over $\mathbf{C}$ has as data:

- a class $\mathrm{Ob}\,\mathbf{M}$ of objects of $\mathbf{M}$, where we write $m \in \mathbf{M}$ for $m \in \mathrm{Ob}\mathbf{M}$;
- for any finite (possibly empty) sequence $(x_1, \dots, x_n)$ of objects in $\mathbf{C}$ and a pair of objects $(m, m')$ of $\mathbf{M}$, a set $\mathbf{M}(x_1, \dots, x_n; m; m')$ of *multimorphisms* from $(x_1, \dots, x_n; m)$ to $m'$;
- for every $m \in \mathbf{M}$, an identity multimorphism $\mathrm{id}_m \in \mathbf{M}(m; m)$; and
- for any $x \in \mathbf{C}$, $m, m', \ell \in \mathbf{M}$, $(x_i)_{i=1}^n \in \mathrm{Ob}(\mathbf{C})^n$, and

$$\big((u_1^1, \dots, u_1^{k_1}), \dots, (u_n^1, \dots, u_n^{k_n})\big) \in \mathrm{Ob}(\mathbf{C})^{k_1+\cdots+k_n},$$

  a composition operation

$$\begin{aligned}&\mathbf{M}(x_1, \dots, x_n; m; m') \times \mathbf{C}(u_1^1, \dots, u_1^{k_1}; x_1) \times \cdots \times \mathbf{C}(u_n^1, \dots, u_n^{k_n}; x_n)\\ &\quad \times \mathbf{M}(u_{n+1}^1, \dots, u_{n+1}^{k_{n+1}}; \ell; m) \longrightarrow \mathbf{M}(u_1^1, \dots, u_{n+1}^{k_{n+1}}; \ell; m').\end{aligned} \tag{8.1.4}$$

This data has to satisfy the (non-marked) associativity and unitality axioms of [AM24, Definitions 4.4 and 6.9].

**Notation 8.29.** From now on, we will often abbreviate sequences of objects using the vector notation similar to [Her00]; e.g., writing $\vec{x}$ for $(x_1, \ldots, x_n)$ and $(\vec{x}_1, \vec{x}_2)$ for $(x_1^1, \ldots, x_1^{n_1}, x_2^1, \ldots, x_2^{n_2})$. We will use this notation for clarity and brevity, when keeping track of the length of the tuples is not required. We also fix the notation for the remainder of this section, letting $\mathbf{C}$ be a multicategory and letting $\mathbf{M}$ be a $\mathbf{C}$-multiactegory.

Recall from [Her00, Chapter 8] that a *pre-universal arrow* for a tuple $\vec{x}$ of objects of $\mathbf{C}$ consists of an object $\otimes(\vec{x})$ of $\mathbf{C}$ together with a multimorphism $\pi \in \mathbf{C}(\vec{x}; \otimes(\vec{x}))$, inducing isomorphisms $\mathbf{C}(\otimes(\vec{x}); y) \xrightarrow{\sim} \mathbf{C}(\vec{x}; y)$. If $\pi$ further induces isomorphisms

$$\mathbf{C}(\vec{z}, \otimes(\vec{x}), \vec{z}'; y) \xrightarrow{\sim} \mathbf{C}(\vec{z}, \vec{x}, \vec{z}'; y),$$

we say that $\pi$ is *universal*. If any sequence of objects of $\mathbf{C}$ is the domain of a universal arrow, we say that $\mathbf{C}$ is *representable*. A *representation* of $\mathbf{C}$ is a choice of a universal arrow from every sequence of objects in $\mathbf{C}$.

**Definition 8.30.** Let $\mathbf{M}$ be a $\mathbf{C}$-multiactegory. A *pre-universal arrow* for a tuple $\vec{x}$ of objects of $\mathbf{C}$ and $m \in \mathbf{M}$ consists of an object $\triangleright(\vec{x}; m)$, together with a multimorphism $\pi \in \mathbf{M}(\vec{x}; m; \triangleright(\vec{x}; m))$, inducing isomorphisms

$$\mathbf{M}(\triangleright(\vec{x}; m); m') \xrightarrow{\sim} \mathbf{M}(\vec{x}; m; m').$$

If $\pi$ further induces isomorphisms

$$\mathbf{M}(\vec{y}, \triangleright(\vec{x}; m); m') \xrightarrow{\sim} \mathbf{M}(\vec{y}, \vec{x}; m; m'),$$

then $\pi$ is said to be *universal*.

**Definition 8.31.** We say that a multimorphism $\rho \in \mathbf{C}(\vec{x}; y)$ is *pre-universal in* $\mathbf{M}$ if it induces isomorphisms $\mathbf{M}(y; m; n) \xrightarrow{\sim} \mathbf{M}(\vec{x}; m; n)$. It is *universal in* $\mathbf{M}$, if it also induces isomorphisms $\mathbf{M}(\vec{z}, \vec{x}, \vec{w}; m; n) \xrightarrow{\sim} \mathbf{M}(\vec{z}, y, \vec{w}; m; n)$.

**Definition 8.32.** For $\mathbf{C}$ representable and a representation $\mathsf{R}$ of $\mathbf{C}$, we say that $\mathbf{M}$ is a *representable with respect to* $\mathsf{R}$, if for every sequence $\vec{x}$ of objects of $\mathbf{C}$ and every object $m \in \mathbf{M}$, there is a universal arrow $\pi_{(\vec{x};m)} \in \mathbf{M}(\vec{x}; m; \triangleright(\vec{x}, m))$, and all morphisms of $\mathsf{R}$ are universal in $\mathbf{M}$.

Recall from [Her00, Proposition 8.5] that universal arrows in a multicategory $\mathbf{C}$ are closed under composition. Further, if every sequence of objects in $\mathbf{C}$ is the domain of a pre-universal arrow and pre-universal arrows are closed under composition, then a pre-universal arrow in $\mathbf{C}$ is universal. A similar result holds for multiactegories.

**Lemma 8.33.**

1. *Composition of universal arrows of* **C** *in* **M** *and universal arrows of* **M** *yields universal arrows in* **M**.
2. *If composition of pre-universal arrows of* **C** *in* **M** *and pre-universal arrows of* **M** *yields pre-universal arrows in* **M**, *then a pre-universal morphism in* **M** *is universal.*
3. *If pre-universal arrows in* **C** *act pre-universally in* **M**, *and pre-universal arrows of* **C** *and* **M** *are closed under composition in* **M**, *then* **M** *is representable.*

*Proof.* For the first part, let $(\pi_0, \rho_1, \ldots, \rho_n, \pi_{n+1})$ be a sequence composable in **M**, with $\pi_0, \pi_{n+1}$ in **M** and $\rho_1, \ldots, \rho_n \in \mathbf{C}$, as indicated in Equation (8.1.4), all of whose elements are universal in **M**. The composite $\pi_0 \circ (\rho_1, \ldots, \rho_n, \pi_{n+1})$ gives rise to the sequence of isomorphisms

$$\begin{aligned}&\mathbf{M}(\vec{x}, \triangleright(\otimes(\vec{y}_1), \ldots, \otimes(\vec{y}_n); \triangleright(\vec{y}_{n+1}, m)), -)\\ &\quad\cong \mathbf{M}(\vec{x}, \otimes(\vec{y}_1), \ldots, \otimes(\vec{y}_n); \triangleright(\vec{y}_{n+1}, m); -) \cong \mathbf{M}(\vec{x}, \ldots, \vec{y}_1, \ldots, \vec{y}_{n+1}; m; -),\end{aligned}$$

proving the universality of the composite.

For the second part, consider the following chain of maps:

$$\begin{aligned}\mathbf{M}(\vec{x}, \vec{y}, \vec{z}; m; -) \to \mathbf{M}(\vec{x}, \otimes(\vec{y}), \vec{z}; m; -) \to \mathbf{M}(\otimes(\vec{x}), \otimes(\vec{y}), \otimes(\vec{z}); m; -)\\ \to \mathbf{M}(\triangleright(\otimes(\vec{x}), \otimes(\vec{y}), \otimes(\vec{z}); m); -),\end{aligned}$$

where the first morphism is induced by the pre-universal arrow $\rho_{\vec{x}}$ of **C**, the second one by the pre-universal arrows $\rho_{\vec{y}}$ and $\rho_{\vec{z}}$ of **C**, and the third is induced by the universal arrow $\pi_{(\vec{y}, \vec{x}, \vec{z}; m)}$. The composite of all three maps is an isomorphism—being induced by a composite of pre-universal arrows—and hence a pre-universal arrow. The same holds for the composite of the latter two maps. These two observations establish the invertibility of the first map.

The third part is established by applying the proof of the second part above to the case of pre-universal arrows of **C**. □

Assume **C** to be representable, let R be a representation of **C**, and suppose that **M** is representable with respect to R. Given a sequence $x_1, \ldots, x_n$ of objects in **C** and an object $m$ in **M**, the set of ways in which we may compose the universal arrows in the given representations to a universal arrow with domain $(x_1, \ldots, x_n; m)$ is canonically in bijection with the set of parenthesisations of the word $x_1 \cdots x_n m$ into subwords of length at least two. Let $\diamond$ and $\diamond'$ be two such parenthesisations. We denote the codomains of the corresponding universal arrows by $\diamond(x_1,\dots,x_n;m)$ and $\diamond'(x_1,\dots,x_n;m)$, and the universal arrows themselves by $\pi_{\diamond(x_1,\dots,x_n;m)}$ and $\pi_{\diamond'(x_1,\dots,x_n;m)}$. The universality of $\pi_{\diamond(x_1,\dots,x_n;m)}$ entails the existence of a unique arrow

$$\alpha^{\diamond,\diamond'}_{\vec{x},m}\colon \diamond(x_1,\dots,x_n;m) \longrightarrow \diamond'(x_1,\dots,x_n;m),$$

satisfying $\alpha^{\diamond,\diamond'}_{\vec{x},m} \circ \pi_{\diamond(\vec{x},m)} = \pi_{\diamond'(\vec{x},m)}$, which is invertible with inverse $\alpha^{\diamond',\diamond}_{\vec{x},m}$.

Further, $\alpha^{\diamond,\diamond''}_{\vec{x},m} = \alpha^{\diamond',\diamond''}_{\vec{x},m} \circ \alpha^{\diamond,\diamond'}_{\vec{x},m}$. Thus, for any sequence $x_1,\dots,x_n,m$ we obtain an indiscrete category[19] whose objects are parenthesisations of $x_1\cdots x_n m$ and $\mathrm{Hom}(\diamond,\diamond') = \{\alpha^{\diamond,\diamond'}_{\vec{x},m}\}$. The arrow $\alpha^{\diamond,\diamond'}_{\vec{x},m}$ is natural in $m$ and $x_i$ for all $i$.

[19] A category is called *indiscrete* if there exists a unique morphism from every object to every other object.

**Remark 8.34.** Recall that given a representable multicategory $\mathbf{C}$ with a fixed representation, there is a monoidal category $\mathscr{C}$ defined by $\mathrm{Ob}\,\mathscr{C} := \mathrm{Ob}\,\mathbf{C}$, $\mathscr{C}(x,y) := \mathbf{C}(x;y)$, and $x \otimes y := \otimes(x,y)$, the codomain of the universal arrow $\pi_{\otimes(x,y)}$ from $(x,y)$ in the representation of $\mathbf{C}$. The associator of $\mathscr{C}$ is defined by the $\alpha^{(wx)y,w(xy)}_{w,x,y}$ described above; one obtains similar arrows for unitality.

The indiscreteness of the category of morphisms $\alpha^{\diamond,\diamond'}$ associated to the parenthesisations of words in $\mathbf{C}$ implies the commutativity of the pentagon and triangle diagrams, establishing the coherence and well-definedness of $\mathscr{C}$.

The indiscrete category described above establishes the well-definedness of the $\mathscr{C}$-module category $\mathscr{M}$ of the following definition.

**Definition 8.35.** Let $\mathbf{M}$ be a representable multiactegory over a representable multicategory $\mathbf{C}$, with fixed representations for both. Let $\mathscr{C}$ be the monoidal category associated to the representation of $\mathbf{C}$. The *$\mathscr{C}$-module category $\mathscr{M}$ associated to the fixed representation of* $\mathbf{M}$ is defined by setting $\mathrm{Ob}\,\mathbf{M} := \mathrm{Ob}\,\mathscr{M}$, $\mathscr{M}(m,n) := \mathbf{M}(m;n)$, and $x \triangleright m := \triangleright(x,m)$, the codomain of the universal arrow $\pi_{\triangleright(x,m)}$. The associator $(x\otimes y)\triangleright m \xrightarrow{\sim} x \triangleright (y \triangleright m)$ is given by the arrow $\alpha^{(xy)m,x(ym)}_{x,y,m}$, and similarly for the unitors.

Having established sufficient analogues of [Her00; Lei04], we continue with analogues of [AHLF18, Section 3]. For $x_i \in \mathscr{C}$ and $m \in \mathscr{M}$, write

$$\varphi\colon x_1 \triangleright \cdots \triangleright x_n \triangleright Tm \longrightarrow T(x_1 \triangleright \cdots \triangleright x_n \triangleright m)$$

for the morphism obtained by repeated application of the lax structure of $T$, where $x_1 \triangleright \cdots \triangleright x_n \triangleright m$ is to be read as $x_1 \triangleright (\cdots \triangleright (x_n \triangleright m))$.

We will also consider the coequaliser diagram

$$T(x_1 \triangleright \cdots \triangleright x_n \triangleright Tm) \underset{T\varphi}{\longrightarrow} T^2(x_1 \triangleright \cdots \triangleright x_n \triangleright m) \underset{\mu}{\longrightarrow} T(x_1 \triangleright \cdots \triangleright x_n \triangleright m) \qquad \text{with } T(x_1 \triangleright \cdots \triangleright x_n \triangleright \nabla_m)\colon T(x_1 \triangleright \cdots \triangleright x_n \triangleright Tm) \longrightarrow T(x_1 \triangleright \cdots \triangleright x_n \triangleright m) \tag{8.1.5}$$

which is reflexive, with $T(x_1 \triangleright \cdots \triangleright x_n \triangleright m) \xrightarrow{T(x_1 \triangleright \cdots \triangleright x_n \triangleright \eta_m)} T(x_1 \triangleright \cdots \triangleright x_n \triangleright Tm)$ as a common section.

**Definition 8.36.** A morphism $f\colon x_1 \triangleright \cdots \triangleright x_n \triangleright m \longrightarrow m'$ is called *$n$-multilinear* if the following diagram commutes:

$$\begin{array}{ccc}
x_1 \triangleright \cdots \triangleright x_n \triangleright Tm & \xrightarrow{\varphi} & T(x_1 \triangleright \cdots \triangleright x_n \triangleright m) \\
 & & \downarrow{\scriptstyle Tf} \\
\downarrow{\scriptstyle x_1 \triangleright \cdots \triangleright x_n \triangleright \nabla_m} & & Tm' \\
 & & \downarrow{\scriptstyle \nabla_{m'}} \\
x_1 \triangleright \cdots \triangleright x_n \triangleright m & \xrightarrow[f]{} & m'
\end{array}$$

**Remark 8.37.** For a monoidal category $\mathscr{C}$, there exists a representable multicategory $\mathbf{C}$, with $\mathbf{C}(x_1, \ldots, x_n; y) := \mathscr{C}(x_1 \otimes (\cdots (x_{n-1} \otimes x_n) \ldots), y)$. Then, a multimorphism $\pi \in \mathbf{C}(x_1, \ldots, x_n; y)$ is universal if and only if it is an isomorphism $\pi\colon x_1 \otimes \cdots \otimes x_n \xrightarrow{\sim} y$.

Similarly, given a $\mathscr{C}$-module category, one obtains a representable multiactegory $\mathbf{M}$ over $\mathbf{C}$ satisfying $\mathbf{M}(x_1, \ldots, x_n; m; n) = \mathscr{M}(x_1 \triangleright \cdots \triangleright x_n \triangleright m, m')$. A universal arrow $\pi \in \mathbf{M}(x_1, \ldots, x_n; m; m')$ is an isomorphism $\pi\colon x_1 \triangleright \cdots \triangleright m \xrightarrow{\sim} m'$.

**Lemma 8.38.** *For a $\mathscr{C}$-module category $\mathscr{M}$, the collection of multilinear morphisms of $\mathbf{M}$ is stable under multicomposition with morphisms of $\mathbf{C}$, forming a $\mathscr{C}$-multiactegory that we denote by $\mathbf{M}^{\varphi}$.*

*Proof.* For $i = 1, \ldots, n$, consider sequences $\vec{z_i}$; morphisms $\otimes(z_i) \xrightarrow{g_i} x_i$, where $\otimes(z_i)$ again is given the rightmost parenthesisation; and a multilinear morphism $f\colon x_1 \triangleright \ldots \triangleright x_n \triangleright m \longrightarrow m'$ of $\mathbf{M}$. The claim follows immediately from the

commutativity of

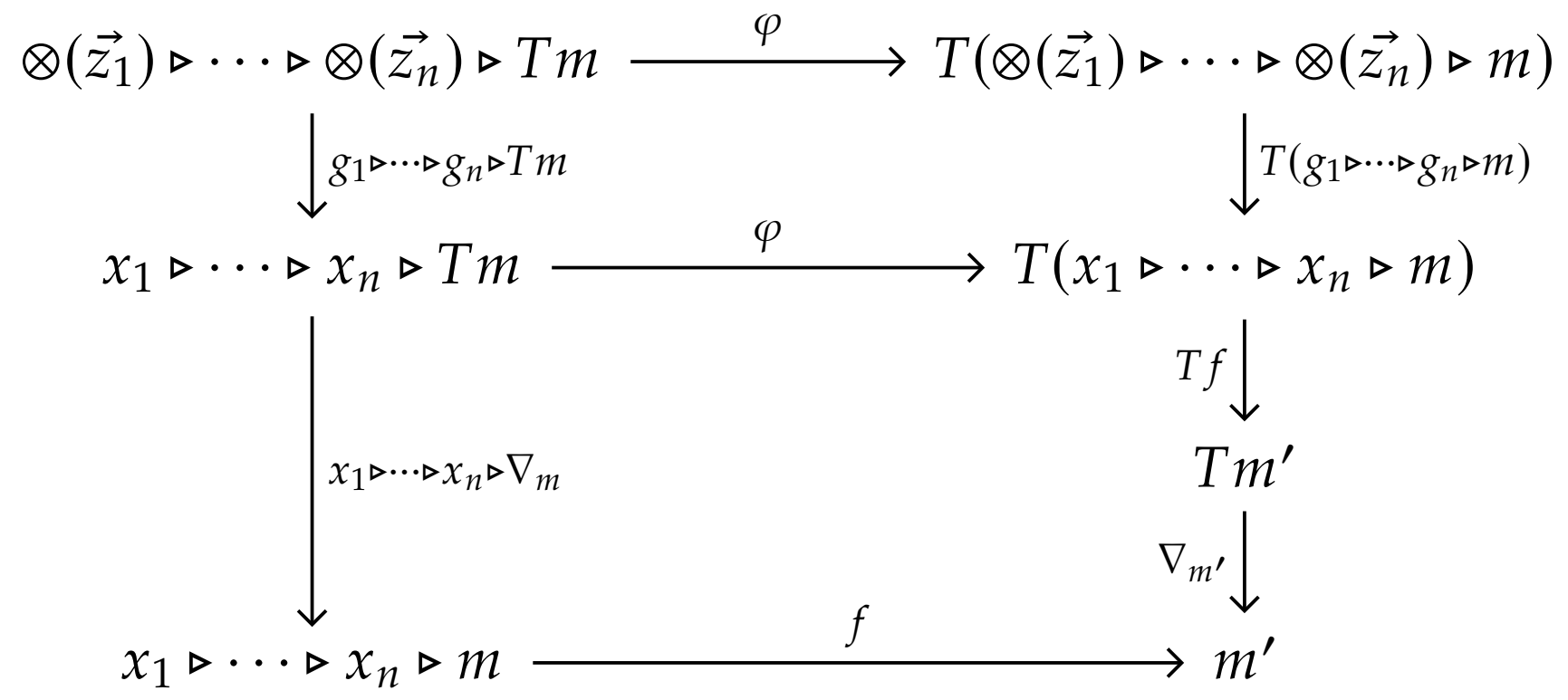

where the top face commutes by naturality of the action of $\mathscr{C}$ in $\mathscr{M}$. □

Note that $\mathbf{M}^{\varphi}$ being representable would prove Theorem 8.25, since then the $\mathscr{C}$-module category can be equipped with coherent associators; see Remark 8.34 and Definition 8.35.

**Lemma 8.39.** *A map $f$ is $n$-multilinear if and only if the map*

$$\underline{f}\colon T(x_1 \triangleright \ldots \triangleright x_n \triangleright m) \xrightarrow{Tf} Tm' \xrightarrow{\nabla_{m'}} m'$$

*coequalises Diagram (8.1.5), and $f$ is pre-universal in $\mathbf{M}^{\varphi}$ if and only if $\underline{f}$ is the (universal) coequaliser of ibid.*

*Proof.* Consider the following diagram:

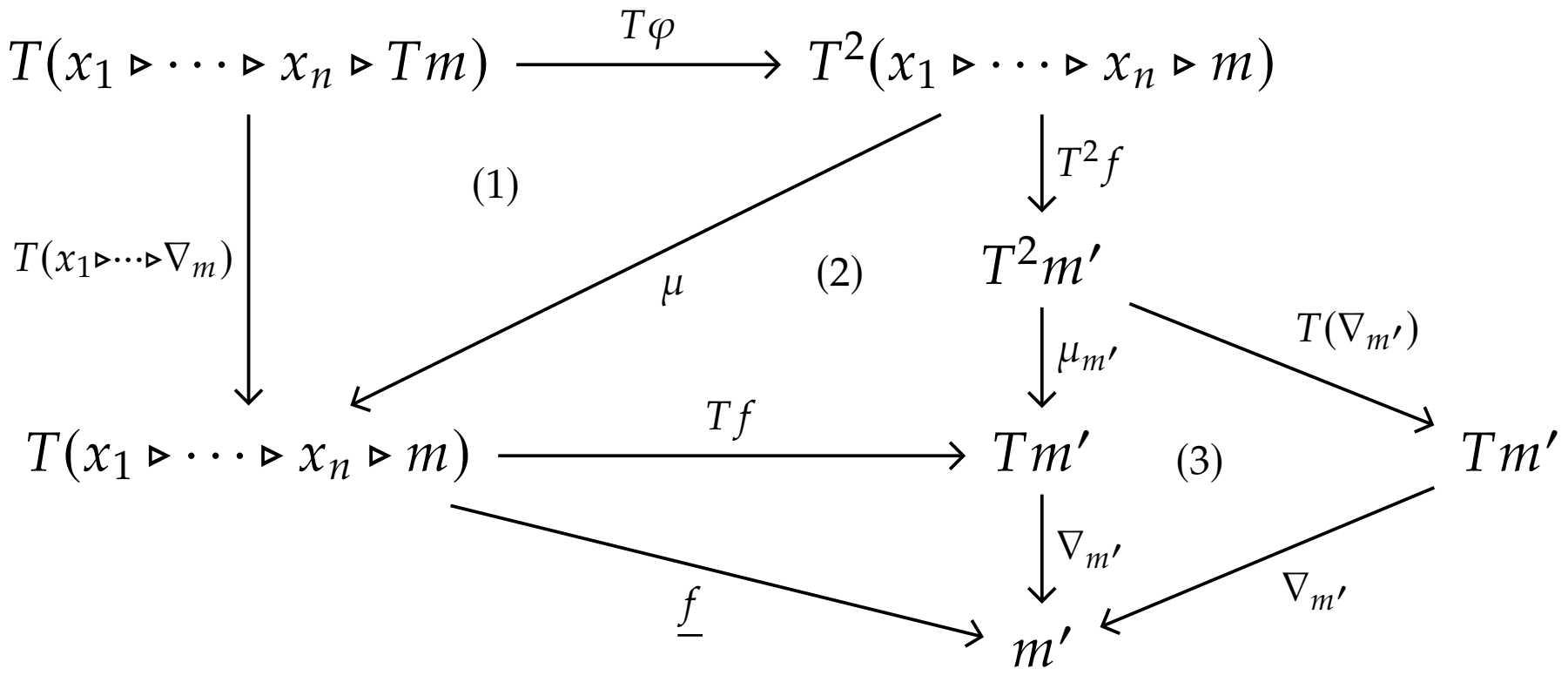

Coequalising the pair (8.1.5) is precisely coequalising face (1) in the above diagram. Thus, since face (2) commutes, it is equivalent to $\nabla_{m'}$ coequalising the square formed jointly by faces (1) and (2).

Observe that, using commutativity of face (3), this square postcomposed by $\nabla_{m'}$ gives precisely the pair of morphisms in $\mathcal{M}^T(T(x_1 \triangleright \cdots \triangleright Tm), m')$ corresponding to the pair of morphisms in $\mathcal{M}(x_1 \triangleright \cdots \triangleright T(m), m')$, defining multilinearity of $f$, under $\mathcal{M}^T(T(-), =) \cong \mathcal{M}(-, =)$. Thus, both conditions in the first statement are equivalent to the commutativity of the diagram.

The second statement is shown analogously to [AHLF18, Lemma 3.5] □

Finally, in order to establish Theorems 8.25 and 8.26, we need an analogue of [AHLF18, Proposition 3.13].

**Proposition 8.40.** *The multicategory* $\mathbf{M}^\varphi$ *is representable.*

*Proof.* Let $(x_1, \ldots, x_n, m')$ be a sequence of objects, with $x_i \in \mathcal{C}$ and $m' \in \mathcal{M}$. The coequalizer of Diagram (8.1.5) defines a pre-universal map in $\mathbf{M}^\varphi$ with domain $(x_1, \ldots, x_n; m')$, by Lemma 8.39. It remains to show that a composition of pre-universal maps is pre-universal. Thus, consider sequences $\vec{z}_1, \ldots, \vec{z}_n$ and $\vec{y}$, universal multimorphisms $h_i \colon \vec{z}_i \longrightarrow x_i$ of $\mathbf{C}$, and universal multimorphisms $f \colon (\vec{y}, m) \longrightarrow m'$ and $g \colon (\vec{x}, m') \longrightarrow m''$.

Behold the following diagram:

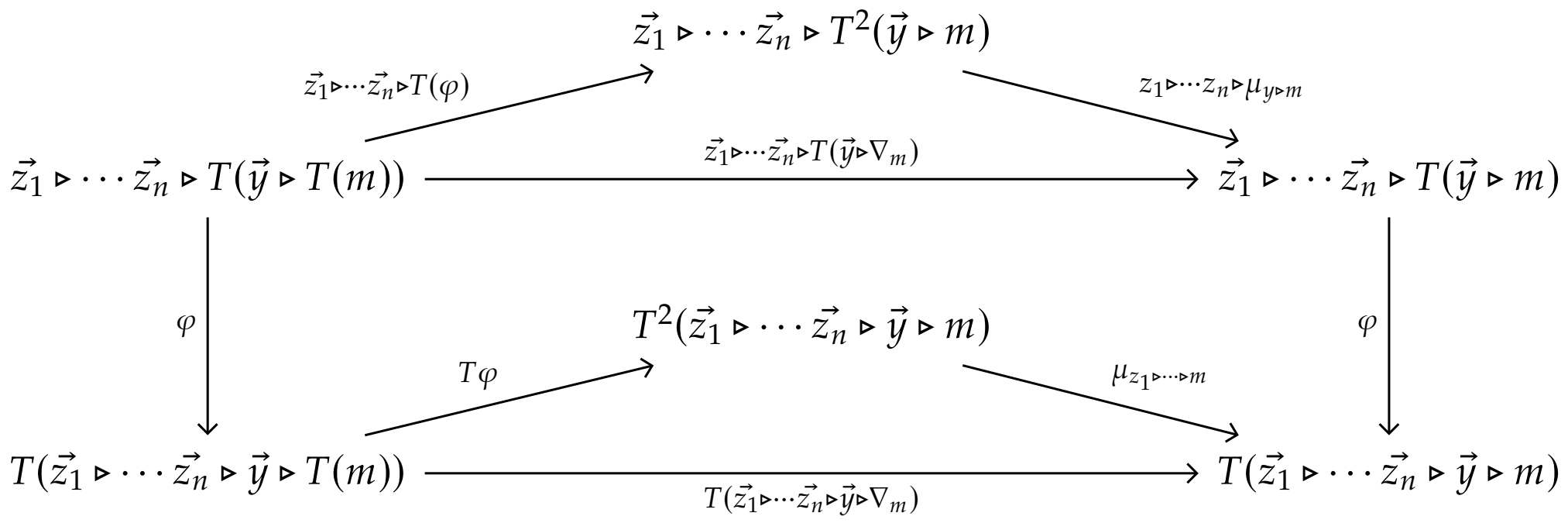

Similarly to the proof of [AHLF18, Proposition 3.13], this is a morphism from the top triangle-shaped diagram to the bottom one, where the top is the action of $(\vec{z}_1, \ldots, \vec{z}_n)$ on the Linton pair for $(\vec{y}, m)$, while the bottom is the Linton pair for $(\vec{z}_1, \ldots, \vec{z}_n, \vec{y}; m)$.

Now, let $X \in \mathbf{M}^\varphi(\vec{z}_1, \ldots, \vec{z}_n, \vec{y}; m; \ell)$. Consider the diagram

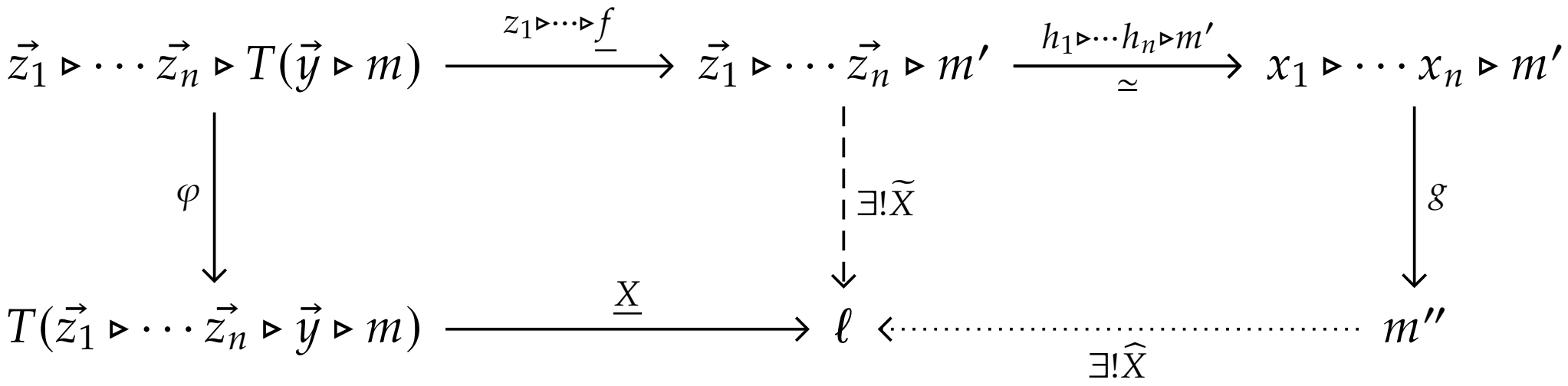

where $\widetilde{X}$ exists since, by Lemma 8.39, $\underline{X}$ coequalises the bottom Linton pair in the penultimate diagram above, and $h_1 \triangleright \cdots \triangleright \underline{f}$ is the coequaliser of the top Linton pair in the same diagram.

Similarly to [AHLF18, Lemma 3.11] one shows that $\widetilde{X}$ is multilinear, and thus so is $h_1^{-1} \triangleright \cdots \triangleright h_n^{-1} \triangleright m'$, since $h_i$ is invertible for all $i$, being a pre-universal, and hence universal, arrow of $\mathbf{C}$. This proves the existence of $\widehat{X}$.

Thus, the multilinear morphism $X$ factors through $g \circ (h_1, \ldots, h_n, f)$, and this factorisation is unique by uniqueness of $\widetilde{X}$ and of $\widehat{X}$. This establishes the pre-universality of the composition $g \circ (h_1, \ldots, h_n, f)$. □

## 8.2 Internal projective and injective objects

**Definition 8.41.** Let $\mathscr{C}$ be an abelian monoidal category and let $\mathscr{M}$ be an abelian $\mathscr{C}$-module category. An object $M \in \mathscr{M}$ is said to be *$\mathscr{C}$-projective* if, for any projective object $P \in \mathscr{C}$, the object $P \triangleright M \in \mathscr{M}$ is projective. Similarly, $M$ is *$\mathscr{C}$-injective* if, for any injective object $I \in \mathscr{C}$, the object $I \triangleright M \in \mathscr{M}$ is injective.

The next result connects the formal definition of a $\mathscr{C}$-projective with the ordinary notion of projective object in an abelian category.

**Proposition 8.42.** *Assume that $\mathscr{C}$ and $\mathscr{M}$ have enough projectives. Let $M$ be a closed object of $\mathscr{M}$. Then $\lfloor M, - \rfloor$ is right exact if and only if $M$ is $\mathscr{C}$-projective.*

*Proof.* If $\lfloor M, - \rfloor$ is right exact and $P$ is a projective object of $\mathscr{C}$. Then $P \triangleright M$ is projective because $\mathscr{M}(P \triangleright M, -) \cong \mathscr{C}(P, \lfloor M, - \rfloor)$, which is a composite of right exact functors, hence itself a right exact functor.

Assume that $M$ is $\mathscr{C}$-projective. Let $\operatorname{colim}_j N_j$ be a finite colimit in $\mathscr{M}$, and $X \in \mathscr{C}$. Since $\mathscr{C}$ has enough projectives, one has $X \cong \operatorname{colim}_i Q_i$, realising $X$ as the colimit of a finite diagram of projectives. We then have

$$\begin{aligned}
\mathscr{C}(X, \lfloor M, \operatorname*{colim}_j N_j \rfloor) &\cong \mathscr{C}(\operatorname*{colim}_i Q_i, \lfloor M, \operatorname*{colim}_j N_j \rfloor) \\
&\cong \lim_i \mathscr{C}(Q_i, \lfloor M, \operatorname*{colim}_j N_j \rfloor) \cong \lim_i \mathscr{M}(Q_i \triangleright M, \operatorname*{colim}_j N_j) \\
&\cong \lim_i \operatorname*{colim}_j \mathscr{M}(Q_i \triangleright M, N_j) \cong \lim_i \operatorname*{colim}_j \mathscr{C}(Q_i, \lfloor M, N_j \rfloor) \\
&\cong \lim_i \mathscr{C}(Q_i, \operatorname*{colim}_j \lfloor M, N_j \rfloor) \cong \mathscr{C}(\operatorname*{colim}_i Q_i, \operatorname*{colim}_j \lfloor M, N_j \rfloor) \\
&\cong \mathscr{C}(X, \operatorname*{colim}_j \lfloor M, N_j \rfloor).
\end{aligned}$$ □

By oppositising Proposition 8.42, we obtain the following result.

**Proposition 8.43.** *Assume that $\mathcal{C}$ and $\mathcal{M}$ have enough injectives. Let $M$ be a coclosed object of $\mathcal{M}$. Then $\lceil M, - \rceil$ is left exact if and only if $M$ is $\mathcal{C}$-injective.*

**Definition 8.44.** Let $\mathcal{C}$ be an abelian monoidal, and $\mathcal{M}$ an abelian $\mathcal{C}$-module category, with both having enough projectives. A $\mathcal{C}$-projective object $M$ is called a *$\mathcal{C}$-projective $\mathcal{C}$-generator* if for any projective object $Q$ in $\mathcal{M}$ there exists a projective object $P$ in $\mathcal{C}$, such that $Q$ is a direct summand of $P \triangleright M$.

Analogously, if $\mathcal{C}$ and $\mathcal{M}$ instead have enough injectives, a $\mathcal{C}$-injective object $M$ is called a *$\mathcal{C}$-injective $\mathcal{C}$-cogenerator* if any injective object $J$ of $\mathcal{M}$ is a direct summand of an object of the form $I \triangleright M$, for an injective object $I$ of $\mathcal{C}$.

**Proposition 8.45.** *Assume that $\mathcal{C}$ and $\mathcal{M}$ have enough projectives. Let $M \in \mathcal{M}$ be a closed $\mathcal{C}$-projective object. If $M$ is a $\mathcal{C}$-generator, then $\lfloor M, - \rfloor$ reflects zero objects.*

*If $\mathcal{C}$-*proj *is a Krull–Schmidt category and every indecomposable projective object of $\mathcal{C}$ is the projective cover of a simple object, then we have that $\lfloor M, - \rfloor$ reflects zero objects if and only if $M$ is a $\mathcal{C}$-generator.*

*Proof.* Assume that $M$ is a $\mathcal{C}$-generator. Let $N$ be a non-zero object of $\mathcal{M}$ and let $Q \twoheadrightarrow N$ be an epimorphism from a projective object in $\mathcal{M}$. Let $P$ be a projective object of $\mathcal{C}$ such that $Q$ is a direct summand of $P \triangleright M$. Then $\mathcal{M}(P \triangleright M, N)$ is non-zero, since $\mathcal{M}(Q, N)$ is a direct summand thereof. Thus

$$0 \neq \mathcal{M}(P \triangleright M, N) \cong \mathcal{C}(P, \lfloor M, N \rfloor),$$

showing that $\lfloor M, N \rfloor$ is not zero.

For the latter statement, let $Q'$ be an indecomposable projective object of $\mathcal{M}$, and let $S$ be its simple top. Since $\lfloor M, S \rfloor \neq 0$, there is some $V \in \mathcal{M}$ such that $\mathcal{C}(V, \lfloor M, S \rfloor) \neq 0$. Let $P' \twoheadrightarrow V$ be an epimorphism from a projective object in $\mathcal{C}$. Then

$$0 \neq \mathcal{C}(V, \lfloor M, S \rfloor) \simeq \mathcal{C}(P', \lfloor M, S \rfloor) \simeq \mathcal{C}(P' \triangleright M, S).$$

Since $M$ is $\mathcal{C}$-projective, the object $P' \triangleright M$ is projective. Thus, $\mathcal{C}(P' \triangleright M, S)$ being non-zero implies that $Q'$ is a direct summand of $P' \triangleright M$. □

More closely following the classical case, one alternatively defines a closed $M \in \mathcal{M}$ as a $\mathcal{C}$-generator if $\lfloor M, - \rfloor$ is faithful, see [DSPS19, Definition 2.21].

Oppositisation of Proposition 8.45 yields a similar variant in terms of coclosed and $\mathcal{C}$-injective objects.

**Proposition 8.46.** *Let $\mathscr{C}$ and $\mathscr{M}$ have enough injectives, and let $M \in \mathscr{M}$ be a coclosed $\mathscr{C}$-injective object. If $M$ is a $\mathscr{C}$-cogenerator, then $\lceil M, - \rceil$ reflects zero objects.*

*If $\mathscr{C}$-inj is a Krull–Schmidt category and every indecomposable injective object of $\mathscr{C}$ is the injective hull of a simple object, then $\lceil M, - \rceil$ reflects zero objects if and only if $M$ is a $\mathscr{C}$-cogenerator.*

Finally, we give a brief account of internally projective and injective objects in a semisimple module category.

Since $\mathscr{M}$ is semisimple, any object of $\mathscr{M}$ is $\mathscr{C}$-projective and $\mathscr{C}$-injective.

**Proposition 8.47.** *Let $\mathscr{C}$ be a monoidal category with enough projectives, and let $\mathscr{M}$ be a semisimple $\mathscr{C}$-module category. An object $M \in \mathscr{M}$ is a $\mathscr{C}$-generator if and only if any simple object $S \in \mathscr{M}$ is a direct summand of an object of the form $P \triangleright M$, for some $P \in \mathscr{C}$-proj.*

*Proof.* Recall that, in a semisimple abelian category, every object is projective. In particular, $P \triangleright M$ is always projective, for $M \in \mathscr{M}$ and $P \in \mathscr{C}$-proj, hence every object in $\mathscr{M}$ is $\mathscr{C}$-projective.

Now, let $M \in \mathscr{M}$ be a $\mathscr{C}$-generator. Since a simple object $S \in \mathscr{M}$ is projective, it immediately follows from $M$ being a $\mathscr{C}$-generator that $S$ is a direct summand of $P \triangleright M$, for some $P \in \mathscr{C}$-proj.

Lastly, assume that every simple is a direct summand of $P \triangleright M$, for some $P \in \mathscr{C}$, and let $Q \in \mathscr{M}$-proj. As $\mathscr{M}$ is semisimple, we have $Q \cong \oplus_i S_i$ for simples $S_i$ in $\mathscr{M}$. Since $S_i$ is a direct summand of $P_i \triangleright M$ for some $P_i \in \mathscr{C}$, one calculates

$$Q \cong \oplus_i S_i \subseteq_\oplus \oplus_i (P_i \triangleright M) \cong (\oplus_i P_i) \triangleright M \overset{\text{def}}{=} P \triangleright M. \qquad \square$$

## 8.3 Reconstruction for lax module endofunctors

The following is our most general module categorical reconstruction result.

**Theorem 8.48.** *Let $\mathscr{C}$ be a monoidal abelian category with enough projectives, $\mathscr{M}$ an abelian $\mathscr{C}$-module category with enough projectives, and assume that $\ell \in \mathscr{M}$ is a closed $\mathscr{C}$-projective $\mathscr{C}$-generator. Then there is an equivalence of $\mathscr{C}$-module categories*

$$\mathscr{M} \simeq \mathscr{C}^{\lfloor \ell, - \triangleright \ell \rfloor},$$

*where $\mathscr{C}^{\lfloor \ell, - \triangleright \ell \rfloor}$ is endowed with the extended $\mathscr{C}$-module structure by means of Linton coequalisers, see Theorem 8.25.*

*Proof.* Following Example 2.50, the functor $- \triangleright \ell$ is a strong $\mathscr{C}$-module functor. By Porism 5.29, its right adjoint $\lfloor \ell, - \rfloor$ is a lax $\mathscr{C}$-module functor. Thus, the resulting monad $\lfloor \ell, - \triangleright \ell \rfloor$ is a right exact lax $\mathscr{C}$-module monad. Proposition 5.35 yields that the comparison functor $\mathscr{C}_{\lfloor \ell, - \triangleright \ell \rfloor} \longrightarrow \mathscr{M}$ is a strong $\mathscr{C}$-module functor. Furthermore, due to Propositions 8.42 and 8.45, $\lfloor \ell, - \rfloor$ is right exact and reflects zero objects, so we are able to apply Theorem 2.92, whence the comparison functor $\mathscr{M} \longrightarrow \mathscr{C}^{\lfloor \ell, - \triangleright \ell \rfloor}$ is an equivalence. Transporting the $\mathscr{C}$-module structure along this equivalence, we obtain an extended $\mathscr{C}$-module structure on $\mathscr{C}^{\lfloor \ell, - \triangleright \ell \rfloor}$, which, by Theorem 8.9, is necessarily that of Theorem 8.25; this proves the result. □

Using Theorem 8.25, we can also formulate a converse.

**Theorem 8.49.** *Let $\mathscr{C}$ be a monoidal abelian category with enough projectives, and let $T$ be a right exact lax $\mathscr{C}$-module monad on $\mathscr{C}$. Then $T1$ is a closed $\mathscr{C}$-projective $\mathscr{C}$-generator in $\mathscr{C}^T$, with the $\mathscr{C}$-module category structure on the latter given by Theorem 8.25. There is a bijection*

$$\{ (\mathscr{M}, \ell) \text{ as in Theorem 8.48} \} \big/ _{\mathscr{M} \simeq \mathscr{N}} \overset{\cong}{\longleftrightarrow} \left\{ \begin{array}{l} \text{Right exact lax } \mathscr{C}\text{-module} \\ \text{monads on } \mathscr{C} \end{array} \right\} \Big/ \mathscr{C}^T \simeq \mathscr{C}^S$$

$$(\mathscr{M}, \ell) \longmapsto \lfloor \ell, - \triangleright \ell \rfloor$$

$$(\mathscr{C}^T, T1) \longleftarrow T$$

*Proof.* The fact that $T1$ is closed follows immediately from the isomorphisms

$$\mathscr{C}^T(- \triangleright T1, =) \overset{(i)}{\cong} \mathscr{C}^T(T(- \otimes 1), =) \simeq \mathscr{C}(- \otimes 1, U^T(=)) \simeq \mathscr{C}(-, U^T(=))$$

showing that $\lfloor T1, - \rfloor \cong U^T$, where $(i)$ follows from Lemma 8.23.

For $P \in \mathscr{C}$-proj projective, $P \triangleright T(1) \cong TP$ is projective by Proposition 8.5. By the same result, $T1$ is a $\mathscr{C}$-projective generator, and $\mathscr{C}^T$ has enough projectives.

For the latter claim, notice that $\mathscr{M} \simeq \mathscr{C}^{\lfloor T1, - \triangleright T1 \rfloor}$ by Theorem 8.48. Lastly, $\mathscr{C}^T \simeq \mathscr{C}^{\lfloor T1, - \triangleright T1 \rfloor}$ holds, since from $\lfloor T1, - \rfloor \cong U^T$ and, by uniqueness of adjoints, $F^T \cong - \triangleright T1$ we deduce that $T \cong \lfloor T1, - \triangleright T1 \rfloor$. □

Using Propositions 8.43 and 8.46 in place of Propositions 8.42 and 8.45, we find dual statements.

**Theorem 8.50.** *Let $\mathscr{C}$ be a monoidal abelian category with enough injectives, $\mathscr{M}$ an abelian $\mathscr{C}$-module category with enough injectives, and assume that $\ell \in \mathscr{M}$ be a coclosed $\mathscr{C}$-injective $\mathscr{C}$-cogenerator.*

*Then there is an equivalence of $\mathscr{C}$-module categories*

$$\mathscr{M} \simeq \mathscr{C}^{\lceil \ell, - \triangleright \ell \rceil},$$

*where $\mathscr{C}^{\lceil \ell, - \triangleright \ell \rceil}$ is endowed with the extended $\mathscr{C}$-module structure of Theorem 8.25.*

**Theorem 8.51.** *Let $\mathscr{C}$ be a monoidal abelian category with enough injectives, and $S$ a left exact oplax $\mathscr{C}$-module comonad on $\mathscr{C}$. Then $\mathscr{C}^S$, endowed with the $\mathscr{C}$-module category structure of Theorem 8.25, is a $\mathscr{C}$-module category, and $S1$ is a coclosed $\mathscr{C}$-injective $\mathscr{C}$-generator. There is a bijection*

$$\{\, (\mathscr{M}, \ell) \text{ as in Theorem 8.50} \,\}\Big/_{\mathscr{M} \simeq \mathscr{N}} \leftrightarrow \left\{ \begin{matrix} \text{Left exact oplax } \mathscr{C}\text{-module} \\ \text{comonads on } \mathscr{C} \end{matrix} \right\} \Big/ \, \mathscr{C}^S \simeq \mathscr{C}^G$$

$$(\mathscr{M}, \ell) \longmapsto \lceil \ell, - \triangleright \ell \rceil$$

$$(\mathscr{C}^S, S1) \longleftarrow S$$

Observe that Theorems 8.48 and 8.50 do not make any finiteness assumptions on $\mathscr{M}$—such conditions are often imposed on it by the existence of a closed $\mathscr{C}$-projective $\mathscr{C}$-generator, or a coclosed $\mathscr{C}$-injective $\mathscr{C}$-cogenerator.

**Proposition 8.52.** *Let $\mathscr{C}$ be a finite abelian monoidal category and suppose $\mathscr{M}$ is an abelian $\mathscr{C}$-module category, such that there exists a closed $\mathscr{C}$-projective $\mathscr{C}$-generator $\ell \in \mathscr{M}$. Then $\mathscr{M}$ is finite abelian.*

*Proof.* By Theorem 8.48, we have a $\mathscr{C}$-module equivalence $\mathscr{M} \simeq \mathscr{C}^{\lfloor \ell, - \triangleright \ell \rfloor}$. Since $- \triangleright \ell$ admits a right adjoint it is right exact. Since $\ell$ is $\mathscr{C}$-projective, the functor $\lfloor \ell, - \rfloor$ is also right exact. Thus $\lfloor \ell, - \triangleright \ell \rfloor$ is a right exact monad on the finite abelian category $\mathscr{C}$. By Proposition 8.5, $\mathscr{C}^{\lfloor \ell, - \triangleright \ell \rfloor}$ is finite abelian. □

Using Proposition 8.3 in place of Proposition 8.5, one obtains the following.

**Proposition 8.53.** *Let $\mathscr{C}$ be a locally finite abelian monoidal category, and let $\mathscr{M}$ be an abelian $\mathscr{C}$-module category with enough injectives, such that there exists a coclosed $\mathsf{Ind}(\mathscr{C})$-injective $\mathsf{Ind}(\mathscr{C})$-generator $\ell \in \mathscr{M}$. Then the full subcategory of compact objects in $\mathscr{M}$ is locally finite abelian.*

In the presence of suitable finiteness conditions, the characterisations of adjoint functors between finite and locally finite abelian categories give sufficient conditions for an object in a module category to be closed or coclosed.

**Proposition 8.54.** *Let $\mathscr{C}$ be a finite abelian monoidal category and let $\mathscr{M}$ be a finite abelian $\mathscr{C}$-module category. Then an object $m \in \mathscr{M}$ is closed if and only if the functor $- \triangleright m$ is right exact.*

*Proof.* This is an immediate consequence of Proposition 2.138. □

Similarly, using Proposition 2.137, one shows the following.

**Proposition 8.55.** *Let $\mathscr{C}$ be a locally finite abelian monoidal category and $\mathscr{M}$ a locally finite abelian $\mathscr{C}$-module category. An object $\ell \in \mathsf{Ind}(\mathscr{M})$ is coclosed with respect to the induced $\mathsf{Ind}(\mathscr{C})$-module structure on $\mathsf{Ind}(\mathscr{M})$ if and only if $- \triangleright \ell \colon \mathscr{C} \longrightarrow \mathsf{Ind}(\mathscr{M})$ is left exact and the induced functor $- \triangleright \ell \colon \mathsf{Ind}(\mathscr{C}) \longrightarrow \mathsf{Ind}(\mathscr{M})$ is quasi-finite.*

### 8.3.1 *Rigid monoidal and (finite) tensor categories*

We obtain an algebraic reconstruction result in case $\mathscr{C}$ is rigid monoidal.

**Theorem 8.56.** *Let $\mathscr{C}$ be a rigid monoidal abelian category with enough projectives, $\mathscr{M}$ an abelian $\mathscr{C}$-module category with enough projectives, and assume that $\ell \in \mathscr{M}$ is a closed $\mathscr{C}$-projective $\mathscr{C}$-generator. Then there is an algebra object $A \in \mathscr{C}$ such that there is an equivalence of $\mathscr{C}$-module categories*

$$\mathsf{mod}_{\mathscr{C}} A \simeq \mathscr{M}.$$

*Proof.* By Theorem 8.48, we have $\mathscr{C}^{\lfloor \ell, - \triangleright \ell \rfloor} \simeq \mathscr{M}$ as $\mathscr{C}$-module categories. By Proposition 2.73, the monad $\lfloor \ell, - \triangleright \ell \rfloor \colon \mathscr{C} \longrightarrow \mathscr{C}$ is a strong $\mathscr{C}$-module functor, as it is lax and $\mathscr{C}$ is rigid. The claim follows by Proposition 5.23. □

**Theorem 8.57.** *Let $\mathscr{C}$ be a rigid monoidal abelian category with enough injectives, $\mathscr{M}$ an abelian $\mathscr{C}$-module category, and assume that $\ell \in \mathscr{M}$ is a coclosed $\mathscr{C}$-injective $\mathscr{C}$-cogenerator. Then there is a coalgebra object $C \in \mathscr{C}$ such that there is an equivalence of $\mathscr{C}$-module categories $\mathsf{comod}_{\mathscr{C}} C \simeq \mathscr{M}$.*

Recall the following result in case that $\mathscr{C}$ is a multitensor category.

**Theorem 8.58** ([EGNO15, Theorem 7.10.1]). *Let $\mathscr{C}$ be a finite multitensor category and let $\mathscr{M}$ be an abelian $\mathscr{C}$-module category, such that*

- *the functor $- \triangleright =$ is exact in the first variable (as $\mathscr{C}$ is rigid, it is always exact in the second variable);*
- *there exists a $\mathscr{C}$-projective $\mathscr{C}$-generator.*

*Then there is an algebra object $A$ in $\mathscr{C}$ such that $\mathscr{M} \simeq \mathsf{mod}_{\mathscr{C}}A$.*

We generalise this statement to the locally finite setting.

**Theorem 8.59.** *Let $\mathscr{C}$ be a multitensor category, and let $\mathscr{M}$ be an abelian $\mathscr{C}$-module category such that the $\mathsf{Ind}(\mathscr{C})$-module category $\mathsf{Ind}(\mathscr{M})$ admits a coclosed $\mathsf{Ind}(\mathscr{C})$-injective $\mathsf{Ind}(\mathscr{C})$-cogenerator $\ell$.*

*Then there is a coalgebra object $C$ in $\mathsf{Ind}(\mathscr{C})$ such that $\mathsf{Ind}(\mathscr{M}) \simeq \mathsf{Comod}_{\mathsf{Ind}(\mathscr{C})}C$. Thus, $\mathscr{M}$ is the category of compact $C$-comodule objects.*

*Proof.* By Theorem 8.50, we have $\mathsf{Ind}(\mathscr{M}) \simeq \mathsf{Ind}(\mathscr{C})^{\lceil \ell, - \triangleright \ell \rceil}$ as $\mathsf{Ind}(\mathscr{C})$-module categories. Observe that the $\mathsf{Ind}(\mathscr{C})$-module category structure on both $\mathsf{Ind}(\mathscr{C})$ and on $\mathsf{Ind}(\mathscr{M})$ is the finitary extension of the respective $\mathscr{C}$-module category structures, following Proposition 2.129. Similarly, $- \triangleright \ell \colon \mathsf{Ind}(\mathscr{C}) \longrightarrow \mathsf{Ind}(\mathscr{M})$ is the extension of $- \triangleright \ell \colon \mathscr{C} \longrightarrow \mathsf{Ind}(\mathscr{M})$. On the ind-completion, the left adjoint $\lceil \ell, - \rceil \colon \mathsf{Ind}(\mathscr{M}) \longrightarrow \mathsf{Ind}(\mathscr{C})$, being finitary and preserving compact objects by Lemma 2.139, restricts to an oplax $\mathscr{C}$-module functor $\lceil \ell, - \rceil \colon \mathscr{M} \longrightarrow \mathscr{C}$.

Since $\mathscr{C}$ is rigid, the restricted $\mathscr{C}$-module functor $\lceil \ell, - \rceil \colon \mathscr{M} \longrightarrow \mathscr{C}$ is in fact a strong $\mathscr{C}$-module functor. Its finitary extension $\lceil \ell, - \rceil \colon \mathsf{Ind}(\mathscr{M}) \longrightarrow \mathsf{Ind}(\mathscr{C})$ is a strong $\mathsf{Ind}(\mathscr{C})$-module functor, and the comonad $\lceil \ell, - \triangleright \ell \rceil$ is a strong $\mathsf{Ind}(\mathscr{C})$-module monad. By Proposition 5.23, it is of the form $- \otimes C$ for a coalgebra object in $\mathsf{Ind}(\mathscr{C})$, and hence

$$\mathsf{Ind}(\mathscr{M}) \simeq \mathsf{Ind}(\mathscr{C})^{\lceil \ell, - \triangleright \ell \rceil} \simeq \mathsf{Comod}_{\mathsf{Ind}(\mathscr{C})}C.$$

By Proposition 8.3, $\mathscr{M}$ consists of compact $C$-comodule objects of $\mathscr{C}$. □

**Corollary 8.60.** *Let $H$ be a Hopf algebra and let $\mathscr{C} = {}^{H}\mathsf{vect}$; in particular, we have a monoidal equivalence $\mathsf{Ind}(\mathscr{C}) \simeq {}^{H}\mathsf{Vect}$. Let $\mathscr{M}$ be an abelian $\mathscr{C}$-module category such that the $\mathsf{Ind}(\mathscr{C})$-module category $\mathsf{Ind}(\mathscr{M})$ admits a coclosed $\mathsf{Ind}(\mathscr{C})$-injective $\mathsf{Ind}(\mathscr{C})$-cogenerator. Then there is an $H$-comodule coalgebra $C$ such that there is an equivalence $\mathsf{Ind}(\mathscr{M}) \simeq \mathsf{Comod}_{H}C$ of $\mathsf{Ind}(\mathscr{C})$-module categories, restricting to a $\mathscr{C}$-module equivalence $\mathscr{M} \simeq \mathsf{comod}_{H}C$.*

## 8.4 An Eilenberg–Watts theorem for lax module monads

The formulation of the bijection of Theorem 8.51 indicates a Morita aspect to the reconstruction theory for module categories, as the equivalence relation imposed on monads strongly resembles—and can be specialised to—Morita equivalence of algebras.

In Theorem 8.72, we characterise the right exact lax $\mathscr{C}$-module functors between Eilenberg–Moore categories for lax $\mathscr{C}$-module monads in terms of suitable bimodule objects in the category of endofunctors of $\mathscr{C}$—this gives a precise meaning to the notion of Morita equivalence of monads, hence Theorem 8.75 yields a bijection

$$\{(\mathcal{M},\ell) \text{ as in Theorem 8.50}\}\big/_{\mathcal{M}\simeq\mathcal{N}} \xleftrightarrow{\simeq} \left\{\begin{matrix}\text{Left exact oplax } \mathscr{C}\text{-mod-}\\ \text{ule comonads on } \mathscr{C}\end{matrix}\right\}\Big/_{\simeq \text{Morita}}$$

$$(\mathcal{M},\ell) \longmapsto [\ell, - \triangleright \ell]$$

$$(\mathcal{M}^T, T1) \longleftarrow T$$

Let us first introduce some vocabulary.

**Definition 8.61.** Let $\mathscr{A}$ be an abelian category and suppose that $T$ and $S$ are right exact monads on $\mathscr{A}$. Viewing them as algebra objects in $\mathsf{Rex}(\mathscr{A},\mathscr{A})$, a *$T$-$S$-biact functor* is an $T$-$S$-bimodule object in $\mathsf{Rex}(\mathscr{A},\mathscr{A})$.

In other words, a biact functor $F\colon \mathscr{A} \longrightarrow \mathscr{A}$ is a triple $(F, F_{\mathsf{la}}, F_{\mathsf{ra}})$, consisting of a right exact endofunctor on $\mathscr{A}$ together with transformations $F_{\mathsf{la}}\colon TF \Longrightarrow F$ and $F_{\mathsf{ra}}\colon FS \Longrightarrow F$, satisfying the natural unitality and associativity axioms, and commuting with each other in the sense that

$$\begin{array}{ccc} TFS & \xrightarrow{TF_{\mathsf{ra}}} & TF \\ {\scriptstyle F_{\mathsf{la}}S}\downarrow & & \downarrow{\scriptstyle F_{\mathsf{la}}} \\ FS & \xrightarrow[F_{\mathsf{ra}}]{} & F \end{array} \tag{8.4.1}$$

commutes. The category of $T$-$S$-biact functors shall be denoted by $T$-$\mathsf{Biact}$-$S$.

**Hypothesis 8.62.** For the rest of this section, let us fix two right exact monads $T$ and $S$ on an abelian category $\mathscr{A}$.

**Remark 8.63.** For a $T$-$S$-biact functor $F$, the left $T$-action $F_{\mathsf{la}}$ endows any object of the form $Fa$, for $a \in \mathscr{A}$, with the structure of a $T$-module by defining $\nabla_{Fa} := F_{\mathsf{la};a}$. Any morphism $f \in \mathscr{A}(a,b)$ lifts to a map $Ff\colon Fa \longrightarrow Fb$ of $T$-modules by naturality of $F_{\mathsf{la}}$. In particular, the functor $F$ factors through the canonical forgetful functor $U^T\colon \mathscr{C}^T \longrightarrow \mathscr{C}$. We write $F = U^T \circ \tilde{F}$.

**Definition 8.64.** Given a $T$-$S$-biact functor $F$, define $\tilde{F} \circ_S -\colon \mathcal{A}^S \longrightarrow \mathcal{A}^T$ by

$$\tilde{F} \circ_S - := \operatorname{coeq}\left( FS \underset{F_{\mathrm{ra}}}{\overset{F\nabla}{\rightrightarrows}} F \right), \tag{8.4.2}$$

where for $x \in \mathcal{A}^S$ we endow $Fx$ and $FSx$ with the $T$-module structures described in Remark 8.63.

The functoriality of Definition 8.64 follows from that of colimits.

**Remark 8.65.** By naturality of $F_{\mathsf{la}}$ and the commutativity condition of Diagram (8.4.1), both morphisms in Equation (8.4.2) are morphisms in $\mathcal{A}^T$; since the coequalisers in $\mathcal{A}^T$ are created by $U^T$, this coequaliser is one in $\mathcal{A}^T$.

**Proposition 8.66.** *Let $T$ and $S$ be right exact monads on an abelian category $\mathcal{A}$. Writing $\varepsilon^T := \varepsilon^{F^T \dashv U^T}$ and $\varepsilon^S := \varepsilon^{F^S \dashv U^S}$, the following assignments extend to an equivalence of categories:*

Notice the similarity of this assignment and Definition 6.1, only that we are now pushing $\Phi$ along two adjunctions.

$$\begin{aligned} \mathsf{Rex}(\mathcal{A}^S, \mathcal{A}^T) &\longleftrightarrow T\text{-}\mathsf{Biact}\text{-}S \\ \Phi &\longmapsto \left(U^T \Phi F^S,\ U^T \varepsilon^T \Phi F^S,\ U^T \Phi \varepsilon^S F^S\right) \\ \tilde{F} \circ_S - &\longleftarrow\!\shortmid F. \end{aligned}$$

*Proof.* Since $F^S$ is right exact and $U^T$ is exact, $U^T \Phi F^S$ is right exact.

It is easy to verify that the transformations $U^S \varepsilon^T \Phi F^T$ and $U^T \Phi \varepsilon^S F^S$ endow $U^T \Phi F^S$ with the structure of a biact functor, and that for a natural transformation $\alpha\colon \Phi \Longrightarrow \Phi'$, the resulting transformation $U^T \alpha F^S$ is a morphism of biact functors. This proves the functoriality of $\mathsf{Rex}(\mathcal{A}^S, \mathcal{A}^T) \longrightarrow T\text{-}\mathsf{Biact}\text{-}S$.

To see that the converse assignment is functorial, observe that given a morphism $f\colon F \Longrightarrow G$ of biact functors, we obtain a morphism of forks in $\mathcal{A}^T$:

$$\begin{array}{ccc} FSx & \underset{F_{\mathrm{ra};x}}{\overset{F\nabla_x}{\rightrightarrows}} & Fx \\ {\scriptstyle f_{Sx}}\downarrow & & \downarrow{\scriptstyle f_x} \\ GSx & \underset{G_{\mathrm{ra};x}}{\overset{G\nabla_x}{\rightrightarrows}} & Gx \end{array}$$

This induces a morphism of coequalisers $\tilde{F} \circ_S - \Longrightarrow \tilde{G} \circ_S -$.

We now prove that these functors define mutually quasi-inverse equivalences. First, observe that $\widetilde{U^T \Phi F^S} = \Phi F^S \colon \mathscr{A} \longrightarrow \mathscr{A}^T$, since $U^T$ is faithful and injective on objects, so it is a monomorphism in the 1-category $\mathsf{Cat}_{\Bbbk}$. By slight abuse of notation, see Remark 8.65, for $x \in \mathscr{A}^S$ we have

$$(\Phi F^S) \circ_S x = \mathrm{coeq}\left( \Phi S^2 x \underset{\Phi S \nabla_x}{\overset{\Phi \mu^S_x}{\rightrightarrows}} \Phi S x \right),$$

which, since $\Phi$ is right exact, is isomorphic to

$$\Phi\left( \mathrm{coeq}\left( S^2 x \underset{S \nabla_x}{\overset{\mu^S_x}{\rightrightarrows}} S x \right)\right).$$

Via the action $\nabla_x \colon Tx \longrightarrow x$, the latter coequaliser is canonically isomorphic to $x$, hence we obtain the following isomorphism, which is natural in $\Phi$:

$$\Phi \hat{\nabla}_x \colon (\Phi F^S) \circ_S x \xrightarrow{\sim} \Phi x.$$

Conversely, let $F$ be a biact functor and $x \in \mathscr{A}$. We have

$$\left( U^T (\tilde{F} \circ_S -) F^S \right)(x) = \mathrm{coeq}\left( F S^2 x \underset{F \nabla_{Sx} = F \mu^S_x}{\overset{F_{\mathsf{ra};Sx}}{\rightrightarrows}} F S x \right)$$

Observe that $F$, being a right $S$-module object, is presented as a coequaliser in the category of such objects, which additionally splits in the underlying monoidal category $\mathsf{Rex}(\mathscr{A}, \mathscr{A})$. This is indicated in the following diagram:

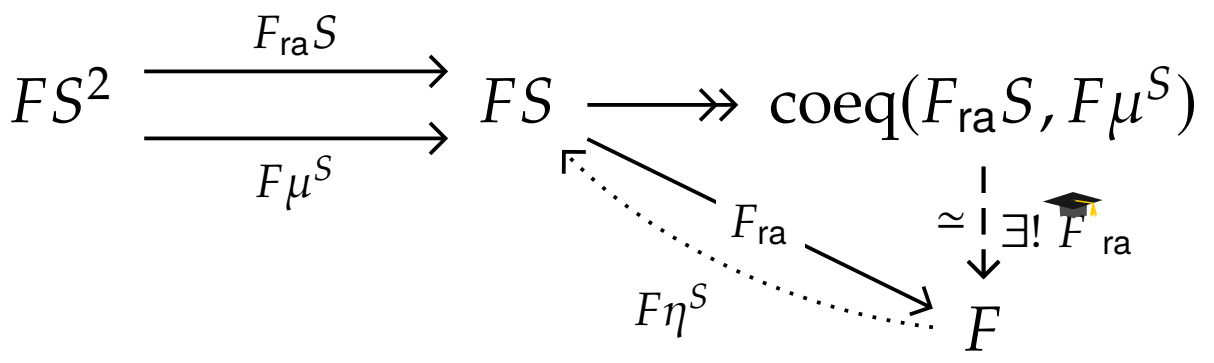

Since $F_{\mathsf{ra}}$ is a morphism of left $T$-module objects, we find that $\hat{F}_{\mathsf{ra}}$ defines an isomorphism of biact functors $U^T (\tilde{F} \circ_S -) F^S \cong F$ that is natural in $F$. ◻

Let us now assume $\mathscr{C}$ and $\mathscr{M}$ to be abelian, $- \triangleright =$ to be right exact in both variables, and $S$ and $T$ to be right exact lax $\mathscr{C}$-module monads on $\mathscr{M}$.

**Definition 8.67.** An $S$-$T$-bimodule object in the category $\mathsf{RexLax}\mathscr{C}\mathsf{Mod}(\mathscr{M}, \mathscr{M})$ of right exact lax $\mathscr{C}$-module endofunctors of $\mathscr{M}$ is called a *lax $\mathscr{C}$-module $S$-$T$-biact functor*. We denote the category of $S$-$T$-biact functors by $S\text{-}\mathsf{Biact}_{\mathscr{C}}\text{-}T$.

In other words, a lax $\mathscr{C}$-module biact functor is a biact functor $(F, F_{\mathsf{la}}, F_{\mathsf{ra}})$, such that $F\colon \mathscr{M} \longrightarrow \mathscr{M}$ is a right exact lax $\mathscr{C}$-module functor and $F_{\mathsf{la}}$ and $F_{\mathsf{ra}}$ are $\mathscr{C}$-module transformations.

**Theorem 8.68.** *The equivalence* $\mathsf{Rex}(\mathscr{M}^T, \mathscr{M}^S) \simeq S\text{-}\mathsf{Biact}\text{-}T$ *of Proposition 8.66 restricts to a faithful functor*

$$\mathscr{C}\mathsf{EW}\colon \mathsf{LaxRex}(\mathscr{M}^T, \mathscr{M}^S) \longrightarrow S\text{-}\mathsf{Biact}_{\mathscr{C}}\text{-}T,$$

*such that the following diagram commutes*:

$$\begin{array}{ccc} \mathsf{LaxRex}(\mathscr{M}^T, \mathscr{M}^S) & \xrightarrow{\mathscr{C}\mathsf{EW}} & S\text{-}\mathsf{Biact}_{\mathscr{C}}\text{-}T \\ \downarrow & & \downarrow \\ \mathsf{Rex}(\mathscr{M}^T, \mathscr{M}^S) & \xrightarrow{\simeq \text{ by } 8.66} & S\text{-}\mathsf{Biact}\text{-}T \end{array} \tag{8.4.3}$$

*The vertical arrows are the respective forgetful functors, and the categories $\mathscr{M}^T$ and $\mathscr{M}^S$ are endowed with the $\mathscr{C}$-module structures of Theorem 8.25.*

*Proof.* By Theorem 8.26, similarly to [AHLF18, Proposition 3.10], the functor $F^T\colon \mathscr{M} \longrightarrow \mathscr{M}^T$ is a strong $\mathscr{C}$-module functor, and so is $F^S$. Thus, by Porism 5.29, the functor $U^S$ is a lax $\mathscr{C}$-module functor, and hence for a lax $\mathscr{C}$-module functor $\Phi\colon \mathscr{M}^T \longrightarrow \mathscr{M}^S$, the composite $U^S \Phi F^T$ is a lax $\mathscr{C}$-module functor. Further, the $S$-$T$-biact structure on $U^S \Phi F^S$ given in Proposition 8.66 assembles to a lax $\mathscr{C}$-module biact functor, and this extends also to $\mathscr{C}$-module biact transformations arising from $\mathscr{C}$-module transformations.

This defines a functor $\mathscr{C}\mathsf{EW}$, such that Diagram (8.4.3) commutes. It is faithful, since so are the remaining functors in that diagram. □

**Lemma 8.69.** *Let $\Phi \in \mathsf{Rex}(\mathscr{M}^T, \mathscr{M}^S)$ and for all $m \in \mathscr{M}$ suppose that*

$$(\Phi T)_{\mathsf{a};m}\colon x \triangleright_{\mathscr{M}} \Phi T m \longrightarrow \Phi T(x \triangleright_{\mathscr{M}} m)$$

*is a lax $\mathscr{C}$-module biact functor structure on the biact functor $U^S \Phi F^T$. Then we have $(\Phi T)_{\mathsf{a}} = \mathscr{C}\mathsf{EW}(\Phi_{\mathsf{a};T(-)})$, where*

$$\Phi_{\mathsf{a};T(-)}\colon x \triangleright_{\mathscr{M}^S} \Phi T(-) \Longrightarrow \Phi(x \triangleright_{\mathscr{M}^T} T(-)) = \Phi(T(x \triangleright_{\mathscr{M}} -))$$

*is the unique natural transformation making the following diagram commute:*

$$\begin{array}{ccccc} x \triangleright_{\mathcal{M}^S} \Phi T- & & & \overset{\exists!\Phi_{\mathsf{a}}}{\dashrightarrow} & \\ \uparrow\!\!\uparrow & & & & \\ S(x \triangleright_{\mathcal{M}} \Phi T-) & \xrightarrow[S(\Phi T)_{\mathsf{a}}]{} & S\Phi T x \triangleright_{\mathcal{M}} - & \xrightarrow[\Phi_{\mathsf{la}} T x \triangleright -]{} & \Phi T x \triangleright_{\mathcal{M}} - \end{array} \tag{8.4.4}$$

*Proof.* In order to verify that $\Phi_{\mathsf{a}}$ is well-defined, observe that the morphism $(\Phi_{\mathsf{la}} T x \triangleright_{\mathcal{M}} -) \circ S(\Phi T)_{\mathsf{a}}$ coequalises $SS_{\mathsf{a};x}\Phi T$ and $Sx \triangleright_{\mathcal{M}} \Phi_{\mathsf{la}} T$:

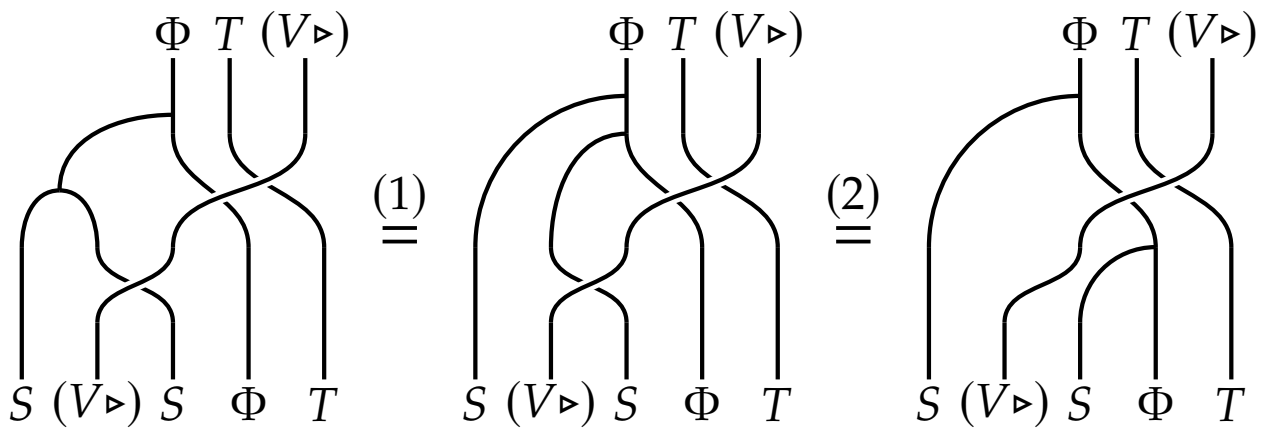

Above, (1) follows by $(\Phi T)_{\mathsf{la}}$ being a left $S$-act structure, and (2) by it being a $\mathcal{C}$-module transformation.

By the definition of the functor $\mathcal{C}\mathsf{EW}$, we have

$$\mathcal{C}\mathsf{EW}(\Phi_{\mathsf{a};T(-)}) = \Phi_{\mathsf{a};T(-)} \circ U^S_{\mathsf{a};x,\Phi T(-)},$$

where $U^S_{\mathsf{a}}$ is the lax $\mathcal{C}$-module functor structure on $U^S$ induced by Porism 5.29.

More generally, for $m \in \mathcal{M}$ the map $U^S_{\mathsf{a};x,m}$ is given by the composite

$$x \triangleright_{\mathcal{M}} m \xrightarrow{\eta^S_{x \triangleright_{\mathcal{M}} m}} S(x \triangleright_{\mathcal{M}} m) \twoheadrightarrow x \triangleright_{\mathcal{M}^S} m,$$

with the latter map being the projection onto the following coequaliser:

$$\begin{array}{ccc} \mathrm{coeq}\Big(T(x \triangleright T^2 m) \rightrightarrows T(x \triangleright Tm)\Big) & = & S(x \triangleright m) \\ \Big\downarrow T(x \triangleright \nabla_m) \qquad\qquad \Big\downarrow T(x \triangleright \nabla_m) & & \Big\downarrow\!\!\downarrow \\ \mathrm{coeq}\Big(T(x \triangleright T^2 m) \rightrightarrows Tx \triangleright Tm\Big) & = & x \triangleright m \end{array}$$

All in all, we find that $\mathcal{C}\mathsf{EW}(\Phi_{\mathsf{a};T(-)})$ is given by the outer path of

$$\begin{array}{ccccccc} & & x \triangleright_{\mathcal{M}^S} \Phi T- & & & \overset{\exists!\Phi_{\mathsf{a}}}{\dashrightarrow} & \\ & & \uparrow\!\!\uparrow & & & & \\ S(x \triangleright_{\mathcal{M}} \Phi T-) & \xrightarrow[\eta_{x \triangleright_{\mathcal{M}} \Phi T-}]{} & S(x \triangleright_{\mathcal{M}} \Phi T-) & \xrightarrow[S(\Phi T)_{\mathsf{a}}]{} & S\Phi T x \triangleright_{\mathcal{M}} - & \xrightarrow[\Phi_{\mathsf{la}} T x \triangleright -]{} & \Phi T x \triangleright_{\mathcal{M}} - \end{array}$$

which, by the unitality of the $S$-act structure $(\Phi T)_{\mathsf{la}}$, equals $(\Phi T)_{\mathsf{a}}$. □

For $m \in \mathcal{M}^T$, let $\Phi_{\mathsf{a};m}$ be the unique map such that

$$
\begin{array}{ccc}
x \triangleright_{\mathcal{M}^S} \Phi T m & \xrightarrow{x \triangleright_{\mathcal{M}^S} \Phi \nabla_m} & x \triangleright_{\mathcal{M}^S} \Phi m \\
\downarrow{\scriptstyle \Phi_{\mathsf{a};Tm}} & & \\
\Phi T(x \triangleright m) & & \vdots\ {\scriptstyle \exists ! \Phi_{\mathsf{a};m}} \\
\downarrow{\scriptstyle =} & & \downarrow \\
\Phi(x \triangleright_{\mathcal{M}^T} Tm) & \xrightarrow[\Phi(x \triangleright_{\mathcal{M}^T} \nabla_m)]{} & \Phi(x \triangleright_{\mathcal{M}^T} m)
\end{array}
\tag{8.4.5}
$$

commutes. By the uniqueness in the defining property of $\Phi_{\mathsf{a};m}$, this defines a natural transformation

$$\Phi_{\mathsf{a}} \colon - \triangleright_{\mathcal{M}^S} \Phi(=) \Longrightarrow \Phi(- \triangleright_{\mathcal{M}^T} =).$$

**Lemma 8.70.** *For $\Phi \in \mathsf{Rex}(\mathcal{M}^T, \mathcal{M}^S)$, the map $\Phi_{\mathsf{a}} \colon - \triangleright_{\mathcal{M}^S} \Phi(=) \Longrightarrow \Phi(- \triangleright_{\mathcal{M}^S} =)$ extended from the morphisms $\Phi_{\mathsf{a};Tm}$ of Diagram (8.4.5) defines a lax $\mathcal{C}$-module structure on $\Phi$.*

*Proof.* In order to simplify notation, we write $\triangleright$ for $\triangleright_{\mathcal{M}}$ and $\blacktriangleright$ for $\triangleright_{\mathcal{M}^T}$ and $\triangleright_{\mathcal{M}^S}$, since each occurrence of either of these symbols is unambiguous with respect to which is used. We also write $x\triangleright$ for $x \triangleright -$, and suppress the horizontal composition symbols in $\mathbb{C}\mathsf{at}_{\Bbbk}$, replacing them simply by concatenation. Hence, every occurrence of $\circ$ in the diagrams that follow is a vertical composition of natural transformations. For example,

$$S(y\triangleright)S(x\triangleright)\Phi T \xrightarrow{S(y\triangleright)(\Phi_{\mathsf{la}}T(x\triangleright)\circ S(\Phi T)_{\mathsf{a},x})} S(y\triangleright)\Phi T(x\triangleright)$$

represents the natural transformation

$$S(y\triangleright S(x\triangleright \Phi T(-))) \xrightarrow{S(y\triangleright S(\Phi T)_{\mathsf{a};y,x})} S(y\triangleright S\Phi T(x\triangleright -)) \xrightarrow{S(y\triangleright \Phi_{\mathsf{la}}T(x\triangleright -))} S(y\triangleright \Phi T(x\triangleright -)).$$

Since $\Phi_{\mathsf{a}}$ is uniquely determined by $\Phi_{\mathsf{a};T(-)}$, it suffices to show that

$$
\begin{array}{ccccc}
y \blacktriangleright \Phi(x \blacktriangleright T) & & \xrightarrow{\Phi_{\mathsf{a};y,x\blacktriangleright T}} & & \Phi(y \blacktriangleright x \blacktriangleright T) \\
\uparrow{\scriptstyle y\blacktriangleright \Phi_{\mathsf{a};x,T}} & \diagdown\!\!\diagdown & & \diagup\!\!\diagup & \\
y \blacktriangleright x \blacktriangleright \Phi T & \xrightarrow{\alpha_1} & y \blacktriangleright \Phi T(x\triangleright) \xrightarrow{\alpha_3} \Phi T(y\triangleright)(x\triangleright) & & \downarrow{\scriptstyle \Phi\alpha_2} \\
\downarrow{\scriptstyle \alpha_2} & & \swarrow{\scriptstyle \alpha_5} & & \\
y \otimes x \blacktriangleright \Phi T & \xrightarrow{\alpha_4} & \Phi T(y \otimes x)\triangleright & = & \Phi((y \otimes x) \blacktriangleright T)
\end{array}
$$

commutes, where we define $\alpha_1 := y \blacktriangleright \Phi_{\mathsf{a},x}$, $\alpha_3 := \Phi_{\mathsf{a},y}(x\triangleright)$, $\alpha_4 := \Phi_{\mathsf{a},y\otimes x}$, and $\alpha_5 := \Phi T(\triangleright_{\mathsf{a},y,x})$. The morphism $\alpha_2$ is obtained from Figure 8.1, in which every marked epimorphism is the coequaliser of the pair preceding it, every dashed arrow is an induced morphism between coequalisers coming from a morphism of diagrams, and the dotted arrows come from the universal property of coequalisers. The notation for the coequaliser $(y, x) \blacktriangleright \Phi T$ is to emphasise the connection with the multiactegorical approach of Section 8.1.3.

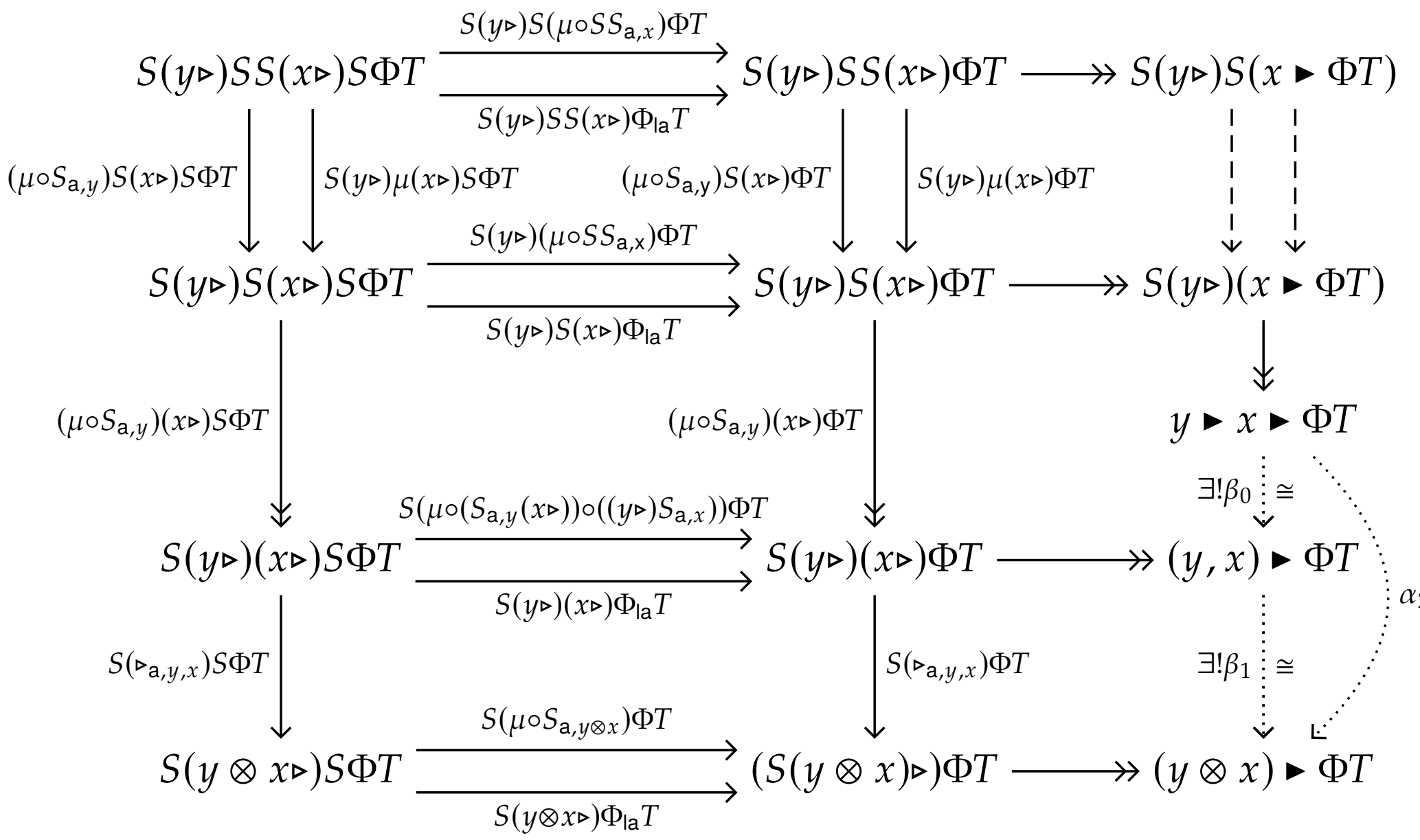

Figure 8.1: Definition of the morphism $\alpha_2$.

We now have the diagram

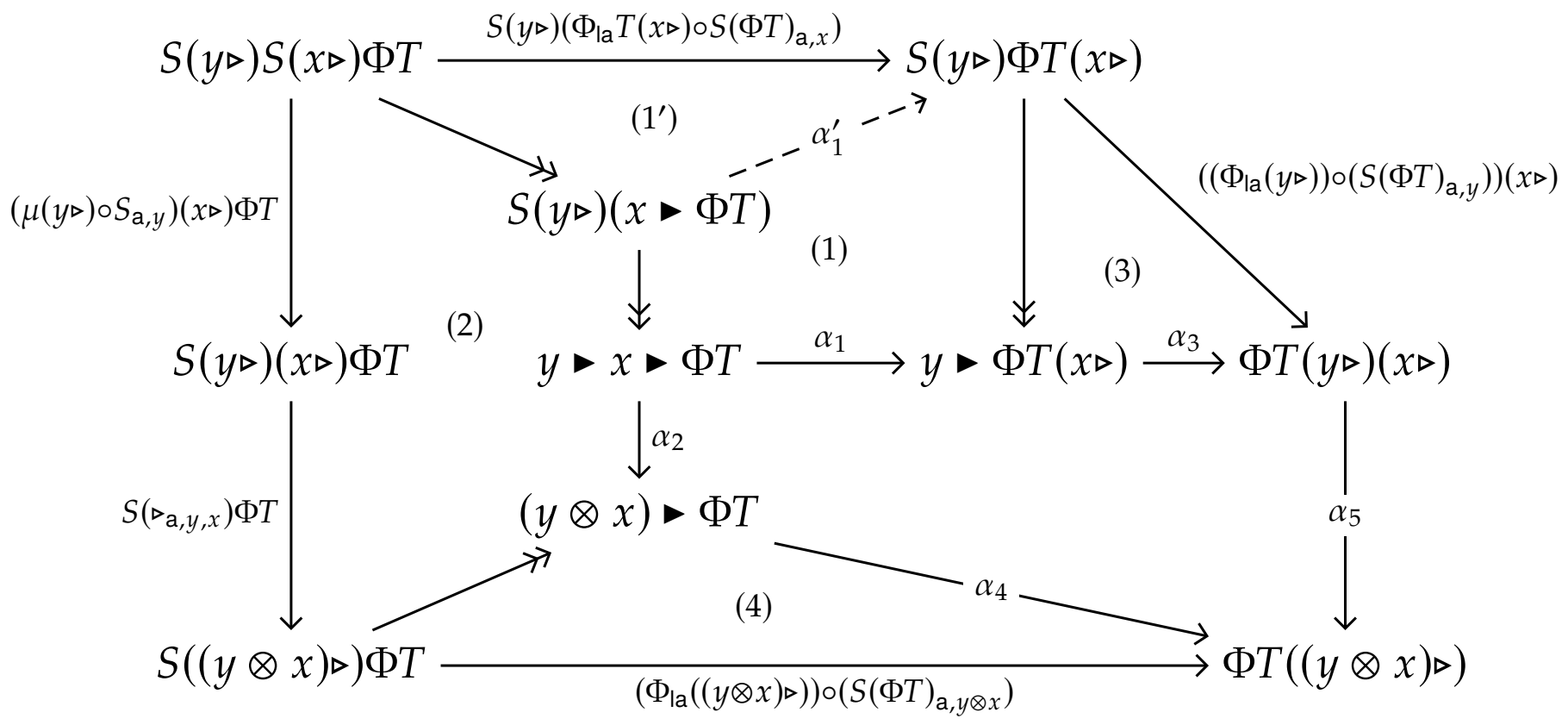

whose labelled faces all commute, and where the morphism decorated by the label of the face it is part of is defined as that making said face commute, via the universal property of coequalisers.

Our aim is to show the commutativity of the inner unlabelled face. Since all of its morphisms are defined by the remaining inner (commutative) faces, and we can reach $y \blacktriangleright x \blacktriangleright \Phi T$ from $S(y\triangleright)S(x\triangleright)\Phi T$ with two epimorphisms, this is implied by the commutativity of the outer face, which follows by the commutativity of the inner faces of Figure 8.2, where

- faces (1), (3), and (4) commute by the interchange law in $\mathbb{C}\mathrm{at}_{\Bbbk}$;
- face (2) commutes by $\Phi_{\mathsf{la}} \colon S\Phi T \Longrightarrow \Phi T$ being a $\mathscr{C}$-module map;
- face (5) commutes by the associativity of action $\Phi_{\mathsf{la}}$;
- face (6) commutes by $\Phi T$ being a lax $\mathscr{C}$-module functor; and
- face (7) commutes by naturality of $\Phi_{\mathsf{la}}$. □

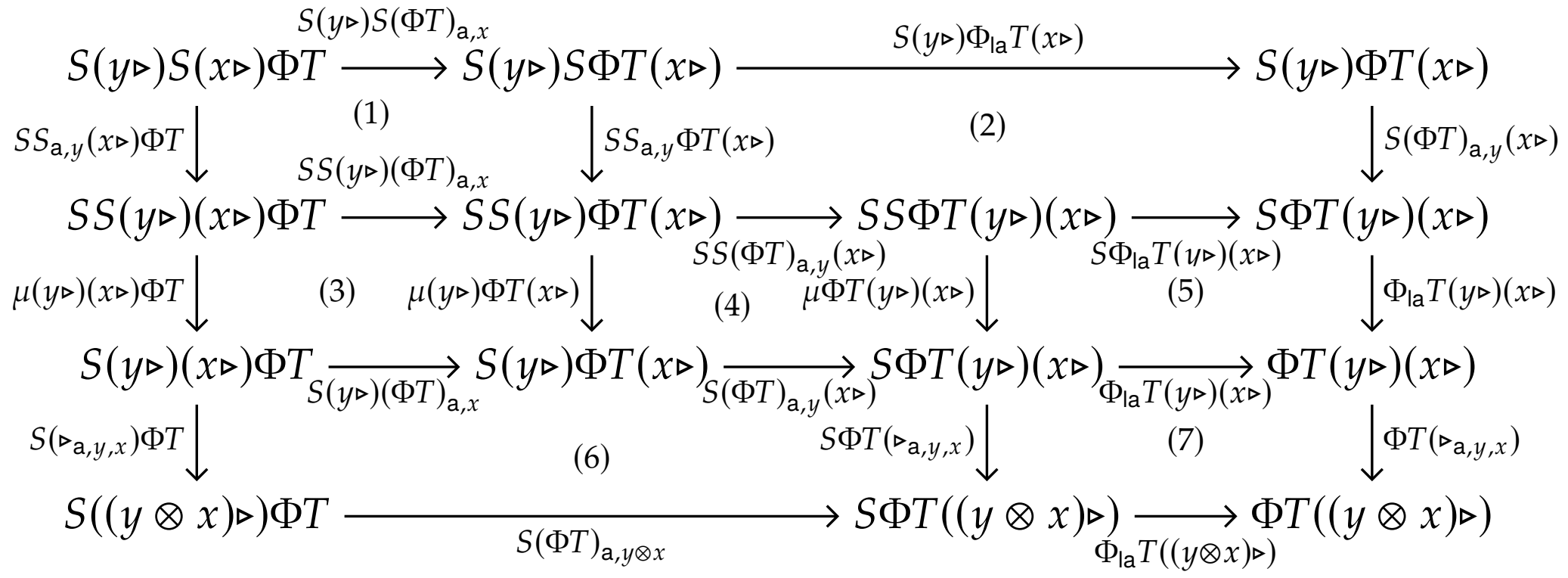

Figure 8.2: The arrow $\Phi_{\mathsf{a}}$ defines a lax $\mathscr{C}$-module structure on $\Phi$.

**Lemma 8.71.** *Let $\Phi, \Phi' \in \mathsf{LaxRex}(\mathscr{M}^T, \mathscr{M}^S)$ and $\phi \colon U^S\Phi F^T \Longrightarrow U^S\Phi' F^T$ be a lax $\mathscr{C}$-module biact transformation. Then*

$$\varphi := \phi \circ_T - \colon \Phi \Longrightarrow \Phi'$$

*is a $\mathscr{C}$-module transformation.*

*Proof.* First, observe that $\varphi$ is determined by $\varphi_{T(-)}$, as indicated by the diagram

$$\begin{array}{ccccc}
\Phi T^2 m & \overset{\Phi T\nabla_m}{\underset{\Phi\mu_m}{\rightrightarrows}} & \Phi T m & \twoheadrightarrow & \Phi m \\
\downarrow{\scriptstyle \varphi_{T^2 m}} & & \downarrow{\scriptstyle \varphi_{Tm}} & & \downarrow{\scriptstyle \exists! \varphi_m} \\
\Phi' T^2 m & \overset{\Phi' T\nabla_m}{\underset{\Phi'\mu_m}{\rightrightarrows}} & \Phi' T m & \twoheadrightarrow & \Phi' m
\end{array}$$

Thus, it suffices to show that for any $m \in \mathscr{M}$ and $x \in \mathscr{C}$, we have that

(8.4.6) $$\Phi'_{\mathsf{a};x,Tm} \circ (x \blacktriangleright \phi_m) = \phi_{x\triangleright m} \circ \Phi_{\mathsf{a},x,T(m)}.$$

Consider the following diagram:

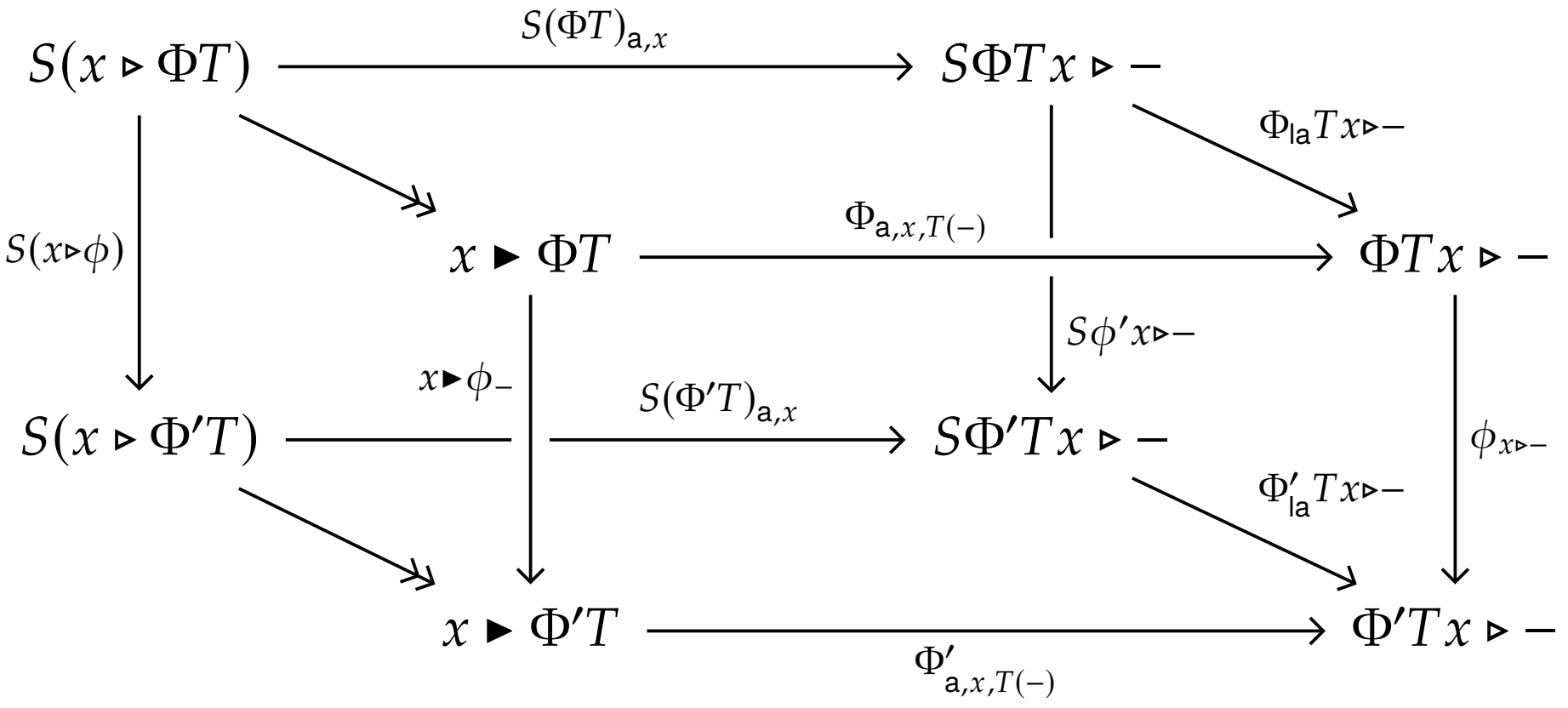

Its top and bottom faces commute by definition of $\Phi_{\mathsf{a}}$, as seen from Diagram (8.4.4); its left face commutes by the definition of $x \blacktriangleright -$; and its right face by $\phi$ being a biact transformation. The back face commutes since $\phi$ is a $\mathscr{C}$-module transformation. The commutativity of the front face is precisely Equation (8.4.6), and since the top-left projection map is an epimorphism, it suffices to verify its commutativity after precomposing with it, which follows easily from the commutativity of the remaining faces. $\square$

**Theorem 8.72.** *The faithful functor*

$$\mathscr{C}\mathsf{EW}\colon \mathsf{LaxRex}(\mathscr{M}^T, \mathscr{M}^S) \longrightarrow S\text{-}\mathsf{Biact}_{\mathscr{C}}\text{-}T$$

*of Theorem 8.68 is an equivalence of categories.*

*Proof.* Fullness follows from Lemma 8.71, so it is left to prove that $\mathscr{C}\mathsf{EW}$ is essentially surjective. Suppose that $F \in T\text{-}\mathsf{Biact}_{\mathscr{C}}\text{-}S$. Since the functor of Proposition 8.66 is essentially surjective, we may assume that $F \cong U^S \Phi F^T$, for some right exact $\Phi\colon \mathscr{M}^T \longrightarrow \mathscr{M}^S$. By Lemma 8.70, $(\Phi, \Phi_{\mathsf{a}}) \in \mathsf{LaxRex}(\mathscr{M}^T, \mathscr{M}^S)$, and $F \cong \mathscr{C}\mathsf{EW}((\Phi, \Phi_{\mathsf{a}}))$ due to Lemma 8.69, giving essential surjectivity. $\square$

**Definition 8.73.** We say that two right exact lax $\mathscr{C}$-module monads $T$ and $S$ on $\mathscr{M}$ are *Morita equivalent* if there are $F \in T\text{-}\mathsf{Biact}_{\mathscr{C}}\text{-}S$ and $G \in S\text{-}\mathsf{Biact}_{\mathscr{C}}\text{-}T$ such that $G \circ_T F \cong \mathrm{Id}_{\mathscr{M}^T}$ and $F \circ_S G \cong \mathrm{Id}_{\mathscr{M}^S}$ as lax $\mathscr{C}$-module biact functors.

**Proposition 8.74.** *Two right exact lax $\mathscr{C}$-module monads $T$ and $S$ on $\mathscr{M}$ are Morita equivalent if and only if there is a $\mathscr{C}$-module equivalence $\mathscr{M}^T \simeq \mathscr{M}^S$, where the Eilenberg–Moore categories are endowed with the extended Linton coequaliser structure of Theorem 8.25.*

*Proof.* This is a direct consequence of Theorem 8.72, using that an equivalence of categories is right exact. □

Combining Proposition 8.74 with Theorem 8.51, we find the following.

**Theorem 8.75.** *Let $\mathscr{C}$ be a monoidal abelian category with enough injectives, and let $T$ be a left exact lax $\mathscr{C}$-module comonad on $\mathscr{C}$. There is a bijection*

$$\left\{ (\mathscr{M}, \ell) \text{ as in Theorem 8.50} \right\} \big/_{\mathscr{M} \simeq \mathscr{N}} \xleftrightarrow{\simeq} \left\{ \begin{matrix} \text{Left exact oplax } \mathscr{C}\text{-module} \\ \text{comonads on } \mathscr{C} \end{matrix} \right\} \Big/ \simeq_{\mathit{Morita}}$$

$$(\mathscr{M}, \ell) \longmapsto \lceil \ell, - \triangleright \ell \rceil$$

$$(\mathscr{M}^T, T1) \longleftarrow T$$

I am big! It's the pictures that got small.

NORMA DESMOND; Sunset Boulevard

# 9 HOPF TRIMODULES

THEOREMS 8.48 TO 8.50 GIVE GENERAL RECONSTRUCTION RESULTS in the non-rigid setting, which has recently seen increased interest in multiple areas, see for example [DSPS19; ALSW21] or [Str24b, Corollary 10.11]. However, lax $\mathscr{C}$-module monads on $\mathscr{C}$ are less accessible than mere algebra objects of $\mathscr{C}$, and may indeed appear not very accessible in general. To show that this need not be the case, we again turn to the Hopf-algebraic setting.

Recall that a ring category—the non-rigid counterpart of a tensor category, see [EGNO15, Sections 4.2 and 5.4]—which admits a fibre functor to $\mathsf{vect}$ is monoidally equivalent to the category of finite-dimensional left $B$-comodules over a (not necessarily Hopf) bialgebra $B$. In this chapter, we show that the category of left exact finitary lax ${}^{B}\mathsf{Vect}$-module endofunctors of ${}^{B}\mathsf{Vect}$ is monoidally equivalent to the category ${}^{B}_{B}\mathsf{Vect}^{B}$ of *Hopf trimodules*[20]: $B$-$B$-bicomodules with an additional left $B$-action that is a $B$-$B$-bicomodule morphism.

[20] The "trimodule" terminology is based on [Sho09] and was suggested by Ulrich Krähmer, in lieu of more unwieldy notations like "(2, 1)-Hopf module".

**Theorem 9.2.** *For a bialgebra $B$ and $\mathscr{V} := {}^{B}\mathsf{Vect}$, there is a monoidal equivalence*

$$ {}^{B}_{B}\mathsf{Vect}^{B} \simeq \mathsf{LexfLax}\mathscr{V}\mathsf{Mod}(\mathscr{V}, \mathscr{V}) $$

*between the category of Hopf trimodules, and the category of left exact finitary lax $\mathscr{V}$-module endofunctors on $\mathscr{V}$.*

Such structures feature prominently in the quasi-bialgebraic generalisation of the fundamental theorem of Hopf modules, see [HN99; Sar17]. This equivalence matches a lax ${}^{B}\mathsf{Vect}$-module monad with a Hopf trimodule algebra: a Hopf trimodule $A$ together with maps $A \square A \longrightarrow A$ and $B \longrightarrow A$ satisfying the usual associativity and unitality axioms for an algebra object.

If $B$ is infinite-dimensional, Theorem 8.50 should yield reconstruction in terms of an *oplax* ${}^{B}\mathsf{Vect}$-module comonad $\lceil X, - \triangleright X \rceil$, using the existence of injective objects in ${}^{B}\mathsf{Vect}$. Since our result gives an algebraic realisation only for the *lax* ${}^{B}\mathsf{Vect}$-module functors, we restrict ourselves to the case where $- \triangleright X$ admits both a left and a right adjoint, by assuming $- \triangleright X$ to be exact. In that case, we obtain an adjunction $\lceil X, - \triangleright X \rceil \dashv \lfloor X, - \triangleright X \rfloor$.

As the monad is right adjoint, the Eilenberg–Moore category of $\lceil X, - \triangleright X \rceil$ is not equivalent to the category of modules over the associated Hopf trimodule algebra $\lfloor X, X \rfloor$, but rather to the category of its *contramodules*: structures extensively studied in the setting of so-called semi-infinite homological algebra, [Pos10]. On the other hand, the Kleisli categories are equivalent, and once again control the ${}^{B}\mathsf{Vect}$-module structure, which can thus be read off directly from the Hopf trimodule algebra. We obtain the following result.

**Theorem 9.23.** *Let $\mathcal{M}$ be a locally finite abelian ${}^{B}\mathsf{vect}$-module category satisfying the assumptions of Theorem 8.50. There exists a Hopf trimodule algebra $A \in {}^{B}_{B}\mathsf{vect}^{B}$, such that there is an equivalence $\mathsf{Ind}(\mathcal{M}) \simeq A\text{-Contramod}$ of ${}^{B}\mathsf{vect}$-module categories, where the ${}^{B}\mathsf{vect}$-module structure on the right-hand side is extended from the category of free $A$-module. This equivalence restricts to $\mathcal{M} \simeq A\text{-contramod}$.*

Since this algebraic realisation of our reconstruction results may at first glance seem difficult to apply in calculations, we give two explicit examples in which we determine a trimodule algebra for a given module category.

The first is Section 9.4, where $B$ is the semigroup algebra for the unique two-element monoid which is not a group. This example is intended to be very similar to [DSPS19, Example 2.20], in that a simple object of $\mathscr{C}$ acts as zero on an indecomposable object. We also show that applying the ordinary, "rigid" reconstruction procedure on this example yields the same algebra object in ${}^{B}\mathsf{vect}$ as that corresponding to the regular action of ${}^{B}\mathsf{vect}$ on itself, and thus fails to classify module categories.

Our second example is given in Section 9.5. Supposing $B$ to be arbitrary, we determine the Hopf trimodule algebra underlying the module functor given by the fibre functor ${}^{B}\mathsf{vect} \longrightarrow \mathsf{vect}$.

Finally, the Morita-theoretic results of Theorem 8.75 give us a precise notion of Morita equivalence for Hopf trimodule algebras with respect to their contramodules in Definition 9.26. This results in a Morita Theorem for ${}^{B}\mathsf{vect}$-module categories, Theorem 9.28.

As another application, we use trimodule reconstruction to give a categorical interpretation of a variant of the fundamental theorem of Hopf modules.

**Proposition 9.38 and Corollary 9.45.** *The functor ${}^{B}\mathsf{Vect} \longrightarrow {}^{B}_{B}\mathsf{Vect}^{B}$ corresponds to the inclusion of strong ${}^{B}\mathsf{Vect}$-module endofunctors in the category of lax ${}^{B}\mathsf{Vect}$-module endofunctors, under the Yoneda lemma and Theorem 9.2. The latter—and* *hence also the former—functor is an equivalence if and only if* ${}^{B}\mathsf{vect}$ *is left rigid, which is the case if and only if* $B$ *has a twisted antipode.*

In Proposition 9.49, we show that the fusion operators of a Hopf monad $T$ on a monoidal category $\mathcal{W}$ can be interpreted as coherence cells for the canonical oplax module structure on $T$. This gives a strong converse to the fact that (op)lax module functors over a rigid category are automatically strong, by providing a distinguished module functor which is strong if and only if the category is rigid. This result, interpreting additional structure of a monad on a monoidal category as morphisms for an (op)lax module functor structure on said monad, is similar to the results of [FLP24, Theorems 3.17 and 3.18], which studies Frobenius monoidal functors rather than Hopf monads.

## 9.1 Bicomodules and their graphical calculus

Before introducing Hopf trimodules, let us first talk about bicomodules over $B$ in the sense of Section 2.6 and their string diagrammatic representation.

**Hypothesis 9.1.** For this chapter, fix a field $\Bbbk$ and a bialgebra $B$ over it. Write $\mathcal{V} := {}^{B}\mathsf{Vect}$ for the category of left $B$-comodules.

After introducing all of the relevant concepts and notation, our first goal is to prove the following theorem.

**Theorem 9.2.** *There is a monoidal equivalence*

$$\begin{aligned} {}^{B}_{B}\mathsf{Vect}^{B} &\longrightarrow \mathsf{LexfLax}\mathcal{V}\mathsf{Mod}(\mathcal{V},\mathcal{V}) \\ X &\longmapsto (X \,\square\, -, \chi) \end{aligned} \tag{9.1.1}$$

*between the category of Hopf trimodules, and the category of left exact finitary lax* $\mathcal{V}$*-module endofunctors on* $\mathcal{V}$*, where* $\chi$ *is the interchange morphism of Definition 9.6, see Equation (9.2.4)*

The quasi-inverse of Equation (9.1.1) evaluates a left exact finitary lax $\mathcal{V}$-module functor at the injective generator $B$. We defer the proof of Theorem 9.2 until Section 9.2.

**Notation 9.3.** A bicomodule over $B$ consists of a vector space $X$, a left coaction $\lambda\colon X \longrightarrow B \otimes_{\Bbbk} X$, and a right coaction $\rho\colon X \longrightarrow X \otimes_{\Bbbk} B$, such that

$$(B \otimes_{\Bbbk} \rho) \circ \lambda = (\lambda \otimes_{\Bbbk} B) \circ \rho.$$

In graphical notation, we write

(9.1.2)

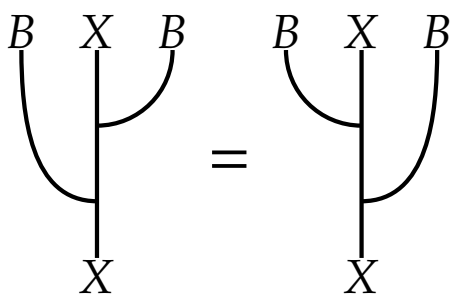

**Example 9.4.** Analogously to Example 2.108, we can form the *cotensor product* of two bicomodules $(X, \lambda^X, \rho^X)$ and $(Y, \lambda^Y, \rho^Y)$ over $B$. This is the equaliser

$$X \,\square\, Y \hookrightarrow X \otimes_{\Bbbk} Y \underset{X \otimes_{\Bbbk} \lambda^Y}{\overset{\rho^X \otimes_{\Bbbk} Y}{\rightrightarrows}} X \otimes_{\Bbbk} B \otimes_{\Bbbk} Y.$$

Our graphical calculus has to differentiate between the tensor product of two bicomodules over $\Bbbk$ and their cotensor product. In particular, we have to indicate which additional transformations are possible with the latter: the equalised actions in the cotensor product will be annotated in grey.

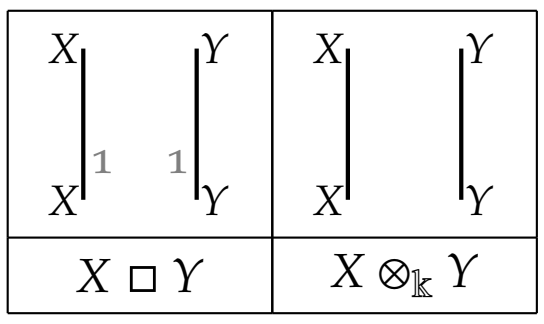

The fact that we may interchange the appropriate left and right actions will be indicated thusly:

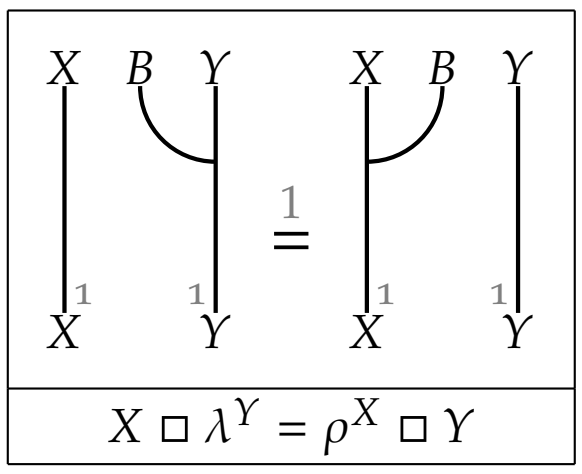

## 9.2 From Hopf trimodules to lax module functors

In this chapter we give a proof of Theorem 9.2.

**Definition 9.5.** A *Hopf trimodule* over $B$ consists of a bicomodule $X$ together with a left $B$-action $\alpha\colon B \otimes_{\Bbbk} X \longrightarrow X$, such that $\alpha$ is a left and right $B$-comodule morphism. We write $X \in {}^B_B\mathsf{Vect}^B$.

Alternatively, Definition 9.5 could impose the conditions that $\lambda$ and $\rho$ are left and right $B$-module morphisms, respectively. This equivalence is easily seen in a string diagrammatic reformulation:

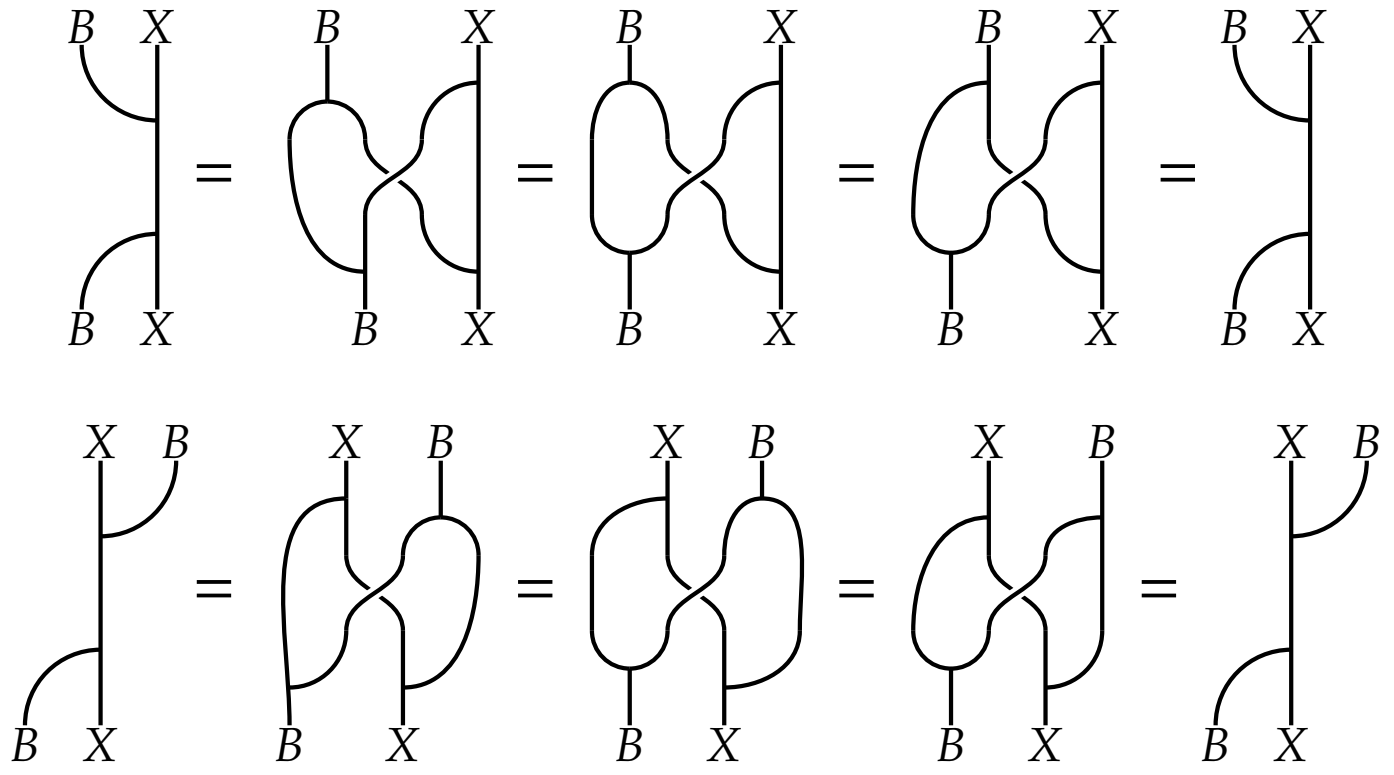

(9.2.1)

**Definition 9.6.** Let $X \in {}^B_B\mathsf{Vect}^B$. For all $M, N \in {}^B\mathsf{Vect}$, define the *interchange morphism* $\chi_{M,N}\colon M \otimes_{\Bbbk} (X \,\square\, N) \longrightarrow X \,\square\, (M \otimes_{\Bbbk} N)$ by

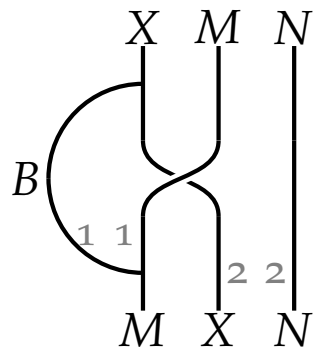

Notice the similarity of the interchange morphism with the Yetter–Drinfeld braiding of Example 2.54.

Observe in particular that

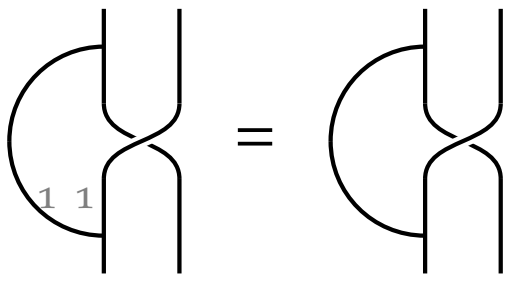

because the image of the coaction $M \longrightarrow B \otimes_{\Bbbk} M$ is contained in the cotensor product $B \,\square\, M$.

**Lemma 9.7.** *The interchange morphism is a well-defined left $B$-comodule morphism.*

*Proof.* To show well-definedness, we need to show that the image of $\chi_{M,N}$ is contained in $X \,\square\, (M \otimes_{\Bbbk} N)$, for all $M, N \in {}^B\mathsf{Vect}$. This follows by the calculation in Figure 9.1.

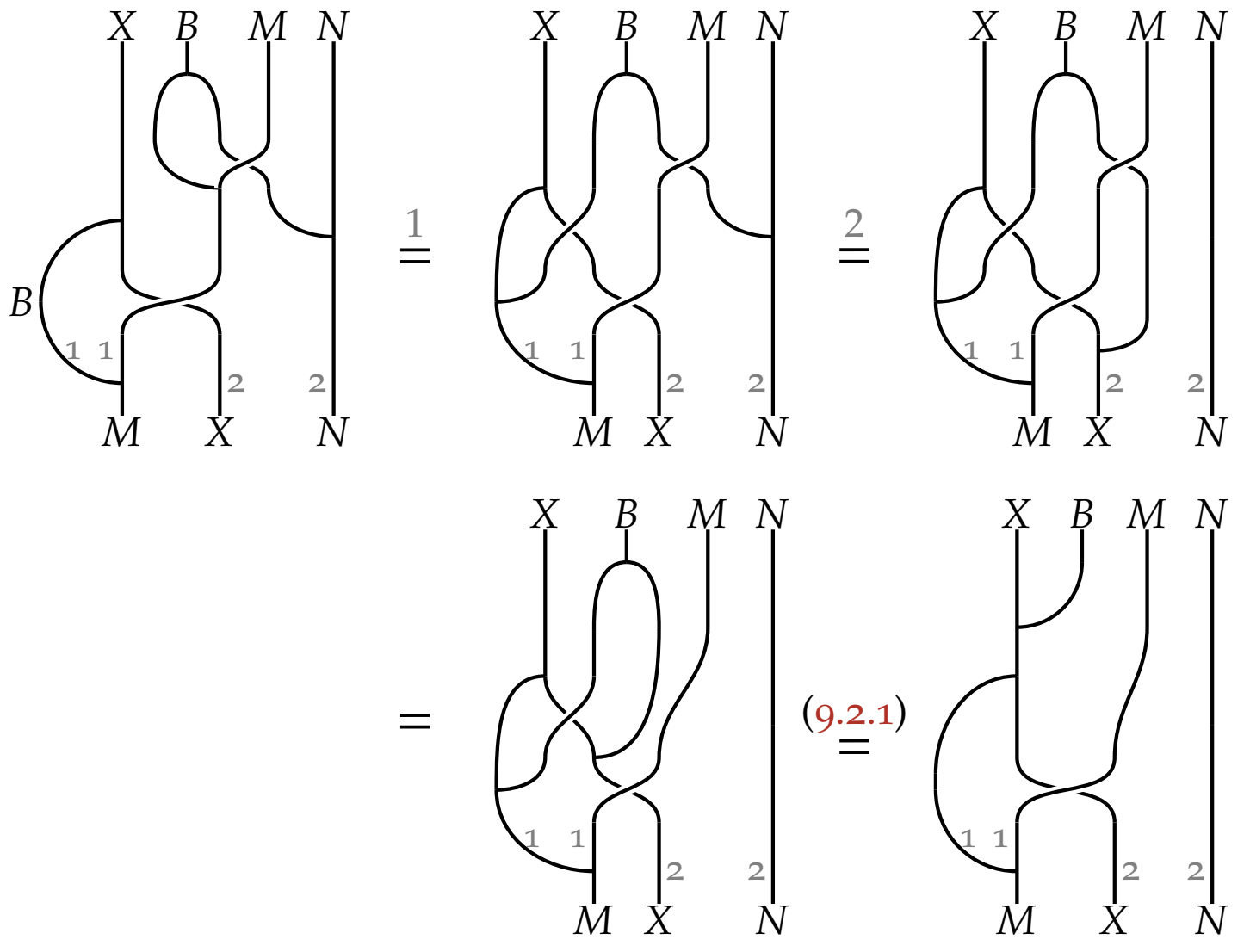

Figure 9.1: The image of $\chi_{M,N}$ is contained in $X \,\square\, (M \otimes_{\Bbbk} N)$.

Graphically, the fact that $\chi$ is a left $B$-comodule morphism means that

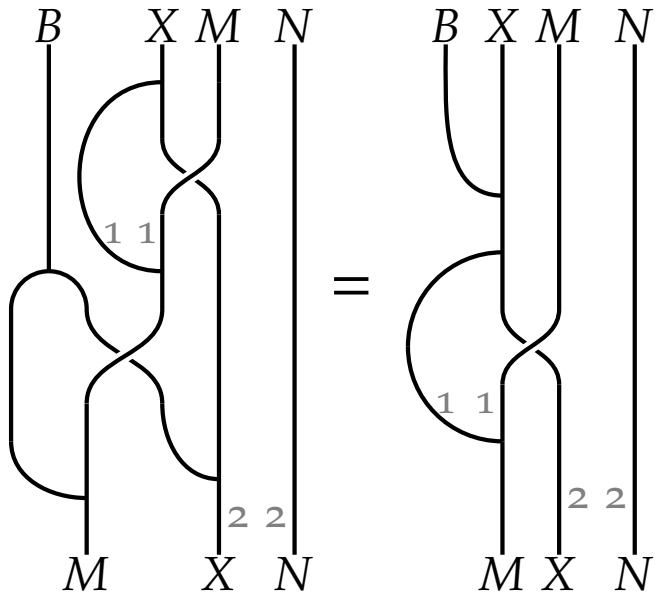

Using this, the calculation in Figure 9.2 finishes the proof. □

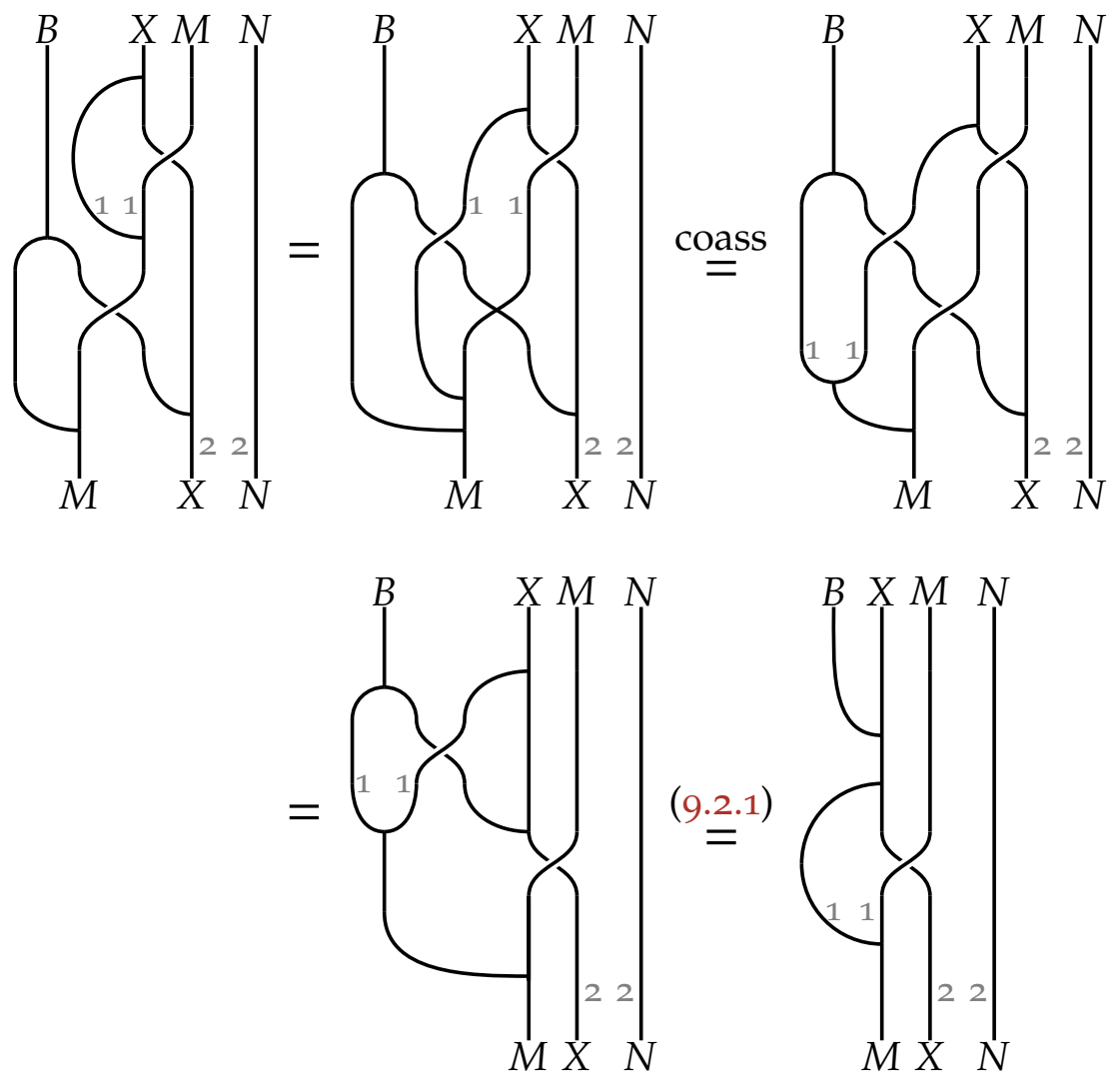

Figure 9.2: The arrow $\chi$ is a morphism of left $B$-comodules.

**Lemma 9.8.** *The interchange morphism is a braiding; i.e., it is natural in both variables and the following diagrams commute for all $M, N, P \in {}^{B}\mathsf{Vect}$:*

$$
\begin{array}{ccc}
M \otimes_{\Bbbk} (N \otimes_{\Bbbk} (X \square P)) & \xrightarrow{\quad M \otimes_{\Bbbk} \chi_{N,P} \quad} & M \otimes_{\Bbbk} (X \square (N \otimes_{\Bbbk} P)) \\
{\scriptstyle \alpha}\downarrow & & \downarrow{\scriptstyle \chi_{M, N \otimes_{\Bbbk} P}} \\
(M \otimes_{\Bbbk} N) \otimes_{\Bbbk} (X \square P) \xrightarrow[\chi_{M \otimes_{\Bbbk} N, P}]{} & X \square ((M \otimes_{\Bbbk} N) \otimes_{\Bbbk} P) & \xrightarrow[\alpha]{} X \square (M \otimes_{\Bbbk} (N \otimes_{\Bbbk} P))
\end{array}
\tag{9.2.2}
$$

$$
\begin{array}{ccc}
\Bbbk \otimes_{\Bbbk} (X \square M) & \xrightarrow{\quad \chi_{\Bbbk, N} \quad} & X \square (\Bbbk \otimes_{\Bbbk} M) \\
& {\scriptstyle \lambda_{X \square M}} \searrow \qquad \swarrow {\scriptstyle X \square \lambda_N} & \\
& X \square N &
\end{array}
\tag{9.2.3}
$$

*Proof.* The interchange morphism is always natural in its second variable. To prove naturality in the first variable, let $f \colon M \longrightarrow M'$ be a left comodule morphism. Then

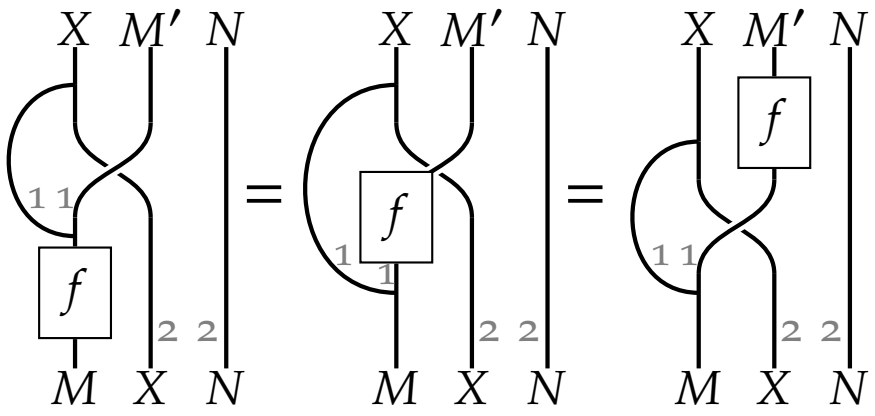

where the first equality follows by $f$ being a left comodule morphism, and the second one is the naturality of the braiding.

Diagram (9.2.2) commuting is equivalent to the following diagram, which is seen to be true by associativity of the action on $X$.

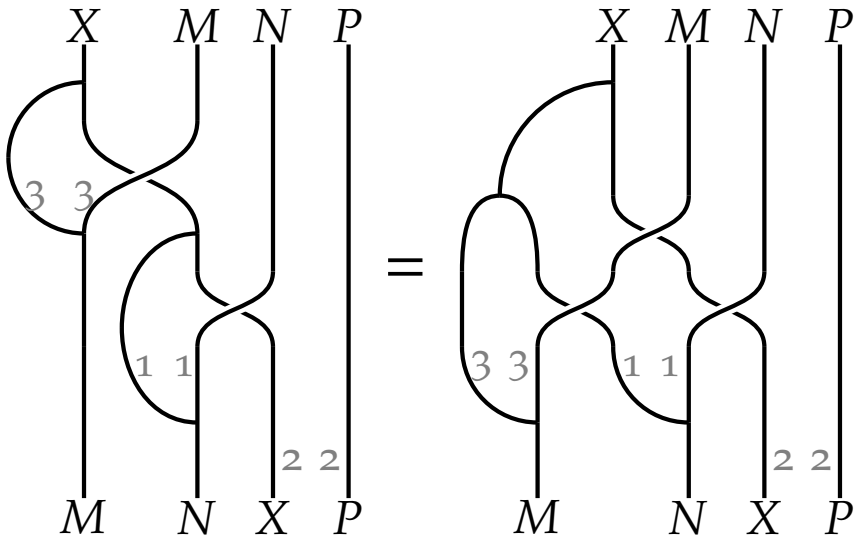

Diagram (9.2.3) follows by the unitality of the action on $X$. □

Lemmas 9.7 and 9.8 taken together say that the well-defined arrow

$$\chi \colon - \otimes_{\Bbbk} (X \square =) \Longrightarrow X \square (- \otimes_{\Bbbk} =)$$

satisfies Diagrams (9.2.2) and (9.2.3). This, in turn, yields the following result.

**Proposition 9.9.** *The pair $(X \,\square\, -, \chi)$ defines a lax* ${}^{B}\mathsf{Vect}$*-module functor.*

We can extend this correspondence to morphisms.

**Lemma 9.10.** *Let $f \in {}^{B}_{B}\mathsf{Vect}^{B}(X, Y)$ be a morphism of Hopf trimodules. Then*

$$f \,\square\, -\colon X \,\square\, - \Longrightarrow Y \,\square\, -$$

*is a* ${}^{B}\mathsf{Vect}$*-module transformation.*

*Proof.* We have to prove that

$$\big(f \,\square\, (M \otimes_{\Bbbk} N)\big) \circ \chi^{X}_{M,N} = \chi^{Y}_{M,N} \circ \big(M \otimes_{\Bbbk} (f \,\square\, N)\big).$$

In our graphical language, this means that

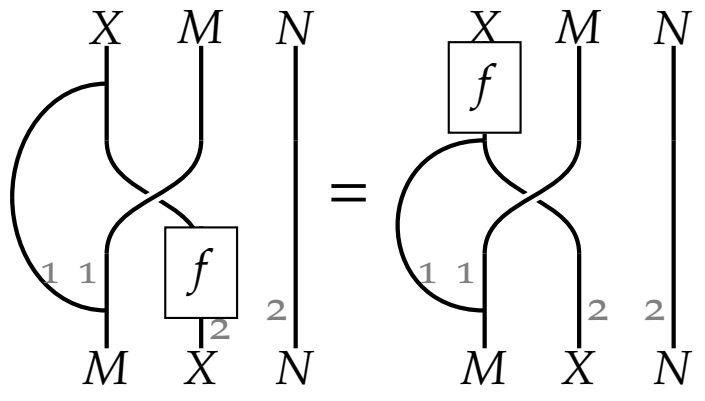

which follows immediately from $f$ being a module morphism. □

Recall the notation $\mathcal{V} := {}^{B}\mathsf{Vect}$. As a result of the previous considerations, there exists a well-defined functor

(9.2.4) $$\Sigma\colon {}^{B}_{B}\mathsf{Vect}^{B} \longrightarrow \mathsf{LexfLax}\mathcal{V}\mathsf{Mod}(\mathcal{V}, \mathcal{V}), \qquad X \longmapsto (X \,\square\, -, \chi),$$

To finish the proof of Theorem 9.2, we have to show that $\Sigma$ is monoidal, as well as an equivalence of categories.

**Lemma 9.11.** *Let $X, Y \in {}^{B}_{B}\mathsf{Vect}^{B}$. Their cotensor product $X \,\square\, Y$ is a Hopf trimodule, where the left action is given diagonally by $b(x \otimes y) := b_{(1)}x \otimes b_{(2)}y$.*

*Proof.* First, we show that the diagonal action of $B$ is well-defined as a map $B \otimes_{\Bbbk} (X \,\square\, Y) \longrightarrow X \,\square\, Y$; i.e., that

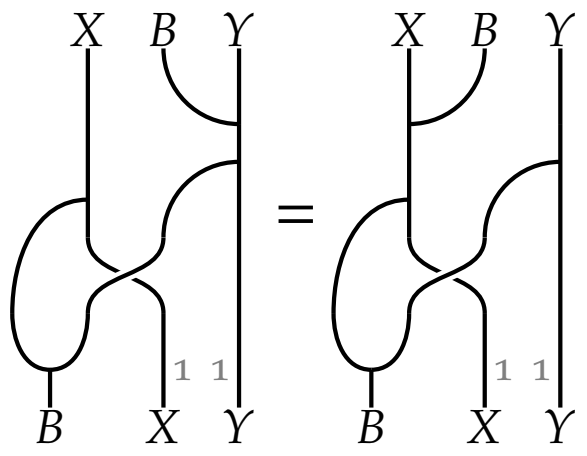

This follows by the calculation in Figure 9.3.

It is left to show that the action is a left and right $B$-comodule morphism. The former case follows by Figure 9.4, and in the latter case one calculates:

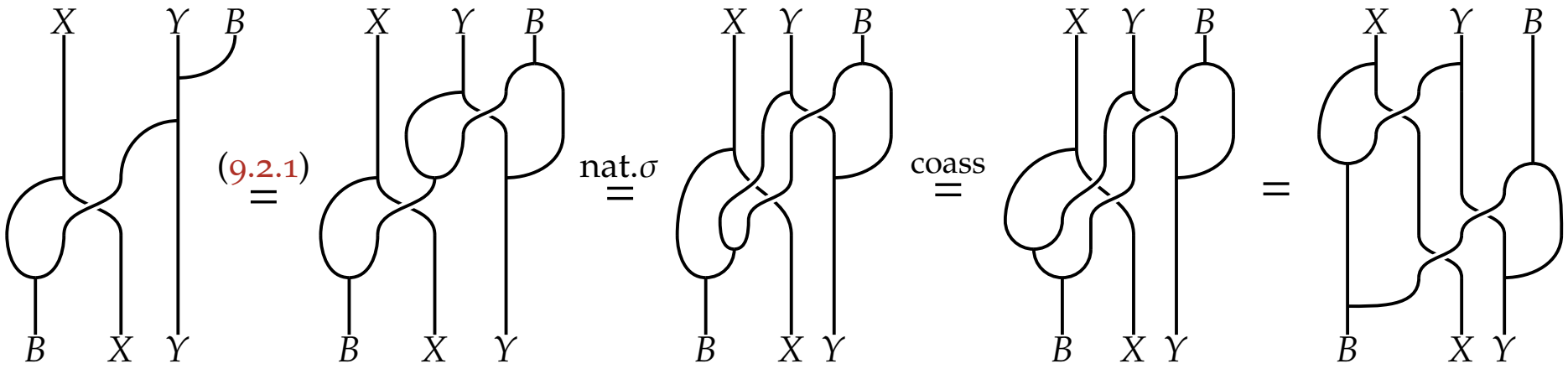

□

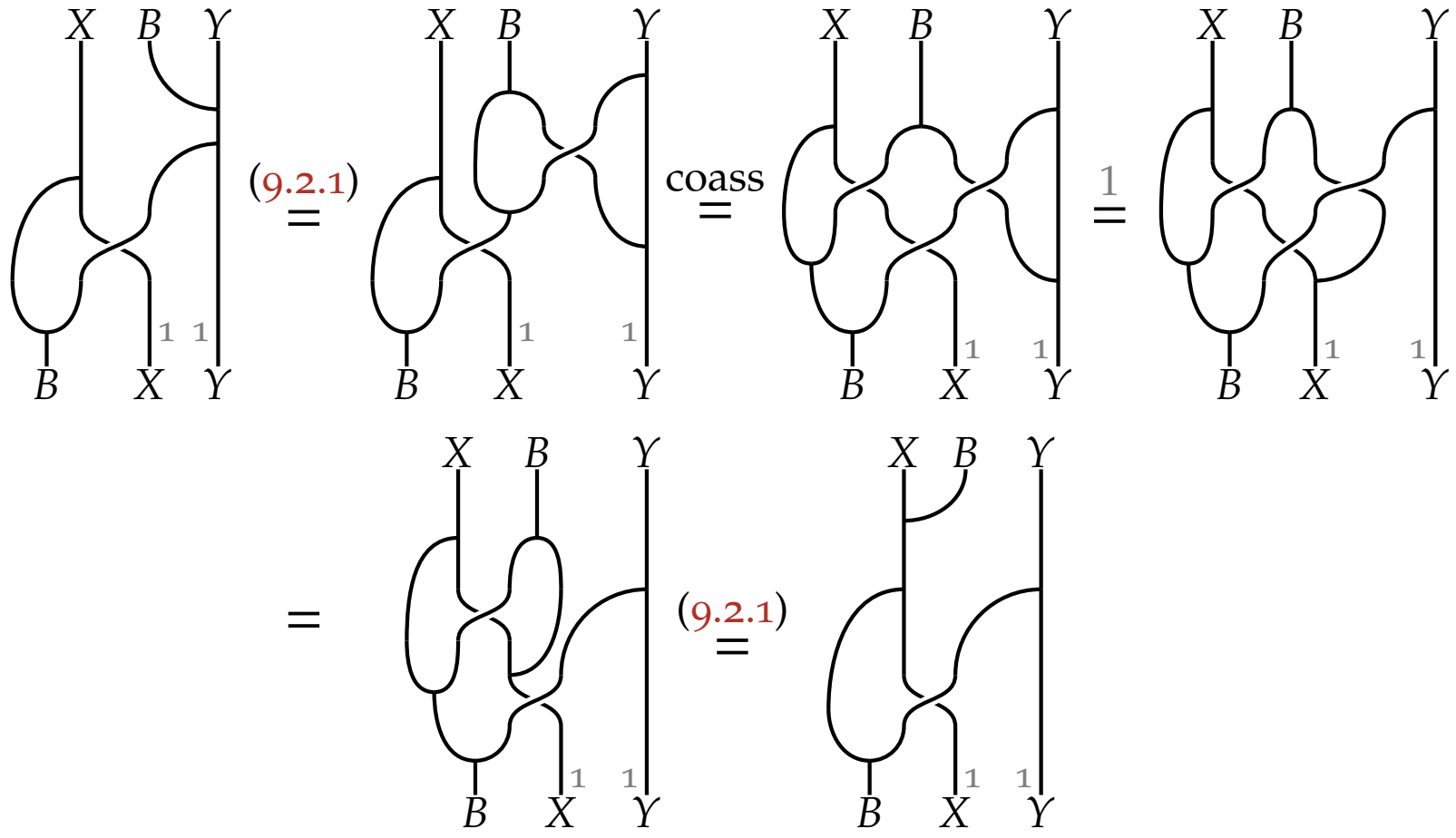

Figure 9.3: The diagonal action is well-defined on the cotensor product.

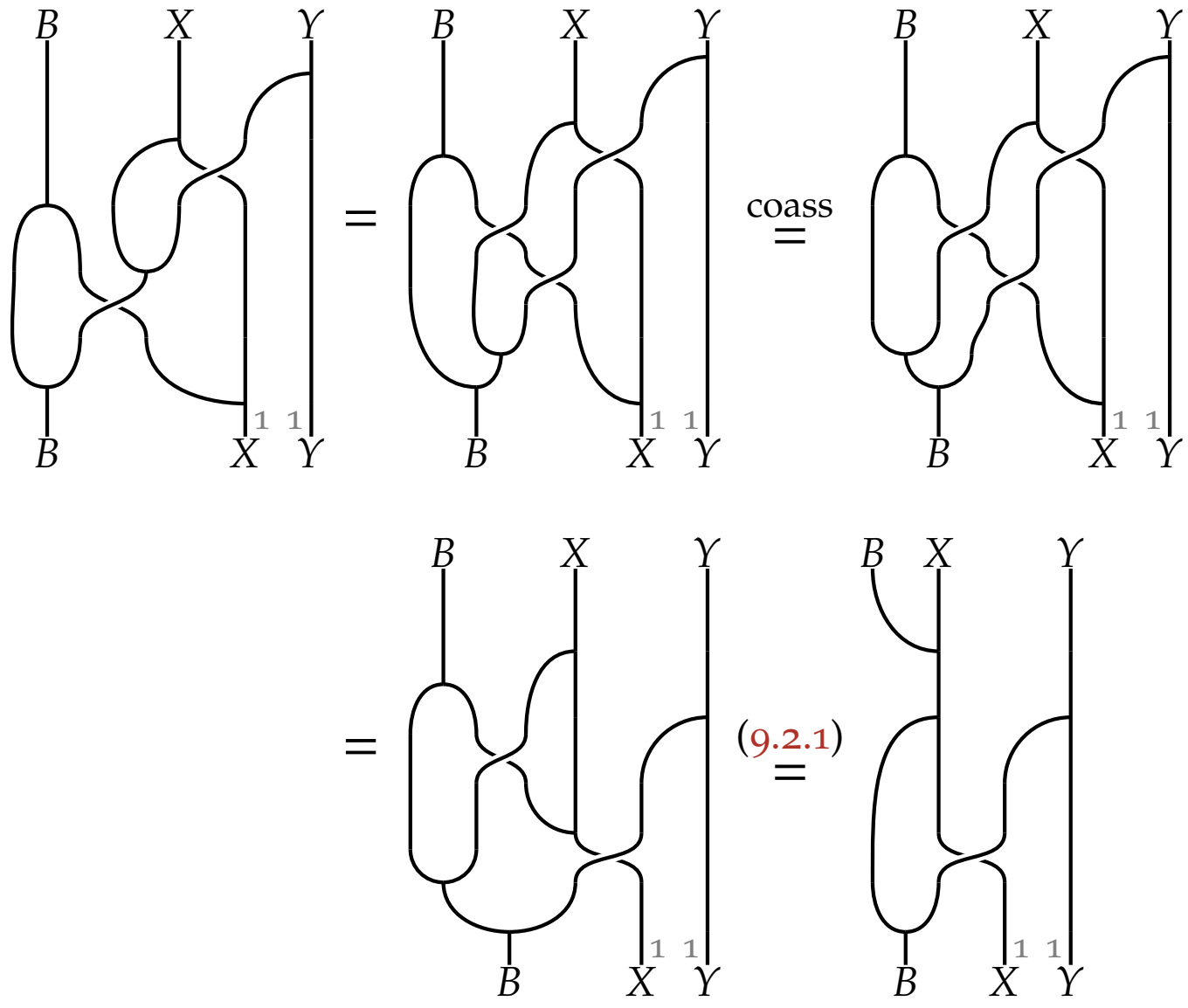

Figure 9.4: The action of the cotensor product is a right $B$-comodule morphism.

**Proposition 9.12.** *Let $\chi^X$ and $\chi^Y$ be as in Definition 9.6. Then*

$$\chi^X \diamond \chi^Y \cong \chi^{X \square Y},$$

*where $\diamond$ is the multiplicative cell for lax module morphisms; i.e.,*

$$(X \square -, \chi^X) \diamond (Y \square -, \chi^Y) := (X \square Y \square -, \chi^X \diamond \chi^Y).$$

*In other words, the functor $\Sigma$ from Equation (9.2.4) is strong monoidal.*

*Proof.* Associativity follows from coassociativity of the coaction and naturality of the braiding:

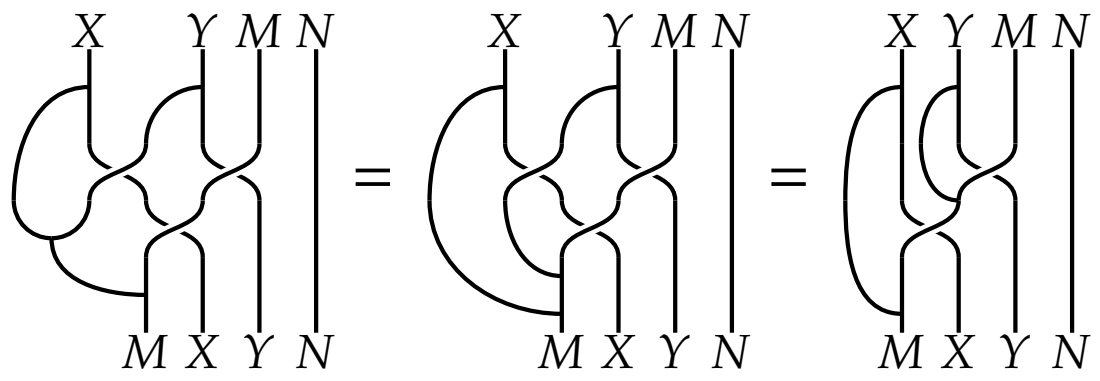

Unitality is immediate. □

To complete the proof of Theorem 9.2, it is left to show that the strong monoidal functor $\Sigma$ from Equation (9.2.4) is an equivalence. We will show that it is fully faithful and essentially surjective.

**Proposition 9.13.** *The functor $\Sigma$ is fully faithful.*

*Proof.* To show that $\Sigma$ is faithful, suppose that $f, g \colon X \longrightarrow Y$ are morphisms of Hopf trimodules, such that $f \square - \cong g \square -$. Then

$$\begin{array}{ccc} X \square B & \overset{f \square B}{\underset{g \square B}{\rightrightarrows}} & Y \square B \\ \cong \downarrow & & \downarrow \cong \\ X & \overset{f}{\underset{g}{\rightrightarrows}} & Y \end{array}$$

commutes, so the result follows.

To see that $\Sigma$ is full, suppose that $\varphi \colon X \square - \Longrightarrow Y \square -$ is a $\mathcal{V}$-module transformation. By Proposition 2.109, we to show that the induced arrow

$$\hat{\varphi}_B \colon X \cong X \square B \xrightarrow{\varphi_B} Y \square B \cong Y$$

is a morphism of Hopf trimodules. For a Hopf trimodule $M \in {}^{B}_{B}\mathsf{Vect}^{B}$, the morphism $\varphi_M : X \,\square\, M \longrightarrow Y \,\square\, M$ and the induced morphism $\hat{\varphi}_B$ can, in string diagrams, be depicted like this:

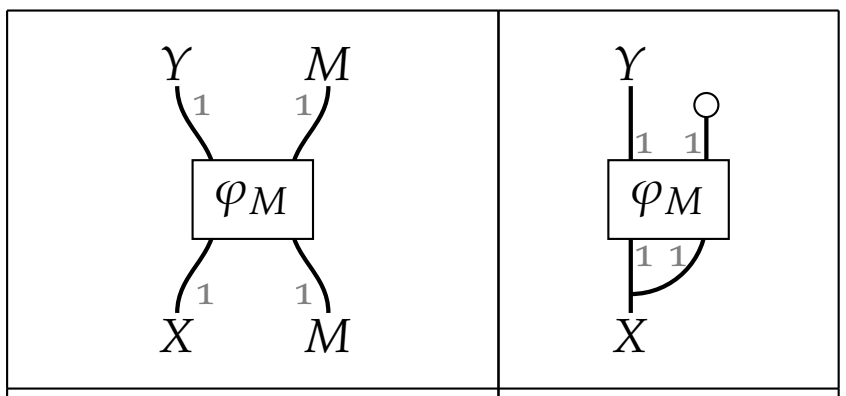

$\varphi : X \,\square\, - \longrightarrow Y \,\square\, -$ | $\hat{\varphi}_B : X \longrightarrow Y$

Figure 9.5 collects several properties of $\chi$ and $\varphi$ in this notation. The arrow $\hat{\varphi}_B$ is a left and right comodule morphism by the calculations in Figure 9.6. Likewise, that $\varphi_B$ is a left module morphism follows from:

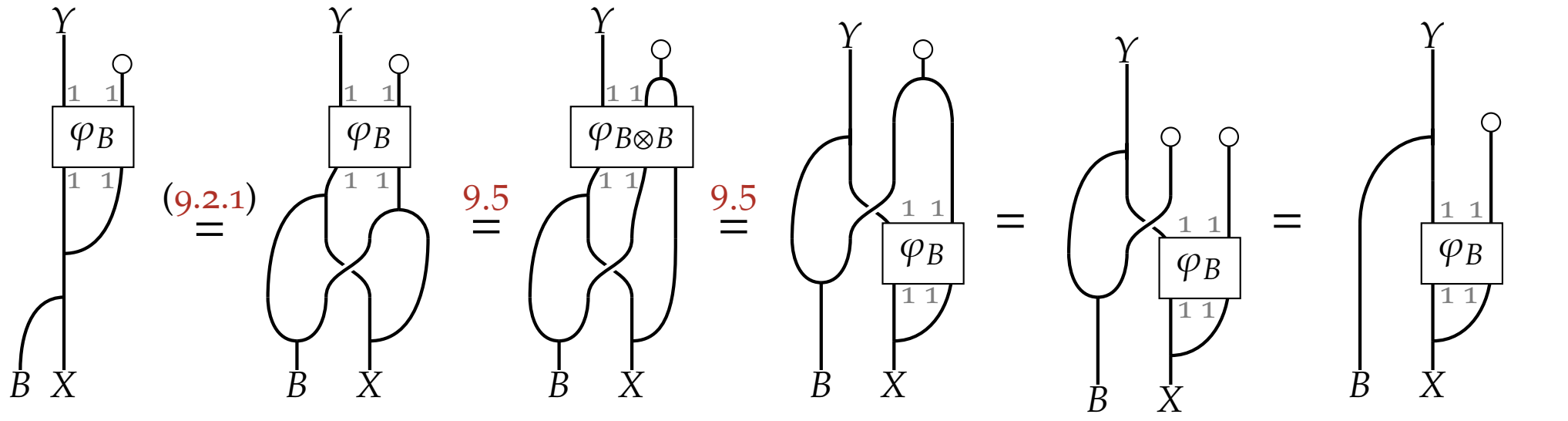

□

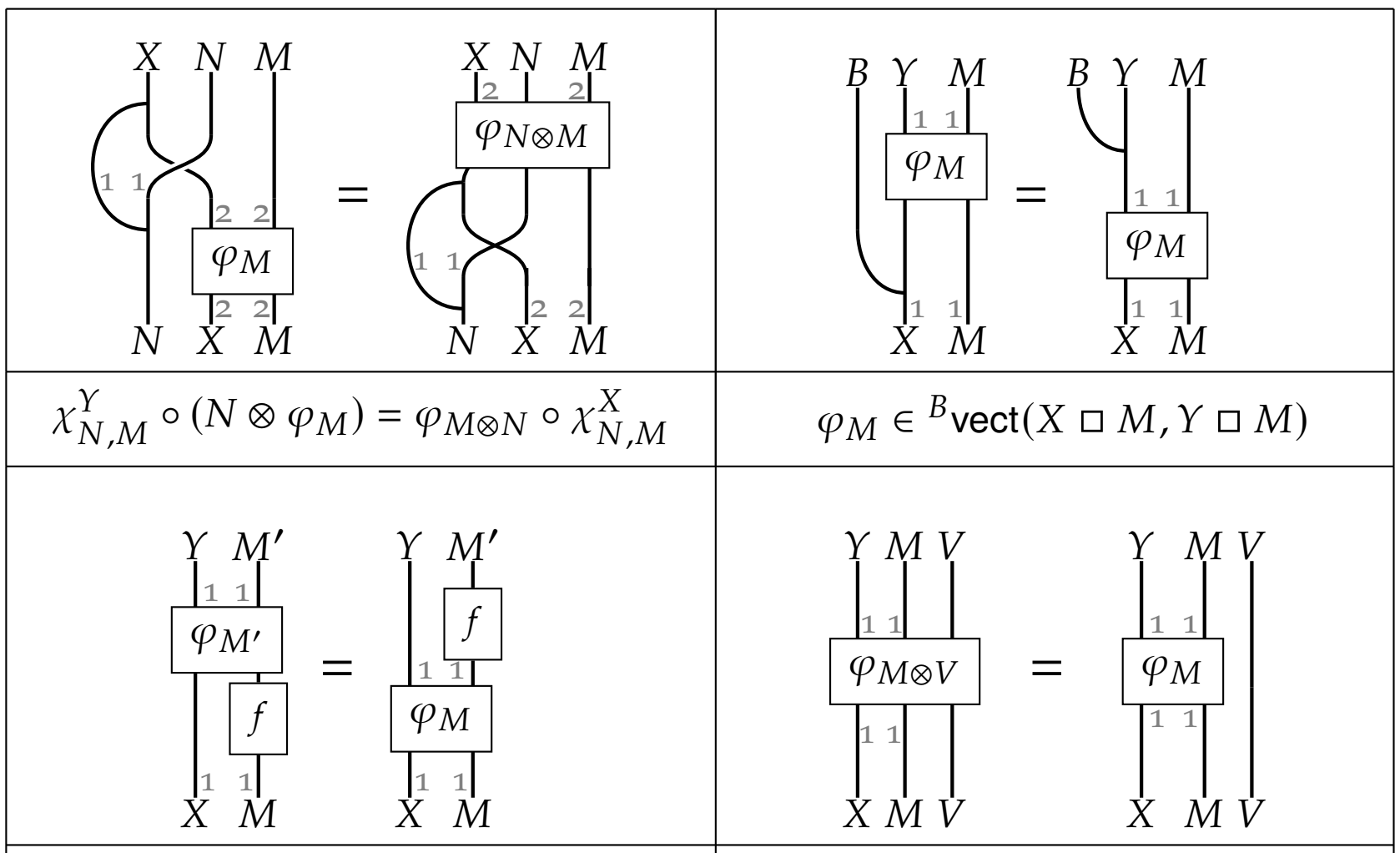

$\chi^{Y}_{N,M} \circ (N \otimes \varphi_M) = \varphi_{M \otimes N} \circ \chi^{X}_{N,M}$ | $\varphi_M \in {}^{B}\mathsf{vect}(X \,\square\, M, Y \,\square\, M)$

$(Y \otimes f) \circ \varphi_M = \varphi_{M'} \circ (X \otimes f)$ | $\forall_{V \in \mathsf{vect}}\ \varphi_{M \otimes V} = \varphi_M \otimes V$

Figure 9.5: Properties of the arrows $\chi$ and $\varphi$.

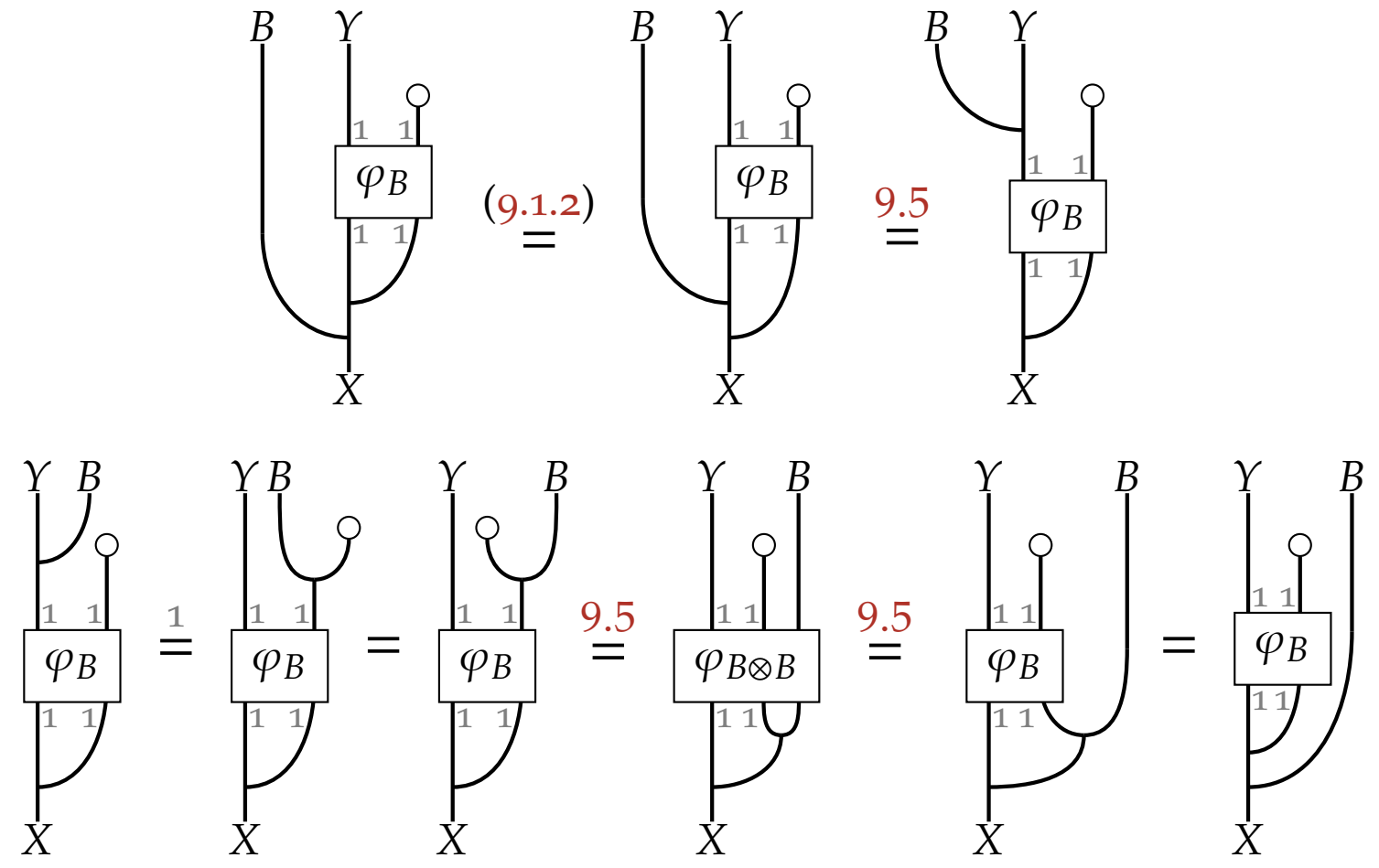

Figure 9.6: Verification that $\hat{\varphi}_B$ is a left and right comodule morphism.

**Proposition 9.14.** *The functor $\Sigma$ is essentially surjective.*

*Proof.* Suppose that $\Phi \in \mathsf{LexfLax}\mathcal{V}\mathsf{Mod}(\mathcal{V}, \mathcal{V})$. Since $\Phi$ is left exact and finitary, by Proposition 2.109 there exists a $B$-$B$-bicomodule $X$, such that $\Phi \cong X \,\square\, -$ as functors. We will show that if $X \,\square\, -$ is a lax module morphism, then $X$ comes equipped with a left module structure, making it into a Hopf trimodule.

Define a $B$-action on $X$ in the following way:

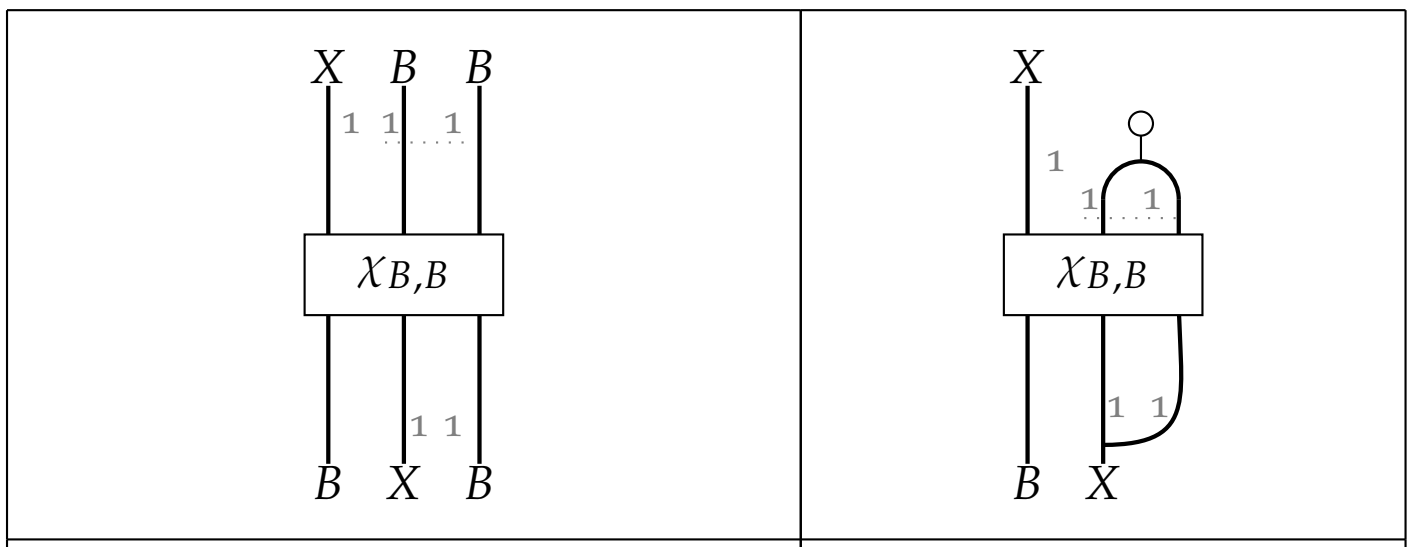

$\chi_{B,B} \colon B \otimes (X \,\square\, B) \longrightarrow X \,\square\, (B \otimes B)$ | $\alpha \colon B \otimes X \longrightarrow X$

Above, we have chosen to represent $\chi_{B,B}$—see Definition 9.6—by a single box. The claim is that this action turns $X$ into a Hopf trimodule.

For brevity, we shall omit the grey action indicators in our string diagrammatic proofs, save when performing the respective transformation.

To show that $\alpha$ is a bicomodule morphism, it suffices to note that the left comodule morphism $\chi_{B,B}$ is also a right comodule morphism, as $\alpha$ is then composed of bicomodule morphisms.

Notice that $\chi$ is a natural transformation of functors which are left exact finitary in each variable, from ${}^{B}\mathsf{Vect} \otimes_{\Bbbk} {}^{B}\mathsf{Vect}$ to ${}^{B}\mathsf{Vect}$. Since $\chi_{B,B}$ is the component of the transformation

$$\chi_{B,-} \colon B \otimes_{\Bbbk} (X \,\square\, -) \Longrightarrow X \,\square\, (B \otimes_{\Bbbk} -)$$

taken at the injective generator of ${}^{B}\mathsf{Vect}$, it follows that it is a right comodule morphism; see Propositions 2.84 and 2.109. Similarly, $\chi_{B,B}$ is the component of $\chi_{-,B}\colon - \otimes_{\Bbbk} (X \mathbin{\square} B) \Longrightarrow X \mathbin{\square} (- \otimes_{\Bbbk} B)$ taken at the injective generator, so $\chi_{B,B}$ also intertwines this additional comodule structure:

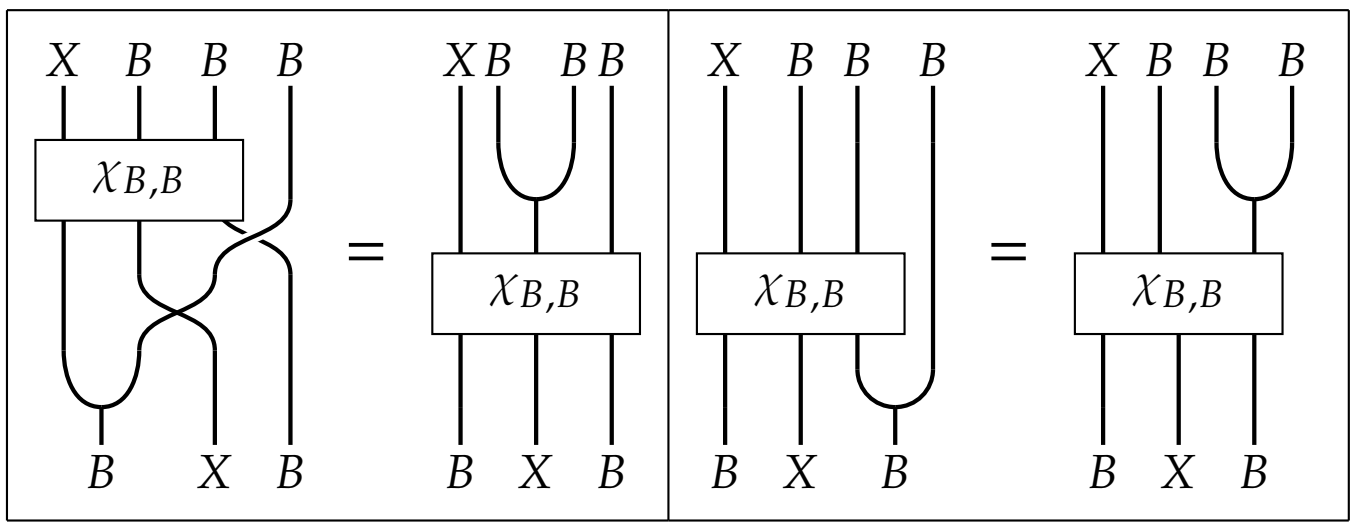

Let us now verify that $\alpha$ is in fact a $B$-module morphism. For the multiplicative identity of $\alpha$, one calculates

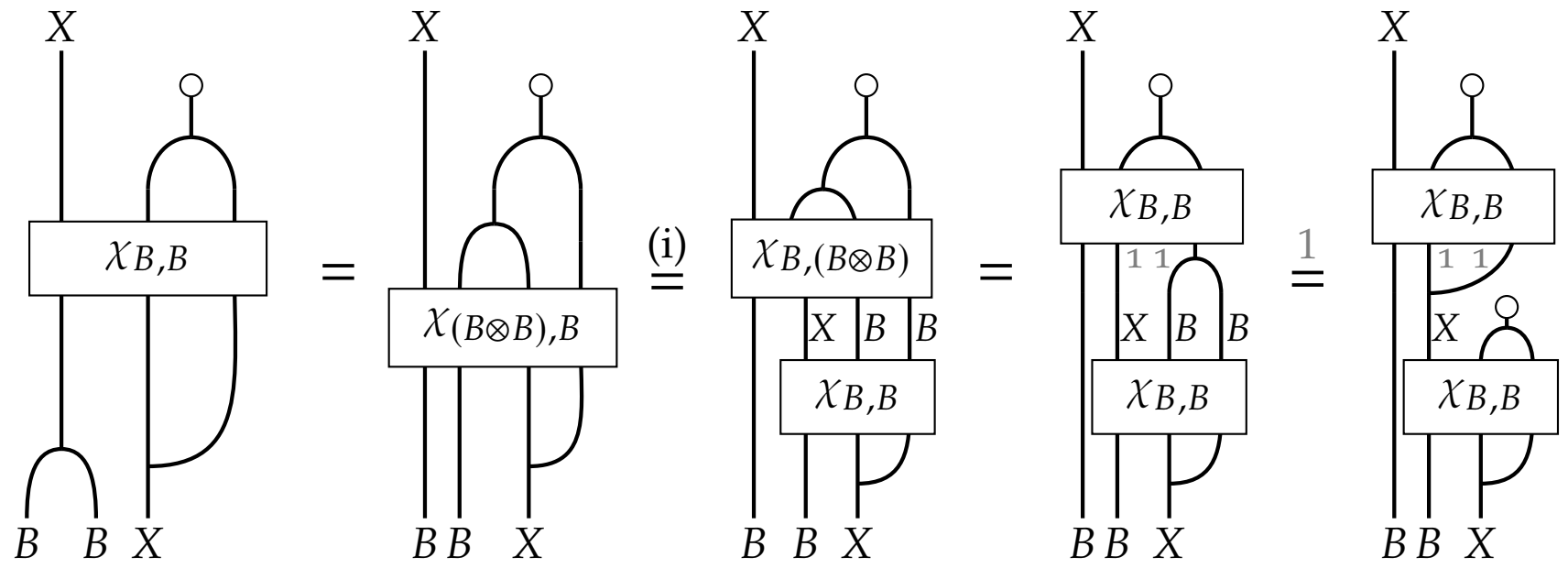

where (i) uses the fact that $(X \mathbin{\square} -, \chi)$ is a ${}^{B}\mathsf{Vect}$-module functor.

The unitality of $\alpha$ follows from the naturality of $\chi$, the unital axiom of $(X \mathbin{\square} -, \chi)$ being a ${}^{B}\mathsf{Vect}$-module functor, and unitality of $B$ and $X$:

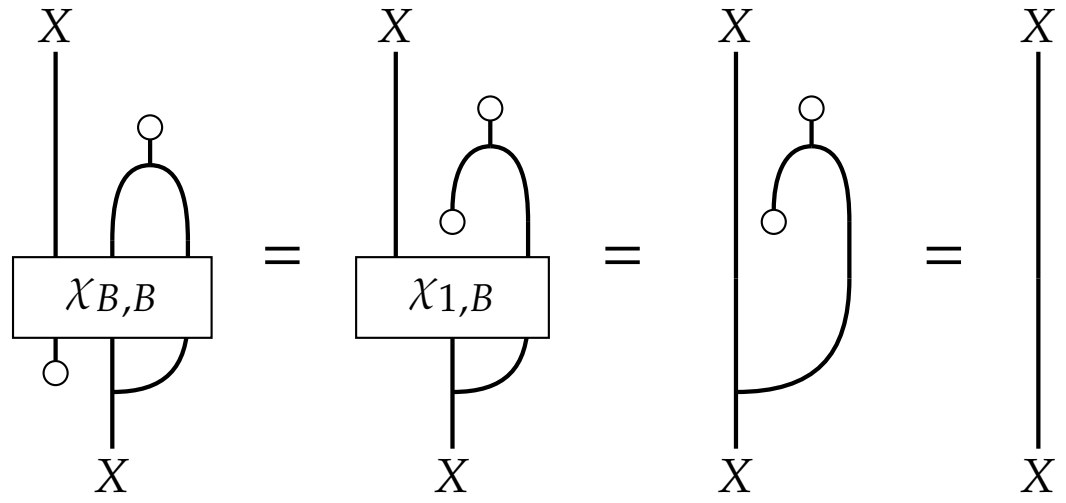

It is left to show that $\Sigma(\alpha) = \chi_{B,B}$. The necessary calculation is given in Figure 9.7, where equation (i) holds because $\chi_{B,B}$ is a right comodule morphism in both variables, and (ii) follows by naturality of $\chi$. $\square$

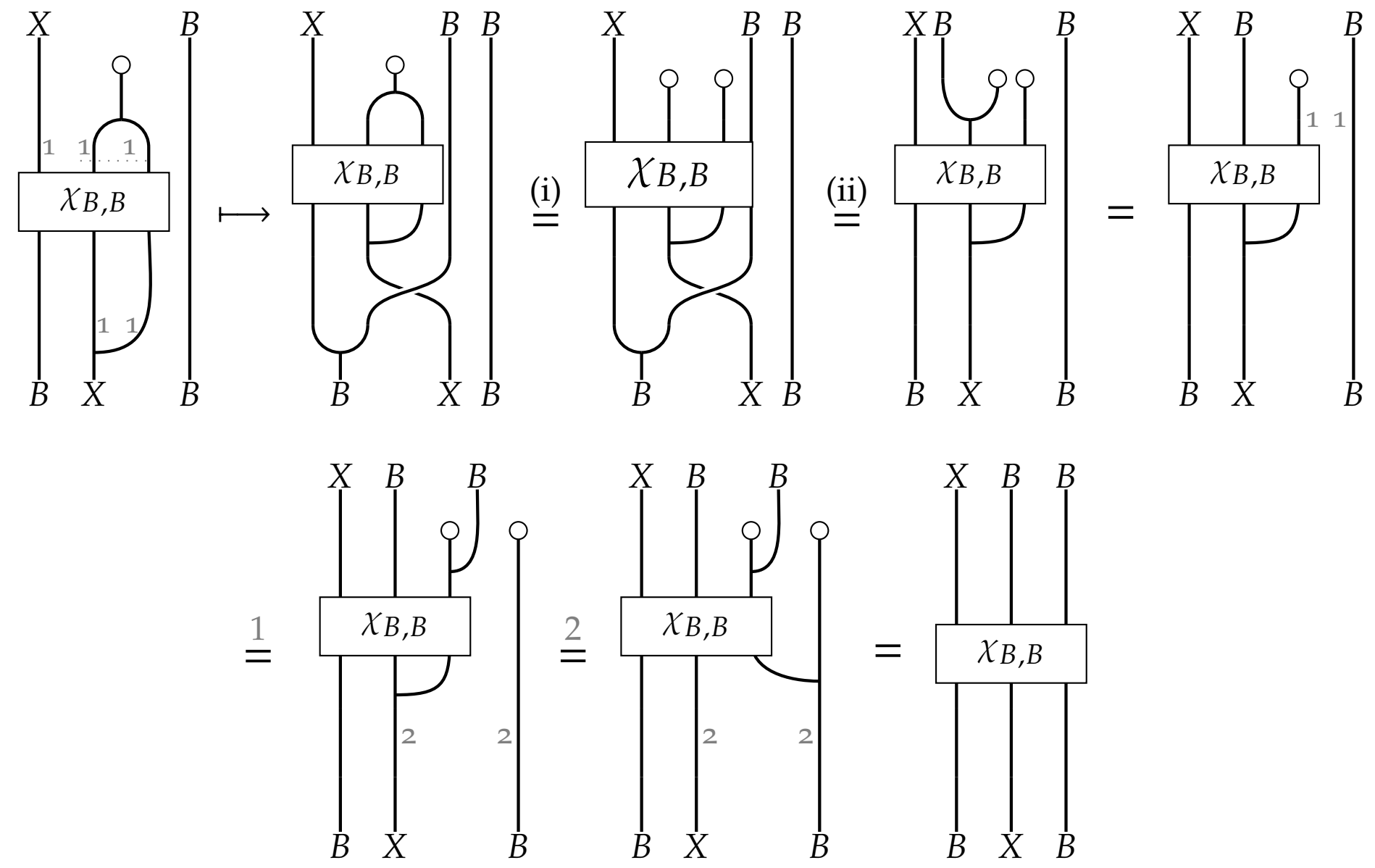

Figure 9.7: Verification that $\Sigma(\alpha) = \chi_{B,B}$.

## 9.3 Contramodule reconstruction over Hopf trimodule algebras

We stress that we are still working under Hypotheses 8.1 and 9.1.

**Definition 9.15.** A *Hopf trimodule algebra* consists of an algebra object in the monoidal category $\left({}^B_B\mathsf{Vect}^B, \square, B\right)$. More explicitly, it consists of a Hopf trimodule $A \in {}^B_B\mathsf{Vect}^B$, together with trimodule morphisms $\mu\colon A \square A \longrightarrow A$ and $\eta\colon B \longrightarrow A$, satisfying associativity and unitality conditions.

**Definition 9.16.** A *left module* over a Hopf trimodule algebra is a left $A$-module object in the ${}^B_B\mathsf{Vect}^B$-module category ${}^B\mathsf{Vect}$ in the sense of Definition 2.98. As such, it consists of a left $B$-comodule $M$ together with a $B$-comodule homomorphism $A \square M \longrightarrow M$, satisfying associativity and unitality conditions with respect to $\mu$ and $\eta$.

Similarly to Definition 9.16 we define right modules over $A$. A module over $A$ is *free* if it is of the form $(A \square N, \mu \square N)$.

We can also dualise this definition: we refer to objects of the category ${}^B_B\mathsf{Vect}_B$ as *Hopf termodules*. A *Hopf termodule coalgebra* is a coalgebra object $C$ in the monoidal category $\left({}^B_B\mathsf{Vect}_B, \otimes_B, B\right)$, and a *comodule over* $C$ is a left $B$-module $N$ together with a $B$-module homomorphism $N \longrightarrow C \otimes_B N$, satisfying usual coaction axioms. Right comodules over $C$ are defined analogously.

**Remark 9.17.** Under the equivalence of Equation (9.2.4) between ${}^B_B\mathsf{Vect}^B$ and $\mathsf{LexfLax}\mathcal{V}\mathsf{Mod}(\mathcal{V}, \mathcal{V})$, the ${}^B_B\mathsf{Vect}^B$-module category $\mathcal{V} = {}^B\mathsf{Vect}$ corresponds to the evaluation action of $\mathsf{LexfLax}\mathcal{V}\mathsf{Mod}(\mathcal{V},\mathcal{V})$ on $\mathcal{V}$. Thus, the category of modules over a Hopf trimodule algebra $A$ is precisely the Eilenberg–Moore category for the monad $A \,\square\, -$. The subcategory of free modules over $A$ is the image of the Kleisli category for this monad under the canonical embedding

$$({}^{B}\mathsf{Vect})_{A\square-} \xrightarrow{\iota} ({}^{B}\mathsf{Vect})^{A\square-}.$$

Let $A \in {}^{B}_{B}\mathsf{Vect}^{B}$ be a Hopf trimodule algebra such that ${}^{B}A$ is quasi-finite, meaning that the corresponding lax ${}^{B}\mathsf{Vect}$-module functor $A \,\square\, -$ has a left adjoint $\mathrm{cohom}_{B-}(A,-)$, which is canonically a comonad by Proposition 2.26.

**Definition 9.18.** A *left contramodule over* $A$ is an object in $({}^{B}\mathsf{Vect})^{\mathrm{cohom}_{B-}(A,-)}$.

Explicitly, this means that a left contramodule is a left $B$-comodule $M$ together with a left $B$-comodule homomorphism $M \longrightarrow \mathrm{cohom}_{B-}(A,M)$, endowing $M$ with the structure of a comodule in ${}^{B}\mathsf{Vect}$ over $\mathrm{cohom}_{B-}(A,-)$.

A contramodule is *free* if it lies in the image of the canonical embedding

$$({}^{B}\mathsf{Vect})_{\mathrm{cohom}_{B-}(A,-)} \xhookrightarrow{\iota} ({}^{B}\mathsf{Vect})^{\mathrm{cohom}_{B-}(A,-)}.$$

In other words, it is of the form $\mathrm{cohom}_{B-}(A,M)$, for a left $B$-comodule $M$, and its coaction is of the form

$$\mathrm{cohom}_{B-}(A,M) \xrightarrow{\mu^*} \mathrm{cohom}_{B-}(A \,\square\, A, M) \cong \mathrm{cohom}_{B-}(A, \mathrm{cohom}_{B-}(A,M)).$$

We denote the category of free contramodules by $A$-ContramodFree.

Similarly, the category $C$-Contramod of left contramodules over a Hopf termodule coalgebra $C$ is the Eilenberg–Moore category $({}_{B}\mathsf{Vect})^{\mathrm{Hom}_{B-}(C,-)}$ for the monad $\mathrm{Hom}_{B-}(C,-)$ right adjoint to the comonad $C \otimes_B -$. The free contramodules are defined analogously to the case of Hopf trimodules.

**Remark 9.19.** Contramodules similar to those of Definition 9.18 were studied extensively in more homological settings, see in particular [Pos10]. Forgetting the left $B$-module structure, a Hopf trimodule algebra yields a so-called *semialgebra*, and by forgetting the left $B$-comodule structure, a Hopf termodule coalgebra yields a so-called *bocs*—bimodule object, coalgebra structure.

**Lemma 9.20.** *Let $A$ be a Hopf trimodule algebra. The action*

$$V \triangleright (A \,\square\, M) := A \,\square\, (V \otimes M), \qquad \textit{for all } V \in \mathsf{Vect} \textit{ and } M \in {}^{B}\mathsf{vect}$$

*endows the category* $A$-ModFree *of free* $A$*-modules with a canonical module category structure over* ${}^B\mathsf{Vect}$*. Similarly,* $A$-ContramodFree *is a module category over* ${}^B\mathsf{Vect}$ *by means of the action* $V \triangleright \mathrm{cohom}_{B-}(A, N) := \mathrm{cohom}_{B-}(A, V \otimes N)$*. These two* ${}^B\mathsf{Vect}$*-module categories are canonically equivalent.*

*Proof.* By definition, we have $A\text{-ModFree} = ({}^B\mathsf{Vect})_{A\square -}$. Since $A$ is a Hopf trimodule algebra, the monad $A \,\square\, -$ is a lax ${}^B\mathsf{Vect}$-module monad, so the Kleisli category $({}^B\mathsf{Vect})_{A\square -}$ is a ${}^B\mathsf{Vect}$-module category, as described in Corollary 5.33, and it is easy to check that the assignments described in the lemma correspond precisely to those of *ibid*.

Similar considerations apply to $A$-ContramodFree, using the equality

$$A\text{-ContramodFree} = ({}^B\mathsf{Vect})_{\mathrm{cohom}_{B-}(A,-)}$$

and the fact that, since $\mathrm{cohom}_{B-}(A, -)$ is the left adjoint of the lax ${}^B\mathsf{Vect}$-module monad $A \,\square\, -$, it is an oplax ${}^B\mathsf{Vect}$-module comonad, by Porism 5.29.

The asserted equivalence between these two ${}^B\mathsf{Vect}$-module categories follows directly from Proposition 5.37. □

Recall from [Tak77, 1.12] that $\mathrm{Cohom}_{B-}(A, -)$ is exact if and only if $A$ is an injective left $B$-comodule.

**Proposition 9.21.** *If* ${}^BA$ *is injective as a left* $B$*-module, then the* ${}^B\mathsf{Vect}$*-module category structure on* $A$-ContramodFree *of Lemma 9.20 extends to a* ${}^B\mathsf{Vect}$*-module category structure on* $A$-Contramod*, in the sense of Theorem 8.25.*

*Proof.* For the comonadic dual of Theorem 8.25, it suffices to notice that the comonad $\mathrm{Cohom}_{B-}(A, -)$ is left exact, which is the case if and only if ${}^BA$ is injective, by [Tak77, 1.12]. □

**Theorem 9.22.** *Let* $B$ *be a bialgebra and let* $\mathcal{M}$ *be a locally finite abelian module category over* ${}^B\mathsf{vect}$*, which admits a* ${}^B\mathsf{Vect}$*-injective* ${}^B\mathsf{Vect}$*-cogenerator* $X \in \mathsf{Ind}(\mathcal{M})$*, such that* $- \triangleright X \colon {}^B\mathsf{vect} \longrightarrow \mathsf{Ind}(\mathcal{M})$ *is exact and quasi-finite. Then there is a Hopf trimodule algebra* $(A, \mu, \eta)$*, such that* ${}^BA$ *is quasi-finite and injective, hence finitely cogenerated injective, such that there is an equivalence of* ${}^B\mathsf{Vect}$*-module categories*

$$\mathsf{Ind}(\mathcal{M}) \simeq A\text{-Contramod},$$

*between* $\mathsf{Ind}(\mathcal{M})$ *and the category of left* $A$*-contramodules. The* ${}^B\mathsf{Vect}$*-module structure on* $A$-Contramod *is extended, in the sense of Definition 8.10, from the category* $A$-ContramodFree *of free left* $A$*-contramodules.*

*This equivalence restricts to an equivalence* $\mathcal{M} \simeq A\text{-contramod}$*, between* $\mathcal{M}$ *and the category of finite-dimensional* $A$*-contramodules.*

*Proof.* Let us again write $\mathcal{V} := {}^{B}\mathsf{Vect}$. Since $- \triangleright X \colon {}^{B}\mathsf{vect} \longrightarrow \mathsf{Ind}(\mathcal{M})$ is right exact, its extension to $\mathsf{Ind}({}^{B}\mathsf{vect}) \cong \mathcal{V}$ admits a right adjoint $\lfloor X, - \rfloor$ by Proposition 2.133. Similarly, since $- \triangleright X \colon {}^{B}\mathsf{vect} \longrightarrow \mathsf{Ind}(\mathcal{M})$ is left exact, so is its extension to $\mathcal{V}$, which is also quasi-finite by assumption. Thus, by Proposition 2.134, the functor $- \triangleright X \colon \mathcal{V} \longrightarrow \mathsf{Ind}(\mathcal{M})$ also has a left adjoint $\lceil X, - \rceil$. We thus obtain an oplax $\mathcal{V}$-module comonad $\lceil X, - \triangleright X \rceil$ on $\mathcal{V}$, and a right adjoint lax $\mathcal{V}$-module monad $\lfloor X, - \triangleright X \rfloor$.

By Theorem 9.2, we have an isomorphism $\lfloor X, - \triangleright X \rfloor \cong \lfloor X, B \triangleright X \rfloor \square -$ of lax $\mathcal{V}$-module monads, where $\lfloor X, B \triangleright X \rfloor$ is a Hopf trimodule algebra. Let us denote $\lfloor X, B \triangleright X \rfloor$ by $A$. Since $A \square -$ admits a left adjoint, ${}^{B}A$ is quasi-finite.

By uniqueness of adjoints, an isomorphism $\lceil X, - \triangleright X \rceil \cong \mathrm{cohom}_{B-}(A, -)$ follows. This is an isomorphism of $\mathcal{V}$-module comonads, since the $\mathcal{V}$-module comonad structure is lifted from the isomorphic oplax $\mathcal{V}$-module monad structures on the respective right adjoints.

By Theorem 8.50, we have an equivalence $\mathsf{Ind}(\mathcal{M}) \simeq \mathcal{V}^{\lceil X, - \triangleright X \rceil}$ of $\mathcal{V}$-module categories. In particular, the functor $\lceil X, - \rceil$ is comonadic and hence exact, see Theorem 2.92. By assumption, so is its right adjoint $- \triangleright X$, whence the comonad $\lceil X, - \triangleright X \rceil$ is exact, which, by [Tak77, 1.12] and since it is isomorphic to $\mathrm{cohom}_{B-}(A, -)$, implies the injectivity of ${}^{B}A$.

The $\mathcal{V}$-module structure on $\mathcal{V}^{\lceil X, - \triangleright X \rceil}$ is extended from $\mathcal{V}_{\lceil X, - \triangleright X \rceil}$. By the isomorphism $\lceil X, - \triangleright X \rceil \cong \mathrm{cohom}_{B-}(A, -)$ of the involved comonads, we have an equivalence of $\mathcal{V}$-module categories

$$\mathcal{V}^{\lceil X, - \triangleright X \rceil} \simeq \mathcal{V}^{\mathrm{cohom}(A, -)},$$

and an equivalence $\mathcal{V}_{\lceil X, - \triangleright X \rceil} \simeq \mathcal{V}_{\mathrm{cohom}(A, -)}$. Further, by definition

$$\mathcal{V}^{\mathrm{cohom}(A, -)} = A\text{-Contramod} \qquad \text{and} \qquad \mathcal{V}_{\mathrm{cohom}(A, -)} = A\text{-ContramodFree},$$

which yields the desired equivalence $\mathsf{Ind}(\mathcal{M}) \simeq A$-Contramod.

The restricted equivalence $\mathcal{M} \simeq A$-contramod follows from Proposition 8.3. In particular, we use the observation that a contramodule is compact if and only if its underlying $B$-comodule is such, if and only if it is finite-dimensional. □

In the finite-dimensional setting, a similar result uses projective objects.

**Theorem 9.23.** *Let $B$ be a finite-dimensional bialgebra and let $\mathcal{M}$ be a finite abelian module category over ${}_{B}\mathsf{vect}$, which admits an object $X \in \mathcal{M}$ such that $X$ is a ${}_{B}\mathsf{vect}$-projective ${}_{B}\mathsf{vect}$-generator; and $- \triangleright X \colon {}_{B}\mathsf{vect} \longrightarrow \mathcal{M}$ is exact.*

*Then there is a finite-dimensional Hopf termodule coalgebra* $C \in {}^{B}_{B}\mathsf{vect}_{B}$, *such that* ${}_{B}C$ *is projective and there is an equivalence of* ${}_{B}\mathsf{vect}$*-module categories*

$$\mathcal{M} \simeq C\text{-contramod}$$

*between* $\mathcal{M}$ *and the category of finite-dimensional left* $C$*-contramodules. The* ${}_{B}\mathsf{vect}$*-module structure on* $C$-contramod *is extended from the free left* $A$*-contramodules.*

*Proof.* This follows analogously to Theorem 9.22, using Theorem 8.48 instead of Theorem 8.50 and Proposition 8.5 instead of Proposition 8.3. □

In the semisimple case, contramodules become equivalent to modules.

**Definition 9.24.** We say that a Hopf trimodule algebra $A$ is *semisimple* if the monad $A \,\square\, -$ is semisimple in the sense of Definition 8.15.

**Proposition 9.25.** *For a semisimple Hopf trimodule algebra* $A \in {}^{B}_{B}\mathsf{Vect}^{B}$, *there is an equivalence* $A$-Mod $\simeq$ $A$-Contramod *of module categories over* ${}^{B}\mathsf{Vect}$.

*Proof.* This is Proposition 8.18 applied to the monad $A \,\square\, -$. □

### 9.3.1 *Morita equivalence for contramodules*

**Definition 9.26.** We say that two Hopf trimodule algebras $A, A' \in {}^{B}_{B}\mathsf{vect}^{B}$, such that ${}^{B}A$ and ${}^{B}A'$ are quasi-finite, are *contraMorita equivalent* if the comonads $\mathrm{cohom}_{B-}(A,-)$ and $\mathrm{cohom}_{B-}(A',-)$ are Morita equivalent.

More explicitly, Definition 9.26 says that $A$ and $A'$ are contraMorita equivalent if there are $M, N \in {}^{B}_{B}\mathsf{Vect}^{B}$ with coaction morphisms in ${}^{B}_{B}\mathsf{Vect}^{B}$

$$\begin{aligned} \mathrm{la}_{M,A}\colon M &\longrightarrow \mathrm{cohom}_{B-}(A,M), & \mathrm{ra}_{N,A}\colon N &\longrightarrow N \,\square_B\, \mathrm{cohom}_{B-}(A,B), \\ \mathrm{la}_{N,A'}\colon N &\longrightarrow \mathrm{cohom}_{B-}(A',N), & \mathrm{ra}_{M,A'}\colon M &\longrightarrow M \,\square_B\, \mathrm{cohom}_{B-}(A',B), \end{aligned}$$

such that

- there is an isomorphism in ${}^{B}_{B}\mathsf{Vect}^{B}$ between the equaliser of

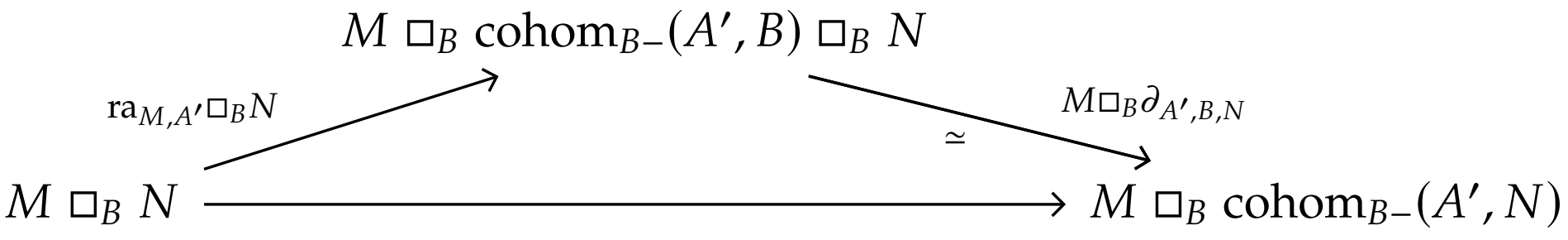

and $A$, intertwining the two coactions $A \longrightarrow A \,\square_B\, \mathrm{cohom}_{B\text{-}}(A, B)$ and $A \longrightarrow \mathrm{cohom}_{B\text{-}}(A, A)$ coming from the algebra structure on $A$, with the coactions on the equaliser induced by $\mathrm{ra}_{N,A}$ and $\mathrm{la}_{M,A}$. Here, $\partial_{A',B,N}$ is the isomorphism of [Tak77, 1.13], which is invertible because ${}^{B}A'$ is injective.

- a similar isomorphism in ${}^{B}_{B}\mathsf{Vect}^{B}$ between a similarly defined equaliser reversing the roles of $M$ and $N$, and $A'$ exists.

**Proposition 9.27.** *For $A, A'$ as in Definition 9.26, the following are equivalent:*

*(i)* *$A$ and $A'$ are contraMorita equivalent,*
*(ii)* $A$-Contramod $\simeq$ $A'$-Contramod *as* ${}^{B}\mathsf{Vect}$*-module categories, and*
*(iii)* $A$-contramod $\simeq$ $A'$-contramod *as* ${}^{B}\mathsf{vect}$*-module categories.*

*Proof.* The equivalence of (i) and (ii) is an instance of Proposition 8.74, for $\mathscr{C} = {}^{B}\mathsf{Vect}$. Since an equivalence of categories preserves compact objects, (iii) follows from (ii). Finally, if $F$ is an equivalence as in (iii), then $\mathsf{Ind}(F)$ is an equivalence as in (ii). □

Combining Proposition 9.27 with Theorem 9.22, we find the following algebraic realisation of Theorem 8.75:

**Theorem 9.28.** *There is a bijection*

$$\{\,(\mathcal{M}, X) \text{ as in Theorem 9.22}\,\}\big/_{\mathcal{M}\,\simeq\,\mathcal{N}} \overset{\cong}{\longleftrightarrow} \left\{\begin{matrix}\text{Left finitely cogenerated}\\ \text{injective Hopf trimodule}\\ \text{algebras over } B\end{matrix}\right\}\Big/ \simeq_{cMorita}.$$

## 9.4 A semisimple example of non-rigid reconstruction

We now give a concrete application of Theorems 8.50 and 9.22, describing module categories for non-rigid categories where the ordinary reconstruction procedure for rigid categories would not yield the correct result.

Let $S = \{e, s\}$ be the commutative two-element monoid which is not a group, with identity element $e$. In particular, $s^2 = s$ is an idempotent. The category $\mathsf{rep}(S)$ of finite-dimensional modules over the monoid algebra $\Bbbk[S]$ is semisimple, with two isoclasses of simple objects given by 1-dimensional representations $\Bbbk_{\mathrm{triv}}$ and $\Bbbk_{\mathrm{grp}}$. In $\Bbbk_{\mathrm{triv}}$, the element $s$ acts by 1, while in $\Bbbk_{\mathrm{grp}}$ it acts by 0. The standard comultiplication on $\Bbbk[S]$, in which $e$ and $s$ are group-like, makes turns $\Bbbk[S]$ into a bialgebra. The category $\mathsf{corep}(S)$ of finite-dimensional comodules over $\Bbbk[S]$ is equivalent to that of $S$-graded finite-dimensional spaces, $\mathsf{vect}^S$. Thus, it too is semisimple of rank two—we denote the simple objects by $\delta_e$ and $\delta_s$.

The monoidal structure on $\mathsf{rep}(S)$ satisfies $\Bbbk_{\mathrm{grp}} \otimes \Bbbk_{\mathrm{grp}} \cong \Bbbk_{\mathrm{grp}}$. In fact, this presentation of the Grothendieck ring $[\mathsf{rep}(S)]$ of $\mathsf{rep}(S)$ determines $\mathsf{rep}(S)$ uniquely as a semisimple monoidal category, see [SCZ23, Proposition 3.5]. Thus, we find monoidal equivalences

$$\mathsf{rep}(S) \simeq \mathsf{corep}(S) \simeq \mathsf{vect}^S,$$

since in the latter case we have $\delta_s \otimes \delta_s \cong \delta_s$.

Let us now focus on the coalgebraic setup of Theorem 8.50, studying the module categories over $\mathsf{corep}(S)$.

A semisimple $\mathsf{corep}(S)$-module category $(\mathcal{M}, \triangleright)$ is indecomposable[21] if and only if the matrix $[\delta_s]_{\mathcal{M}}$ describing the action of $\delta_s$ on the Grothendieck ring $[\mathcal{M}]$ in the basis of simple objects for $\mathcal{M}$ is indecomposable. But $[\delta_s]$ is idempotent, so it must be the $1 \times 1$-matrix $(1)$, and, as a category, we must have $\mathcal{M} \simeq \mathsf{vect}$. We then either have $\delta_s \triangleright - \cong \mathrm{Id}_{\mathcal{M}}$, or $\delta_s \triangleright - \cong 0$. One can show that any two module categories satisfying one of these conditions are equivalent, by showing that the second cohomology group $H^2(S, \Bbbk^\times)$ is trivial, similarly to [EGNO15, Example 7.4.10]. Indeed, $\psi \in Z^2(S, \Bbbk^\times)$ must satisfy the relation $\psi(e,s) = \psi(e,e) = \psi(s,e)$, and thus can be written as $\psi(x,y) = \sigma(x) - \sigma(xy) + \sigma(y)$, by setting $\sigma(e) := \psi(e,e)$ and $\sigma(s) := \psi(s,s)$. The condition $\delta_s \triangleright - \cong \mathrm{Id}_{\mathcal{M}}$ is satisfied by the fibre functor $\mathrm{Fr}\colon \mathsf{corep}(S) \longrightarrow \mathsf{vect}$.

[21] A $\mathscr{C}$-module category $\mathcal{M}$ is called *indecomposable* if it does not split as a direct sum of $\mathscr{C}$-module subcategories.

Assume instead that we have an indecomposable module category $\mathcal{M}$ satisfying $\delta_s \triangleright - \cong 0$. Let $X \in \mathcal{M}$ be a simple object. One can verify that $X$ satisfies the assumptions of Theorem 8.50. Given $V \in \mathsf{corep}(S)$, we have

$$\lceil X, V \triangleright X \rceil_s = \mathsf{corep}(S)(\lceil X, V \triangleright X \rceil, \delta_s) \cong \mathsf{corep}(S)(V \triangleright X, \delta_s \triangleright X),$$

and similarly for $\delta_e$.

Using $\delta_s \triangleright X \cong 0$ and $\delta_e \triangleright X \cong X$, we find that $\lceil X, \delta_e \triangleright X \rceil \cong \delta_e$ and $\lceil X, \delta_s \triangleright X \rceil = 0$. To describe the bicomodule corresponding to the functor $\lceil X, - \triangleright X \rceil$ under the correspondence of Proposition 2.109, we observe that a bicomodule can be identified with an $S \times S$-graded space, and the cotensor product is given by

$$(V \mathbin{\square} W)_{x,y} = \bigoplus_{z \in S} V_{x,z} \otimes W_{z,y}.$$

Using this, we find that $\lceil X, - \triangleright X \rceil \cong \delta_{e,e} \,\square\, -$, where $\delta_{e,e}$ is the one-dimensional space concentrated in degree $(e, e)$. The unit Hopf trimodule homomorphism $\Bbbk[S] \longrightarrow \lceil X, \Bbbk[S] \triangleright X \rceil$ must thus have the submodule $\Bbbk\{s\}$ as its kernel, so the reconstructing Hopf trimodule $\lceil X, \Bbbk[S] \triangleright X \rceil$ is given by $\Bbbk[S]/\Bbbk\{s\}$. As a bicomodule, it is indeed concentrated in degree $(e, e)$, and as a module it is isomorphic to $\Bbbk_{\mathrm{grp}}$. The composition

$$\lceil X, \delta_e \triangleright X \rceil \otimes \lceil X, \delta_e \triangleright X \rceil \longrightarrow \lceil X, \delta_e \triangleright X \rceil$$

sending $\mathrm{id}_X \otimes \mathrm{id}_X$ to $\mathrm{id}_X$ corresponds to the algebra structure

$$\begin{aligned} \Bbbk[S]/\Bbbk\{s\} \,\square\, \Bbbk[S]/\Bbbk\{s\} \simeq \Bbbk[S]/\Bbbk\{s\} \otimes \Bbbk[S]/\Bbbk\{s\} &\longrightarrow \Bbbk[S]/\Bbbk\{s\} \\ (e + \Bbbk\{s\}) \otimes (e + \Bbbk\{s\}) &\longmapsto (e + \Bbbk\{s\}), \end{aligned}$$

and this multiplication defines a Hopf trimodule algebra.

To check that the category $\Bbbk[S]/\Bbbk\{s\}$-Contramod satisfies the properties we have assumed of $\mathcal{M}$ previously, note that the free module $\Bbbk[S]/\Bbbk\{s\} \,\square\, B$ is one-dimensional, hence semisimple, showing that the category of modules over $\Bbbk[S]/\Bbbk\{s\}$ is semisimple of rank one. Thus, by Proposition 9.25 we have

$$\Bbbk[S]/\Bbbk\{s\}\text{-Contramod} \simeq \Bbbk[S]/\Bbbk\{s\}\text{-Mod}.$$

Further, $\delta_s \triangleright \Bbbk[S]/\Bbbk\{s\} = \Bbbk[S]/\Bbbk\{s\} \,\square\, \delta_s = 0$, which corresponds to $\delta_s \triangleright X = 0$.

If we instead follow the reconstruction procedure described in [EGNO15, Chapter 7], the reconstructing algebra object would be the object $A$ of $\mathsf{corep}(S)$ representing the presheaf $\mathcal{M}(- \triangleright X, X)$. We see that we must have $A \cong \delta_e$, and the algebra structure on $\delta_e$ is unique. However, $\delta_e \cong 1_{\mathsf{corep}(S)}$ and so

$$\mathrm{Mod}(\delta_e) \simeq \mathsf{corep}(S) \not\simeq \mathcal{M},$$

showing that the reconstruction procedure for module categories over rigid monoidal categories fails in this non-rigid case. It also illustrates that this failure can be observed already in the semisimple case.

## 9.5 A Hopf trimodule algebra constructing the fibre functor

This section discusses the general construction of a Hopf trimodule algebra, which covers the remaining indecomposable semisimple module category over the category $\mathsf{corep}(S)$ considered in Section 9.4.

**Proposition 9.29.** *The $\Bbbk$-vector space $B \otimes_{\Bbbk} B$ forms a Hopf trimodule with coactions*

$$b \otimes c \longmapsto b_{(1)} \otimes b_{(2)} \otimes c, \qquad c \otimes b \longmapsto c \otimes b_{(1)} \otimes b_{(2)},$$

*and action*

$$x \otimes b \otimes c \longmapsto x_{(1)} b \otimes x_{(2)} b.$$

*We shall denote it by $B \bullet B$.*

*Proof.* It is clear that $B \bullet B$ is a bicomodule, as both of the coactions are given by that of $B$ on the respective Let.

Thus, we first verify that the action is a morphism of left and right comodules; i.e., that $B \bullet B$ is in ${}_B\mathsf{Vect}^B$ and ${}_B^B\mathsf{Vect}$. The former follows by Figure 9.8, and the latter is due to the calculation in Figure 9.9. □

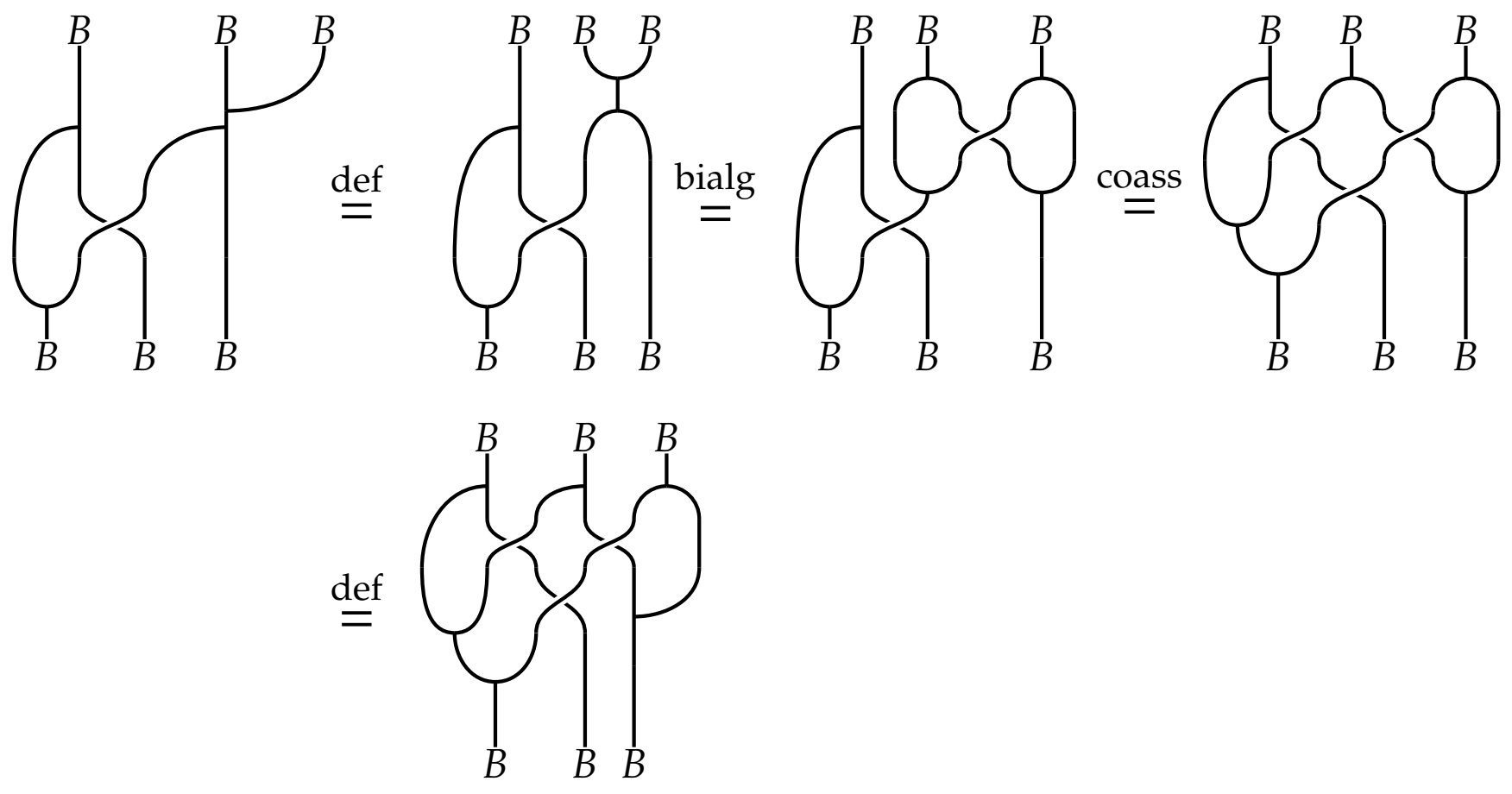

Figure 9.8: The action of $B \bullet B$ is a morphism of left comodules.

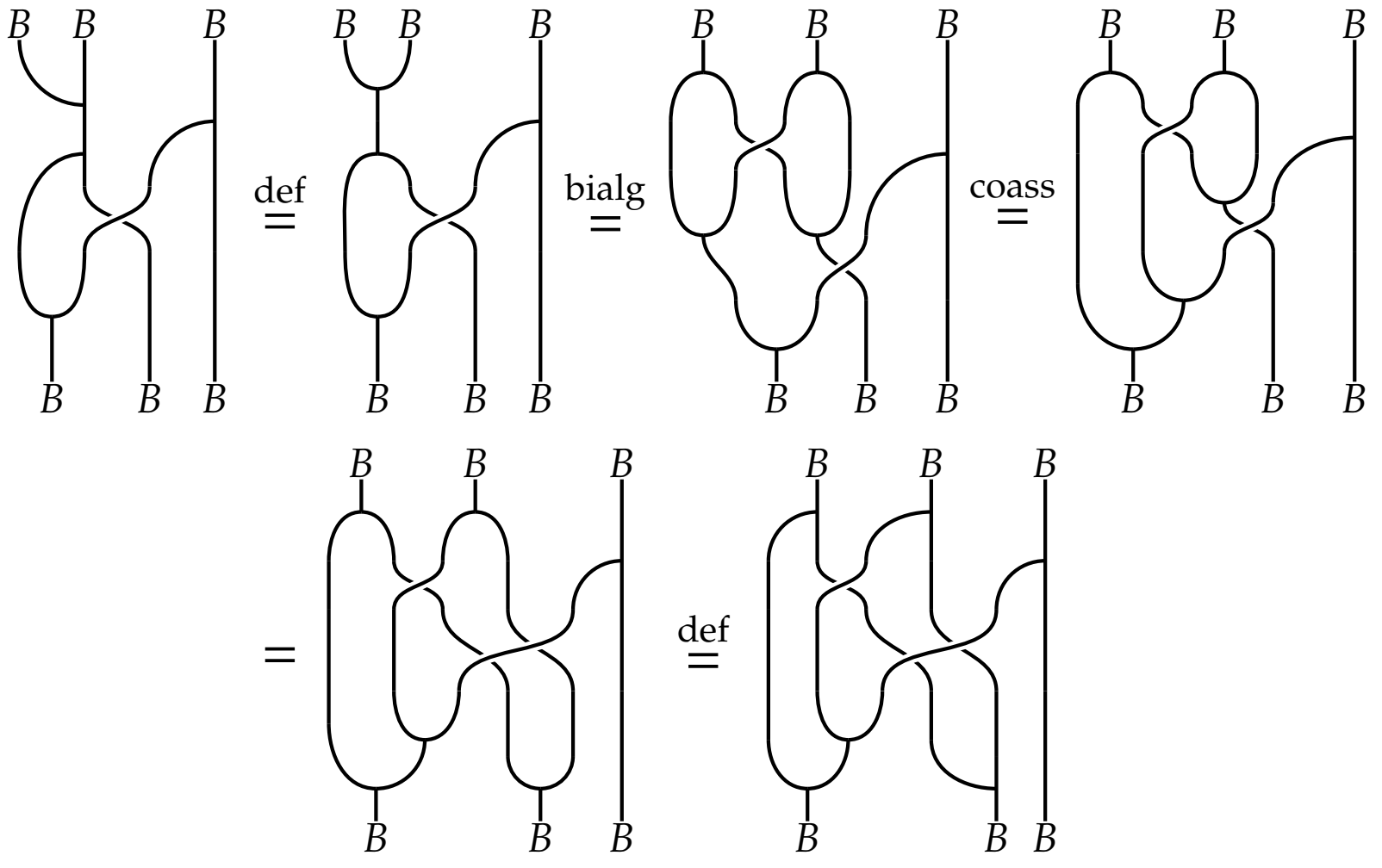

Figure 9.9: The action of $B \bullet B$ is a morphism of right comodules.

**Proposition 9.30.** *The Hopf trimodule $B \bullet B$ from Proposition 9.29 is an algebra object in ${}^{B}_{B}\mathsf{Vect}^{B}$.*

*Proof.* We define the multiplication and unit of $B \bullet B$ by

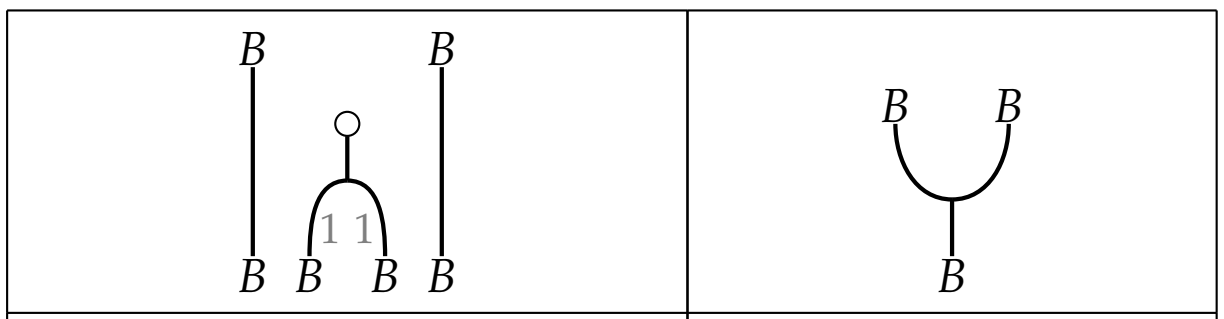

| $\mu\colon (B \bullet B) \,\square\, (B \bullet B) \longrightarrow B \bullet B$ | $\eta = \Delta\colon B \longrightarrow B \bullet B$ |
|---|---|

where we have denote the counit $\varepsilon$ of $B$ with a small white dot, and have employed the same numbering system for representing the cotensor product as in Section 9.1. It is clear that $\mu$ is associative, left unitality follows from and right unitality is similar.

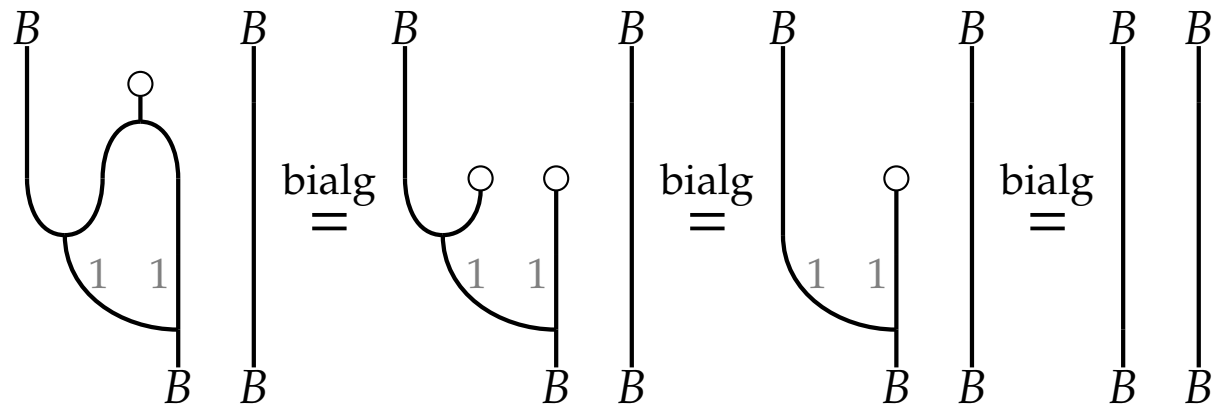

One immediately has that $\mu$ is a ${}^{B}\mathsf{Vect}^{B}$-morphism, and for $\eta = \Delta$ this follows from coassociativity. To finish the proof we need only show that both are morphisms of left $B$-modules, which follows from Figure 9.10 and the following calculation:

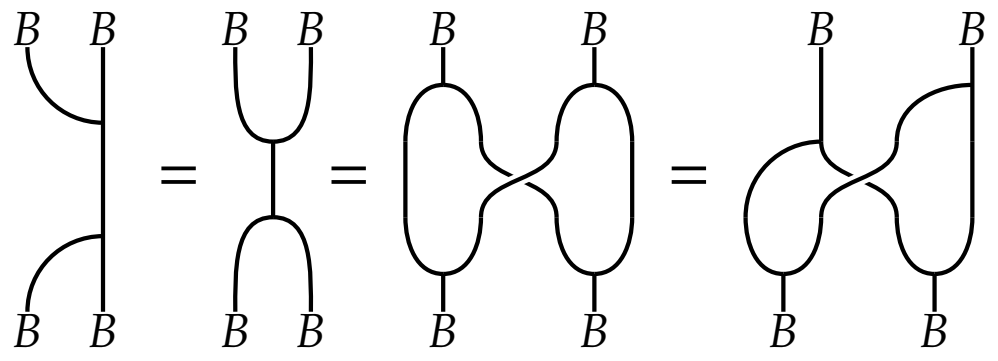

□

**Definition 9.31.** Let $J\colon (B \bullet B)\text{-ModFree} \longrightarrow \mathsf{Vect}$ be the functor given by $J((B \bullet B) \,\square\, M) = M$ and by local maps

$$J_{M,N}\colon (B \bullet B)\text{-ModFree}\big((B \bullet B) \,\square\, M, (B \bullet B) \,\square\, N\big) \longrightarrow \mathrm{Hom}_{\Bbbk}(M, N)$$
$$\sigma \longmapsto (\varepsilon \otimes \varepsilon \,\square\, \mathrm{id}_N) \circ (\sigma(1 \otimes 1 \otimes -))$$

**Lemma 9.32.** *Definition 9.31 yields an equivalence of ${}^{B}\mathsf{Vect}$-module categories.*

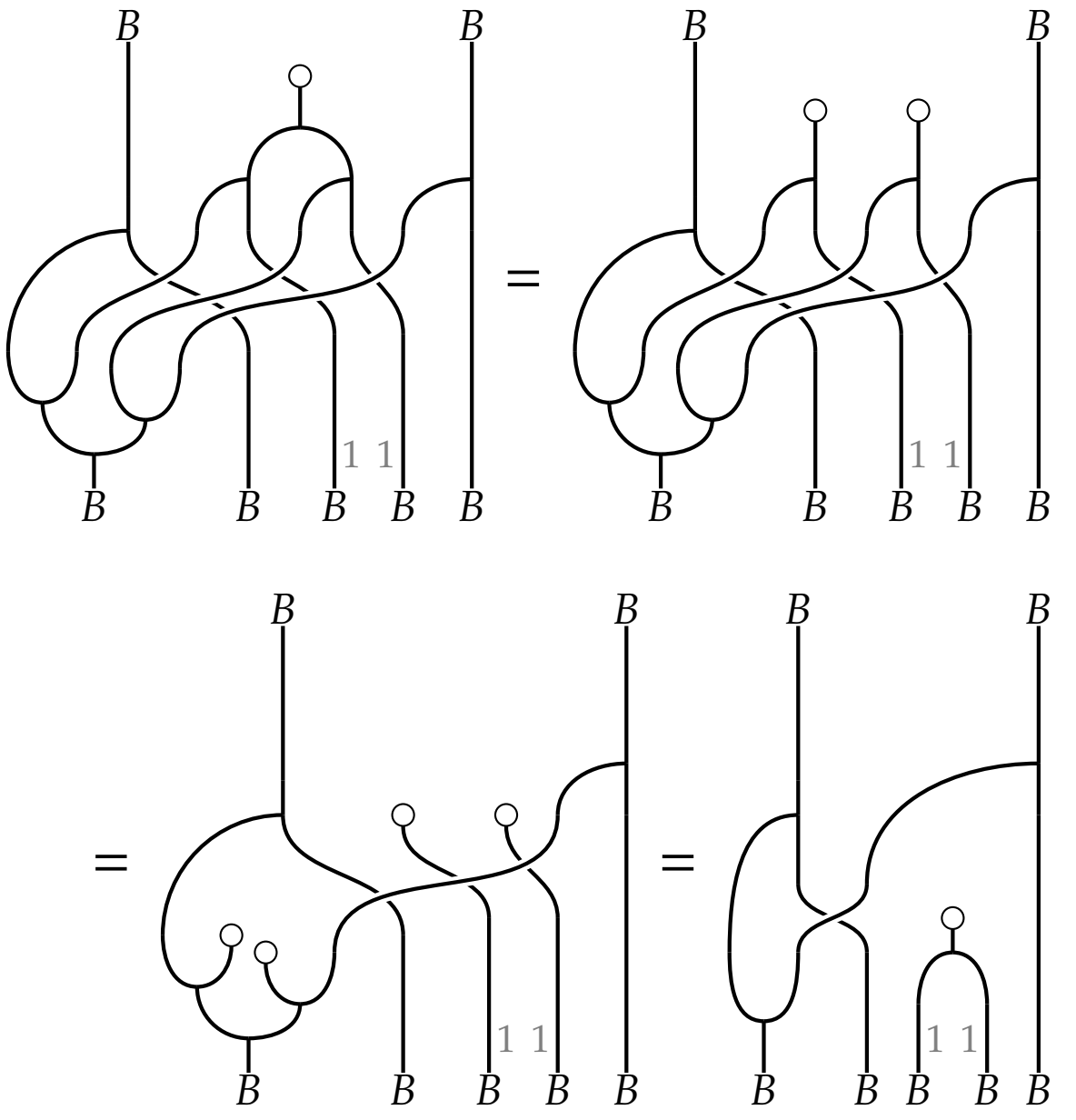

Figure 9.10: The unit $\eta = \Delta$is a morphism of left $B$-modules.

*Proof.* First, assume that $J$ is a functor; we will show it is an equivalence of $({}^{B}\mathsf{Vect})$-module categories. Let $V \in \mathsf{Vect}$. Essential surjectivity follows by

$$V \cong \Bbbk^{\dim V} \cong J((B \bullet B) \,\square\, \Bbbk_{\mathrm{triv}}^{\dim V}).$$

The fact that $J$ is fully faithful follows by $J_{M,N}$ being defined as the composite of the isomorphisms

$$(B \bullet B)\text{-ModFree}\big((B \bullet B) \,\square\, M, (B \bullet B) \,\square\, N\big) \xrightarrow{-\circ\eta} {}^{B}\mathsf{Vect}(M, (B \bullet B) \,\square\, N)$$
$$\xrightarrow{(B\otimes\varepsilon\square N)\circ-} {}^{B}\mathsf{Vect}(M, B \otimes N) \xrightarrow{(\varepsilon\otimes N)\circ-} \mathsf{vect}(M, N).$$

The first morphism is an isomorphism of the form $\mathcal{A}^{T}(TX, TY) \xrightarrow{\sim} \mathcal{A}(X, TY)$, for the monad $T := (B \bullet B) \,\square\, -$ on the category $\mathcal{A} := {}^{B}\mathsf{Vect}$, see Remark 8.4. The second isomorphism is simply postcomposition with the isomorphism $(B \bullet B) \,\square\, N = B \otimes B \,\square\, N \xrightarrow{B\otimes\varepsilon\square N} B \otimes N$, since $B \,\square\, N \xrightarrow{\varepsilon\square N} N$ is an isomorphism. The third is an isomorphism $\mathcal{B}^{C}(CX, CY) \xrightarrow{\sim} \mathcal{B}(CX, Y)$, for the comonad $K := B \otimes -$ on the category $\mathcal{B} := \mathsf{Vect}$.

Lastly, $J$ is a $({}^{B}\mathsf{Vect})$-module functor: for $W \in {}^{B}\mathsf{Vect}$, we have

$$W \triangleright J(B \bullet B \,\square\, M) = W \otimes M = J(B \bullet B \,\square\, (W \otimes M)) = J(W \triangleright B \bullet B \,\square\, M).$$

Thus, if $J$ defines a functor, it is an equivalence of $({}^{B}\mathsf{Vect})$-module categories.

It is left to verify the functoriality of $J$. Let

$$\sigma \in (B \bullet B)\text{-ModFree}((B \bullet B) \square M, (B \bullet B) \square N),$$
$$\tau \in (B \bullet B)\text{-ModFree}((B \bullet B) \square N, (B \bullet B) \square W).$$

We represent these morphisms graphically by

| | |
|---|---|
| 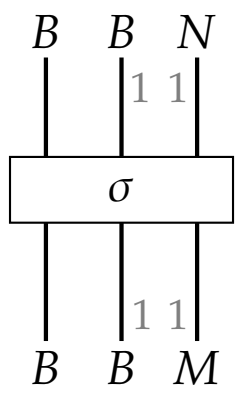  | 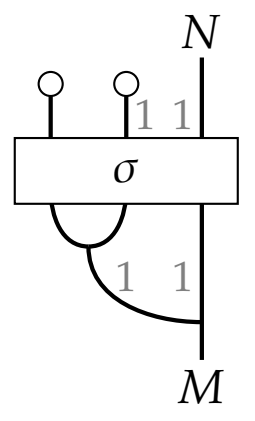  |
| $\sigma \in (B \bullet B)\text{-ModFree}((B \bullet B) \square M, (B \bullet B) \square N)$ | $J(\sigma) \in \mathrm{Hom}_{\Bbbk}(M, N)$ |

Functoriality of $J$ follows by Figure 9.11; (i) and (ii) hold due to

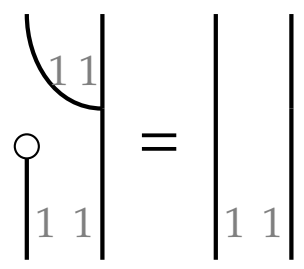

Equality (iii) is satisfied by $\tau$ being a $(B \bullet B)$-module morphism, (iv) follows by the following calculation, using that $B$ is a bialgebra:

$$\varepsilon(b_{(1)}b_{(2)}) \otimes b_{(3)} = \varepsilon(b_{(1)})\varepsilon(b_{(2)}) \otimes b_{(3)} = \varepsilon(b_{(1)}\varepsilon(b_{(2)})) \otimes b_{(3)} = \varepsilon(b_{(1)}) \otimes b_{(2)} = b,$$

and (v) holds by an analogous computation. □

**Corollary 9.33.** *The monad* $(B \bullet B) \square -$ *is semisimple, and there is an equivalence of* ${}^{B}\mathsf{Vect}$*-module categories* $(B \bullet B)\text{-Mod} \simeq \mathsf{Vect}$.

## 9.6 THE FUNDAMENTAL THEOREM OF HOPF MODULES

**Definition 9.34** ([Mon93, Proposition 1.5.10])**.** Let $B$ be a bialgebra. A *twisted antipode* for $B$ is an antipode for the bialgebra $B^{\mathrm{cop}}$.[22]

[22] Equivalently, a twisted antipode is an antipode for $B^{\mathrm{op}}$.

**Remark 9.35.** A Hopf algebra admits a twisted antipode $\tilde{S}$ if and only if its own antipode $S$ is invertible, in which case we have $S^{-1} = \tilde{S}$. Notably, a finite-dimensional Hopf algebra always admits a twisted antipode.

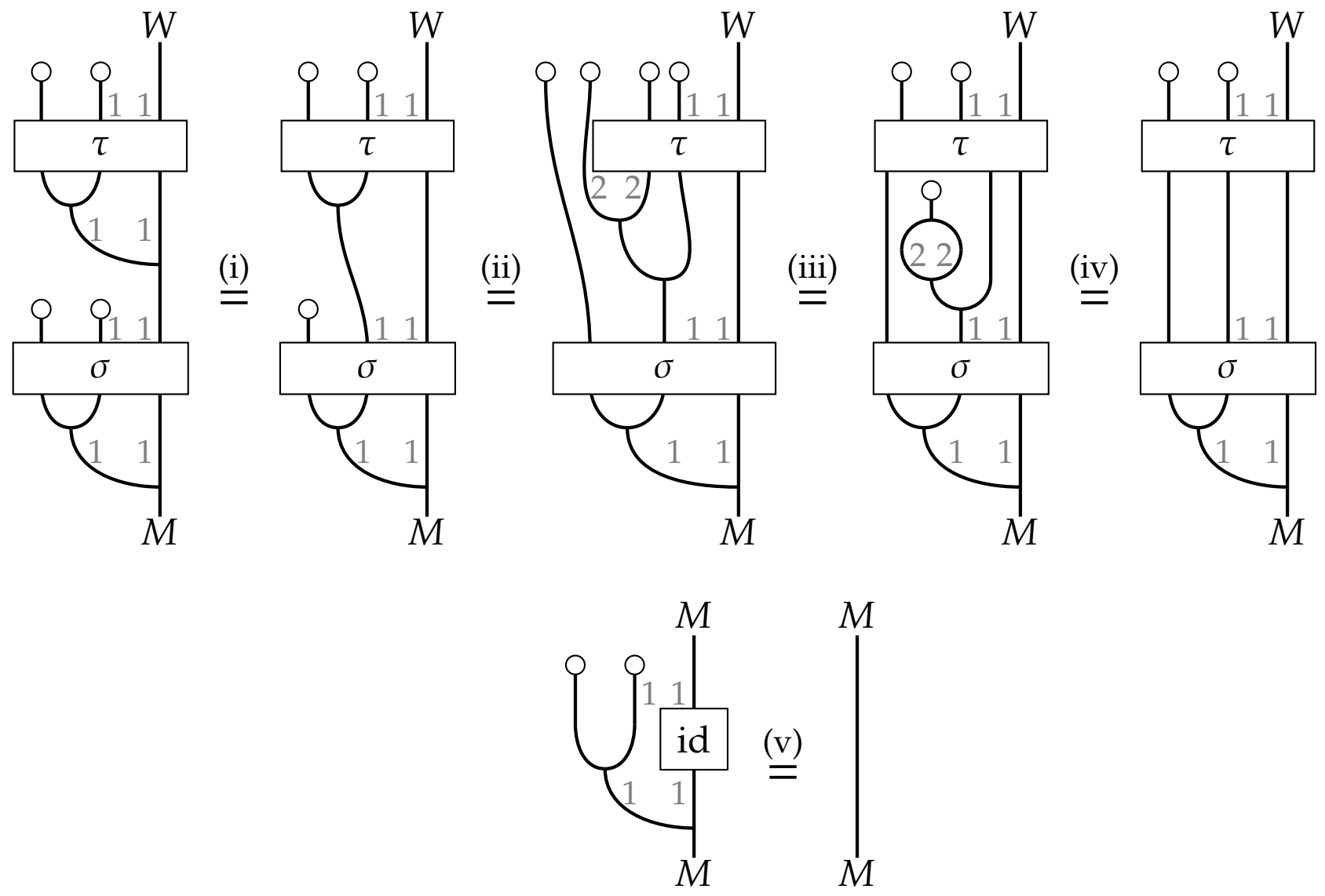

Figure 9.11: The assignment $J$ is functorial.

Lemma 9.36 below is a well-known result, see for example [EGNO15, Theorem 5.4.1.(2)]. We emphasise that while its former part is seemingly opposite to the claim of [EGNO15, Theorem 5.4.1.(2)], identifying the presence of an antipode with right, rather than left, rigidity, this is because we consider left, rather than right, comodules. The latter claim of Lemma 9.36 follows from the former by using the equivalence $({}^{B}\mathsf{Vect})^{\mathrm{rev}} \simeq {}^{B^{\mathrm{op}}}\mathsf{Vect}$.

**Lemma 9.36.** *Let $B$ be a bialgebra. Then ${}^{B}\mathsf{vect}$ is right rigid if and only if $B$ admits an antipode. Similarly, ${}^{B}\mathsf{vect}$ is left rigid if and only if $B$ admits a twisted antipode.*

The next theorem was shown in the setting of quasi-bialgebras in [HN99, Proposition 3.11] and [Sar17, Theorem 3.9], as a generalisation of the structure theorem for Hopf modules of Larson and Sweedler [LS69].

**Theorem 9.37.** *Let $B$ be a bialgebra. The functor $- \otimes B \colon {}_{B}\mathsf{Vect} \longrightarrow {}_{B}\mathsf{Vect}^{B}_{B}$ is an equivalence if and only if $B$ admits an antipode.*

In the rest of this section we show that a variant of Theorem 9.37 can be deduced from the following general statement.

**Proposition 9.38.** *Let $\mathscr{C}$ be a left closed monoidal category such that every left $\mathscr{C}$-module endofunctor of $\mathscr{C}$ is strong—in other words,that the monoidal embedding*

$$\mathsf{Str}\mathscr{C}\mathsf{Mod}(\mathscr{C}, \mathscr{C}) \hookrightarrow \mathsf{Lax}\mathscr{C}\mathsf{Mod}(\mathscr{C}, \mathscr{C})$$

*is an equivalence. Then $\mathscr{C}$ is left rigid.*

*Proof.* Since $\mathscr{C}$ is left closed, for every $x \in \mathscr{C}$ there is a right adjoint $\lfloor x, - \rfloor$ to the strong $\mathscr{C}$-module endofunctor $- \otimes x$. By Porism 5.29, the functor $\lfloor x, - \rfloor$ is a lax $\mathscr{C}$-module endofunctor of $\mathscr{C}$, and thus, by assumption, it is a strong module endofunctor. Proposition 2.51 yields an object ${}^{\vee}x$ such that there is an isomorphism of $\mathscr{C}$-module functors $\lfloor x, - \rfloor \cong - \otimes {}^{\vee}x$. Thus we obtain an adjunction of $\mathscr{C}$-module functors $- \otimes x \dashv - \otimes {}^{\vee}x$, which by [Kel72, pp. 102–103] shows that $({}^{\vee}x, x)$ is a dual pair in $\mathscr{C}$. It follows that $\mathscr{C}$ is left rigid. □

**Remark 9.39.** Given a rigid monoidal category $\mathscr{C}$, notice that every lax $\mathscr{C}$-module endofunctor on $\mathscr{C}$ is automatically strong by Proposition 2.73. Together with Proposition 9.38, this yields a characterisation when a (left) closed monoidal category is (left) rigid, solely based on its lax $\mathscr{C}$-module functors. This is a special property of module functors, as in Section 3.1.2 we have already seen that ordinary adjunctions do not suffice for such a characterisation.

We now define the functor $B \otimes -$ analogous to $- \otimes B$ of Theorem 9.37.

**Definition 9.40.** Define the functor $B \otimes - \colon {}^{B}\mathsf{Vect} \longrightarrow {}^{B}_{B}\mathsf{Vect}^{B}$ by sending a left $B$-comodule $M$ to the trimodule $B \otimes M$ whose left and right coaction, as well as left action, is given by

$$b \otimes m \longmapsto b_{(1)} m_{(-1)} \otimes b_{(2)} \otimes m_{(0)}; \quad b \otimes m \longmapsto b_{(1)} \otimes m \otimes b_{(2)}; \quad a \otimes b \otimes m \longmapsto ab \otimes m;$$

where we extends to morphisms in the natural way. It is easy to verify that $B \otimes M$ defines a Hopf trimodule and also that $B \otimes -$ defines a functor.

**Remark 9.41.** Comparing Definition 9.40 with Theorem 9.37, we remark on two differences. First, the domain of the functor we study is the category of comodules over $B$, rather than modules. This is merely to maintain consistence with Section 9.2, and a variant for modules can be formulated without greater difficulty. Second, the functor we study endows a left comodule with an additional comodule structure on the right, and a module structure on the left, rather than on the right. We will see that it is related to the functor $\mathsf{Vect} \longrightarrow {}_{B}\mathsf{Vect}^{B}$, which is known to be an equivalence if $B$ admits a twisted antipode; see for example [Mon93, Section 1.9.4].

We now give a categorical interpretation of the functor of Definition 9.40. Recall the notation $\mathscr{V} := {}^{B}\mathsf{Vect}$ of Section 9.2.

**Proposition 9.42.** *The diagram*

$$\begin{array}{ccc} \mathsf{Str}\mathcal{V}\mathsf{Mod}(\mathcal{V},\mathcal{V}) & \hookrightarrow & \mathsf{LexfLax}\mathcal{V}\mathsf{Mod}(\mathcal{V},\mathcal{V}) \\ \downarrow \simeq & & \downarrow \simeq \\ \mathcal{V}^{\mathrm{rev}} & \hookrightarrow & {}^{B}_{B}\mathsf{Vect}^{B} \end{array}$$

*commutes up to a monoidal natural isomorphism. Its upper horizontal arrow is the natural inclusion functor, its left vertical arrow is the equivalence of Proposition 2.51, its right vertical arrows is the equivalence of Theorem 9.2, and its lower horizontal arrow is the functor of Definition 9.40.*

*Proof.* Chasing a functor $F$ in the diagram, we find the following:

$$\begin{array}{ccccccc} F & \longmapsto & & & & & F \\ \downarrow & & & & & & \downarrow \\ F\Bbbk_{\mathrm{triv}} & \longmapsto & B \otimes F\Bbbk_{\mathrm{triv}} & \xrightarrow{\simeq} & F(B \otimes \Bbbk_{\mathrm{triv}}) & \xrightarrow{\simeq} & FB \end{array}$$

It is easy to verify that the isomorphism indicated in the bottom row of this chase is natural in $F$, and also monoidal. □

**Corollary 9.43.** *If the functor $B \otimes -$ of Definition 9.40 is an equivalence, then $B$ admits a twisted antipode.*

*Proof.* By Proposition 9.42, $B\otimes-$ is an equivalence if and only if the embedding

$$\mathsf{Str}\mathcal{V}\mathsf{Mod}(\mathcal{V},\mathcal{V}) \hookrightarrow \mathsf{LexfLax}\mathcal{V}\mathsf{Mod}(\mathcal{V},\mathcal{V})$$

is an equivalence, so assume that it is.

For $M \in {}^{B}\mathsf{vect}$, the right adjoint of $- \otimes_{\Bbbk} M$ is finitary by Proposition 2.134, since $- \otimes_{\Bbbk} M$ sends finite-dimensional comodules to finite-dimensional comodules. The right adjoint is then a finitary left exact lax $\mathcal{V}$-module functor. By assumption it is strong, and thus of the form $- \otimes {}^{\vee}M$, for ${}^{\vee}M \in \mathcal{V}$. Since the fibre functor $U\colon \mathcal{V} \longrightarrow \mathsf{Vect}$ is strong monoidal, we have ${}^{\vee}M \cong M^{*}$. This shows that ${}^{B}\mathsf{vect}$ is left rigid; by Lemma 9.36, $B$ admits a twisted antipode. □

**Proposition 9.44.** *If a bialgebra $B$ admits a twisted antipode, then the functor $B\otimes-$ of Definition 9.40 is an equivalence.*

*Proof.* Since the functor $\mathcal{V} \hookrightarrow \mathsf{LexfLax}\mathcal{V}\mathsf{Mod}(\mathcal{V}, \mathcal{V})$ is always fully faithful, so is the functor $B \otimes -$ by Proposition 9.42. It remains to show that $B \otimes -$ is essentially surjective.

Following [Mon93, Section 1.9.4], for $X \in {}_B\mathsf{Vect}^B$ the composite map

$$\Gamma\colon X \xrightarrow{\Delta_X} X \otimes B \xrightarrow{\mathrm{t} \otimes B} X^{\mathrm{co}B} \otimes B \xrightarrow{\sigma_{X^{\mathrm{co}B}, B}} B \otimes X^{\mathrm{co}B}$$

is an isomorphism in ${}_B\mathsf{Vect}^B$, where $X^{\mathrm{co}B} = \{\, m \in X \mid \Delta^r_X(m) = m \otimes 1 \,\}$ and $\mathrm{t}$ is the projection diagrammatically represented by

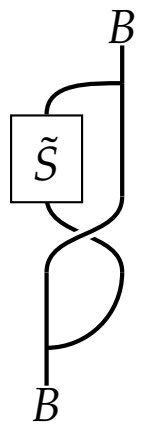

The inverse of $\Gamma$ is given by the restriction $\nabla_X|_{\mathrm{co}B}\colon B \otimes X^{\mathrm{co}B} \longrightarrow X$ of the left $B$-action on $X$.

For $X \in {}^B\mathsf{Vect}_B$, the space $X^{\mathrm{co}B}$ of coinvariants is a left $B$-subcomodule of $X$, since $X$ is a bicomodule. It is easy to check that $\nabla_X|_{\mathrm{co}B}\colon B \otimes X^{\mathrm{co}B} \longrightarrow X$ is a left $B$-comodule morphism, and thus a morphism in ${}^B_B\mathsf{Vect}^B$, making the functor $B \otimes -$ essentially surjective. □

**Corollary 9.45.** *A bialgebra $B$ admits a twisted antipode if and only if the functor $B \otimes -$ of Definition 9.40 is an equivalence.*

**Remark 9.46.** Here we see that it is crucial that Section 9.2 identifies Hopf trimodules with *finitary* lax module functors. Given an infinite-dimensional comodule $M$ over $B$, the functor $\mathrm{Hom}_{\Bbbk}(M, -)\colon {}^B\mathsf{Vect} \longrightarrow {}^B\mathsf{Vect}$ is endowed with lax ${}^B\mathsf{Vect}$-structure by virtue of Theorem 5.28, but is not finitary itself.

## 9.7 Fusion operators for Hopf monads

Similarly to our study of Theorem 9.37 in Section 9.6, the results of Bruguières, Lack, and Virelizier in [BLV11] can be used to characterise rigidity of a monoidal category.

**Lemma 9.47.** *Let $F\colon \mathcal{C} \rightleftarrows \mathcal{D} :U$ be an oplax monoidal adjunction. The strong monoidal structure of $U$ turns $\mathcal{C}$ into a $\mathcal{D}$-module category, by defining*

$$- \triangleright = \;:=\; U(-) \otimes =.$$

A variant of this result using "Doi–Koppinen data" appears in [KKS15, Section 5.5].

*With respect to this $\mathfrak{D}$-module structure, the bimonad $T := UF$ on $\mathscr{C}$ becomes an oplax $\mathfrak{D}$-module monad. The coherence morphism is given by*

$$T_{\mathsf{a};d,c}\colon T(d \triangleright c) = T(Ud \otimes c) \xrightarrow{T_{2;Ud,c}} TUd \otimes Tc \xrightarrow{U\varepsilon_d \otimes Tc} Ud \otimes Tc, \tag{9.7.1}$$

*for all $c \in \mathscr{C}$ and $d \in \mathfrak{D}$. In particular, for $x \in \mathscr{C}$, the coherence morphism $T_{\mathsf{a};Fx,c}$ is precisely the right fusion operator for $T$ at $(x, c)$. In other words, $T_{\mathsf{a};Fx,c} = T_{\mathsf{rf};x,c}$.*

*Proof.* The first statement is a well-known result on transport of structure. For $y \in \mathfrak{D}$, the $\mathfrak{D}$-module structure on $U$ is given by the morphisms

$$U_{\mathsf{m}}\colon d \triangleright Uy = Ud \otimes Uy \xrightarrow{\sim} U(d \otimes y) = U(d \triangleright y).$$

For the second statement, observe that by Porism 5.29, $F$ inherits an oplax $\mathfrak{D}$-module functor structure from the strong $\mathfrak{D}$-module structure on $U$. The oplax $\mathfrak{D}$-module functor on $T$ is given by the upper part of the following diagram, where the lower part is included clarity:

$$\begin{array}{ccccc}
UF(Ud \otimes c) & \xrightarrow{UF_{2;Ud,c}} & U(FUd \otimes Fc) & & \\
\| & & \downarrow{\scriptstyle U(\varepsilon_d \otimes Fc)} & & \\
T(d \triangleright c) & & U(d \otimes Fc) & \xrightarrow{U_{\mathsf{m};d,Fc}} & Ud \otimes UFc \\
\| & & \| & & \| \\
UF(d \triangleright c) & \longrightarrow & U(d \triangleright Fc) & & Ud \otimes Tc \\
 & & \downarrow & & \| \\
 & & d \triangleright UFc & = & d \triangleright Tc
\end{array}$$

To see that this morphism equals that defined in Equation (9.7.1), observe that the following diagram of functors commutes by naturality of $U_2$:

$$\begin{array}{ccc}
U(FU \otimes F) & \xrightarrow{U_{2;FU,F}} & UFU \otimes UF \\
{\scriptstyle U(\varepsilon \otimes F)}\downarrow & & \downarrow{\scriptstyle U\varepsilon \otimes UF} \\
U(\mathrm{Id}_{\mathfrak{D}} \otimes F) & \xrightarrow[U_{2;\mathrm{Id},F}]{} & U \otimes UF
\end{array}$$

Since by definition we have that $\mu = U\varepsilon F$, substituting $d = Fx$ in Equation (9.7.1) yields $T_{\mathsf{a};Fx,c} = T_{\mathsf{rf};x,c}$, which proves the result. □

**Lemma 9.48.** *Let $\mathscr{C}$ be a category, $f, g\colon x \longrightarrow y$ a pair of parallel morphisms in $\mathscr{C}$, such that the coequaliser $\operatorname{coeq}(f, g)$ exists. If $G, H\colon \mathscr{C} \longrightarrow \mathfrak{D}$ is a pair of functors preserving that coequaliser, $\alpha\colon G \Longrightarrow H$ a natural transformation, then $\alpha_x$ and $\alpha_y$ being invertible implies that so is $\alpha_{\operatorname{coeq}(f,g)}$.*

*Proof.* Let $c := \operatorname{coeq}(f, g)$ and write $p\colon y \twoheadrightarrow c$ for the canonical projection. We are in the following situation:

$$\begin{array}{ccccc} Gx & \overset{Gf}{\underset{Gg}{\rightrightarrows}} & Gy & \overset{Gp}{\twoheadrightarrow} & Gc \\ {\scriptstyle\alpha_x}\downarrow & & \downarrow{\scriptstyle\alpha_y} & & \\ Hx & \overset{Hf}{\underset{Hg}{\rightrightarrows}} & Hy & \overset{Hp}{\twoheadrightarrow} & Hc \end{array}$$

Then we calculate

$$Hp \circ \alpha_y \circ Gf = Hp \circ Hf \circ \alpha_x = Hp \circ Hg \circ \alpha_x = Hp \circ Gg \circ \alpha_x,$$

which, by the universal property of $Gc$, implies the existence of a unique arrow $!\colon Gc \longrightarrow Hc$ that is equal to $\alpha_c$. Since $\alpha_x$ and $\alpha_y$ are invertible we also obtain a unique arrow $Hc \longrightarrow Gc$, which is easily seen to be inverse to $\alpha_c$. □

**Proposition 9.49.** *The right fusion operators of bimonad $T$ on a monoidal category $\mathscr{C}$ are invertible if and only if the coherence morphisms for the oplax $\mathscr{C}^T$-module structure on $T$ of Lemma 9.47 are invertible.*

*In other words, $T$ is right Hopf if and only if it is a strong $\mathscr{C}^T$-module monad.*

*Proof.* Since, by Lemma 9.47, the right fusion operators for $T$ are a special case of the coherence morphisms for the $\mathscr{C}^T$-module functor structure on $T$, it is immediate that $T$ being a strong $\mathscr{C}^T$-module monad implies invertibility of the right fusion operators.

To show the converse, assume that the right fusion operators are invertible, in particular, by Lemma 9.47, the morphisms $T_{\mathsf{a};Fx,c}$ are invertible. Recall any $T$-algebra $(x, \alpha) \in \mathscr{C}^T$ is the coequaliser of the diagram

$$T^2x \overset{\mu_x}{\underset{T\alpha}{\rightrightarrows}} Tx \overset{\alpha}{\longrightarrow} x$$

in $\mathscr{C}^T$. In other words, $x$ is a $U$-split coequaliser of morphisms from $FTx$ to $Fx$, since $Fx$ is the free $T$-module on $x$.

Consider the natural transformation $T_{\mathsf{a};-,c} \colon T(- \triangleright c) \longrightarrow - \triangleright Tc$. We have

$$T(- \triangleright c) = T(U(-) \otimes c) = (- \otimes c) \circ U,$$

which means that $T(- \triangleright c)$ preserves the coequaliser $\operatorname{coeq}(T\alpha, \mu_x)$ by first sending it to a split coequaliser under $U$, and then preserving it, since split coequalisers are absolute.

In a similar fashion, $- \triangleright Tc = U(-) \otimes Tc = (- \otimes Tc) \circ U$ also preserves this coequaliser. The invertibility of $T_{\mathsf{a};x,c}$ now follows from the assumed invertibility of $T_{\mathsf{rf};x,c} = T_{\mathsf{a};Fx,c}$ and $T_{\mathsf{rf};Tx,c} = T_{\mathsf{a};FTx,c}$, by Lemma 9.48. □

# LIST OF FIGURES

Hiermit versichere ich, dass ich die vorliegende Arbeit ohne unzulässige Hilfe Dritter und ohne Benutzung anderer als der angegebenen Hilfsmittel angefertigt habe; die aus fremden Quellen direkt oder indirekt übernommenen Gedanken sind als solche kenntlich gemacht. Die Arbeit wurde bisher weder im Inland noch im Ausland in gleicher oder ähnlicher Form einer anderen Prüfungsbehörde vorgelegt.

Tony Zorman
Dresden, 01. April 2025